\pdfoutput=1 
\documentclass[a4paper]{amsbook}

\newcommand{\uref}[1]{{\upshape\eqref{#1}}} 
\newcommand{\lref}[2]{\ref{#2}}

\newcommand{\chnotconv}{Notation and conventions}

\numberwithin{equation}{chapter}

\newcounter{ctr} \numberwithin{ctr}{section}

\usepackage[pdfpagelabels, breaklinks=true]{hyperref}

\usepackage{amssymb}

\usepackage{graphicx}

\usepackage[all]{xy} 
\usepackage[mathscr]{eucal} 
\usepackage{bbm} 
\usepackage{mathtools} 

\usepackage{amsmidx}

\usepackage{etoolbox}

\usepackage{tikz}
\usetikzlibrary{arrows,positioning}

\newcommand{\breakflow}{\medbreak}
\newcommand{\aftereq}{\medbreak}
\newcommand{\afterenum}{\medbreak}
\newcommand{\afterall}{\medbreak}

\newcommand{\thmpre}[4]{#1\noindent\refstepcounter{ctr}{#3#2 \thectr.}#4}
\newcommand{\thmpreit}[2]{\thmpre{#1}{#2}{\bf}{ \itshape}}
\newcommand{\thmpreup}[2]{\thmpre{#1}{#2}{\bf}{ }}
\newcommand{\thmpost}{\ignorespacesafterend}

\newcommand{\athmpre}[4]{#1\noindent{#3#2.}#4}

\newenvironment{thm}[ 1][\afterall]{\thmpreit{#1}{Theorem}}{\thmpost}
\newenvironment{lem}[ 1][\afterall]{\thmpreit{#1}{Lemma}}{\thmpost}
\newenvironment{prp}[ 1][\afterall]{\thmpreit{#1}{Proposition}}{\thmpost}
\newenvironment{cor}[ 1][\afterall]{\thmpreit{#1}{Corollary}}{\thmpost}

\newenvironment{dfn}[ 1][\afterall]{\thmpreup{#1}{Definition}}{\thmpost}
\newenvironment{nrmk}[ 1][\afterall]{\thmpre{#1}{Remark}{\bf}{ }}{\thmpost}
\newenvironment{example}[ 1][\afterall]{\thmpre{#1}{Example}{\bf}{ }}{\thmpost}
\newenvironment{rmk}[ 1][\afterall]{\breakflow}{\thmpost}

\newenvironment{dfn*}[ 1][\afterall]{\athmpre{#1}{Definition}{\bf}{ }}{\thmpost}
\newenvironment{example*}[ 1][\afterall]{\athmpre{#1}{Example}{\bf}{ }}{\thmpost}
\newenvironment{nrmk*}[ 1][\afterall]{\athmpre{#1}{Remark}{\bf}{ }}{\thmpost}
\newenvironment{cor*}[ 1][\afterall]{\athmpre{#1}{Corollary}{\bf}{ \itshape}}{\thmpost}

\newcommand{\pf}[ 1][\medbreak]{#1\noindent\textit{Proof}. }
\newcommand{\quod}{\qed}

\makeatletter
\DeclareRobustCommand\widecheck[1]{{\mathpalette\@widecheck{#1}}}
\def\@widecheck#1#2{%
    \setbox\z@\hbox{\m@th$#1#2$}%
    \setbox\tw@\hbox{\m@th$#1%
       \widehat{%
          \vrule\@width\z@\@height\ht\z@
          \vrule\@height\z@\@width\wd\z@}$}%
    \dp\tw@-\ht\z@
    \@tempdima\ht\z@ \advance\@tempdima2\ht\tw@ \divide\@tempdima\thr@@
    \setbox\tw@\hbox{%
       \raise\@tempdima\hbox{\scalebox{1}[-1]{\lower\@tempdima\box
\tw@}}}%
    {\ooalign{\box\tw@ \cr \box\z@}}}
\makeatother

\newcommand{\bF}{\mathbb{F}}
\newcommand{\bG}{\mathbb{G}}
\newcommand{\bR}{\mathbb{R}}
\newcommand{\bZ}{\mathbb{Z}}

\newcommand{\bU}{\mathbb{U}}

\newcommand{\cC}{\mathcal{C}}
\newcommand{\cE}{\mathcal{E}}
\newcommand{\cF}{\mathcal{F}}
\newcommand{\cG}{\mathcal{G}}
\newcommand{\cI}{\mathcal{I}}
\newcommand{\cJ}{\mathcal{J}}
\newcommand{\cM}{\mathcal{M}}
\newcommand{\cN}{\mathcal{N}}
\newcommand{\cO}{\mathcal{O}}

\newcommand{\cQ}{\mathcal{Q}}

\newcommand{\cU}{\mathcal{U}}
\newcommand{\cV}{\mathcal{V}}

\newcommand{\fm}{\mathfrak{m}}
\newcommand{\fp}{\mathfrak{p}}
\newcommand{\fq}{\mathfrak{q}}

\newcommand{\Fq}{{\mathbb{F}_{\!q}}}
\newcommand{\Fp}{{\mathbb{F}_{\!p}}}

\newcommand{\Gm}{{\mathbb{G}_{\textup{m}}}}

\newcommand{\unit}{\mathbbm{1}} 
\newcommand{\unic}{\bU} 

\newcommand{\nC}{A}
\newcommand{\nCalt}{B}

\newcommand{\nmC}{M}
\newcommand{\nmCalt}{N}

\newcommand{\nuC}{z}
\newcommand{\nuCalt}{\zeta}

\newcommand{\neC}{a}
\newcommand{\nmeC}{m}
\newcommand{\nmeCalt}{n}

\newcommand{\ngC}{\Lambda}
\newcommand{\ngS}{S} 

\newcommand{\ngB}{B}

\newcommand{\ngCalt}{{\Lambda'}}
\newcommand{\ngSalt}{S'} 

\newcommand{\ngOF}{\cO_F}
\newcommand{\ngF}{{F}}

\newcommand{\nsOF}{D}
\newcommand{\nsF}{D^\circ}

\newcommand{\nsOFX}{\nsOF\times X}
\newcommand{\nsFX}{\nsF\times X}

\newcommand{\ibase}{\mathfrak{a}}
\newcommand{\icoef}{\mathfrak{b}}

\newcommand{\nfA}{F^\natural}
\newcommand{\nOA}{\omega^\natural}

\newcommand{\nSpA}{C^\circ}

\DeclareMathOperator{\Res}{Res}

\newcommand{\fmF}{\mathfrak{m}}
\newcommand{\rami}{\mathfrak{e}}

\newcommand{\isosign}{\smash{\raisebox{-0.65ex}{\ensuremath{\sim}}}} 
\newcommand{\tisosign}{\scalebox{0.9}{\ensuremath{\sim}}} 
\newcommand{\bisosign}{\scalebox{0.9}[0.9]{\ensuremath{\sim}}} 
\newcommand{\visosign}{\scalebox{1.0}{\ensuremath{\wr}}} 
\newcommand{\ltviso}{\visosign\!}
\newcommand{\rtviso}{\!\visosign}

\DeclareMathOperator*{\colim}{colim}
\DeclareMathOperator{\cone}{cone}
\DeclareMathOperator{\coker}{coker}

\DeclareMathOperator{\rank}{rank}
\DeclareMathOperator{\tr}{tr}

\DeclareMathOperator{\Aut}{Aut}
\DeclareMathOperator{\End}{End}

\DeclareMathOperator{\Hom}{Hom}
\DeclareMathOperator{\Spec}{Spec}
\DeclareMathOperator{\Frac}{Frac}

\DeclareMathOperator{\Ext}{Ext}
\DeclareMathOperator{\Tor}{Tor}
\DeclareMathOperator{\RHom}{RHom}
\DeclareMathOperator{\Rlim}{Rlim}
\DeclareMathOperator{\RGamma}{R\Gamma}
\DeclareMathOperator{\RGammac}{R\Gamma_{\!c}}
\DeclareMathOperator{\RGammag}{R\Gamma_{\!\mathit{g}}}
\DeclareMathOperator{\RvGamma}{R\widecheck{\Gamma}}
\DeclareMathOperator{\RcGamma}{R\widehat{\Gamma}}
\DeclareMathOperator{\cGamma}{\Gamma_{\!a}}
\DeclareMathOperator{\RH}{RH}

\DeclareMathOperator{\Lie}{Lie}
\DeclareMathOperator{\iHom}{\mathscr{H}\mathrm{om}}
\DeclareMathOperator{\siHom}{\underline{\mathscr{H}\mathrm{om}}}

\DeclareMathOperator{\res}{res}

\DeclareMathOperator{\Sht}{Sht}

\newcommand{\uD}{\mathrm{D}}
\newcommand{\uH}{\mathrm{H}}
\newcommand{\uK}{\mathrm{K}}
\newcommand{\uR}{\mathrm{R}}
\newcommand*\Der{\mathop{}\!\mathbin\nabla}

\newcommand{\complot}{\mathbin{\widehat{\otimes}}}
\newcommand{\indcot}{\mathbin{\widecheck{\otimes}}}

\makeatletter
\newcommand*{\relrelbarsep}{.386ex}
\newcommand*{\relrelbar}{%
  \mathrel{%
    \mathpalette\@relrelbar\relrelbarsep
  }%
}
\newcommand*{\@relrelbar}[2]{%
  \raise#2\hbox to 0pt{$\m@th#1\relbar$\hss}%
  \lower#2\hbox{$\m@th#1\relbar$}%
}
\providecommand*{\rightrightarrowsfill@}{%
  \arrowfill@\relrelbar\relrelbar\rightrightarrows
}
\providecommand*{\shtuka}[2]{%
  \ext@arrow 0359\rightrightarrowsfill@{#2}{#1}%
}
\makeatother

\newcommand{\stacks}[1]{\href{http://stacks.math.columbia.edu/tag/#1}{#1}}

\makeatletter
\def\@idxitem{\par\addvspace{5\p@ \@plus 5\p@ \@minus 3\p@}\hangindent 40\p@}
\def\subitem{\par\hangindent 40\p@ \hspace*{20\p@}}
\def\subsubitem{\par\hangindent 40\p@ \hspace*{30\p@}}

\patchcmd\theindex{\indexname}{\indexname}{}{}
\makeatother

\makeindex{idx}
\makeindex{nidx}

\begin{document}
\title{Shtuka cohomology and \\ special values of Goss $L$-functions}
\author{M. Mornev}
\address{ETH Z\"urich -- D-MATH, R\"amistrasse 101, 8092 Z\"urich, Switzerland}
\email{maxim.mornev@math.ethz.ch} 
\frontmatter
\maketitle
\tableofcontents
\mainmatter

\chapter*{Introduction}
\setcounter{section}{0}

Assuming everywhere good reduction
we generalize the class number formula of Taelman \cite{ltsv} to Drinfeld
modules over arbitrary coefficient rings. In order to prove this formula we
develop a theory of shtukas and their cohomology.

\section{A class number formula for Drinfeld modules}

Fix a finite field $\Fq$. In the following all morphisms, fiber and tensor
products will be over $\Fq$ unless indicated otherwise.
Let $C$ be a smooth projective connected curve over $\Fq$. Fix a closed point
$\infty \in C$. The $\Fq$-algebra
\begin{equation*}
A = \Gamma(C-\{\infty\}, \cO_C)
\end{equation*}
will be called the \emph{coefficient ring}.
%
Fix an $A$-algebra $R$ (the base ring). We denote $\iota\colon A \to R$
the natural map.



Consider the group scheme $\bG_a$ over $R$. It is well-known that every
$\Fq$-linear endomorphism of $\bG_a$ can be uniquely written in the form
of a \emph{$\tau$-polynomial}
\begin{equation*}
r_0 + r_1 \tau + \dotsc + r_n \tau^n
\end{equation*}
where $r_0, \dotsc, r_n \in R$, $r_n \ne 0$ and $\tau\colon \bG_a \to \bG_a$
denotes the $q$-Frobenius.

Recall that an $A$-module scheme is an abelian group scheme equipped with an
action of $A$.

\begin{dfn*}
A Drinfeld module over $R$ with coefficients in $A$ is an $A$-module scheme $E$
over $R$ which has the following properties:
\begin{enumerate}
\item The underlying
additive group
scheme of $E$ is Zariski-locally
isomorphic to $\bG_a$.


\item For every element $a \in A$ the induced endomorphism of the Lie algebra scheme
$\Lie_E$ is the multiplication by $\iota(a)$.

\item There exists an element $a \in A$ such that locally on $\Spec R$
the action of $a$ on $E$ is given by a $\tau$-polynomial of positive degree
and with top coefficient a unit.
\end{enumerate}
\end{dfn*}

\begin{example*}
Let $A = \Fq[t]$, $R = \Fq[\theta]$ and let $\iota\colon A \to R$ be the
isomorphism which sends $t$ to $\theta$. An example of a Drinfeld
$A$-module over $R$ is the \emph{Carlitz module $E$}. Its underlying additive group
scheme is $\bG_a$. The action of $t \in A$ on $E$ is given by the $\tau$-polynomial
\begin{equation*}
\theta + \tau.
\end{equation*}
The Frobenius $\tau$ induces the zero endomorphism on $\Lie_{E}$.
Hence $t$ acts on $\Lie_{E}$ as the multiplication by $\theta = \iota(t)$.
It follows that
the condition (2) holds for $E$. The conditions (1) and (3) are clear.
\end{example*}

\breakflow
From now on we assume that $R$ is a domain and that it is finite flat over $A$.
It follows that $R$ is a Dedekind domain of finite type over $\Fq$.
We denote $K$ the fraction
field of $R$. 
%
The generic fiber of a Drinfeld module over $R$ is a Drinfeld module over $K$.
However, not every Drinfeld module over $K$ extends to a Drinfeld module over
$R$. The ones which do are said to have \emph{good reduction everywhere}.

Drinfeld modules 
behave in a way similar to elliptic curves. For the latter the
role of the coefficient ring $A$ is played by $\bZ$. Given an elliptic curve $E$
over a number field
and a prime $(p) \subset \bZ$ one can consider its $p$-adic Tate module
\begin{equation*}
T_p E = \Hom_{\bZ}(\mathbb{Q}_p/\mathbb{Z}_p, E(\overline{\mathbb{Q}})).
\end{equation*}
Much in the same way for a Drinfeld module $E$ over the function field $K$ and a
prime $\fp \subset A$ one has the $\fp$-adic Tate module
\begin{equation*}
T_\fp E = \Hom_A(F_\fp/A_\fp, E(K^{\textup{sep}}))
\end{equation*}
where $K^{\textup{sep}}$ denotes a separable closure of $K$,
$A_\fp$ is the completion of $A$ at $\fp$ and $F_\fp$ is the field
of fractions of $A_\fp$.

The Tate module $T_\fp E$ is a finitely generated free
$A_\fp$-module. Its rank can be any integer greater than zero.
This rank does not depend on $\fp$ and is called the \emph{rank of $E$}.
The Tate module
$T_\fp E$ is naturally a continuous representation of the Galois group
$G(K^{\textrm{sep}}/K)$ 
unramified at almost all primes $\fm \subset R$.
%
%
 
Let us fix a Drinfeld module $E$ over $R$.
%
Let $\fp \subset A$ and $\fm \subset R$ be primes such that $\fp \ne
\iota^{-1}(\fm)$.
Since $E$ is defined over $R$
the Tate module $T_\fp E$ is unramified at $\fm$.
It thus makes sense to consider the inverse characteristic polynomial of
the geometric Frobenius element at $\fm$ acting on $T_\fp E$. This polynomial
has coefficients in the fraction field of $A$ and is independent of the
choice of $\fp$. We denote it $P_\fm(T)$.

\begin{dfn*} Let $F_\infty$ be the local field of the curve $C$ at $\infty$.
We define $L(E^*,0) \in F_\infty$ by the formula
\begin{equation*}
L(E^*, 0) = \prod_{\fm} \frac{1}{P_\fm(1)}
\end{equation*}
where the product ranges
over all primes $\fm \subset R$.
\end{dfn*}

\breakflow
It is not difficult to show that this product converges.
The resulting element
$L(E^*,0) \in F_\infty$ is indeed a value of a certain function,
the Goss $L$-function of the strictly compatible family of Galois representations given
by the Tate modules $T_\fp E$.
We use the notation $L(E^*,0)$ instead of $L(E,0)$ since the usual Goss
$L$-function of $E$ is given by the family of dual Tate modules $(T_\fp E)^*$.
For the experts we remark that our $L(E^*,0)$ coincides with the special value
considered by Taelman in \cite{ltsv}.


\begin{example*}
Let us examine the Carlitz module $E$.
For every prime $\fm \subset R = \Fq[\theta]$ there exists a unique monic irreducible
polynomial $f \in \Fq[t]$ such that
\begin{equation*}
\fm = \iota(f) R.
\end{equation*}
Let $g \in \Fq[t]$ be another monic irreducible polynomial, $g \ne f$.
The Tate module $T_g E$ is free of rank $1$.
One can show that
the \emph{arithmetic} Frobenius at $\fm$ acts on $T_g E$ by multiplication by $f$.
Therefore the inverse charateristic polynomial of the \emph{geometric} Frobenius at
$\fm$ is
\begin{equation*}
P_\fm(T) = 1 - f^{-1} T.
\end{equation*}
We conclude that the special value $L(E^*,0)$ is given by the Euler product
\begin{equation*}
L(E^*, 0) =
\prod_{f}
\frac{1}{1 - \dfrac{1}{f}}
\end{equation*}
where $f \in \Fq[t]$ runs over monic irreducible polynomials.
Expanding this product we get the formula
\begin{equation*}
L(E^*, 0) =
\sum_h
\frac{1}{h}
\end{equation*}
where $h \in \Fq[t]$ runs over all monic polynomials.
\end{example*}

\breakflow
Let $K_\infty = R \otimes_A F_\infty$.
%
The exponential map of the Drinfeld module $E$ is the unique map
\begin{equation*}
\exp\colon
\Lie_E(K_\infty) \to E(K_\infty)
\end{equation*}
satisfying the
following conditions:
\begin{enumerate}
\item $\exp$ is a homomorphism of $A$-modules,

\item $\exp$ is an analytic function with derivative $1$ at zero in the following sense.
Fix an $\Fq$-linear isomorphism of group schemes $E \cong \bG_a$ defined over
$K_\infty$.
It
identifies $E(K_\infty)$ with $K_\infty$ while its
differential identifies $\Lie_E(K_\infty)$ with $K_\infty$.
We demand that the resulting map $\exp\colon K_\infty \to K_\infty$
is given by an everywhere convergent power series
\begin{equation*}
\exp(z) = z + a_1 z^q + a_2 z^{q^2} + \dotsc
\end{equation*}
with coefficients in $K_\infty$.
\end{enumerate}

\begin{dfn*} The \emph{complex of units} of $E$ is the $A$-module complex
\begin{equation*}
U_E =
\Big[ \Lie_E(K_\infty) \xrightarrow{\quad\exp\quad}
\frac{E(K_\infty)}{E(R)} \Big]
\end{equation*}
where $\Lie_E(K_\infty)$ is placed in degree $0$.\end{dfn*}

\breakflow
By construction
\begin{align*}
\uH^0(U_E) &= \exp^{-1} E(R), \\
\uH^1(U_E) &= \frac{E(K_\infty)}{\exp(K_\infty) + E(R)}.
\end{align*}
The cohomology modules of $U_E$ have individual names:
$\uH^0(U_E)$ is the module of units and $\uH^1(U_E)$ is the class module.

\begin{dfn*}
The \emph{regulator}
\begin{equation*}
\rho\colon F_\infty \otimes_A U_E \to \Lie_E(K_\infty)[0]
\end{equation*}
is the $F_\infty$-linear extension of the morphism $U_E \to \Lie_E(K_\infty)[0]$
given by the identity in degree zero.
\end{dfn*}

\afterall\noindent\refstepcounter{ctr}\textbf{Theorem \thectr}\label{lennyut}
(Taelman \cite{ltut}).
\textit{The $A$-module complex $U_E$
is perfect and the regulator $\rho$ is a quasi-isomorphism.}

\breakflow
%
%
%
%
%
%
This theorem is usually stated in a different form:
$\uH^1(U_E)$ is finite and $\uH^0(U_E)$ is a lattice in
$\Lie_E(K_\infty)$.
Recall that an $A$-submodule $\Lambda$ in a finite-dimensional $F_\infty$-vector
space $V$ is called a \emph{lattice} if one of the following equivalent conditions is
satisfied:
\begin{itemize}
\item $\Lambda$ is discrete and cocompact.

\item The natural map $F_\infty\otimes_A \Lambda \to V$ is an isomorphism.
\end{itemize}
A lattice is automatically a finitely generated projective $A$-module.

One may interpret Taelman's theorem
as saying that the complex 
$U_E$ is a lattice in $\Lie_E(K_\infty)$ in a
derived sense.
$\Lie_E(K_\infty)$ contains one more natural lattice: the
integral Lie algebra $\Lie_E(R)$. 
We would like
to determine 
their relative position.

Since $U_E$ is perfect
the theory of Knudsen-Mumford determinants \cite{kmdet}
provides us with
an invertible $A$-module
$\det_A U_E$
and an isomorphism of one-dimensional $F_\infty$-vector spaces
\begin{equation*}
\det\nolimits_{F_\infty}(\rho)\colon
F_\infty \otimes_A \det\nolimits_A U_E \xrightarrow{\isosign}
\det\nolimits_{F_\infty} \Lie_E(K_\infty).
\end{equation*}
The 
vector space
$\det\nolimits_{F_\infty} \Lie_E(K_\infty) = F_\infty \otimes_A \det\nolimits_A \Lie_E(R)$
contains a canonical lattice $\det_A \Lie_E(R)$. 
Now we are ready to state our main result.
\begin{thm}\label{mainthm}%
The image of $\det_A U_E$ under $\det_{F_\infty}(\rho)$ is
\begin{equation*}
L(E^*,0) \cdot \det\nolimits_A \Lie_E(R).
\end{equation*}%
\end{thm}

\begin{nrmk*}
Theorem \ref{mainthm} implies that the 
$A$-modules $\det_A U_E$ and $\det_A \Lie_E(R)$ are isomorphic.
This is by no means immediate if the class group
of the coefficient ring $A$ is not zero. In fact, it was not
previously known apart from the trivial case $\mathrm{Pic}\,\, A = 0$.\end{nrmk*}

\begin{nrmk*}
It is important to realize that Theorem \ref{mainthm} gives
a formula for $L(E^*,0)$ in terms of the lattices $U_E$ and $\Lie_E(R)$.
A~priori the relation
\begin{equation*}
\det\nolimits_{F_\infty}(\rho)\Big(\det\nolimits_A U_E\Big) = x \cdot \det\nolimits_A \Lie_E(R)
\end{equation*}
determines $x \in F_\infty$ up to a unit of $A$. However one can prove
that $L(E^*,0)$ is a $1$-unit in $F_\infty$. 
The only unit of $A$ which is also a $1$-unit of $F_\infty$ is the element $1$.
So the relation above determines a $1$-unit $x$ uniquely.
\end{nrmk*}

\begin{nrmk*}
The statement of Theorem \ref{mainthm} goes back to the fundamental work of
Taelman \cite{ltsv} where he established a formula for
$L(E^*,0)$ under assumption that the coefficient ring $A$ is $\Fq[t]$.
Unlike
our Theorem \ref{mainthm} the result of Taelman applies to Drinfeld modules with
arbitrary reduction type.

Fang \cite{fang} extended the result of Taelman to Anderson modules \cite{and}
which are a higher-dimensional generalization of Drinfeld modules. 
He also considered the coefficient ring $A = \Fq[t]$ only.

Debry \cite{debry} was the first to generalize this formula to coefficient rings
$A$ different from $\Fq[t]$. His result applies to coefficient
rings with trivial class group.
In our Theorem \ref{mainthm} the coefficient ring $A$ can be arbitrary but
the Drinfeld module $E$ is assumed to have good reduction everywhere.
\end{nrmk*}

\begin{nrmk*} One can describe the image of $\det_A U_E$ under
$\det_{F_\infty}(\rho)$
as follows. Theorem \ref{lennyut} implies that the
$A$-submodule
$\exp^{-1} E(R)$ of $\Lie_E(K_\infty)$
is a lattice so that its top exterior power
$\det\nolimits_A\exp^{-1} E(R)$
is an invertible $A$-submodule in the determinant of $\Lie_E(K_\infty)$.
The image of $\det_A U_E$ 
is the $A$-submodule 
\begin{equation*}
I \cdot \det\nolimits_A \exp^{-1} E(R)
\end{equation*}
where $I$ is the $0$-th Fitting ideal of the class module $\uH^1(U_E)$.
This ideal has the following explicit description.
The $A$-module $\uH^1(U_E)$ can be written as a finite direct sum
$\bigoplus_n A/I_n$
where $I_n \subset A$ are ideals. The $0$-th Fitting ideal of $\uH^1(U_E)$ is
\begin{equation*}
I = \prod_n I_n.
\end{equation*}
\end{nrmk*}

\begin{example*} Let us show how Theorem \ref{mainthm} works for the
Carlitz module $E$. In this case $F_\infty = \Fq(\!(t^{-1})\!)$ and $K_\infty =
\Fq(\!(\theta^{-1})\!)$. The exponential map $\exp\colon K_\infty \to K_\infty$
of the Carlitz module admits a local inverse around $0$, the Carlitz logarithm
map. It is given by the power series
\begin{equation*}
\log z = z - \frac{z^q}{\theta^q - \theta} +
\frac{z^{q^2}}{(\theta^q - \theta)(\theta^{q^2} - \theta)} - \dotsc
\end{equation*}
The series converges for $z$ such that $|z| \leqslant q^{\frac{q}{q-1}}$. In
particular the image of the exponential map contains the unit ball
$\Fq[[\theta^{-1}]] \subset K_\infty$. As a consequence the class module
\begin{equation*}
\uH^1(U_E) = \frac{\Fq(\!(\theta^{-1})\!)}{\exp(K_\infty) + \Fq[\theta]}
\end{equation*}
is zero.

The $F_\infty$-vector space $\Lie_E(K_\infty)$ is of dimension $1$. Hence
Theorem \ref{lennyut} implies that $\uH^0(U_E)$ is a free $A$-module of rank $1$.
By construction the element
\begin{equation*}
\widetilde{\pi} = \log(1)
\end{equation*}
belongs to $\uH^0(U_E) = \exp^{-1} E(R)$. A priori it generates an $A$-submodule
of finite index. 
However it is easy to show that a nonconstant element of $A$ can not divide $1 \in E(R)$.
Therefore
\begin{equation*}
\uH^0(U_E) = A \cdot \widetilde{\pi}.
\end{equation*}

The $A$-module
$\uH^0(U_E)$ coincides with its determinant since it is free of rank $1$.
Now $\uH^1(U_E) = 0$ so Theorem \ref{mainthm} implies that
\begin{equation*}
L(E^*,0) = 
\alpha \iota^{-1}(\widetilde{\pi})
\end{equation*}
for some $\alpha \in A^\times$.
Here $\iota\colon F_\infty \cong K_\infty$ is the natural isomorphism.
As observed before, one can show that $L(E^*,0)$
is a $1$-unit of $F_\infty$. Since $\widetilde{\pi}$ is a $1$-unit by
construction it follows that
\begin{equation*}
L(E^*,0) = 
\iota^{-1}(\widetilde{\pi}).
\end{equation*}
Expanding the definitions we obtain a formula
\begin{equation*}
\sum_{\substack{h \in \Fq[t] \\ \textup{monic}}}
\frac{1}{h}
=
1 - \frac{1}{t^q - t} +
\frac{1}{(t^q - t)(t^{q^2} - t)} - \dotsc
\end{equation*}
in $\Fq(\!(t^{-1})\!)$.
Observe that the series on the right hand side converges much faster than the
series
on the left hand side. This formula was discovered by Carlitz \cite{carlitz} in the
1930-ies.
\end{example*}

\begin{nrmk*}
The complex of units $U_E$ has an interesting analog
in the context of number fields.
For the moment, let $K$ be a number field and let $\cO_K \subset K$ be
its ring of integers. Set $K_\infty = \mathbb{R} \otimes_\mathbb{Q} K$ and
consider the complex
\begin{equation*}
U_K = \Big[
\Lie_{\Gm}(K_\infty)
\xrightarrow{\quad\exp\quad}
\frac{\Gm(K_\infty)}{\Gm(\cO_K)} \Big]
\end{equation*}
where $\exp$ is the exponential of the Lie group
$\mathbb{G}_\textup{m}(K_\infty)$.
In this setting the regulator
\begin{equation*}
\rho\colon \mathbb{R}\otimes_\bZ U_K \to
\Lie_{\Gm}(K_\infty)
\end{equation*}
is the $\mathbb{R}$-linear extension of the morphism
$U_K \to \Lie_{\Gm}(K_\infty)[0]$
given by the identity in degree zero.
The Dirichlet's unit theorem for $K$ is equivalent to the following statement:

\afterall\noindent
\textbf{Theorem 1.0}
(Dirichlet). \textit{The $\mathbb{Z}$-module complex $U_K$
is perfect and the cone of $\rho$ is quasi-\hspace{0pt}isomorphic to
$\mathbb{R}[0]$.}\thmpost
\end{nrmk*}

\section{Overview of the proof}

To prepare the ground for the proof of Theorem \ref{mainthm}
we develop a theory of shtukas and their cohomology.
While retaining some 
features of
the works of Taelman \cite{ltsv} and Fang \cite{fang}, our approach differs from them in an essential
way. Certain aspects of this approach were envisaged by Taelman in
\cite{ltcarlitz}. The central idea of using Anderson trace formula \cite{trace}
to study special values of shtukas is due to V.~Lafforgue \cite{valeurs}.
In general, the ideas of Anderson \cite{and,trace} play an important role in
this text. 
Our cohomology theory for shtukas was heavily motivated by the works
of B\"ockle-Pink \cite{bp} and V.~Lafforgue \cite{valeurs}.
The notion of a nilpotent $\tau$-sheaf from \cite{bp} figures prominently in it.
Our intellectual debt to Drinfeld \cite{ell,nc} is obvious.

\begin{nrmk*}
To avoid confusion we should stress that the definitions in this section are
simplified for expository purposes.\end{nrmk*}

\breakflow
We begin with an overview of shtuka theory relevant to the proof of
Theorem~\ref{mainthm}. Let us first describe the setting. 
The finite flat $A$-algebra $R$ is a Dedekind domain of finite type over $\Fq$.
To such an algebra $R$ one can functorially associate a smooth connected
projective curve $X$ over $\Fq$ together with an open embedding $\Spec R \subset X$.
%
Consider the scheme 
$\Spec A \times X$. 
Let $\tau\colon \Spec A \times X \to \Spec A \times X$ be the endomorphism which
acts as the identity on $A$ and as the $q$-Frobenius on~$X$.

\begin{dfn*} A \emph{shtuka} $\cM$ on $\Spec A \times X$ is given by a diagram
\begin{equation*}
\cM_0 \shtuka{i}{j} \cM_1
\end{equation*}
where $\cM_0$, $\cM_1$ are coherent sheaves
on $\Spec A \times X$ and
\begin{align*}
i\colon &\cM_0 \to \cM_1,\\
j\colon &\cM_0 \to \tau_\ast\cM_1
\end{align*}
are morphisms of coherent sheaves. Morphisms of shtukas are given by morphisms
of underlying coherent sheaves which commute with $i$ and $j$ (cf. Definition
\lref{shtukadefs}{dfnshtuka}). 
\end{dfn*}


\breakflow
An example of a shtuka 
is the \emph{unit shtuka}
\begin{equation*}
\unit = \Big[\cO \shtuka{1}{\,\,\tau^\#\,\,} \cO\Big]
\end{equation*}
where $\cO$ is the structure sheaf of $\Spec A \times X$ and $\tau^\#\colon \cO
\to \tau_\ast \cO$ is the map defined by the endomorphism $\tau$.
Shtukas on $\Spec A \times X$ form an abelian category.

\begin{dfn*} Let $\cM$ be a shtuka on $A \times X$.
The \emph{cohomology complex of $\cM$} is the $A$-module complex
\begin{equation*}
\RGamma(\cM) = \RHom(\unit, \cM)
\end{equation*}
where $\RHom$ on the right hand side is computed in the derived category of
shtukas.\end{dfn*}

\breakflow
Shtuka cohomology can be computed in terms of
coherent cohomology: for every shtuka
\begin{equation*}
\cM = \Big[\cM_0 \shtuka{i}{j} \cM_1\Big]
\end{equation*}
there exists a natural distinguished triangle
\begin{equation}\label{introtri}
\RGamma(\cM) \to
\RGamma(\cM_0) \xrightarrow{\,\,i-j\,\,}
\RGamma(\cM_1) \to [1].
\end{equation}
We write $\RGamma(-)$ instead of $\RGamma(\Spec A \times X, -)$ to
improve legibility.

\begin{dfn*} 
We define the
\emph{linearization functor}\index{nidx}{shtukas!$\Der$, linearization}
$\Der$ from the category of shtuka to itself in the following way:
\begin{equation*}
\Der\Big[ \cM_0 \shtuka{i}{j} \cM_1 \Big] =
\Big[ \cM_0 \shtuka{i}{0} \cM_1 \Big].
\end{equation*}
\end{dfn*}

\breakflow
The cohomology of $\Der\cM$ is often easier to compute
than the cohomology of $\cM$. Even though the complexes $\RGamma(\cM)$
and $\RGamma(\Der\cM)$ are usually quite different, there is a
subtle link between them.
%

The sheaves $\cM_0$ and $\cM_1$ are coherent by assumption.
As a consequence
the $A$-module complexes $\RGamma(\cM_0)$ and $\RGamma(\cM_1)$ are
perfect. 
The distinguished triangle \eqref{introtri} now implies that $\RGamma(\cM)$
and $\RGamma(\Der\cM)$ are perfect.
So we can apply the theory of
Knudsen-Mumford determinants to $\RGamma(\cM)$ and $\RGamma(\Der\cM)$.

\begin{dfn*} We define the \emph{$\zeta$-isomorphism}
\begin{equation*}
\zeta_\cM \colon
\det\nolimits_A\RGamma(\cM) \xrightarrow{\isosign}
\det\nolimits_A\RGamma(\Der\cM)
\end{equation*}
as the composition of the isomorphisms
\begin{equation*}
\det\nolimits_A^{\phantom{1}}\RGamma(\cM) \xrightarrow{\isosign}
\det\nolimits_A^{\phantom{1}}\RGamma(\cM_0) \otimes_A^{\phantom{1}}
\det\nolimits_A^{-1}\RGamma(\cM_1) \xleftarrow{\isosign}
\det\nolimits_A^{\phantom{1}}\RGamma(\Der\cM)
\end{equation*}
induced by the natural distinguished triangles \eqref{introtri}
for $\cM$ and $\Der\cM$.%
\end{dfn*}

\breakflow%
The $\zeta$-isomorphisms are named by analogy with the
$\zeta$-elements of Kato \cite{kato}.
Such isomorphisms first appeared in the article \cite{valeurs} of V.~Lafforgue.

Now we are in position to describe the main steps in the proof of Theorem~\ref{mainthm}.
To a Drinfeld module $E$ over $R$ we associate a certain shtuka called
\emph{a~model of $E$}. Its precise definition is a bit technical and is not
necessary to understand the following.
The construction of a model proceeds roughly as follows.
We take the Anderson motive of $E$ and dualize it to obtain a shtuka on $\Spec A
\otimes R$. We then extend it to $\Spec A \times X$ using the functor of
extension by zero from the theory of B\"ockle-Pink \cite{bp}.

Every Drinfeld module $E$ admits many different shtuka models.
However all the models share important properties.
The underlying coherent sheaves of a model are locally free. Their rank
coincides with the rank of $E$. 
The cohomology of a model captures important arithmetic
invariants of $E$.

\begin{thm}\label{introthma} For every shtuka model $\cM$ of $E$ there are natural
quasi-\hspace{0pt}isomorphisms
\begin{equation*}
\RGamma(\cM) \xrightarrow{\isosign} U_E[-1], \quad
\RGamma(\Der\cM) \xrightarrow{\isosign} \Lie_E(R)[-1].
\end{equation*}\end{thm}

%

\breakflow
Recall that
the complex of units $U_E$ is defined in terms of an analytic map,
the exponential of $E$. Theorem \ref{introthma} provides an algebraic
description of this complex.
One important application of it is the following:

\begin{cor*} 
$\det_A U_E \cong
\det_A\Lie_E(R)$.\end{cor*}

\pf Indeed 
we have a $\zeta$-isomorphism
$\zeta_\cM\colon
\det_A\RGamma(\cM) \cong
\det_A\RGamma(\Der\cM)$. \quod

\begin{nrmk*}%
The second quasi-isomorphism in Theorem \ref{introthma} is easy to construct.
In contrast there is no obvious natural map between the complexes $\RGamma(\cM)$
and $U_E[-1]$. The construction of the quasi-isomorphism $\RGamma(\cM) \cong
U_E[-1]$ is rather intricate.
The proof of Theorem \ref{introthma} uses Hochschild cohomology of $A$
as the main computational tool.
\end{nrmk*}

\begin{nrmk*} Taelman \cite{ltcarlitz} established Theorem \ref{introthma} for
the Carlitz module $E$ over an arbitrary finite flat $A$-algebra $R$, $A =
\Fq[t]$. He
constructed shtuka models of $E$ in an ad hoc manner. His result was generalized
by Fang \cite{fang} to Anderson modules with coefficients in $A = \Fq[t]$.
The construction of shtuka models in \cite{fang} is also ad hoc.
\end{nrmk*}

\begin{nrmk*} It is necessary to mention that our proof of Theorem
\ref{introthma} extends without change to arbitrary Anderson $A$-modules,
including the non-uniformizable ones.
In this text we limit the exposition to Drinfeld modules
since other important parts of the theory still depend on their special
properties.\end{nrmk*}

\begin{nrmk*} Our proof of Theorem \ref{introthma} was inspired by the article
\cite{and} of Anderson. In \cite[\S2]{and} he proves a vanishing statement for
$\textup{Ext}^1$ which in retrospect can be viewed as a statement on cohomology
of certain shtukas related to Drinfeld modules.\end{nrmk*}

\breakflow
As we mentioned above the cohomology complexes of a shtuka $\cM$ and its
linearization $\Der\cM$ are quite different in general. The $\zeta$-isomorphism
$\zeta_\cM$
relates their determinants. A more direct link is given by the
$\emph{regulator}$
\begin{equation*}
\rho_\cM\colon F_\infty\otimes_A \RGamma(\cM) \xrightarrow{\isosign} F_\infty \otimes_A
\RGamma(\Der\cM).
\end{equation*}
The regulator is a quasi-isomorphism and is natural in $\cM$.
It is defined for \emph{elliptic shtukas}, 
a natural class of shtukas 
which generalize shtuka models of Drinfeld modules.

The central result about the regulator is the trace formula which expresses
$\zeta_\cM$ in terms of $\rho_\cM$ and an explicit numerical invariant $L(\cM)
\in F_\infty$.
%
This invariant is a product of
local factors, one for each prime $\fm \subset R$. The local factor at $\fm$
depends only on the restriction of $\cM$ to $A \otimes R/\fm$.

\afterall\noindent\refstepcounter{ctr}{\bf Theorem \thectr}\label{introthmb}
(Trace formula)
{\itshape If $\cM$ is an elliptic shtuka satisfying certain
technical condition\footnote{At present we can only prove the theorem under
a technical condition on the cohomology of sheaves underlying $\cM$
(cf. Theorem \lref{globell}{elltrace}).
It is enough for the proof of the class number formula. We expect that the
trace formula holds without this condition.} then}
\begin{equation*}
\zeta_\cM = L(\cM) \cdot \det\nolimits_{F_\infty}(\rho_\cM)
\end{equation*}
{\itshape as maps from $F_\infty \otimes_A \det_A \RGamma(\cM)$ to
$F_\infty \otimes_A \det_A \RGamma(\Der\cM)$.}

\breakflow
The following important lemma is very easy to prove:

\afterall\noindent\textbf{Lemma.}
\textit{If $\cM$ is a shtuka model of $E$ then $L(\cM) = L(E^*,0)$.}


\begin{nrmk*}Theorem \ref{introthmb} is basically the trace formula of Anderson
\cite{trace} in disguise.\end{nrmk*}

\begin{nrmk*}In general the invariant $L(\cM) \in F_\infty$ is transcendental
over $A \subset F_\infty$. Its inherent complexity reflects in the construction of
the regulator making it rather involved. By contrast the definition of
the regulator (Definition \lref{reg}{defellreg}) is simple.\end{nrmk*}


\breakflow
Now we have almost all the tools to prove Theorem \ref{mainthm}.
Thanks to Theorem~\ref{introthma} the
shtuka-theoretic regulator $\rho_\cM$ of a model $\cM$ induces a
quasi-isomorphism
\begin{equation*}
F_\infty \otimes_A U_E \to 
\Lie_E(K_\infty)[0].
\end{equation*}
However there is no a priori reason for it to coincide with
the arithmetic regulator
\begin{equation*}
\rho_E\colon F_\infty \otimes_A U_E \to \Lie_E(K_\infty)[0]
\end{equation*}
which is defined purely in terms of the Drinfeld module $E$. 

\begin{thm}\label{introthmc} Let $\cM$ be a shtuka model of $E$.
The quasi-isomorphisms of Theorem~\ref{introthma} identify the shtuka-theoretic
regulator $\rho_\cM$ with the shifted arithmetic regulator
$\rho_E[-1]$.\end{thm}
\begin{nrmk*}The only proof of Theorem \ref{introthmc} which we have at
the moment is rather technical, and is based on explicit computations
in the case $A = \Fq[t]$.\end{nrmk*}

\breakflow
As we observed above, shtuka models $\cM$ of $E$ exist and have the property
that $L(\cM) = L(E^*,0)$. Hence
Theorems \ref{introthma}, \ref{introthmb} and \ref{introthmc} imply Theorem
\ref{mainthm} for $E$.
 
\begin{nrmk*}
Our theory of shtukas is very sensitive to 
reduction properties of Drinfeld modules. Its extension to the bad reduction
case is not at all straightforward and may be difficult. Such an extension is a
subject of current research.
We also work on an extension of our theory to Anderson modules \cite{and}.
\end{nrmk*}

\section{Acknowledgements}

I am very grateful to Marina Prokhorova who inspired me to become a mathematician.
It is my pleasure to thank my advisor Lenny Taelman who introduced me to Drinfeld
modules and shtukas, proposed the class number formula as a research project
and has offered continuous guidance and support.

I want to thank Fabrizio Andreatta, Yuri Bilu, Peter Bruin, Bas Edixhoven,
Lance Gurney, Arno Kret, Richard Pink, Pavel Solomatin, Peter Stevenhagen and
Dima Sustretov for very useful conversations. I am indebted to Gebhard B\"ockle
and Urs Hartl who reviewed the first version of this text.
Many thanks to Amina and Joost van Rossum for their support.

Most of this text was written
when I was a PhD student at the universities of Leiden
and Milan and a guest researcher at the University of Amsterdam. The final editing
has been done when I was a post-doc at ETH Z\"urich.
I would like to thank the universities of Amsterdam, Leiden, Milan and ETH
Z\"urich for hospitality, the ALGANT consortium, ETH Z\"urich and Leiden University
for their generous financial support.

\chapter*{Leitfaden}

\noindent
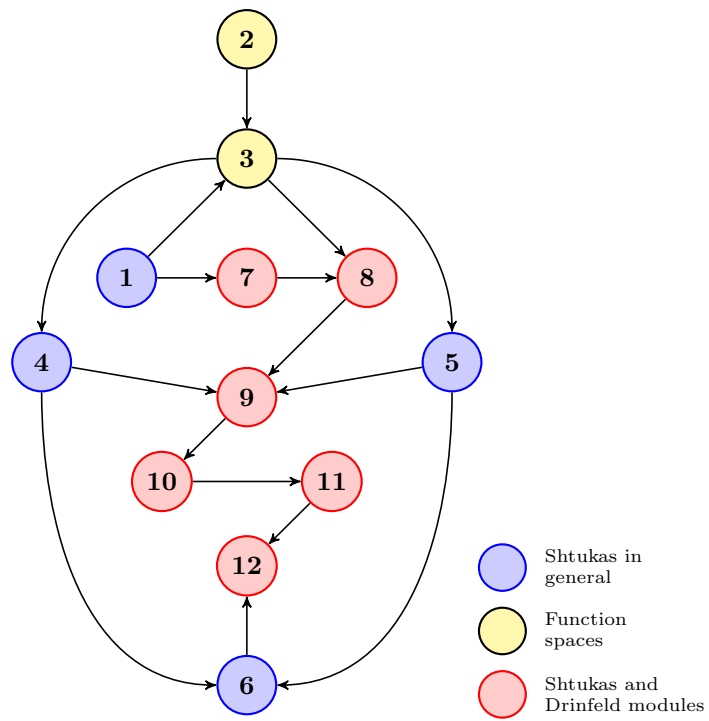
\hspace{9em}
\begin{tikzpicture}[node distance=4.5em,thick,circle,minimum size=2.2em]
\tikzstyle{zchap}=[draw=blue!100,fill=blue!20]
\tikzstyle{echap}=[draw=red!100,fill=red!20]
\tikzstyle{schap}=[draw=blue!100,fill=blue!20]
\tikzstyle{dchap}=[draw=red!100,fill=red!20]
\tikzstyle{achap}=[draw,fill=yellow!40]
\tikzstyle{dep}=[->,>=stealth',semithick]

\node[achap] (2) {\textbf{2}};
\node[achap, below of=2] (3) {\textbf{3}};
\node[dchap, below of=3] (7) {\textbf{7}};
\node[zchap, left of=7] (1) {\textbf{1}};
\node[dchap, right of=7] (8) {\textbf{8}};
\node[dchap, below of=7] (9) {\textbf{9}};
\node[schap, below left of=1] (4) {\textbf{4}};
\node[schap, below right of=8] (5) {\textbf{5}};
\node[dchap, below left of=9] (10) {\textbf{10}};
\node[dchap, below right of=9] (11) {\textbf{11}};
\node[echap, below left of=11] (12) {\textbf{12}};
\node[schap, below of=12] (6) {\textbf{6}};

\draw[dep] (1) -- (3);
\draw[dep] (2) -- (3);
\draw[dep] (3) to [out=180,in=90] (4);
\draw[dep] (3) to [out=0,in=90] (5);
\draw[dep] (4) to [out=270,in=180] (6);
\draw[dep] (5) to [out=270,in=0] (6);

\draw[dep] (1) -- (7);
\draw[dep] (3) -- (8);
\draw[dep] (7) -- (8);
\draw[dep] (4) -- (9);
\draw[dep] (5) -- (9);
\draw[dep] (8) -- (9);
\draw[dep] (9) -- (10);
\draw[dep] (10) -- (11);
\draw[dep] (11) -- (12);
\draw[dep] (6) -- (12);

\tikzstyle{lab}=[text width=8em]
\tikzstyle{cir}=[circle,thick,minimum size=2.2em]
\scriptsize

\pgftransformshift{\pgfpoint{12em}{-25em}};
\path 
(0, 0) node[cir, draw=blue!100,fill=blue!20] {}
(0, -3em) node[cir, draw,fill=yellow!40] (F) {}
(0, -6em) node[cir, draw=red!100,fill=red!20] {}

(6em, 0) node[lab] {Shtukas in \\ general}
(6em, -3em) node[lab] {Function\\spaces}
(6em, -6em) node[lab] {Shtukas and \\ Drinfeld modules};

\end{tikzpicture}

\vspace{1em}
\begin{center}
\begin{tabular}{rl}
\textbf{1} & Shtukas \\
\textbf{2} & Topological vector spaces over finite fields \\
\textbf{3} & Topological rings and modules \\
\textbf{4} & Cohomology of shtukas \\
\textbf{5} & Regulator theory \\
\textbf{6} & Trace formula \\
\textbf{7} & The motive of a Drinfeld module \\
\textbf{8} & The motive and the Hom shtuka \\
\textbf{9} & Local models \\
\textbf{10} & Change of coefficients \\
\textbf{11} & Regulators of local models \\
\textbf{12} & Global models and the class number formula
\end{tabular}
\end{center}

\chapter*{Notation and conventions}
\label{ch:lfprod}
\setcounter{section}{0}
\renewcommand*{\theHsection}{chnotconv.\the\value{section}}

\section{Stacks Project}

We use the Stacks Project \cite{stacks} as a reference for the theory of
schemes, commutative and homological algebra. We follow the conventions,
the terminology and the notation of the Stacks Project, with two amendments:
\begin{itemize}
\item
We generally write distinguished triangles as $A \to B \to C \to [1]$
omitting the object $A$ at the last place.

\item\index{idx}{complex!concentrated in degrees $[a,b]$}%
We say that a complex $A$ is \emph{concentrated in degrees} $[a,b]$
if $H^n(A) = 0$ whenever $n \not\in[a,b]$.
\end{itemize}

Since the order and numeration of items in the
Stacks Project is subject to change we refer to them by tags as explained at the
page \url{http://stacks.math.columbia.edu/tags}. A reference to a tag has the
form [$wxyz$] where ``$wxyz$'' is a combination of four letters and numbers.
The corresponding item of the Stacks Project
can be accessed by the URI
\url{http://stacks.math.columbia.edu/tag/}$wxyz$.

\section{Lattices}
\label{sec:lattices}
In Chapter~\ref{chapter:trm} and Chapters~\ref{chapter:locmod}--\ref{chapter:cnf}
we will use the notion of a lattice in a module. 

\begin{dfn}\index{idx}{lattice}%
Let $R_0 \to R$ be a ring homomorphism, $M$ an $R$-module and $M_0 \subset M$
an $R_0$-submodule. We say that $M_0$ is an $R_0$-\emph{lattice} in $M$ if the
natural map $R \otimes_{R_0} M_0 \to M$ is an isomorphism.\end{dfn}

\section{Ground field}

Throughout the text we fix a finite field $\Fq$. Correspondingly the letter $q$
stands for its cardinality. Apart from Chapter \ref{chapter:shtukas} we work
over $\Fq$. The tensor product $\otimes$ and the fiber product $\times$ without
subscripts mean the products over $\Fq$.

\section{Finite products of local fields}
\label{sec:lfprod}

In our context a local field always means a local field containing $\Fq$.
Let $F = \prod_{i=1}^n F_i$ be a finite product of local fields.
It will be convenient for us to treat such products in a uniform way independent
of $n$. To do that we set up some notation and terminology.

Observe that $F$ is a locally compact $\Fq$-algebra. It has a compact open
subalgebra $\cO_F = \prod_{i=1}^n \cO_{F_i}$ which we call the \emph{ring of
integers} of $F$. We call an element $\pi\in\cO_F$ a \emph{uniformizer} if
its projection to every $\cO_{F_i}$ is a uniformizer. Observe that $F =
\cO_F[\pi^{-1}]$. By a slight abuse of notation we denote $\fm_F \subset \cO_F$
the Jacobson radical of $\cO_F$. It is the cartesian product of maximal ideals
$\fm_{F_i} \subset \cO_{F_i}$. Every uniformizer generates $\fm_F$.
If $F$ is not a single local field then $\fm_F$ is not maximal.

\section{Mapping fiber}
\label{sec:mappingfiber}

\begin{dfn}\label{mappingfiber} Let $f\colon A \to B$ be a morphism in an
abelian category. The \emph{mapping fiber of $f$} is the complex
\begin{equation*}
\big[ A \xrightarrow{f} B \big]
\end{equation*}
where $A$ is placed in degree $0$ and $B$ in degree $1$. It coincides
with $\cone(f)[-1]$ up to sign.

We extend this definition to a morphism $f\colon A \to B$ of complexes in an
abelian category in the following way.
The mapping fiber complex
\begin{equation*}
\big[ A \xrightarrow{f} B \big]
\end{equation*}
has the object $A^n \oplus B^{n-1}$ in degree $n$ and the
differential is given by the matrix
\begin{equation*}
\begin{pmatrix}
d_A & 0 \\
f & -d_B
\end{pmatrix}
\end{equation*}
where $d_A$ and $d_B$ are differentials of $A$ respectively $B$.

Alternatively one can describe the mapping fiber complex as the total complex of the
double complex
\begin{equation*}
\xymatrix{
& \\
A^{n+1} \ar[r]^f \ar[u] & B^{n+1} \ar[u] \\
A^n \ar[r]^f \ar[u]^{d_A} & B^n \ar[u]_{-d_B} \\
\ar[u] & \ar[u]
}
\end{equation*}
where $A^n$ is placed in bidegree $(n,0)$ and $B^n$ in bidegree $(n,1)$.

Denoting the
mapping fiber complex $C$ we get a natural distinguished triangle
\begin{equation*}
C \xrightarrow{\,\,p\,\,} A \xrightarrow{\,\,f\,\,} B \xrightarrow{\,-i\,} C[1]
\end{equation*}
where $p$ is the natural projection and $i$ is the natural embedding.
The sign change for $i$ is necessary to make the triangle distinguished.
%
\end{dfn}

%

\chapter{Shtukas}
\label{chapter:shtukas}
\label{ch:shtukadefs}
\label{ch:shtukacoh}
\label{ch:nilp}
\label{ch:genhomsht}

\numberwithin{section}{chapter}


In this chapter we present a theory of shtuka cohomology together with some
supplementary constructions. By itself, shtuka
cohomology is nothing new. It usually appears in the form of explicit
complexes such as the one of Theorem \lref{shtukacoh}{shtaffcoh}
or the one of Theorem \lref{globcoh}{shtglobartcoh}.
By contrast the point of view we take in
this chapter is rather abstract. Given a scheme $X$ and an endomorphism $\tau$
we define an abelian category of shtukas on $(X,\tau)$, prove that it has enough injectives
and define a shtuka cohomology functor as the right derived functor of a certain
global sections functor.
This theory is developed for an arbitrary scheme $X$ over $\Spec \bZ$ and an arbitrary
endomorphism $\tau$. Assumptions on $X$ or $\tau$ are neither necessary nor
will they make the theory simpler.
%
%
%
%
%
%
Our treatment of shtuka cohomology was inspired by the article \cite{valeurs} of
V. Lafforgue and the book \cite{bp} of G. B\"ockle and R. Pink.

The general theory of shtuka cohomology occupies the first nine sections of
this chapter. Section~\ref{sec:shtnilp} introduces the important notion of nilpotence
borrowed from the theory of B\"ockle-Pink \cite{bp}.
The construction of $\zeta$-isomorphisms in Section~\ref{sec:shtlin} is due to
V.~Lafforgue \cite{valeurs}.
The material of Section~\ref{sec:taupoly} is well-known. In Section~\ref{sec:genhomsht}
we study a Hom shtuka construction.
Theorem \ref{homshtcohrhom} of that section relates the
cohomology of the Hom shtuka to $\RHom$ in the category of left modules
over a $\tau$-polynomial ring.
This result is of central importance to our
computations of shtuka cohomology in the context of Drinfeld modules.

In reading this chapter a certain degree of familiarity with derived
categories will be beneficial. 

\section{Basic definitions}

\begin{dfn}\label{tauring}\index{idx}{$\tau$!$\tau$-ring}%
A $\tau$-ring is a pair $(R,\tau)$
consisting of a ring $R$ and a ring
endomorphism $\tau\colon R \to R$. A morphism of $\tau$-rings $f\colon (R,\tau)
\to (S,\sigma)$ is a ring homomorphism $f\colon R \to S$ such that $f \tau =
\sigma f$.

A $\tau$-scheme\index{idx}{$\tau$!$\tau$-scheme}
is a pair $(X,\tau)$ consisting of a scheme $X$ and
an endomorphism $\tau\colon X \to X$. A morphism of $\tau$-schemes
$f\colon (X,\tau) \to (Y,\sigma)$ is a morphism of schemes $f\colon X \to Y$
such that $f \tau = \sigma f$. \end{dfn}

\breakflow
As we never work with more than one
$\tau$-ring structure on a given ring $R$ we speak of a $\tau$-ring $R$ instead
of $(R,\tau)$ and reserve the letter $\tau$ to denote the corresponding ring
endomorphism. The same applies to $\tau$-schemes.

A typical example of a $\tau$-scheme appearing in this text is the following.
Let $\Fq$ be a finite field with $q$ elements, $A$ an $\Fq$-algebra and $X$ a
smooth projective curve over $\Fq$. We equip the product $\Spec A \times_{\Fq}
X$ with the $\tau$-scheme structure given by the endomorphism which acts as the
identity on $\Spec A$ and as the $q$-Frobenius on $X$.

\begin{dfn}\label{dfnshtuka}\index{idx}{shtuka}Let $X$ be a $\tau$-scheme.
An $\cO_X$-module \emph{shtuka} is a diagram
\begin{equation*}
\cM_0 \shtuka{i}{j}\cM_1
\end{equation*}
where $\cM_0$, $\cM_1$ are $\cO_X$-modules and
\begin{align*}
i&\colon \cM_0 \to \cM_1, \\
j&\colon \cM_0 \to \tau_\ast\cM_1
\end{align*}
are morphisms of $\cO_X$-modules.
A shtuka is called \emph{quasi-coherent} if $\cM_0$ and $\cM_1$ are
quasi-coherent $\cO_X$-modules. It is called \emph{locally free} if $\cM_0$
and $\cM_1$ are locally free $\cO_X$-modules \underline{\textbf{of finite rank}}.

Let $\cM$ and $\cN$ be $\cO_X$-module shtukas given by diagrams
\begin{equation*}
\cM = \Big[\cM_0 \shtuka{i_M}{j_M} \cM_1\Big], \quad
\cN = \Big[\cN_0 \shtuka{i_N}{j_N} \cN_1\Big].
\end{equation*}
A morphism from $\cM$ to $\cN$ is a pair
$(f_0,f_1)$ where $f_n\colon \cM_n \to \cN_n$ are $\cO_X$-module morphisms such
that the diagrams
\begin{equation*}
\xymatrix{
\cM_0 \ar[d]^{i_M} \ar[r]^{f_0} & \cN_0 \ar[d]^{i_N} \\
\cM_1 \ar[r]^{f_1} & \cN_1
}\quad
\xymatrix{
\cM_0 \ar[d]^{j_M} \ar[r]^{f_0} & \cN_0 \ar[d]^{j_N} \\
\tau_\ast \cM_1 \ar[r]^{\tau_\ast(f_1)} & \tau_\ast \cN_1
}
\end{equation*}
commute.
\end{dfn}

\breakflow
Our definition of a shtuka differs from the ones present
in the literature in that we assume no restriction on $\cM_0$, $\cM_1$, $i$,
$j$, $X$ and even $\tau$. This definition is the most convenient one for our purposes. We
work with arbitrary $\cO_X$-modules instead of just the quasi-coherent ones
to make our definition of shtuka cohomology compatible with the
cohomology of coherent sheaves. The latter relies on resolutions by 
injective $\cO_X$-modules which are not quasi-coherent in general.

\section{The category of shtukas}

Let $X$ be a $\tau$-scheme. In the following we denote $\Sht \cO_X$ the category
of $\cO_X$-module shtukas. Strictly speaking $\Sht\cO_X$ depends
not only on $\cO_X$ but also on the endomorphism $\tau$. We drop $\tau$ from the
notation since we never work with more than one $\tau$-structure on a given
scheme $X$.
In this section we establish basic properties of the category
$\Sht\cO_X$.

\begin{lem}\label{shtmoradj} Let $X$ be a $\tau$-scheme. 
Let $\cM$, $\cN$ be
$\cO_X$-module shtukas defined by diagrams
\begin{equation*}
\cM = \Big[ \cM_0 \shtuka{i_M}{j_M} \cM_1 \Big], \quad
\cN = \Big[\cN_0 \shtuka{i_N}{j_N} \cN_1 \Big].
\end{equation*}
Denote
$j_M^a\colon \tau^\ast \cM_0 \to \cM_1$, 
$j_N^a\colon \tau^\ast \cN_0 \to \cN_1$
the adjoints of 
$j_M\colon \cM_0 \to \tau_\ast\cM_1$, 
$j_N\colon \cN_0 \to \tau_\ast\cN_1$
respectively.

Let $f_0\colon \cM_0 \to \cN_0$ and $f_1\colon \cM_1 \to \cN_1$ be morphisms
of $\cO_X$-modules.
The pair $(f_0,f_1)$ is a morphism of shtukas if and only if the squares
\begin{equation*}
\xymatrix{
\cM_0 \ar[d]_{i_M} \ar[r]^{f_0} & \cN_0 \ar[d]^{i_N} \\
\cM_1 \ar[r]^{f_1} & \cN_1
}\quad
\xymatrix{
\tau^\ast\cM_0 \ar[d]_{j_M^a} \ar[r]^{\tau^\ast(f_0)} & \tau^\ast\cN_0 \ar[d]^{j_N^a} \\
\cM_1 \ar[r]^{f_1} & \cN_1
}
\end{equation*}
are commutative. \quod\end{lem}

\begin{dfn}\label{defalphabeta}%
Let $X$ be a $\tau$-scheme.
We define functors from $\Sht\cO_X$ to $\cO_X$-modules:
$\alpha_\ast [ \cM_0 \shtuka{}{} \cM_1 ] = \cM_0$ and 
$\beta_\ast [ \cM_0 \shtuka{}{} \cM_1 ] = \cM_1$.%
\end{dfn}

\afterall
\begin{prp} Let $X$ be a $\tau$-scheme.
\begin{enumerate}
\item $\Sht \cO_X$ is an abelian category.

\item The functors $\alpha_\ast$ and $\beta_\ast$ are exact.

\item A sequence
$\cM \to \cM' \to \cM''$
of $\cO_X$-module shtukas is exact if and only if the induced
sequences
\begin{equation*}
\alpha_\ast\cM \to \alpha_\ast\cM' \to \alpha_\ast\cM'', \quad
\beta_\ast\cM \to \beta_\ast\cM' \to \beta_\ast\cM''.
\end{equation*}
are exact.\quod
\end{enumerate}\end{prp}

\pf $\Sht\cO_X$ is clearly an additive category. As the functor $\tau_\ast$ is left
exact it is straightforward to show that kernels in $\Sht\cO_X$ exist and
commute with $\alpha_\ast$, $\beta_\ast$. In a similar way Lemma
\ref{shtmoradj} and the fact that $\tau^\ast$ is right exact imply that cokernels
exist and commute with $\alpha_\ast$, $\beta_\ast$. A morphism of shtukas
$f\colon \cM \to \cN$ is an isomorphism if and only if $\alpha_\ast(f)$
and $\beta_\ast(f)$ are isomorphisms. Therefore $\Sht\cO_X$ is an abelian category.
(2) and (3) are clear.\quod

\begin{dfn} Let $X$ be a $\tau$-scheme. We define functors $\alpha^\ast$,
$\beta^\ast$ from the category of
$\cO_X$-modules to $\Sht\cO_X$:
\begin{equation*}
\alpha^\ast \cF = \Big[ \cF \shtuka{(1,0)}{(0,\eta)} \cF \oplus \tau^\ast \cF
\Big], \quad
\beta^\ast \cF = \Big[ 0 \shtuka{}{} \cF \Big].
\end{equation*}
Here $\eta\colon \cF \to \tau_\ast \tau^\ast \cF$ is the adjunction unit.
\end{dfn}

%
%
%

\begin{lem} $\alpha^\ast$ is left adjoint to $\alpha_\ast$ and $\beta^\ast$ is
left adjoint to $\beta_\ast$.\end{lem}

\pf The first adjunction follows from Lemma \ref{shtmoradj}. The second
adjunction is clear. \quod
 
\breakflow
The following Theorem is of fundamental importance to our treatment of
shtuka cohomology.
Recall that an object $U$ of an abelian category is called a generator if
for every nonzero morphism $f\colon A \to B$ there is a morphism $g\colon U \to
A$ such that the composition $f \circ g$ is nonzero.

\begin{thm}\label{shtcatgen} Let $X$ be a $\tau$-scheme.
\begin{enumerate}
\item $\Sht\cO_X$ has all colimits and filtered
colimits are exact.


\item $\Sht\cO_X$ admits a generator.
\end{enumerate}\end{thm}

\breakflow
It is a fundamental result of Grothendieck \cite{tohoku} that every abelian
category
satisfying (1) and (2)
has enough injective
objects.

\breakflow\noindent\textit{Proof of Theorem \ref{shtcatgen}.}
(1) Taking the direct sum of underlying $\cO_X$-modules one concludes that
$\Sht\cO_X$ has arbitrary direct sums. As it is abelian it follows that it has all
colimits. By construction the functors $\alpha_\ast$ and $\beta_\ast$ commute
with colimits. Applying $\alpha_\ast$ and $\beta_\ast$ to a colimit of
$\cO_X$-module shtukas we deduce that filtered colimits are exact in $\Sht\cO_X$
from the fact that they are exact in the category of $\cO_X$-modules.

(2) 
Consider the $\cO_X$-module
\begin{equation*}
U = \bigoplus_{V \subset X} (\iota_{V})_! \cO_V
\end{equation*}
where $V \subset X$ runs over all open subsets and $\iota_V\colon V
\hookrightarrow X$ denotes the corresponding open embedding. It is easy to see
that $U$ is a generator of the category of $\cO_X$-modules.

We claim that $\alpha^\ast U \oplus \beta^\ast U$ is a generator of
$\Sht\cO_X$.
Let $f\colon \cM \to \cN$ be a morphism of $\cO_X$-module shtukas. If $f \ne 0$
then either $\alpha_\ast f$ or $\beta_\ast f$ is nonzero, say the first one. As
$U$ is a generator there exists a morphism $g\colon U \to \alpha_\ast\cM$ such
that $\alpha_\ast f \circ g \ne 0$. As a consequence the composition of the adjoint
$g^a\colon \alpha^\ast U \to \cM$ and $f$ is
nonzero.\quod

\breakflow
Our treatment of shtuka cohomology relies on the notion of a K-injective
complex. Recall that a complex $C$ of objects in an abelian category is called
K-injective if every morphism from an acyclic complex to $C$ is zero up to
homotopy. A bounded below complex of injective objects is K-injective. In
general K-injective objects play the role of injective resolutions for unbounded
complexes. The reader who does not want to bother with unbounded complexes can
safely replace K-injective complexes with bounded below complexes of
injective objects in all the statements of this chapter. However unbounded
complexes are used in some proofs.


\begin{cor}\label{shtenoughinj} Let $X$ be a $\tau$-scheme. The category
$\Sht\cO_X$ has enough injectives. Every complex of $\cO_X$-module shtukas has a
K-injective resolution.\end{cor}

\pf By [\stacks{079I}] it follows from Theorem \ref{shtcatgen}. \quod


%
%

\section{Injective shtukas}

If $\cI$ is an injective shtuka then, as we demonstrate below, $\beta_\ast\cI$
is an injective sheaf of modules. On the contrary $\alpha_\ast\cI$ need not be
injective. Nevertheless we will show that it is good enough to compute
derived pushforwards.


\begin{lem}\label{omegainjsht}%
If $\cI$ is a K-injective complex of
$\cO_X$-module shtukas over a $\tau$-scheme $X$ then $\beta_\ast \cI$ is a
K-injective complex of $\cO_X$-modules.%
\end{lem}

\pf Immediate since $\beta_\ast$ admits an exact left adjoint $\beta^\ast$. \quod

%
%

\breakflow
In the following $\uK(\cO_X)$ stands for the homotopy category of
$\cO_X$-module complexes and $\uD(\cO_X)$ for the derived category.

Recall that a complex $\cF$ of $\cO_X$-modules on a ringed space $X$ is called
K-flat if the functor $\cF \otimes_{\cO_X} -$ preserves quasi-isomorphisms.
A bounded above complex of flat $\cO_X$-modules is K-flat. Spaltenstein
\cite{spalt} proved that every complex of $\cO_X$-modules has a K-flat
resolution.

\begin{lem}\label{shtalphaadjflat}%
Let $X$ be a $\tau$-scheme. If $\cF$ is a K-flat complex of
$\cO_X$-modules and $\cI$ a K-injective complex of $\cO_X$-module shtukas then
$\Hom_{\uK(\cO_X)}(\cF, \alpha_\ast \cI) =
\Hom_{\uD(\cO_X)}(\cF, \alpha_\ast \cI)$.%
\end{lem}

\pf Assume that $\cF$ is acyclic. By adjunction
\begin{equation*}
\Hom_{\uK(\cO_X)}(\cF,\alpha_\ast\cI) = \Hom_{\uK(\Sht \cO_X)}(\alpha^\ast\cF,\cI)
\end{equation*}
where
$\uK(\Sht\cO_X)$ is the homotopy category of $\cO_X$-module shtukas.

The complex $\tau^\ast\cF$ is acyclic since $\cF$ is K-flat. As a consequence
$\alpha^\ast\cF$ is acyclic
and the Hom on the right side of the equation is zero.
Now let $\cF$ be an arbitrary K-flat
complex and $f\colon \cF' \to \cF$ 
a quasi-isomorphism of K-flat complexes. The cone of $f$ is K-flat and acyclic.
Applying $\Hom_{\uK(\cO_X)}(-,\alpha_\ast\cI)$ to a distinguished triangle
extending $f$ we deduce that every map $g\colon \cF' \to \alpha_\ast\cI$ in
$\uK(\cO_X)$
factors through $\cF$. As every $\cO_X$-module complex admits a K-flat
resolution [\stacks{06YF}] we conclude that $\Hom_{\uK(\cO_X)}(\cF,\alpha_\ast\cI) =
\Hom_{\uD(\cO_X)}(\cF,\alpha_\ast\cI)$.\quod
 


\begin{lem}\label{alphainjsht}%
Let $X$ be a $\tau$-scheme and $f\colon X \to Y$ a morphism of
schemes. If $\cI$ is a K-injective complex of $\cO_X$-module shtukas
then the natural map $f_\ast\alpha_\ast \cI \to \uR f_\ast\alpha_\ast
\cI$ is a quasi-isomorphism.%
\end{lem}

\pf Pick a K-injective resolution $\iota\colon \alpha_\ast\cI \to \cJ$ and let
$C$ be the cone of $\iota$ so that we have a distinguished triangle
\begin{equation*}
\alpha_\ast\cI \xrightarrow{\iota} \cJ \to C \to [1]
\end{equation*}
in $\uK(\cO_X)$. 
We need to prove that $f_\ast(\iota)$ is a quasi-isomorphism or equivalently
that $f_\ast C$ is acyclic. Let $\cF$ be a K-flat $\cO_Y$-module complex.
Applying the functors 
$\Hom_{\uK(\cO_X)}(f^\ast\cF,-)$ and $\Hom_{\uD(\cO_X)}(f^\ast\cF,-)$
to the triangle above we get a morphism of long exact sequences
\begin{equation*}
\xymatrix{
\Hom_{\uK(\cO_X)}(f^\ast\cF,\alpha_\ast\cI) \ar[r] \ar[d] & \Hom_{\uD(\cO_X)}(f^\ast\cF, \alpha_\ast\cI) \ar[d] \\
\Hom_{\uK(\cO_X)}(f^\ast\cF,\cJ) \ar[r] \ar[d] & \Hom_{\uD(\cO_X)}(f^\ast\cF,\cJ) \ar[d] \\
\Hom_{\uK(\cO_X)}(f^\ast\cF, C) \ar[r] \ar[d] & \Hom_{\uD(\cO_X)}(f^\ast\cF, C) \ar[d] \\
\vdots & \vdots
}
\end{equation*}
The complex $f^\ast\cF$ is K-flat so
the top horizontal 
arrow in this diagram is an isomorphism by Lemma \ref{shtalphaadjflat}.
The middle horizontal 
arrow is an isomorphism since $\cJ$ is K-injective.
Thus the five lemma shows that
the bottom horizontal 
arrow is an isomorphism. As $C$ is acyclic we deduce that
\begin{equation*}
0 = \Hom_{\uD(\cO_X)}(f^\ast\cF,C) = \Hom_{\uK(\cO_X)}(f^\ast\cF,C) =
\Hom_{\uK(\cO_Y)}(\cF,f_\ast C)
\end{equation*}
for an arbitrary K-flat complex $\cF$.
Since the complex $f_\ast C$ admits a K-flat resolution $\cF \to f_\ast C$ we
conclude that $f_\ast C$ is acyclic.
\quod

\section{Cohomology of shtukas}

We work over a fixed $\tau$-scheme $X$.

\begin{dfn}
The \emph{ring of invariants} $\cO_X(X)^{\tau=1}$ is 
$\{ s \mid \tau(s) = s \} \subset \cO_X(X)$.\end{dfn}

\breakflow
The category of $\cO_X$-module shtukas is
$\cO_X(X)^{\tau=1}$-linear by construction.

\begin{dfn} 
The \emph{unit shtuka} $\unit_X$ is defined by the diagram
\begin{equation*}
\cO_X \shtuka{\,\,\,1\,\,\,}{\tau^\sharp} \cO_X
\end{equation*}
where $\tau^\sharp\colon \cO_X \to \tau_\ast \cO_X$ is the homomorphism of
sheaves of rings determined by $\tau$.
\end{dfn}

\begin{dfn}\index{idx}{shtuka cohomology}\index{nidx}{shtuka cohomology!RGamma@$\RGamma$, cohomology complex}%
We define the \emph{cohomology functor} $\RGamma(X,-)$ from
the derived category of $\Sht\cO_X$ to the derived category of
$\cO_X(X)^{\tau=1}$-modules as follows:
\begin{equation*}
\RGamma(X,\,\cM) = \RHom(\unit_X, \,\cM).
\end{equation*}
We call $\RGamma(X,\,\cM)$ the
\emph{cohomology complex} of $\cM$ or simply the cohomology of $\cM$.
The $n$-th cohomology module of $\RGamma(X,\,\cM)$ is denoted
$\uH^n(X,\,\cM)$.
By construction $\uH^n(X,\,\cM) = \Ext^n(\unit_X,\,\cM)$.%
\end{dfn}

\breakflow
Let $\cM$ be an $\cO_X$-module shtuka given by a diagram
\begin{equation*}
\cM_0 \shtuka{i}{j} \cM_1\Big.
\end{equation*}
The arrows of $\cM$ determine natural maps
\begin{equation*}
i,j\colon \Gamma(X,\,\cM_0) \to \Gamma(X,\,\cM_1)
\end{equation*}
with the same source and target.
In the case of $j$ we identify $\Gamma(X,\,\tau_\ast\cM_1)$ with $\Gamma(X,\,\cM_1)$
using the fact that $\tau^{-1} X = X$.
Observe that $j$ is only
$\cO_X(X)^{\tau=1}$-linear since the natural identification
$\Gamma(X,\,\cM_1) = \Gamma(X,\,\tau_\ast\cM_1)$
is $\tau$-linear. 

\begin{prp}\label{shtcohzerodesc}%
$\uH^0(X,\,\cM) = \{\, s \in \Gamma(X,\,\cM_0) \mid i(s) = j(s) \,\}$.%
\end{prp}

\pf A morphism $f\colon \unit_X \to \cM$ is a pair of maps
$f_0\colon \cO_X \to \cM_0$, 
$f_1\colon \cO_X \to \cM_1$
such that
$i \circ f_0 = f_1$ and 
$j \circ f_0 = \tau_\ast(f_1) \circ \tau^\sharp$.
The pair $(f_0,f_1)$ is determined by the section $s = f_0(1)$ of
$\Gamma(X,\cM_0)$ which satisfies the equation $i(s) = j(s)$.\quod

\section{Canonical triangle}

The constructions of this section are due to
V.~Lafforgue \cite[Section 4]{valeurs}.

\begin{dfn}We define a morphism
$\delta\colon \beta^\ast\cO_X \to \alpha^\ast\cO_X$
by the diagram
\begin{equation*}
\xymatrix{
0 \ar[rr] \ar@<-.75ex>[d] \ar@<.75ex>[d] &&
\cO_X \ar@<-.75ex>[d]_{(1,0)} \ar@<.75ex>[d]^{(0,\tau^\sharp)} \\
\cO_X \ar[rr]^{(1,-1)} && \cO_X\oplus\cO_X
}
\end{equation*}\end{dfn}

\begin{lem}\label{shtdij}%
For every $\cO_X$-module shtuka $\cM$ the diagram
\begin{equation*}
\xymatrix{
\Hom(\alpha^\ast\cO_X,\cM) \ar[r]^{\delta \circ -} \ar@{=}[d] &
\Hom(\beta^\ast\cO_X,\cM) \ar@{=}[d] \\
\Gamma(X,\,\alpha_\ast\cM) \ar[r]^{i-j} & \Gamma(X,\,\beta_\ast\cM)
}
\end{equation*}
is commutative.%
\quod\end{lem}

\begin{lem}\label{deralphabeta}%
For every complex of $\cO_X$-module shtukas the maps
\begin{align*}
\RHom(\alpha^\ast\cO_X, \cM) &\to \RGamma(X, \,\alpha_\ast\cM) \\
\RHom(\beta^\ast\cO_X, \cM) &\to \RGamma(X, \,\beta_\ast\cM)
\end{align*}
induced by the natural identifications
\begin{align*}
\Hom(\alpha^\ast\cO_X, \cM) &= \Gamma(X, \,\alpha_\ast\cM) \\
\Hom(\beta^\ast\cO_X, \cM) &= \Gamma(X, \,\beta_\ast\cM)
\end{align*}
are quasi-\hspace{0pt}isomorphisms.\end{lem}

\pf Without loss of generality we assume that $\cM$ is K-injective.
The result for $\beta$ is then immediate since $\beta_\ast\cM$ is K-injective by Lemma~\ref{omegainjsht}.
Applying Lemma~\ref{alphainjsht} to the map $X \to \Spec\bZ$
we conclude that $\Gamma(X,\,\alpha_\ast\cM) = \RGamma(X,\,\alpha_\ast\cM)$.\quod

\breakflow
Let $\cM$ be a complex of $\cO_X$-module shtukas.
The arrows of the shtukas of $\cM$ determine natural maps
\begin{equation*}
i,j\colon \RGamma(X,\alpha_\ast\cM) \to \RGamma(X,\beta_\ast\cM),
\end{equation*}
in the following way.
The first map is induced by the $i$-arrows.
The $j$-arrows induce a map
$\RGamma(X,\alpha_\ast\cM) \to \RGamma(X,\tau_\ast\beta_\ast\cM)$.
Taking its composition with the natural map
$\RGamma(X,\tau_\ast\beta_\ast\cM) \to \RGamma(X,\uR\tau_\ast\beta_\ast\cM)$
and using the identity
$\RGamma(X,\beta_\ast\cM) = \RGamma(X,\uR\tau_\ast\beta_\ast\cM)$
we get a map of the desired form.

\begin{lem}\label{shtdijder}%
For every complex of $\cO_X$-module shtukas $\cM$ the square
\begin{equation*}
\xymatrix{
\RHom(\alpha^\ast\cO_X, \,\cM) \ar[rr]^{\delta \circ -} \ar[d]_{\ltviso} &&
\RHom(\beta^\ast\cO_X, \,\cM) \ar[d]^{\rtviso} \\
\RGamma(X,\,\alpha_\ast\cM) \ar[rr]^{i-j} &&
\RGamma(X,\,\beta_\ast\cM)
}
\end{equation*}
is commutative.\end{lem}

\pf Without loss of generality we may assume that $\cM$ is K-injective.
The result then follows from Lemmas \ref{deralphabeta} and \ref{shtdij}.\quod

%
%
%
%

\begin{dfn}We define a morphism
$\rho\colon\alpha^\ast \cO_X \to \unit_X$
by the diagram
\begin{equation*}
\xymatrix{
\cO_X \ar[rr]^1 \ar@<-.75ex>[d]_{(1,0)} \ar@<.75ex>[d]^{(0,\tau^\sharp)} &&
\cO_X \ar@<-.75ex>[d]_{1} \ar@<.75ex>[d]^{\tau^\sharp} \\
\cO_X\oplus\cO_X \ar[rr]^{\sum} && \cO_X
}
\end{equation*}\end{dfn}

\begin{dfn}\label{defunitcomp}%
We denote $\unic_X = \cone\big(\beta^\ast\cO_X \xrightarrow{-\delta} \alpha^\ast\cO_X\big)$.
We define a map $\pi\colon\unic_X \to \beta^\ast\cO_X[1]$ by the identity in degree $-1$
and a map $\widetilde{\rho}\colon\unic_X \to \unit_X[0]$ 
by $\rho$ in degree $0$.\end{dfn}

\begin{lem}$\widetilde{\rho}$ is a quasi-isomorphism.\quod\end{lem}

%

\begin{prp}\label{unittriok}%
The triangle
\begin{equation}\label{unittri}
\unit_X[-1] \xrightarrow{\,-\pi\widetilde{\rho}^{-1}[-1]\,}
\beta^\ast\cO_X[0] \xrightarrow{\quad \delta\quad}
\alpha^\ast\cO_X[0] \xrightarrow{\quad \rho\quad}
\unit_X[0]
\end{equation}
in $\uD(\Sht\cO_X)$ is distinguished.\end{prp}

\breakflow
The following definition is central to our theory:

\begin{dfn}\label{defcantri}\index{idx}{shtuka cohomology!canonical triangle}%
Let $\cM$ be a complex of $\cO_X$-module shtukas.
We define the \emph{canonical distinguished triangle}
\begin{equation*}
\RGamma(X,\,\cM) \xrightarrow{\phantom{\,\, i-j\,\,}}
\RGamma(X,\,\alpha_\ast\cM) \xrightarrow{\,\, i-j\,\,}
\RGamma(X,\,\beta_\ast\cM)\xrightarrow{\phantom{\,\, i-j\,\,}}
[1]
\end{equation*}
by applying $\RHom(-,\cM)$ to the triangle
\eqref{unittri} and using the identifications
of Lemma \ref{shtdijder}.\end{dfn}

\afterall\noindent\textit{Proof of Proposition \ref{unittriok}.}
Let $\iota\colon\unic_X \to \alpha^\ast\cO_X[0]$ be the map
given by the identity in degree $0$.
By definition of $\unic_X$ the triangle
\begin{equation*}
\beta^\ast\cO_X[0] \xrightarrow{\,\,-\delta\,\,}
\alpha^\ast\cO_X[0] \xrightarrow{\,\,\iota\,\,}
\unic_X \xrightarrow{\,\,-\pi\,\,}
\beta^\ast\cO_X[1]
\end{equation*}
is distinguished [\stacks{014P}]. Hence so is the isomorphic
triangle
\begin{equation*}
\beta^\ast\cO_X[0] \xrightarrow{\,\,-\delta\,\,}
\alpha^\ast\cO_X[0] \xrightarrow{\,\,\rho\,\,}
\unit_X[0] \xrightarrow{\,\,-\pi\widetilde{\rho}^{-1}\,\,}
\beta^\ast\cO_X[1]
\end{equation*}
Rotating it we obtain a triangle
\begin{equation*}
\unit_X[-1] \xrightarrow{\,\pi\widetilde{\rho}^{-1}[-1]\,}
\beta^\ast\cO_X[0] \xrightarrow{\,\,-\delta\,\,}
\alpha^\ast\cO_X[0] \xrightarrow{\,\,\rho\,\,}
\unit_X[0]
\end{equation*}
which is isomorphic to \eqref{unittri}.\quod

\section{Associated complex}
\label{sec:assoccompl}

It will often be convenient for us to view the functor $\RGamma$ on shtukas not
as the derived functor of $\Hom(\unit_X,-)$
but as the derived functor of the so-called \emph{associated complex} functor
which we now introduce

\begin{dfn}\label{assoccompl}%
\index{idx}{shtuka cohomology!associated complex}\index{nidx}{shtuka cohomology!cGamma@$\cGamma$, associated complex}%
Let $\cM$ be a complex of $\cO_X$-module shtukas.
We define
\begin{equation*}
\cGamma(X,\,\cM) = \Hom^\bullet(\unic_X,\,\cM)
\end{equation*}
and call $\cGamma(X,\,\cM)$ the \emph{associated complex} of $\cM$.
Here $\Hom^\bullet$ is the Hom complex [\stacks{09L9}].
We define a natural morphism
$\cGamma(X,\,\cM) \to \RGamma(X,\,\cM)$
as the composition
\begin{equation*}
\Hom^\bullet(\unic_X,\,\cM)\xrightarrow{\phantom{\,\widetilde{q} \circ -\,}}
\RHom(\unic_X,\,\cM) \xrightarrow{\,\widetilde{\rho}^{-1} \circ -\,}
\RHom(\unit_X,\,\cM)
\end{equation*}
where $\widetilde{\rho}$ is the quasi-\hspace{0pt}isomorphism
of Definition \ref{defunitcomp}.\end{dfn}


\begin{example*}
Let an $\cO_X$-module shtuka $\cM$ be given by a diagram
\begin{equation*}
\cM_0 \shtuka{i}{j} \cM_1.
\end{equation*}
Regarding $\cM$ as a complex of shtukas concentrated in degree $0$ we have
\begin{equation*}
\cGamma(X,\,\cM) = \Big[ \Gamma(X,\,\cM_0) \xrightarrow{i-j} \Gamma(X,\,\cM_1) \Big].
\end{equation*}
The square brackets denote the mapping fiber complex of Chapter ``\chnotconv''.
\end{example*}

\begin{lem}\label{assoccompldesc}%
Let $\cM$ be a
complex of $\cO_X$-module shtukas. Denote
$\cM^n$ the shtuka in degree $n$. Then $\cGamma(X,\,\cM)$
coincides with the total complex of the double complex
\begin{equation*}
\xymatrix{
&& \\
\Gamma(X,\,\alpha_\ast\cM^{n+1}) \ar[rr]^{(-1)^{n+1}(i-j)} \ar[u] &&
\Gamma(X,\,\beta_\ast\cM^{n+1}) \ar[u] \\
\Gamma(X,\,\alpha_\ast\cM^{n}) \ar[rr]^{(-1)^{n}(i-j)} \ar[u] &&
\Gamma(X,\,\beta_\ast\cM^{n}) \ar[u] \\
\ar[u] && \ar[u]
}
\end{equation*}
The object $\Gamma(X,\,\alpha_\ast\cM^n)$
is placed in the bidegree $(n,0)$ while $\Gamma(X,\,\beta_\ast\cM^n)$ is in the
bidegree $(n,1)$.
The vertical maps are the differentials of $\Gamma(X,\,\alpha_\ast\cM)$
respectively $\Gamma(X,\,\beta_\ast\cM)$.
The maps $i$ and $j$ are induced by the arrows of
the shtukas $\cM^n$.\end{lem}

\pf In view of Lemma \ref{shtdij} it follows from the definition of $\Hom^\bullet$
[\stacks{09L9}].\quod

\begin{dfn}\label{defassoctri}%
Let $\cM$ be a complex of $\cO_X$-module shtukas.
We define a triangle
\begin{equation}\label{assoccompltri}
\cGamma(X,\,\cM) \xrightarrow{\,\,p\,\,}
\Gamma(X,\,\alpha_\ast \cM) \xrightarrow{\,i-j\,}
\Gamma(X,\,\beta_\ast\cM) \xrightarrow{\,-q\,}
\cGamma(X,\,\cM)[1]
\end{equation}
as follows. Denote $\cM^n$ the shtuka in degree $n$.
According to Lemma \ref{assoccompldesc} the object of $\cGamma(X,\cM)$ in
degree $n$ is
\begin{equation*}
\Gamma(X,\alpha_\ast\cM^n) \oplus \Gamma(X,\beta_\ast\cM^{n-1}).
\end{equation*}
The morphism $p$ is the natural projection. 
The morphism $q$ is defined by the formula $b \mapsto
(0,(-1)^n b)$ in degree $n$.
The maps $i$ and $j$ are induced by the arrows of the shtukas $\cM^n$.%
\end{dfn}

\begin{prp}\label{assocmor}%
For every complex of $\cO_X$-module shtukas $\cM$ 
the following holds:
\begin{enumerate}
\item The triangle \eqref{assoccompltri} is distinguished.

\item\label{assocmordiag}%
The natural diagram
\begin{equation*}%
\xymatrix{
\cGamma(X,\,\cM) \ar[r]^{p\,\,} \ar[d] &
\Gamma(X,\,\alpha_\ast\cM) \ar[r]^{i-j} \ar[d] &
\Gamma(X,\,\beta_\ast\cM) \ar[r]^{\quad\,\,\,\,-q} \ar[d] & [1] \\
\RGamma(X,\,\cM) \ar[r] & \RGamma(X,\alpha_\ast\cM) \ar[r]^{i-j} &
\RGamma(X,\,\beta_\ast\cM)\ar[r] & [1]
}
\end{equation*}%
is a morphism of triangles. 
Here
the bottom row is the canonical triangle of Definition~\ref{defcantri},
the first vertical map is the one of Definition~\ref{assoccompl},
the second and the third are the natural ones.%

\item If $\cM$ is K-injective then the diagram \eqref{assocmordiag} is an isomorphism.%
\end{enumerate}%
\end{prp}%

\begin{rmk}%
In particular the triangle \eqref{assoccompltri} gives an explicit description
of the canonical triangle of Definition~\ref{defcantri} in the case when
$\cM$ is K-injective.\end{rmk}

\afterall\noindent\textit{Proof of Proposition~\ref{assocmor}.}
Consider a commutative diagram
\begin{equation}\label{unictrimor}
\vcenter{\vbox{\xymatrix{
\unic_X[-1] \ar[rr]^{-\pi} \ar[d]_{\widetilde{\rho}[-1]}^{\rtviso} &&
\beta^\ast\cO_X[0] \ar[r]^{\delta} \ar@{=}[d] &
\alpha^\ast\cO_X[0] \ar[r]^{\iota} \ar@{=}[d] &
\unic_X \ar[d]_{\ltviso}^{\widetilde{\rho}} \\
\unit_X[-1] \ar[rr]^{-\pi\widetilde{\rho}^{-1}[-1]} &&
\beta^\ast\cO_X[0] \ar[r]^{\delta} &
\alpha^\ast\cO_X[0] \ar[r]^{\rho} &
\unit_X[0]
}}}
\end{equation}
The map $\iota$ is given by the identity in degree $0$.
The bottom triangle is distinguished by Proposition \ref{unittriok}.
Hence so is the top one.

We next apply $\Hom^\bullet(-,\cM)$ to the top row.
We use the map given by $(-1)^n$ in degree $n$
to identify $\Hom^\bullet(-[-1],\cM)$ with
$\Hom^\bullet(-,\cM)[1]$.
A straightforward computation shows that the resulting
triangle is \eqref{assoccompltri}. In particular this triangle is
distinguished.

Applying $\RHom(-,\cM)$ to \eqref{unictrimor} 
and using Lemma \ref{shtdijder}
we obtain an isomorphism from the canonical distinguished triangle of $\cM$
to a distinguished triangle
\begin{equation*}
\RHom(\unic_X,\,\cM) \xrightarrow{\phantom{\,i-j\,}}
\RGamma(X,\,\alpha_\ast\cM) \xrightarrow{\,i-j\,}
\RGamma(X,\,\beta_\ast\cM) \xrightarrow{\phantom{\,i-j\,}}
[1].
\end{equation*}
Composing the inverse of this isomorphism with the map of triangles
induced by the natural morphism
$\Hom^\bullet(-,\cM) \to \RHom(-,\cM)$ we get the diagram
\eqref{assocmordiag}.
It is therefore a morphism of triangles.

If $\cM$ is K-injective then the third vertical arrow
in \eqref{assocmordiag} is a quasi-\hspace{0pt}isomorphism
by Lemma~\ref{omegainjsht}. Applying Lemma~\ref{alphainjsht} to the
natural map $X \to \Spec\bZ$ we conclude that the second
vertical arrow of \eqref{assocmordiag} is a quasi-\hspace{0pt}isomorphism.
Hence so is the first one, and the diagram \eqref{assocmordiag}
is an isomorphism of triangles.\quod

\section{Pushforward}

We work with a fixed morphism of $\tau$-schemes $f\colon X \to Y$.

\begin{dfn}%
Let $\cM$ be an $\cO_X$-module shtuka given by a diagram
\begin{equation*}
\cM_0 \shtuka{i}{j} \cM_1.
\end{equation*}
Define 
\begin{equation*}
f_\ast \cM = \Big[ f_\ast \cM_0 \shtuka{f_\ast i}{f_\ast j} f_\ast \cM_1 \Big].
\end{equation*}
Here we use the natural isomorphism $f_\ast \tau_\ast \cM_1 = \tau_\ast f_\ast
\cM_1$ to interpret $f_\ast j$ as a map to $\tau_\ast f_\ast \cM_1$.\end{dfn}

\begin{dfn}%
The functor $f_\ast$ on the category of $\cO_X$-module shtukas is
left exact. We define $\uR f_\ast$ as its right derived functor.\end{dfn}

\begin{lem}\label{alphaomegapush}%
%
The natural maps 
$\alpha_\ast\uR f_\ast \to \uR f_\ast \alpha_\ast$
and
$\beta_\ast\uR f_\ast \to \uR f_\ast \beta_\ast$
are quasi-\hspace{0pt}isomorphisms.%
\end{lem}

\pf Let $\cM$ be a K-injective complex of $\cO_X$-module shtukas.
Lemma \ref{alphainjsht} shows that the natural map
$f_\ast\alpha_\ast\cM \to \uR f_\ast\alpha_\ast \cM$ is a quasi-isomorphism.
However $f_\ast\alpha_\ast\cM = \alpha_\ast f_\ast\cM = \alpha_\ast\uR f_\ast\cM$
and we get the result for $\alpha_\ast$.
The result for $\beta_\ast$ follows in a similar way
since the complex $\beta_\ast\cI$ is K-injective
by Lemma \ref{omegainjsht}.\quod

\breakflow
Let $\cM$ be a complex of $\cO_X$-module shtukas.
We have a natural isomorphism
$\Hom(\unit_X,\cM) = \Hom(\unit_Y,\,f_\ast \cM)$. It induces a
map $\RGamma(X,\cM) \to \RGamma(Y,\,\uR f_\ast \cM)$. Furthermore
we have a natural quasi-isomorphism
\begin{equation*}
\RGamma(X,\,\alpha_\ast\cM) \xrightarrow{\,\isosign\,}
\RGamma(Y,\,\uR f_\ast \alpha_\ast \cM) \xrightarrow{\,\isosign\,}
\RGamma(Y,\,\alpha_\ast \uR f_\ast \cM)
\end{equation*}
where the second arrow is the quasi-\hspace{0pt}isomorphism of Lemma \ref{alphaomegapush}.
In a similar way we have a natural quasi-\hspace{0pt}isomorphism
$\RGamma(X,\,\beta_\ast\cM) \xrightarrow{\,\isosign\,}\RGamma(Y,\,\beta_\ast \uR f_\ast \cM)$.

\begin{prp}\label{shtpushmor}%
For every complex $\cM$ of $\cO_X$-module shtukas
the natural map
$\RGamma(X,\,\cM) \to \RGamma(Y,\,\uR f_\ast\cM)$
is a quasi-\hspace{0pt}isomorphism.
Moreover the natural diagram
\begin{equation}\label{shtpushmordiag}
\vcenter{\vbox{
\xymatrix{
\RGamma(X,\,\cM) \ar[r] \ar[d] & \RGamma(X,\,\alpha_\ast\cM) \ar[r]^{i-j} \ar[d] &
\RGamma(X,\,\beta_\ast\cM) \ar[r] \ar[d] & [1] \\
\RGamma(Y,\,\uR f_\ast \cM) \ar[r] & \RGamma(Y,\,\alpha_\ast \uR f_\ast \cM) \ar[r]^{i-j} &
\RGamma(Y,\,\beta_\ast \uR f_\ast \cM) \ar[r] & [1] \\
}}}
\end{equation}
is an isomorphism of distinguished triangles.\end{prp}

\pf Without loss of generality we assume that $\cM$ is K-injective
so that $f_\ast\cM = \uR f_\ast\cM$.
Let $\cI$ be a K-injective resolution of $f_\ast\cM$.
The map $f_\ast\cM \to \cI$ induces a morphism of
triangles \eqref{assoccompltri}:
\begin{equation}\label{fpushmor}
\vcenter{\vbox{\xymatrix{
\cGamma(Y,\,f_\ast \cM) \ar[r] \ar[d] &
\Gamma(Y,\,\alpha_\ast f_\ast \cM) \ar[r]^{i-j} \ar[d] &
\Gamma(Y,\,\beta_\ast f_\ast \cM) \ar[r] \ar[d] & [1] \\
\cGamma(Y,\,\cI) \ar[r] & \Gamma(Y,\,\alpha_\ast\cI) \ar[r]^{i-j} &
\Gamma(Y,\,\beta_\ast\cI) \ar[r] & [1]
}}}
\end{equation}
The top triangle coincides with the
triangle \eqref{assoccompltri} for $\cGamma(X,\,\cM)$.
Proposition \ref{assocmor} then identifies the diagram \eqref{fpushmor}
with the diagram \eqref{shtpushmordiag}. Whence the result. \quod

\section{Pullback}

Let $f\colon X \to Y$ be a morphism of $\tau$-schemes 
and let
\begin{equation*}
\cM = \Big[\cM_0 \shtuka{i}{j} \cM_1\Big]
\end{equation*}
be an $\cO_Y$-module shtuka. 

\begin{dfn}We set
\begin{equation*}
f^\ast \cM = \Big[ f^\ast \cM_0 \shtuka{f^\ast i}{\mu\circ f^\ast j} f^\ast \cM_1 \Big]
\end{equation*}
where $\mu$ is
the base change map $f^\ast \tau_\ast \cM_1 \to \tau_\ast f^\ast \cM_1$ arising
from the commutative square
\begin{equation*}
\xymatrix{
X \ar[r]^f \ar[d]^\tau & Y \ar[d]^\tau \\
X \ar[r]^f & Y.
}
\end{equation*}
We will use the following notation:
\begin{itemize}
\item For an affine $\tau$-scheme $X = \Spec R$ 
we denote
$\cM(R)$ the $R$-module shtuka
\begin{equation*}
\Gamma(X, f^\ast \cM_0) \shtuka{f^\ast i}{\mu \circ f^\ast j} \Gamma(X, f^\ast \cM_1)
\end{equation*}
arising by pullback of $\cM$ along $f\colon X \to Y$.

\item To make the expressions more legible we will generally write $\RGamma(X,\cM)$
in place of $\RGamma(X, f^\ast \cM)$.\end{itemize}\end{dfn}

\begin{lem}
There exists a unique adjunction
\begin{equation*}
\Hom_{\Sht \cO_X}(f^\ast -, -) \cong \Hom_{\Sht \cO_Y}(-,f_\ast -)
\end{equation*}
which is compatible with the adjunction
\begin{equation*}
\Hom_{\cO_X}(f^\ast -, -) \cong \Hom_{\cO_Y}(-,f_\ast -)
\end{equation*}
through the natural maps given by functors $\alpha_\ast$ and $\beta_\ast$.
\quod
\end{lem}

\begin{dfn}
We define the \emph{pullback map}
\begin{equation*}
\RGamma(Y,\cM) \xrightarrow{\quad f^\ast \quad} \RGamma(X, f^\ast\cM)
\end{equation*}
in the following way.
Let $\eta\colon \cM \to f_\ast f^\ast \cM$ be the adjunction unit. Taking its
composition with the natural map $f_\ast f^\ast \cM \to \uR f_\ast f^\ast \cM$
and applying $\RGamma(Y, -)$ we obtain a map from $\RGamma(Y,\cM)$ to
$\RGamma(Y,\uR f_\ast f^\ast \cM)$. Proposition \ref{shtpushmor} identifies $\RGamma(Y,\uR
f_\ast f^\ast \cM)$ with $\RGamma(X,f^\ast \cM)$. The resulting
map from $\RGamma(Y,\cM)$ to $\RGamma(X, f^\ast\cM)$ is the pullback map.
\end{dfn}

\breakflow
Observe that $\alpha_\ast f^\ast \cM = f^\ast \alpha_\ast\cM$ and
$\beta_\ast f^\ast \cM = f^\ast \beta_\ast \cM$ by construction. 
So we have natural pullback maps 
$\RGamma(Y, \alpha_\ast \cM) \to \RGamma(X, \alpha_\ast f^\ast \cM)$
and $\RGamma(Y,\beta_\ast \cM) \to \RGamma(X, \beta_\ast f^\ast \cM)$.

\begin{prp}\label{shtpullmor}
For every $\cO_Y$-module shtuka $\cM$ the natural diagram
\begin{equation}\label{shtpullmordiag}
\vcenter{\vbox{\xymatrix{
\RGamma(Y,\,\cM) \ar[r] \ar[d] &
\RGamma(Y,\,\alpha_\ast\cM) \ar[r]^{i-j} \ar[d] &
\RGamma(Y,\,\beta_\ast\cM) \ar[r] \ar[d] & [1] \\
\RGamma(X,\,f^\ast \cM) \ar[r] &
\RGamma(X,\,\alpha_\ast f^\ast\cM) \ar[r]^{i-j} &
\RGamma(X,\,\beta_\ast f^\ast\cM) \ar[r] & [1]
}}}
\end{equation}
is a morphism of distinguished triangles.\end{prp}

\pf The natural map $\cM \to \uR f_\ast f^\ast \cM$ induces a morphism of
distinguished triangles
\begin{equation}\label{shtpullmorstage1}
\vcenter{\vbox{\xymatrix@C=0.835em{
\RGamma(Y,\,\cM) \ar[rr] \ar[d] &&
\RGamma(Y,\,\alpha_\ast\cM) \ar[rr]^{i-j} \ar[d] &&
\RGamma(Y,\,\beta_\ast\cM) \ar[d] \ar[r] & [1] \\
\RGamma(Y,\,\uR f_\ast f^\ast \cM) \ar[rr] &&
\RGamma(Y,\,\alpha_\ast \uR f_\ast f^\ast\cM) \ar[rr]^{i-j} &&
\RGamma(Y,\,\beta_\ast \uR f_\ast f^\ast\cM) \ar[r] & [1]
}}}
\end{equation}
At the same time Proposition \ref{shtpushmor} states that the natural diagram
\begin{equation}\label{shtpullmorstage2}
\vcenter{\vbox{\xymatrix@C=0.835em{
\RGamma(X,\,f^\ast\cM) \ar[rr] \ar[d]^{\rtviso} &&
\RGamma(X,\,\alpha_\ast f^\ast\cM) \ar[rr]^{i-j} \ar[d]^{\rtviso} &&
\RGamma(X,\,\beta_\ast f^\ast\cM) \ar[r] \ar[d]^{\rtviso} & [1] \\
\RGamma(Y,\,\uR f_\ast f^\ast \cM) \ar[rr] &&
\RGamma(Y,\,\alpha_\ast \uR f_\ast f^\ast \cM) \ar[rr]^{i-j} &&
\RGamma(Y,\,\beta_\ast \uR f_\ast f^\ast \cM) \ar[r] & [1] \\
}}}
\end{equation}
is an isomorphism of distinguished triangles. A quick inspection shows that
the composition of \eqref{shtpullmorstage1} and the inverse of
\eqref{shtpullmorstage2} gives the diagram \eqref{shtpullmordiag}. \quod


\section{Shtukas over affine schemes}
\label{sec:shtaffcoh}

It follows from Definition \ref{dfnshtuka} that a quasi-coherent shtuka $\cM$ on an affine
$\tau$-scheme $X = \Spec R$ is given by a diagram
\begin{equation*}
M_0 \shtuka{i}{j} M_1
\end{equation*}
where $M_0$, $M_1$ are $R$-modules, $i\colon M_0 \to M_1$ an $R$-module
homomorphism and $j\colon M_0 \to M_1$ a $\tau$-linear $R$-module homomorphism:
for all $r \in R$ and $m \in M_0$ one has $j(r m) = \tau(r) j(m)$. The
associated complex of $\cM$ is
\begin{equation*}
\cGamma(X, \cM) = \Big[ M_0 \xrightarrow{i-j} M_1 \Big].
\end{equation*}
We will show that this complex computes the cohomology of $\cM$.

\begin{thm}\label{shtaffcoh}%
If $\cM$ is a quasi-coherent shtuka over an affine
$\tau$-scheme $X$ then the natural map
$\cGamma(X,\,\cM) \to \RGamma(X,\,\cM)$ is a quasi-\hspace{0pt}isomorphism.\end{thm}

\pf By Proposition \ref{assocmor} the natural map in question extends to a morphism
of distinguished triangles
\begin{equation*}
\xymatrix{
\cGamma(X,\,\cM) \ar[r] \ar[d] & \Gamma(X,\,\cM_0) \ar[r]^{i-j}
\ar[d] & \Gamma(X,\,\cM_1) \ar[r] \ar[d] & [1] \\
\RGamma(X,\,\cM) \ar[r] & \RGamma(X,\,\cM_0) \ar[r]^{i-j}
& \RGamma(X,\,\cM_1) \ar[r] & [1]
}
\end{equation*}
As $\cM_0$ and $\cM_1$ are quasi-\hspace{0pt}coherent $\cO_X$-modules over an affine scheme
$X$ the complexes $\RGamma(X,\,\cM_0)$ and $\RGamma(X,\,\cM_1)$ are
concentrated in degree $0$ [\stacks{01XB}]. Hence the second and third
vertical maps in the diagram above are quasi-\hspace{0pt}isomorphism. It follows that so is
the first map. \quod

\breakflow
To make the expressions more legible
we will often write $\RGamma(R,\,\cM)$ instead of
$\RGamma(\Spec R,\,\cM)$. If there is no ambiguity in the choice of $R$ then we
further shorten it to $\RGamma(\cM)$. For a quasi-coherent shtuka
$\cM$ we identify $\RGamma(\cM)$ with $\cGamma(X,\,\cM)$ using the
Theorem above.

\section{Nilpotence}
\label{sec:shtnilp}

The notion of nilpotence for shtukas is crucial to this work. It first appeared
in the book of B\"ockle-Pink \cite{bp} in the context of $\tau$-sheaves.

\begin{dfn}\label{defnilp}\index{idx}{shtuka!nilpotent}%
Let $X$ be a $\tau$-scheme. An $\cO_X$-module shtuka
\begin{equation*}
\cM_0 \shtuka{i}{j} \cM_1
\end{equation*}
is called \emph{nilpotent} if $i$ is an isomorphism and the composition
\begin{equation}\label{nilpdiag}
\cM_0 \xrightarrow{\,\mu\,} \tau_\ast \cM_0 \xrightarrow{\,\tau_\ast \mu\,}
\tau_\ast^2 \cM_0 \to \dotsc \to \tau_\ast^n \cM_0, \quad
\mu = \tau_\ast(i^{-1}) \circ j,
\end{equation}
is zero for some $n \geqslant 1$.\end{dfn}

\begin{prp}Let $f\colon X \to Y$ be a morphism of $\tau$-schemes and let $\cM$
be an $\cO_Y$-module shtuka. If $\cM$ is nilpotent then $f^\ast \cM$ is
nilpotent.\end{prp}

\pf Without loss of generality we assume that
\begin{equation*}
\cM = \Big[ \cM_0 \shtuka{1}{j} \cM_0 \Big].
\end{equation*}
Let $j^a\colon \tau^\ast\cM_0 \to \cM_0$ be the adjoint of $j$.
Using the naturality of the adjunction $\tau^\ast \leftrightarrow \tau_\ast$
it is easy to show that \eqref{nilpdiag} is zero
if and only if the composition
\begin{equation*}
\tau^{\ast n}(\cM_0) \xrightarrow{\,\,\tau^{\ast (n-1)}(j^a)\,\,}
\tau^{\ast (n-1)}(\cM_0) \to \dotsc \to
\tau^\ast(\cM_0) \xrightarrow{\,\,j^a\,\,} \cM_0
\end{equation*}
is zero. Taking the pullback by $f$
and using the equality $\tau \circ f = f \circ \tau$
we get the result. \quod



\begin{prp}\label{nilpvanish}%
Let $\cM$ be an $\cO_X$-module shtuka over a $\tau$-scheme $X$. If $\cM$ is
nilpotent then $\RGamma(X,\,\cM) = 0$.\end{prp}

\pf 
Without loss of generality we assume that $\cM$ is given by a diagram
\begin{equation*}
\cM_0 \shtuka{1}{j} \cM_0.
\end{equation*}
As $\cM$ is nilpotent one easily deduces that the endomorphism of
$\RGamma(\cM_0)$ induced by $j$ is nilpotent. As a consequence the
endomorphism $1-j$ is a quasi-\hspace{0pt}isomorphism.
Now the canonical trinangle of $\RGamma(\cM)$ shows that $\RGamma(\cM)$
is the mapping fiber of $1 - j$ 
so the result follows.\quod
%
%

\breakflow
The following proposition is our main tool to deduce vanishing of cohomology.

\begin{prp}\label{nilpcomp} Let $R$ be a Noetherian $\tau$-ring complete with
respect to an ideal $I \subset R$. Assume that $\tau(I) \subset I$ so that
$\tau$ descends to the quotient $R/I$. Let $\cM$
be a locally free $R$-module
shtuka. If $\cM(R/I)$ is nilpotent then
the following holds:
\begin{enumerate}
\item $\RGamma(\cM) = 0$.

\item For every $n > 0$ the shtuka $\cM(R/I^n)$ is nilpotent.
\end{enumerate}
\end{prp}

\pf 
Suppose that $\cM$ is given by a diagram
\begin{equation*}
M_0 \shtuka{i}{j} M_1.
\end{equation*}
The ideal $I$ is in the Jacobson radical of $R$ so Nakayama's lemma
implies that $i$ is surjective. Now $i$ splits since $M_1$ 
is projective
and one more application of Nakayama's lemma shows that
$i$ is an isomorphism. 

The endomorphism $i^{-1} j$ of $M_0$ preserves the filtration by powers of $I$.
Furthermore $(i^{-1} j)^m M_0 \subset IM_0$ for some $m \geqslant 0$ since
$\cM(R/I)$ is nilpotent. As a consequence $(i^{-1} j)^{mn} M_0 \subset I^n M_0$
and we get (2). Moreover (2) implies that
$1 - i^{-1} j$ is an isomorphism modulo every power of $I$.
Since $M_0$ is $I$-adically complete we deduce that $1 - i^{-1} j$ is an
isomorphism. As $i$ is an isomorphism the claim (1) now follows from Theorem
\ref{shtaffcoh}.\quod

\section{\texorpdfstring{The linearization functor and $\zeta$-isomorphisms}{The Lie
functor and zeta-isomorphisms}}
\label{sec:shtlin}

The constructions of this section are due to V.~Lafforgue \cite{valeurs} but the
terminology and the notation is our own.
The notion of a $\zeta$-isomorphism is at the heart of our approach to the class
number formula.

\begin{dfn} Let $X$ be a $\tau$-scheme. We define
the \emph{linearization functor} 
$\Der$ from $\Sht\cO_X$ to $\Sht\cO_X$ as follows:
\begin{equation*}
\Der\Big[ \cM_0 \shtuka{i}{j} \cM_1 \Big] =
\Big[ \cM_0 \shtuka{i}{0} \cM_1 \Big].
\end{equation*}
We say that an $\cO_X$-module shtuka $[ \cM_0 \shtuka{i}{j} \cM_1 ]$ is
\emph{linear}\index{idx}{shtuka!linear} if $j = 0$.\end{dfn}

%
%

\breakflow
The complex $\RGamma(X,\Der\cM)$ is often easier to compute than
$\RGamma(X,\cM)$.
Even though the complexes
$\RGamma(X,\cM)$ and $\RGamma(X,\Der\cM)$ are very different in general,
a link beween them exists under some natural assumptions on $\cM$.

Fix a subring $A \subset \cO_X(X)^{\tau=1}$ and assume that
$\RGamma(X,\cM)$ is a perfect complex of $A$-modules.
The theory of Knudsen-Mumford \cite{kmdet}
associates to $\RGamma(X,\cM)$ an invertible $A$-module\footnote{Strictly speaking the
determinant is a pair $(L,\alpha)$ consisting of an invertible $A$-module $L$
and a continuous function $\alpha\colon \Spec A \to \bZ$. This function is not
important for the following discussion so we ignore it.} $\det_A \RGamma(X,\cM)$.
This determinant is functorial in quasi-\hspace{0pt}isomorphisms.
If the $A$-module complexes
$\RGamma(X,\,\cM)$, $\RGamma(X,\,\cM_0)$ and $\RGamma(X,\,\cM_1)$
are bounded and their cohomology modules are perfect 
then the canonical distinguished triangle
\begin{equation*}
\RGamma(X,\cM) \to \RGamma(X,\cM_0) \xrightarrow{i-j} \RGamma(X,\cM_1) \to [1].
\end{equation*}
determines a natural $A$-module isomorphism
\begin{equation*}
\det\nolimits_A \RGamma(X,\,\cM) \xrightarrow{\isosign}
\det\nolimits_A \RGamma(X,\,\cM_0) \otimes_A \det\nolimits_A^{-1} \RGamma(X,\,\cM_1)
\end{equation*}
\cite[Corollary 2 after Theorem 2]{kmdet}.

\begin{dfn}\label{defzeta}\index{idx}{shtuka cohomology!$\zeta$-isomorphism}%
Let $X$ be a $\tau$-scheme and let $A \subset \cO_X(X)^{\tau=1}$ be
a subring. Let $\cM$ be an $\cO_X$-module shtuka given by a diagram
\begin{equation*}
\cM_0 \shtuka{\,\,i\,\,}{\,\,j\,\,} \cM_1.
\end{equation*}
We say that the \emph{$\zeta$-isomorphism is defined for $\cM$}
if the complexes of $A$-modules
$\RGamma(X,\,\cM)$, $\RGamma(X,\,\Der\cM)$,
$\RGamma(X,\,\cM_0)$ and $\RGamma(X,\,\cM_1)$
are bounded and their cohomology modules are perfect.
Under this assumption
we define the \emph{$\zeta$-isomorphism}
\begin{equation*}
\zeta_\cM\colon \det\nolimits_A \RGamma(X,\cM) \xrightarrow{\isosign}
\det\nolimits_A\RGamma(X,\Der\cM)
\end{equation*}
as the composition
\begin{equation*}
\det\nolimits_A^{\phantom{1}} \!\RGamma(X,\cM) \xrightarrow{\isosign}
\det\nolimits_A^{\phantom{1}} \!\RGamma(X,\cM_0) \otimes_A^{\phantom{1}}
\det\nolimits_A^{-1} \!\RGamma(X,\cM_1) \xleftarrow{\isosign}
\det\nolimits_A^{\phantom{1}} \!\RGamma(X,\Der\!\cM)
\end{equation*}
of isomorphisms determined by the canonical triangles
\begin{align*}
\RGamma(X,\cM) &\to \RGamma(X,\cM_0) \xrightarrow{i-j} \RGamma(X,\cM_1) \to [1],
\\
\RGamma(X,\Der\cM) &\to \RGamma(X,\cM_0) \xrightarrow{\,\,\,i\,\,\,\,} \RGamma(X,\cM_1) \to [1].
\end{align*}
\end{dfn}


\section{\texorpdfstring{$\tau$-polynomials}{Tau-polynomials}}
\label{sec:taupoly}

In this section we work with a fixed $\tau$-ring $R$.

\begin{dfn} We define the ring $R\{\tau\}$ as follows.\index{idx}{$\tau$!the ring of $\tau$-polynomials}\index{nidx}{$R\{\tau\}$, the ring of $\tau$-polynomials}
Its elements are formal polynomials
$r_0 + r_1 \tau + r_2 \tau^2 + \dotsc r_n \tau^n$,
$r_0, \dotsc r_n \in R$, $n \geqslant 0$.
The multiplication in $R\{\tau\}$ is subject to the following identity: for
every $r \in R$ we have
\begin{equation*}
\tau \cdot r = r^\tau \cdot \tau
\end{equation*}
where $r^\tau = \tau(r)$ is the image of $r$ under $\tau\colon R \to
R$.\end{dfn}

\breakflow
Unlike all the other rings in this text the ring $R\{\tau\}$ is not commutative
in general. Still it is associative and has the multiplicative unit $1$. Left
$R\{\tau\}$-modules are directly related to $R$-module shtukas.

\begin{dfn}\label{rmodassoc}Let $M$ be a left $R\{\tau\}$-module. The
\emph{$R$-module shtuka associated to $M$} is
\begin{equation*}
M \shtuka{1}{\tau} M.
\end{equation*}
Here $\tau\colon M \to M$ is the $\tau$-multiplication map. It is tautologically
$\tau$-linear so that the diagram above indeed defines a shtuka.\end{dfn}

%

\breakflow
%
Let $M$ be a left $R\{\tau\}$-module. In this section we
denote $a\colon R\{\tau\} \otimes_R M \to M$ the map which sends a tensor $\varphi
\otimes m$ to $\varphi\cdot m$.  The letter $a$ stands for ``action''.

\begin{lem}\label{rtaumodresb} If $M$ is a left $R\{\tau\}$-module then the
sequence of left $R\{\tau\}$-modules
\begin{equation*}
0 \to R\{\tau\}\tau \otimes_R M \xrightarrow{\,\,d\,\,}
R\{\tau\} \otimes_R M \xrightarrow{\,\,a\,\,} M \to 0
\end{equation*}
is exact. Here
$d(\varphi \tau \otimes m) = \varphi \otimes \tau \cdot m - \varphi \tau \otimes m$. \end{lem}

\pf It is clear that $a \circ d = 0$ and $a$ is surjective.
Let us verify the injectivity of $d$.
The modules $R\{\tau\}$, $R\{\tau\}\tau$ carry filtrations by degree of
$\tau$-polynomials. The map $d$ is compatible with the induced filtrations on $R\{\tau\}\tau\otimes_R M$ and
$R\{\tau\}\otimes_R M$ and is injective on subquotients. It is therefore injecitve.

Let us verify the exactness of the sequence at $R\{\tau\}\otimes_R M$. Consider
the quotient of $R \{\tau\}\otimes_R M$ by the image of $d$.
In this quotient we have the identity
$r \tau^{n+1} \otimes m \equiv r \tau^n \otimes \tau m$
for all $r \in R$, $m \in M$, $n \geqslant 0$.
As a consequence $\varphi \otimes m
\equiv \varphi \cdot m$ for every $\varphi \in R\{\tau\}$.
Hence every element $y \in R\{\tau\}\otimes_R M$ is equivalent to $1
\otimes a(y)$. If $a(y) = 0$ then $y \equiv 0$
or in other words $y$ is in the image of $d$.
\quod

\begin{nrmk}\label{taumod} Let $M$ be an $R$-module. We denote $\tau^\ast M$ the $R$-module
$R^\tau \otimes_R M$
where $R^\tau$ is $R$ with the $R$-algebra structure given by the homomorphism $\tau\colon
R \to R$. 
We write the elements of $\tau^\ast M$ as sums of pure tensors $r \otimes m$,
$r \in R^\tau$, $m \in M$. If $r, r_1 \in R$ and $m \in M$ then
\begin{equation*}
r \otimes r_1 m = r \tau(r_1) \otimes m.
\end{equation*}
The ring $R$ acts on $\tau^\ast M = R^\tau \otimes_R M$ via the factor $R^\tau =
R$. If $r, r_1 \in R$ and $m \in M$ then
\begin{equation*}
r_1 \cdot (r \otimes m) = r_1 r \otimes m.
\end{equation*}\end{nrmk}

\begin{lem}\label{rtaumodpulldesc} Let $M$ be an $R$-module.
The maps
\begin{equation*}
R\{\tau\}\tau\otimes_R M \to R\{\tau\} \otimes_R \tau^\ast M, \quad
\varphi\tau \otimes m \mapsto \varphi \otimes (1 \otimes m)
\end{equation*}
and
\begin{equation*}
R\{\tau\} \otimes_R \tau^\ast M \to R\{\tau\}\tau\otimes_R M, \quad
\varphi \otimes (r \otimes m) \mapsto \varphi r\tau \otimes m
\end{equation*}
are mutually inverse isomorphisms of left $R\{\tau\}$-modules.
%
%
\quod\end{lem}

%

\breakflow
\begin{prp}\label{rtaumodres} If $M$ is a left $R\{\tau\}$-module then the
sequence of left $R\{\tau\}$-modules
\begin{equation*}
0 \to R\{\tau\}\otimes_R \tau^\ast M \xrightarrow{\,\,1\otimes \tau^a - \eta\,\,}
R\{\tau\} \otimes_R M \xrightarrow{\,\,a\,\,} M \to 0
\end{equation*}
is exact. Here $\tau^a\colon \tau^\ast M \to M$ is the adjoint of the
$\tau$-multiplication map $M \to \tau_\ast M$ and $\eta$ is the map given by the
formula
\begin{equation*}
\eta\colon \varphi \otimes (r\otimes m) \mapsto \varphi r \tau \otimes m.
\end{equation*}
\end{prp}

\pf Using the isomorphism $R\{\tau\}\tau\otimes_R M \cong R\{\tau\}\otimes_R
\tau^\ast M$ of Lemma \ref{rtaumodpulldesc} we rewrite the short sequence
in question as
\begin{equation*}
0 \to R\{\tau\}\tau\otimes_R M \xrightarrow{\,\,d\,\,}
R\{\tau\} \otimes_R M \xrightarrow{\,\,a\,\,} M \to 0.
\end{equation*}
An easy computation shows that
\begin{equation*}
d(\varphi \tau \otimes m) = \varphi \otimes \tau \cdot m - \varphi \tau \otimes m.
\end{equation*}
The result thus follows from Lemma \ref{rtaumodresb}. \quod

\section{The Hom shtuka in an affine setting}
\label{sec:genhomsht}

Let $R$ be a $\tau$-ring and let $M$ and $N$ be $R$-module shtukas.
In this section we introduce the Hom shtuka $\iHom_R(M,N)$.
To some extent it behaves
like an internal Hom in the category of $R$-module shtukas. It is literally the internal
Hom for shtukas which come from left $R\{\tau\}$-modules.
Even if both $M$ and $N$ are left $R\{\tau\}$-modules, $\iHom_R(M,N)$ is in
general a genuine shtuka which does not come from a left $R\{\tau\}$-module.
Apart from the Drinfeld construction of Chapter \ref{ch:drinconstr} the $\iHom$
construction is the main source of nontrivial shtukas in the present work.

\begin{dfn}\label{defhomsht}\index{idx}{shtuka!Hom}\index{nidx}{shtukas!$\iHom$}%
Let $M$ and $N$ be $R$-module shtukas given by diagrams
\begin{equation*}
M = \Big[ M_0 \shtuka{\,\,i_M\,\,}{\,\,j_M\,\,} M_1 \Big], \quad
N = \Big[ N_0 \shtuka{\,\,i_N\,\,}{\,\,j_N\,\,} N_1 \Big].
\end{equation*}
The \emph{Hom shtuka} $\iHom_R(M,N)$ is defined by the diagram
\begin{equation*}
\Hom_R(M_1, N_0) \shtuka{\quad\,\, i\,\,\quad}{j} \Hom_R(\tau^\ast M_0, N_1)
\end{equation*}
where
\begin{align*}
i(f) &= i_N \circ f \circ j_M^a, \\
j(f) &= (j_N \circ f)^a \circ \tau^\ast(i_M),
\end{align*}
$j^a_M$ is the $\tau$-adjoint of $j_M$
and
$(-)^a$ denotes the adjunction
$\Hom_R(M_1, \tau_\ast N_1) = \Hom_R(\tau^\ast M_1, N_1)$.\end{dfn}

\breakflow
The adjunction $\Hom_R(M_1, \tau_\ast N_1) = \Hom_R(\tau^\ast M_1, N_1)$
is $\tau$-linear
so the diagram above indeed defines a shtuka.

If $M$ and $N$ are left $R\{\tau\}$-modules then $\iHom_R(M,N)$ means
$\iHom_R$ applied to the $R$-module shtukas associated to $M$ and $N$ as in
Definition~\ref{rmodassoc}. The Hom shtukas we work with are typically of this sort.
We will also need $\iHom_R(M,N)$ in the case when $M$ is an $R$-module
shtuka which does not come from a left $R\{\tau\}$-module
(cf. Section~\lref{locmod}{sec:deflocmod}).

In the rest of this section we use the notation of Remark \ref{taumod}
for the elements of $\tau^\ast M$.

\begin{lem}\label{rtaumodhomsht}%
If $M$ and $N$ are left $R\{\tau\}$-modules then
$\iHom_R(M,N)$ is represented by the diagram
\begin{gather*}
\Hom_R(M,N) \shtuka{\quad i\quad}{j} \Hom_R(\tau^\ast M,N), \\
i(f) = f \circ \tau_M^a, \quad
j(f) = \tau_N^a \circ \tau^\ast(f)
\end{gather*}
where
$\tau_M^a$ and $\tau_N^a$ 
are the adjoints of the $\tau$-multiplication maps. In other words
\begin{equation*}
i(f)\colon r \otimes m \mapsto f(r \tau \cdot m), \quad
j(f)\colon r \otimes m \mapsto r \tau \cdot f(m).\quod
\end{equation*}\end{lem}


\breakflow
Our goal is to describe the cohomology of the shtuka $\iHom_R(M,N)$
in the case when $M$ and $N$ are left $R\{\tau\}$-modules.

\begin{prp}\label{ihomzerocoh}%
If $M$ and $N$ are left $R\{\tau\}$-modules then
\begin{equation*}
\Hom_{R\{\tau\}}(M,N) =
\uH^0(\iHom_R(M,N))
\end{equation*}
as abelian subgroups of $\Hom_R(M,N)$.%
\end{prp}

\pf Let $i$ and $j$ be the arrows of $\iHom_R(M,N)$. Let $f \in \Hom_R(M,N)$.
Lemma \ref{rtaumodhomsht} implies that $i(f) = j(f)$ if and only if $f$ commutes
with $\tau$. \quod

\breakflow
If $M$ is projective as an $R$-module then the functor
$N \mapsto \iHom_R(M,N)$
on the category of left $R\{\tau\}$-modules
is exact.
So the isomorphism of Proposition~\ref{ihomzerocoh}
induces a natural map
$\RHom_{R\{\tau\}}(M,N) \to \RGamma(\iHom_R(M,N))$.

\begin{thm}\label{homshtcohrhom}%
Let $M$ and $N$ be left $R\{\tau\}$-modules.
If $M$ is projective as an $R$-module then the natural map
$\RHom_{R\{\tau\}}(M,N) \to \RGamma(\iHom_R(M,N))$
is a quasi-\hspace{0pt}isomorphism
\end{thm}

\pf By Proposition \ref{rtaumodres} we have a short exact sequence
\begin{equation}\label{rtmres}
0 \to R\{\tau\}\otimes_R \tau^\ast M
\xrightarrow{\,\,1 \otimes \tau^a_M - \eta\,\,}
R\{\tau\}\otimes_R M \xrightarrow{\,\,a\,\,} M \to 0.
\end{equation}
If $M$ is a projective $R$-module then so is $\tau^\ast M$. As a consequence
$R\{\tau\}\otimes_R M$ and $R\{\tau\}\otimes_R \tau^\ast M$ are projective left
$R\{\tau\}$-modules. Thus \eqref{rtmres} is a projective resolution of $M$ as a
left $R\{\tau\}$-module. Applying $\Hom_{R\{\tau\}}(-,N)$ to (3) we conclude that
\begin{equation}\label{rtmcompl}
\RHom_{R\{\tau\}}(M,N) =
\Big[ \Hom_R(M,N) \xrightarrow{(1\otimes\tau^a_M)^*-\eta^*} \Hom_R(\tau^\ast M,N) \Big].
\end{equation}
where * indicates the induced maps.
Consider the shtuka
\begin{equation*}
\iHom_R(M,N) = \Big[ \Hom_R(M,N) \shtuka{i}{j} \Hom_R(\tau^\ast M, N) \Big].
\end{equation*}
According to Lemma \ref{rtaumodhomsht} the maps $i$ and $(1 \otimes \tau^a_M)^*$ coincide. 
A straightforward computation shows that $\eta^* = j$.
Therefore \eqref{rtmcompl} computes $\RGamma(\iHom_R(M,N))$ by Theorem \ref{shtaffcoh}. \quod

\section{A global variant of the Hom sthuka}
\label{sec:globhomsht}

Let $X$ be a $\tau$-scheme and let $\cM$ and $\cN$ be $\cO_X$-module shtukas.
In this section define the $\cO_X$-module Hom shtuka $\siHom_X(\cM,\cN)$,
a global variant of the construction of the previous section.
This construction is important to Chapter~\ref{chapter:cnf}.

Given $\cO_X$-modules $\cE$ and $\cF$ we denote $\siHom_X(\cE,\cF)$
the sheaf Hom [\stacks{01CM}].

\begin{dfn}\label{defglobhomsht}\index{idx}{shtuka!Hom}\index{nidx}{shtukas!$\siHom$}%
Let $\cM$ and $\cN$ be $\cO_X$-module shtukas given by diagrams
\begin{equation*}
\cM = \Big[ \cM_0 \shtuka{\,\,i_\cM\,\,}{\,\,j_\cM\,\,} \cM_1 \Big], \quad
\cN = \Big[ \cN_0 \shtuka{\,\,i_\cN\,\,}{\,\,j_\cN\,\,} \cN_1 \Big].
\end{equation*}
The \emph{Hom shtuka} $\siHom_X(\cM,\cN)$ is defined by the diagram
\begin{equation*}
\siHom_X(\cM_1, \cN_0) \shtuka{\quad\,\, i\,\,\quad}{j} \siHom_X(\tau^\ast \cM_0, \cN_1).
\end{equation*}
The arrow $i$ is the compostion
\begin{equation*}
\siHom_X(\cM_1, \cN_0) \xrightarrow{\,\,i_\cN \circ - \,\,}
\siHom_X(\cM_1, \cN_1) \xrightarrow{\,\, - \circ j_\cM^a\,\,}
\siHom_X(\tau^\ast\cM_0, \cN_1)
\end{equation*}
where $j^a_\cM$ is the $\tau$-adjoint of $j_\cM$.
The arrow $j$ is defined by the diagram
\begin{equation*}
\xymatrix{
\siHom_X(\cM_1, \cN_0) \ar[d]_{j_\cN \circ -} \ar[rr]^j && \tau_\ast\siHom_X(\tau^\ast\cM_0, \cN_1) \\
\siHom_X(\cM_1, \tau_\ast \cN_1) \ar[rr]^{\isosign} &&
\tau_\ast\siHom_X(\tau^\ast\cM_1, \cN_1), \ar[u]_{\tau_\ast(h)}
}
\end{equation*}
the unlabelled arrow is the natural adjunction isomorphism
and $h$ is given by composition with $\tau^\ast(i_\cM)\colon \tau^\ast\cM_0 \to \tau^\ast\cM_1$.
\end{dfn}

\begin{lem}%
Let $X = \Spec R$ be an affine $\tau$-scheme and let $M$ and $N$
be $R$-module shtukas. Denote $\cM$ and $\cN$ the $\cO_X$-module
shtukas corresponding to $M$ and $N$ respectively.
If the underlying $R$-modules of $M$ are of finite type
then $\siHom_X(\cM,\cN)$ is the $\cO_X$-module shtuka
associated to $\iHom_R(M,N)$.\quod\end{lem}

\section{\texorpdfstring{Special $\tau$-structure}{Special tau-structure}}
\label{sec:globcoh}

The last section of this chapter is devoted to
a technical result which will only be used
in Section~\ref{sec:globartreg}.
Let $X$ be a $\tau$-scheme.
We assume that $\tau$
acts as the identity on the underlying topological space. The most important case
for our applications is $X = \Spec\ngC\times_\Fq X_0$ where $X_0$ is an $\Fq$-scheme,
$\ngC$ a finite artinian $\Fq$-algebra and $\tau\colon X \to X$ the endomorphism
which acts as the identity on $\ngC$ and as the $q$-Frobenius on $X_0$.

Let $\cF$ be an $\cO_X$-module. Since $\tau$ acts as the identity on the underlying
topological space we can identify $\cF$ and $\tau_\ast\cF$ as sheaves of abelian groups.
We can thus make the following construction.
Let $\cM$ be a shtuka on $X$ given by a diagram
\begin{equation*}
\cM_0 \shtuka{\,\,i\,\,}{j} \cM_1.
\end{equation*}
We associate to it a complex of $\cO_X(X)^{\tau=1}$-module sheaves
\begin{equation}\label{artshcomplex}\tag{$\ast$}
\big[ \cM_0 \xrightarrow{\,\,i-j\,\,} \cM_1 \big]
\end{equation}
on the underlying topological space of $X$.
We identify $\tau_\ast\cM_1$ and $\cM_1$ as sheaves 
and view $j\colon\cM_0 \to \tau_\ast\cM_1$ as a morphism from $\cM_0$ to $\cM_1$.

Our goal is to show that the sheaf cohomology of \eqref{artshcomplex}
coincides with the shtuka cohomology of $\cM$.
It is useful to first extend the definition of
\eqref{artshcomplex} to complexes of shtukas.

\begin{dfn}Let $\cM$ be a complex of $\cO_X$-module shtukas.
We define a complex of $\cO_X(X)^{\tau=1}$-module sheaves $\cG_a(\cM)$
as the total complex of the double complex
\begin{equation*}
\xymatrix{
&& \\
\alpha_\ast\cM^{n+1} \ar[rr]^{(-1)^{n+1}(i-j)} \ar[u] &&
\beta_\ast\cM^{n+1} \ar[u] \\
\alpha_\ast\cM^{n} \ar[rr]^{(-1)^{n}(i-j)} \ar[u] &&
\beta_\ast\cM^{n} \ar[u] \\
\ar[u] && \ar[u]
}
\end{equation*}
The vertical maps are the differentials of $\alpha_\ast\cM$
respectively $\beta_\ast\cM$. The object $\alpha_\ast\cM^n$
is placed in the bidegree $(n,0)$ while $\beta_\ast\cM^n$ is in the
bidegree $(n,1)$.\end{dfn}

\begin{dfn}\label{artshtri}%
We define a natural triangle
\begin{equation*}
\cG_a(\cM) \xrightarrow{\,\,p\,\,}
\alpha_\ast \cM \xrightarrow{\,i-j\,}
\beta_\ast\cM \xrightarrow{\,-q\,}
[1]
\end{equation*}
as follows. The object of $\cG_a(\cM)$ in degree $n$ is
$\alpha_\ast\cM^n \oplus \beta_\ast\cM^{n-1}$.
The morphism $p$ is the projection to the first factor.
The morphism $q$ is the injection to the second factor
multiplied by $(-1)^n$.\end{dfn}

\begin{example*}%
If $\cM$ is a single shtuka placed in degree $0$ then the triangle above is
the canonical distinguished triangle
\begin{equation*}
\big[\cM_0 \xrightarrow{\,i-j\,} \cM_1\big] \xrightarrow{\phantom{\,i-j\,}}
\cM_0 \xrightarrow{\,i-j\,} \cM_1 \xrightarrow{\phantom{\,i-j\,}} [1]
\end{equation*}
of the mapping fiber complex.\end{example*}

\begin{lem}The triangle of Definition \ref{artshtri} is distinguished.\end{lem}

\pf The sequence
\begin{equation*}
0 \to \beta_\ast\cM[-1] \xrightarrow{\,q[-1]\,}
\cG_a(\cM) \xrightarrow{\,\,\,\,p\,\,\,\,}
\alpha_\ast\cM \to 0
\end{equation*}
is exact and is termwise split. Such a sequence determines a distinguished
triangle in the following way. 
Let $r$ be the splitting of $q[-1]$ 
given by the formula $(a,b)
\mapsto (-1)^n b$ in degree $n+1$ and let $s$ be the splitting of $p$ 
given by the formula $a \mapsto (a,0)$.
Let $f = r \circ d \circ s$ where $d$ is the differential of
$\cG_a(\cM)$. The triangle
\begin{equation*}
\beta_\ast\cM[-1] \xrightarrow{q[-1]}
\cG_a(\cM) \xrightarrow{\,\,p\,\,}
\alpha_\ast\cM \xrightarrow{\,\,f\,\,}
\beta_\ast\cM
\end{equation*}
is distinguished [\stacks{014Q}]. An easy computation reveals that
$f = i - j$. Rotating this triangle we conclude that \eqref{assoccompltri}
is distinguished.\quod


\breakflow
Lemma \ref{assoccompldesc} implies that
$\cGamma(X,\,\cM) = \Gamma(X,\,\cG_a(\cM))$.
Moreover
the functor $\cG_a$ is exact on the level of homotopy categories.
So taking the derived
functors we obtain a canonical morphism
$\RGamma(X, \,\cM) \to
\RGamma(X, \,\cG_a(\cM))$.

\begin{thm}\label{shtglobartcoh}%
Assume that $\tau\colon X \to X$
acts as the identity on the underlying topological space.
For every complex of $\cO_X$-module shtukas $\cM$
the natural map
$\RGamma(X, \,\cM) \to \RGamma(X,\,\cG_a(\cM))$
is a quasi-\hspace{0pt}isomorphism. 
Furthermore it extends to an isomorphism
of distinguished triangles
\begin{equation*}
\xymatrix{
\RGamma(X, \,\cM) \ar[r] \ar[d]^{\rtviso}
& \RGamma(X, \,\alpha_\ast\cM) \ar[r]^{i-j} \ar@{=}[d]
& \RGamma(X, \,\beta_\ast\cM) \ar[r] \ar@{=}[d]
& [1] \\
\RGamma(X, \,\cG_a(\cM)) \ar[r] 
& \RGamma(X, \,\alpha_\ast\cM) \ar[r]^{i-j}
& \RGamma(X, \,\beta_\ast\cM) \ar[r]
& [1] \\
}
\end{equation*}
Here the top row is the canonical triangle
and the bottom row is obtained by applying $\RGamma(X,-)$
to the triangle of Definition \ref{artshtri}.\end{thm}

\breakflow
In a way this result justifies the use
of notation $\RGamma$ for shtuka cohomology. With some caution one may
think of a shtuka as a two term complex \eqref{artshcomplex}.
Shtuka cohomology is then the usual sheaf cohomology of this complex.


\afterall\noindent\textit{Proof of Theorem \ref{shtglobartcoh}. }%
Without loss of generality we can assume that the complex $\cM$ is K-injective.
Consider the diagram
\begin{equation*}
\xymatrix{
\Gamma(X,\,\cG_a(\cM)) \ar[r] \ar[d] & \Gamma(X,\,\alpha_\ast\cM) \ar[r]^{i-j} \ar[d] &
\Gamma(X,\,\beta_\ast\cM) \ar[r] \ar[d] & [1] \\
\RGamma(X,\,\cG_a(\cM)) \ar[r] & \RGamma(X,\,\alpha_\ast\cM) \ar[r]^{i-j} &
\RGamma(X,\,\beta_\ast\cM) \ar[r] & [1]
}
\end{equation*}
where the rows are obtained by applying $\Gamma(X,-)$ respectively
$\RGamma(X,-)$ to the distinguished triangle of Definition \ref{artshtri}
and the vertical arrows are the natural morphisms $\Gamma(X,-)\to\RGamma(X,-)$.
By construction the top row coincides with the distinguished triangle
\begin{equation*}
\cGamma(X,\,\cM) \to \Gamma(X,\,\alpha_\ast\cM)\xrightarrow{i-j} \Gamma(X,\,\beta_\ast\cM) \to [1]
\end{equation*}
of Definition \ref{defassoctri}. The diagram above is a morphism of distinguished
triangles by naturality. 
Now $\cM$ is K-injective so Proposition~\ref{assocmor}
identifies the top row with the canonical distinguished
triangle
\begin{equation*}
\RGamma(X,\,\cM) \to \RGamma(X,\,\alpha_\ast\cM)\xrightarrow{i-j} \RGamma(X,\,\beta_\ast\cM) \to [1].
\end{equation*}
and the second and third vertical arrows with the identity morphisms.
It follows that
the first vertical arrow is a quasi-\hspace{0pt}isomorphism.\quod

\chapter{Topological vector spaces over finite fields}
\label{chapter:tvs}
\label{ch:pont}
\label{ch:veccompl}
\label{ch:complot}
\label{ch:indcot}
\label{ch:contfn}
\label{ch:bndfn}
\label{ch:bndlc}
\label{ch:germ}

In this chapter we present some results on topological vector spaces over
finite fields.
%
%
%
The base field $\Fq$ is fixed throughout the chapter.
It is assumed to carry the discrete topology.
In the following a subspace of a vector space always means an $\Fq$-vector
subspace, not an arbitrary topological subspace. 
As is usual in the theory of topological groups all our locally compact
topological vector spaces are assumed to be Hausdorff.
A topological vector space is said to be \emph{linearly topologized} if every
open neighbourhood of zero contains an open subspace.
We mainly study linearly topologized Hausdorff vector spaces.
Throughout this chapter we abbreviate ``linearly topologized Hausdorff'' as
``\emph{lth}''.

%

\section{Overview}

In our computations of shtuka cohomology we will extensively use various spaces
of continuous $\Fq$-linear functions and germs of such functions. This chapter is
devoted to their study.

Let $V$ and $W$ be locally compact topological $\Fq$-vector spaces. 
We consider the following function spaces:
\begin{itemize}
\item the space
$c(V,W)$ of continuous $\Fq$-linear maps from $V$ to $W$,

\item the space $b(V,W)$ of bounded continuous
$\Fq$-linear maps, i.e. the maps which have image in a compact subspace,

\item the space $a(V,W)$ of locally constant bounded $\Fq$-linear maps.
\end{itemize}
The function spaces are equipped with topologies which make them into complete
vector spaces. These topologies are carefully chosen to suit the 
applications.
The space $c(V,W)$ carries the compact-open topology. The topologies on its
dense subspaces $a(V,W)$ and $b(V,W)$ are finer than the induced ones.

An important object related to the function spaces is the space of germs $g(V,W)$.
Its elements are equivalence classes of continuous $\Fq$-linear maps from $V$ to
$W$. Two such maps are equivalent if they restrict to the same map on an
open subspace.
The main property of $g(V,W)$ is invariance under local isomorphisms on the
source $V$ and the target $W$.
This property will be indispensable for cohomological computations of
Chapter
\ref{ch:infbasecoh} among others.

To describe the structure of the function spaces we employ two topological
tensor products:
\begin{itemize}
\item the completed tensor product
$V \complot W = \lim_{U, Y} V/U \otimes W/Y$.

\item the ind-complete tensor product
$V \indcot W = \lim_{U, Y} (V \otimes W) / (U \otimes Y)$.
\end{itemize}
Here $U \subset V$ and $Y \subset W$ run over all open subspaces.
The tensor products $\complot$ and $\indcot$
are closely related: we will show that the natural commutative square
\begin{equation*}
\xymatrix{
V \indcot W \ar[r] \ar[d] & V\complot W^\#\ar[d] \\
V^\# \complot W \ar[r] & V \complot W
}
\end{equation*}
is cartesian in the category of topological vector spaces
(Proposition \ref{indcotdescr}).
Here $(-)^\#$ denotes the same vector space taken with the discrete topology.

Using the tensor products above we will construct natural topological isomorphisms
\begin{equation*}
\begin{array}{r@{}ll}
c(V,W)\,\,&\cong V^* \complot W &\textrm{(Proposition \ref{contcomplot})}, \\
b(V,W)\,\,&\cong (V^*)^\#\complot W &\textrm{(Proposition \ref{bndcomplot})}, \\
a(V,W)\,\,&\cong V^* \indcot W &\textrm{(Proposition \ref{bndlcot})}
\end{array}
\end{equation*}
where $V^*$ is the continuous $\Fq$-linear dual of $V$.

Function spaces, germ spaces and the topological tensor products
$\complot$, $\indcot$ figure prominently in this text.
Past this chapter some degree of familiarity with
them will be assumed.

Much of the material in this chapter is largely well-known. However a reservation
should be made about the ind-complete tensor product $V \indcot W$, the function
spaces $a(V,W)$, $b(V,W)$ and the germ space $g(V,W)$. While these constructions
are very natural and should have certainly appeared before, we are not aware of
a reference for them in the literature.

Last but not least we should acknowledge our intellectual debt to G.~W.~Anderson.
The inspiration for this chapter
comes from his article \cite{and}, specifically from \S 2 of that text where he
uses function spaces to compute what in retrospect is the cohomology of certain
shtukas
associated to Drinfeld modules.

\section{Examples}

To lend this discussion more substance let us give some examples. We begin with
a few examples of locally compact vector spaces:
\begin{itemize}
\item $A = \Fq[t]$ with the discrete topology,

\item $F = \Fq(\!(t^{-1})\!)$ with the locally compact topology,

\item the compact open subspace $\cO_F = \Fq[[t^{-1}]]$ in $F$.
\end{itemize}
The quotient $F/A$ is naturally a compact $\Fq$-vector space.

Let $K$ be a local field containing $\Fq$. An example of a function space which
is particularly relevant to our study is
$c(F/A, K)$, the space of continuous $\Fq$-linear maps from $F/A$ to $K$.
Proposition \ref{contcomplot} provides us with a topological isomorphism
\begin{equation*}
(F/A)^* \complot K \cong c(F/A, K).
\end{equation*}
Let us show that
this isomorphism gives a rather hands-on description of $c(F/A,K)$. Let
$\Omega = \Fq[t]\,dt $ be the module of K\"ahler differentials of $A$ over $\Fq$ equipped with
the discrete topology and let
\begin{equation*}
\res\colon \Omega \otimes_A F \to \Fq, \quad
\sum_n a_n t^n \cdot dt \mapsto -a_{-1}
\end{equation*}
be the residue map at infinity.
The map
\begin{equation*}
\Omega \to (F/A)^*, \quad
\omega \mapsto [x \mapsto \res(x\omega)]
\end{equation*}
is easily shown to be
an isomorphism of topological vector spaces. As a consequence the isomorphism
$\Omega \complot K \cong c(F/A,K)$ of Proposition \ref{contcomplot} identifies
$c(F/A,K)$ with the topological vector space of formal series \`a la Tate:
\begin{equation*}
K\langle t \rangle\, dt = \Big\{\sum_{n \geqslant 0} \alpha_n t^n dt \in K[[t]] dt \,\,\Big|
\lim_{n \to \infty} \alpha_n = 0 \Big\}.
\end{equation*}
The topology on this space is induced by the
norm
\begin{equation*}
\Big|\sum_{n \geqslant 0} \alpha_n t^n \cdot dt \Big| = \sup_{n \geqslant 0} |\alpha_n|.
\end{equation*}
A series $\sum_{n\geqslant 0} \alpha_n t^n dt$ corresponds to the continuous
function
\begin{equation*}
F/A \to K, \quad
x \mapsto \sum_{n \geqslant 0} \alpha_n \res(x t^n dt).
\end{equation*}
This function maps $t^{-n}$ to $-\alpha_{n-1}$.

From this discussion one easily deduces that $g(F/A,K)$, the space of germs of
continuous functions from $F/A$ to $K$, is
isomorphic to the quotient
\begin{equation*}
\frac{\Omega \complot K}{\Omega \otimes K}.
\end{equation*}
Such quotients arise naturally in the cohomological computations of the
subsequent chapters. The utility of $g(F/A, K)$ is that it gives them an 
accessible interpretation.

\section{Basic properties}

\begin{lem}\label{opensplit}Every open embedding of topological vector spaces
is continuously split.\end{lem}

\pf Let $j\colon U \hookrightarrow V$ be an open embedding. The quotient
topology on $V/U$ is discrete. So every $\Fq$-linear splitting $i\colon V/U
\to V$ of the quotient map is continuous and the map $j\oplus i\colon U \oplus
V/U \to V$ is a continuous bijection. If $U' \subset U$ is an open subset and
$Y' \subset V/U$ any subset then the image of $U' \oplus Y'$ in $V$ is a union
of translates of $U'$ whence open. Thus $j\oplus i$ is a topological
isomorphism. \quod

\begin{cor}\label{openlift}Let $V, W$ be topological vector spaces, $U \subset
V$ an open subspace. Every continuous $\Fq$-linear map $f\colon U \to W$
admits an extension to $V$ with image $f(U)$.\end{cor}

\pf Take a splitting $V = U \oplus Y$ and extend $f$ to $Y$ by zero.
\quod

\begin{lem}\label{lthaus}A topological vector space is Hausdorff if
and only if its zero point is closed.\quod\end{lem}



\begin{lem} Let $V$ be an lth vector space, $V' \subset V$ a closed subspace.
The quotient topology on $V/V'$ is lth. \quod\end{lem}


\begin{lem} The category of lth vector spaces and continuous
$\Fq$-linear maps is additive and has arbitrary limits. The limits commute
with forgetful functors to $\Fq$-vector spaces (without topology) and to
topological abelian groups. \quod\end{lem}


\section{Locally compact vector spaces}

%
%
%
%
%
%

\begin{prp}\label{lclth}
Every locally compact vector space is lth and contains a compact open subspace.
\end{prp}

\pf%
Let $V$ be a locally compact vector space.
We first assume that $\Fq = \bF_p$ for a prime $p$.
In effect we work
with a locally compact $p$-torsion abelian group~$V$.
If $V$ is connected then its Pontrjagin dual is trivial since the $p$-torsion
subgroup of $\mathbb{C}^\times$ is disconnected. Hence $V$ is itself trivial.
The connected component of $0$ for a general $V$ is a connected locally compact
subgroup so it is trivial. Translating by elements of $V$ we conclude that $V$
is totally disconnected. Now a theorem of van~Dantzig
\cite[Ch. II, Theorem 7.7, p. 62]{hr}
states that every open neighbourhood of $0 \in V$ contains a compact open
subgroup.
 
Next let $\Fq$ be an arbitrary finite field.
If $U \subset V$ is an open subgroup then the subgroup
\begin{equation*}
U' = \bigcap_{\alpha \in \Fq^\times} \alpha U
\end{equation*}
is open 
and stable under the $\Fq$-action.\quod

\begin{lem}\label{compopen} Let $V$ be a locally compact vector space. For every
compact subset $K \subset V$ there exists a compact open subspace of $V$
containing $K$.\end{lem}

\pf Let $U \subset V$ be a compact open subspace which exists by Proposition
\ref{lclth}. The quotient $V/U$ is discrete so the image of $K$ in it is finite.
Let $K' \subset U/V$ be a finite $\Fq$-vector subspace containing the image of
$K$.  The preimage of $K'$ in $V$ is compact open and contains $K$ by
construction.\quod

%

\begin{dfn}\label{defdual}\index{nidx}{topological vector spaces and modules!$V^*$, continuous dual}%
For a locally compact vector space $V$ we define $V^*$ to be the
space of continuous $\Fq$-linear functions $V \to \Fq$ equipped with the
compact-open topology.\end{dfn}

\breakflow
In our context the duality theorem of Pontrjagin takes the following form:

\begin{thm}\label{pontdual} Let $V$ be a locally compact vector space.
\begin{enumerate}
\item $V^*$ is locally compact. Moreover:
\begin{enumerate}
\item $V$ is discrete if and only if $V^*$ is compact.

\item $V$ is compact if and only if $V^*$ is discrete.
\end{enumerate}

\item The natural map $V \to V^{**}$ is a topological isomorphism.
\end{enumerate}\end{thm}

\pf 
Let $\bF_{\!p} \subset \Fq$ be the prime subfield.
Using the trace map $\tr\colon \Fq \to \bF_{\!p}$ one easily reduces the problem to
the case $\Fq = \bF_{\!p}$.
In this case
we can invoke the usual Pontrjagin duality 
as follows. A choice of a primitive $p$-th root of unity determines a
topological isomorphism of $(\Fp,+)$ and the $p$-torsion subgroup
$\mu_p(\mathbb{C}) \subset \mathbb{C}^\times$. Every character of $V$ as
a locally compact $p$-torsion abelian group has the image in $\mu_p(\mathbb{C})$.
Thus the chosen isomorphism $(\Fp,+) \cong \mu_p(\mathbb{C})$ identifies $V^*$
with the Pontrjagin dual of $V$. The result is now clear. \quod

\section{Completion}\index{nidx}{topological vector spaces and modules!$\widehat{V}$, completion}

Let $V$ be an lth vector space. The completion of $V$ is the lth vector space
\begin{equation*}
\widehat{V} = \lim_U V/U
\end{equation*}
where $U \subset V$ runs over all open subspaces. It is enough to take the limit
over a fundamental system of open subspaces. $V$ is called complete if
the natural map $V \to \widehat{V}$ is a topological isomorphism. Every
continuous $\Fq$-linear map from $V$ to a complete lth vector space
factors uniquely over $\widehat{V}$. 

\begin{lem}\label{complopen} If $U$ is an open subspace in an lth
vector space $V$ then the
natural sequence
$0 \to \widehat{U} \to \widehat{V} \to V/U \to 0$
is exact. In particular the natural map $\widehat{U} \to \widehat{V}$ is an open
embedding. \quod\end{lem}


\begin{lem}\label{complmap} Let $V$ be an lth vector space.
\begin{enumerate}
\item $\widehat{V}$ is complete.

\item The natural map $V \to \widehat{V}$ is injective with dense image.

\item If $V \to \widehat{V}$ is bijective then it is a topological
isomorphism. \quod
\end{enumerate}\end{lem}



\begin{lem}A locally compact vector space is complete.\end{lem}

\pf Let $V$ be a compact vector space. The natural map $V \to \widehat{V}$ is
injective with dense image so a topological isomorphism. Thus a compact 
space is complete. If a space $V$ admits a complete open subspace then it is
complete. In particular every locally compact space is complete. \quod

%

\begin{prp}
\label{complcover}%
\textit{Let $\cU = \{U_i\}$ be a covering of an lth vector space
$V$ by open subspaces. If for every $U_i, U_j \in \cU$ there exists
$U_k \in \cU$ such that $U_i + U_j \subset U_k$ then $\{\widehat{U}_i\}$ covers
$\widehat{V}$.}\end{prp}

\begin{nrmk*}%
The draft version of this proposition was wrong.
Many thanks to Hendrik Lenstra for the correction.\end{nrmk*}

\afterall\noindent\textit{Proof of Proposition~\ref{complcover}.}
According to Lemma \ref{complopen} the natural maps $\widehat{U}_i \to
\widehat{V}$ are open embeddings. Let $W$ be the union of $\widehat{U}_i$
inside $\widehat{V}$. 
The assumption on $\cU$ implies that $W$ is an $\Fq$-vector subspace.
By construction $W$ contains $V$ and so is dense.
Being an open $\Fq$-vector subspace it is automatically closed hence coincides
with $\widehat{V}$.\quod

\section{Completed tensor product}
\label{sec:complot}

Recall that according to our convention a tensor product $\otimes$ without
subscript means a tensor product over $\Fq$.

\begin{dfn}\index{idx}{topological tensor product!tensor product topology}%
\index{nidx}{topological tensor products!$\otimes_{\textup{c}}$, tensor product topology}%
Let $V, W$ be lth vector spaces. We define the \emph{tensor product topology} on
$V \otimes W$ by the fundamental system of subspaces $U \otimes W + V \otimes
Y$ where $U \subset V, Y \subset W$ run over all open subspaces. The
$\Fq$-vector space $V \otimes W$ equipped with this topology is denoted $V
\otimes_{\textup{c}} W$. We reserve the tensor product $V \otimes W$ without
decorations to indicate the corresponding vector space without
topology.\end{dfn}

\begin{rmk} In general a continuous bilinear map $U \times V \to W$ does not 
induce a continuous map $U \otimes_{\textup{c}} V \to W$.\end{rmk}

\begin{lem}\label{otlth} If $V$ and $W$ are lth vector spaces then $V
\otimes_{\textup{c}} W$ is lth.\end{lem}

\pf We need to prove that $V \otimes_{\textup{c}} W$ is Hausdorff. According to
Lemma \ref{lthaus} it suffices to prove that $0 \in V \otimes_{\textup{c}} W$ is
closed. Assume $W$ is discrete. If $U \subset V$ is an open subspace then $U
\otimes W \subset V \otimes_{\textup{c}} W$ is open and hence closed. Letting
$U$ run over all open subspaces of $V$ we conclude that $0 = \bigcap_U U \otimes
W$ is closed. Now let $W$ be arbitrary. Fix an open subspace $Y \subset W$. The
natural map $V \otimes_{\textup{c}} W \to V \otimes_{\textup{c}} W/Y$ is
continuous. As the latter space is Hausdorff it follows that $V \otimes Y$ is
closed. The intersection $\bigcap_Y V \otimes Y = 0$ is thus closed.\quod

\begin{dfn}\label{defcomplot}\index{idx}{topological tensor product!completed}\index{nidx}{topological tensor products!$\complot$, completed}%
Let $V, W$ be lth vector spaces. We define the \emph{completed tensor product}
$V\complot W$ as the completion of $V \otimes_{\textup{c}} W$. In other words
\begin{equation*}
V \complot W = \lim_{U, Y} V/U \otimes W/Y
\end{equation*}
where $U \subset V$, $Y \subset W$ run over all open subspaces and the tensor
products in the limit diagram are equipped with the discrete topology. If $f\colon V_1
\to V_2$ and $g\colon W_1 \to W_2$ are continuous $\Fq$-linear maps then $f
\complot g\colon V_1 \complot W_1 \to V_2 \complot W_2$ is defined as the
completion of $f \otimes g$.\end{dfn}


\begin{rmk}A completed tensor product of two compact spaces is compact.
However a completed tensor product of an infinite discrete and an infinite
compact space is never locally compact.\end{rmk}

\begin{dfn}\index{nidx}{topological vector spaces and modules!$V^\#$, discretization}%
For a vector space $V$ we define $V^\#$ to be this space
equipped with the discrete topology.\end{dfn}

\begin{prp}\label{discrcomplot} Let $V, W$ be lth vector spaces. Consider the
natural map
\begin{equation*}
\iota\colon V^\#\complot W \to V \complot W.
\end{equation*}
\begin{enumerate}
\item The map $\iota$ is injective with dense image.

\item If $V$ is complete and $W$ compact then $\iota$ is a bijection.\end{enumerate}\end{prp}

\pf For every open subspace $Y \subset W$ let $\iota_Y\colon V
\otimes_{\textup{c}} W/Y \to \lim_U (V/U \otimes_{\textup{c}} W/Y)$ be
the completion map.
At the level of $\Fq$-vector spaces without topology
the map $\iota$ is the limit of $\iota_Y$ over all open
subspaces $Y \subset W$.

(1) The space $V \otimes_{\textup{c}} W/Y$ is Hausdorff by Lemma \ref{otlth}.
Hence $\iota_Y$ is injective by Lemma \ref{complmap}. It follows that $\iota$ is
injective.
The density statement is clear.

(2) The space $W/Y$ is finite since $W$ is compact.
Hence the natural map
\begin{equation*}
\lim_U (V/U \otimes W/Y) \to \lim_U(V/U) \otimes W/Y
\end{equation*}
is an isomorphism of $\Fq$-vector spaces without topology. Since $V$ is complete we
conclude that $\iota_Y$ is bijective. As a consequence
$\iota$ is a bijection.
\quod

\section{Ind-complete tensor product}
\label{sec:indcot}


In this section we introduce a different topology on $V \otimes W$ which is
better for some purposes than the usual tensor product topology.

\begin{dfn}\label{indcotdef}\index{idx}{topological tensor product!ind-tensor product topology}%
\index{nidx}{topological tensor products!$\otimes_{\textup{ic}}$, ind-tensor product topology}%
Let $V, W$ be lth vector spaces. We define the \emph{ind-tensor product
topology} on $V \otimes W$ by the fundamental system of open subspaces $U
\otimes Y$ where $U \subset V$, $Y \subset W$ run over all open subspaces. We
denote $V \otimes_{\textup{ic}} W$ the tensor product $V \otimes W$ equipped
with this topology.\end{dfn}

\begin{rmk} One can prove that a continuous bilinear map $U \times V \to W$
induces a continuous map $U \otimes_{\textup{ic}} V \to W$.
On the downside the ind-tensor product topology has some counterintuitive
properties. For example $V \otimes_{\textup{ic}} \Fq = V^\#$
so $\Fq$ is not a tensor unit for $\otimes_{\textup{ic}}$.
\end{rmk}

%

\begin{lem}\label{indcottodiscrot} If $V, W$ are lth vector spaces then the
natural map $V \otimes_{\textup{ic}} W \to V^\# \otimes_{\textup{c}} W$ is
continuous.\end{lem}

\pf Indeed if $Y \subset W$ is an open subspace then $V \otimes Y$ is open both
in $V \otimes_{\textup{ic}} W$ and in $V^\# \otimes_{\textup{c}} W$. The result
follows since the subspaces $V \otimes Y$ form a fundamental system in $V^\#
\otimes_{\textup{c}} W$.\quod

\begin{lem} If $V, W$ are lth vector spaces then $V \otimes_{\textup{ic}} W$ is
lth.\end{lem}

\pf The space $V^\# \otimes_{\textup{c}} W$ is Hausdorff by Lemma \ref{otlth}.
As the natural bijection $V \otimes_{\textup{ic}} W \to V^\#
\otimes_{\textup{c}} W$ is continuous the point $0 \in V \otimes_{\textup{ic}}
W$ is closed. So $V \otimes_{\textup{ic}} W$ is Hausdorff by Lemma
\ref{lthaus}.\quod

\begin{lem}\label{indcotmapcont} If $f\colon V_1 \to V_2$ and $g\colon W_1 \to
W_2$ are continuous $\Fq$-linear maps of lth vector spaces then
the map $f \otimes g\colon V_1 \otimes_{\textup{ic}} W_1 \to V_2
\otimes_{\textup{ic}} W_2$ is continuous. \quod\end{lem}


\begin{dfn}\label{defindcot}\index{idx}{topological tensor product!ind-complete}%
\index{nidx}{topological tensor products!$\indcot$, Angry Bird}%
Let $V, W$ be lth vector spaces. We define the \emph{ind-complete tensor product}
$V \indcot W$ as the completion of $V \otimes_{\textup{ic}} W$. If $f\colon V_1
\to V_2$ and $g\colon W_1 \to W_2$ are continuous $\Fq$-linear maps of lth
vector spaces then $f \indcot g\colon V_1 \indcot W_1 \to V_2 \indcot W_2$ is
defined as the completion of $f \otimes g$.\end{dfn}

\breakflow
Lemma \ref{indcottodiscrot} equips us with natural maps $V \indcot W \to
V^\#\complot W$ and $V \indcot W \to V \complot W^\#$.

\begin{prp}\label{indcotdescr} If $V, W$ are lth vector spaces then
the natural square
\begin{equation*}
\xymatrix{
V \indcot W \ar[r] \ar[d]&  V \complot W^\#  \ar[d] \\
V^\# \complot W \ar[r] & V \complot W
}
\end{equation*}
is cartesian in the category of topological vector spaces.\end{prp}

\pf The proposition claims that the map
\begin{equation*}
f\colon V \indcot W \to (V^\# \complot W) \times_{V \complot W} (V \complot W^\#)
\end{equation*}
defined by the diagram above is a topological isomorphism. 

Given open subspaces $U \subset V$, $Y \subset W$ let us temporarily
denote
\begin{align*}
[ U, Y ] &= V/U \otimes_{\textup{c}} W/Y, \\
\langle U, Y \rangle &= (V^\# \otimes_{\textup{c}} W/Y) \times_{[U,Y]} (V/U
\otimes_{\textup{c}} W^\#).
\end{align*}
As limits commute with limits the target of the map $f$ is
$\lim_{U,Y} \langle U, Y \rangle$
where $U, Y$ range over all open subspaces. Hence $f$ is defined by the
natural projections
\begin{equation*}
f_{U,Y} \colon V \otimes_{\textup{ic}} W \to \langle U, Y \rangle.
\end{equation*}
A choice of splittings $V \cong U \oplus V/U$, $W \cong Y \oplus W/Y$ induces
isomorphisms
\begin{align*}
\langle U, Y \rangle & \cong
(U^\# \otimes_{\textup{c}} W/Y) \times (V/U \otimes_{\textup{c}} Y^\#) \times [U,Y], \\
V \otimes_{\textup{ic}} W & \cong (U \otimes_{\textup{ic}} Y) \times 
(U \otimes_{\textup{ic}} W/Y) \times
(V/U \otimes_{\textup{ic}} W) \times [U,Y]
\end{align*}
which identify $f_{U,Y}$ with the projection to the last three factors.
Hence $f_{U,Y}$ is onto with the kernel $U \otimes_{\textup{ic}}
Y$. The resulting continuous bijection
\begin{equation*}
(V \otimes_{\textup{ic}} W)/(U \otimes_{\textup{ic}} Y) \to \langle U, Y \rangle
\end{equation*}
is a topological isomorphism since its target and source are both discrete.
Taking the limit over all $U, Y$ we deduce the desired result.\quod

\begin{cor}\label{indcotinj} If $V, W$ are lth vector spaces then
the natural maps $V \indcot W \to V^\# \complot W$ and $V \indcot W \to V
\complot W^\#$ are injective.\end{cor}

\pf According to Proposition \ref{discrcomplot} both natural maps $V^\#\complot
W \to V \complot W$ and $V \complot W^\# \to V \complot W$ are injective. So the
result follows from Proposition \ref{indcotdescr}.\quod

\begin{cor}\label{indcotdiscrcomplot} If $V$ is a compact vector space and $W$ a
complete lth vector space then the natural map $V \indcot W \to V^\# \complot W$ is a
continuous bijection.\end{cor}

\pf Indeed $V \complot W^\# \to V \complot W$ is a bijection by Proposition
\ref{discrcomplot} (2) whence the result follows from Proposition
\ref{indcotdescr}. \quod

\section{Continuous functions}

\begin{dfn}\label{defcontfn}\index{idx}{function spaces and germ spaces!continuous functions}\index{nidx}{function spaces and germ spaces!$c(V,W)$, continuous functions}%
Let $V, W$ be topological vector spaces. We define $c(V,W)$ to be
the space of continuous $\Fq$-linear maps $V \to W$ equipped with the
compact-open topology.\end{dfn}

\begin{lem} Let $V$ be a topological vector space. If $W$ is an lth vector space
then so is $c(V,W)$. \quod\end{lem}

\begin{lem} If an $\Fq$-linear map $V \to W$ is continuous then so are the
induced natural transformations $c(W,-) \to c(V,-)$ and $c(-,V) \to c(-,W)$.
\quod\end{lem}

\begin{rmk}
The natural map $c(U \otimes_{\textup{c}} V, W) \to c(U, c(V,W))$
is not surjective in general and so does not define an adjunction of
$-\otimes_{\textup{c}} V$ and $c(V,-)$.\end{rmk}

\begin{dfn}\label{defevmap} Let $V, W$ be topological vector spaces. We denote
\begin{equation*}
\sigma_{V,W}\colon V^* \otimes W \to c(V,W)
\end{equation*}
the map which sends a tensor $f
\otimes w$ to the function $v \mapsto f(v)w$.\end{dfn}

\begin{prp}\label{contcomplot}For every locally compact vector space $V$ and
a complete lth vector space $W$ there exists a unique topological
isomorphism
\begin{equation*}
V^* \complot W \xrightarrow{\,\isosign\,} c(V,W)
\end{equation*}
extending $\sigma_{V,W}$ on $V^* \otimes W$.\end{prp}

\pf We split the proof in two steps.

\breakflow
\textbf{Step 1.}
\textit{$V$ is compact and $W$ is discrete.}

The space $V^*$ is discrete by
Theorem \ref{pontdual} whence $V^* \complot W = V^* \otimes_{\textup{c}} W$ is discrete.
As a consequence
\begin{equation*}
V^* \complot W = \bigcup_{W' \subset W} V^* \otimes_{\textup{c}} W'
\end{equation*}
where $W' \subset W$ ranges over all finite subspaces.

The space $c(V,W)$ is discrete since $V \subset V$ is compact and $\{0\}
\subset W$ is open. As $V$ is compact and $W$ is discrete the image of every
continuous $\Fq$-linear map $V \to W$ is finite. Therefore
\begin{equation*}
c(V,W) = \bigcup_{W' \subset W} c(V,W')
\end{equation*}
where $W' \subset W$ again ranges over all finite subspaces. Hence it is enough
to consider the case when $W$ is finite. This case instantly reduces to the case
$W \cong \Fq$ which is clear.

\breakflow
\textbf{Step 2.}
\textit{$V$ is locally compact and $W$ is complete lth.}

For an open subspace $Y \subset W$ let $q_Y\colon W \to W/Y$ be the quotient
map. For a compact open subspace $U \subset V$
let $\rho_{U,Y}\colon c(V,W) \to c(U,W/Y)$ be the map given by restriction to
$U$ and composition with $q_Y$.
The subspace $\ker\rho_{U,Y} \subset c(V,W)$ consists of functions which send
the compact subset $U$ to the open subset $Y$. As a consequence it is open.
Lemma \ref{compopen} implies that the collection of all the subspaces
$\ker\rho_{U,Y}$ is a fundamental system of open subspaces in $c(V,W)$.
Every $\rho_{U,Y}$ is surjective by Corollary \ref{openlift}.
Therefore the limit map
\begin{equation*}
\rho\colon c(V,W) \to \lim_{U,Y} c(U,W/Y).
\end{equation*}
defined by the $\rho_{U,Y}$ is the completion map $c(V,W) \to c(V,W)^{\widehat{}}$.
The map $\rho$ is bijective since $V$ is the union of its compact open
subspaces and $W = \lim_Y W/Y$. As $\rho$ is the completion map of $c(V,W)$ it
is in fact a topological isomorphism.

Let $\rho_U\colon V^* \to U^*$ be the restriction map.
Arguing as in the case of $\rho_{U,Y}$ above we conclude that the collection of
subspaces $\ker\rho_U$ is a fundamental system in $V^*$ and that every map
$\rho_U$ is surjective.  As a consequence the limit map
\begin{equation*}
\psi\colon V^* \otimes_{\textup{c}} W \to \lim_{U,Y}(U^*\otimes W/Y).
\end{equation*}
defined by the $\rho_U \otimes q_Y$ is the completion map of $V^*
\otimes_{\textup{c}} W$.
Altogether we obtain a commutative diagram
\begin{equation*}
\xymatrix{
V^* \otimes_{\textup{c}} W \ar[rr]^{\sigma_{V,W}} \ar[d]^\psi && c(V,W) \ar[d]^\rho \\
\lim_{U,Y} (U^* \otimes_{\textup{c}} W/Y) \ar[rr]^{\lim \sigma_{U,W/Y}} && \lim_{U,Y}
c(U,W/Y).
}
\end{equation*}
The bottom arrow is a topological isomorphism as a limit of topological
isomorphisms $\sigma_{U,W/Y}$. Hence $\sigma_{V,W}$ factors through a
topological isomorphism $V^*\complot W \to c(V,W)$. As $V^* \otimes W$ is dense
in $V^* \complot W$ this is the topological isomorphism $\sigma$ we
need to construct.\quod

\section{Bounded functions}

\begin{dfn}\label{defbndfn}\index{idx}{function spaces and germ spaces!bounded functions}\index{nidx}{function spaces and germ spaces!$b(V,W)$, bounded functions}%
Let $V, W$ be locally compact vector spaces. A continuous
$\Fq$-linear map $f\colon V \to W$ is said to be bounded if its image is
contained in a compact subspace. We define $b(V,W) \subset c(V,W)$ to be the
space of all bounded maps. The topology on $b(V,W)$ is given by the fundamental
system of subspaces $c(V,Y)$ where $Y \subset W$ ranges over all compact open
subspaces.\end{dfn}

\begin{lem} The inclusion $b(V,W) \subset c(V,W)$ is continuous. \quod\end{lem}


\begin{lem} If an $\Fq$-linear map $V \to W$ is continuous then so
are the induced natural transformations $b(W,-) \to b(V,-)$ and $b(-,V) \to
b(-,W)$. \quod\end{lem}

\breakflow
Observe that the map $\sigma_{V,W}\colon V^* \otimes W \to c(V,W)$ of Definition
\ref{defevmap} has image in $b(V,W)$.

\begin{prp}\label{bndcomplot} For every pair $V, W$ of locally compact vector spaces
there exists a unique topological isomorphism
\begin{equation*}
(V^*)^\# \complot W \xrightarrow{\,\isosign\,} b(V,W)
\end{equation*}
extending $\sigma_{V,W}$ on $V^* \otimes W$.\end{prp}

\pf The composition
\begin{equation*}
(V^*)^\#\complot W \xrightarrow{\iota} V^* \complot W
\xrightarrow{\sigma} c(V,W)
\end{equation*}
of the natural inclusion $\iota$ and the topological isomorphism $\sigma$ of
Proposition \ref{contcomplot} is continuous and restricts to $\sigma_{V,W}$ on
$V^* \otimes W$. Hence our claim follows if $\sigma\iota$ is a homeomorphism
onto $b(V,W)$.

Let $\cU$ be the family of all compact open subspaces of $W$. $\cU$ is a
fundamental system in $W$ by Proposition \ref{lclth}. It covers $W$ by
Lemma \ref{compopen}. Since $(V^*)^\#$ is discrete it follows that the family
\begin{equation*}
\{ (V^*)^\# \otimes Y \mid Y \in \cU \}
\end{equation*}
is a fundamental system which covers $(V^*)^\#\otimes_{\textup{c}} W$. As a
consequence the family
\begin{equation*}
\{ (V^*)^\# \complot Y \mid Y \in \cU \}
\end{equation*}
is a fundamental system in $(V^*)^\#\complot W$. It covers $(V^*)^\#\complot W$
by Proposition \ref{complcover}.

As $Y \in \cU$ is compact the map $\iota$ identifies $(V^*)^\# \complot Y$ with
$(V^*)\complot Y$ by Proposition \ref{discrcomplot} (2). The map $\sigma$ sends
the latter subspace isomorphically onto $c(V,Y) \subset b(V,W)$. The subspaces
$c(V,Y)$ form a fundamental system in $b(V,W)$. This system covers $b(V,W)$ as a
consequence of Lemma \ref{compopen}. We conclude that $\iota\sigma$ is a
homeomorphism onto $b(V,W)$.
%
\quod

\begin{cor}%
Let $V$ and $W$ be locally compact vector spaces.
If $V$ is compact then the inclusion $b(V,W) \hookrightarrow c(V,W)$
is a topological isomorphism.\end{cor}

\pf Indeed $V^*$ is discrete by Theorem~\ref{pontdual} so the natural map
$(V^*)^\#\complot W \to V^*\complot W$ is a topological isomorphism.\quod

\section{Bounded locally constant functions}

\begin{dfn}\label{defbndlc}\index{idx}{function spaces and germ spaces!bounded locally constant functions}\index{nidx}{function spaces and germ spaces!$a(V,W)$, bounded locally constant functions}%
Let $V, W$ be locally compact vector spaces. A continuous
$\Fq$-linear map $f\colon V \to W$ is called bounded locally constant if it is
bounded and its kernel is open. We define $a(V,W) \subset b(V,W)$ to be the
space of all bounded locally constant maps. The space $a(V,W)$ is equipped with
the minimal topology such that the inclusions $a(V,W) \subset b(V,W)$ and
$a(V,W) \subset c(V,W^\#)$ are continuous.\end{dfn}

\breakflow
One can describe $a(V,W)$ set-theoretically as the intersection
\begin{equation*}
a(V,W) = b(V,W) \cap c(V,W^\#) \subset c(V,W).
\end{equation*}
By construction the topology on $a(V,W)$ is that of the fibre product
\begin{equation*}
b(V,W) \times_{c(V,W)} c(V,W^\#).
\end{equation*}

\begin{lem} If an $\Fq$-linear map $V \to W$ is continuous then so are the
induced natural transformations $a(W,-) \to a(V,-)$ and $a(-,V) \to
a(-,W)$. \quod\end{lem}

\breakflow
Observe that the map $\sigma_{V,W}\colon V^* \otimes W \to c(V,W)$ of Definition
\ref{defevmap} has image in $a(V,W)$.

\begin{prp}\label{bndlcot} For every pair of locally compact vector spaces $V, W$ 
there exists a unique topological isomorphism
\begin{equation*}
V^* \indcot W \xrightarrow{\,\isosign\,} a(V,W)
\end{equation*}
extending $\sigma_{V,W}$ on $V^* \otimes W$.\end{prp}

\pf Consider the commutative diagram
\begin{equation}\label{bndlcotmap}\tag{$\ast$}
\vcenter{\vbox{\xymatrix{
(V^*)^\# \complot W \ar[r] \ar[d] & V^* \complot W \ar[d] & V^* \complot W^\# \ar[l] \ar[d] \\
b(V,W) \ar[r] & c(V,W) & c(V,W^\#) \ar[l]
}}}
\end{equation}
where the horizontal arrows are the natural maps, the left vertical arrow is the
topological isomorphism of Proposition \ref{bndcomplot} and the other two
vertical arrows are the topological isomorphisms provided by Propositon
\ref{contcomplot}. According to Proposition \ref{indcotdescr}
\begin{equation*}
V^* \indcot W = ((V^*)^\# \complot W) \times_{V^* \complot W} (V^* \complot W^\#).
\end{equation*}
At the same time
\begin{equation*}
a(V,W) = b(V,W) \times_{c(V,W)} c(V,W^\#).
\end{equation*}
Thus \eqref{bndlcotmap} defines a topological isomorphism $V^*
\indcot W \to a(V,W)$. It extends $\sigma_{V,W}$ since the vertical maps in
\eqref{bndlcotmap} do so.\quod

\begin{cor}%
Let $V$ and $W$ be locally compact vector spaces.
If $W$ is discrete then the inclusion $a(V,W) \hookrightarrow b(V,W)$
is a topological isomorphism.\end{cor}

\pf The map $V^*\complot W^\# \to V^*\complot W$ is a topological isomorphism.
Proposition~\ref{indcotdescr} implies that the map $V^*\indcot W \to (V^*)^\#\complot W$
is a topological isomorphism. Whence the result.\quod

\section{Germ spaces}
\label{sec:germ}

\begin{dfn}\label{defgerm}\index{idx}{function spaces and germ spaces!the space of germs}\index{nidx}{function spaces and germ spaces!$g(V,W)$, germs}%
Let $V, W$ be lth vector spaces. The $\Fq$-vector
\emph{space of germs} $g(V,W)$
is
\begin{equation*}
g(V,W) = \colim_U c(U,W)
\end{equation*}
where $U \subset V$ runs over all open subspaces and the transition maps are
restrictions. We do not equip $g(V,W)$ with a topology. The image of $f\in
c(U,W)$ under the natural map $c(U,W) \to g(V,W)$ is called the germ of
$f$.\end{dfn}

\breakflow
An element of $g(V,W)$ can be represented by a pair $(U,f)$ where $U
\subset V$ is an open subspace and $f\colon U \to W$ a continuous $\Fq$-linear
map. Two such pairs $(U_1,f_1)$, $(U_2, f_2)$ represent the same element of
$g(V,W)$ if there exists an open subspace $U \subset U_1 \cap U_2$ such that
$f_1|_U = f_2|_U$.

%
\begin{prp}\label{germses}%
The natural sequence
\begin{equation*}
0 \to a(V,W) \to b(V,W) \to g(V,W) \to 0
\end{equation*}
is exact
for all locally compact vector spaces $V$, $W$.\end{prp}

\pf The sequence is clearly left exact. We need to prove that $b(V,W) \to
g(V,W)$ is surjective.  Let $U \subset V$ be an open subspace. As $V$ is locally
compact there exists a  compact open subspace $U'\subset U$. According to
Corollary \ref{openlift} the restriction map $c(U,W) \to c(U',W)$ is onto.
Furthermore $c(U',W) = b(U',W)$ since $U'$ is compact. The map $b(V,W) \to
b(U',W)$ is surjective by Corollary \ref{openlift} again. Hence the map $b(V,W) \to
g(V,W)$ is surjective.\quod

\begin{dfn}\label{deflociso}%
Let $V, W$ be lth vector spaces. A continuous $\Fq$-linear map
$f\colon V \to W$ is called a \emph{local isomorphism} if there exists an open
subspace $U \subset V$ such that $f(U) \subset W$ is open and the restriction
$f|_U\colon U \to f(U)$ is a topological isomorphism.\end{dfn}

\begin{prp}\label{germiso} If $f\colon V \to W$ is a local isomorphism of lth
vector spaces
then the induced natural transformations $g(W,-) \to g(V,-)$ and $g(-,V) \to
g(-,W)$ are isomorphisms. \quod\end{prp}



\begin{lem}\label{germquotses} For every pair of locally compact vector spaces
$V, W$ the natural map $V\indcot W \to V^\#\complot W$ extends to a
natural short exact sequence
\begin{equation*}
0 \to V \indcot W \to V^\#\complot W \to g(V^*,W) \to 0.
\end{equation*}\end{lem}

\pf Consider the short exact sequence
\begin{equation*}
0 \to a(V^*,W) \to b(V^*,W) \to g(V^*,W) \to 0
\end{equation*}
of Proposition \ref{germses}. The isomorphisms
\begin{align*}
b(V^*,W) &\cong (V^{**})^\#\complot W, \\
a(V^*,W) &\cong V^{**} \indcot W
\end{align*}
of Propositions \ref{bndcomplot}, \ref{bndlcot}
and Pontrjagin duality of Theorem \ref{pontdual} transform it to
\begin{equation*}
0 \to V \indcot W \to V^\#\complot W \to g(V^*,W) \to 0
\end{equation*}
and the result follows. \quod

\begin{prp}\label{germquotiso} Let $f_V\colon V_1 \to V_2$ and $f_W\colon W_1
\to W_2$ be continuous $\Fq$-linear maps of locally compact vector spaces.
If $f_V^*$ and $f_W$ are local isomorphisms then the induced map
\begin{equation*}
\frac{V_1^\# \complot W_1}{V_1 \indcot W_1} \to \frac{V_2^\#\complot W_2}{V_2\indcot
W_2}
\end{equation*}
is an isomorphism.\end{prp}

\pf The induced map $g(V_1^*,W_1) \to g(V_2^*,W_2)$ is an isomorphism by
Proposition \ref{germiso}. Hence the result follows from Lemma
\ref{germquotses}. \quod

%

\chapter{Topological rings and modules}
\label{chapter:trm}
\label{ch:algcompl}
\label{ch:bounded}
\label{ch:compideal}
\label{ch:ringot}
\label{ch:otalgprops}
\label{ch:otalgexamples}
\label{ch:ottaustruct}
\label{ch:funcmod}
\label{ch:resdual}
\label{ch:chdom}
\label{ch:topmodcompl}

In this chapter we use the tensor product and function space constructions of
Chapter \ref{chapter:tvs} to produce and study rings and modules over them.
%
%

We keep the conventions of
Chapter \ref{chapter:tvs}.
In particular we continue using the acronym ``lth'' and assume all locally
compact vector spaces to be Hausdorff.  We work with topological algebras over the
fixed field $\Fq$ and with modules over such algebras.
In this chapter an algebra (without further qualifications) means an $\Fq$-algebra.

A topological $\Fq$-algebra is a topological $\Fq$-vector space $\nC$ equipped with a
continuous multiplication map $\nC \times \nC \to \nC$ which makes $\nC$ into a
commutative associative unital $\Fq$-algebra. A topological module $\nmC$ over a
topological algebra $\nC$ is a topological vector space $\nmC$ equipped with a
continuous $\nC$-action map $\nC \times \nmC \to \nmC$ which makes $\nmC$ into an $\nC$-module.



\section{Overview}

Let $\nC, \nCalt$ be locally compact $\Fq$-algebras. Typical examples of such
algebras relevant to our applications are
\begin{itemize}
\item the discrete algebra $\Fq[t]$,

\item the locally compact algebra $\Fq(\!(t^{-1})\!)$,

\item the compact algebra
$\Fq[[t^{-1}]]$.
\end{itemize}
The first goal of this chapter is to
equip the tensor products $\nC\indcot \nCalt$ and $\nC \complot \nCalt$ with topological
$\Fq$-algebra structures compatible with the dense subalgebra $\nC \otimes \nCalt$. In
the case of $\nC \complot \nCalt$ it can be done only under certain assumptions on $\nC$,
$\nCalt$ (cf. Example~\ref{algcomplotfail}). To handle this difficulty we work
out some preliminaries 
in Section~\ref{sec:otprelim}.
The rings $\nC \indcot \nCalt$ and $\nC \complot \nCalt$ play a prominent role
in this work. Some degree of familiarity
with them will be assumed in the subsequent chapters.
We discuss examples of such rings in Section \ref{sec:otalgexamples}.

Let $\nmC$ be a locally compact $\nC$-module and $\nmCalt$ a locally compact $\nCalt$-module.
Another important goal of this chapter is to equip the function spaces $a(\nmC,\nmCalt)$,
$b(\nmC,\nmCalt)$, $c(\nmC,\nmCalt)$ and the germ space $g(\nmC,\nmCalt)$ with natural actions of tensor
product rings:
\begin{itemize}
\item an $\nC\indcot \nCalt$-module structure on $a(\nmC,\nmCalt)$ and $g(\nmC,\nmCalt)$,

\item an $\nC^\#\complot \nCalt$-module structure on $b(\nmC,\nmCalt)$,

\item an $\nC\complot \nCalt$-module structure on $c(\nmC,\nmCalt)$, under certain assumptions.
\end{itemize}

In Section \ref{sec:ottaustruct} 
we fix $\tau$-ring structures on $\nC \indcot \nCalt$ and $\nC \complot \nCalt$ 
to facilitate applications in the context of shtukas.
We also fix the structures
of left modules over suitable $\tau$-polynomial rings
on the function spaces and the germ spaces.
As a result one can use them
as arguments for the Hom shtuka construction of Section~\lref{genhomsht}{sec:genhomsht}.
A Hom shtuka with a function space or a germ space argument is one of the main
constructions of this text.

In Section \ref{sec:resdual} we study $a(\nmC,\nmCalt)$, $b(\nmC,\nmCalt)$ and $c(\nmC,\nmCalt)$ as modules in one
particular case which is central to our applications. Let $C$ be a smooth
projective connected curve over $\Fq$. Fix a closed point $\infty \in C$ and set $A =
\Gamma(C - \{\infty\},\cO_C)$. The local field $F$ of $C$ at
$\infty$ is a locally compact $\Fq$-algebra. Its ring of integers
$\cO_F \subset F$ is a compact open subalgebra while $A \subset F$ is a discrete
cocompact subalgebra. Let $\omega_A$ be the module of K\"ahler differentials
of $A$ over $\Fq$ and let $\omega_F = \omega_A \otimes_A F$.
Serre duality for $C$ implies important results for module
structures on $a(\nmC,\nmCalt)$, $b(\nmC,\nmCalt), c(\nmC,\nmCalt)$ where $\nmC$ is one of
the spaces
\begin{equation*}
F, \quad F/A, \quad F/\cO_F.
\end{equation*}
For example $a(F, \nCalt)$ is topologically isomorphic to
$\omega_F \indcot \nCalt$, a free $F \indcot \nCalt$-module of rank $1$,
while $b(F/A,\nCalt) \cong \omega_A \complot \nCalt$ is a locally free $A \complot \nCalt$-module of rank~$1$.

The content of Section \ref{sec:otprelim} 
is well-known. The same applies to the rest of the chapter modulo the reservations
we made
on the tensor product $\indcot$, the function spaces $a(\nmC,\nmCalt)$, $b(\nmC,\nmCalt)$ and the
germ space $g(\nmC,\nmCalt)$ in Chapter \ref{chapter:tvs}.

As is the case for Chapter \ref{chapter:tvs}, this chapter was inspired by and
owes much to Anderson's work \cite{and}, especially to \S 2 of that text where
Anderson uses function spaces to compute certain $\textup{Ext}$'s for modules over
$\tau$-polynomial rings.

\section{Preliminaries}\index{nidx}{topological vector spaces and modules!$\widehat{V}$, completion}
\label{sec:otprelim}
\label{sec:algcompl}
\label{sec:boundedness}
\label{sec:compideal}

\begin{lem}\label{topmod}Let $U, V$ and $W$ be topological vector spaces.
A bilinear map $\mu\colon U \times V \to W$ is continuous if and only if
\begin{enumerate}
\item For every $u \in U$ the map $\mu(u,-)\colon V \to W$ is continuous.

\item For every $v \in V$ the map $\mu(-,v)\colon U \to W$ is continuous.

\item The map $\mu$ is continuous at $(0,0)$.
\end{enumerate}\end{lem}

\pf Assume (1), (2), (3). Let $(u,v) \in U \times V$ and let $u' \in U$,
$v' \in V$. By bilinearity of $\mu$ we have
\begin{equation*}
\mu(u + u', v + v') = \mu(u,v) + \mu(u,v') + \mu(u',v) + \mu(u',v').
\end{equation*}
From (1), (2), (3) it follows that the map $(u',v') \mapsto \mu(u + u', v + v')$
is continuous at $(0,0)$. Hence $\mu$ is continuous. \quod
 
\begin{lem} Let $U, V$ and $W$ be lth vector spaces, $\mu\colon U \times V \to W$
a bilinear map. If $\mu$ is continuous then there exists a unique bilinear continuous
map $\widehat{\mu}\colon \widehat{U} \times \widehat{V} \to \widehat{W}$
extending $\mu$.\end{lem}

\pf \cite[III, p.50, Theorem 1]{btg}.\quod

\begin{lem}\label{complalg} Let $\nC$ be an lth vector space equipped with an
$\Fq$-algebra structure.
\begin{enumerate}
\item $\widehat{\nC}$ admits at most one structure of a topological $\Fq$-algebra such that the
natural map $\nC \to \widehat{\nC}$ is a homomorphism.

\item Such a structure exists if and only if the multiplication map
$\mu\colon \nC \times \nC \to \nC$ is continuous.
\end{enumerate}\end{lem}

\pf (1) follows from the fact that $\nC \times \nC$ is dense in $\widehat{\nC} \times
\widehat{\nC}$.
(2) If $\mu$ is continuous then we get a continuous multiplication
$\widehat{\mu}\colon \widehat{\nC} \times \widehat{\nC} \to \widehat{\nC}$
by completion. Conversely if there exists a $\widehat{\mu}$ such that the
natural map $\nC \to \widehat{\nC}$ is a homomorphism then $\mu$ is
continuous since every open subset of $\nC$ is a preimage of an open subset
in $\widehat{\nC}$.\quod

\begin{lem}\label{compllattice}%
Let $\nmC$ be an lth module of finite type over an
lth algebra $\nC$. If $\nmC$ is a topological direct summand of $\nC^{\oplus n}$ for some $n$
then the natural map $\nmC \otimes_{\nC} \widehat{\nC} \to \widehat{\nmC}$ is an
$\widehat{\nC}$-module isomorphism.\end{lem}

\pf The claim reduces to the case when $\nmC = \nC^{\oplus n}$ and this case is
clear.\quod

%


\breakflow
We will need the notion of boundedness from topological ring theory.

\begin{dfn} Let $\nC$ be a topological algebra and $\nmC$ a topological $\nC$-module.
\begin{enumerate}
\item If $V \subset \nC$,
$K \subset \nmC$ are subsets then $V \cdot K$ denotes the set of products
\begin{equation*}
\{ \neC \cdot \nmeC \mid \neC \in V, \nmeC \in K \} \subset \nmC.
\end{equation*}
\item A subset $K \subset \nmC$ is called \emph{bounded} if for every open subspace $U
\subset \nmC$ there exists an open subspace $V \subset \nC$ such that $V \cdot K
\subset U$.
%
\end{enumerate}\end{dfn}


\begin{lem}\label{compbounded} Let $\nmC$ be an lth module over an lth algebra $\nC$.
Every compact subset of $\nmC$ is bounded.
\end{lem}

\pf Let $K$ be a compact subset and $\nmeC \in K$ a point. The preimage of $U$ in $\nC
\times \nmC$ under the multiplication map contains the point $(0,\nmeC)$. As this
preimage is open there exist open subspaces $V_\nmeC \subset \nC$, $Y_\nmeC \subset \nmC$ such
that $V_\nmeC \cdot (\nmeC + Y_\nmeC) \subset U$.  The translates $\nmeC + Y_\nmeC$ cover $K$. As
$K$ is compact we can choose finitely many such translates. Let $V$ be the
intersection of the corresponding $V_\nmeC$. By construction $V \cdot K \subset
U$.\quod


\begin{lem}\label{compideal} Let $\nC$ be a compact algebra and $\nmC$ an lth 
$\nC$-module. Every compact open subspace $U \subset \nmC$ contains an open
submodule.\end{lem}

\pf Given $\neC \in \nC$ let $\mu_\neC\colon \nmC \to \nmC$ denote the multiplication by $\neC$
map. As $U$ is compact and open Lemma \ref{compbounded} provides us with an open
subspace $V \subset \nC$ such that $V \cdot U \subset U$. Consider the
intersection
\begin{equation*}
U' = U \cap \bigcap_\neC \mu_\neC^{-1} U
\end{equation*}
where $\neC$ runs over a set of representatives of the classes in $\nC/V$. The
subspace $U'$ is open since $\nC/V$ is finite. By construction $\nC \cdot U'
\subset U$ hence the submodule generated by $U'$ is contained in
$U$.\quod

\section{Tensor products}

Let $\nC, \nCalt$ be locally compact algebras. Even though $\nC \otimes_{\textup{c}} \nCalt$
is both an $\Fq$-algebra and an lth vector space it is not necessarily an lth
algebra since its multiplication need not be continuous in the tensor product
topology.

%
%

\begin{example}\label{algcomplotfail} Let us demonstrate that the multiplication
on $\Fq(\!(z)\!) \otimes_{\textup{c}} \Fq(\!(z)\!)$ is not continuous.
Denote $F = \Fq(\!(z)\!)$ temporarily. We have
\begin{equation*}
\lim_{n \to \infty} z^{-n} \otimes z^n = 0 = \lim_{n \to \infty} z^n \otimes z^{-n}
\end{equation*}
in $F \otimes_{\textup{c}} F$. As a consequence
\begin{equation*}
\lim_{n \to \infty} (z^{-n} \otimes z^n, z^n \otimes z^{-n}) = 0
\end{equation*}
in $(F \otimes_{\textup{c}} F) \times (F \otimes_{\textup{c}} F)$.
The multiplication maps
this sequence to the constant sequence $1 \otimes 1$. Hence it is not
continuous.
\end{example}

\breakflow
We would like to give a sufficient condition for $\nC\complot \nCalt$ to carry a natural
topological algebra structure. To do it we need some preparation.


%

%
%

\begin{lem}\label{ltbmodot}%
Let $\nmC$ be an lth module over an lth algebra $\nC$.
If $\nmC$ is bounded and admits a fundamental system of open submodules then the
map $\mu\colon \nC \otimes_{\textup{c}} \nmC \to \nmC$ induced by the
$\nC$-multiplication on $\nmC$ is continuous.\end{lem}

\pf Let $U \subset \nmC$ be an open submodule. As $\nmC$ is bounded there exists an
open subspace $V \subset \nC$ such that $V \cdot \nmC \subset U$.  Therefore $\mu(V
\otimes \nmC + \nC \otimes U) \subset U$. \qed

\begin{lem}\label{modcontot}%
Let $\nC, \nCalt$ be lth algebras, $\nmC$ an lth $\nC$-module
and $\nmCalt$ an lth $\nCalt$-module. If the map $\nC \otimes_{\textup{c}} \nmC \to \nmC$ induced
by the $\nC$-multiplication on $\nmC$ is continuous then the $\nC \otimes_{\textup{c}}
\nCalt$-multiplication on $\nmC \otimes_{\textup{c}} \nmCalt$ is continuous.\end{lem}

\pf The $\nC \otimes_{\textup{c}} \nCalt$-multiplication on $\nmC \otimes_{\textup{c}} \nmCalt$
is a bilinear map which satisfies the conditions (1) and (2) of Lemma~\ref{topmod}.
It remains to show that the condition (3) holds.
Let $U_\nmC \subset \nmC$, $U_\nmCalt \subset \nmCalt$ be open subspaces. 
As the map $\nC \otimes_{\textup{c}} \nmC \to \nmC$ is continuous there are open subspaces
$U_\nC \subset \nC$ and $V_\nmC \subset \nmC$ such that $U_\nC \cdot \nmC \subset U_\nmC$ and
$\nC \cdot V_\nmC \subset U_\nmC$. By continuity of $\nCalt$-multiplication on $\nmCalt$
there exist open subspaces $U_\nCalt \subset \nCalt$ and $V_\nmCalt \subset \nmCalt$ such that
$U_\nCalt \cdot V_\nmCalt \subset U_\nmCalt$. We then
have
\begin{gather*}
(U_\nC \otimes \nCalt + \nC \otimes U_\nCalt) \cdot (V_\nmC \otimes \nmCalt + \nmC \otimes V_\nmCalt) 
\subset \\
(U_\nC \cdot V_\nmC) \otimes \nmCalt + (U_\nC \cdot \nmC) \otimes \nmCalt +
(\nC \cdot V_\nmC) \otimes \nmCalt + \nmC \otimes (U_\nCalt \cdot V_\nmCalt)
\subset \\
U_\nmC \otimes \nmCalt + \nmC \otimes U_\nmCalt
\end{gather*}
and the result follows.\quod

\begin{prp}\label{algcomplot}\index{idx}{topological tensor product!completed}\index{nidx}{topological tensor products!$\complot$, completed}%
Let $\nC, \nCalt$ be lth algebras. %
If $\nC$ is compact or discrete %
then $\nC \complot \nCalt$ admits a %
unique structure of a topological algebra such that
the natural map $\nC \otimes \nCalt \to \nC \complot \nCalt$ is a homomorphism.\end{prp}

\pf We first prove that the multiplication map $\nC \otimes_{\textup{c}} \nC \to \nC$
is continuous. This is clear if $\nC$ is discrete. If $\nC$ is compact then it is
bounded by Lemma \ref{compbounded} and admits a fundamental system of open
$\nC$-submodules by Lemma \ref{compideal}. Hence the map $\nC
\otimes_{\textup{c}} \nC \to \nC$ is continuous by Lemma \ref{ltbmodot}.
Now Lemma \ref{modcontot} implies that 
the multiplication map $(\nC \otimes_{\textup{c}} \nCalt) \times (\nC
\otimes_{\textup{c}} \nCalt) \to (\nC \otimes_{\textup{c}} \nCalt)$ is continuous.
Whence the result follows from Lemma \ref{complalg}. \quod

\afterall
Let $\nC, \nCalt$ be lth algebras, $\nmC$ an lth $\nC$-module and $\nmCalt$ an lth $\nCalt$-module.
The topological vector space $\nmC \complot \nmCalt$ comes equipped with a natural action of $\nC \otimes \nCalt$
by functoriality of $\complot$. Recall that according to our convention 
the ring $\nC \otimes \nCalt$ carries no topology so its action is not supposed to be continuous.
Nonetheless we can prove the following result.

\begin{prp}\label{modcomplot}%
Let $\nC, \nCalt$ be lth algebras, $\nmC$ an lth $\nC$-module
and $\nmCalt$ an lth $\nCalt$-module.
Assume that one of the following conditions hold:
\begin{enumerate}
\item $\nC$ and $\nmC$ are discrete,

\item $\nC$ and $\nmC$ are compact.
\end{enumerate}
Then the $\nC \otimes \nCalt$-module structure on $\nmC \complot \nmCalt$ extends
to a unique topological $\nC\complot \nCalt$-module structure.
\end{prp}

\pf The unicity is clear. If either (1) or (2) holds then the same argument
as in the proof of Proposition \ref{algcomplot} shows that the multiplication
map $(\nC \otimes_{\textup{c}} \nCalt) \times (\nmC \otimes_{\textup{c}} \nmCalt) \to
(\nmC \otimes_{\textup{c}} \nmCalt)$ is continuous. Taking its completion
we get a topological $\nC\complot \nCalt$-module structure on $\nmC \complot \nmCalt$.
It is compatible with the natural $\nC \otimes \nCalt$-module structure on the
dense subspace $\nmC \otimes \nmCalt$. Whence the result. \quod

\breakflow
As we observed above the completed tensor product $\nC \complot \nCalt$ of locally
compact algebras $\nC$ and $\nCalt$ need not carry a natural lth algebra
structure. The ind-complete tensor product $\nC \indcot \nCalt$ does not suffer from
such a problem.

\begin{lem}\label{indcotmodcont}If $\nC, \nCalt$ are lth algebras, $\nmC$ an lth
$\nC$-module and $\nmCalt$ an lth $\nCalt$-module then the multiplication map
$(\nC \otimes_{\textup{ic}} \nCalt) \times (\nmC \otimes_{\textup{ic}} \nmCalt) \to \nmC
\otimes_{\textup{ic}} \nmCalt$ is continuous.\end{lem}

\pf The conditions (1) and (2) of Lemma~\ref{topmod} are clear and the
condition (3) follows instantly from the definition of the ind-tensor product topology
(Definition \lref{indcot}{indcotdef}).\qed

\begin{prp}\index{idx}{topological tensor product!ind-complete}\index{nidx}{topological tensor products!$\indcot$, Angry Bird}%
If $\nC, \nCalt$ are locally compact algebras then there exists a unique
lth algebra structure on $\nC \indcot \nCalt$ such that the natural map $\nC \otimes \nCalt
\to \nC \indcot \nCalt$ is a homomorphism. \quod\end{prp}

\begin{prp}\label{indcotmod}If $\nC, \nCalt$ are locally compact algebras, $\nmC$ a locally compact
$\nC$-module and $\nmCalt$ a locally compact $\nCalt$-module then the natural $\nC \otimes \nCalt$-module
structure on $\nmC \indcot \nmCalt$ extends to a unique topological
$\nC \indcot \nCalt$-module structure.\quod\end{prp}

\breakflow
Next we study some natural maps of tensor product algebras.

\begin{prp}\label{ringdiscrcomplot} Let $\nC, \nCalt$ be lth algebras.
Assume that either $\nC$ or $\nCalt$ is compact or discrete. 
\begin{enumerate}
\item  The natural map $\iota\colon \nC \complot \nCalt^\# \to \nC \complot \nCalt$ is an injective
$\Fq$-algebra homomorphism.

\item If $\nC$ is compact then $\iota$ is an $\Fq$-algebra
isomorphism.
\end{enumerate}\end{prp}

\pf The map $\iota$ is injective by Proposition \lref{complot}{discrcomplot} (1).
It is a homomorphism of $\Fq$-algebras since it is the completion
of the continuous homomorphism 
$\nC\otimes_{\textup{c}} \nCalt^\# \to \nC \otimes_{\textup{c}} \nCalt$.
If $\nC$ is compact then $\iota$ is bijective by Proposition
\lref{complot}{discrcomplot} (2).\quod

\begin{prp}\label{ringindcotdiscrcomplot} Let $\nC, \nCalt$ be locally compact
algebras.
\begin{enumerate}
\item The natural map $\iota\colon \nC \indcot \nCalt \to \nC^\#\complot \nCalt$ is an injective
$\Fq$-algebra homomorphism.

\item If $\nC$ is compact then $\iota$ is an isomorphism.
\end{enumerate}\end{prp}

\pf Follows from Corollaries \lref{indcot}{indcotinj}
and \lref{indcot}{indcotdiscrcomplot}.\quod


\begin{prp}\label{indcotcomplot} If $\nC, \nCalt$ are compact algebras then the natural
map $\nC \indcot \nCalt \to \nC \complot \nCalt$ is a $\Fq$-algebra isomorphism.\end{prp}

\pf{}The map in question is the composition of continuous homomorphisms $\nC \indcot \nCalt
\to \nC \complot \nCalt^\# \to \nC \complot \nCalt$. The first one is an isomorphism by
Proposition \ref{ringindcotdiscrcomplot} while the second one is an isomorphism
by Proposition \ref{ringdiscrcomplot} (2). \quod


\begin{rmk} The natural maps of Propositions \ref{ringdiscrcomplot} and
\ref{ringindcotdiscrcomplot} are always continuous. However their inverses,
if they exist, are not continuous in general. Similarly the natural map of
Proposition \ref{indcotcomplot} is a continuous bijection whose inverse is not
continuous in general.\end{rmk}

\section{Examples of tensor product algebras}
\label{sec:otalgexamples}
We are now in position to discuss some examples of tensor product algebras.
Consider the locally compact algebras
\begin{equation*}
\begin{array}{r@{}lr@{}l}
F &{}= \Fq(\!(\nuC)\!), &\cO_{F} &{}=\Fq[[\nuC]] \\
K &{}= \Fq(\!(\nuCalt)\!), &\cO_{K} &{}=\Fq[[\nuCalt]].
\end{array}
\end{equation*}
We have
\begin{align*}
\cO_{F}\indcot\cO_{K} &= \Fq[[\nuC,\nuCalt]], \\
\cO_{F}\complot\cO_{K} &= \Fq[[\nuC,\nuCalt]]
\end{align*}
as abstract $\Fq$-algebras.
However the topologies on $\cO_{F}\indcot\cO_{K}$ and
$\cO_{F}\complot\cO_{K}$ are different. The topology on
$\cO_{F}\indcot\cO_{K}$ is given by the powers of the ideal $(\nuC\nuCalt)$ while the
topology on $\cO_{F}\complot\cO_{K}$ is given by the powers of the ideal
$(\nuC,\nuCalt)$. In particular $\cO_{F}\complot\cO_{K}$ is compact while
$\cO_{F}\indcot\cO_{K}$ is not even locally compact.

The completed tensor products with a discrete factor look as follows:
\begin{align*}
\cO_{F}^\#\complot\cO_{K} &= \Fq[[\nuC,\nuCalt]], \\
\cO_{F}^\#\complot K &= \Fq[[\nuC]](\!(\nuCalt)\!), \\
F^\#\complot\cO_{K} &= \Fq(\!(\nuC)\!)[[\nuCalt]], \\
F^\#\complot K &= \Fq(\!(\nuC)\!)(\!(\nuCalt)\!).
\end{align*}
The topologies on $\cO_{F}^\#\complot\cO_{K}$ and $F^\#\complot\cO_{K}$ are given by
powers of the ideals $(\nuCalt)$. The topologies on $\cO_{F}^\#\complot K$ and
$F^\#\complot K$ are determined by open subalgebras $\cO_{F}^\#\complot\cO_{K}$ and
$F^\#\complot\cO_{K}$ respectively.

The ind-complete tensor products with a compact factor has the following form:
\begin{equation*}
F \indcot \cO_{K} = \Fq[[\nuCalt]](\!(\nuC)\!).
\end{equation*}
Its topology is defined by the open subalgebra $\cO_{F} \indcot\cO_{K}$.
The ind-complete tensor product of $F$ and $K$ is
\begin{equation*}
F\indcot K =
\Fq[[\nuC,\nuCalt]] [(\nuC\nuCalt)^{-1}]
\end{equation*}
with the topology given by the open subalgebra $\cO_{F}\indcot\cO_{K}$. As we
demonstrated in Example \ref{algcomplotfail} the tensor product
$F \complot K$ makes no sense as an lth algebra.

Another important example is $F \complot \bF_q[x]$.
It is topologically isomorphic to the algebra of Tate series
\begin{equation*}
F \langle x \rangle =
\Big\{ \sum_{n\geqslant 0} \alpha_n x^n \in F[[x]] \,\,\Big| \lim_{n \to \infty}
\alpha_n = 0 \Big\}.
\end{equation*}
For the sake of completeness let us describe the algebra $\cO_{F}\complot K$.
It can be identified with the algebra of power series
\begin{equation*}
\sum_{n \in \bZ} \alpha_n \nuCalt^n, \,\, \alpha_n \in \cO_{F}, \,\, \lim_{n \to -\infty}
\alpha_n = 0.
\end{equation*}
For every nonzero ideal $I \subset \cO_{F}$ and every integer $m \in \bZ$
the subspace
\begin{equation*}
\Big\{ \sum_{n \in \bZ} \alpha_n \nuCalt^n \mid 
\alpha_n \in I \textup{ for all }n \leqslant m \Big\} \subset \cO_{F}\complot K
\end{equation*}
is open. Such subspaces form a fundamental system.

%
%

\section{Algebraic properties}
\label{sec:otalgprops}

In this section we study the properties of $\nC \indcot \nCalt$, $\nC \complot \nCalt$ as
commutative rings without topology. We are primarily interested in the case when
$\nC$ and $\nCalt$ are finite products of local fields or the rings of integers in such
finite products.
We begin with some localization properties.

\begin{prp}\label{discrcomplotloc}If $\nC$ is an $\Fq$-algebra, $K$ a finite
product of local fields and $\nuCalt \in \cO_{K}$ a uniformizer then
$(\nC^\#\complot\cO_{K})[\nuCalt^{-1}] = \nC^\#\complot K$.\end{prp}

\pf The family of open subspaces $\{\nC^\# \otimes_{\textup{c}} \nuCalt^{-n} \cO_{K}\}_{n
\geqslant 0}$ covers $\nC^\# \otimes_{\textup{c}} K$. Proposition
\lref{veccompl}{complcover} implies that the family $\{\nC^\#\complot \nuCalt^{-n} \cO_{K}\}_{n
\geqslant 0}$ covers $\nC^\#\complot K$. Multiplication by $\nuCalt^n$ maps
$\nC^\#\complot \nuCalt^{-n} \cO_{K}$ bijectively onto $\nC^\#
\complot \cO_{K}$ since the same is true with $\otimes_{\textup{c}}$ in place of
$\complot$. The claim is now clear. \quod

\begin{prp}\label{compindcotloc}If $\nC$ is a compact $\Fq$-algebra, $K$ a
finite product of local fields and $\nuCalt \in K$ a uniformizer then $(\nC \indcot
\cO_{K})[\nuCalt^{-1}] = \nC \indcot K$.\end{prp}

\pf As $\nC$ is compact the natural maps $\nC \indcot \cO_{K} \to \nC^\#\complot\cO_{K}$
and $\nC \indcot K \to \nC^\# \complot K$ are isomorphisms by Proposition
\ref{ringindcotdiscrcomplot}. Hence the claim follows from Proposition
\ref{discrcomplotloc}. \quod 

\begin{prp}\label{indcotloc}If $F$, $K$ are finite products of local fields with
uniformizers $\nuC \in \cO_{F}$, $\nuCalt \in \cO_{K}$ then
$(\cO_{F} \indcot \cO_{K})[(\nuC\nuCalt)^{-1}] = F \indcot K$.\end{prp}

\pf The family of open subspaces $\{\nuC^{-n}\cO_{F} \otimes_{\textup{ic}}
\nuCalt^{-n}\cO_{K}\}_{n \geqslant 0}$ covers $F \otimes_{\textup{ic}} K$.
Its completion $\{\nuC^{-n}\cO_{F} \indcot \nuCalt^{-n}\cO_{K}\}_{n \geqslant 1}$ covers
$F \indcot K$ by Proposition \lref{veccompl}{complcover}. Multiplication by
$(\nuC\nuCalt)^n$ maps $\nuC^{-n} \cO_{F} \indcot \nuCalt^{-n} \cO_{K}$ bijectively onto $\cO_{F}
\indcot \cO_{K}$ since the same is true with $\otimes_{\textup{ic}}$ in place of
$\indcot$.  The claim is now clear.
\quod

\breakflow
Next we study quotients of tensor product algebras.
Observe that an ideal $I \subset \cO_{K}$ is open if and only if it projects to
nonzero ideals in all factors of $\cO_{K} = \prod_{i=1}^n \cO_{K_i}$. If $\nC$ is an
$\Fq$-algebra and $I \subset \cO_{K}$ an open ideal then we have a natural map
$\nC^\#\complot\cO_{K} \to \nC \otimes \cO_{K}/I$.

\begin{prp}\label{discrcomplotquot}%
Let $\nC$ be an $\Fq$-algebra and $K$ a finite product of local fields. If $I
\subset \cO_{K}$ is an open ideal then
%
%
the following holds:
\begin{enumerate}
\item The sequence
$0 \to \nC^\#\complot I \to \nC^\#\complot \cO_{K} \to \nC \otimes \cO_{K}/I \to 0$
is exact.

\item The natural map $(\nC^\#\complot\cO_{K})\otimes_{\cO_{K}} I \to
\nC^\#\complot\cO_{K}$ is injective with image $\nC^\#\complot I$.
\end{enumerate}\end{prp}

\pf (1) Indeed the sequence $0 \to \nC^\#\otimes_{\textup{c}} I \to
\nC^\#\otimes_{\textup{c}} \cO_{K} \to \nC^\#\otimes_{\textup{c}} \cO_{K}/I \to 0$ is
clearly exact and the first map in it is an open embedding. Hence the result
follows from Lemma \lref{veccompl}{complopen}.
(2) Observe that $I$ is a free $\cO_{K}$-module of rank $1$. Let $x \in I$ be a
generator.
Multiplication by $x$ identifies $\nC^\#\otimes_{\textup{c}} \cO_{K}$
with $\nC^\#\otimes_{\textup{c}} I$. Taking completion we get the result.
\quod


\breakflow
If $I \subset \cO_{K}$ is an open ideal then the quotient $\cO_{K}/I$ is discrete.
As a consequence $\nC \otimes_{\textup{ic}} \cO_{K}/I$ is discrete. Taking the
completion of $\nC \otimes_{\textup{ic}} \cO_{K} \to \nC \otimes_{\textup{ic}}
\cO_{K}/I$ we get a natural map $\nC \indcot \cO_{K} \to \nC \otimes \cO_{K}/I$.

\begin{prp}\label{indcotquot}
Let $\nC$ be an lth $\Fq$-algebra, $K$ a finite product of local fields.
If $I \subset \cO_{K}$ is an open ideal then
then the following holds:
\begin{enumerate}
\item The sequence $0 \to \nC \indcot I \to \nC \indcot \cO_{K} \to \nC \otimes \cO_{K}/I
\to 0$ is exact.

\item The natural map $(\nC \indcot \cO_{K}) \otimes_{\cO_{K}} I \to \nC \indcot \cO_{K}$
is injective with image $\nC \indcot I$.
\end{enumerate}\end{prp}

\pf Follows by the same argument as Proposition \ref{discrcomplotquot}. \quod

\breakflow
Finally we discuss some structural properties of tensor product algebras.

\begin{prp}\label{discrcomplotnoether} Let $\nC$ be a noetherian $\Fq$-algebra and $K$ a finite
product of local fields.
\begin{enumerate}
\item $\nC^\#\complot\cO_{K}$ is noetherian and complete with respect to the ideal
$\nC^\#\complot\fm_{K}$.

\item $\nC^\#\complot K$ is noetherian.
\end{enumerate}\end{prp}

\pf Without loss of generality we assume that $K$ is a local field. In this case
$K \cong k(\!(\nuCalt)\!)$ for some finite field extension $k$ of $\Fq$. (1) By
definition of the completed tensor product
\begin{equation*}
\nC^\#\complot\cO_{K} = \lim_{n \geqslant 1} \nC^\# \otimes \cO_{K}/\fm_{K}^n.
\end{equation*}
Therefore $\nC^\#\complot \cO_{K}$ is the completion of the ring $(\nC\otimes k)[\nuCalt]$
at the ideal $(\nuCalt)$. The ring $\nC\otimes k$ is of finite type over $\nC$ and so is
noetherian. Thus $(\nC\otimes k)[\nuCalt]$ is noetherian and so is its completion
$\nC^\#\complot \cO_{K}$. The completion of the ideal $(\nuCalt) \subset (\nC \otimes k)[\nuCalt]$
is clearly $\nC^\#\complot\fm_{K}$ so $\nC^\#\complot\cO_{K}$ is complete with respect
to $\nC^\#\complot\fm_{K}$. (2) follows from (1) in view of Proposition
\ref{discrcomplotloc}. \quod

\begin{prp}\label{complotstruct} Let $F$, $K$ be finite products of local fields,
$\nuC \in \cO_{F}$ and $\nuCalt \in \cO_{K}$ uniformizers.

\begin{enumerate}
\item $\cO_{F} \complot \cO_{K}$ is a finite product of complete regular
2-dimensional local rings.

\item The maximal ideals of $\cO_{F} \complot \cO_{K}$ are precisely the prime ideals
containing $\nuC$ and $\nuCalt$.
\end{enumerate}\end{prp}

\pf It is enough to assume that $F$ and $K$ are local fields.  In this case $F
\cong k_1 [\![\nuC]\!]$ and $K \cong k_2 [\![\nuCalt]\!]$ for some finite field
extensions $k_1$ and $k_2$ of $\Fq$. Therefore
\begin{equation*}
\cO_{F} \complot \cO_{K} = \lim_{n,m \geqslant 0} (k_1 \otimes k_2) [\nuC, \nuCalt]/(\nuC^n,
\nuCalt^m) = (k_1 \otimes k_2)[\![\nuC, \nuCalt]\!].
\end{equation*}
Observe that $k_1 \otimes k_2$ is a finite product of finite fields. (1) and (2)
are now clear. \quod

\section{Function spaces as modules}
\label{sec:funcspmod}

Let $\nC$, $\nCalt$ be locally compact algebras, $\nmC$ a locally compact $\nC$-module and $\nmCalt$
a locally compact $\nCalt$-module. The function spaces $a(\nmC,\nmCalt)$, $b(\nmC,\nmCalt)$ and $c(\nmC,\nmCalt)$
carry an action of $\nC$ on the right and $\nCalt$ on the left by
functoriality.
Since $\nC$ is commutative we get a natural $\nC \otimes \nCalt$-module structure.
In this section we show that it extends canonically to a structure of an
$\nC \indcot \nCalt$-module on $a(\nmC,\nmCalt)$ and an $\nC^\#\complot \nCalt$-module on
$b(\nmC,\nmCalt)$. Under certain assumptions on $\nC$ and $\nmC$ we show that it
extends to a structure of an $\nC\complot \nCalt$-module on $c(\nmC,\nmCalt)$.

\begin{dfn}\index{nidx}{topological vector spaces and modules!$V^*$, continuous dual}
Let $\nC$ be an lth algebra and $\nmC$ an lth $\nC$-module. We equip
$\nmC^*$, the continuous $\Fq$-linear dual of $\nmC$, with the $\nC$-module structure
given by the action of $\nC$ on $\nmC$.\end{dfn}

\begin{lem} If $\nmC$ is locally compact then the $\nC$-action map $\nC \times \nmC^* \to
\nmC^*$ is continuous.\end{lem}

\pf We will deduce the result from Lemma~\ref{topmod}. To do it we need to
check the following conditions:
\begin{enumerate}
\item For every $\neC \in \nC$ the induced map $\nmC^* \to \nmC^*$, $f \mapsto f \neC$ is continuous.

\item For every $f \in \nmC^*$ the induced map $\nC \to \nmC^*$, $\neC \mapsto f \neC$ is
continuous.

\item The $\nC$-action map $\nC \times \nmC^* \to \nmC^*$ is continuous at $(0,0)$.
\end{enumerate}
The condition (1) holds by functoriality of $\nmC^*$.
Let us check (3).
For a compact open subset $U \subset \nmC$ let $[U] \subset \nmC^*$ be the subspace of
functions which vanish on $U$. By Lemma
\ref{compbounded} there exists an open $\Fq$-vector subspace $V \subset \nC$
such that $V \cdot U \subset U$. As a consequence $V \cdot [U] \subset [U]$ so
the condition (3) holds.
Given $f \in \nmC^*$ there exists a compact open $U \subset \nmC$ such that $f$
vanishes on $U$. As before we can find an open $V \subset \nC$ with the property
that $V \cdot U \subset U$. Hence $V \cdot f \subset [U]$ and the condition (2)
holds as well.\qed

%
%

\begin{prp}\label{bndlcmod}%
\index{idx}{function spaces and germ spaces!bounded locally constant functions}%
\index{nidx}{function spaces and germ spaces!$a(V,W)$, bounded locally constant functions}%
Let $\nC, \nCalt$ be locally compact algebras, $\nmC$ a
locally compact $\nC$-module and $\nmCalt$ a locally compact $\nCalt$-module.
\begin{enumerate}
\item The $\nC \otimes \nCalt$-module structure on $a(\nmC,\nmCalt)$ extends uniquely to a topological
$\nC\indcot \nCalt$-module structure.

\item The map $\nmC^* \otimes \nmCalt \to a(\nmC,\nmCalt)$, $f \otimes \nmeCalt \mapsto (\nmeC \mapsto f(\nmeC)\nmeCalt)$
extends uniquely to a topological isomorphism
$\nmC^* \indcot \nmCalt \cong a(\nmC,\nmCalt)$ of $\nC \indcot \nCalt$-modules.\end{enumerate}\end{prp}

\pf The map $\nmC^* \otimes \nmCalt \to a(\nmC,\nmCalt)$ above extends to a topological isomorphism
$\nmC^* \indcot \nmCalt \cong a(\nmC,\nmCalt)$ of $\Fq$-vector spaces by Proposition \lref{bndlc}{bndlcot}.
It is $\nC \otimes \nCalt$-linear by naturality. The $\nC \otimes \nCalt$-module structure
on $\nmC^* \indcot \nmCalt$ extends canonically to a topological $\nC\indcot \nCalt$-module structure
by Proposition \ref{indcotmod} and we get the result. \quod




%

\begin{lem}\label{bndlcmoddesc}%
Let $\nC, \nCalt$ be locally compact algebras and $\nmC$ a
locally compact $\nC$-module.
If $\nmC^*$ is a topological direct summand of a free
$\nC$-module of finite rank then the natural map
$\nmC^* \otimes_{\nC} (\nC \indcot \nCalt) \to a(\nmC,\nCalt)$
is an $\nC \indcot \nCalt$-module isomorphism.\end{lem}

\pf By Proposition \ref{bndlcmod} we can identify $a(\nmC,\nCalt)$ with the completion of
an lth $\nC \otimes_{\textup{ic}} \nCalt$-module $\nmC^* \otimes_{\textup{ic}} \nCalt$.
Due to our assumption on $\nmC^*$ the module $\nmC^* \otimes_{\textup{ic}} \nCalt$ is a
topological direct summand of $(\nC \otimes_{\textup{ic}} \nCalt)^{\oplus n}$ for some
$n$. So the result is a consequence of Lemma \ref{compllattice}.\quod


\begin{prp}\label{bndmod}%
\index{idx}{function spaces and germ spaces!bounded functions}%
\index{nidx}{function spaces and germ spaces!$b(V,W)$, bounded functions}%
Let $\nC, \nCalt$ be locally compact algebras, $\nmC$ a locally
compact $\nC$-module and $\nmCalt$ a locally compact $\nCalt$-module.
\begin{enumerate}
\item The $\nC \otimes \nCalt$-module structure on $b(\nmC,\nmCalt)$ extends uniquely
to a topological $\nC^\#\complot \nCalt$-module structure.

\item The natural map $\nmC^* \otimes \nmCalt \to b(\nmC,\nmCalt)$,
$f \otimes \nmeCalt \mapsto (\nmeC \mapsto f(\nmeC) \nmeCalt)$ extends uniquely
to a topological isomorphism $(\nmC^*)^\#\complot \nmCalt \cong b(\nmC,\nmCalt)$
of $\nC^\#\complot \nCalt$-modules.\end{enumerate}\end{prp}

\pf Follows from Propositions \lref{bndfn}{bndcomplot} and \ref{modcomplot}.\quod

\begin{lem}\label{bndmoddesc}%
Let $\nC, \nCalt$ be locally compact algebras and $\nmC$ a
locally compact $\nC$-module. If $\nmC^*$ is projective of finite type as an
$\nC$-module without topology then the natural map
$\nmC^* \otimes_{\nC} (\nC^\# \complot \nCalt) \to b(\nmC,\nCalt)$
is an $\nC^\#\complot \nCalt$-module isomorphism.\end{lem}

\pf Follows from Proposition \ref{bndmod} and Lemma \ref{compllattice}. \quod 

\begin{prp}\label{contmod}%
\index{idx}{function spaces and germ spaces!continuous functions}%
\index{nidx}{function spaces and germ spaces!$c(V,W)$, continuous functions}%
Let $\nC, \nCalt$ be lth algebras, $\nmC$ an lth $\nC$-module and
$\nmCalt$ an lth $\nCalt$-module. Assume that either of the following conditions hold:
\begin{enumerate}
\item $\nC$ is discrete and $\nmC$ is compact.

\item $\nC$ is compact and $\nmC$ is discrete.
\end{enumerate}
Then the following is true:
\begin{enumerate}
\item The $\nC \otimes \nCalt$-module structure on $c(\nmC,\nmCalt)$ extends uniquely to
a topological $\nC\complot \nCalt$-module structure.

\item The natural map $\nmC^* \otimes \nmCalt \to c(\nmC,\nmCalt)$,
$f \otimes \nmeCalt \mapsto (\nmeC \mapsto f(\nmeC) \nmeCalt)$ extends uniquely to a topological
isomorphism $\nmC^* \complot \nmCalt \cong c(\nmC,\nmCalt)$ of
$\nC \complot \nCalt$-modules.
\end{enumerate}\end{prp}

\pf Follows from Pontrjagin duality (Theorem \lref{pont}{pontdual}),
Propositions \lref{contfn}{contcomplot} and \ref{modcomplot}.\quod

\begin{lem}\label{contmoddesc}Let $\nC, \nCalt$ be lth algebras and $\nmC$ an lth $\nC$-module.
Assume the following:
\begin{enumerate}
\item $\nC$ is either discrete or compact.

\item $\nmC^*$ is a topological direct summand of $\nC^{\oplus n}$ for some $n \geqslant 0$.
\end{enumerate}
Then the following holds:
\begin{enumerate}
\item If $\nC$ is discrete then $\nmC$ is compact. If $\nC$ is compact then $\nmC$ is discrete.
In particular the $\nC\otimes \nCalt$-module structure on $c(\nmC,\nCalt)$ extends uniquely to
a topological $\nC \complot \nCalt$-module structure.

\item The natural map $\nmC^* \otimes_{\nC} (\nC\complot \nCalt) \to c(\nmC,\nCalt)$ is an $\nC\complot \nCalt$-module
isomorphism.
\end{enumerate}\end{lem}

\pf Follows from Pontrjagin duality,
Proposition~\ref{contmod} and Lemma~\ref{compllattice}.\quod

\section{Germ spaces as modules}
\label{sec:germspmod}

\begin{dfn}\index{idx}{function spaces and germ spaces!the space of germs}\index{nidx}{function spaces and germ spaces!$g(V,W)$, germs}%
Let $\nC, \nCalt$ be locally compact $\Fq$-algebras, $\nmC$ a locally compact
$\nC$-module and $\nmCalt$ a locally compact $\nCalt$-module. We equip the germ space
$g(\nmC,\nmCalt)$ with an $\nC \indcot \nCalt$-module structure in the following way. Consider
the short exact sequence of Proposition \lref{germ}{germses}:
\begin{equation*}
0 \to a(\nmC,\nmCalt) \to b(\nmC,\nmCalt) \to g(\nmC,\nmCalt) \to 0.
\end{equation*}
The first arrow in this sequence is an $\nC\indcot \nCalt$-module homomorphism by
construction.
We equip $g(\nmC,\nmCalt)$ with the resulting $\nC\indcot \nCalt$-module structure.\end{dfn}

\section{\texorpdfstring{$\tau$-ring and $\tau$-module structures}
{Tau-ring and tau-module structures}}
\label{sec:ottaustruct}

Recall from Definition \lref{shtukadefs}{tauring} that a $\tau$-ring is a
ring $R$ equipped with an endomorphism $\tau\colon R \to R$. We would like to
fix $\tau$-ring structures on algebras of the form $\nC \indcot \nCalt$, $\nC \complot
\nCalt$.

\begin{dfn}\label{xottauring} Let $\nC, \nCalt$ be locally compact $\Fq$-algebras.
Let $\sigma\colon \nCalt \to \nCalt$ be the $q$-power map.
\begin{enumerate}
\item We equip $\nC \indcot \nCalt$ with the $\tau$-ring structure given by the
endomorphism $1 \indcot \sigma$.

\item Assuming $\nC \complot \nCalt$ admits the natural topological $\Fq$-algebra structure we equip
it with the $\tau$-ring structure given by the endomorphism $1 \complot \sigma$.
\end{enumerate}
\end{dfn}

\begin{lem}\label{xotfuncspmod} Let $\nC, \nCalt$ be locally compact $\Fq$-algebras and
$\nmC$ a locally compact $\nC$-module. Let $\sigma\colon \nCalt \to \nCalt$ be the $q$-power
map.
If $a(\nmC,\nCalt)$, $b(\nmC,\nCalt)$, $c(\nmC,\nCalt)$ and $g(\nmC,\nCalt)$ are equipped with endomorphisms $\tau$ given
by composition with $\sigma$ then the following is true:
\begin{enumerate}
\item $a(\nmC,\nCalt)$ is a left $\nC \indcot \nCalt\{\tau\}$-module.

\item $b(\nmC,\nCalt)$ is a left $\nC^\#\complot \nCalt\{\tau\}$-module.

\item $g(\nmC,\nCalt)$ is a left $\nC\indcot \nCalt\{\tau\}$-module.

\item $c(\nmC,\nCalt)$ is a left $\nC\complot \nCalt\{\tau\}$-module provided
the assumptions of Proposition~\ref{contmod} on $\nC$ and $\nmC$ hold.
\end{enumerate}
In all cases the $\tau$-ring structures are as in Definition \ref{xottauring}.
\end{lem}
 
\begin{dfn}\label{xotfuncspdef}Under assumptions of Lemma
\ref{xotfuncspmod} we equip the spaces $a(\nmC,\nCalt)$, $b(\nmC,\nCalt)$, $c(\nmC,\nCalt)$ and $g(\nmC,\nCalt)$ with the
$\tau$-module structures as described above. From now on we work
with only these $\tau$-module structures.\end{dfn}

\breakflow\noindent
\textit{Proof of Lemma \ref{xotfuncspmod}.}
(1) Let $f \in a(\nmC,\nCalt)$, $x \in \nC \indcot \nCalt$. We need to prove that
\begin{equation*}
\sigma \circ (x \cdot f) = \tau(x) \cdot (\sigma \circ f).
\end{equation*}
This is clear if $x \in \nC\otimes \nCalt$. As $\nC \otimes \nCalt \subset \nC \indcot \nCalt$ is
dense and $\tau\colon \nC\indcot \nCalt \to \nC\indcot \nCalt$ is continuous the general statement
follows. (2) and (4) follow in the same manner. (3) follows from (1), (2) and the short
exact sequence of Proposition \lref{germ}{germses}.\quod

\section{Residue and duality}
\label{sec:resdual}

In this section we show that in one special case the function spaces $a(\nmC,\nmCalt)$,
$b(\nmC,\nmCalt)$ and $c(\nmC,\nmCalt)$ have particularly nice module structures. It is exactly the case which 
appears in our applications.

Let $C$ be a smooth projective 
connected curve over $\Fq$.
Fix a closed point $\infty \in C$. Let $F$ be the local field of $C$ at
$\infty$, $\cO_F \subset F$ the ring of integers and
$A = \Gamma(C - \{\infty\},\,\cO_C)$
where $\cO_C$ is the structure sheaf of $C$. The
natural topology on $F$ makes it into a locally compact $\Fq$-algebra with a
compact open subalgebra $\cO_F \subset F$ and a discrete cocompact subalgebra $A
\subset F$.


Let $\Omega_C$ be the sheaf of K\"ahler differentials of $C$ over $\Fq$.
We use the following notation:
\begin{equation*}
\omega_A = \Gamma(\Spec A,\,\Omega_C), \quad
\omega_{\cO_F} = \Gamma(\Spec\cO_F, \,\Omega_C), \quad
\omega_F = \Gamma(\Spec F,\,\Omega_C).
\end{equation*}
The $F$-module $\omega_F$ carries a natural locally compact
topology with $\omega_{\cO_F} \subset \omega_F$ a compact open $\cO_F$-submodule
and $\omega_A \subset \omega_F$ a discrete cocompact $A$-submodule.
It comes equipped with a residue map $\omega_F \to k$ where $k$
is the residue field at $\infty$. We denote
$\res\colon \omega_F \to \Fq$ its composition with the trace map
$\tr\colon k \to \Fq$.
In our study we need the following
duality theorem for $\res$:

\begin{thm}\label{resdual}%
The pairing
$\omega_F \times F \to \Fq$, $(\eta,x) \mapsto\res(x\eta)$
induces the following topological isomorphisms:
\begin{equation*}
\omega_A \cong (F/A)^*, \quad
\omega_{\cO_F} \cong (F/\cO_F)^*, \quad
\omega_F \cong F^*.
\end{equation*}
\end{thm}

\pf The result is well-known. Still we sketch a proof for the
reader's convenience.
Let 
$\cO_C(1)$ be the Serre twist of $\cO_C$ by the divisor $\infty$.
Let $n \in \bZ$. A \v{C}ech computation shows that
\begin{align*}
\RGamma(C,\,\cO_C(n)) &= \Big[ A \oplus z^{-n} \cO_F \to F \Big], \\
\RGamma(C,\,\Omega_C(-n)) &= \Big[ \omega_A \oplus z^n \omega_{\cO_F} \to
\omega_F \Big]
\end{align*}
where $z \in \cO_F$ is a uniformizer and the differentials send $(x,y)$ to $x - y$.
The residue pairing $\omega_F \times F \to \Fq$ induces the following perfect pairings:
\begin{align*}
\uH^1(C,\,\Omega_C(-n)) \times \uH^0(C,\,\cO(n)) &\to \Fq, \\
\uH^0(C,\,\Omega_C(-n)) \times \uH^1(C,\,\cO(n)) &\to \Fq.
\end{align*}
Using the explicit descriptions of cohomology groups provided by the
complexes above we rewrite these pairings as
\begin{align}
\label{respairinga}
\frac{\omega_F}{\omega_A + z^n\omega_{\cO_F}} \times (A \cap z^{-n} \cO_F) &\to \Fq, \\
\label{respairingb}
[\omega_A \cap z^n\omega_{\cO_F}] \times \frac{F}{A + z^{-n} \cO_F} &\to \Fq.
\end{align}
The open subspaces $(A + z^{-n} \cO_F)/A$ form a fundamental system which covers $F/A$.
Taking the limit of \eqref{respairingb} as $n \to -\infty$ we conclude that the
residue pairing induces a topological isomorphism $\omega_A \cong (F/A)^*$.
It remains to deduce the topological isomorphisms $\omega_F \cong F^*$ and
$\omega_{\cO_F} \cong (F/\cO_F)^*$.

Let us denote $\rho\colon\omega_F \to \Hom_\Fq(F,\Fq)$ the map defined by the
residue pairing. A priori we do not even know whether its image is contained in
$F^* \subset \Hom_\Fq(F,\Fq)$.
First we prove that $\rho$ sends $z^n \omega_{\cO_F}$
to $(F/z^{-n} \cO_F)^* \subset F^*$. As $\rho$ is $F$-linear it is enough to
treat the case $n = 0$. In this case \eqref{respairinga} implies that
$\res(\eta) = 0$ for every $\eta \in \omega_{\cO_F}$. Hence $\res(x\eta) = 0$ for
all $x \in \cO_F$ and $\eta \in \omega_{\cO_F}$. We conclude that $\rho(\eta) \in
(F/\cO_F)^* \subset F^*$.

Our next step is to prove that for every $n \in \bZ$ the induced map
\begin{equation}\label{respairingsubquot}
\rho\colon \frac{z^n \omega_{\cO_F}}{z^{n+1}\omega_{\cO_F}} \to
\Big(\frac{z^{-(n+1)} \cO_F}{z^{-n}\cO_F}\Big)^*
\end{equation}
is injective. Since $\rho(\eta)(x) = \res(x\eta)$ it is enough to prove this for
a single $n \in \bZ$. As the divisor $\infty \in C$ is ample there exists an
$n \gg 0$ such that $\uH^1(C,\,\cO_C(n)) = 0$. Now \eqref{respairingb} implies that
$\omega_A \cap z^n \omega_{\cO_F} = 0$. If $\eta \in z^n \omega_{\cO_F}$ is such
that $\res(z^{-(n+1)} x \eta) = 0$ for any $x \in \cO_F^\times$ then the pairing
\eqref{respairinga} implies that $\eta \in \omega_A + z^{n+1} \omega_{\cO_F}$.
Since $\omega_A \cap z^n \omega_{\cO_F} = 0$ we conclude that $\eta \in
z^{n+1} \omega_{\cO_F}$. Whence \eqref{respairingsubquot} is injective.

At the same time \eqref{respairingsubquot} is a morphism of one-dimensional
$\cO_F/z$-vector spaces. It is therefore an isomorphism. We conclude that for
every $n > 0$ the induced map
\begin{equation*}
\rho\colon\frac{\omega_{\cO_F}}{z^n\omega_{\cO_F}} \to
\Big(\frac{z^{-n}\cO_F}{\cO_F}\Big)^*
\end{equation*}
is an isomorphism. As $\omega_{\cO_F}$ is complete it follows that
$\rho\colon\omega_{\cO_F} \to (F/\cO_F)^*$ is a topological isomorphism. Since
the open subspaces
$z^n \omega_{\cO_F}$ cover $\omega_F$ we deduce that
$\rho\colon\omega_F \to F^*$ is a topological isomorphism. \quod

\begin{cor}\label{funcmoddesc}%
Let $\nmCalt$ be a locally compact module over a locally compact algebra
$\nCalt$.
The map
$\omega \otimes \nmeCalt \mapsto (x \mapsto \res(x \omega) \nmeCalt)$
extends uniquely to the following topological isomorphisms:
\begin{itemize}
\item an isomorphism
$\omega_A \otimes \nmCalt \xrightarrow{\isosign} a(F/A,\nmCalt)$
of $A \otimes \nCalt$-modules.

\item an isomorphism
$\omega_A \complot \nmCalt \xrightarrow{\isosign} b(F/A,\nmCalt)$
of $A \complot \nCalt$-modules.

\item an isomorphism
$\omega_F \indcot \nmCalt \xrightarrow{\isosign} a(F,\nmCalt)$
of $F \indcot \nCalt$-modules.

\item an isomorphism
$\omega_F^\#\complot \nmCalt \xrightarrow{\isosign} b(F,\nmCalt)$
of $F^\# \complot \nCalt$-modules.

\item an isomorphism
$\omega_{\cO_F}\indcot \nmCalt \xrightarrow{\isosign} a(F/\cO_F,\nmCalt)$
of $\cO_F \indcot \nCalt$-modules.

\item an isomorphism
$\omega_{\cO_F}\complot \nmCalt \xrightarrow{\isosign} c(F/\cO_F,\nmCalt)$
of $\cO_F \complot \nCalt$-modules.
\end{itemize}\end{cor}

\begin{rmk}Note that $c(F/A,\nmCalt) = b(F/A,\nmCalt)$ and $b(F/\cO_F,\nmCalt) = a(F/\cO_F,\nmCalt)$.\end{rmk}

\afterall\noindent\textit{Proof of Corollary \ref{funcmoddesc}.} In view of Theorem
\ref{resdual} it follows from Propositions \ref{bndlcmod}, \ref{bndmod} and
\ref{contmod}.\quod

\breakflow
Let $R_0 \to R$ be a ring homomorphism and $\nmC$ an $R$-module.
Recall that an $R_0$-submodule $\nmC_0 \subset \nmC$
is called a lattice if the natural map $R \otimes_{R_0} \nmC_0 \to \nmC$
is an isomorphism (Section~\ref{sec:lattices} in the chapter ``\chnotconv'').

\begin{cor}\label{funcomegalat}%
Let $\nCalt$ be a locally compact algebra.
\begin{enumerate}
\item $\omega_A$ is an $A$-lattice in $a(F/A,\nCalt)$ and $b(F/A,\nCalt)$.

\item $\omega_F$ is an $F$-lattice in $a(F,\nCalt)$ and $b(F,\nCalt)$.

\item $\omega_{\cO_F}$ is an $\cO_F$-lattice in $a(F/\cO_F,\nCalt)$
and $c(F/\cO_F,\nCalt)$.\end{enumerate}\end{cor}

\pf Theorem \ref{resdual} identifies $\omega_A$ with
$(F/A)^* \subset a(F/A,\nCalt) \subset b(F/A,\nCalt)$.
As $\omega_A$ is a locally free module over a discrete ring $A$
the result (1) follows from Lemmas \ref{bndlcmoddesc}
and \ref{bndmoddesc}. Similarly
(2) and (3) follow from Lemmas \ref{bndlcmoddesc}, \ref{bndmoddesc} and \ref{contmoddesc}.\quod

\begin{cor}\label{funcmodstruct}%
Let $\nCalt$ be a locally compact algebra.
\begin{enumerate}
\item $a(F/A,\nCalt)$ is a locally free $A \otimes \nCalt$-module of rank $1$.

\item $b(F/A,\nCalt)$ is a locally free $A \complot \nCalt$-module of rank $1$.

\item $a(F,\nCalt)$ is a free $F \indcot \nCalt$-module of rank $1$.

\item $b(F,\nCalt)$ is a free $F^\#\complot \nCalt$-module of rank $1$.

\item $a(F/\cO_F,\nCalt)$ is a free $\cO_F\indcot \nCalt$-module of rank $1$.

\item $c(F/\cO_F,\nCalt)$ is a free $\cO_F\complot \nCalt$-module of rank $1$.
\end{enumerate}\end{cor}

\pf Follows from Corollary \ref{funcomegalat} since
$\Omega_C$ is a locally free sheaf of rank $1$ on $C$.\quod

\begin{cor}\label{funclat}%
Let $\nCalt$ be a locally compact algebra.
\begin{enumerate}
\item $a(F/A,\nCalt)$ is an $A \otimes \nCalt$-lattice
in $b(F/A,\nCalt)$ and $a(F,\nCalt)$.

\item $a(F/\cO_F,\nCalt)$ is an $\cO_F\indcot \nCalt$-lattice
in $a(F,\nCalt)$ and $c(F/\cO_F,\nCalt)$.

\item $a(F,\nCalt)$ is an $F\indcot \nCalt$-latticein $b(F,\nCalt)$.

\item $b(F/A,\nCalt)$ is an $A\complot \nCalt$-lattice in $b(F,\nCalt)$.
\end{enumerate}\end{cor}

\pf Follows from Corollary \ref{funcomegalat}.\quod

\chapter{Cohomology of shtukas}
\label{chapter:shtcoh}
\label{ch:germcoh}
\label{ch:gcmp}
\label{ch:globcoh}
\label{ch:lcmp}
\label{ch:csc}
\label{ch:locglobcmp}
\label{ch:cechcoh}
\label{ch:complcech}
\label{ch:globcoeffch}
\label{ch:globzeta}

Fix a locally compact noetherian $\Fq$-algebra $\ngC$ and a smooth projective curve
$X$ over $\Fq$. We call $\ngC$ the \emph{coefficient ring}
and $X$ the \emph{base curve}. Set $\ngS = \Spec \ngC$ and
consider the product $\ngS \times X$.
We equip $\ngS \times X$ with the
$\tau$-scheme structure given by the endomorphism which acts as the identity on
$\ngS$ and as the $q$-power map on $X$.

In this chapter we study the cohomology of locally free shtukas on $\ngS \times X$.
The basis of our approach is a \v{C}ech method described in
Section~\ref{sec:cechcoh}. We also introduce and study a few supplementary
constructions.
The reader should be warned that some of them
will not reappear until the last chapter of the book. They are the
compactly supported cohomology functor, the global germ map and the
local-global compatibility theorem.

The \v{C}ech method presented here
involves a choice of additional data, the points ``at infinity'' on $X$.
So let us fix finitely many closed points $x_1, \dotsc, x_n$ of $X$.
The complement of $\{x_1,\dotsc,x_n\}$ in $X$ is an affine subscheme
which we denote $Y = \Spec R$.
The product
of the local fields of $X$ at $x_1,\dotsc,x_n$ is denoted $K$
with 
$\cO_K \subset K$ the ring
of integers and $\fm_K \subset \cO_K$ the Jacobson
radical. By construction $\Spec\cO_K/\fm_K = \{x_1,\dotsc,x_n\}\subset X$.
The natural topology on $\cO_K$ makes it into a compact open $\Fq$-subalgebra
of a locally compact $\Fq$-algebra $K$.

We begin with a
study of shtuka cohomology around the points ``at infinity''.
Let $\cM$ be a locally free shtuka on $\ngC \indcot \cO_K$.
In Section~\ref{sec:germcoh} we introduce the germ cohomology complex
$\RGammag(\ngC \indcot K,\,\cM)$. As suggested by the notation it depends only on
the restriction of $\cM$ to $\ngC \indcot K$. The germ cohomology is modelled on
the germ spaces of Section \lref{germ}{sec:germ}.
With some degree of caution it can be regarded as
compactly supported cohomology
for shtukas on $\ngC \indcot K$ with respect to the compactification given by the
ring $\ngC \indcot \cO_K$.

The germ cohomology is related to the usual cohomology
via the local germ map
\begin{equation*}
\RGamma(\ngC \indcot\cO_K,\,\cM) \xrightarrow{\,\,\isosign\,\,} \RGammag(\ngC \indcot K,\,\cM)
\end{equation*}
which we construct in Section \ref{sec:locgermmap}.
This map is defined only if
$\cM(\ngC \otimes \cO_K/\fm_K)$ is nilpotent and is always a quasi-\hspace{0pt}isomorphism.
The local germ map will play a role in Chapters \ref{chapter:mothomsht} -- \ref{chapter:cnf}.

Starting from Section \ref{sec:cechcoh} we shift the focus to
the global setting.
Let $\cM$ be a locally free shtuka on $\ngS \times X$.
We introduce a
\v{C}ech method which computes $\RGamma(\ngS \times X,\,\cM)$ as the mapping fiber
\begin{equation*}
\Big[ \RGamma(\ngC \otimes R,\,\cM) \oplus
\RGamma(\ngC^\#\complot\cO_K,\,\cM) \xrightarrow{\,\,\textrm{difference}\,\,}
\RGamma(\ngC^\#\complot K,\,\cM) \Big].
\end{equation*}
This is our main tool to handle the cohomology of shtukas on $\ngS\times X$.

The compactly supported cohomology functor $\RGammac(\ngS \times Y,\,\cM)$ is
introduced in Section \ref{sec:cscoh}. As suggested by the notation
$\RGammac(\ngS \times Y,\,\cM)$
depends only on the restriction of $\cM$ to $\ngS \times Y$. It comes equipped with a
natural map $\RGamma(\ngS \times X,\,\cM) \to \RGammac(\ngS \times Y,\,\cM)$. We prove
that this map is a quasi-isomorphism if $\cM(\ngC \otimes \cO_K/\fm_K)$ is
nilpotent.
One can interpet the nilpotence condition as saying that
$\cM$ is an extension by zero of a shtuka on $\ngS \times Y$.
We use $\RGammac(\ngS \times Y,\,\cM)$ to construct the global germ map
\begin{equation*}
\RGamma(\ngS \times X,\,\cM) \to \RGammag(\ngC \indcot K,\,\cM).
\end{equation*}
Similarly to its local counterpart the global germ map is defined under
assumption that $\cM(\ngC \otimes \cO_K/\fm_K)$ is nilpotent.
However it is not a quasi-isomorphism in general.

Section \ref{sec:locglobcmp} is devoted to the proof of 
Theorem \ref{locglobcmp}, a compatibility statement for the local and global
germ maps.
This statement will be used in Chapter \ref{chapter:cnf} in the proof of the
class number formula.

In Section \ref{sec:complcech} we present an advanced version of the
\v{C}ech method for shtukas on $\Spec\ngOF \times X$ where $\ngOF$ is the ring of
integers of a local field $\ngF$. 
This method enables us to prove the following: if $\cM(\ngOF/\fm_\ngF \otimes R)$
is nilpotent then the natural map
\begin{equation*}
\RGamma(\Spec\ngOF \times X,\,\cM) \to \RGamma(\ngOF\complot\cO_K,\,\cM)
\end{equation*}
is a quasi-\hspace{0pt}isomorphism (Theorem \ref{intconc}).
Informally speaking, the
cohomology of $\cM$ concentrates on $\ngOF\complot\cO_K$. This phenomenon 
is important to the theory of regulator developed in
Chapters \ref{chapter:reg} and \ref{chapter:trace}.

Finally in Sections \ref{sec:globcoeffch} and \ref{sec:globzeta} we study how
the change of the coefficient ring $S$ reflects on shtuka cohomology and 
$\zeta$-isomorphisms. The results of these sections will be used
in Chapters \ref{chapter:trace} and \ref{chapter:cnf}.



The cohomology functors 
in this chapter
are typically given by mapping fibers 
\begin{equation*}
\Big[\RGamma(\mathcal{X},\,\cM)  \to \RGamma(\mathcal{X}',\,\cM)\Big]
\end{equation*}
where $\mathcal{X}$ and $\mathcal{X}'$ are affine schemes
(for the mapping fiber see Definition \ref{mappingfiber} in Chapter ``\chnotconv'').
A~word of warning about them is necessary.
The complexes $\RGamma(\mathcal{X},\,\cM)$ and $\RGamma(\mathcal{X}',\,\cM)$ are well-defined only as
objects in the derived category. As a consequence one gets a problem
with functoriality.
If $\cM \to \cN$ is a morphism of shtukas then the induced maps
\begin{equation}\tag{$\ast$}\label{conetrouble}
\begin{array}{l@{}l}
\RGamma(\mathcal{X},\,\cM)\,\,&\to \RGamma(\mathcal{X},\,\cN), \\
\RGamma(\mathcal{X}',\,\cM)\,\,&\to \RGamma(\mathcal{X}',\,\cN)
\end{array}
\end{equation}
do not determine a unique morphism of the
mapping fibers. The reason is that this morphism 
depends on the choice of non-derived representatives for the maps
\eqref{conetrouble}. 

We solve this problem in the following way.
Since the schemes $\mathcal{X}$ and $\mathcal{X}'$ are affine Theorem \lref{shtukacoh}{shtaffcoh} provides us with
canonical non-derived representatives for the complexes $\RGamma(\mathcal{X},\,\cM)$ and
$\RGamma(\mathcal{X}',\,\cM)$, namely the associated complexes $\cGamma(\mathcal{X},\,\cM)$ and
$\cGamma(\mathcal{X}',\,\cM)$.
The mapping fiber construction is
functorial on the level of the non-derived category of complexes.
Following the convention of Section \ref{ch:shtukacoh}.\ref{sec:shtaffcoh} we identify
$\RGamma(\mathcal{X},\,\cM)$ with $\cGamma(\mathcal{X},\,\cM)$ and $\RGamma(\mathcal{X}',\,\cM)$ with
$\cGamma(\mathcal{X}',\,\cM)$. We thus regain the functoriality.


\section{Germ cohomology}
\label{sec:germcoh}

In this section we fix locally compact $\Fq$-algebras $\ngC$ and $\ngB$.
Following the conventions of Section
\lref{ottaustruct}{sec:ottaustruct} we equip
$\ngC \indcot \ngB$ and $\ngC^\#\complot \ngB$ with the $\tau$-ring structures given by the
endomorphisms which act as identity on $\ngC$ and as the $q$-power map on $\ngB$.

\begin{dfn}\label{defgermcoh}\index{idx}{shtuka cohomology!germ cohomology}\index{nidx}{shtuka cohomology!zRGammag@$\RGammag$, germ cohomology}%
Let $\cM$ be a quasi-coherent shtuka on
$\ngC \indcot \ngB$. The \emph{germ cohomology} complex of $\cM$ is the
$\ngC$-module complex
\begin{equation*}
\RGammag(\ngC \indcot \ngB, \cM) =
\Big[ \RGamma(\ngC\indcot \ngB,\,\cM) \to \RGamma(\ngC^\#\complot \ngB,\,\cM) \Big].
\end{equation*}
The differential in this complex is induced by the natural map $\ngC \indcot \ngB \to
\ngC^\#\complot \ngB$
which is the completion of the continuous bijection
$\ngC \otimes_{\textup{ic}} \ngB \to \ngC^\#\otimes_{\textup{c}} \ngB$.
The $n$-th cohomology group of $\RGammag(\ngC \indcot \ngB, \cM)$ is denoted
$\uH^n_g(\ngC \indcot \ngB,\,\cM)$.
\end{dfn}

\begin{prp}\label{locfreegermcoh} If $\cM$ is a locally free shtuka on $\ngC
\indcot \ngB$ then the natural map
\begin{equation*}
\RGammag(\ngC \indcot \ngB, \,\cM) \to
\RGamma\!\Big(\ngC \indcot \ngB, \,\tfrac{\cM(\ngC^\#\complot \ngB)}{\cM(\ngC \indcot \ngB)}\Big)[-1]
\end{equation*}
is a quasi-isomorphism.\end{prp}
 
\pf Tensoring the short exact sequence
\begin{equation*}
0 \to \ngC \indcot \ngB \to \ngC^\#\complot \ngB \to \frac{\ngC^\#\complot \ngB}{\ngC \indcot \ngB} \to
0
\end{equation*}
with a locally free $\ngC \indcot \ngB$-module of finite rank we get a short exact
sequence. As $\cM$ is locally free the claim follows. \quod

\begin{prp}\label{germcohlociso} Let $f\colon \ngC_1 \to \ngC_2$ and $g\colon \ngB_1 \to
\ngB_2$ be continuous homomorphisms of locally compact $\Fq$-algebras. Let $\cM$ be
a locally free shtuka on $\ngC_1 \indcot \ngB_1$. If $f^*$ and $g$ are local
isomorphisms of topological $\Fq$-vector spaces then the natural map
\begin{equation*}
\RGammag(\ngC_1 \indcot \ngB_1, \,\cM) \to \RGammag(\ngC_2 \indcot \ngB_2, \,\cM)
\end{equation*}
induced by $f \indcot g$ is a quasi-isomorphism.\end{prp}

\pf By Proposition \lref{germ}{germquotiso} the maps $f$ and $g$ induce a bijection
\begin{equation*}
\frac{\ngC_1^\#\complot \ngB_1}{\ngC_1\indcot \ngB_1} \cong
\frac{\ngC_2^\#\complot \ngB_2}{\ngC_2\indcot \ngB_2}.
\end{equation*}
As the shtuka $\cM$ is locally free it follows that $f \indcot g$ induces an
isomorphism of shtukas
\begin{equation*}
\frac{\cM(\ngC_1^\#\complot \ngB_1)}{\cM(\ngC_1\indcot \ngB_1)} \cong
\frac{\cM(\ngC_2^\#\complot \ngB_2)}{\cM(\ngC_2\indcot \ngB_2)}.
\end{equation*}
The result now follows from Proposition \ref{locfreegermcoh}. \quod

\section{Local germ map}
\label{sec:locgermmap}

Fix a noetherian locally compact $\Fq$-algebra $\ngC$. In applications this algebra will
usually be a local field. Let $K$ be a finite product of local
fields, $\cO_K \subset K$ the ring of integers and $\fm_K \subset \cO_K$ the Jacobson radical.

The ideal $\fm_K \subset \cO_K$ is open so that we have a natural map
$\ngC^\#\complot\cO_K \to \ngC \otimes \cO_K/\fm_K$. Taking the completion of $\ngC
\otimes_{\textup{ic}} \cO_K \to \ngC \otimes_{\textup{ic}} \cO_K/\fm_K$ we get a
natural map $\ngC \indcot \cO_K \to \ngC \otimes \cO_K/\fm_K$ since $\ngC
\otimes_{\textup{ic}} \cO_K/\fm_K$ is discrete.

\begin{prp}\label{loccs} Let $\cM$ be a locally free shtuka on $\ngC^\# \complot
\cO_K$. If $\cM(\ngC\otimes\cO_K/\fm_K)$ is nilpotent then
$\RGamma(\ngC^\#\complot\cO_K, \,\cM) = 0$.\end{prp}

\pf According to Proposition \lref{otalgprops}{discrcomplotnoether} the ring
$\ngC^\#\complot\cO_K$ is noetherian and complete with respect to the ideal
$\ngC^\#\complot\fm_K$. By Proposition \lref{otalgprops}{discrcomplotquot} the
natural map $\ngC^\#\complot\cO_K \to \ngC \otimes \cO_K/\fm_K$ is surjective with
kernel $\ngC^\#\complot\fm_K$.
So the result follows from Proposition~\lref{nilp}{nilpcomp}. \quod

\begin{lem}\label{rgammaok} If $\cM$ is a locally free shtuka on $\ngC\indcot\cO_K$
then the natural map $\RGammag(\ngC\indcot\cO_K, \,\cM) \to \RGammag(\ngC\indcot K, \,\cM)$
is a quasi-isomorphism.\end{lem}

\pf The inclusion $\cO_K \hookrightarrow K$ is a local isomorphism of topological
$\Fq$-vector spaces. So the result is a consequence of Proposition
\ref{germcohlociso}. \quod

\breakflow
Let $\cM$ be a locally free shtuka on $\ngC \indcot \cO_K$. 
According to Definition \ref{defgermcoh}
\begin{equation*}
\RGammag(\ngC \indcot \cO_K, \cM) =
\Big[ \RGamma(\ngC \indcot \cO_K,\,\cM) \to \RGamma(\ngC^\#\complot\cO_K,\,\cM) \Big].
\end{equation*}
The projection to the first argument of the mapping fiber 
defines a natural map
\begin{equation*}
\RGammag(\ngC\indcot\cO_K,\,\cM) \to \RGamma(\ngC \indcot \cO_K, \cM).
\end{equation*}
Taking its composition with the quasi-isomorphism
$\RGammag(\ngC\indcot K,\,\cM) \cong \RGammag(\ngC\indcot\cO_K,\,\cM)$ of Lemma
\ref{rgammaok} we obtain a map
\begin{equation}\label{locinvmap}
\RGammag(\ngC\indcot K,\,\cM) \to \RGamma(\ngC \indcot \cO_K,\,\cM).
\end{equation}

\begin{lem}\label{locinvmapiso} Let $\cM$ be a locally free shtuka on
$\ngC \indcot \cO_K$. If $\cM(\ngC \otimes \cO_K/\fm_K)$ is nilpotent then the
natural map \eqref{locinvmap} is a quasi-isomorphism.\end{lem}

\pf By construction the natural map $\RGammag(\ngC\indcot\cO_K,\,\cM) \to
\RGamma(\ngC\indcot\cO_K,\,\cM)$ extends to a distinguished triangle
\begin{equation*}
\RGammag(\ngC\indcot\cO_K,\,\cM) \to \RGamma(\ngC\indcot\cO_K,\,\cM) \to
\RGamma(\ngC^\#\complot\cO_K,\,\cM) \to [1].
\end{equation*}
Together with the quasi-isomorphism $\RGammag(\ngC\indcot\cO_K,\,\cM) \cong
\RGammag(\ngC\indcot K,\,\cM)$ it gives us a distinguished triangle
\begin{equation*}
\RGammag(\ngC \indcot K,\,\cM) \to \RGamma(\ngC\indcot \cO_K,\,\cM) \to
\RGamma(\ngC^\#\complot\cO_K,\,\cM) \to [1].
\end{equation*}
Proposition \ref{loccs} shows that $\RGamma(\ngC^\#\complot\cO_K,\,\cM) = 0$, so the
result follows. \quod

\begin{dfn}\index{idx}{shtuka cohomology!local germ map}\label{defloccmp}%
Let $\cM$ be a locally free shtuka on $\ngC \indcot \cO_K$ such that
$\cM(\ngC\otimes\cO_K/\fm_K)$ is nilpotent. By Lemma \ref{loccs} the natural map
\eqref{locinvmap} is a quasi-\hspace{0pt}isomorphism. We define \emph{the local germ map}
\begin{equation*}
\RGamma(\ngC\indcot \cO_K,\,\cM) \xrightarrow{\,\,\isosign\,\,} \RGammag(\ngC\indcot K,\,\cM)
\end{equation*}
as its inverse. The adjective ``local'' signifies that it involves a shtuka
defined over a semil-local ring $\cO_K$. Observe that the local germ map is a
quasi-isomorphism by construction.\end{dfn}

\begin{prp}\label{loccoh} Let $\cM$ be a locally free shtuka on $\ngC\indcot\cO_K$.
If $\cM(\ngC\otimes\cO_K/\fm_K)$ is nilpotent then $\RGamma(\ngC\indcot\cO_K,\,\cM)$ and
$\RGammag(\ngC \indcot K,\,\cM)$ are concentrated in degree~$1$.\end{prp}

\pf The complex $\RGamma(\ngC\indcot\cO_K,\,\cM)$ is concentrated in degrees $0$ and
$1$ since $\Spec(\ngC \indcot \cO_K)$ is affine. The complex $\RGammag(\ngC \indcot K,\,\cM)$
is concentrated in degrees $1$ and $2$ by Proposition \ref{locfreegermcoh}. As
these complexes are quasi-isomorphic via the local germ map the conclusion
follows.\quod

\begin{prp}\label{loccmpdesc} Let $\cM = [\cM_0 \shtuka{i}{j} \cM_1]$ be a
locally free shtuka on $\ngC\indcot\cO_K$ and let $x \in \cM_1$.
Assume that $\cM(\ngC \otimes \cO_K/\fm_K)$ is nilpotent.
\begin{enumerate}
\item There exists a unique $y \in \cM_0(\ngC^\#\complot\cO_K)$ such that $(i-j)(y)
= x$.

\item Consider the composition
\begin{equation*}
\uH^1(\ngC \indcot\cO_K,\,\cM) \xrightarrow{\quad\textup{local}\quad}
\uH^1_g(\ngC \indcot K,\,\cM) \xrightarrow{\,\,\isosign\,\,}
\uH^0\Big(\ngC \indcot K, \tfrac{\cM(\ngC^\#\complot K)}{\cM(\ngC\indcot K)}\Big)
\end{equation*}
of the local germ map and the natural isomorphism of Proposition
\ref{locfreegermcoh}. This composition sends the class of $x$ to the image of
$y$ in the quotient $\cM_0(\ngC^\#\complot K)/\cM_0(\ngC\indcot K)$.
\end{enumerate}\end{prp}

%

\pf (1) $\RGamma(\ngC^\#\complot\cO_K,\,\cM)$ is represented by the complex
\begin{equation*}
\Big[ \cM_0(\ngC^\#\complot\cO_K) \xrightarrow{\,\,i-j\,\,}
\cM_1(\ngC^\#\complot\cO_K) \Big].
\end{equation*}
By Proposition \ref{loccs} $\RGamma(\ngC^\#\complot\cO_K,\,\cM) = 0$. So the map
$i-j$ in the complex above is a bijection and (1) follows.

(2) Consider the maps
\begin{align*}
\uH^1_g(\ngC\indcot\cO_K,\,\cM) &\to \uH^1(\ngC\indcot\cO_K,\,\cM),\\
\uH^1_g(\ngC\indcot\cO_K,\,\cM) &\to
\uH^0\Big(\ngC\indcot\cO_K,\tfrac{\cM(\ngC^\#\complot\cO_K)}{\cM(\ngC\indcot\cO_K)}\Big)
\end{align*}
determined by the natural maps of complexes
\begin{align*}
\RGammag(\ngC\indcot\cO_K,\,\cM) &\to
\RGamma(\ngC\indcot\cO_K,\,\cM), \\
\RGammag(\ngC\indcot\cO_K,\,\cM) &\to
\RGamma\!\Big(\ngC\indcot\cO_K,\tfrac{\cM(\ngC^\#\complot\cO_K)}{\cM(\ngC\indcot\cO_K)}\Big)[-1].
\end{align*}
In order to prove (2) it is enough to produce a cohomology class $h \in
\uH^1_g(\ngC\indcot\cO_K,\,\cM)$ which maps to the class of $x$ in
$\uH^1(\ngC\indcot\cO_K,\,\cM)$ and to the class of $y$ in
$\cM_0(\ngC^\#\complot\cO_K)/\cM_0(\ngC \indcot\cO_K)$.

By definition $\RGammag(\ngC\indcot\cO_K,\,\cM)$ is represented by the total complex
of the double complex
\begin{equation*}
\xymatrix{
\cM_1(\ngC\indcot\cO_K) \ar[r] & \cM_1(\ngC^\#\complot\cO_K) \\
\cM_0(\ngC\indcot\cO_K) \ar[u]^{i-j} \ar[r] & \cM_0(\ngC^\#\complot\cO_K) \ar[u]_{j-i}
}
\end{equation*}
The element $(x,y) \in \cM_1(\ngC\indcot\cO_K) \oplus \cM_0(\ngC^\#\complot\cO_K)$
is a $1$-cocyle in the total complex since $x + (j-i)(y) = 0$ by definition
of $y$. By construction $(x,y)$ maps to the class of $x$ in
$\uH^1(\ngC\indcot\cO_K,\,\cM)$ and to the class of $y$ in the quotient
$\cM_0(\ngC^\#\complot\cO_K)/\cM_0(\ngC\indcot\cO_K)$. Thus (2) follows. \quod

\section{\texorpdfstring{\v{C}ech cohomology}{Cech cohomology}}
\label{sec:cechcoh}

In this section we work with a smooth projective curve $X$ over $\Fq$ and
a coefficient algebra $\ngC$ as
described in the introduction. We set $\ngS = \Spec\ngC$.
Our goal is to present a \v{C}ech method for
computing the cohomology of shtukas on $\ngS \times X$.

The coefficient algebra $\ngC$ is assumed to be \emph{noetherian}.
As usual the $\tau$-structures on the tensor product rings $\ngC \otimes R$,
$\ngC^\#\complot\cO_K$, $\ngC^\#\complot K$ and on the scheme
$\ngS \times X = \Spec\ngC \times X$ are given by endomorphisms which act
as the identity on $\ngC$ and as the $q$-power map on the other factor.

\begin{dfn}\label{defcechcoh}\index{idx}{shtuka cohomology!\v{C}ech cohomology}\index{nidx}{shtuka cohomology!RGammaA@$\RvGamma$, \v{C}ech cohomology}%
Let $\cM$ be a quasi-coherent shtuka on $\ngS \times X$.
The \emph{\v{C}ech chomology complex} of $\cM$ is the $\ngC$-module complex
\begin{equation*}
\RvGamma(\ngS\times X,\,\cM) =
\Big[ \RGamma(\ngC\otimes R, \,\cM) \oplus
\RGamma(\ngC^\#\complot \cO_K, \,\cM) \to \RGamma(\ngC^\#\complot K,\,\cM) \Big]
\end{equation*}
where the unlabelled morphism
is the difference of the natural maps.
The $n$-th cohomology group of this complex is
denoted $\check{\uH}^n(\ngS\times X, \,\cM)$.%
\end{dfn}


%
%
%
\begin{lem}\label{cechsquare}The natural commutative square
\begin{equation*}
\xymatrix{
\Spec(\ngC^\#\complot K) \ar[r]^{\iota'} \ar[d]^{f'} & \Spec(\ngC^\#\complot \cO_K) \ar[d]^f \\
\Spec(\ngC \otimes R) \ar[r]^\iota & \ngS\times X
}
\end{equation*}
is cartesian. Furthermore $\Spec(\ngC \otimes R)$ and $\Spec(\ngC^\#\complot\cO_K)$
form a flat covering of $\ngS\times X$.\end{lem}

\pf 
Proposition \lref{otalgprops}{discrcomplotloc} implies that the square is
cartesian.
The complement of $\Spec(\ngC \otimes R)$ in $\ngS\times X$
is $\Spec(\ngC \otimes \cO_K/\fm_K)$ so the images of $\Spec(\ngC \otimes R)$ and
$\Spec(\ngC^\#\complot\cO_K)$ cover $\ngS\times X$. 
It remains to prove that $\Spec(\ngC^\#\complot\cO_K)$
is flat over $\ngS\times X$.

Pick an affine
open subscheme $\Spec R' \subset X$ which contains $\Spec\cO_K/\fm_K$. 
Shrinking $\Spec R'$ if necessary we can find an element $r' \in R'$ which
is a uniformizer of $\cO_K$. By Proposition
\lref{otalgprops}{discrcomplotnoether}
the ring $\ngC^\#\complot\cO_K$ is complete with respect to the ideal
$\ngC^\#\complot\fm_K$.
This ideal is generated by $\fm_K$ according to
Proposition \lref{otalgprops}{discrcomplotquot}.
As $r'$ is a generator of $\fm_K$ we deduce that 
$\ngC^\#\complot\cO_K$ is the completion of $\ngC \otimes R'$ with respect to $\ngC
\otimes r'R'$.
Now the fact that $\ngC \otimes R'$ is noetherian implies that
$\ngC^\#\complot\cO_K$ is flat over $\ngC \otimes R'$ and therefore over $\ngS\times X$.\quod

\breakflow
Let $\cF$ be a quasi-coherent sheaf on $\ngS\times X$. We define a complex of sheaves
on $\ngS\times X$:
\begin{equation*}
\cC(\cF) =
\Big[\iota_\ast\iota^\ast\cF \oplus f_\ast f^\ast \cF \to
g_\ast g^\ast \cF\Big]
\end{equation*}
where $g\colon \Spec(\ngC^\#\complot K) \to \ngS\times X$ is the natural map
and the differential is the difference of the natural maps as in the definition
of $\RvGamma$. The sum of adjunction units provides us with a natural morphism
$\cF[0] \to \cC(\cF)$.


\begin{lem}\label{fcechlemb} If $\cF$ is a quasi-coherent sheaf on $\ngS\times X$
then the natural map $\cF[0] \to \cC(\cF)$ is a quasi-isomorphism.\end{lem}

\pf We first show that
natural sequence
\begin{equation}\label{cechoses}
0 \to \cO_{\ngS\times X} \to \iota_\ast\iota^\ast\cO_{\ngS\times X} \oplus f_\ast
f^\ast \cO_{\ngS\times X} \to g_\ast g^\ast \cO_{\ngS\times X} \to 0
\end{equation}
is  exact.
As the
commutative diagram of Lemma \ref{cechsquare} is cartesian and
the morphism
$f\colon \Spec(\ngC^\#\complot\cO_K) \to \ngS\times X$
is affine the pullback of
\eqref{cechoses} to $\Spec(\ngC^\#\complot\cO_K)$ is
\begin{equation*}
0 \to \ngC^\#\complot\cO_K \xrightarrow{(1,\iota')} (\ngC^\#\complot\cO_K)\oplus(\ngC^\#\complot K)
\xrightarrow{(\iota',-1)} \ngC^\#\complot K \to 0.
\end{equation*}
This sequence is clearly exact. The same argument shows that the pullback
of \eqref{cechoses} to $\Spec(\ngC \otimes R)$ is exact. As $\Spec(\ngC\otimes
R)$ and $\Spec(\ngC^\#\complot\cO_K)$ form a flat covering of $\ngS\times X$ it
follows that \eqref{cechoses} is exact.

Now let $\cF$ be a quasi-coherent sheaf on the scheme $\ngS\times X$. Consider the morphism
$g\colon\Spec(\ngC^\#\complot K) \to \ngS\times X$.
As $\ngS\times X$ is
separated over $\Fq$ the morphism $g$ is affine. Thus
the natural map
\begin{equation*}
\big(g_\ast \cO_{\Spec\ngC^\#\complot K}\big)\otimes_{\cO_{\ngS\times X}} \cF \to g_\ast g^\ast \cF
\end{equation*}
is an isomorphism. The same argument applies to the maps $\iota$ and $f$. We
conclude that
\begin{equation*}
\cC(\cO_{\ngS\times X}) \otimes_{\cO_{\ngS\times X}} \cF =
\cC(\cF).
\end{equation*}
Consider the distinguished triangle
\begin{equation*}
\cO_{\ngS\times X}[0] \to \cC(\cO_{\ngS\times X}) \to C \to [1]
\end{equation*}
extending the natural quasi-isomorphism $\cO_{\ngS\times X}[0] \to
\cC(\cO_{\ngS\times X})$. Applying the functor $-\otimes_{\cO_{\ngC\times
X}} \cF$ we obtain a distinguished triangle
\begin{equation*}
\cF[0] \to \cC(\cF) \to C \otimes_{\cO_{\ngS\times X}} \cF \to [1]
\end{equation*}
where the first arrow is the natural map $\cF[0] \to \cC(\cF)$.
By construction $C$ is a bounded acyclic complex of flat
$\cO_{\ngS\times X}$-modules. Hence the complex $C \otimes_{\cO_{\ngS\times X}} \cF$ is
acyclic and the first arrow in the triangle above is a quasi-isomorphism. \quod

\begin{dfn}\label{fcechmor} Let $\cF$ be a quasi-coherent sheaf on $\ngS\times X$.
\begin{enumerate}
\item
The \emph{\v{C}ech cohomology complex} of $\cF$ is the $\ngC$-module complex
\begin{equation*}
\RvGamma(\ngS\times X,\,\cF) = \Gamma(\ngS\times X,\,\cC(\cF)).
\end{equation*}

\item
We define a natural map
$\RvGamma(\ngS\times X,\,\cF) \to \RGamma(\ngS\times X,\,\cF)$
as the composition 
\begin{equation*}
\Gamma(\ngS\times X,\,\cC(\cF)) \to
\RGamma(\ngS\times X,\,\cC(\cF)) \xleftarrow{\,\isosign\,}
\RGamma(\ngS\times X,\,\cF)
\end{equation*}
of the natural map $\Gamma \to \RGamma$ and 
the quasi-isomorphism provided by Lemma \ref{fcechlemb}.
\end{enumerate}
More explicitly
\begin{equation*}
\RvGamma(\ngS\times X,\,\cF) =
\Big[ \Gamma(\ngC \otimes R,\,\cF) \oplus
\Gamma(\ngC^\#\complot \cO_K,\,\cF) \to \Gamma(\ngC^\#\complot K,\,\cF) \Big].
\end{equation*}
The differential is as in the definition of $\RvGamma$ for shtukas.%
\end{dfn}

\breakflow
To make the expressions
in the rest of the section more legible
we will generally omit the argument $\ngS\times X$ of
the functors $\Gamma$, $\RvGamma$ and $\RGamma$ for quasi-coherent sheaves and
shtukas. The same applies to the associated complex functor $\cGamma$.

\begin{thm}\label{fcechlema}%
The natural map
$\RvGamma(\cF) \to \RGamma(\cF)$
is a quasi-\hspace{0pt}isomorphism
for every quasi-coherent sheaf $\cF$.
\end{thm}

\pf By construction $\cC(\cF)$ sits in a distinguished triangle
\begin{equation*}
\cC(\cF) \to (\iota_\ast \iota^\ast\cF \oplus f_\ast f^\ast \cF)[0] \to g_\ast
g^\ast\cF[0] \to [1].
\end{equation*}
Applying $\Gamma$ and $\RGamma$ we obtain a
morphism of distinguished triangles
\begin{equation*}
\xymatrix{
\Gamma(\cC(\cF)) \ar[r] \ar[d] &
\Gamma(\iota_\ast \iota^\ast \cF \oplus f_\ast f^\ast \cF)[0]
\ar[r] \ar[d] &
\Gamma(g_\ast g^\ast \cF)[0] \ar[r] \ar[d] & [1] \\
\RGamma(\cC(\cF)) \ar[r] &
\RGamma(\iota_\ast \iota^\ast \cF \oplus f_\ast f^\ast \cF)
\ar[r] &
\RGamma(g_\ast g^\ast \cF) \ar[r] & [1]
}
\end{equation*}
We will prove that the second and third vertical arrows in this diagram
are quasi-isomorphisms. It follows that the first arrow is a
quasi-isomorphism and so the lemma is proven.

Consider the third vertical arrow. The map $g$ is affine so that
$g_\ast g^\ast \cF [0] = \uR g_\ast g^\ast \cF$. Hence $\RGamma(g_\ast g^\ast
\cF) = \RGamma(\uR g_\ast g^\ast \cF) = \RGamma(\ngC^\#\complot K, \cF)$.
As $\Spec(\ngC^\#\complot K)$ is affine the natural map $\Gamma(\ngC^\#\complot K,
\cF)[0] \to \RGamma(\ngC^\#\complot K,\cF)$ is a quasi-isomorphism.
Hence the third vertical  map in the diagram above is a quasi-isomorphism.  The
maps $\iota$ and $f$ are also affine whence the same argument shows that the
second vertical map is a quasi-isomorphism. \quod

\breakflow
Let $\cM$ be a quasi-coherent shtuka on $\ngS\times X$. Define a complex of
shtukas on $\ngS\times X$:
\begin{equation*}
\cC(\cM) =
\Big[\iota_\ast\iota^\ast\cM \oplus f_\ast f^\ast \cM \to
g_\ast g^\ast \cM\Big].
\end{equation*}
Here the differential is the difference of the natural maps as in the definition
of $\RvGamma(\ngS\times X,\,\cM)$.
The sum of the adjunction units gives a natural morphism $\cM[0] \to \cC(\cM)$.
%

%
%
%
%
%
%
%

\begin{lem} If $\cM$ is a quasi-coherent shtuka on $\ngS\times X$ then
$\RvGamma(\cM) = \cGamma(\cC(\cM))$.\end{lem}

\pf Let $f\colon \cN \to \cN'$ be a morphism of shtukas and let
$C = [\cN \to \cN']$ be its mapping fiber. The associated complex functor
$\cGamma$ is defined in such a way that $\cGamma(C)$ is the mapping fiber
of $\cGamma(f)$. Applying this observation to $C = \cC(\cM)$ we get the result.\quod


\begin{lem}\label{cechlemb} If $\cM$ is a quasi-coherent shtuka on $\ngS\times X$
then the natural map $\cM[0] \to \cC(\cM)$ is a quasi-isomorphism.
\end{lem}

\pf Follows instantly from Lemma \ref{fcechlemb}. \quod

\begin{dfn}\label{cechmor}%
We define a natural map
$\RvGamma(\cM) \to \RGamma(\cM)$
as the composition
\begin{equation*}
\cGamma(\cC(\cM)) \to
\RGamma(\cC(\cM)) \xleftarrow{\,\isosign\,}
\RGamma(\cM)
\end{equation*}
of the natural map $\cGamma \to \RGamma$ and the quasi-isomorphism provided by
Lemma \ref{cechlemb}.\end{dfn}

\begin{thm}\label{cechcoh}%
%
For every quasi-coherent shtuka $\cM$ on $\ngS\times X$ the natural map
$\RvGamma(\cM) \to \RGamma(\cM)$ is a quasi-\hspace{0pt}isomorphism.\end{thm}

\pf According to Proposition~\ref{assocmor} the natural diagram
\begin{equation*}
\xymatrix{
\cGamma(\cC(\cM)) \ar[r] \ar[d] &
\Gamma(\cC(\cM_0)) \ar[r]^{i-j} \ar[d] &
\Gamma(\cC(\cM_1)) \ar[r] \ar[d] & [1] \\
\RGamma(\cC(\cM)) \ar[r] &
\RGamma(\cC(\cM_0)) \ar[r]^{i-j} &
\RGamma(\cC(\cM_1)) \ar[r] & [1]
}\end{equation*}
is a morphism of distinguished triangles. Furthermore the canonical
triangles are natural. Hence
the quasi-\hspace{0pt}isomorphism $\cM[0] \to \cC(\cM)$
induces an isomorphism of distinguished triangles
\begin{equation*}
\xymatrix{
\RGamma(\cC(\cM)) \ar[r] &
\RGamma(\cC(\cM_0)) \ar[r]^{i-j} &
\RGamma(\cC(\cM_1)) \ar[r] & [1] \\
\RGamma(\cM) \ar[r] \ar[u]^{\ltviso} &
\RGamma(\cM_0) \ar[r]^{i-j} \ar[u]^{\ltviso} &
\RGamma(\cM_1) \ar[r] \ar[u]_{\rtviso} & [1]
}\end{equation*}
Taking the composition of the two diagrams above we get a morphism of distinguished
triangles
\begin{equation*}
\xymatrix{
\RvGamma(\cM) \ar[r] \ar[d] &
\RvGamma(\cM_0) \ar[r]^{i-j} \ar[d] &
\RvGamma(\cM_1) \ar[r] \ar[d] & [1] \\
\RGamma(\cM) \ar[r] &
\RGamma(\cM_0) \ar[r]^{i-j} &
\RGamma(\cM_1) \ar[r] & [1]
}\end{equation*}
The second and third vertical arrows are quasi-\hspace{0pt}isomorphisms
by Theorem \ref{fcechlema}. Whence the result.\quod

\section{Compactly supported cohomology}
\label{sec:cscoh}

We continue using the notation and the conventions of Section
\ref{sec:cechcoh}.
%
In this section we assume that the coefficient algebra $\ngC$ carries a structure
of a locally compact $\Fq$-algebra. A typical example of $\ngC$ relevant to our
applications is the discrete algebra $\Fq[t]$ and the locally compact algebra
$\Fq(\!(t^{-1})\!)$. We denote $Y = \Spec R$ as in the introduction.

\begin{dfn}[\aftereq]\label{defcsc}\index{idx}{shtuka cohomology!compactly supported cohomology}\index{nidx}{shtuka cohomology!RGammac@$\RGammac$, compactly supported cohomology}%
Let $\cM$ be a quasi-coherent shtuka on $\ngS\times Y$.
The \emph{compactly supported cohomology complex} of $\cM$ is the $\ngC$-module
complex
\begin{equation*}
\RGammac(\ngS\times Y,\,\cM) =
\Big[\RGamma(\ngC \otimes R,\,\cM) \to \RGamma(\ngC^\#\complot K,\,\cM) \Big].
\end{equation*}
Here the differential is induced by the natural inclusion
$\ngC\otimes R \to \ngC^\#\complot K$.
The $n$-th cohomology group of $\RGammac(\ngS\times Y,\,\cM)$ is denoted
$\uH^n_c(\ngS\times Y,\,\cM)$.\end{dfn}

\begin{prp}\label{locfreecsc}%
If $\cM$ is a locally free shtuka on $\ngS\times Y$ then the natural map
\begin{equation*}
\RGammac(\ngS\times Y,\,\cM) \to
\RGamma\!\Big(\ngS\times Y,\,\tfrac{\cM(\ngC^\#\complot K)}{\cM(\ngC\otimes R)}\Big)[-1]
\end{equation*}
is a quasi-\hspace{0pt}isomorphism.\quod%
\end{prp}

\begin{dfn}\label{csctogerm}%
We define a map $\RGammac(\ngS\times Y,\,\cM) \to \RGammag(\ngC\indcot K,\,\cM)$
by the diagram
\begin{equation*}
\xymatrix{
\big[ \RGamma(\ngC\otimes R,\,\cM) \ar[r] \ar[d] & \RGamma(\ngC^\#\complot K,\,\cM) \ar@{=}[d] \big] \\
\big[ \RGamma(\ngC\indcot K,\,\cM) \ar[r] & \RGamma(\ngC^\#\complot K,\,\cM) \big]
}
\end{equation*}\end{dfn}

\begin{dfn}\label{csctocoh}%
Let $\cM$ be a quasi-coherent shtuka on $\ngS\times X$.
We define a map
$\RGammac(\ngS\times Y,\,\cM) \to \RGamma(\ngS\times X,\,\cM)$
as the composition
\begin{equation*}
\RGammac(\ngS\times Y,\,\cM) \xrightarrow{\textup{embedding}}
\RvGamma(\ngS\times X,\,\cM) \xrightarrow{\textup{Thm. }\ref{cechcoh}}
\RGamma(\ngS\times X,\,\cM).
\end{equation*}\end{dfn}
%

\begin{prp}\label{bcohmap}%
Let $\cM$ be a locally free shtuka on $\ngS\times X$.
If $\cM(\ngC \otimes \cO_K/\fm_K)$ is nilpotent then
$\RGammac(\ngS\times Y,\,\cM) \to \RGamma(\ngS\times X,\,\cM)$
is a quasi-\hspace{0pt}isomorphism.\end{prp}

\breakflow
The condition that $\cM(\ngC\otimes \cO_K/\fm_K)$ is nilpotent may be interpreted
as saying that $\cM$ is an extension by zero of a shtuka on the open $\tau$-subscheme
$\ngS\times Y \subset \ngS\times X$.
 
\breakflow\noindent
\textit{Proof of Proposition \ref{bcohmap}.} %
The map in question
extends to a distinguished triangle
\begin{equation*}
\RGammac(\ngS\times Y,\,\cM) \to \RvGamma(\ngS\times X,\,\cM) \to
\RGamma(\ngC^\#\complot\cO_K,\,\cM) \to [1].
\end{equation*}
The result follows since
$\RGamma(\ngC^\#\complot\cO_K,\,\cM) = 0$ by Proposition \ref{loccs}. \quod

\begin{dfn}\label{defgcmp}\index{idx}{shtuka cohomology!global germ map}%
Let $\cM$ be a locally free shtuka on $\ngS\times X$ such that $\cM(\ngC \otimes \cO_K/\fm_K)$
is nilpotent. The \emph{global germ map} is defined as the composition
\begin{equation*}
\RGamma(\ngS\times X,\,\cM) \xleftarrow{\,\isosign\,} \RGammac(\ngS\times Y,\,\cM)
\xrightarrow{\textup{Def. }\ref{csctogerm}} \RGammag(\ngC \indcot K,\,\cM)
\end{equation*}
in $\uD(\ngC)$ 
where the first arrow is the quasi-isomorphism of Definition~\ref{csctocoh}.
The adjective ``global'' indicates that this map involves a shtuka on the
whole $\ngS\times X$ as opposed to $\ngC \indcot \cO_K$.\end{dfn}

\section{Local-global compatibility}
\label{sec:locglobcmp}

We keep the conventions and the notation of Section~\ref{sec:cechcoh}.
Our coefficient algebra will be either a local field $\ngF$
or its ring of integers $\ngOF$.
We denote $\nsOF = \Spec\ngOF$ and $\nsF = \Spec\ngF$.
The letter ``$D$'' stands for a disk.

Let $\cM$ be a locally free shtuka on $\nsFX$ such that $\cM(\ngF\otimes \cO_K/\fm_K)$ is nilpotent.
We have two maps from the cohomology of $\cM$ to the germ cohomology
$\RGammag(\ngF\indcot K,\,\cM)$:
\begin{itemize}
\item the local germ map
of Definition \ref{defloccmp},

\item the global germ map
of Definition \ref{defgcmp}.
\end{itemize}
They form a square in the derived category of $\ngF$-vector spaces:
\begin{equation}\label{locglobcmpsq}
\vcenter{\vbox{\xymatrix{
\RGamma(\nsFX,\,\cM) \ar[d] \ar[rrr]^{\textrm{global}} &&&
\RGammag(\ngF\indcot K,\,\cM) \\
\RGamma(\ngF\indcot \cO_K,\,\cM) \ar[rrr]^{\textrm{local}}_{\bisosign} &&&
\RGammag(\ngF\indcot K,\,\cM) \ar@{=}[u]
}}}
\end{equation}
The left arrow in this square is the pullback map.

The definitions of the
local and the global germ map have nothing in common so there is no a priori
reason for \eqref{locglobcmpsq} to be commutative. Nevertheless we will prove
that \eqref{locglobcmpsq} commutes under the assumption that $\cM$
extends to a locally free shtuka on $\nsOFX$.

\begin{thm}\label{locglobcmp}If $\cM$ is a locally free shtuka on
$\nsOFX$ such that $\cM(\ngF\otimes \cO_K/\fm_K)$ is nilpotent then the
square \eqref{locglobcmpsq} is commutative.%
\end{thm}

\breakflow
Later in this chapter we will show that the left arrow in \eqref{locglobcmpsq}
is a quasi-isomorphism provided $\cM(\cO_F/\fm_F \otimes R)$ is nilpotent
(Theorem \ref{intconc}). The local germ map is a quasi-isomorphism by
construction. The commutativity of \eqref{locglobcmpsq} then implies that the
global germ map is a quasi-isomorphism, a property which is not evident
from its definition.

\afterall\noindent%
\textit{Proof of Theorem \ref{locglobcmp}.}
Take $\uH^1$ of \eqref{locglobcmpsq} and extend it to the left as follows:
\begin{equation}\label{extcmp}
\vcenter{\vbox{\xymatrix{
\uH^1(\nsOFX,\,\cM) \ar@{-->}[r] \ar@{-->}[d]
& \uH^1(\nsFX,\,\cM) \ar[r] \ar[d] &
\uH^1_g(\ngF\indcot K,\,\cM) \\
\uH^1(\ngOF \indcot \cO_K,\,\cM) \ar@{-->}[r]
& \uH^1(\ngF\indcot \cO_K,\,\cM) \ar[r] &
\uH^1_g(\ngF\indcot K,\,\cM) \ar@{=}[u]
}}}
\end{equation}
The three additional maps are the pullback morphisms. We proceed to prove that the
outer rectangle of \eqref{extcmp} commutes. 

Theorem \ref{cechcoh} equips us with natural isomorphisms
\begin{align*}
\uH^1(\nsOFX,\,\cM) &\cong \check{\uH}^1(\nsOFX,\,\cM), \\
\uH^1(\nsFX,\,\cM) &\cong \check{\uH}^1(\nsFX,\,\cM),
\end{align*}
while Proposition \ref{locfreegermcoh} provides a natural isomorphism
\begin{equation*}
\uH^1_g(\ngF\indcot K,\,\cM) \cong
\uH^0(\ngF\indcot K,\,\cQ)
\end{equation*}
where
\begin{equation*}
\cQ =\frac{\cM(\ngF^\#\complot K)}{\cM(\ngF \indcot K)}.
\end{equation*}
Using them we rewrite \eqref{extcmp} as
\begin{equation}\label{extcmpcech}
\vcenter{\vbox{\xymatrix{
\check{\uH}^1(\nsOFX,\,\cM) \ar[r] \ar[d]
& \check{\uH}^1(\nsFX,\,\cM) \ar[r] &
\uH^0(\ngF\indcot K,\,\cQ) \\
\uH^1(\ngOF \indcot \cO_K,\,\cM) \ar[r]
& \uH^1(\ngF \indcot \cO_K,\,\cM) \ar[r] &
\uH^0(\ngF \indcot K,\,\cQ) \ar@{=}[u].
}}}
\end{equation}
The middle arrow is omitted since it is not easy to describe in terms of \v{C}ech
cohomology.

Let the shtuka $\cM$ be given by a diagram
\begin{equation*}
\cM_0 \shtuka{\,\,i\,\,}{\,\,j\,\,} \cM_1.
\end{equation*}
%
Let $\ngC$ be either $\ngOF$ or $\ngF$.
By definition $\RvGamma(\Spec\ngC\times X,\,\cM)$ is
the total complex of the double complex
\begin{equation*}
\xymatrix{
\cM_1(\ngC \otimes R) \oplus \cM_1(\ngC^\#\complot\cO_K)
\ar[rr]^{\quad\quad\quad\textup{difference}}
&& \cM_1(\ngC^\#\complot K) \\
\cM_0(\ngC \otimes R) \oplus \cM_0(\ngC^\#\complot\cO_K) \ar[u]^{i-j}
\ar[rr]^{\quad\quad\quad\textup{difference}}
&& \cM_0(\ngC^\#\complot K) \ar[u]_{j-i}
}
\end{equation*}
So a cohomology class in $\check{\uH}^1(\Spec\ngC\times X,\,\cM)$ is represented by
a triple
\begin{equation*}
(a,x,b) \in \cM_1(\ngC \otimes R) \oplus \cM_1(\ngC^\#\complot\cO_K) \oplus
\cM_0(\ngC^\#\complot K)
\end{equation*}
satisfying $a - x + (j-i)(b) = 0$.

Fix a cohomology class $h \in \check{\uH}^1(\nsOFX,\,\cM)$. We want to
compute its image under the composition
\begin{equation*}
\check{\uH}^1(\nsOFX,\,\cM) \to
\check{\uH}^1(\nsFX,\,\cM) \to
\uH^0(\ngF\indcot K,\,\cQ)
\end{equation*}
of the two top arrows in \eqref{extcmpcech}. Let $(a,x,b)$ be a triple
representing $h$.
The image of $h$ in $\check{\uH}^1(\nsFX,\,\cM)$ is represented by the same
triple $(a,x,b)$.
From Definition \ref{defgcmp} it follows that the map
$\check{\uH}^1(\nsFX,\,\cM) \to \uH^1_g(\ngF\indcot K,\,\cM)$
of \eqref{extcmpcech} is a composition
\begin{equation}\label{h1gcmp}
\check{\uH}^1(\nsFX,\,\cM)
\xleftarrow{\,\isosign\,}
\uH^1_c(\ngF\otimes R,\,\cM) \to
\uH^0(\ngF\indcot K,\,\cQ)
\end{equation}
By construction $\RGammac(\ngF\otimes R,\,\cM)$ is the total complex of the double
complex
\begin{equation*}
\xymatrix{
\cM_1(\ngF \otimes R) \ar[r] & \cM_1(\ngF^\#\complot K) \\
\cM_0(\ngF \otimes R) \ar[u]^{i-j} \ar[r] & \cM_0(\ngF^\#\complot K) \ar[u]_{j-i}
}
\end{equation*}
So a cohomology class in $\uH^1_c(\ngF\otimes R,\,\cM)$ is represented by a pair
\begin{equation*}
(a',b') \in \cM_1(\ngF\otimes R) \oplus \cM_0(\ngF^\#\complot K)
\end{equation*}
such that $a' + (j-i)(b') = 0$.

The isomorphism $\uH^1_c(\ngF\otimes R,\,\cM) \cong \check{\uH}^1(\nsFX,\,\cM)$
of \eqref{h1gcmp} sends $(a',b')$ to $(a',0,b')$.
Thus in order to compute the image of $h$ in
$\uH^1_c(\ngF\otimes R,\,\cM)$ we need to replace $(a,x,b)$ with a
cohomologous triple of the form $(a',0,b')$.
A triple is a coboundary if and only if it has the form
\begin{equation*}
\big((i-j)(a'), (i-j)(y), a' - y\big)
\end{equation*}
where $a' \in \cM_0(\ngF\otimes R)$ and $y \in \cM_0(\ngF^\#\complot\cO_K)$.
By Proposition \ref{loccmpdesc} (1) there is a unique $y \in \cM_0(\ngF^\#\complot\cO_K)$ such
that $(i-j)(y) = x$.
The triple $(0,x,-y)$ is then a coboundary so that
$(a,x,b)$ is cohomologous to $(a, 0, b + y)$ and the image
of $h$ in $\uH^1_c(\ngF\otimes R,\,\cM)$ is represented by $(a, b + y)$.
 
We are finally ready to compute the image of
$h \in \check{\uH}^1(\nsOFX,\,\cM)$
under the composition \eqref{h1gcmp}:
\begin{equation*}
\check{\uH}^1(\nsFX,\,\cM)
\xleftarrow{\,\isosign\,}
\uH^1_c(\ngF\otimes R,\,\cM) \to
\uH^0(\ngF \indcot K,\,\cQ)
\end{equation*}
By definition of $\cQ$ we have
\begin{equation*}
\RGamma(\ngF \indcot K,\,\cQ) =
\Big[
\tfrac{\cM_0(\ngF^\# \complot K)}{\cM_0(\ngF\indcot K)} \xrightarrow{\,\,i - j\,\,}
\tfrac{\cM_1(\ngF^\#\complot K)}{\cM_1(\ngF\indcot K)}
\Big].
\end{equation*}
The second arrow in \eqref{h1gcmp} sends a pair $(a',b')$ representing a class
in $\uH^1_c(\ngF \otimes R,\,\cM)$ to the equivalence class $[b']$ of $b'$ in the quotient
$\cM_0(\ngF^\#\complot K)/\cM_0(\ngF\indcot K)$.
Above we demonstrated that
the image of $h$ in $\uH^1_c(\ngF \otimes R,\,\cM)$ is represented by the pair $(a,
b + y)$. Hence the image of $h$ in $\uH^0(\ngF\indcot K,\,\cQ)$ is given by the
equivalence class $[b + y]$.

The key observation in this proof is that $[b + y] = [y]$. Indeed
the left arrow in the natural commutative square
\begin{equation*}
\xymatrix{
\ngOF \indcot K \ar[d]\ar[r] & \ngF \indcot K \ar[d] \\
\ngOF^\#\complot K \ar[r] & \ngF^\# \complot K.
}
\end{equation*}
is an isomorphism by Proposition
\lref{ringot}{ringindcotdiscrcomplot}. Hence the homomorphism
$\ngOF^\# \complot K \to \ngF^\# \complot K$ factors through $\ngF \indcot K$ and the
natural map
\begin{equation*}
\cM_0(\ngOF^\#\complot K) \to 
\frac{\cM_0(\ngF^\#\complot K)}{\cM_0(\ngF\indcot K)}
\end{equation*}
is zero.
As $b \in \cM_0(\ngOF^\#\complot K)$ by construction we conclude that $[b + y] = [y]$.

So far we have demonstrated the following. Let
$h \in \check{\uH}^1(\nsOFX,\,\cM)$ be a cohomology class.
If $h$ is represented by a triple
\begin{equation*}
(a,x,b) \in \cM_1(\ngOF \otimes R) \oplus \cM_1(\ngOF^\#\complot\cO_K) \oplus
\cM_0(\ngOF^\#\complot K)
\end{equation*}
then the image of $h$ under the composition
\begin{equation*}
\check{\uH}^1(\nsOFX,\,\cM) \to
\check{\uH}^1(\nsFX,\,\cM) \to
\uH^0(\ngF\indcot K,\,\cQ)
\end{equation*}
of the two top arrows in the square \eqref{extcmpcech} is given by the
equivalence class
\begin{equation*}
[y] \in \frac{\cM_0(\ngF^\#\complot K)}{\cM_0(\ngF\indcot K)}
\end{equation*}
where $y \in \cM_0(\ngF^\#\complot\cO_K)$
is the unique element satisfying $(i-j)(y) = x$.
We are now in position to prove that the square \eqref{extcmpcech} is
commutative.

The cohomology classes in $\uH^1(\ngOF\indcot\cO_K,\,\cM)$ are represented by
elements of $\cM_1(\ngOF\indcot\cO_K)$.
By Proposition \lref{ringot}{ringindcotdiscrcomplot}
the natural map
$\ngOF \indcot \cO_K \to \ngOF^\# \complot \cO_K$
is an isomorphism.
Hence we can identify $\cM_1(\ngOF\indcot\cO_K)$ with
$\cM_1(\ngOF^\#\complot\cO_K)$. The left arrow
$\check{\uH}^1(\nsOFX,\,\cM) \to \uH^1(\ngOF\indcot\cO_K,\,\cM)$
of the square \eqref{extcmpcech} sends the cocycle $(a,x,b)$ to $x \in
\cM_1(\ngOF^\#\complot\cO_K) = \cM_1(\ngOF \indcot \cO_K)$. Now Proposition
\ref{loccmpdesc} (2) implies that the image of $x$ under the composition of the
two bottom arrows in \eqref{extcmpcech} is $[y]$. Therefore
the square \eqref{extcmpcech} is commutative.

We deduce that the outer rectangle of \eqref{extcmp} is commutative.
Since the $\ngF$-linear extension
\begin{equation*}
\ngF \otimes_{\ngOF} \uH^1(\nsOFX,\,\cM) \to \uH^1(\nsFX,\,\cM)
\end{equation*}
of the top horizontal map in this square is an isomorphism the right square of
\eqref{extcmp} is commutative too. By Proposition \ref{loccoh} the complex
$\RGammag(\ngF\indcot K,\,\cM)$ is concentrated in degree $1$. Therefore
commutativity of the right square of \eqref{extcmp} implies commutativity of the
main diagram \eqref{locglobcmpsq} in the derived category of $\ngF$-vector spaces.\quod

\section{\texorpdfstring{Completed \v{C}ech cohomology}{Completed Cech cohomology}}
\label{sec:complcech}

We keep the notation and the conventions of Sections \ref{sec:cechcoh} and \ref{sec:locglobcmp}.
%
In this section
we present a refined version of the \v{C}ech
method for computing the cohomology of \emph{coherent} shtukas on $\nsOFX$.
In essense it is the \v{C}ech method of Section~\ref{sec:cechcoh} developed in
the setting of formal schemes over $\nsOF = \Spec\ngOF$.
%

\begin{dfn}
Let $\cM$ be a coherent shtuka on $\nsOFX$. The
\emph{completed \v{C}ech chomology complex} of $\cM$ is
\begin{equation*}
\RcGamma(\nsOFX,\,\cM) =
\Big[ \RGamma(\ngOF\complot R,\,\cM) \oplus \RGamma(\ngOF\complot\cO_K,\,\cM) \to
\RGamma(\ngOF\complot K,\,\cM) \Big].
\end{equation*}
Here the differential is the difference of the natural maps.%
\end{dfn}

\breakflow
Recall that the \v{C}ech complex of $\cM$ is
\begin{equation*}
\RvGamma(\nsOFX,\,\cM) = 
\Big[ \RGamma(\ngOF\otimes R,\,\cM) \oplus \RGamma(\ngOF^\#\complot\cO_K,\,\cM)\to
\RGamma(\ngOF^\#\complot K,\,\cM) \Big].
\end{equation*}
We thus have a natural map
$\RvGamma(\nsOFX,\,\cM) \to \RcGamma(\nsOFX,\,\cM)$.
\begin{thm}\label{complcechqi}%
Let $\cM$ be a shtuka on $\nsOFX$. If $\cM$ is coherent then
the natural map
$\RvGamma(\nsOFX,\,\cM) \to \RcGamma(\nsOFX,\,\cM)$
is a quasi-\hspace{0pt}isomorphism.%
\end{thm}

\breakflow
Rather than using completed \v{C}ech cohomology directly we will
rely on 
Theorem~\ref{intconc}
which captures
a ``cohomology concentration'' phenomenon for locally free shtukas on $\nsOFX$.
It will play a role in Chapter~\ref{chapter:trace}.


\begin{thm}\label{intconc}%
Let $\cM$ be a locally free shtuka on $\nsOFX$.
If $\cM(\cO_F/\fm_F \otimes R)$ is nilpotent then the natural
map $\RGamma(\nsOFX,\,\cM) \to \RGamma(\cO_F \complot \cO_K,\,\cM)$ is a
quasi-\hspace{0pt}isomorphism.\end{thm}


\pf Due to Theorem \ref{complcechqi} it is enough to prove that the natural map
\begin{equation*}
\RcGamma(\nsOFX,\,\cM) \to \RGamma(\cO_F \complot\cO_K,\,\cM)
\end{equation*}
is a quasi-isomorphism. By definition
\begin{equation*}
\RcGamma(\nsOFX,\,\cM) =
\Big[ \RGamma(\cO_F\complot R,\,\cM) \oplus \RGamma(\cO_F\complot\cO_K,\,\cM) \to
\RGamma(\cO_F\complot K,\,\cM) \Big].
\end{equation*}
Hence it is enough to show that the complexes $\RGamma(\cO_F\complot R,\cM)$ and
$\RGamma(\cO_F\complot K,\cM)$ are acyclic.

According to Proposition \lref{otalgprops}{discrcomplotnoether} the ring
$\cO_F\complot R$ is noetherian and complete with respect to the ideal $\fm_F \complot
R$. By Proposition \lref{otalgprops}{discrcomplotquot} the natural map
$\cO_F\complot R \to \cO_F/\fm_F \otimes R$ is surjective with kernel $\fm
\complot R$. 
Thus $\RGamma(\cO_F\complot R,\cM) = 0$ by Proposition \lref{nilp}{nilpcomp}.
Applying the same argument to $\cO_F\complot K^\#$ we deduce that
$\RGamma(\cO_F\complot K^\#,\cM) = 0$. 
The natural map $\cO_F\complot K^\#\to \cO_F\complot K$ is an isomorphism by
Proposition \lref{ringot}{ringdiscrcomplot} whence $\RGamma(\cO_F\complot K,\cM)
= 0$.\quod

\breakflow
We now turn to the proof of Theorem~\ref{complcechqi}. We will derive it from
a similar statement for coherent sheaves.

Let $\cF$ be a quasi-coherent sheaf on $\nsOFX$.
Recall that
\begin{equation*}
\RvGamma(\nsOFX,\,\cF) =
\Big[ \Gamma(\cO_F \otimes R,\,\cF) \oplus \Gamma(\cO_F^\#\complot \cO_K,\,\cF)\to
\Gamma(\cO_F^\#\complot K,\,\cF) \Big].
\end{equation*}
We set
\begin{equation*}
\RcGamma(\nsOFX,\,\cF) =
\Big[\Gamma(\cO_F\complot R,\,\cF) \oplus \Gamma(\cO_F\complot\cO_K,\,\cF) \to
\Gamma(\cO_F\complot K,\,\cF) \Big]
\end{equation*}
with the same differentials as in the definition of $\RcGamma$ for shtukas.
To improve the legibility we will generally omit the argument $\nsOFX$ of
the functors $\RvGamma$ and $\RcGamma$.
By construction we have a natural map
$\RvGamma(\cF) \to\RcGamma(\cF)$.


For technical reasons it will be more convenient for us to work with
different presentations of the complexes
$\RvGamma(\cF)$ and $\RcGamma(\cF)$.
We define the complexes
\begin{align*}
B(\cF) &=
\Big[ \Gamma(\cO_F \otimes R, \cF) \oplus \Gamma(\cO_F\indcot \cO_K, \cF)\to
\Gamma(\cO_F\indcot K, \cF) \Big], \\
\widehat{B}(\cF) &=
\Big[\Gamma(\cO_F\complot R, \cF) \oplus \Gamma(\cO_F\complot\cO_K^\#,\cF) \to
\Gamma(\cO_F\complot K^\#,\cF) \Big]
\end{align*}
with the same differentials as
$\RvGamma(\cF)$ and $\RcGamma(\cF)$.

\begin{lem}\label{ccechbiso}The natural map $B(\cF) \to \RvGamma(\cF)$ is an
isomorphism.\end{lem}

\pf Follows from Proposition \lref{otalgprops}{ringindcotdiscrcomplot} since $\cO_F$ is
compact. \quod

\begin{lem}\label{ccechbcompliso}The natural map
$\widehat{B}(\cF) \to \RcGamma(\cF)$ is an isomorphism.\end{lem}

\pf Follows from Proposition \lref{otalgprops}{ringdiscrcomplot} since $\cO_F$ is compact.
\quod

\begin{lem}\label{fcechfull} For every quasi-coherent sheaf $\cF$ on
$\nsOFX$ there exists a natural quasi-isomorphism
$B(\cF) \xrightarrow{\isosign} \RGamma(\nsOFX,\,\cF)$.\end{lem}

\pf The natural map $B(\cF) \to \RvGamma(\cF)$ is an isomorphism by Lemma
\ref{ccechbiso}.
So the result is a consequence of Theorem \ref{fcechlema}.\quod

\breakflow
\begin{lem}\label{ccechfiniso}Let $\cF$ be a quasi-coherent sheaf on $\nsOFX$.
If $\fm_F^n \cF = 0$ for some $n \gg 0$ then the natural map
$B(\cF) \to \widehat{B}(\cF)$ is an isomorphism.\end{lem}

\pf 
Consider the natural diagram
\begin{equation*}
\Gamma(\cO_F \indcot K, \cF) \to
\Gamma(\cO_F \complot K^\#, \cF) \to
\Gamma(\cO_F/\fm^n_F \otimes K, \cF).
\end{equation*}
By Proposition \lref{otalgprops}{discrcomplotquot} the second arrow in this
diagram is the reduction modulo $\fm^n_F$. The composite arrow is the reduction
modulo $\fm^n_F$ by Propostion \lref{otalgprops}{indcotquot}.
Both arrows are isomorphisms since $\fm_F^n \cF = 0$. Hence so is the first
arrow. The same argument shows that the natural maps
$\Gamma(\cO_F\indcot\cO_K,\cF) \to \Gamma(\cO_F\complot\cO_K^\#,\cF)$ and
$\Gamma(\cO_F\otimes R,\cF) \to \Gamma(\cO_F\complot R,\cF)$ are
isomorphisms.\quod

%
%
%

\begin{lem}\label{ccechlimiso}If $\cF$ is a coherent sheaf on $\nsOFX$
then the natural map $\widehat{B}(\cF) \to \lim_n \widehat{B}(\cF/\fm_F^n)$ is
an isomorphism. \end{lem}

\pf
Proposition
\lref{otalgprops}{discrcomplotnoether} shows that
the ring $\cO_F\complot K^\#$ is noetherian and complete with respect to the
ideal $\fm_F\complot K^\#$. According to Proposition
\lref{otalgprops}{discrcomplotquot} this ideal is generated by $\fm_F$.
The $\cO_F\complot K^\#$-module
$\Gamma(\cO_F\complot K^\#, \cF)$ is finitely generated. As a consequence it is
complete with respect to $\fm_F (\cO_F\complot K^\#)$. Hence the natural map
\begin{equation*}
\Gamma(\cO_F\complot K^\#,\cF) \to
\lim\nolimits_n \Gamma(\cO_F\complot K^\#,\cF/\fm^n_F) =
\lim\nolimits_n \Gamma(\cO_F\complot K,\cF)/\fm^n_F
\end{equation*}
is an isomorphism. 
The same argument applies to $\cO_F\complot R$ and
$\cO_F\complot\cO_K^\#$.\quod

\breakflow
The next two lemmas use the derived limit functor $\Rlim$ 
for abelian groups. We use the Stacks Project [\stacks{07KV}] as a
reference for $\Rlim$.

\begin{lem}\label{ccechrlimqi}If $\cF$ is a quasi-coherent sheaf on $\nsOFX$
then the natural map
\begin{equation*}
\lim\nolimits_n \widehat{B}(\cF/\fm_F^n) \to \Rlim_n\widehat{B}(\cF/\fm^n)
\end{equation*}
is a quasi-isomorphism.\end{lem}

\pf Let us denote
\begin{align*}
A_n &= \Gamma(\cO_F\complot R,\cF/\fm_F^n),\\
B_n &= \Gamma(\cO_F\complot\cO_K^\#,\cF/\fm_F^n),\\
C_n &= \Gamma(\cO_F\complot K^\#,\cF/\fm_F^n).
\end{align*}
The natural map in question extends to a morphism of distinguished triangles
\begin{equation*}
\xymatrix{
\lim_n \widehat{B}(\cF/\fm_F^n) \ar[r] \ar[d] &
\lim_n A_n \oplus B_n \ar[r] \ar[d] &
\lim_n C_n \ar[r] \ar[d] &
[1] \\
\Rlim_n \widehat{B}(\cF/\fm_F^n) \ar[r] &
\Rlim_n A_n \oplus B_n \ar[r] &
\Rlim_n C_n \ar[r] &
[1]
}
\end{equation*}
So in order to show that the first vertical arrow is a quasi-isomorphism it is
enough to prove that so are the second and the third vertical arrows.

The transition maps in the projective system
$\{B_n\}_{n\geqslant 1}$ are surjective by construction. Hence this system
satisfies the Mittag-Leffler condition [\stacks{02N0}]. As a consequence $\uR^1\!\lim_n B_n
= 0$ [\stacks{07KW}]. The natural map $\lim_n B_n \to \Rlim_n
B_n$ is thus a quasi-isomorphism. The same argument applies to $\{A_n\}$ and
$\{C_n\}$. \quod

\begin{lem}\label{ccechhomlim} If $\cF$ is a coherent sheaf on $\nsOFX$
then the natural map
\begin{equation*}
\uH^i(\widehat{B}(\cF)) \to \lim\nolimits_n \uH^i(\widehat{B}(\cF/\fm_F^n))
\end{equation*}
is an isomorphism for every $i$. \end{lem}

\pf Lemma \ref{ccechlimiso} implies that the map
\begin{equation*}
\uH^i(\widehat{B}(\cF)) \to \uH^i(\lim\nolimits_n \widehat{B}(\cF/\fm_F^n))
\end{equation*}
is an isomorphism for every $i$. At the same time the map
\begin{equation*}
\uH^i(\lim\nolimits_n \widehat{B}(\cF/\fm_F^n)) \to \uH^i(\Rlim_n\widehat{B}(\cF/\fm_F^n))
\end{equation*}
is an isomorphism by Lemma \ref{ccechrlimqi}.
The cohomology group $\uH^i(\Rlim_n \widehat{B}(\cF/\fm_F^n))$ sits in a
natural short exact sequence [\stacks{07KY}]
\begin{multline*}
0\! \to\! \uR^1\!\lim\nolimits_n \uH^{i-1}(\widehat{B}(\cF/\fm_F^n))\! \to
\!\uH^i(\Rlim_n \widehat{B}(\cF/\fm_F^n))\! \to
\!\lim\nolimits_n \uH^i(\widehat{B}(\cF/\fm_F^n))\! \to\! 0.
\end{multline*}
We thus need to prove that the first term in this sequence vanishes. 

Lemma \ref{fcechfull} provides us with natural isomorphisms
\begin{equation*}
\uH^{i-1}(\nsOFX,\,\cF/\fm_F^n) \xrightarrow{\isosign} \uH^{i-1}(\widehat{B}(\cF/\fm_F^n)).
\end{equation*}
As $\nsOFX$ is proper over $\cO_F$ it follows that 
$\uH^{i-1}(\widehat{B}(\cF/\fm_F^n))$ is a finitely generated
$\cO_F/\fm_F^n$-module for every $n$. Thus the image of $\uH^{i-1}(\widehat{B}(\cF/\fm_F^m))$ in
$\uH^{i-1}(\widehat{B}(\cF/\fm_F^n))$ is independent of $m$ for $m \gg n$.
In other words the projective system
\begin{equation*}
\{\uH^{i-1}(\widehat{B}(\cF/\fm_F^n))\}_{n\geqslant 1}
\end{equation*}
satisfies the
Mittag-Leffler condition [\stacks{02N0}]. Hence its first derived
limit $\uR^1\!\lim$ is zero [\stacks{07KW}]. \quod

\begin{lem}\label{fccomplqi} If $\cF$ is a coherent sheaf on $\nsOFX$
then the natural map $\RvGamma(\cF) \to \RcGamma(\cF)$ is a quasi-isomorphism.\end{lem}

\pf In view of Lemmas \ref{ccechbiso} and \ref{ccechbcompliso} it is enough to
prove that the natural map $B(\cF) \to \widehat{B}(\cF)$ is a quasi-isomorphism.
Let $i \in \bZ$. We have a natural commutative diagram
\begin{equation*}
\xymatrix{
\uH^i(B(\cF)) \ar[r] \ar[d] & \uH^i(\widehat{B}(\cF)) \ar[d] \\
\lim_n \uH^i(B(\cF/\fm_F^n)) \ar[r] &
\lim_n \uH^i(\widehat{B}(\cF/\fm_F^n)).
}
\end{equation*}
The right arrow is an isomorphism by Lemma \ref{ccechhomlim} while
the bottom arrow is an isomorphism by Lemma \ref{ccechfiniso}.
Thus in order to prove that the top arrow is an isomorphism it is enough to show
that the left arrow is so.
This arrow fits into a natural commutative square
\begin{equation*}
\xymatrix{
\uH^i(\nsOFX, \cF) \ar[r] \ar[d]&
\uH^i(B(\cF)) \ar[d] \\
\lim_n \uH^i(\nsOFX,\cF/\fm_F^n) \ar[r] &
\lim_n \uH^i(B(\cF/\fm_F^n))
}
\end{equation*}
where the horizontal arrows are the natural isomorphisms of Lemma
\ref{fcechfull}. According to the Theorem on formal functions [\stacks{02OC}] 
the left arrow in this square is an isomorphism.
Whence the result follows.
\quod

\afterall\noindent\textit{Proof of Theorem \ref{complcechqi}. }%
Let the shtuka $\cM$ be given by a diagram
\begin{equation*}
\cM_0 \shtuka{i}{j} \cM_1.
\end{equation*}
We have a natural morphism of distinguished triangles
\begin{equation*}
\xymatrix{
\RvGamma(\cM) \ar[r] \ar[d] &
\RvGamma(\cM_0) \ar[r]^{i-j} \ar[d] &
\RvGamma(\cM_1) \ar[r] \ar[d] & [1] \\
\RcGamma(\cM) \ar[r] &
\RcGamma(\cM_0) \ar[r]^{i-j} &
\RcGamma(\cM_1) \ar[r] & [1]
}
\end{equation*}
The second and third vertical arrows are quasi-isomorphisms by Lemma
\ref{fccomplqi} so we are done. \quod

\section{Change of coefficients}
\label{sec:globcoeffch}

Fix a noetherian $\Fq$-algebra $\ngC$. In this section we study how the cohomology of
shtukas on $\Spec\ngC \times X$ changes under the pullback to $\Spec\ngCalt \times X$
where $\ngCalt$ is an $\ngC$-algebra. 
%
%
We denote $\ngS = \Spec\ngC$ and $\ngSalt = \Spec\ngCalt$.

\begin{dfn}\label{shtpullcohadj}%
Let $\cM$ be an $\cO_{\ngS\times X}$-module shtuka.
We define a natural morphism
\begin{equation*}
\RGamma(\ngS\times X,\,\cM) \otimes_\ngC^{\mathbf L} \ngCalt \to \RGamma(\ngSalt \times X,\,\cM)
\end{equation*}
by extension of scalars of the pullback morphism
$\RGamma(\ngS\times X,\,\cM) \to \RGamma(\ngSalt\times X,\,\cM)$.
Given an $\cO_{\ngS\times X}$-module $\cE$
we define a natural morphism
$\RGamma(\ngS\times X, \cE) \otimes_\ngC^{\mathbf L} \ngCalt \to \RGamma(\ngSalt \times X, \cE)$
in the same way.%
\end{dfn}

\begin{lem}\label{coeffchangetrimor}%
Let $\ngCalt$ be an $\ngC$-algebra.
If $\cM$ is an $\cO_{\ngS\times X}$-module shtuka given by a diagram
\begin{equation*}
\cM_0 \shtuka{i}{j} \cM_1
\end{equation*}
then the natural diagram
\begin{equation}\label{coeffchangetri}
\vcenter{\vbox{\xymatrix{
\RGamma(\ngS\times X, \cM) \otimes_\ngC^{\mathbf L} \ngCalt \ar[r] \ar[d] &
\RGamma(\ngSalt \times X, \cM) \ar[d] \\
\RGamma(\ngS\times X,\cM_0) \otimes_\ngC^{\mathbf L} \ngCalt \ar[r] \ar[d]_{i-j} &
\RGamma(\ngSalt \times X,\cM_0) \ar[d]^{i-j}  \\
\RGamma(\ngS\times X,\cM_1) \otimes_\ngC^{\mathbf L} \ngCalt \ar[r] \ar[d] &
\RGamma(\ngSalt \times X,\cM_1) \ar[d] \\
[1] & [1]
}}}
\end{equation}
is a morphism of distinguished triangles.
Here the left column 
is the image under $-\otimes_\ngC^{\mathbf L} \ngCalt$ of the
canonical triangle for $\cM$  and
the right column 
is the canonical triangle 
for the pullback of $\cM$ to $\ngSalt \times X$.%
\end{lem}

\pf Follows from Proposition \ref{shtpullmor}.\quod


\begin{prp}\label{coeffchangetriiso}
Let $\cM$ be an $\cO_{\ngS\times X}$-module shtuka given by a diagram
\begin{equation*}
\cM_0 \shtuka{i}{j} \cM_1.
\end{equation*}
If $\cM_0$, $\cM_1$ are coherent and flat over $\ngC$ then the following holds:
\begin{enumerate}
\item 
$\RGamma(\ngS\times X, \cM)$ is a perfect $\ngC$-module complex.

\item For every $\ngC$-algebra $\ngCalt$ the natural map
\begin{equation*}
\RGamma(\ngS\times X,\cM) \otimes_\ngC^{\mathbf L} \ngCalt \to \RGamma(\ngSalt\times X,\cM)
\end{equation*}
is a quasi-\hspace{0pt}isomorphism. Moreover the diagram
\eqref{coeffchangetri} is an isomorphism of distinguished triangles.
\end{enumerate}
\end{prp}

\pf Since $\cM_0$ is coherent and flat over $\ngC$
the base change theorem for coherent cohomology [\stacks{07VK}] shows that
$\RGamma(\ngS\times X,\,\cM_0)$ is a perfect $\ngC$-module complex and the natural map
$\RGamma(\ngS\times X,\,\cM_0) \otimes_\ngC^{\mathbf L} \ngCalt \to \RGamma(\ngSalt \times X,\,\cM_0)$
is a quasi-\hspace{0pt}isomorphism. The same applies to $\cM_1$.
As $\RGamma(\ngS\times X,\cM)$ fits to a distinguished triangle
\begin{equation*}
\RGamma(\ngS\times X,\cM) \to
\RGamma(\ngS\times X,\cM_0) \xrightarrow{\,\,i-j\,\,}
\RGamma(\ngS\times X,\cM_1) \to [1]
\end{equation*}
we conclude that it is a perfect $\ngC$-module complex.
Finally Lemma \ref{coeffchangetrimor} implies that \eqref{coeffchangetri} is an
isomorphism of distinguished triangles. \quod

\section{\texorpdfstring{$\zeta$-isomorphisms}{Zeta-isomorphisms}}
\label{sec:globzeta}

Let $\ngC$ be a noetherian $\Fq$-algebra and let $\ngS = \Spec\ngC$.
In this section we study $\zeta$-isomorphisms for shtukas over $\ngS\times X$.
We will prove that under suitable conditions they are stable under change of $\ngC$.

Let $\cM$ be an $\cO_{\ngS\times X}$-module shtuka given by a diagram
\begin{equation*}
\cM_0 \shtuka{i}{j} \cM_1.
\end{equation*}
Assume that
the $\ngC$-module complexes
$\RGamma(\cM)$, $\RGamma(\Der\cM)$, $\RGamma(\cM_0)$ and $\RGamma(\cM_1)$
are bounded with perfect cohomology modules.
In this situation we have
a $\zeta$-isomorphism
\begin{equation*}
\zeta_\cM\colon
\det\nolimits_\ngC \RGamma(\cM) \xrightarrow{\isosign}
\det\nolimits_\ngC \RGamma(\Der\cM).
\end{equation*}
It
is the composition
\begin{equation*}
\det\nolimits_\ngC^{\phantom{1}} \RGamma(\cM) \xrightarrow{\isosign}
\det\nolimits_\ngC^{\phantom{1}} \RGamma(\cM_0) \otimes_\ngC^{\phantom{1}}
\det\nolimits_\ngC^{-1} \RGamma(\cM_1) \xrightarrow{\isosign}
\det\nolimits_\ngC^{\phantom{1}} \RGamma(\Der\cM)
\end{equation*}
of isomorphisms induced by the canonical triangles
\begin{align*}
&\RGamma(\cM) \to \RGamma(\cM_0) \xrightarrow{\,\,i-j\,\,}
\RGamma(\cM_1) \to [1], \\
&\RGamma(\Der\cM) \to \RGamma(\cM_0) \xrightarrow{\,\,i\,\,}
\RGamma(\cM_1) \to [1].
\end{align*}

\begin{prp}\label{regzetadefined}%
If $\ngC$ is regular then the $\zeta$-isomorphism
is defined for every coherent shtuka on $\ngS\times X$.\end{prp}

\breakflow
We will also need $\zeta$-isomorphisms for coefficient rings $\ngC$ which are
not regular. The example of such an $\ngC$ relevant to our study is a local
artinian ring which is not a field.

\afterall\noindent
\textit{Proof of Proposition \ref{regzetadefined}.}
Suppose that a coherent shtuka $\cM$ is given by a diagram
\begin{equation*}
\cM_0 \shtuka{i}{j} \cM_1.
\end{equation*}
Grothendieck vanishing theorem [\stacks{02UZ}] shows that the $\ngC$-modules
$\uH^n(\cM_0)$ and $\uH^n(\cM_1)$ are zero for $n > 1$.
By [\stacks{02O5}] they are finitely
generated for $n = 0, 1$.
Thus $\uH^n(\cM)$ and $\uH^n(\Der\cM)$ are zero for $n \not\in\{0,1,2\}$ and
finitely generated for $n = 0,1,2$.
The ring $\ngC$ has finite global
dimension since it is regular [\stacks{00O7}]. Whence $\uH^n(\cM)$,
$\uH^n(\Der\cM)$, $\uH^n(\cM_0)$ and $\uH^n(\cM_1)$ are perfect $\ngC$-modules.\quod

\begin{prp}\label{globcoeffchzeta}%
Let $\ngCalt$ be an $\ngC$-algebra and let $\ngSalt = \Spec\ngCalt$.
Let $\cM$ be a shtuka on $\ngS\times X$ given by a diagram
\begin{equation*}
\cM_0 \shtuka{i}{j} \cM_1.
\end{equation*}
Let $\cM_{\ngSalt}$ be the pullback of $\cM$ to $\ngSalt \times X$.
Assume that
\begin{enumerate}
\item the $\zeta$-isomorphisms are defined for $\cM$ and $\cM_{\ngSalt}$,

\item $\cM_0$ and $\cM_1$ are coherent and flat over $\ngC$.
\end{enumerate}
Then the following holds:
\begin{enumerate}
\item The natural maps
\begin{align*}
&\RGamma(\cM)\otimes_\ngC^{\mathbf L} \ngCalt \to \RGamma(\cM_{\ngSalt}), \\
&\RGamma(\Der\cM)\otimes_\ngC^{\mathbf L} \ngCalt \to \RGamma(\Der\cM_{\ngSalt})
\end{align*}
of Definition \ref{shtpullcohadj} are quasi-\hspace{0pt}isomorphisms.

\item The natural square
\begin{equation*}
\xymatrix{
\det\nolimits_\ngC \RGamma(\cM) \otimes_\ngC \ngCalt
\ar[rr]^{\zeta_\cM \otimes 1} \ar[d]_{\ltviso} &&
\det\nolimits_\ngC \RGamma(\Der\cM)\otimes_\ngC \ngCalt \ar[d]^{\rtviso} \\
\det\nolimits_\ngCalt \RGamma(\cM_{\ngSalt}) \ar[rr]^{{\zeta_\cM}_{\ngSalt}} &&
\det\nolimits_\ngCalt \RGamma(\Der\cM_{\ngSalt})
}
\end{equation*}
is commutative. Here the vertical arrows are induced by the quasi-\hspace{0pt}isomorphisms
of (1).\end{enumerate}
\end{prp}


\pf The natural isomorphisms of determinants
\begin{align*}
&\det\nolimits_\ngC^{\phantom{1}} \RGamma(\cM) \to
\det\nolimits_\ngC^{\phantom{1}} \RGamma(\cM_0) \otimes_\ngC^{\phantom{1}}
\det\nolimits_\ngC^{-1} \RGamma(\cM_1),\\
&\det\nolimits_\ngC^{\phantom{1}} \RGamma(\Der\cM) \to
\det\nolimits_\ngC^{\phantom{1}} \RGamma(\cM_0) \otimes_\ngC^{\phantom{1}}
\det\nolimits_\ngC^{-1} \RGamma(\cM_1)
\end{align*}
induced by the triangles of $\cM$ and $\Der\cM$ are stable under
the pullback to $\ngCalt$ by construction (see the proof of Corollary 2 after Theorem 2
in \cite{kmdet}). So the result follows from Proposition \ref{coeffchangetriiso}. \quod

\chapter{Regulator theory}
\label{chapter:reg}
\label{ch:locloccoh}
\label{ch:artnilp}
\label{ch:artreg}
\label{ch:reg}

Let $F$ be a local field and let $K$ be a finite product of local fields.
As usual we assume $F$ and $K$ to contain $\Fq$.
We denote $\cO_F \subset F$ and $\cO_K \subset K$ the rings of
integers, $\fmF$ the maximal ideal of $\cO_F$ and $\fm_K$ the Jacobson radical
of $\cO_K$. We omit the subscript $F$ for the ideal $\fm \subset \cO_F$ to
improve the legibility.

We mainly work with $\cO_F\complot\cO_K$-module shtukas. In
agreement with the conventions of Section \lref{ottaustruct}{sec:ottaustruct} the
$\tau$-structure on $\cO_F\complot\cO_K$ is given by the endomorphism which acts
as the identity on $\cO_F$ and as the $q$-Frobenius on $\cO_K$.

The aim of this chapter is to construct for a certain class of shtukas $\cM$ on
$\cO_F\complot\cO_K$ a natural quasi-isomorphism
\begin{equation*}
\RGamma(\cM) \xrightarrow{\quad\rho\quad}
\RGamma(\Der\cM)
\end{equation*}
called the \emph{regulator}.
We do it for \emph{elliptic shtukas}, a class of shtukas
generalizing the 
models of Drinfeld modules in the sense of Chapter \ref{chapter:locmod}.
The key definitions and results of this chapter are as follows:
\begin{itemize}
\item Theorem \ref{llcoh} 
describes the cohomology of a certain class of shtukas on $\cO_F\complot\cO_K$.
It applies in particular to elliptic shtukas.

\item Definition \ref{defell} introduces elliptic shtukas.

%
\item Definition \ref{defellreg} introduces the regulator for elliptic shtukas.

\item The existence and unicity of the regulator is affirmed by Theorem
\ref{ellregexist}.
%
%
\end{itemize}

It is easy to characterize the regulator as a 
natural transformation of functors on the category of elliptic shtukas
(see Definition \ref{defellreg}).
However
its construction
is a bit involved. 

The content of this chapter is new save for the preliminary Sections
\ref{sec:elltopprelim} and \ref{sec:ellalgprelim}.
Lemma~\ref{ofmodfgcrit} should have certainly appeared before
but we are not aware of a reference for it.
Our search for a shtuka-theoretic regulator
was motivated by the article \cite{valeurs} of V.~Lafforgue.

\afterall\noindent\textbf{Remark.} %
After the work on this text was finished the author found out
that V.~Lafforgue has constructed a map \cite[Lemme 4.8]{valeurs} which
resembles in some respects
the regulator isomorphism of this chapter.
It is an interesting question whether the two maps are actually the same.


\section{Topological preliminaries}
\label{sec:elltopprelim}

In this section we give a topological criterion for an $\cO_F$-module to be
finitely generated. We will
use it to prove that cohomology modules of certain shtukas are finitely
generated. 

\begin{lem}\label{ofmodfgcrit} Let $M$ be a compact Hausdorff $\cO_F$-module.
The following are equivalent:
\begin{enumerate}
\item $M/\fmF$ is finite as a set.

\item $M$ is finitely generated as an $\cO_F$-module without topology.
\end{enumerate}\end{lem}

\pf (1) $\Rightarrow$ (2).
Let $z$ be a uniformizer of $\cO_F$. The submodule $\fmF M \subset M$ is the
image of $M$ under multiplication by $z$ so it is closed, compact and Hausdorff.
Furthermore multiplication by $z$ defines a surjective map
$M/\fmF \to (\fmF M)/\fmF$ so that $(\fmF M)/\fmF$ is finite. By induction we
conclude that the submodules $\fmF^n M \subset M$ are closed and of finite index,
hence open.

Let us show that the open submodules $\fmF^n M$ form a fundamental system of
neighbourhoods of zero. 
Proposition \lref{pont}{lclth} and Lemma \lref{compideal}{compideal} imply that $M$ admits a
fundamental system of open $\cO_F$-submodules. If $U \subset M$ is an
open submodule then $M/U$ is finite so there exists an $n > 0$ such
that $z^n$ acts by zero on $M/U$. Hence $\fmF^n M \subset U$.

Now $M$ is a compact $\cO_F$-module so it is complete as a topological
$\Fq$-vector space. As the submodules $\fmF^n M$ form a fundamental system
we conclude that $M = \lim_{n > 0} M/\fmF^n$.

Let $r$ be the dimension of $M/\fmF$ as an $\cO_F/\fmF$-vector space.
For every $n > 0$ let $H_n$ be the set of surjective $\cO_F$-linear maps from
$\cO_F^{\oplus r}$ to $M/\fmF^n$.
The sets $H_n$ form a projective system $H_*$
in a natural way. Every point of the limit of $H_*$ defines a
continuous morphism from $\cO_F^{\oplus r}$ to $M$. Such a morphism
is surjective since it has dense image by construction and its domain
$\cO_F^{\oplus r}$ is compact.

By definition the set $H_1$ is nonempty.
Nakayama's lemma implies that 
the transition maps $H_{n+1} \to H_n$ are surjective for all $n > 0$.
Therefore the projective system $H_*$ has a nonempty limit.\quod

\breakflow
It is worth mentioning that Lemma \ref{ofmodfgcrit} works for any
nonarchimedean local field $F$ and more generally for any local noetherian
ring $\cO_F$ with finite residue field. 
Indeed one can show that a (locally) compact Hausdorff $\cO_F$-module $M$ admits
a fundamental system of open submodules and the rest of the argument applies
essentially as is.

\section{Algebraic preliminaries}
\label{sec:ellalgprelim}

Let us review some elementary algebraic properties of the ring
$\cO_F\complot\cO_K$.

\begin{lem}\label{openideal}%
%
An ideal $\ibase \subset \cO_K$ is open if and only if it is a free $\cO_K$-module of
rank $1$.%
\end{lem}

\pf By definition $\cO_K$ is a finite product of complete discrete valuation
rings. The ideal $\ibase$ is open if and only if it projects to a nonzero ideal in
every factor. Whence the result. \quod

\breakflow
Let $\ibase \subset \cO_K$ be an open ideal. We will often use
natural homomorphisms
\begin{align*}
f_\ibase&\colon \cO_F\complot\cO_K \to \cO_F\otimes\cO_K/\ibase, \\
g_\ibase&\colon \cO_F\complot\cO_K \to F \otimes \cO_K/\ibase.
\end{align*}
The homomorphism $f_\ibase$ is
the completion of the natural map $\cO_F\otimes_{\textup{c}} \cO_K \to
\cO_F\otimes_{\textup{c}} \cO_K/\ibase$. We use the fact that $\cO_K/\ibase$ is
finite to identify $\cO_F\complot(\cO_K/\ibase)$ with $\cO_F\otimes\cO_K/\ibase$.
%
The homomorphism $g_\ibase$ is
the composition of $f_\ibase$ with the natural map
$\cO_F\otimes\cO_K/\ibase \to F \otimes\cO_K/\ibase$. By construction $g_\ibase$
factors over $F\otimes_{\cO_F}(\cO_F\complot\cO_K)$. 
The notation $f_\ibase$ and $g_\ibase$ will not be used.
The same constructions apply to a nonzero ideal $\icoef \subset \cO_F$.

By definition
$\cO_F = k[[z]]$
for a finite field extension $k$ of $\bF_q$ and a uniformizer $z \in \cO_F$.
In a similar way
\begin{equation*}
\cO_K = \prod_{i=1}^d k_i [[\zeta_i]]
\end{equation*}
where $k_i$ are finite field extensions of $\Fq$ and $\zeta_i$ are uniformizers
of the factors of $\cO_K$. As a consequence
\begin{equation*}
\cO_F\complot\cO_K =
\prod_{i=1}^d (k \otimes k_i)[[z,\zeta_i]]
\end{equation*}
is a finite product of power series rings in two variables. With this
observation in mind the following lemmas become obvious.

\begin{lem}\label{basecocart} Let $\ibase \subset\cO_K$ be an open ideal.
\begin{enumerate}
\item The ideal $\ibase\cdot(\cO_F\complot\cO_K) = \cO_F\complot \ibase$ is a free
$\cO_F\complot\cO_K$-module of rank $1$.

\item The sequence $0 \to \ibase\cdot(\cO_F\complot\cO_K) \to \cO_F\complot\cO_K \to
\cO_F\otimes\cO_K/\ibase \to 0$ is exact.\quod\end{enumerate}\end{lem}


\begin{lem}\label{coeffcocart} Let $\icoef \subset \cO_F$ be an open ideal.
\begin{enumerate}
\item The ideal $\icoef \cdot (\cO_F\complot\cO_K) = \icoef \complot\cO_K$ is a free module of rank $1$.

\item The sequence $0 \to \icoef\cdot(\cO_F\complot\cO_K) \to
\cO_F\complot\cO_K\to \cO_F/\icoef \otimes \cO_K \to 0$ is exact.\quod
\end{enumerate}\end{lem}

%
%

\breakflow
The next lemma will only be used in the proof of Proposition \ref{elldelta}.

\begin{lem}\label{ellintersect}
If $\ibase \subset \cO_K$ and $\icoef \subset \cO_F$ are open ideals then the natural
sequence
\begin{equation*}
0 \to
\frac{\cO_F\complot \ibase}{\icoef \complot \ibase} \oplus \frac{\icoef\complot\cO_K}{\icoef \complot \ibase}
\to \frac{\cO_F\complot\cO_K}{\icoef \complot \ibase} \to
\frac{\cO_F\complot\cO_K}{\icoef \complot \cO_K + \cO_F\complot \ibase} \to 0
\end{equation*}
is exact.\end{lem}

\pf Follows since $\icoef \complot \ibase = \icoef\complot\cO_K \cap \cO_F\complot \ibase$. \quod

\section{Cohomology with artinian coefficients}
\label{sec:ellnilplemma}

Fix a finite $\Fq$-algebra $\ngC$ which is a local artinian ring.
Let $\fm \subset \ngC$ be the maximal ideal. 
In this section we work with the ring $\ngC \otimes \cO_K$.
We equip it with the $\tau$-ring
structure given by the endomorphism
which acts as the identity on $\ngC$ and as the
$q$-Frobenius on $\cO_K$. 


\begin{lem}\label{artmodnilp}%
Let $\cM$ be a locally free shtuka on $\ngC\otimes \cO_K$.
If $\cM(\ngC/\fm \otimes K)$ is nilpotent then
$\cM(\ngC\otimes K)$ is nilpotent.\end{lem}

\pf The ring $\ngC \otimes K$ is noetherian and complete with respect to the
ideal $\fm \otimes K$. As the ideal $\fm$ is nilpotent the result follows from
Proposition \lref{nilp}{nilpcomp}.\quod 

\begin{thm}\label{artnilpcoh}%
Let $\cM$ be a locally free shtuka on  $\ngC\otimes \cO_K$.
If $\cM(\ngC/\fm \otimes K)$ is nilpotent then:
\begin{enumerate}
\item $\uH^0(\cM) = 0$,

\item $\uH^1(\cM)$ is a free $\ngC$-module of finite rank.
\end{enumerate}\end{thm}

\pf (1) Since $\cM$ is locally free the natural
map $\uH^0(\ngC\otimes\cO_K,\,\cM) \to \uH^0(\ngC \otimes K,\,\cM)$ injective.
However $\cM(\ngC \otimes K)$ is nilpotent by Lemma \ref{artmodnilp}
so $\uH^0(\ngC \otimes K,\cM) = 0$ by Propostion \lref{nilp}{nilpvanish}.

(2) First let us prove that $\uH^1(\cM)$ is a flat $\ngC$-module.
Suppose that $\cM$ is given by a diagram
\begin{equation*}
M_0 \shtuka{i}{j} M_1.
\end{equation*}
The cohomology complex $\RGamma(\cM)$ is represented by the associated complex
\begin{equation*}
\cGamma(\cM) = [ M_0 \xrightarrow{\,\,i-j\,\,} M_1 ].
\end{equation*}
It is a comlex of flat $\ngC$-modules since $\cM$ is locally free by assumption.
As a consequence
\begin{equation*}
\cGamma(\cM) \otimes^{\mathbf L}_\ngC \ngC/\fm = \cGamma(\cM) \otimes_\ngC \ngC/\fm.
\end{equation*}
However $\cGamma(\cM) \otimes_\ngC \ngC/\fm$ is the complex representing
$\RGamma(\ngC/\fm \otimes \cO_K,\,\cM)$. Applying the argument (1) above to the
shtuka $\cM(\ngC/\fm \otimes \cO_K)$ we deduce that
$\RGamma(\ngC/\fm \otimes\cO_K,\,\cM)$ is concentrated in degree $1$.
Hence $\uH^1(\cM) \otimes^{\mathbf L}_\ngC \ngC/\fm$ is concentrated in degree $0$
or in other words $\Tor_n(\uH^1(\cM),\,\ngC/\fm) = 0$ for $n > 0$. Therefore
$\uH^1(\cM)$ is a flat $\ngC$-module [\stacks{051K}].

Next we prove that $\uH^1(\cM)$ is finitely generated.
The $\cO_K$-modules $M_0$ and $M_1$ are finitely generated by assumption. 
They carry a natural compact Hausdorff topology given by the powers of the ideal $\fm_K$.
We would like to prove that the map $(i-j)\colon M_0 \to M_1$ is open.
Since $M_1$ is a compact $\cO_K$-module it then follows that $M_1/(i-j)M_0 =
\uH^1(\cM)$ is a finite set.

Consider the locally compact $K$-vector spaces $V_0 = M_0 \otimes_{\cO_K} K$ and
$V_1 = M_1 \otimes_{\cO_K} K$. By Lemma \ref{artmodnilp} the shtuka $\cM(\ngC \otimes K)$
is nilpotent whence $i\colon V_0 \to V_1$ is an isomorphism and the
endomorphism $i^{-1}j$ of $V_0$ is nilpotent. The isomorphism $i^{-1}\colon V_1 \to V_0$
is continuous by $K$-linearity. The map $j\colon V_0 \to V_1$ is continuous
since it is a Frobenius-linear morphism of finite-dimensional $K$-vector
spaces. As a consequence $i^{-1}j\colon V_0 \to V_0$ is continuous. Since it is
nilpotent we conclude that the endomorphism $(1-i^{-1} j)^{-1}$ is continuous.
Therefore $1 - i^{-1} j$ is open. However
\begin{equation*}
i - j = i \circ (1 - i^{-1} j).
\end{equation*}
So $(i-j)\colon V_0 \to V_1$ is a composition of open maps. We conclude that
$(i-j)\colon M_0 \to M_1$ is open.\quod

\section{Finiteness of cohomology}
\label{sec:llcoh}


In this section we prove that under certain natural conditions the cohomology
groups of shtukas over $\cO_F\complot\cO_K$ are finitely generated free
$\cO_F$-modules.

\begin{lem}\label{ztriangle}For every
locally free shtuka $\cM$ on $\cO_F\complot\cO_K$ the natural map
$\RGamma(\cM) \otimes_{\cO_F}^{\mathbf{L}} \cO_F/\fmF \to \RGamma(\cO_F/\fmF\otimes\cO_K,\,\cM)$
is a quasi-\hspace{0pt}isomorphism.\end{lem}

\pf Let $z \in \cO_F$ be a uniformizer. Observe that
$\RGamma(\cM)\otimes_{\cO_F}^{\mathbf{L}} \cO_F/\fmF$ is the cone of
the multiplication map $z\colon \RGamma(\cM)\to\RGamma(\cM)$.
Applying $\RGamma(-)$ to the short exact sequence
$0 \to \cM \xrightarrow{z} \cM \to \cM/z \to 0$
we conclude that the natural map $\RGamma(\cM)\otimes_{\cO_F}^{\mathbf{L}}\cO_F/\fmF \to \RGamma(\cM/z)$ 
is a quasi-\hspace{0pt}isomorphism.
By Lemma \ref{coeffcocart}
the quotient $\cM/z$ is the restriction of $\cM$ to $\cO_F/\fmF \otimes \cO_K$.
Whence the result.\quod

\begin{prp}\label{llcohfin}Let $\cM$ be a locally free shtuka on $\cO_F\complot\cO_K$.
The following are equivalent:
\begin{enumerate}
\item $\RGamma(\cO_F/\fmF\otimes\cO_K,\,\cM)$ is a prerfect complex
of $\cO_F/\fmF$-vector spaces.

\item $\RGamma(\cM)$ is a perfect complex of $\cO_F$-modules.\end{enumerate}\end{prp}

\pf (2) $\Rightarrow$ (1) follows from Lemma \ref{ztriangle}. (1) $\Rightarrow$ (2).
Suppose that $\cM$ is given by a diagram
\begin{equation*}
M_0 \shtuka{i}{j} M_1.
\end{equation*}
The $\cO_F\complot\cO_K$-modules
$M_0$, $M_1$ come equipped with a
canonical topology given by powers of the ideal
\begin{equation*}
\ibase = \fmF\complot\cO_K + \cO_F\complot\fm_K.
\end{equation*}
The map $i$ is continuous in this topology since it is $\cO_F\complot\cO_K$-linear.
The partial Frobenius $\tau\colon \cO_F\complot\cO_K \to \cO_F\complot\cO_K$
maps $\ibase$ to itself. As $j$ is $\tau$-linear it follows that $j$ is continuous.
Hence $i-j$ is continuous.

Consider the complex
\begin{equation*}
\big[M_0 \xrightarrow{\,\,i-j\,\,} M_1\big].
\end{equation*}
The ring $\cO_F\complot\cO_K$ is compact with respect to the $\ibase$-adic topology.
Therefore the finitely generated
$\cO_F\complot\cO_K$-modules $M_0$, $M_1$ are compact Hausdorff.
It follows that the image of $M_0$ in $M_1$ is closed and the quotient topology on
$\uH^1(\cM)$ is compact Hausdorff. So is the subspace topology on $\uH^0(\cM)$.

The natural map $\cO_F\to\cO_F\complot\cO_K$ is continuous and the differential
$i-j$ in the complex above is $\cO_F$-linear. As a consequence $\uH^0(\cM)$ and
$\uH^1(\cM)$ are topological $\cO_F$-modules. 
Now $\uH^0(\cM)$ and $\uH^1(\cM)$ are compact Hausdorff so
Lemma \ref{ofmodfgcrit} shows that they are finitely generated.\quod

\breakflow
The first main result of this chapter is the following:

\begin{thm}\label{llcoh} Let $\cM$ be a locally free shtuka on $\cO_F\complot\cO_K$. If
$\cM(\cO_F/\fmF\otimes K)$ is nilpotent then the following holds:
\begin{enumerate}
\item $\uH^0(\cM) = 0$,

\item $\uH^1(\cM)$ is a finitely generated free
$\cO_F$-module.\end{enumerate}\end{thm}

\pf The shtuka $\cM(\cO_F/\fmF \otimes K)$ is nilpotent
so Theorem \ref{artnilpcoh} shows that
\begin{enumerate}
\item[(i)] $\uH^0(\cO_F/\fmF\otimes \cO_K, \,\cM) = 0$,

\item[(ii)] $\uH^1(\cO_F/\fmF\otimes\cO_K, \,\cM)$ is finite as a set.\end{enumerate}
Using Proposition \ref{llcohfin} we deduce that $\uH^0(\cM)$
and $\uH^1(\cM)$ are finitely generated $\cO_F$-modules.
At the same time (i) in combination with Lemma \ref{ztriangle} implies
that $\uH^0(\cM)$ is divisible and
$\uH^1(\cM)$ is torsion-free. So the result follows.\quod

%

\section{Artinian regulators}
\label{sec:artreg}

Let us fix a finite $\Fq$-algebra $\ngC$ which is a local artinian ring.
We denote $\fm \subset \ngC$ the maximal ideal.
In this section we work over the ring $\ngC \otimes \cO_K$.
We equip it with the $\tau$-ring
structure given by the endomorphism
which acts as the identity on $\ngC$ and as the
$q$-Frobenius on $\cO_K$. 

We study locally free shtukas on $\ngC \otimes \cO_K$ which
restrict to nilpotent shtukas on $\ngC/\fm \otimes K$.
Under certain conditions we will define a regulator map for such shtukas,
the \emph{artinian regulator}.

\begin{lem} Let $\cM$ be a locally free shtuka on $\ngC \otimes \cO_K$ given by a diagram
\begin{equation*}
M_0 \shtuka{i}{j} M_1.
\end{equation*}
If $\cM(\ngC/\fm \otimes K)$ is nilpotent then the map $i\colon M_0 \otimes_{\cO_K}
K \to M_1 \otimes_{\cO_K} K$ is an isomorphism.\end{lem}

\pf By Lemma \ref{artmodnilp} the shtuka $\cM(\ngC \otimes K)$ is nilpotent. The
$i$-map of such a shtuka is an isomorphism by definition. \quod

\begin{dfn}\label{defartreg}%
Let a locally free shtuka $\cM$ on $\ngC
\otimes \cO_K$ be given by a diagram
\begin{equation*}
M_0 \shtuka{i}{j} M_1.
\end{equation*}
Suppose that $\cM(\ngC/\fm \otimes K)$ is nilpotent.
We say that the artinian regulator is defined for $\cM$ if the endomorphism
$i^{-1} j$ of the $\ngC$-module $M_0 \otimes_{\cO_K} K$ preserves the submodule $M_0$.
In this case we define 
the \emph{artinian regulator}
$\rho_\cM \colon \cGamma(\cM) \to \cGamma(\Der\cM)$
by the diagram
\begin{equation*}
\xymatrix{
[ M_0 \ar[r]^{i-j} \ar[d]_{1 - i^{-1} j} & M_1 \ar[d]^1 ] \\
[ M_0 \ar[r]^{i} & M_1 ].
}
\end{equation*}\end{dfn}

\breakflow
In a moment we will give a sufficient condition for the regulator to be
defined. Before that let us study its properties.

\begin{lem}\label{artregprop}%
The regulator of Definition \ref{defartreg} has the following properties.
\begin{enumerate}
\item $\rho_\cM$ is natural in $\cM$.

\item $\rho_\cM$ is an isomorphism.
\end{enumerate}\end{lem}

\pf (1) Clear. (2) Suppose that $\cM$ is given by a diagram
\begin{equation*}
M_0 \shtuka{i}{j} M_1.
\end{equation*}
The shtuka $\cM(\ngC \otimes K)$ is nilpotent by Lemma \ref{artmodnilp} whence
the endomorphism $i^{-1} j$ of $M_0 \otimes_{\cO_K} K$ is nilpotent. As a
consequence the endomorphism $1 - i^{-1} j$ is in fact an automorphism.
\quod

%
%

\begin{dfn} Let $\ibase \subset \cO_K$ be an ideal. Let a quasi-coherent shtuka
$\cM$ on $\ngC\otimes\cO_K$ be given by a diagram
\begin{equation*}
M_0 \shtuka{i}{j} M_1.
\end{equation*}
We define
\begin{equation*}
\ibase\cM =
\Big[ \ibase M_0 \shtuka{i}{j} \ibase M_1 \Big].
\end{equation*}\end{dfn}

\begin{rmk}
According to Lemma~\ref{openideal}
an ideal $\ibase \subset \cO_K$ is open if and only if
it is free of rank $1$.\end{rmk}

\begin{lem}\label{arttwist} Let $\cM$ be a locally free shtuka on
$\ngC\otimes\cO_K$ such that $\cM(\ngC/\fm \otimes K)$ is nilpotent.
Let $\ibase \subset \cO_K$ be an open ideal. 
\begin{enumerate}
\item $\ibase\cM$ is a locally free shtuka which restricts to a nilpotent shtuka on
$\ngC/\fm \otimes K$.

\item If the regulator is defined for $\cM$ then it is defined for $\ibase\cM$.
\end{enumerate}\end{lem}

\pf (1) Clear.
(2)
Suppose that $\cM$ is given by a diagram
\begin{equation*}
M_0 \shtuka{i}{j} M_1.
\end{equation*}
We know that $j(M_0) \subset i(M_0)$. As a consequence $j(\ibase M_0) \subset i(\ibase^q
M_0)$ which implies that $\ibase M_0$ is invariant under $i^{-1} j$. \quod

\begin{lem}\label{artcohles} Let $\ibase \subset \cO_K$ be an open ideal and let
$\cM$ be a locally free shtuka on $\ngC \otimes \cO_K$. If $\cM(\ngC/\fm\otimes K)$ is
nilpotent then the short exact sequence of shtukas
\begin{equation*}
0 \to \ibase \cM \to \cM \to \cM/\ibase \to 0
\end{equation*}
induces a long exact sequence of cohomology
\begin{equation}\label{basedeltaredseq}
0 \to \uH^0(\cM/\ibase) \xrightarrow{\,\,\delta\,\,}
\uH^1(\ibase\cM) \xrightarrow{\quad}
\uH^1(\cM) \xrightarrow{\textup{red.}}
\uH^1(\cM/\ibase) \to 0.
\end{equation}\end{lem}

\breakflow
The exact sequence \eqref{basedeltaredseq} will play an important role in our
theory. It will mainly appear as the sequence
\eqref{ellcoeffdotseq} for elliptic shtukas on $\cO_F\complot\cO_K$.
 
\afterall\noindent\textit{Proof of Lemma \ref{artcohles}.} 
By Lemma \ref{arttwist} $\ibase\cM$ is a locally free shtuka on
$\ngC\otimes\cO_K$ whose restriction to $\ngC/\fm\otimes K$ is nilpotent.
Hence $\uH^0(\ibase\cM)$ and $\uH^0(\cM)$ vanish by Theorem \ref{artnilpcoh}. \quod

\breakflow
If a shtuka $\cM$ on $\ngC \otimes \cO_K$ is linear then $\RGamma(\cM)$ is
represented by a complex with an $\cO_K$-linear differential. As a consequence
$\RGamma(\cM)$ carries a natural action of $\cO_K$.

\begin{lem}\label{artlierediso} 
Let $\cM$ be a locally free shtuka on $\ngC \otimes \cO_K$ such that $\cM(\ngC/\fm \otimes
K)$ is nilpotent. Let $\ibase \subset \cO_K$ be an open ideal.
If $\cM$ is linear then the following are equivalent:
\begin{enumerate}
\item $\ibase \cdot \uH^1(\cM) = 0$.

\item The map
$\uH^0(\cM/\ibase) \xrightarrow{\delta}
\uH^1(\ibase\cM)$ in \eqref{basedeltaredseq} is an isomorphism.

\item The map
$\uH^1(\cM) \xrightarrow{\textup{red.}}
\uH^1(\cM/\ibase)$
in \eqref{basedeltaredseq} is an
isomorphism.\end{enumerate}\end{lem}

\pf Suppose that $\cM$ is given by a diagram
\begin{equation*}
M_0 \shtuka{i}{0} M_1.
\end{equation*}
By definition $\uH^1(\cM) = M_1/i(M_0)$. Hence the following are equivalent:
\begin{enumerate}
\item $\ibase \cdot \uH^1(\cM) = 0$.

\item[(1$'$)] $\ibase M_1 \subset i(M_0)$.

\item[(1$''$)] The natural map $\uH^1(\ibase\cM) \to \uH^1(\cM)$ is zero.
\end{enumerate}
By Lemma \ref{artcohles} either of the conditions (2) or (3) is equivalent to
(1$''$). \quod

\begin{lem}\label{artrediso} 
Let $\cM$ be a locally free shtuka on $\ngC \otimes \cO_K$ such that $\cM(\ngC/\fm \otimes
K)$ is nilpotent. Suppose that the regulator is defined for $\cM$.
Let $\ibase \subset \cO_K$ be an open ideal.
The following are equivalent:
\begin{enumerate}
\item $\ibase \cdot \uH^1(\Der\cM) = 0$.

\item The map
$\uH^0(\cM/\ibase) \xrightarrow{\delta}
\uH^1(\ibase\cM)$ in \eqref{basedeltaredseq} is an isomorphism.

\item The map
$\uH^1(\cM) \xrightarrow{\textup{red.}}
\uH^1(\cM/\ibase)$
in \eqref{basedeltaredseq} is an
isomorphism.\end{enumerate}\end{lem}

\pf 
Lemma \ref{arttwist} tells that the regulator is defined for the shtuka
$\ibase\cM$ as well. By Lemma \ref{artregprop} regulators are natural.
We thus get a commutative diagram of complexes
\begin{equation*}
\xymatrix{
0 \ar[r] & \cGamma(\ibase\cM) \ar[r] \ar[d]_{\rho_{\ibase\cM}}
& \cGamma(\cM) \ar[r] \ar[d]_{\rho_{\cM}}
&\cGamma(\cM/\ibase) \ar[r] \ar[d]^{\rho_{\cM/\ibase}} &
0 \\
0 \ar[r] & \cGamma(\Der \ibase\cM) \ar[r]
& \cGamma(\Der\cM) \ar[r]
& \cGamma(\Der\cM/\ibase) \ar[r] & 0 
}
\end{equation*}
where $\rho_{\cM/\ibase}$ is induced by $\rho_{\cM}$. The regulators $\rho_{\cM}$ and
$\rho_{\ibase\cM}$ are isomorphisms by Lemma \ref{artregprop}.
As a consequence $\rho_{\cM/\ibase}$ is an isomorphism.
Taking $\uH^1$ of the diagram above we conclude that the following are equivalent:
\begin{enumerate}
\item[(3)] The reduction map $\uH^1(\cM) \to \uH^1(\cM/\ibase)$ is an isomorphism.

\item[(3$'$)] The reduction map $\uH^1(\Der\cM) \to \uH^1(\Der\cM/\ibase)$ is an isomorphism.
\end{enumerate}
According to Lemma \ref{artlierediso} the condition (1) is equivalent to (3$'$).
Hence the condition (1) is equivalent to (3).
The conditions (2) and (3) are equivalent by Lemma \ref{artcohles}. \quod


%
%

\begin{prp}\label{artregexist}
Let $\cM$ be a locally free shtuka on $\ngC \otimes \cO_K$ such that $\cM(\ngC/\fm
\otimes K)$ is nilpotent.
Let $\ibase \subset \cO_K$ be an open ideal.
Assume that
\begin{enumerate}
\item $\ibase \cdot \uH^1(\Der\cM) = 0$,

\item $\cM/\ibase$ is linear.
\end{enumerate}
Then the following holds:
\begin{enumerate}
\item The regulator is defined for $\cM$.

\item The map
$\uH^0(\cM/\ibase) \xrightarrow{\delta} \uH^1(\ibase\cM)$
in \eqref{basedeltaredseq} is an isomorphism.

\item The map
$\uH^1(\cM) \xrightarrow{\textup{red.}} \uH^1(\cM/\ibase)$
in \eqref{basedeltaredseq} is an isomorphism.
\end{enumerate}\end{prp}

\pf
Suppse that $\cM$ is given by a diagram
\begin{equation*}
M_0 \shtuka{i}{j} M_1.
\end{equation*}
Now $j(M_0) \subset \ibase M_1$ by assumption (2) and $\ibase M_1 \subset i(M_0)$ by
assumption (1). Hence $M_0$ is preserved by $i^{-1} j$. We conclude that the
regulator is defined for $\cM$. In view of this fact the results (2) and (3)
follow from the assumption (1) by Lemma \ref{artrediso}. \quod

\begin{prp}\label{artregred}
Let $\cM$ be a locally free shtuka on $\ngC \otimes \cO_K$ such that
$\cM(\ngC/\fm \otimes K)$ is nilpotent.
Let $\ibase \subset \cO_K$ be an open ideal. Assume that $\cM/\ibase$ is linear.
If the regulator is defined for $\cM$
then the diagram
\begin{equation*}
\xymatrix{
\uH^1(\cM) \ar[r]^{\textup{red.}} \ar[d]_{\rho_{\cM}} &
\uH^1(\cM/\ibase) \ar[d]^1 \\
\uH^1(\Der\cM) \ar[r]^{\textup{red.}} &
\uH^1(\Der\cM/\ibase).
}
\end{equation*}
is commutative.\end{prp}

\pf The regulator is defined for $\ibase\cM$ by Lemma \ref{arttwist}. Since the regulators
are natural by Lemma \ref{artregprop} we obtain a commutative diagram of complexes
\begin{equation*}
\xymatrix{
0 \ar[r] & \cGamma(\ibase\cM) \ar[r] \ar[d]_{\rho_{\ibase\cM}}
& \cGamma(\cM) \ar[r] \ar[d]_{\rho_{\cM}}
&\cGamma(\cM/\ibase) \ar[r] \ar[d]^{\rho_{\cM/\ibase}} &
0 \\
0 \ar[r] & \cGamma(\Der \ibase\cM) \ar[r]
& \cGamma(\Der\cM) \ar[r]
& \cGamma(\Der\cM/\ibase) \ar[r] & 0 
}
\end{equation*}
where $\rho_{\cM/\ibase}$ is the morphism induced by $\rho_{\cM}$. The regulator
$\rho_{\cM}$ is given by the identity map in degree $1$. So the same is true for
$\rho_{\cM/\ibase}$. Taking $\uH^1$ of the diagram above we get the result. \quod

%
%

\begin{prp}\label{artregdelta}
Let $\cM$ be a locally free shtuka on $\ngC \otimes \cO_K$ such that
$\cM(\ngC/\fm \otimes K)$ is nilpotent.
Let $\ibase \subset \cO_K$ be an open ideal. Assume that
\begin{enumerate}
\item $\ibase \cdot \uH^1(\Der\cM) = 0$,

\item $\cM/\ibase^2$ is linear.
\end{enumerate}
Then the regulator is defined for $\ibase\cM$ and the square
\begin{equation*}
\xymatrix{
\uH^0(\cM/\ibase) \ar[r]^{\delta} \ar[d]_{1} &
\uH^1(\ibase\cM) \ar[d]^{\rho_{\ibase\cM}} \\
\uH^0(\Der\cM/\ibase) \ar[r]^{\delta} &
\uH^1(\Der \ibase\cM)
}
\end{equation*}
is commutative. Here the maps $\delta$ are the boundary homomorphisms of the
long exact
sequence \eqref{basedeltaredseq}.
\end{prp}

\pf The regulator is defined for $\cM$ by Proposition \ref{artregexist}.
Hence it is defined for $\ibase\cM$ by Lemma \ref{arttwist}.
We then have a diagram of complexes
\begin{equation}\label{artdeltacommdiag}
\vcenter{\vbox{\xymatrix{
0 \ar[r] & \cGamma(\ibase\cM) \ar[r] \ar[d]_{\rho_{\ibase\cM}}
& \cGamma(\cM) \ar[r] \ar[d]_{\rho_{\cM}}
&\cGamma(\cM/\ibase) \ar[r] \ar[d]^1 &
0 \\
0 \ar[r] & \cGamma(\Der \ibase\cM) \ar[r]
& \cGamma(\Der\cM) \ar[r]
& \cGamma(\Der\cM/\ibase) \ar[r] & 0 
}}}
\end{equation}
with short exact rows.
We will
show that the right square commutes. The result then follows since
\eqref{artdeltacommdiag}
induces a morphism of long exact cohomology sequences.

The right square of \eqref{artdeltacommdiag}
commutes in degree $1$ since
$\rho_{\cM}$ is given by the identity map in degree $1$. We thus need to show
commutativity in degree $0$.
Suppose that $\cM$ is given by a diagram
\begin{equation*}
M_0 \shtuka{i}{j} M_1.
\end{equation*}
To show commutativity of \eqref{artdeltacommdiag} is is enough to prove that
$i^{-1} j(M_0) \subset \ibase M_0$. By assumption (2) we have
$j(M_0) \subset \ibase^2 M_1$. Assumption (1) implies that $\ibase^2 M_1 \subset i(\ibase
M_0)$. Hence $i^{-1} j(M_0) \subset \ibase M_0$. \quod


\section{Elliptic shtukas}
\label{sec:defell}

Starting from 
this section we work over the ring $\cO_F\complot\cO_K$.
Recall that an ideal of $\ibase \subset \cO_K$ is open if and only if
it is free of rank $1$ (Lemma \ref{openideal}).


%

\begin{dfn}\label{defell}\index{idx}{shtuka!elliptic}\index{idx}{ramification ideal!elliptic shtuka@of an elliptic shtuka}\index{nidx}{$\rami$, ramification ideal}%
Let $\cM$ be a shtuka on $\cO_F\complot\cO_K$ and let $\rami\subset\cO_K$ be an open ideal.
We say that $\cM$ is an \emph{elliptic shtuka of
ramification ideal $\rami$} if the following holds:
\begin{enumerate}
\setcounter{enumi}{-1}
\renewcommand{\theenumi}{E\arabic{enumi}}
\item\label{ell0}
$\cM$ is locally free.

\item\label{ell1}
$\cM(\cO_F/\fmF\otimes K)$ is nilpotent.

\item\label{ell2}
$\cM(F\otimes\cO_K/\fm_K)$ is nilpotent.

\item\label{ell3}
$\fmF \cdot \uH^1(\Der\cM) = \rami\cdot\uH^1(\Der\cM)$.

\item\label{ell4}
$\cM(\cO_F\otimes\cO_K/\rami)$ is linear.
\end{enumerate}\end{dfn}


\begin{rmk}\noindent The ramification ideal $\rami\subset\cO_K$ is fixed throughout the rest of the
chapter.
In the following we speak simply of elliptic shtukas instead of elliptic shtukas
of ramification ideal~$\rami$.\end{rmk}

\begin{example*} Let $\cO_F = \Fq[[z]]$, $\cO_K = \Fq[[\zeta]]$.
We have
$\cO_F\complot\cO_K = \Fq[[z,\zeta]]$.
The endomorphism $\tau$ of this ring preserves $\Fq[[z]]$ and sends $\zeta$ to
$\zeta^q$.
%
Our ramification ideal $\rami \subset
\Fq[[\zeta]]$ 
will be
the ideal generated by $\zeta$.

Consider the shtuka
\begin{equation*}
\cM =
\Big[\,\,
\Fq[[z,\zeta]] \shtuka{\quad\zeta-z\quad}{z\zeta\cdot\tau} \Fq[[z,\zeta]]
\,\,\Big].
\end{equation*}
In fact $\cM$ is (a part of) a model of the Carlitz module. We claim that $\cM$
is an elliptic shtuka of ramification ideal $\rami$. 
Indeed $\cM$ is locally free by construction. Furthermore
\begin{align*}
\cM(\cO_F/\fmF \otimes K) &=
\Big[\,\, \Fq(\!(\zeta)\!) \shtuka{\quad\zeta\quad}{0} \Fq(\!(\zeta)\!)
\,\,\Big] \\
\cM(F \otimes \cO_K/\fm_K) &=
\Big[\,\, \Fq(\!(z)\!) \shtuka{\,\,\,\,\,-z\,\,\,\,\,}{0} \Fq(\!(z)\!)
\,\,\Big].
\end{align*}
Hence the restrictions of $\cM$ to $\cO_F/\fmF \otimes K$ and
$F\otimes\cO_K/\fm_K$ are nilpotent.
The cohomology of $\Der\cM$ is easy to compute:
\begin{equation*}
\uH^1(\Der\cM) = \Fq[[z,\zeta]]/(\zeta-z) = \Fq[[\zeta]]. 
\end{equation*}
The element $z \in \Fq[[z]]$ acts on $\uH^1(\Der\cM)$ by
multiplication by $\zeta$.
So $\fmF \cdot \uH^1(\Der\cM) = \rami \cdot \uH^1(\Der\cM)$.
Finally the linearity condition holds since
\begin{equation*}
\cM(\cO_F \otimes\cO_K/\rami) =
\Big[\,\, \Fq[[z]] \shtuka{\quad-z\quad}{0} \Fq[[z]]
\,\,\Big].
\end{equation*}
%
\end{example*}

\begin{prp}\label{ellder} If a shtuka $\cM$ is elliptic then so is $\Der\cM$.\end{prp}

\pf Indeed the functor $\Der$ commutes with arbitrary restrictions and preserves
nilpotence so that $\Der\cM$ satisfies the conditions \eqref{ell1} and \eqref{ell2}
of Definition \ref{defell}. The condition \eqref{ell3} is tautologically satisfied and
\eqref{ell4} follows since the shtuka $\Der\cM$ is already linear.\quod

\begin{thm}\label{ellcoh} If $\cM$ is an elliptic shtuka then the following
holds:
\begin{enumerate}
\item $\uH^0(\cM) = 0$

\item $\uH^1(\cM)$ is a finitely generated free $\cO_F$-module.
%
\end{enumerate}
\end{thm}

\pf Indeed $\cM$ is locally free by \eqref{ell0} and $\cM(\cO_F/\fmF \otimes K)$ is nilpotent by \eqref{ell1} so the result
follows from Theorem \ref{llcoh}.
\quod

%

\section{Twists and quotients}
\label{sec:elltwist}

\begin{dfn}\label{defelltwist}%
Let $I \subset \cO_F\complot\cO_K$ be a $\tau$-invariant ideal.
Let $\cM$ be a shtuka on $\cO_F\complot\cO_K$
given by a diagram
\begin{equation*}
M_0 \shtuka{i}{j} M_1.
\end{equation*}
We define
\begin{equation*}
I\cM =
\Big[ IM_0 \shtuka{i}{j} IM_1 \Big].
\end{equation*}
We call $I\cM$ the \emph{twist} of $\cM$ by $I$.
The shtuka $I\cM$ comes equipped with a natural embedding $I\cM \hookrightarrow \cM$.
We denote the quotient $\cM/I$.

\afterall\noindent
\textbf{Warning.} 
Given invariant ideals
$I, J \subset\cO_F\complot\cO_K$ we denote $I\cM/J$ the quotient
$(I\cM)/(IJ\cM)$. In other words we first do the twist by $I$ and then take the
quotient by $J$.\end{dfn}

\breakflow
We will use the following invariant ideals:
\begin{itemize}
%
%
%
%
\item $\fmF^n\complot\cO_K$ for $n \geqslant 0$.

\item $\cO_F\complot \ibase$ for $\ibase \subset \cO_K$ an open ideal.

\item $\fmF^n\complot\cO_K+\cO_F\complot \ibase$ for $\fmF^n$ and $\ibase$ as above.

\item $\fmF^n \complot \ibase$ for $\fmF^n$ and $\ibase$ as above.
\end{itemize}
To simplify the notation we will write $\fmF^n \cM$ instead of
$(\fmF^n\complot\cO_K) \cM$. The same applies to $\ibase$, $\fmF^n + \ibase$ and the
quotients by the ideals of these three types. The twist of $\cM$ by the last
ideal will be denoted $\fmF^n \ibase \cM$ and the quotient will be denoted
$\cM/\fmF^n \ibase$.

\begin{lem}\label{ellquot}
Let $n \geqslant 0$ and let $\ibase \subset\cO_K$ be an open ideal.
If $\cM$ is a shtuka on $\cO_F\complot\cO_K$ then
\begin{align*}
\cM/\fmF^n &= \cM(\cO_F/\fmF^n\otimes\cO_K), \\
\cM/\ibase &= \cM(\cO_F\otimes\cO_K/\ibase), \\
\cM/(\fmF^n+\ibase) &= \cM(\cO_F/\fmF^n \otimes\cO_K/\ibase).
\end{align*}\end{lem}

\pf By Lemma \ref{coeffcocart} we have $\cO_F/\fmF^n \otimes_{\cO_F}
(\cO_F\complot\cO_K) = \cO_F/\fmF^n \otimes \cO_K$ so the first formula holds.
In a similar way Lemma \ref{basecocart} implies the second formula.  The last
formula follows from the first two. \quod


\begin{prp}\label{elltwist}If $\cM$ is an elliptic shtuka then so is $\fmF\cM$.\end{prp}

\pf Indeed the shtukas $\cM$ and $\fmF\cM$ are isomorphic.\quod


\begin{lem}\label{nilpbasequot}%
Let $\cM$ be a locally free shtuka on $\cO_F\complot\cO_K$ such that
$\cM(F\otimes\cO_K/\fm_K)$ is nilpotent. For every open ideal $\ibase \subset\cO_K$
the following holds:
\begin{enumerate}
\item $\cM(F\otimes\cO_K/\ibase)$ is nilpotent,

\item $\uH^0(\cM/\ibase) = 0$.
\end{enumerate}\end{lem}

\pf (1) It is enough to assume that $K$ is a single local field.
In this case $\ibase = \fm_K^n$ for some $n \geqslant 0$.
The ring $F\otimes\cO_K/\fm_K^n$ is
noetherian and complete with respect to the $\tau$-invariant ideal
$F\otimes\fm_K/\fm_K^n$.
Since $\cM(F\otimes\cO_K/\fm_K)$ is nilpotent Proposition
\lref{nilp}{nilpcomp} implies that $\cM(F\otimes\cO_K/\fm_K^n)$ is
nilpotent.

(2) Lemma \ref{ellquot} shows that $\cM/\ibase = \cM(\cO_F\otimes\cO_K/\ibase)$.
The natural map
\begin{equation*}
\uH^0(\cO_F\otimes\cO_K/\ibase,\,\cM) \to \uH^0(F\otimes\cO_K/\ibase,\,\cM)
\end{equation*}
is injective since $\cM$ is locally free. However the shtuka
$\cM(F\otimes\cO_K/\ibase)$ is nilpotent by (1). So Proposition
\lref{nilp}{nilpvanish} shows that $\uH^0(F\otimes\cO_K/\ibase,\,\cM) =0$.
\quod

\begin{lem}\label{ellnilptwist}%
Let $\cM$ be a shtuka on $\cO_F\complot\cO_K$.
If $\cM(F \otimes \cO_K/\fm_K)$ is nilpotent then
$(\rami\cM)(F\otimes \cO_K/\fm_K)$ is nilpotent.\end{lem}

\pf
We have a short exact sequence of shtukas
\begin{equation*}
0 \to \rami\cM/\fm_K \rami\cM \to \cM/\fm_K \rami \cM \to \cM/\fm_K\cM \to 0.
\end{equation*}
Using Lemma \ref{ellquot} we rewrite it as follows:
\begin{equation*}
0 \to (\rami\cM)(\cO_F\otimes\cO_K/\fm_K) \to
\cM(\cO_F\otimes\cO_K/\fm_K \rami) \to
\cM(\cO_F\otimes\cO_K/\fm_K) \to 0.
\end{equation*}
Localizing at a uniformizer of $\cO_F$ we get a short exact sequence
\begin{equation*}
0 \to (\rami\cM)(F\otimes\cO_K/\fm_K) \to
\cM(F\otimes\cO_K/\fm_K \rami) \to
\cM(F\otimes\cO_K/\fm_K) \to 0.
\end{equation*}
The third shtuka is nilpotent by assumption
while the second shtuka is nilpotent by Lemma
\ref{nilpbasequot}. Hence the first shtuka is nilpotent. \quod

\begin{prp}\label{elltwistvanish}%
Let $\cM$ be a shtuka on $\cO_F\complot\cO_K$.
For every $n \geqslant 0$ the shtuka $(\rami^n\cM)/\rami^n$ is linear.\end{prp}

\pf Suppose that $\cM$ is given by a diagram
\begin{equation*}
M_0 \shtuka{i}{j} M_1.
\end{equation*}
We need to prove that $j(\rami^n M_0) \subset \rami^{2n}M_1$.
The endomorphism $\tau$ of $\cO_F\complot\cO_K$ sends
$\rami$ to $\rami^q$. Since
$j$ is $\tau$-linear it follows that
$j(\rami^n M_0) \subset \rami^{n q} M_1$.
The result follows since $q > 1$. \quod

\begin{prp}\label{ellbasetwist}%
If $\cM$ is an elliptic shtuka
then so is $\rami\cM$.\end{prp}

\pf
%
Let us verify the conditions of Definition \ref{defell} for $\rami\cM$.
\begin{enumerate}
\setcounter{enumi}{-1}
\renewcommand{\theenumi}{E\arabic{enumi}}
\item 
Lemma \ref{basecocart} shows that $\cO_F\complot \rami$ is a
free $\cO_F\complot\cO_K$-module of rank $1$. Hence $\rami\cM$ is a locally free
shtuka. 

\item The shtukas $\cM$ and $\rami\cM$ coincide on
$(\cO_F\complot\cO_K)\otimes_{\cO_K} K$
so the restriction of $\rami\cM$ to $\cO_F/\fm_F\otimes K$
is nilpotent.

\item Follows by Lemma \ref{ellnilptwist}.

\item Consider the short exact sequence of shtukas
\begin{equation*}
0 \to \Der \rami \cM \to \Der\cM \to \Der\cM/\rami \to 0.
\end{equation*}
The module $\uH^0(\Der\cM/\rami)$ vanishes by Lemma \ref{nilpbasequot}. Taking
the cohomology sequence we conclude that the natural map $\uH^1(\Der \rami\cM)
\to \uH^1(\Der\cM)$ is injective. This map is $\cO_F\otimes\cO_K$-linear by
construction. Now the image of $\uH^1(\Der \rami\cM)$ in $\uH^1(\Der\cM)$ is
$\rami\cdot \uH^1(\Der\cM)$ by definition. Therefore
\begin{equation*}
\fmF \cdot \uH^1(\Der \rami\cM) =
\fmF \rami \cdot \uH^1(\Der\cM) =
\rami \fmF \cdot \uH^1(\Der\cM) =
\rami \rami \cdot \uH^1(\Der\cM) =
\rami \cdot \uH^1(\Der \rami\cM).
\end{equation*}

\item Indeed the shtuka $(\rami\cM)/\rami$ is
linear according to Proposition \ref{elltwistvanish}.
\quod
\end{enumerate}

\section{Filtration on cohomology}
\label{sec:ellcohfilt}

An elliptic shtuka $\cM$ carries a natural filtration by elliptic subshtukas
$\rami^n\cM$.
In this section we would like to describe the induced filtration
on $\uH^1(\cM)$. 
If the elliptic shtuka $\cM$ is linear then
\begin{equation*}
\fmF \cdot \uH^1(\cM) = \rami \cdot \uH^1(\cM) = \uH^1(\rami\cM)
\end{equation*}
by the condition \eqref{ell3}. As a consequence the filtration on $\uH^1(\cM)$ induced by
$\rami^n\cM$ is the filtration by powers of $\fmF$. Our goal is to prove that
the same is true without the linearity assumption. 

%

\begin{lem}\label{ellmfquot}
Let $\cM$ be an elliptic shtuka and let $n \geqslant 0$. The shtuka $\cN = \cM/\fmF^n$
has the following properties:
\begin{enumerate}
\item $\cN$ is a locally free $\cO_F/\fmF^n \otimes\cO_K$-shtuka.

\item $\cN(\cO_F/\fmF \otimes K)$ is nilpotent.

\item $\rami^n \cdot \uH^1(\Der\cN) = 0$.
\end{enumerate}\end{lem}

\pf By Lemma \ref{ellquot} we have $\cN = \cM(\cO_F/\fmF^n \otimes \cO_K)$.
So (1) and (2) follow since $\cM$ is an elliptic shtuka.

The natural map
$\uH^1(\Der\cM) \to \uH^1(\Der\cN)$
is surjective 
and $\cO_F\otimes\cO_K$-linear.
According to the condition \eqref{ell3} we have
$\rami\cdot\uH^1(\Der\cM) = \fmF\cdot\uH^1(\Der\cM)$.
As a consequence
$\rami^n \cdot \uH^1(\Der\cN) = \fmF^n\cdot\uH^1(\Der\cN)$.
However $\fmF^n$ acts on this module by zero since $\Der\cN$ is a shtuka on
$\cO_F/\fmF^n\otimes\cO_K$. We thus get (3). \quod

\begin{prp}\label{ellcoefftwistcoh}%
Let $\cM$ be an elliptic shtuka and let $n \geqslant 0$.
\begin{enumerate}
\item The sequence
$0 \to \uH^1(\fmF^n\cM) \to \uH^1(\cM) \to \uH^1(\cM/\fmF^n) \to 0$
is exact.

\item The image of $\uH^1(\fmF^n\cM)$ in $\uH^1(\cM)$ is $\fmF^n\uH^1(\cM)$.
\end{enumerate}\end{prp}

\pf (1) By Lemma \ref{ellmfquot} the quotient $\cM/\fmF^n$ is a
locally free $\cO_F/\fmF^n \otimes\cO_K$-shtuka whose restriction to
$\cO_F/\fmF\otimes K$ is nilpotent. Hence $\uH^0(\cM/\fmF^n) = 0$ by Theorem \ref{artnilpcoh}.
Taking the cohomology sequence of the short exact sequence
$0 \to \fmF^n\cM \to \cM \to \cM/\fmF^n \to 0$
we get the result. (2) is clear.\quod


\begin{lem}\label{ellmfquotcoh}%
Let $\cM$ be an elliptic shtuka and let $n
\geqslant 0$. If $\cM/\rami^n$ is linear then the shtuka
$\cN = \cM/\fmF^n$ has the following properties:
\begin{enumerate}
\item The artinian regulator is defined for $\cN$ (Definition \ref{defartreg}).

\item The reduction map $\uH^1(\cN) \to \uH^1(\cN/\rami^n)$ is an isomorphism.
\end{enumerate}\end{lem}

\pf We claim that the shtuka $\cN$ has the following properties:
\begin{enumerate}
\item[(i)] $\cN$ is a locally free $\cO_F/\fmF^n \otimes \cO_K$-module shtuka whose
restriction to $\cO_F/\fmF \otimes K$ is niltpotent.

\item[(ii)] $\rami^n \cdot \uH^1(\Der\cN) = 0$.

\item[(iii)] $\cN/\rami^n$ is linear.
\end{enumerate}
Indeed Lemma \ref{ellmfquot} implies (i) and (ii) while (iii) follows since
the quotient $\cM/\rami^n$ is linear by assumption.
We then apply Proposition \ref{artregexist} to $\cN$
with $\ngC = \cO_F/\fmF^n$ and $\ibase = \rami^n$ and conclude that $\cN$ has the
properties (1) and (2).\quod

\begin{lem}\label{ellbasequotred}%
If $\cM$ is an elliptic shtuka then the reduction map $\uH^1(\cM/\rami) \to
\uH^1(\cM/(\fmF + \rami))$ is an isomorphism.\end{lem}

\pf The short exact sequence of shtukas $0 \to (\fmF\cM)/\rami \to \cM/\rami \to
\cM/(\fmF + \rami) \to 0$ induces an exact sequence of cohomology
\begin{equation*}
\uH^1(\fmF\cM/\rami) \to \uH^1(\cM/\rami) \to \uH^1(\cM/(\fmF+\rami)) \to 0.
\end{equation*}
So to prove the lemma it is enough to demonstrate that the first map in this
sequence is zero.

The shtukas $\cM/\rami$ and $(\fmF\cM)/\rami$ are linear since $\cM$ is elliptic.
So we can assume without loss of generality that $\cM$ is itself linear.
The natural map $\uH^1(\cM) \to \uH^1(\cM/\rami)$ is $\cO_F\otimes\cO_K$-linear
and surjective. Furthermore $\fmF \cdot \uH^1(\cM) = \rami \cdot \uH^1(\cM)$
by definition of $\rami$. As a consequence
\begin{equation*}
\fmF \cdot \uH^1(\cM/\rami) = \rami \cdot \uH^1(\cM/\rami).
\end{equation*}
However $\cM/\rami$ is a linear shtuka on $\cO_F\otimes\cO_K/\rami$ so $\rami$ acts on
$\uH^1(\cM/\rami)$ by zero. Thus $\fmF \cdot \uH^1(\cM/\rami) = 0$ which implies
that the natural map $\uH^1(\fmF\cM/\rami) \to \uH^1(\cM/\rami)$ is zero. \quod


\begin{thm}\label{ellcohfilt}%
Let $\cM$ be an elliptic shtuka and let $n \geqslant 0$.
\begin{enumerate}
\item The sequence
$0 \to \uH^1(\rami^n\cM) \to \uH^1(\cM) \to \uH^1(\cM/\rami^n) \to 0$
is exact.

\item The image of $\uH^1(\rami^n\cM)$ in $\uH^1(\cM)$ is $\fmF^n\uH^1(\cM)$.
\end{enumerate}\end{thm}

\breakflow
So as we claimed at the beginning of this section
the filtration on $\uH^1(\cM)$ induced by the subshtukas $\rami^n\cM$
is the filtration by powers of $\fmF$.

\afterall\noindent
\textit{Proof of Theorem \ref{ellcohfilt}.}
(1) Lemma \ref{nilpbasequot} claims that $\uH^0(\cM/\rami^n) = 0$.
So
the natural sequence above is exact.
(2) By Proposition \ref{ellbasetwist} the shtuka $\rami^n\cM$ is elliptic.
It is thus enough to treat the case $n = 1$.
Consider the natural commutative square
\begin{equation*}
\xymatrix{
\uH^1(\cM) \ar[r] \ar[d] & \uH^1(\cM/\rami) \ar[d] \\
\uH^1(\cM/\fmF) \ar[r] & \uH^1(\cM/\fmF + \rami)
}
\end{equation*}
We just demonstrated that the top arrow in this square is surjective with kernel
$\uH^1(\rami\cM)$. According to Proposition \ref{ellcoefftwistcoh} the left
arrow is surjective with kernel $\fmF \uH^1(\cM)$. The right arrow is an
isomorphism by Lemma \ref{ellbasequotred}. Since $\cM/\rami$ is linear by \eqref{ell4}
Lemma \ref{ellmfquotcoh} shows that the bottom arrow is an isomorphism.
So the result follows. \quod

%

\section{Overview of regulators}
\label{sec:ellregoverview}

Now as our discussion of elliptic shtukas has gained some substance we can give
an overview of the regulator theory which will follow.
Let $\cM$ be an elliptic shtuka. If $\cM$ is linear then one tautologically has
a natural isomorphism
$\rho\colon\uH^1(\cM) \to \uH^1(\Der\cM)$
induced by the identity of the shtukas $\cM$ and $\Der\cM$, the \emph{regulator}
of $\cM$. We would like to
extend it to all elliptic shtukas $\cM$. A rough idea is to
approximate $\cM$ with linear pieces.
\begin{dfn*} A natural transformation
\begin{equation*}
\uH^1(\cM) \xrightarrow{\quad\rho\quad} \uH^1(\Der\cM)
\end{equation*}
of functors
on the category of elliptic shtukas is called a \emph{regulator} if for every
$\cM$ such that $\cM/\rami^{2n}$ is linear the diagram
\begin{equation*}
\xymatrix{
\uH^1(\cM) \ar[rr]^{\rho} \ar[d]_{\textup{red.}} &&\uH^1(\Der\cM)
\ar[d]^{\textup{red.}} \\
\uH^1(\cM/\rami^n) \ar[rr]^1 && \uH^1(\Der\cM/\rami^n)
}
\end{equation*}
is commutative.
\end{dfn*}

\breakflow
Note that we demand the quotient $\cM/\rami^{2n}$ to be
linear, not just $\cM/\rami^n$.
Even though the square above makes sense if merely $\cM/\rami^n$ is linear
the regulator will fail to exist if one demands all such squares to commute.

A regulator $\rho$ will send the submodule $\uH^1(\rami^n\cM)
\subset \uH^1(\cM)$ to $\uH^1(\Der\rami^n\cM) \subset \uH^1(\Der\cM)$ since
it is natural.
By Proposition \ref{elltwistvanish}
the subquotients
$(\rami^{2n}\cM)/\rami^{2n} = \rami^{2n}\cM/\rami^{4n}\cM$ are linear for all $n \geqslant 0$.
Hence the regulator is an isomorphism and is unique.
However the existence of the regulator is a different
question altogether.
Its construction occupies the rest of the chapter.

%
%
In Sections \ref{sec:ellcohsubquot} and \ref{sec:ellfiltiso} we present some
auxillary results on cohomology of elliptic shtukas and their subquotients.
These sections are of a technical nature.
The level of technicality reaches its high point in Section
\ref{sec:elltautreg} where
we study the tautological regulators. They are the basic building blocks
for the regulator on $\uH^1(\cM)$.
Already in Section \ref{sec:ellfpreg} the statements and proofs become
much more natural.
The construction of the regulator on $\uH^1(\cM)$ in Section
\ref{sec:ellreg} is actually quite simple. Besides this construction, the main
results of Section \ref{sec:ellreg} are Theorem \ref{ellregcmp} which gives a criterion
for an $\cO_F$-linear map $f\colon\uH^1(\cM) \to \uH^1(\Der\cM)$ to coincide
with the regulator, and Theorem \ref{ellregart} which serves as a link with the trace
formula of Chapter \ref{chapter:trace}.


%
%

\section{Cohomology of subquotients}
\label{sec:ellcohsubquot}

In this section we present
several technical propositions which will be used in the construction
of the regulator map.
First we use Theorem \ref{ellcohfilt} to derive two exact sequences of
cohomology modules.

\begin{prp}\label{ellcoeffdotseqprp}
Let $n \geqslant 0$.
If $\cM$ is an elliptic shtuka
then the short exact sequence of shtukas
\begin{equation*}
0 \to \rami^n\cM/\fmF^n \to \cM/\fmF^n \to \cM/(\fmF^n + \rami^n) \to 0
\end{equation*}
induces an exact sequence of cohomology modules
\begin{multline}\label{ellcoeffdotseq}
0 \to \uH^0(\cM/(\fmF^n + \rami^n)) \xrightarrow{\,\,\delta\,\,}
\uH^1(\rami^n\cM/\fmF^n) \xrightarrow{\,\,0\,\,} \\
\xrightarrow{\,\,\phantom{0}\,\,}
\uH^1(\cM/\fmF^n) \xrightarrow{\,\textup{red.}\,}
\uH^1(\cM/(\fmF^n + \rami^n)) \to 0.
\end{multline}
The middle map in this sequence is zero and the adjacent maps are isomorphisms.
\end{prp}

\pf According to Lemma \ref{ellmfquot} the shtuka $\cM/\fmF^n$ is a locally free
$\cO_F/\fmF^n \otimes \cO_K$-module shtuka whose restriction to
$\cO_F/\fmF \otimes K$ is nilpotent. Theorem \ref{artnilpcoh} implies
that $\uH^0(\cM/\fmF^n) = 0$.
Therefore the sequence \eqref{ellcoeffdotseq} is exact. To prove the result it
is enough to show that the middle map in this sequence is zero. This map fits
into a natural commutative square
\begin{equation*}
\xymatrix{
\uH^1(\rami^n\cM) \ar@{-->}[r] \ar@{-->}[d] & \uH^1(\cM) \ar@{-->}[d] \\
\uH^1(\rami^n\cM/\fmF^n) \ar[r] & \uH^1(\cM/\fmF^n)
}
\end{equation*}
By Proposition \ref{ellcoefftwistcoh} the vertical maps in this square are
surjections and the kernel of the right map is $\fmF^n\uH^1(\cM)$. However
Theorem \ref{ellcohfilt} shows that the image of the top map is
$\fmF^n\uH^1(\cM)$ whence the composition of the top and the right maps is zero.
Since the left map is surjective we conclude that the bottom map is zero. \quod

\begin{prp}\label{ellbasedotseqprp}
Let $n \geqslant 0$.
If $\cM$ is an elliptic shtuka
then the short exact sequence of shtukas
\begin{equation*}
0 \to \fmF^n \cM/\rami^n \to \cM/\rami^n \to \cM/(\fmF^n + \rami^n) \to 0
\end{equation*}
induces an exact sequence of cohomology modules
\begin{multline}\label{ellbasedotseq}
0 \to \uH^0(\cM/(\fmF^n + \rami^n)) \xrightarrow{\,\,\delta\,\,}
\uH^1(\fmF^n\cM/\rami^n) \xrightarrow{\,\,0\,\,} \\
\xrightarrow{\,\,\phantom{0}\,\,}
\uH^1(\cM/\rami^n) \xrightarrow{\,\textup{red.}\,}
\uH^1(\cM/(\fmF^n + \rami^n)) \to 0.
\end{multline}
The middle map in this sequence is zero and the adjacent maps are isomorphisms.
\end{prp}

\pf The cohomology module 
$\uH^0(\cM/\rami^n)$ vanishes
according to Lemma \ref{nilpbasequot}.
Thus the sequence \eqref{ellbasedotseq} is exact. To prove the proposition
it is enough to show that the middle map in \eqref{ellbasedotseq} is zero. This
map fits into a natural commutative square
\begin{equation*}
\xymatrix{
\uH^1(\fmF^n\cM) \ar@{-->}[r] \ar@{-->}[d] & \uH^1(\cM) \ar@{-->}[d] \\
\uH^1(\fmF^n\cM/\rami^n) \ar[r] & \uH^1(\cM/\rami^n).
}
\end{equation*}
Theorem \ref{ellcohfilt} shows that the vertical maps are surjective
and that the kernel of the right map is $\fmF^n\uH^1(\cM)$. So the composition
of the top and the right maps is zero. As the left map 
is surjective we conclude that the bottom map is zero. \quod

\breakflow
Next we study the boundary homomorphisms $\delta$ of the sequences
\eqref{ellcoeffdotseq} and \eqref{ellbasedotseq}.

\begin{prp}\label{elldelta}%
Let $n \geqslant 0$. If $\cM$ is an elliptic shtuka then the natural
diagram
\begin{equation*}
\xymatrix{
\uH^0(\cM/(\fmF^n + \rami^n)) \ar[d]_\delta \ar[r]^\delta &
\uH^1(\fmF^n\cM/\rami^n) \ar[d] \\
\uH^1(\rami^n\cM/\fmF^n) \ar[r] & \uH^1(\cM/\fmF^n\rami^n)
}
\end{equation*}
is \emph{anticommutative}. Here the vertical $\delta$ is the boundary
homomorphism of \eqref{ellcoeffdotseq} while the horizontal $\delta$ is the
boundary homomorphism of \eqref{ellbasedotseq}.\end{prp}

\pf To improve legibility we will write $\ibase$ for $\rami^n$ and $\icoef$ for $\fmF^n$.
Since the shtuka $\cM$ is locally free Lemma \ref{ellintersect} implies that the
natural sequence
\begin{equation*}
0 \to \ibase\cM/\icoef \oplus \icoef\cM/\ibase \to \cM/\ibase\icoef \to \cM/(\ibase+\icoef) \to 0
\end{equation*}
is exact. Now consider the natural diagram
\begin{equation*}
\xymatrix{
0 \ar[r] & \ibase\cM/\icoef \oplus \icoef\cM/\ibase \ar[d]^1 \ar[r] &
\cM/\ibase\icoef \ar[r] \ar[d]^{(\textup{red.},\,\textup{red.})} &
\cM/(\ibase+\icoef) \ar[r] \ar[d]^{\Delta} & 0 \\
0 \ar[r] & \ibase\cM/\icoef \oplus \icoef\cM/\ibase \ar[r] &
\cM/\ibase \oplus \cM/\icoef \ar[r] & (\cM/\ibase+\icoef)^{\oplus 2} \ar[r] & 0
}
\end{equation*}
Here $\Delta$ means the diagonal map. Observe that the lower row is the direct
sum of the short exact sequences
\begin{align*}
&0 \to \icoef\cM/\ibase \to \cM/\icoef \to \cM/(\ibase+\icoef) \to 0, \\
&0 \to \ibase\cM/\icoef \to \cM/\ibase \to \cM/(\ibase+\icoef) \to 0
\end{align*}
which give rise to the cohomology sequences \eqref{ellcoeffdotseq} and
\eqref{ellbasedotseq} respectively.
The diagram above is clearly commutative. So it induces a morphism of long
exact sequences. A part of it looks like this:
\begin{equation*}
\xymatrix{
\uH^0(\cM/(\ibase\!+\!\icoef)) \ar[d]^\Delta \ar[r] &
\uH^1(\ibase\cM/\!\icoef) \!\oplus\! \uH^1(\icoef\cM/\ibase) \ar[r] \ar[d]^1 &
\uH^1(\cM/\ibase\icoef) \ar[d]^{(\textup{red.},\,\textup{red.})} \\
\uH^0(\cM/(\ibase\!+\!\icoef))^{\oplus 2} \ar[r]^{\delta\oplus\delta\quad\,\,\,\,} &
\uH^1(\ibase\cM/\!\icoef) \!\oplus\! \uH^1(\icoef\cM/\ibase) \ar[r] &
\uH^1(\cM/\!\icoef) \!\oplus\! \uH^1(\cM/\ibase)
}
\end{equation*}
As a consequence the boundary homomorphism
$\uH^0(\cM/(\ibase+\icoef)) \to \uH^1(\ibase\cM/\icoef) \oplus \uH^1(\icoef\cM/\ibase)$
in the top row coincides with $(\delta,\delta)$. Since the composition of the
adjacent homomorphisms
\begin{equation*}
\uH^0(\cM/(\ibase+\icoef)) \xrightarrow{\,(\delta,\delta)\,}
\uH^1(\ibase\cM/\icoef) \oplus \uH^1(\icoef\cM/\ibase) \xrightarrow{\quad}
\uH^1(\cM/\ibase\icoef)
\end{equation*}
is zero we get the result. \quod

\begin{prp}\label{ellcrossinj}
Let $n \geqslant 0$.
If $\cM$ is an elliptic shtuka then the natural maps
\begin{align*}
&\uH^1(\rami^n\cM/\fmF^n) \to \uH^1(\cM/\fmF^n\rami^n), \\
&\uH^1(\fmF^n\cM/\rami^n) \to \uH^1(\cM/\fmF^n\rami^n)
\end{align*}
are injective.\end{prp}

\pf Consider the short exact sequence
$0 \to \rami^n\cM/\fmF^n \to \cM/\fmF^n\rami^n \to \cM/\rami^n \to 0$. By Lemma
\ref{nilpbasequot} the module $\uH^0(\cM/\rami^n)$ vanishes. Hence the natural map
$\uH^1(\rami^n\cM/\fmF^n) \to \uH^1(\cM/\fmF^n\rami^n)$ is injective. In a
similar way Theorem \ref{artnilpcoh} implies that the map
$\uH^1(\fmF^n\cM/\rami^n) \to \uH^1(\cM/\fmF^n\rami^n)$ is injective. \quod

\section{Natural isomorphisms on cohomology}
\label{sec:ellfiltiso}

Let $\cM$ be an elliptic shtuka. According to Theorem \ref{ellcohfilt} the
images of $\uH^1(\fmF\cM)$ and $\uH^1(\rami\cM)$ inside $\uH^1(\cM)$ are equal
to $\fmF\uH^1(\cM)$.

\begin{dfn}\label{defellfiltstep}
We define the natural isomorphism
\begin{equation*}
\gamma_\cM\colon\uH^1(\fmF\cM) \to \uH^1(\rami\cM) 
\end{equation*}
as the unique map which makes the triangle
\begin{equation*}
\xymatrix{
\uH^1(\fmF\cM) \ar@{-->}[rr]^{\gamma_\cM}_{\bisosign}
\ar[dr] && \uH^1(\rami\cM) \ar[dl] \\
& \uH^1(\cM) &
}
\end{equation*}
commutative.\end{dfn}

\breakflow
To simplify the expressions we denote $\gamma^n_\cM$ the composition
\begin{equation*}
\uH^1(\fmF^n \cM) \xrightarrow{\,\gamma_{\fmF^{n-1}\cM}\,}
\uH^1(\rami\fmF^{n-1}\cM) \xrightarrow{\,\gamma_{\rami\fmF^{n-2}\cM}\,} \dotsc
\xrightarrow{\,\gamma_{\rami^{n-1}\cM}\,}
\uH^1(\rami^n\cM).
\end{equation*}
We write $\gamma^n$ instead of $\gamma^n_\cM$ if the corresponding shtuka $\cM$
is clear from the context.
The map $\gamma^n_\cM$ will only
be used in Sections \ref{sec:elltautreg} and \ref{sec:ellfpreg}. 

\begin{lem}
If $n, k \geqslant 0$ are integers then the composition
\begin{equation*}
\uH^1(\fmF^{n+k}\cM) \xrightarrow{\,\gamma^n_{\fmF^k\cM}\,}
\uH^1(\rami^n \fmF^k\cM) \xrightarrow{\,\gamma^k_{\rami^n\cM}\,}
\uH^1(\rami^{n+k}\cM)
\end{equation*}
is equal to $\gamma^{n+k}_\cM$.\quod\end{lem}


\breakflow
Our goal is to relate
$\gamma^n$ to the maps which appear in the natural sequences
\eqref{ellcoeffdotseq} and \eqref{ellbasedotseq}.

\begin{lem}\label{ellbargamma}
Let $\cM$ be an elliptic shtuka. For every $n \geqslant 0$ there exists a unique
natural isomorphism
$\overline{\gamma}_n\colon \uH^1(\fmF^n\cM/\rami^n) \xrightarrow{\isosign}
\uH^1(\rami^n\cM/\fmF^n)$
such that the diagram
\begin{equation*}
\xymatrix{
\uH^1(\fmF^n\cM) \ar[r]^{\gamma^n}_{\bisosign} \ar[d]_{\textup{red.}} &
\uH^1(\rami^n\cM) \ar[d]^{\textup{red.}} \\
\uH^1(\fmF^n\cM/\rami^n) \ar[r]^{\overline{\gamma}_n}_{\bisosign} &
\uH^1(\rami^n\cM/\fmF^n) 
}
\end{equation*}
is commutative.\end{lem}

\breakflow
The map $\overline{\gamma}_n$ will only be used in Section \ref{sec:elltautreg}.
The same remark applies to the related map $\overline{\varepsilon}_n$ which we
introduce below.

\afterall\noindent
\textit{Proof of Lemma \ref{ellbargamma}}. By Theorem \ref{ellcohfilt} the left
reduction map is surjective with kernel $\fmF^n\uH^1(\fmF^n\cM)$.
Similarly Proposition \ref{ellcoefftwistcoh} shows that the right
reduction map is surjective with kernel $\fmF^n\uH^1(\rami^n\cM)$. Since
$\gamma^n$ is an isomorphism the result follows. \quod

\begin{prp}\label{ellgammadelta}
Let $\cM$ be an elliptic shtuka. For every $n \geqslant 0$ the
maps $\delta$ of \eqref{ellcoeffdotseq}, \eqref{ellbasedotseq} are isomorphisms
and the
triangle
\begin{equation*}
\xymatrix{
& \uH^0(\cM/\fmF^n+\rami^n) \ar[dl]_\delta \ar[dr]^\delta & \\
\uH^1(\fmF^n\cM/\rami^n) \ar[rr]^{\overline{\gamma}_n} &&
\uH^1(\rami^n\cM/\fmF^n)
}
\end{equation*}
is \emph{anticommutative}. 
\end{prp}

\pf The maps $\delta$ are isomorphisms by Propositions \ref{ellcoeffdotseqprp}
and \ref{ellbasedotseqprp}. Next, consider the diagram
\begin{equation*}
\xymatrix{
& \uH^0(\cM/\fmF^n+\rami^n) \ar[dl]_\delta \ar[dr]^\delta & \\
\uH^1(\fmF^n\cM/\rami^n) \ar[rr]^{\overline{\gamma}_n}
\ar[dr] &&
\uH^1(\rami^n\cM/\fmF^n) \ar[dl] \\
& \uH^1(\cM/\fmF^n\rami^n) &
}
\end{equation*}
The square in this diagram commutes by Proposition \ref{elldelta}. The bottom
diagonal maps are injective by Proposition \ref{ellcrossinj}. Hence it is
enough to prove that the bottom triangle commutes. By definition of
$\overline{\gamma}_n$ we need to show the commutativity of the triangle
\begin{equation*}
\xymatrix{
\uH^1(\fmF^n\cM) \ar[rr]^{\gamma^n}
\ar[dr] &&
\uH^1(\rami^n\cM) \ar[dl] \\
& \uH^1(\cM/\fmF^n\rami^n) &
}
\end{equation*}
However the two diagonal maps factor over the natural maps to $\uH^1(\cM)$.
Hence the commutativity follows from the definition of $\gamma^n$.\quod

\begin{lem} Let $\cM$ be an elliptic shtuka. For every $n \geqslant 0$ there
exists a unique natural isomorphism
$\overline{\varepsilon}_n\colon
\uH^1(\rami^n\cM/\fmF^n) \xrightarrow{\isosign}
\uH^1(\rami^n\cM/\rami^n)$
such that the square
\begin{equation*}
\xymatrix{
\uH^1(\rami^n\cM) \ar@{=}[r] \ar[d]_{\textup{red.}} &
\uH^1(\rami^n\cM) \ar[d]^{\textup{red.}} \\
\uH^1(\rami^n\cM/\fmF^n) \ar[r]^{\overline{\varepsilon}_n}_{\bisosign} &
\uH^1(\rami^n\cM/\rami^n)
}
\end{equation*}
is commutative.\end{lem}

\pf Indeed the left reduction map is surjective with kernel
$\fmF^n\uH^1(\rami^n\cM)$ by Proposition \ref{ellcoefftwistcoh} while the right
reduction map is surjective with the same kernel by Theorem
\ref{ellcohfilt}.\quod

\begin{lem}\label{ellgammaepsilon}
Let $\cM$ be an elliptic shtuka. For every $n \geqslant 0$ the square
\begin{equation*}
\xymatrix{
\uH^1(\fmF^n\cM) \ar[rr]^{\gamma^n} \ar[d]_{\textup{red.}} &&
\uH^1(\rami^n\cM) \ar[d]^{\textup{red.}} \\
\uH^1(\fmF^n\cM/\rami^n) \ar[rr]^{\overline{\varepsilon}_n \circ
\overline{\gamma}_n} && \uH^1(\rami^n\cM/\rami^n)
}
\end{equation*}
is commutative.\end{lem}

\pf Follows instantly from the definitions of $\overline{\gamma}_n$ and
$\overline{\varepsilon}_n$.\quod

\begin{prp}\label{ellepsilonred} Let $\cM$ be an elliptic shtuka and let $n
\geqslant 0$.
\begin{enumerate}
\item The reduction maps
\begin{align*}
&\uH^1(\rami^n\cM/\fmF^n) \xrightarrow{\,\textup{red.}\,}
\uH^1(\rami^n\cM/\fmF^n + \rami^n), \\
&\uH^1(\rami^n\cM/\rami^n) \xrightarrow{\,\textup{red.}\,}
\uH^1(\rami^n\cM/\fmF^n + \rami^n)
\end{align*}
are isomorphisms.

\item
The diagram
\begin{equation*}
\xymatrix{
\uH^1(\rami^n\cM/\fmF^n) \ar[rr]^{\overline{\varepsilon}_n}
\ar[dr]_{\textup{red.}} &&
\uH^1(\rami^n\cM/\rami^n) \ar[dl]^{\textup{red.}} \\
& \uH^1(\rami^n\cM/\fmF^n+\rami^n)
}
\end{equation*}
is commutative.\end{enumerate}\end{prp}

\pf (1) follows from Propositions \ref{ellcoeffdotseqprp} and
\ref{ellbasedotseqprp}. (2) Indeed the outer square in the diagram
\begin{equation*}
\xymatrix{
& \uH^1(\rami^n\cM) \ar[dl]_{\textup{red.}} \ar[dr]^{\textup{red.}}& \\
\uH^1(\rami^n\cM/\fmF^n) \ar[rr]^{\overline{\varepsilon}_n}
\ar[dr]_{\textup{red.}} &&
\uH^1(\rami^n\cM/\rami^n) \ar[dl]^{\textup{red.}} \\
& \uH^1(\rami^n\cM/\fmF^n+\rami^n)
}
\end{equation*}
commutes by definition of $\overline{\varepsilon}_n$.
The top left reduction map is surjective
by Proposition
\ref{ellcoefftwistcoh}
So the result follows. \quod

\section{Tautological regulators}
\label{sec:elltautreg}

\begin{lem}\label{ellfprediso} Let $d \geqslant 1$. If $\cM$ is an elliptic
shtuka then the reduction map $\uH^1(\cM) \to \uH^1(\cM/\rami^d)$ induces an
isomorphism $\uH^1(\cM)/\fmF^d\to
\uH^1(\cM/\rami^d)$.
\end{lem}

\pf Indeed Theorem \ref{ellcohfilt} states that the reduction map is surjective
with kernel $\fmF^d \uH^1(\cM)$. \quod

\breakflow
Let $\cM$ be an elliptic shtuka and let $d \geqslant 1$. If $\cM/\rami^d$ is
linear then the shtukas $\cM/\rami^d$ and $\Der\cM/\rami^d$ coincide
tautologically. We thus get a natural isomorphism 
\begin{equation*}
\uH^1(\cM/\rami^d) \xrightarrow{\quad 1\quad} \uH^1(\Der\cM/\rami^d).
\end{equation*}
Using it 
we will now define a natural isomorphism $\uH^1(\cM)/\fmF^d
\xrightarrow{\isosign} \uH^1(\Der\cM)/\fmF^d$.

\begin{dfn}\label{elldeftautreg}
Let $d \geqslant 1$ and let $\cM$ be an elliptic shtuka such that
$\cM/\rami^d$ is linear. We define the map $\overline{\rho}_d$ by the
commutative diagram
\begin{equation*}
\xymatrix{
\uH^1(\cM)/\fmF^d \ar[rr]^{\overline{\rho}_d} \ar[d]_{\textup{red.}}^{\rtviso} &&
\uH^1(\Der\cM)/\fmF^d \ar[d]^{\textup{red.}}_{\ltviso} \\
\uH^1(\cM/\rami^d) \ar[rr]^1 && \uH^1(\Der\cM/\rami^d).
}
\end{equation*}
Here the vertical maps are the isomorphisms of Lemma \ref{ellfprediso}.
We call $\overline{\rho}_d$ the \emph{tautological regulator of order $d$.}
By construction $\overline{\rho}_d$ is a natural isomorphism.\end{dfn}

\breakflow
For the duration of this section and Section \ref{sec:ellfpreg} we
fix a uniformizer $z \in \cO_F$. As it will be shown in Section \ref{sec:ellfpreg}
our results do
not depend on the choice of $z$. However this choice simplifies the exposition.
With the uniformizer $z$ fixed we have for every elliptic shtuka $\cM$ a natural
isomorphism $\varpi\colon\cM \to \fmF\cM$.

\begin{dfn}\label{defellslideiso} Let $\cM$ be an elliptic shtuka. We define a
natural isomorphism
\begin{equation*}
\uH^1(\cM) \xrightarrow{\quad s\quad} \uH^1(\rami\cM)
\end{equation*}
as the composition
$\uH^1(\cM) \xrightarrow{\varpi}
\uH^1(\fmF\cM) \xrightarrow{\gamma}
\uH^1(\rami\cM)$.
Here $\gamma$ is the natural isomorphism of Definition \ref{defellfiltstep}.
We call $s$ the \emph{sliding isomorphism}.\end{dfn}

\breakflow
Apart from this section the sliding isomorphism $s$ will only be used in Section
\ref{sec:ellfpreg}.

\begin{lem}\label{ellslidecomp}
Let $\cM$ be an elliptic shtuka. For every $n \geqslant 0$ the
natural diagram
\begin{equation*}
\xymatrix{
\uH^1(\cM) \ar[d]_{\varpi^n} \ar[rd]^{s^n} & \\
\uH^1(\fmF^n\cM) \ar[r]^{\gamma^n} & \uH^1(\rami^n\cM)
}
\end{equation*}
is commutative.\end{lem}

\pf This diagram commutes for $n = 0$. Assuming that
it commutes for some $n$ we prove that it does so for $n + 1$.
Consider the natural diagram
\begin{equation*}
\xymatrix{
\uH^1(\cM) \ar[r]^{\varpi^n} \ar[dr]_{\varpi^{n+1}} &
\uH^1(\fmF^n \cM) \ar[d]^{\varpi} \ar[rr]^{\gamma^n_\cM} &&
\uH^1(\rami^n\cM) \ar[d]^{\varpi} \ar[drr]^{s} \\
& \uH^1(\fmF^{n+1}\!\cM) \ar[rr]^{\gamma^n_{\fmF\cM}} &&
\uH^1(\fmF\rami^n\cM) \ar[rr]^{\gamma_{\rami^n\cM}} &&
\uH^1(\rami^{n+1}\!\cM)
}
\end{equation*}
%
The left triangle in this diagram commutes by definition of $\varpi$, the square
commutes by naturality of $\gamma$ and the right triangle commutes by definition of
$s$. By assumption $s^n = \gamma^n(\cM) \circ \varpi^n$. However
\begin{equation*}
\gamma_{\rami^n\cM} \circ 
\gamma^{n}_{\fmF\cM} = \gamma^{n+1}_{\cM}
\end{equation*}
by definition of $\gamma$ so the result follows.\quod


\breakflow
Recall that the natural isomorphism $\overline{\rho}_d$ is defined only under
assumption that $\cM/\rami^{d}$ is linear.
We would like to extend it to all elliptic shtukas. To that end we will prove
that the diagram
\begin{equation*}
\xymatrix{
\uH^1(\cM)/\fmF^d \ar[rr]^{\overline{\rho}_d} \ar[d]_{s^d}^{\rtviso}
&&\uH^1(\Der\cM)/\fmF^d \ar[d]^{s^d}_{\ltviso} \\
\uH^1(\rami^d\cM)/\fmF^d \ar[rr]^{\overline{\rho}_d}
&&\uH^1(\Der\rami^d\cM)/\fmF^d
}
\end{equation*}
commutes 
provided the shtuka $\cM/\rami^{2d}$ is linear.
The proof is a bit technical. We split it into a chain of auxillary lemmas.
In the following the integer $d$ and the elliptic shtuka $\cM$ will be fixed.
To improve the legibility we will write $\ibase$ in place of $\rami^d$.

\begin{lem}\label{ellfpregtwistvanish}
The shtuka $(\ibase\cM)/\ibase$ is linear.\end{lem}

\pf Follows instantly from Proposition \ref{elltwistvanish}. \quod

\breakflow
Consider the shtuka $\ibase\cM/\fmF^d$. According to Lemma \ref{ellmfquot}
it is a locally free shtuka on $\cO_F/\fmF^d \otimes \cO_K$ whose restriction to
$\cO_F/\fmF^d \otimes K$ is nilpotent.
In Section \ref{sec:artreg} we equipped the shtukas of this kind with a natural
isomorphism
\begin{equation*}
\rho\colon \uH^1(\ibase\cM/\fmF^d) \to \uH^1(\Der \ibase\cM/\fmF^d)
\end{equation*}
called the artinian regulator.
It is defined only under certain conditions.

\begin{lem}\label{ellfpregartregexist}
If $\cM/\ibase$ is linear then the artinian regulator
is defined for $\ibase\cM/\fmF^d$.
\end{lem}

\pf Indeed $(\ibase\cM)/\ibase$ is linear by Lemma \ref{ellfpregtwistvanish} whence the
result follows from Lemma \ref{ellmfquotcoh} applied to the shtuka $\ibase\cM$. \quod

%

\breakflow
In Section \ref{sec:ellfiltiso} we introduced natural isomorphisms
\begin{equation*}
\uH^1(\fmF^d\cM/\ibase) \xrightarrow{\quad\overline{\gamma}_d\quad}
\uH^1(\ibase\cM/\fmF^d) \xrightarrow{\quad\overline{\varepsilon}_d\quad}
\uH^1(\ibase\cM/\ibase).
\end{equation*}
In the following we drop the indices $d$ for legibility. Our next step is to
study how the artinian regulator $\rho$ of the shtuka $\ibase\cM/\fmF^d$ interacts with
$\overline{\gamma}$ and $\overline{\varepsilon}$.

\begin{lem}\label{ellfpregbeta}
If $\cM/\ibase$ is linear then the diagram
\begin{equation*}
\xymatrix{
\uH^1(\ibase\cM/\fmF^d) \ar[r]^{\rho\,\,\,} \ar[d]_{\overline{\varepsilon}}^{\rtviso} &
\uH^1(\Der \ibase\cM/\fmF^d) \ar[d]^{\overline{\varepsilon}}_{\ltviso} \\
\uH^1(\ibase\cM/\ibase) \ar[r]^{1\,\,\,\,} &
\uH^1(\Der \ibase\cM/\ibase)
}
\end{equation*}
is commutative.
\end{lem}

\pf Consider the diagram
\begin{equation*}
\xymatrix{
\uH^1(\ibase\cM/\fmF^d) \ar[r]^{\rho} \ar[d]_{\textup{red.}}^{\rtviso} &
\uH^1(\Der \ibase\cM/\fmF^d) \ar[d]^{\textup{red.}}_{\ltviso} \\
\uH^1(\ibase\cM/(\fmF^d + \ibase)) \ar[r]^{1\,\,\,\,} &
\uH^1(\Der \ibase\cM/(\fmF^d + \ibase)) \\
\uH^1(\ibase\cM/\ibase) \ar[r]^{1\,\,\,\,} \ar[u]^{\textup{red.}}_{\rtviso} &
\uH^1(\Der \ibase\cM/\ibase) \ar[u]_{\textup{red.}}^{\ltviso}
}
\end{equation*}
The bottom square is clearly commutative.
Applying Proposition \ref{artregred} to $\ibase\cM/\fmF^d$ we conclude that
the top square is commutative.
Now according to Proposition \ref{ellepsilonred} the isomorphism $\overline{\varepsilon}$
is the composition of reduction isomorphisms
\begin{equation*}
\uH^1(\ibase\cM/\fmF^d) \xrightarrow{\isosign}
\uH^1(\ibase\cM/(\fmF^d + \ibase)) \xleftarrow{\isosign}
\uH^1(\ibase\cM/\ibase).
\end{equation*}
Applying the same Proposition to $\Der \ibase\cM/\fmF^d$
we get the result. \quod

\begin{lem}\label{ellfpregalpha}
If $\cM/\ibase^2$ is linear then the diagram
\begin{equation*}
\xymatrix{
\uH^1(\fmF^d\cM/\ibase) \ar[r]^{1\,\,\,\,} \ar[d]_{\overline{\gamma}}^{\rtviso} &
\uH^1(\Der\fmF^d\cM/\ibase) \ar[d]^{\overline{\gamma}}_{\ltviso} \\
\uH^1(\ibase\cM/\fmF^d) \ar[r]^{\rho\,\,\,} &
\uH^1(\Der \ibase\cM/\fmF^d)
}
\end{equation*}
is commutative.
\end{lem}

\pf Consider the boundary homomorphisms
\begin{align*}
\uH^0(\cM/(\fmF^d + \ibase)) &\xrightarrow{\,\,\delta\,\,}
\uH^1(\ibase \cM/\fmF^d), \\
\uH^0(\cM/(\fmF^d + \ibase)) &\xrightarrow{\,\,\delta\,\,}
\uH^1(\fmF^d \cM/\ibase)
\end{align*}
of the exact sequences \eqref{ellcoeffdotseq} and \eqref{ellbasedotseq}.
Using them 
we construct a diagram
\begin{equation}\label{ellfpregalphasplit}
\vcenter{\vbox{\xymatrix{
\uH^1(\fmF^d\cM/\ibase) \ar[r]^{1\,\,\,\,} &
\uH^1(\Der\fmF^d\cM/\ibase) \\
\uH^0(\cM/(\fmF^d + \ibase)) \ar[u]^{\delta}_{\rtviso}
\ar[r]^{1\,\,\,\,} \ar[d]_{\delta}^{\rtviso} &
\uH^0(\Der\cM/(\fmF^d + \ibase)) \ar[u]^{\delta}_{\rtviso}
\ar[d]_{\delta}^{\rtviso} \\
\uH^1(\ibase\cM/\fmF^d) \ar[r]^{\rho\,\,\,} &
\uH^1(\Der \ibase\cM/\fmF^d) 
}}}
\end{equation}
The top square is clearly commutative. We apply Proposition
\ref{artregdelta} to deduce the commutativity of the bottom square.
To use it we need to verify that the following conditions hold for the shtuka
$\cN = \cM/\fmF^d$:
\begin{enumerate}
\item $\ibase \cdot \uH^1(\Der\cN) = 0$.

\item $\cN/\ibase^2$ is linear.
\end{enumerate}
Lemma \ref{ellmfquot} implies (1) while (2) follows since $\cM/\ibase^2$ is linear.
We conclude that \eqref{ellfpregalphasplit} is commutative.
Now Proposition \ref{ellgammadelta} shows that the maps $\delta$ are
isomorphisms and that the composition
\begin{equation*}
\uH^1(\fmF^d\cM/\ibase) \xrightarrow{\,\,\delta^{-1}\,\,}
\uH^0(\cM/(\fmF^d+\ibase)) \xrightarrow{\,\,\delta\,\,}
\uH^1(\ibase\cM/\fmF^d)
\end{equation*}
is equal to $-\overline{\gamma}$. The same Proposition applies to $\Der\cM$ as
well.
So the commutativity of
\eqref{ellfpregalphasplit} implies our result. \quod


\breakflow
We are finally ready to obtain 
the key result of this section.

\begin{prp}\label{ellfpregtautslide}
Let $d \geqslant 1$ and
let $\cM$ be an elliptic shtuka. If $\cM/\rami^{2d}$ is linear then the square
\begin{equation*}
\xymatrix{
\uH^1(\cM)/\fmF^d \ar[rr]^{\overline{\rho}_d} \ar[d]_{s^d}^{\rtviso}
&&\uH^1(\Der\cM)/\fmF^d \ar[d]^{s^d}_{\ltviso} \\
\uH^1(\rami^d\cM)/\fmF^d \ar[rr]^{\overline{\rho}_d}
&&\uH^1(\Der \rami^d\cM)/\fmF^d.
}
\end{equation*}
is commutative.\end{prp}

\pf We continue to use the notation $\ibase$ for $\rami^d$.
By Lemma \ref{ellfpregtwistvanish} the shtuka $(\ibase\cM)/\ibase$ is
linear so that the square makes sense.
We proceed by repeated splitting of this square till the problem is reduced
to its core. 

Using Lemma \ref{ellslidecomp} we split the square as follows:
\begin{equation}\label{ellfpregsplit1}
\vcenter{\vbox{\xymatrix{
\uH^1(\cM)/\fmF^d \ar[rr]^{\overline{\rho}_d} \ar[d]_{\varpi^d}^{\rtviso}
&&\uH^1(\Der\cM)/\fmF^d \ar[d]^{\varpi^d}_{\ltviso} \\
\uH^1(\fmF^d\cM)/\fmF^d \ar[rr]^{\overline{\rho}_d} \ar[d]_{\gamma^d}^{\rtviso}
&&\uH^1(\Der\fmF^d\cM)/\fmF^d \ar[d]^{\gamma^d}_{\ltviso} \\
\uH^1(\ibase\cM)/\fmF^d \ar[rr]^{\overline{\rho}_d}
&&\uH^1(\Der \ibase\cM)/\fmF^d.
}}}
\end{equation}
The top square commutes by functoriality of $\overline{\rho}_d$ so we
concentrate on the bottom square. It is necessary to split this square even
further.

Recall that $\overline{\rho}_d$ is defined as the composition
\begin{equation*}
\xymatrix{
\uH^1(\cM)/\fmF^d \ar[r]^{\textup{red.}} &
\uH^1(\cM/\ibase) \ar[r]^{1\,\,\,\,} &
\uH^1(\Der\cM/\ibase) &
\uH^1(\Der\cM)/\fmF^d \ar[l]_{\textup{red.}}
}
\end{equation*}
So we can rewrite the bottom square as follows:
\begin{equation*}
\xymatrix{
\uH^1(\fmF^d\!\cM)/\fmF^d \ar[r]^{\textup{red.}} \ar[d]_{\gamma^d}^{\rtviso}&
\uH^1(\fmF^d\!\cM/\ibase) \ar[r]^{1\,\,\,\,}
\ar[d]_{(\overline{\varepsilon} \circ \overline{\gamma})_\cM}^{\rtviso} &
\uH^1(\Der\fmF^d\!\cM/\ibase)
\ar[d]^{(\overline{\varepsilon} \circ \overline{\gamma})_{\Der\cM}}_{\ltviso} &
\uH^1(\Der\fmF^d\!\cM)/\fmF^d \ar[l]_{\,\,\textup{red.}}
\ar[d]^{\gamma^d}_{\ltviso}
\\
\uH^1(\ibase\cM)/\fmF^d \ar[r]^{\textup{red.}} &
\uH^1(\ibase\cM/\ibase) \ar[r]^{1\,\,\,\,} &
\uH^1(\Der \ibase\cM/\ibase) &
\uH^1(\Der \ibase\cM)/\fmF^d \ar[l]_{\,\,\textup{red.}}
}
\end{equation*}
Here the notation gets a bit confusing, so let us elaborate on it.
The compositions
$(\overline{\varepsilon}\circ\overline{\gamma})_\cM$ and
$(\overline{\varepsilon}\circ\overline{\gamma})_{\Der\cM}$
have the same source and target. 
So one would expect that the middle square
comutes tautologically. 
However the actual situation is more complicated.
The maps
$(\overline{\varepsilon}\circ\overline{\gamma})_\cM$ and
$(\overline{\varepsilon}\circ\overline{\gamma})_{\Der\cM}$
are defined in terms of data which
comes from completely
different shtukas
$\cM$ and $\Der\cM$.
We add the subscripts $\cM$ and $\Der\cM$ to emphasise this fact.

The left and right squares in the diagram above commute by Lemma
\ref{ellgammaepsilon}.
Let us consider the middle square. We split it in two:
\begin{equation*}
\xymatrix{
\uH^1(\fmF^d\cM/\ibase) \ar[d]_{\overline{\gamma}}^{\rtviso} \ar[r]^{1\,\,\,\,} &
\uH^1(\Der\fmF^d\cM/\ibase) \ar[d]^{\overline{\gamma}}_{\ltviso} \\
\uH^1(\ibase\cM/\fmF^d) \ar[r]^{\rho\,\,\,} \ar[d]_{\overline{\varepsilon}}^{\rtviso} &
\uH^1(\Der \ibase\cM/\fmF^d) \ar[d]^{\overline{\varepsilon}}_{\ltviso}\\
\uH^1(\ibase\cM/\ibase) \ar[r]^{1\,\,\,\,} &
\uH^1(\Der \ibase\cM/\ibase).
}
\end{equation*}
Here $\rho$ is the artinian regulator of the shtuka $\ibase\cM/\fmF^d$ which is
defined by Lemma \ref{ellfpregartregexist}. The bottom square commutes by Lemma
\ref{ellfpregbeta}. Since $\cM/\ibase^2$ is linear the top square commutes by Lemma
\ref{ellfpregalpha}. So we are done. \quod

\section{Regulators of finite order}
\label{sec:ellfpreg}

\begin{dfn}\label{defellfpreg} Let $d \geqslant 1$. 
An $\cO_F/\fmF^d$-linear natural transformation
\begin{equation*}
\rho_d\colon
\uH^1(\cM)/\fmF^d \xrightarrow{\quad}
\uH^1(\Der\cM)/\fmF^d
\end{equation*}
of functors
on the category of elliptic shtukas is called a \emph{regulator of order $d$} if
the following holds:
\begin{enumerate}
\item If $\cM$ is an elliptic shtuka such that $\cM/\rami^{2d}$ is linear then
$\rho_d$ coincides with the tautological regulator $\overline{\rho}_d$ of
Definition \ref{elldeftautreg}.

\item For every elliptic shtuka $\cM$ the natural diagram
\begin{equation*}
\xymatrix{
\uH^1(\cM)/\fmF^d \ar[rr]^{\rho_d} \ar[d]_{s^d}^{\rtviso} &&
\uH^1(\Der\cM)/\fmF^d \ar[d]^{s^d}_{\ltviso} \\
\uH^1(\rami^d\cM)/\fmF^d \ar[rr]^{\rho_d} &&
\uH^1(\Der\rami^d\cM)/\fmF^d
}
\end{equation*}
is commutative.
\end{enumerate}\end{dfn}

\breakflow
The exponent $2d$ in the condition (1) does not look natural. Indeed the
tautological regulator is defined even if $\cM/\rami^d$ is linear. However with
the exponent $d$ the regulators will 
fail to exist.

It is worth noting that the definition of the regulator for the quotient
$\uH^1(\cM)/\fmF^d$ is actually more complicated than the definition for the
module $\uH^1(\cM)$ itself (Definition \ref{defellreg}). In the latter case one
does not need the condition (2).

\begin{prp}\label{ellfpregexist}
Let $d \geqslant 1$.
A regulator of order $d$ exists, 
is unique and is an isomorphism.\end{prp}
%

\pf Let $\cM$ be an elliptic shtuka. According to Proposition
\ref{elltwistvanish} the shtuka $(\rami^{2d}\cM)/\rami^{2d}$ is linear.
In particular we have the tautological regulator $\overline{\rho}_d$ for
$\rami^{2d}\cM$.
We define the regulator $\rho_d$ for $\cM$ by the commutative diagram
\begin{equation*}
\xymatrix{
\uH^1(\cM)/\fmF^d \ar[rr]^{\rho_d} \ar[d]_{s^{2d}}^{\rtviso} &&
\uH^1(\Der\cM)/\fmF^d \ar[d]^{s^{2d}}_{\ltviso} \\
\uH^1(\rami^{2d}\cM)/\fmF^d \ar[rr]^{\overline{\rho}_d} &&
\uH^1(\Der\rami^{2d}\cM)/\fmF^d
}
\end{equation*}
Due to condition (1) of Definition \ref{defellfpreg} this diagram should commute
for any regulator of order $d$. We thus get the unicity of $\rho_d$.

Let us prove that the map $\rho_d$ we just defined is a regulator.
It is a natural $\cO_F/\fmF^d$-linear isomorphism since the maps $s^{2d}$ and
$\overline{\rho}_d$ are so.
If $\cM/\rami^{2d}$ is itself linear then the diagram
\begin{equation*}
\xymatrix{
\uH^1(\cM)/\fmF^d \ar[rr]^{\overline{\rho}_d} \ar[d]_{s^{2d}}^{\rtviso} &&
\uH^1(\Der\cM)/\fmF^d \ar[d]^{s^{2d}}_{\ltviso} \\
\uH^1(\rami^{2d}\cM)/\fmF^d \ar[rr]^{\overline{\rho}_d} &&
\uH^1(\Der\rami^{2d}\cM)/\fmF^d
}
\end{equation*}
commutes by Proposition \ref{ellfpregtautslide}. Hence the condition (1) of
Definition \ref{defellfpreg} is satisfied. To check the condition (2) consider
the diagram
\begin{equation*}
\xymatrix{
\uH^1(\cM)/\fmF^d \ar[rr]^{\rho_d} \ar[d]_{s^d}^{\rtviso} &&
\uH^1(\Der\cM)/\fmF^d \ar[d]^{s^d}_{\ltviso} \\
\uH^1(\rami^d\cM)/\fmF^d \ar[rr]^{\rho_d} \ar[d]_{s^d}^{\rtviso} &&
\uH^1(\Der\rami^d\cM)/\fmF^d \ar[d]^{s^d}_{\ltviso} \\
\uH^1(\rami^{2d}\cM)/\fmF^d \ar[rr]^{\overline{\rho}_d} \ar[d]_{s^d}^{\rtviso} &&
\uH^1(\Der\rami^{2d}\cM)/\fmF^d \ar[d]^{s^d}_{\ltviso} \\
\uH^1(\rami^{3d}\cM)/\fmF^d \ar[rr]^{\overline{\rho}_d} &&
\uH^1(\Der\rami^{3d}\cM)/\fmF^d
}
\end{equation*}
We need to prove that the top square commutes.
The rectangles of height $2$ commute by definition of $\rho_d$. 
The bottom square commutes by Proposition \ref{ellfpregtautslide}. As a
consequence the middle square commutes which implies the commutativity of the
top square. \quod

\begin{prp}\label{ellfpregindep} The regulator $\rho_d$ does not depend on the
choice of the uniformizer $z \in \cO_F$ in the definition of the sliding
isomorphism $s$ (Definition \ref{defellslideiso}).\end{prp}

\pf We will show that $\rho_d$ satisfies the condition (2) of Definition
\ref{defellfpreg} with any choice of $z$.
According to Lemma \ref{ellslidecomp} the sliding isomorphism $s^d$ is the composition
\begin{equation*}
\uH^1(\cM) \xrightarrow{\,\,\varpi^d\,\,}
\uH^1(\fmF^d \cM) \xrightarrow{\,\,\gamma^d\,\,}
\uH^1(\rami^d\cM)
\end{equation*}
Here $\gamma^d$ does not depend on the choice of $z$. The natural isomorpism
$\varpi\colon\cM \to \fmF \cM$ is the unique map whose composition with the
natural embedding $\fmF\cM \hookrightarrow\cM$ is the multiplication by $z$.

Now let $u\in\cO_F^\times$. The regulator $\rho_d$ commutes with
multiplication by $u^d$ since it is $\cO_F/\fm^d$-linear. However
$u\varpi$ is the natural isomorphism $\cM \to \fmF\cM$ with the choice of
uniformizer $u z$. We conclude that $\rho_d$ satisfies the condition (2) of
Definition \ref{defellfpreg} with the uniformizer $u z$ as well. \quod


\breakflow
To construct the regulator isomorphism on the entire module $\uH^1(\cM)$ we would like
to take the limit of regulators $\rho_d$ for $d \to \infty$. To do it we first
need to show that these regulators agree.

\begin{prp}\label{ellfpregcomp} If $k \geqslant d \geqslant 1$ are integers then
the regulators $\rho_d$ and $\rho_k$ coincide modulo $\fmF^d$.\end{prp}

\pf Set $n = 2dk$.
Let $\cM$ be an elliptic shtuka. According to Proposition
\ref{elltwistvanish} the shtuka $(\rami^{n}\cM)/\rami^{n}$ is linear.
The tautological regulators $\overline{\rho}_d$ and $\overline{\rho}_{dk}$ are
defined for $\rami^{n}\cM$.
Since $d$ and $dk$ divide $n$
the diagrams
\begin{equation*}
\xymatrix{
\uH^1(\cM)/\fmF^d \ar[r]^{\rho_d} \ar[d]^{s^{n}}_{\ltviso}
& \uH^1(\Der\cM)/\fmF^d \ar[d]^{s^{n}}_{\ltviso} \\
\uH^1(\rami^{n}\cM)/\fmF^d \ar[r]^{\overline{\rho}_d} &
\uH^1(\Der\rami^{n}\cM)/\fmF^d
}
\quad 
\xymatrix{
\uH^1(\cM)/\fmF^{dk} \ar[r]^{\rho_{dk}} \ar[d]^{s^{n}}_{\ltviso}
& \uH^1(\Der\cM)/\fmF^{dk} \ar[d]^{s^{n}}_{\ltviso} \\
\uH^1(\rami^{n}\cM)/\fmF^{dk} \ar[r]^{\overline{\rho}_{dk}} &
\uH^1(\Der\rami^{n}\cM)/\fmF^{dk}
} 
\end{equation*}
commute by definition of $\rho_d$ and $\rho_{dk}$.
However $\overline{\rho}_d \equiv \overline{\rho}_{dk} \pmod{\fmF^d}$ by
construction. We conclude that $\rho_d \equiv \rho_{dk} \pmod{\fmF^d}$. Applying
the same argument with $d$ and $k$ interchanged we deduce that
$\rho_k \equiv \rho_{dk} \equiv \rho_d \pmod{\fmF^d}.$\quod


\section{Regulators}
\label{sec:ellreg}

%

\begin{dfn}\label{defellreg}\index{idx}{regulator!elliptic shtuka@of an elliptic shtuka}%
An $\cO_F$-linear natural transformation
\begin{equation*}
\uH^1(\cM) \xrightarrow{\quad\quad\rho\quad\quad} \uH^1(\Der\cM)
\end{equation*}
of functors on
the category of elliptic shtukas is called a \emph{regulator} if
for every $\cM$ such that $\cM/\rami^{2n}$ is linear the diagram
\begin{equation*}
\xymatrix{
\uH^1(\cM) \ar[rr]^{\rho} \ar[d]_{\textup{red.}} &&\uH^1(\Der\cM)
\ar[d]^{\textup{red.}} \\
\uH^1(\cM/\rami^n) \ar[rr]^1 && \uH^1(\Der\cM/\rami^n)
}
\end{equation*}
is commutative.\end{dfn}

\begin{lem}\label{ofmodmapzero}%
Let $f\colon M \to N$ be a morphism of
$\cO_F$-modules. If $M$ and $N$ are finitely generated free then the following
are equivalent:
\begin{enumerate}
\item $f = 0$.

\item For every $d > 0$ there exists an $n \geqslant 0$ such that
$f(\fmF^n M) \subset \fmF^{n+d} N$.\quod
\end{enumerate}\end{lem}

%
%

\begin{thm}\label{ellregexist}%
The regulator exists, is unique and is an isomorphism.%
%
\end{thm}

\pf Let $\cM$ be an elliptic shtuka.
By Proposition \ref{elltwistvanish} the shtuka $(\rami^{2d}\cM)/\rami^{2d}$
is linear for every $d \geqslant 1$.
Theorem \ref{ellcohfilt} shows that:
\begin{itemize}
\item $\uH^1(\rami^{2d}\cM) = \fmF^{2d} \uH^1(\cM)$ as submodules of $\uH^1(\cM)$,

\item the kernel of the reduction map
$\uH^1(\Der\rami^{2d}\cM) \to \uH^1(\Der(\rami^{2d}\cM)/\rami^{2d})$ is
the submodule $\fmF^{4d}\uH^1(\Der\cM)$.
\end{itemize}
So the unicity follows from Lemma \ref{ofmodmapzero}.

Now let us construct the regulator and prove that it is an isomorphism. According to
Proposition \ref{ellfpregexist} for every $d \geqslant 1$ there exists a unique
regulator of order $d$
\begin{equation*}
\rho_d\colon \uH^1(\cM)/\fmF^d \to \uH^1(\Der\cM)/\fmF^d.
\end{equation*}
It is a natural $\cO_F/\fmF^d$-linear isomorphism.
The regulators of different orders are compatible by Proposition
\ref{ellfpregcomp} and do not depend on the auxillary choice of a uniformizer $z
\in \cO_F$ by Proposition \ref{ellfpregindep}.
Now take their limit for $d \to \infty$. Since $\uH^1(\cM)$ and $\uH^1(\Der\cM)$
are finitely generated $\cO_F$-modules we get a natural $\cO_F$-linear
isomorphism
$\rho\colon \uH^1(\cM) \to \uH^1(\Der\cM)$.
It satisfies the condition of Definition \ref{defellreg} since every
$\rho_d$ satisfies it up to order $d$.\quod

\begin{thm}\label{ellregcmp}%
Let $\cM$ be an elliptic shtuka and let $f\colon
\uH^1(\cM) \to \uH^1(\Der\cM)$ be an $\cO_F$-linear map.
\begin{enumerate}
\item For every $n \geqslant 0$ the map $f$ sends the submodule
$\uH^1(\rami^n\cM)$ to $\uH^1(\Der\rami^n\cM)$.

%
\item If for every $d > 0$ there exists an $n \geqslant 0$ such that the
shtuka $\rami^n\cM/\rami^{2d}$ is linear and the diagram
\begin{equation*}
\xymatrix{
\uH^1(\rami^n\cM) \ar[rr]^{f} \ar[d]_{\textup{red.}} &&\uH^1(\Der\rami^n\cM) \ar[d]^{\textup{red.}} \\
\uH^1(\rami^n\cM/\rami^d) \ar[rr]^1 && \uH^1(\Der\rami^n\cM/\rami^d)
}
\end{equation*}
is commutative then $f$ coincides with the regulator $\rho$ of $\cM$.%
\end{enumerate}
\end{thm}

\pf Follows from Lemma \ref{ofmodmapzero} in view of Theorem \ref{ellcohfilt}.\quod


\breakflow
The next theorem relates the regulator $\rho$ to the artinian regulator of
Section~\ref{sec:artreg}.
It is essential to the proof of the trace formula in
Chapter~\ref{chapter:trace}.
It can also be used as an alternative definition of the regulator.

\begin{thm}\label{ellregart}%
Let $\cM$ be an elliptic shtuka of ramification ideal $\rami$ and let $n > 0$.
If $\cM/\rami^{2n}$ is linear then the following holds:
\begin{enumerate}
\item The artinian regulator $\rho_{\cM/\fmF^n}$
is defined for the shtuka $\cM/\fmF^n$.

\item The diagram
\begin{equation*}
\xymatrix{
\uH^1(\cM) \ar[rr]^{\rho} \ar[d]_{\textup{red.}} &&\uH^1(\Der\cM)
\ar[d]^{\textup{red.}} \\
\uH^1(\cM/\fm^n) \ar[rr]^{\rho_{\cM/\fm^n}} && \uH^1(\Der\cM/\fm^n)
}
\end{equation*}
is commutative.
\end{enumerate}
\end{thm}

\pf (1) follows from Lemma \ref{ellmfquotcoh} since $\cM/\rami^n$ is linear. Let
us concentrate on (2).
Consider the diagram
\begin{equation*}
\xymatrix{
\uH^1(\cM) \ar[rr]^{\rho} \ar[d]_{\textup{red.}}
&& \uH^1(\Der\cM) \ar[d]^{\textup{red.}}
\\ \uH^1(\cM/\fmF^n) \ar[rr]^{\rho_{\cM/\fmF^d}} \ar[d]_{\textup{red.}}^{\rtviso}
&& \uH^1(\Der\cM/\fmF^n) \ar[d]^{\textup{red.}}_{\ltviso}
\\ \uH^1(\cM/(\fmF^n+\rami^n)) \ar[rr]^{1\,\,\,\,}
&& \uH^1(\Der\cM/(\fmF^n+\rami^n))
}
\end{equation*}
The bottom square in this diagram commutes by Proposition \ref{artregred}.
The vertical arrows in the bottom square are isomorphisms by Proposition
\ref{ellcoeffdotseqprp}. Hence the top square of this diagram commutes if and
only if the outer rectangle commutes. Now we can split the outer rectangle in
two other squares:
\begin{equation*}
\xymatrix{
\uH^1(\cM) \ar[rr]^{\rho} \ar[d]_{\textup{red.}}
&& \uH^1(\Der\cM) \ar[d]^{\textup{red.}}
\\ \uH^1(\cM/\rami^n) \ar[rr]^{1\,\,\,\,} \ar[d]_{\textup{red.}}^{\rtviso}
&& \uH^1(\Der\cM/\rami^n) \ar[d]^{\textup{red.}}_{\ltviso}
\\ \uH^1(\cM/(\fmF^n+\rami^n)) \ar[rr]^{1\,\,\,\,}
&& \uH^1(\Der\cM/(\fmF^n+\rami^n))
}
\end{equation*}
The bottom square commutes tautologically.
As the shtuka $\cM/\rami^{2n}$ is linear the top square 
commutes by definition of $\rho$. Hence we are done. \quod

\chapter{Trace formula}
\label{chapter:trace}
\label{ch:lftwistcoh}
\label{ch:globell}
\label{ch:euprod}

Let $X$ be a smooth projective curve over $\Fq$. As in Chapter \ref{chapter:shtcoh}
we fix an open dense affine subscheme $\Spec R \subset X$. Its complement consists
finitely many closed points. 
We denote $K$
the product of the local fields of $X$ at these points, 
$\cO_K \subset K$ the ring of integers and $\fm_K \subset \cO_K$
the Jacobson radical.

%
Fix a local field $\ngF$ containing $\Fq$. Let $\ngOF\subset \ngF$ be the ring of
integers and $\fm_F \subset\ngOF$ the maximal ideal.
We set $\nsOF = \Spec\ngOF$.
In this chapter we mainly work with the $\tau$-scheme $\nsOFX$.
As usual the endomorphism $\tau\colon\nsOFX\to\nsOFX$
acts as the identity on the left hand side
of the product $\times$ and as the $q$-Frobenius on the right hand side.
The same applies to the schemes
of the form $\Spec\ngC\times X$ for an $\Fq$-algebra $\ngC$
and to tensor product rings.

Let $\rami\subset\fm_K$ be an open ideal.
We say that a shtuka $\cM$ on $\nsOFX$ is \emph{elliptic of ramification ideal $\rami$}
if it has the following properties:
\begin{enumerate}
\item $\cM$ is locally free,

\item $\cM(\ngOF/\fm_F \otimes R)$ is nilpotent,

\item $\cM(\ngOF\complot\cO_K)$ is an elliptic shtuka of ramification ideal $\rami$
in the sense of Definition \lref{reg}{defell}.\end{enumerate}
Using the theory of Chapter \ref{chapter:reg} we will construct
for every elliptic shtuka $\cM$ 
a natural quasi-isomorphism
\begin{equation*}
\rho_\cM\colon \RGamma(\cM) \xrightarrow{\,\,\isosign\,\,} \RGamma(\Der\cM),
\end{equation*}
the \emph{regulator} of $\cM$.
We will also define a numerical invariant $L(\cM) \in \ngOF^\times$.
This invariant is given by an infinite product of local factors
\begin{equation*}
\prod_\fm L(\cM(\ngOF\otimes R/\fm))^{-1}
\end{equation*}
where $\fm$ runs over the maximal ideals of $R$.

The main result of this
chapter is Theorem  \ref{elltrace}
which states that
\begin{equation*}
\zeta_\cM = L(\cM) \cdot \det\nolimits_{\ngOF}(\rho_\cM)
\end{equation*}
under a certain technical condition on $\cM$.
Here $\zeta_\cM$ is the $\zeta$-isomorphism of $\cM$ in the sense of Definition
\lref{shtukacoh}{defzeta}.
We call this expression the trace formula for regulators of elliptic shtukas.
%
Theorem \ref{elltrace} is 
based on
Anderson's trace
formula \cite{trace} 
in the form 
given to it by B\"ockle-Pink \cite{bp} and V.~Lafforgue \cite{valeurs}.
The statement of
Theorem \ref{elltrace}
has its roots in the article \cite{valeurs}
of V.~Lafforgue as well.


\section{Preliminaries}

In this section we prove an auxillary statement on cohomology of coherent 
sheaves over $\ngS \times X$ where $\ngS$ is the spectrum of a local noetherian
ring $\ngC$. We denote $\fm_\ngC \subset \ngC$ the maximal ideal.

\begin{lem}\label{lfcohchar}%
Let $\cF$ be a coherent sheaf on $\ngS\times X$.
If $\cF$ is $\Lambda$-flat then the following are equivalent:
\begin{enumerate}
\item $\uH^0(\Spec(\ngC/\fm_\ngC)\times X,\,\cF) = 0$.

\item $\uH^0(\cF) = 0$ and $\uH^1(\cF)$ is a free $\ngC$-module of finite rank.%
\end{enumerate}\end{lem}

\pf The base change theorem for coherent cohomology [\stacks{07VK}] states that
$\RGamma(\cF)$ is a perfect complex of $\ngC$-modules and
\begin{equation*}
\RGamma(\cF) \otimes_{\ngC}^{\mathbf L} \ngC/\fm_\ngC =
\RGamma(\Spec(\ngC/\fm_\ngC)\times X, \,\cF).
\end{equation*}
Hence (2) $\Rightarrow$ (1). Let us prove the other direction.
By base change we know that $\RGamma(\cF)\otimes_{\ngC}^{\mathbf L} \ngC/\fm_\ngC$
is concentrated in degree $1$. As $\RGamma(\cF)$ is a perfect complex it follows
[\stacks{068V}] that $\RGamma(\cF)$ has Tor amplitude $[1,1]$.
Now [\stacks{0658}] shows that $\RGamma(\cF)$ is quasi-\hspace{0pt}isomorphic
to a finitely generated free $\Lambda$-module placed in degree $1$.\quod

\section{Euler products in the artinian case}

We work with a finite $\Fq$-algebra $\ngC$ which is a local artinian ring.
%
As before we denote $\fm_\ngC \subset \ngC$ the maximal ideal.

\begin{lem}\label{euartfieldnilp}
Let $k$ be a finite field extension of $\Fq$. Let $\cM$ be a locally free shtuka
on $\ngC\otimes k$ given by a diagram
\begin{equation*}
M_0 \shtuka{i}{j} M_1.
\end{equation*}
If $\cM(\ngC/\fm_\ngC \otimes k)$ is nilpotent then the following holds:
\begin{enumerate}
\item $M_0$ is a free $\ngC$-module of finite rank,

\item $i\colon M_0 \to M_1$ is an isomorphism,

\item $i^{-1} j \colon M_0 \to M_0$ is a $\ngC$-linear nilpotent endomorphism.
\end{enumerate}\end{lem}

\pf (1) is clear. 
The ring $\ngC \otimes k$ is noetherian and complete with respect to the
nilpotent ideal $\fm_\ngC \otimes R$ so Proposition \lref{nilp}{nilpcomp}
implies that $\cM$ is nilpotent. (2) and (3) follow by definition of a
nilpotent shtuka. \quod

\begin{dfn}\label{defeuart}
Let $k$ be a finite field extension of $\Fq$. Let $\cM$ be a locally free shtuka
on $\ngC \otimes k$ given by a diagram
\begin{equation*}
M_0 \shtuka{i}{j} M_1.
\end{equation*}
Assuming that $\cM(\ngC/\fm_\ngC \otimes k)$ is nilpotent we define
\begin{equation*}
L(\cM) =
\det\nolimits_\ngC(1 - i^{-1} j \mid M_0) \in \ngC^\times.
\end{equation*}
\end{dfn}

\begin{lem}\label{euartnilp}
Let $\cM$ be a locally free shtuka on $\ngC \otimes R$. If $\cM(\ngC/\fm_\ngC \otimes
R)$ is nilpotent then $\cM$ itself is nilpotent.\end{lem}

\pf The ring $\ngC \otimes R$ is noetherian and complete with respect to the
nilpotent ideal $\fm_\ngC \otimes R$. Whence Proposition \lref{nilp}{nilpcomp}
implies that $\cM$ is nilpotent. \quod


\begin{lem}\label{euartok}%
Let $\cM$ be a locally free shtuka on $\ngC \otimes R$.
If $\cM(\ngC/\fm_\ngC \otimes R)$ is nilpotent then 
%
for almost all maximal ideals
$\fm \subset R$ we have $L(\cM(\ngC \otimes R/\fm)) = 1$.
\end{lem}

\pf
Suppose that $\cM$ is given by a diagram
\begin{equation*}
M_0 \shtuka{i}{j} M_1.
\end{equation*}
Let $\fm \subset R$ be a maximal ideal.
According to \cite[Lemma 8.1.3]{bp} we have
\begin{equation*}
\det\nolimits_\ngC (1 - i^{-1} j \mid M_0/\fm) = 
\det\nolimits_{\ngC \otimes R/\fm}\big(1 - (i^{-1} j)^d \mid M_0/\fm\big)
\end{equation*}
where $d$ is the degree of the finite field $R/\fm$ over $\Fq$.
However $i^{-1} j\colon M_0 \to M_0$ is a nilpotent endomorphism
by Lemma \ref{euartnilp}.
Hence $(i^{-1} j)^d = 0$ for $d \gg 0$ and we get the result. \quod

\begin{dfn}Let $\cM$ be a locally free shtuka on $\ngC\otimes R$.
Assuming that $\cM(\ngC/\fm_\ngC \otimes R)$ is nilpotent we define
\begin{equation*}
L(\cM) = \prod_\fm L(\cM(\ngC \otimes R/\fm))^{-1} \in \ngC^\times
\end{equation*}
where $\fm \subset R$ ranges over the maximal ideals.
This product is well-defined by Lemma \ref{euartok}.%
\end{dfn}

\begin{lem}\label{euartone}%
Let $\cM$ be a locally free shtuka on $\ngC\otimes R$ such that
$\cM(\ngC/\fm_\ngC \otimes R)$ is nilpotent. If $\ngC$ is a field then $L(\cM) =
1$.\end{lem}

\pf If $V$ is a vector space over a field and $N$ a nilpotent
endomorphism of $V$ then $\det(1 - N\,|\,V) = 1$. Hence for every maximal
ideal $\fm \subset R$ the invariant $L(\cM(\ngC \otimes R/\fm))$ is equal to $1$.\quod

\breakflow
Let $\ngS = \Spec\ngC$
and let $\cM$ be a locally free shtuka on $\ngS\times X$ such that
$\cM(\ngC/\fm_\ngC\otimes R)$ is nilpotent. 
We write $L(\cM)$ in place of $L(\cM(\ngC\otimes R))$ to make
the formulas more legible.
Note that the closed points of $X$ in the complement of $\Spec R$
are not taken into account.

\section{Anderson's trace formula}

We continue working with a finite $\Fq$-algebra $\ngC$ which is a local artinian
ring. 
As usual $\fm_\ngC\subset\ngC$ stands for the maximal ideal
and $\ngS$ denotes the spectrum of $\ngC$.
%
%
Our trace formula for the shtuka-theoretic regulator is based on
the following narrow variant of Anderson's trace formula:

\begin{thm}\label{nilptrace}%
Let $\cM$ be a locally free shtuka on $\ngS\times X$ given by a diagram
\begin{equation*}
\cE \shtuka{1}{j} \cE.
\end{equation*}
Suppose that
\begin{enumerate}
%
\item $\uH^0(\Spec(\ngC/\fm_\ngC)\times X,\,\cE) = 0$,

\item $\cM$ is nilpotent,

\item $\cM(\ngC\otimes \cO_K/\fm_K)$ is linear.
\end{enumerate}
Then 
%
$\uH^1(\cE)$ is a free $\ngC$-module of finite rank
and
$L(\cM) = \det\nolimits_\ngC(1 - j \mid \uH^1(\cE))$.%
\end{thm}

\breakflow
We will deduce it from Proposition 3.2 in the article
\cite{valeurs} of V.~Lafforgue.
However our setting differs slightly from Lafforgue's.
In \cite[Section 3]{valeurs}
the coefficient ring (denoted $A$) is a power series ring
while we work with a local artinian ring $\ngC$.
The next lemma helps to bridge this gap.
In the following $\cI \subset \cO_{\ngS\times X}$ stands for the unique invertible ideal
sheaf such that
\begin{equation*}
\cI(\ngC\otimes R) = \ngC\otimes R,\quad
\cI(\ngC\otimes \cO_K) = \ngC\otimes\fm_K.
\end{equation*}

\begin{lem}\label{lftwistcohgood}%
Let $\cE$ be a locally free sheaf on $\ngS\times X$.
%
If $\uH^0(\Spec(\ngC/\fm_\ngC)\times X,\,\cE) = 0$
then for every $n \geqslant 0$ the following holds:
\begin{enumerate}
%
\item $\uH^0(\cI^n\cE) = 0$ and $\uH^1(\cI^n\cE)$ is a free $\ngC$-module of finite
rank.

\item The natural map $\uH^1(\cI^n\cE) \to \uH^1(\cE)$ is a split surjection.%
\end{enumerate}\end{lem}

\pf 
Consider the short exact sequence of coherent sheaves
$0 \to \cI^n \cE \to \cE \to \cE/\cI^n \cE \to 0$.
The assumption 
implies that $\uH^0(\Spec(\ngC/\fm_\ngC)\times X,\,\cI^n\cE) = 0$
so (1) follows from Lemma \ref{lfcohchar}.
Since $\cE/\cI^n\cE$ is supported at a closed affine subscheme of $\ngS\times X$
the long exact cohomology sequence implies
that the maps $\uH^1(\cI^n\cE) \to \uH^1(\cE)$ are onto.
They split since $\uH^1(\cE)$ is free.\quod

\breakflow
We also need a lemma on shtukas:

\begin{lem}\label{nilptwistdet}%
Let $\cM$ be a locally free shtuka on $\ngS\times X$ given by a diagram
\begin{equation*}
\cE \shtuka{1}{j} \cE.
\end{equation*}
Suppose that
\begin{enumerate}
%
\item $\uH^0(\Spec(\ngC/\fm_\ngC)\times X,\,\cE) = 0$,

\item $\cM(\ngC\otimes \cO_K/\fm_K)$ is linear.
\end{enumerate}
Then for every $n \geqslant 0$ the following holds:
\begin{enumerate}
\item $\uH^1(\cI^n\cE)$ is a free $\ngC$-module of finite rank,

\item $\det_\ngC(1 - j \mid \uH^1(\cI^n\cE)) = \det_\ngC(1 - j\mid\uH^1(\cE))$.
\end{enumerate}\end{lem}

\pf Part (1) is immediate from Lemma \ref{lftwistcohgood}.
The same lemma implies that the short exact sequence
of coherent sheaves
$0 \to \cI^{n+1} \cE \to \cI^n \cE \to (\cI^n\cE)/\cI \to 0$
induces a short exact sequence of cohomology modules
\begin{equation*}
0 \to \uH^0(\cI^n\cE/\cI) \to \uH^1(\cI^{n+1} \cE) \to \uH^1(\cI^n\cE) \to 0
\end{equation*}
with $\uH^0(\cI^n\cE/\cI)$ a free $\ngC$-module of finite rank. As a consequence
\begin{equation*}
\det\nolimits_\ngC\!\big(1 - j \mid \uH^1(\cI^{n+1} \cE\big) =
\det\nolimits_\ngC\!\big(1 - j \mid \uH^0(\cI^n\cE/\cI)\big) \cdot
\det\nolimits_\ngC\!\big(1 - j \mid \uH^1(\cI^{n} \cE)\big).
\end{equation*}
By assumption $j(\cE) \subset \cI\cE$. Since $\tau(\cI) \subset \cI^q$ and
$j$ is $\tau$-linear we conclude that $j(\cI^n\cE) \subset \cI^{qn + 1}\cE$.
Hence $j$ is zero on the quotient $\cI^n\cE/\cI$ and we get the result.\quod

\afterall\noindent
\textit{Proof of Theorem \ref{nilptrace}}.
Let $\Omega_{\ngS\times X}$ be the canonical sheaf of $\ngS\times X$ over $\ngS$
and let
\begin{align*}
\cV &= \uH^0(\Spec\ngC\otimes R, \,\iHom(\cE,\,\Omega_{\ngS\times X})), \\
\cV_t &= \uH^0(\ngS\times X, \,\iHom(\cI^t\cE,\,\Omega_{\ngS\times X})).
\end{align*}
Grothendieck-Serre duality identifies $\cV_t$ with the $(-1)$-st cohomology
module of the complex
\begin{equation*}
\RHom_\ngC(\RGamma(\cI^t\cE),\,\ngC).
\end{equation*}
However Lemma \ref{lftwistcohgood} shows that $\uH^0(\cI^t\cE) = 0$ and
$\uH^1(\cI^t\cE)$ is a free $\ngC$-module of finite rank. Hence
Grothendieck-Serre duality gives a natural isomorphism
\begin{equation*}
\Hom_\ngC(\uH^1(\cI^t\cE),\,\ngC) \xrightarrow{\isosign} \cV_t.
\end{equation*}
In particular every $\cV_t$ is a free $\ngC$-module of finite rank.
By Lemma \ref{lftwistcohgood} the natural maps
$\uH^1(\cI^{t+1}\cE)  \to \uH^1(\cI^t\cE)$
are split surjections whence the
inclusions $\cV_t \subset \cV_{t+1}$ are split.

For every $t \geqslant 0$
the endomorphism $j$ of $\cE$ induces 
an endomorphism of $\uH^1(\cI^t \cE)$ and 
Grothendieck-Serre duality identifies it
with the Cartier-linear endomorphisms $\kappa_\cV\colon \cV_t \to \cV_t$
which are used in Section 3 of \cite{valeurs}.

As we explained above the $\ngC$-modules $\cV_t$ are free of finite rank and the
inclusions $\cV_t \subset \cV_{t+1}$ are split. So the argument 
of \cite[Proposition 3.2]{valeurs} applies with only one minor change.
Namely, one
needs to ensure that the auxillary locally free sheaves $\cF$ on $\ngS\times X$
constructed in the course of the proof have the property that $\uH^0(\cF) = 0$
and $\uH^1(\cF)$ is a free $\ngC$-module of finite rank. This can be done
as follows. The dual of the ideal sheaf $\cI$ is ample relative to $\ngS$.
Hence $\uH^0(\Spec(\ngC/\fm_\ngC)\times X, \,\cI^n\cF)$ is zero for $n \gg 0$.
Lemma \ref{lfcohchar} implies that $\cI^n\cF$ has the desired property.

The rest of the argument works without change. It shows that for $t \gg 0$ we
have an equality of power series in
$\ngC[[T]]$:
\begin{equation*}
\det\nolimits_\ngC(1 - T \kappa_\cV \mid \cV_t) = 
\prod_\fm
\det\nolimits_\ngC\!\big(1 - T j \mid \cE(\ngC\otimes R/\fm)\big)^{-1}
\end{equation*}
where $\fm \subset R$ runs over the maximal ideals.

Now recall the endomorphism $j\colon\cE \to \tau_\ast\cE$ is assumed to be
nilpotent. Furthermore the maximal ideal of $\ngC$ is nilpotent too.
As a consequence \cite[Lemma 8.1.4]{bp} implies that only finitely
many factors in the product on the right hand side are different from $1$.
Thus we can evaluate the equality above at $T = 1$ and conclude that
\begin{equation*}
\det\nolimits_\ngC(1 - \kappa_\cV \mid \cV_t) =
\prod_\fm
\det\nolimits_\ngC\!\big(1 - j \mid \cE(\ngC\otimes R/\fm)\big)^{-1} = L(\cM).
\end{equation*}
Therefore
$\det\nolimits_\ngC(1 - j \mid \uH^1(\cI^t \cE)) = L(\cM)$.
Lemma \ref{nilptwistdet} shows that this equality holds already for $t = 0$.\quod

\section{Artinian regulators}
\label{sec:globartreg}

We keep the assumption that $\ngC$ is a finite $\Fq$-algebra which is local
artinian.

\begin{lem}%
Let $\cM$ be a locally free shtuka on $\ngS\times X$ given by a diagram
\begin{equation*}
\cM_0 \shtuka{i}{j} \cM_1.
\end{equation*}
If $\cM(\ngC/\fm_\ngC \otimes R)$ is nilpotent then
the restriction of $i$ to $\ngC\otimes R$
is an isomorphism.\end{lem}

\pf Indeed $\cM(\ngC \otimes R)$ is nilpotent by Lemma \ref{euartnilp}, and the
$i$-arrow of such a shtuka is an isomorphism by definition. \quod

\begin{lem}\label{globartregij}%
Let $\cM$ be a locally free shtuka on $\ngS\times X$ given by a diagram
\begin{equation*}
\cM_0 \shtuka{i}{j} \cM_1.
\end{equation*}
Suppose that the shtuka $\cM(\ngC/\fm_\ngC\otimes R)$ is nilpotent.
If the endomorphism $i^{-1} j$ of $\cM_0(\Spec\ngC\otimes K)$
preserves the submodule $\cM_0(\Spec\ngC\otimes\cO_K)$
then the following holds:
\begin{enumerate}
\item The $\tau$-linear endomorphism $i^{-1} j$ of
$\cM_0|_{\Spec\ngC\otimes R}$ extends to a unique
$\tau$-linear endomorphism of $\cM_0$.

\item The endomorphism $i^{-1}j$ of $\cM_0$ is nilpotent.\end{enumerate}%
\end{lem}

\pf (1)
%
According to Lemma \ref{cechsquare} the fibre product of the schemes
$\Spec(\ngC\otimes R)$ and $\Spec(\ngC\otimes\cO_K)$ over $\ngS\times X$
is $\Spec(\ngC\otimes K)$. Furthermore the closed subscheme
$\Spec(\ngC\otimes \cO_K/\fm_K)$ of $\ngS\times X$ is defined locally
by a principal ideal sheaf.
Hence Beauville-\hspace{0pt}Laszlo glueing theorem [\stacks{0BP2}]
implies that there exists a unique morphism
$(i^{-1}j)^a\colon \tau^\ast\cM_0 \to \cM_0$
restricting to the adjoints
\begin{align*}
(i^{-1}j)^a &\colon \tau^\ast \cM_0|_{\Spec\ngC \otimes R} \to \cM_0|_{\Spec\ngC\otimes R},\\
(i^{-1}j)^a &\colon \tau^\ast \cM_0|_{\Spec\ngC \otimes \cO_K} \to \cM_0|_{\Spec\ngC\otimes\cO_K}
\end{align*}
of the endomorphisms $i^{-1} j$ on $\ngC \otimes R$ respectively $\ngC\otimes\cO_K$.
(2) follows since $i^{-1} j$ is nilpotent on $\ngC\otimes R$.\quod

\breakflow
Observe that $\tau$ acts as the identity on the underlying topological
space of $\ngS\times X$.
We are thus in position to apply the constructions of Section~\ref{sec:globcoh}:
given a shtuka
\begin{equation*}
\cM = \Big[ \cM_0\shtuka{\,\,i\,\,}{j}\cM_1 \Big]
\end{equation*}
on the scheme $\ngS\times X$ we have the associated complex
of $\Lambda$-module sheaves
\begin{equation*}
\cG_a(\cM) = \big[ \cM_0 \xrightarrow{\,\,i-j\,\,} \cM_1 \big]
\end{equation*}
on the underlying topological space of $\ngS\times X$.

\begin{dfn}\label{defglobartreg}%
Let $\cM$ be a locally free shtuka on $\ngS\times X$ given by a diagram
\begin{equation*}
\cM_0 \shtuka{i}{j} \cM_1.
\end{equation*}
Suppose that $\cM(\ngC/\fm_\ngC \otimes R)$ is nilpotent. We say that the artinian
regulator is defined for $\cM$ if the $\tau$-linear endomorphism $i^{-1} j$ of
$\cM_0(\Spec\ngC\otimes K)$ preserves the submodule $\cM_0(\Spec\ngC\otimes\cO_K)$.
In this case we define the \emph{artinian regulator}
$\rho_\cM\colon \cG_a(\cM) \to \cG_a(\Der\cM)$
by the diagram
\begin{equation*}
\xymatrix{
\big[\cM_0 \ar[r]^{i-j} \ar[d]_{1-i^{-1} j} & \cM_1\big] \ar[d]^1 \\
\big[\cM_0 \ar[r]^i & \cM_1\big]
}
\end{equation*}
Here $i^{-1} j\colon \cM_0 \to \cM_0$ is the $\tau$-linear endomorphism
constructed in Lemma~\ref{globartregij}.%
\end{dfn}

\breakflow
Theorem~\ref{shtglobartcoh} identifies $\RGamma(\ngS\times X,\,\cM)$ and
$\RGamma(\ngS\times X,\,\cG_a(\cM))$. So the artinian regulator $\rho_\cM$
induces a morphism $\RGamma(\cM) \to \RGamma(\Der\cM)$
on shtuka cohomology. We also denote it $\rho_\cM$.

\begin{lem}\label{globartregcmp}%
Let $\cM$ be a locally free shtuka on $\ngS\times X$ such that
$\cM(\ngC/\fm_\ngC\otimes R)$ is nilpotent. If the artinian regulator
is defined for $\cM$ then it is defined for $\cM(\ngC\otimes\cO_K)$
and the diagram
\begin{equation*}
\xymatrix{
\RGamma(\cM) \ar[rr]^{\rho_\cM} \ar[d] && \RGamma(\Der\cM) \ar[d] \\
\RGamma(\ngC\otimes\cO_K,\,\cM) \ar[rr]^{\rho_{\cM(\ngC\otimes\cO_K)}} &&
\RGamma(\ngC\otimes\cO_K,\,\Der\cM)
}
\end{equation*}
is commutative.%
\end{lem}

\pf Suppose that $\cM$ is given by a diagram
\begin{equation*}
\cM_0 \shtuka{i}{j} \cM_1.
\end{equation*}
By assumption the endomorphism $i^{-1}j$ of $\cM_0(\Spec\ngC\otimes K)$
preserves the submodule $\cM_0(\Spec\ngC\otimes\cO_K)$.
Hence the artinian regulator in the sense of Definition \ref{defartreg}
is defined for $\cM(\ngC\otimes\cO_K)$.

Let $U$ be an affine open neighbourhood of the closed
subscheme $\Spec\cO_K/\fm_K$ inside $X$.
Applying the same construction as for the shtukas on $\ngS\times X$
we get an artinian regulator $\rho_U$ for the restriction of $\cM$ to $\ngS\times U$.
We now claim that the squares in the diagram
\begin{equation}\label{globartregcmpdiag}
\vcenter{\vbox{\xymatrix{
\RGamma(\ngS\times X,\,\cM) \ar[rr]^{\rho_\cM} \ar[d] && \RGamma(\ngS\times X,\,\Der\cM) \ar[d] \\
\RGamma(\ngS\times U,\,\cM) \ar[rr]^{\rho_U} \ar[d] && \RGamma(\ngS\times U,\,\Der\cM) \ar[d] \\
\RGamma(\ngC\otimes\cO_K,\,\cM) \ar[rr]^{\rho_{\cM(\ngC\otimes\cO_K)}} &&
\RGamma(\ngC\otimes\cO_K,\,\Der\cM)
}}}
\end{equation}
are commutative.
For the top square it follows from Theorem \ref{shtglobartcoh}
and the definition of $\rho$, $\rho_U$.
For the bottom square we argue as follows.

The higher cohomology of quasi-\hspace{0pt}coherent sheaves on the affine
scheme $\ngS\times U$ is zero. As a consequence
\begin{align*}
\RGamma(\ngS\times U,\,\cG_a(\cM)) &= \Gamma(\ngS\times U,\,\cG_a(\cM)), \\
\RGamma(\ngS\times U,\,\cG_a(\Der\cM)) &= \Gamma(\ngS\times U,\,\cG_a(\Der\cM)).
\end{align*}
The complexes on the right hand side are the associated complexes
of $\cM$ respectively $\Der\cM$.
Denoting $M_0 = \cM_0(\ngS\times U)$ and $M_1 = \cM_1(\ngS\times U)$
we conclude that $\rho_U$ comes from
a morphism of the associated complexes
given by the diagram
\begin{equation*}
\xymatrix{
[ M_0 \ar[r]^{i-j} \ar[d]_{1 - i^{-1} j} & M_1 \ar[d]^1 ] \\
[ M_0 \ar[r]^{i} & M_1 ].
}
\end{equation*}
Therefore the bottom square of \eqref{globartregcmpdiag}
commutes by definition of artinian regulators
for shtukas on $\ngC\otimes\cO_K$ (Definition \ref{defartreg}).\quod

\begin{lem}\label{globartreg}%
Let $\cM$ be a locally free shtuka on $\ngS\times X$ given by the diagram
\begin{equation*}
\cM_0 \shtuka{i}{j} \cM_1.
\end{equation*}
Suppose that $\cM(\ngC/\fm_\ngC \otimes R)$ is nilpotent.
If the artinian regulator is defined for $\cM$ then the diagram
\begin{equation*}
\xymatrix{
\RGamma(\cM) \ar[r] \ar[d]_{\rho_\cM} &
\RGamma(\cM_0) \ar[r]^{i-j} \ar[d]_{1 - i^{-1} j} &
\RGamma(\cM_1) \ar[r] \ar[d]^1 & [1] \\
\RGamma(\Der\cM) \ar[r] &
\RGamma(\cM_0) \ar[r]^i &
\RGamma(\cM_1) \ar[r] & [1]
}
\end{equation*}
is a morphism of canonical triangles.%
\end{lem}

\pf Consider a diagram of sheaf complexes
\begin{equation}\label{globartregmordiag}
\vcenter{\vbox{\xymatrix{
\cG_a(\cM) \ar[r] \ar[d]_{\rho_\cM} & \cM_0[0] \ar[r]^{i-j} \ar[d]_{1-i^{-1}j} & \cM_1[0] \ar[r] \ar[d]^1 & [1] \\
\cG_a(\Der\cM) \ar[r] & \cM_0[0] \ar[r]^i & \cM_1[0] \ar[r] & [1]
}}}
\end{equation}
The rows are the distinguished triangles of the mapping fiber
complexes $\cG_a(\cM)$ respectively $\cG_a(\Der\cM)$. By construction
\eqref{globartregmordiag} is a morphism of distinguished triangles.
Now Theorem \ref{shtglobartcoh} shows that applying the sheaf cohomology
functor $\RGamma(\ngS\times X, -)$ to \eqref{globartregmordiag} we get
exactly the diagram in the statement of this lemma.\quod

\section{Trace formula for artinian regulators}
\label{sec:arttrace}

We continue working with a finite $\Fq$-algebra $\ngC$ which is a local artinian ring.
As before $\fm_\ngC\subset\ngC$ denotes the maximal ideal and
$\ngS$ stands for the spectrum of $\ngC$.

\begin{prp}\label{globartcoh}%
Let $\cM$ be a locally free shtuka on $\ngS\times X$.
If $\cM(\ngC/\fm_\ngC \otimes R)$ is nilpotent then the following holds:
\begin{enumerate}
\item The natural map
$\RGamma(\cM) \to \RGamma(\ngC \otimes \cO_K,\,\cM)$
is a quasi-isomorphism.

\item The complex $\RGamma(\cM)$ is concentrated in degree $1$.

\item $\uH^1(\cM)$ is a free $\ngC$-module of finite rank.
\end{enumerate}\end{prp}

\pf (1) Consider the \v{C}ech complex
\begin{equation*}
\RvGamma(\cM) =
\Big[
\RGamma(\ngC \otimes R,\,\cM) \oplus
\RGamma(\ngC \otimes \cO_K,\,\cM) \to
\RGamma(\ngC \otimes K,\,\cM) \Big].
\end{equation*}
According to
Theorem \lref{cechcoh}{cechcoh} there is a natural quasi-isomorphism
$\RGamma(\cM) \xrightarrow{\isosign} \RvGamma(\cM)$.
Its composition with the natural map
$\RvGamma(\cM) \to \RGamma(\ngC\otimes \cO_K,\,\cM)$ is the pullback map
$\RGamma(\cM) \to \RGamma(\ngC\otimes\cO_K,\,\cM)$. So to prove that this pullback
map is a quasi-isomorphism it is enough to demonstrate that
\begin{equation*}
\RGamma(\ngC\otimes R,\,\cM) = 0, \quad
\RGamma(\ngC\otimes K,\,\cM) = 0.
\end{equation*}
Now the ring $\ngC \otimes R$ is noetherian and complete with respect to the ideal
$\fm_\ngC \otimes R$. As the shtuka $\cM(\ngC/\fm_\ngC \otimes R)$ is nilpotent
Proposition \lref{nilp}{nilpcomp} implies that $\RGamma(\ngC \otimes R,\,\cM) = 0$.
The shtuka $\cM(\ngC/\fm_\ngC \otimes K)$ is nilpotent since nilpotence is preserved under
pullbacks.
So the same argument shows that $\RGamma(\ngC \otimes K,\,\cM) = 0$. Whence the result.
In view of (1) the results (2) and (3) follow from Theorem \lref{artnilp}{artnilpcoh}.\quod


\breakflow
Let $\cM$ be a locally free shtuka on $\ngS\times X$
such that
$\cM(\ngC/\fm_\ngC \otimes R)$ is nilpotent.
Our goal is to compare the artinian regulator $\rho_\cM$ with the
$\zeta$-isomorphism of~$\cM$. This isomorphism is
constructed in the following way
(see Definition \lref{shtukacoh}{defzeta}).
Suppose that 
$\cM$ is given by a diagram
\begin{equation*}
\cM_0 \shtuka{i}{j} \cM_1.
\end{equation*}
We make the additional assumption that
the $\ngC$-modules
$\uH^n(\cM_0)$,
$\uH^n(\cM_1)$,
$\uH^n(\cM)$,
$\uH^n(\Der\cM)$,
$n \geqslant 0$,
are finitely generated free.
The $\zeta$-isomorphism
\begin{equation*}
\zeta_\cM\colon
\det\nolimits_\ngC\RGamma(\cM) \xrightarrow{\,\isosign\,}
\det\nolimits_\ngC\RGamma(\Der\cM)
\end{equation*}
is then defined as the composition of the isomorphisms
\begin{equation*}
\det\nolimits_\ngC^{\phantom{1}}\RGamma(\cM) \xrightarrow{\isosign}
\det\nolimits_\ngC^{\phantom{1}}\RGamma(\cM_0) \otimes_\ngC^{\phantom{1}}
\det\nolimits_\ngC^{-1}\RGamma(\cM_1) \xleftarrow{\isosign}
\det\nolimits_\ngC^{\phantom{1}}\RGamma(\Der\cM)
\end{equation*}
induced by the canonical triangles
\begin{align*}
&\RGamma(\cM) \to \RGamma(\cM_0) \xrightarrow{\,\,i-j\,\,}
\RGamma(\cM_1) \to [1], \\
&\RGamma(\Der\cM) \to \RGamma(\cM_0) \xrightarrow{\,\,i\,\,}
\RGamma(\cM_1) \to [1].
\end{align*}
 
The freeness hypothesis on the
cohomology 
is necessary to make the construction
of $\zeta_\cM$ work.
The reason is that
in the theory of Knudsen-Mumford \cite[Corollary 2 after Theorem 2]{kmdet}
a distinguished triangle
$A \to B \to C \to [1]$
of perfect $\ngC$-module complexes
determines a canonical isomorphisms
\begin{equation*}
\det\nolimits_\ngC^{\phantom{1}} A \xrightarrow{\,\isosign\,}
\det\nolimits_\ngC^{\phantom{1}} B \otimes_\ngC^{\phantom{1}}
\det\nolimits_\ngC^{-1} C
\end{equation*}
only if the cohomology modules of the
complexes $A$, $B$ and $C$ are themselves perfect. 
A module over the local artinian ring $\ngC$ is perfect if and only if it
is free of finite rank while the cohomology modules of a perfect $\ngC$-module
complex can be arbitrary finitely generated $\ngC$-modules.


\begin{thm}\label{arttrace}%
Let $\cM$ be a locally free shtuka on $\ngS\times X$ given by a diagram
\begin{equation*}
\cM_0 \shtuka{i}{j} \cM_1.
\end{equation*}
Assume that
\begin{enumerate}
\item $\uH^0(\Spec(\ngC/\fm_\ngC)\times X,\,\cM_*) = 0$ for $*\in\{0,1\}$,

\item $\cM(\ngC/\fm_\ngC \otimes R)$ is nilpotent,

\item There exists an open ideal $\ibase \subset \fm_K$ such that
$\ibase \cdot \uH^1(\ngC \otimes \cO_K, \Der\cM) = 0$
and $\cM(\ngC \otimes \cO_K/\ibase^2)$ is linear.
\end{enumerate}
Then the following holds:
\begin{enumerate}
\item The artinian regulator $\rho_\cM$ is defined for $\cM$.

\item The $\zeta$-isomorphism $\zeta_\cM$ is defined for $\cM$.

\item $\zeta_\cM = L(\cM) \cdot \det_\ngC(\rho_\cM)$.
\end{enumerate}
\end{thm}

\pf (1) According to Proposition \lref{artreg}{artregexist}
the assumptions (2) and (3) imply that the artinian regulator
is defined for $\cM(\ngC\otimes\cO_K)$ or equivalently that
the endomorphism $i^{-1} j$ of $\cM_0(\Spec\Lambda\otimes K)$
preserves the submodule $\cM_0(\Spec\Lambda\otimes\cO_K)$.
Hence the artinian regulator is defined for $\cM$.

(2) %
%
By assumption (1) the module $\uH^0(\Spec(\ngC/\fm_\ngC)\times X,\,\cM_0)$ is zero.
Hence Lemma \ref{lfcohchar} shows that
$\uH^0(\cM_0) = 0$ and $\uH^1(\cM_0)$ is a free $\ngC$-module
of finite rank. The same applies to $\cM_1$.
Since $\cM(\ngC/\fm_\ngC\otimes R)$ is nilpotent by assumption (2)
Proposition \ref{globartcoh} demonstrates that the cohomology modules of $\RGamma(\cM)$
and $\RGamma(\Der\cM)$ are also perfect. 

(3)
Lemma \ref{globartreg} implies that
$\det\nolimits_\ngC(\rho_\cM) = 
\det\nolimits_\ngC \big(1 - i^{-1} j \,\,\big|\,\, \RGamma(\cM_0)\big) \cdot
\zeta_\cM$.
However $\uH^0(\cM_0) = 0$ and $\uH^1(\cM_0)$ is a free $\ngC$-module.
Therefore
\begin{equation*}
\zeta_\cM =
\det\nolimits_\ngC \big(1 - i^{-1} j \,\,\big|\,\, \uH^1(\cM_0)\big) \cdot
\det\nolimits_\ngC(\rho_\cM).
\end{equation*}
We thus need to compute
$\det\nolimits_\ngC \big(1 - i^{-1} j \,\,\big|\,\, \uH^1(\cM_0)\big)$.
Consider a locally free shtuka
\begin{equation*}
\cN = \Big[ \cM_0 \shtuka{1}{\,\,i^{-1} j\,\,} \cM_0 \Big].
\end{equation*}
This shtuka has the following properties:
\begin{enumerate}
\item[(i)] $\uH^0(\Spec(\ngC/\fm_\ngC)\times X, \,\cM_0) = 0$,

\item[(ii)] $\cN$ is nilpotent,

\item[(iii)] $\cN(\ngC \otimes \cO_K/\fm_K)$ is linear.
\end{enumerate}
Indeed (i) is just the assumption (1) and (ii) holds
by Lemma \ref{globartregij}. Let us prove that (iii) also holds.
We temporarily denote
\begin{equation*}
M_0 = \cM_0(\Spec\ngC\otimes\cO_K),\quad
M_1 = \cM_1(\Spec\ngC\otimes\cO_K).
\end{equation*}
According to the assumption (3) we have
$\ibase M_1 \subset i(M_0)$.
Hence
$\ibase^2 M_1 \subset i(\ibase M_0)$.
At the same time the assumption (3) implies that
$j(M_0) \subset \ibase^2 M_1$.
As a consequence $j(M_0) \subset i(\ibase M_0)$.
We conclude that
the shtuka $\cN(\ngC\otimes\cO_K/\ibase)$ is linear.
However $\ibase \subset \fm_K$ by assumption so we get
the property (iii).

We now apply Theorem \ref{nilptrace} to $\cN$ and conclude that
\begin{equation*}
\det\nolimits_\ngC \big(1 - i^{-1} j \,\,\big|\,\, \uH^1(\cM_0)\big) =
L(\cN).
\end{equation*}
To get the result it remains to observe that
$L(\cN) = L(\cM)$ by construction.\quod 

%
%
%
%

\section{A lemma on determinants}
\label{sec:detlemma}

Our aim for the moment is to prove a technical statement on determinants
of finite-dimensional $F$-vector spaces.
In this section we omit the subscript $F$ of $\det$ and $\otimes$ in order to
simlify the notation.

In the theory of Knudsen-Mumford \cite{kmdet}
a distinguished triangle
$A \to B \to C \to [1]$ 
of bounded complexes with perfect cohomology modules 
determines a natural isomorphism
$\det A \xrightarrow{\isosign} \det B \otimes \det\nolimits^{-1} C$
\cite[Corollary 2 after Theorem 2]{kmdet}. In particular a short exact
sequence 
of finite-dimensional $F$-vector spaces
gives rise to such a natural isomorphism.

\begin{lem}\label{detlemma}%
Consider a diagram of finite-dimensional $F$-vector spaces:
\begin{equation*}
\xymatrix{
& & B \ar@{=}[r]\ar[d] & B\ar[d] & \\
0 \ar[r] & A \ar@{=}[d] \ar[r]^p & B \oplus A_0 \ar[r]^q \ar[d] & B \oplus A_1 \ar[r]\ar[d] & 0 \\
0 \ar[r] & A \ar[r]^f & A_0 \ar[r]^g & A_1 \ar[r] & 0 \\
}
\end{equation*}
Assume the following:
\begin{enumerate}
\item The rows are short exact sequences.

\item The vertical arrows in the second and third column are natural
inclusions respectively projections.

\item The diagram is commutative.
\end{enumerate}
Then the natural square of determinants
\begin{equation*}
\xymatrix{
\det A \ar[r]^{\tisosign\quad\quad\quad\,\,\quad\quad} \ar@{=}[d] &
\det(B \oplus A_0) \otimes \det\nolimits^{-1}(B \oplus A_1) \ar[d]^{\rtviso} \\
\det A \ar[r]^{\tisosign\quad\quad\quad} & \det A_0 \otimes \det\nolimits^{-1} A_1
}
\end{equation*}
commutes up to the sign $(-1)^{nm}$ where $n = \dim A$, $m = \dim B$.
The right vertical arrow in this square is induced by the
natural isomorphism $\det B \otimes \det^{-1} B \xrightarrow{\isosign} F$.
\end{lem}

\breakflow
Morally this lemma is a special case of \cite[Proposition 1 (ii)]{kmdet}. There is,
however, a gap between morality and reality which one has to fill
with a sound argument.

\afterall\noindent
\textit{Proof of Lemma \ref{detlemma}.}
Let $s\colon A_1 \to A_0$ be a section of the surjection $g\colon A_0 \to A_1$.
Pick an element $a_1 \in \det A_1$.
It is easy to show that the natural isomorphism
$\det A \xrightarrow{\isosign} \det A_0 \otimes \det^{-1} A_1$
sends $a \in \det A$ to the element
\begin{equation}\label{detlemmaf1}
\big(f(a) \wedge s(a_1)\big) \otimes a_1^*
\end{equation}
where $a_1^*\colon \det A_1 \to F$ is the unique linear map such that
$a_1^*(a_1) = 1$. This element depends neither on the choice of $s$ nor on the
choice of $a_1$. We would like to obtain a similar explicit formula for the map
determined by the short exact sequence
\begin{equation*}
0 \to A \xrightarrow{\,\,p\,\,} B \oplus A_0 \xrightarrow{\,\,q\,\,} B \oplus A_1
\to 0.
\end{equation*}
The commutativity of the diagram implies that
the map $q$
is given by a matrix
\begin{equation*}
\begin{pmatrix}
1 & h \\
0 & g
\end{pmatrix}
\end{equation*}
where $h\colon A_0 \to B$ is a certain map. Similarly
\begin{equation*}
p = \begin{pmatrix}
\delta \\ f
\end{pmatrix}
\end{equation*}
where $\delta\colon A \to B$ is a certain map.

Let $s\colon A_1 \to A_0$ be a section of $g$. A quick computation shows that
the map
\begin{equation*}
t = \begin{pmatrix}
1 & - hs \\ 0 & s
\end{pmatrix}\colon B \oplus A_1 \to B \oplus A_0
\end{equation*}
is a section of $q$. Next, fix elements $b \in \det B$ and $a_1 \in\det A_1$.
Let $m = \dim B$ and $d_1 = \dim A_1$. By a slight abuse of notation
we write $b \wedge a_1 \in \Lambda^{m + d_1}(B \oplus A_1)$ for the element of
$\det(B \oplus A_1)$ defined by $b$, $a_1$.
For every $a \in \det A$
we need to compute
\begin{equation*}
p(a) \wedge t(b \wedge a_1).
\end{equation*}
Suppose that
\begin{equation*}
a = a^1 \wedge \dotsc \wedge a^n, \quad
b = b^1 \wedge \dotsc \wedge b^m, \quad
a_1 = a_1^1 \wedge \dotsc \wedge a_1^{d_1}.
\end{equation*}
In this case
\begin{equation*}
t(b \wedge a_1) =
b_1 \wedge \dotsc \wedge b^m \wedge \big(s(a_1^1) - h s(a_1^1)\big) \wedge
\dotsc \wedge \big(s(a_1^{d_1}) - h s(a_1^{d_1})\big) = b \wedge s(a_1).
\end{equation*}
Furthermore
\begin{equation*}
p(a) = \big(\delta(a^1) + f(a^1)\big) \wedge \dotsc\wedge
\big(\delta(a^n) + f(a^n)\big)
\end{equation*}
so that
\begin{equation*}
p(a) \wedge t(b \wedge a_1) =
f(a) \wedge b \wedge s(a_1) =
(-1)^{nm} b \wedge f(a) \wedge s(a_1).
\end{equation*}
We conclude that the natural isomorphism
$\det A \xrightarrow{\isosign} \det(B \oplus A_0) \otimes
\det^{-1}(B \oplus A_1)$
sends the element $a \in \det A$ to
\begin{equation*}
(-1)^{nm} b \otimes \big(f(a) \wedge s(a_1)\big) \otimes b^* \otimes a_1^*
\end{equation*}
where $b^*\colon \det B \to F$ and $a_1^*\colon \det A_1 \to F$ are the unique
elements such that $b^*(b) = 1$, $a_1^*(a_1) = 1$.
Comparing this formula with \eqref{detlemmaf1} we get the result. \quod

\section{\texorpdfstring{A functoriality statement for $\zeta$-isomorphisms}{A functoriality statement for zeta-isomorphisms}}


We work with a coefficient ring $\ngF$ which is a local field.
As before we denote $\nsF = \Spec\ngF$.

\begin{lem}%
Let $\cM$ be a shtuka on $\nsF\times X$.
If $\cM$ is coherent then the $\zeta$-isomorphism is defined for
$\cM$.\end{lem}

\pf The ring $\ngF$ is regular so the result is a consequence of
Proposition~\lref{globzeta}{regzetadefined}.\quod 

%
 
\begin{lem}\label{zetadecomp}%
Let $0 \to \cN \to \cM \to \cQ \to 0$ be a short exact sequence
of coherent shtukas on $\nsF\times X$.
Suppose that $\cN$, $\cM$ and $\cQ$ are given by diagrams
\begin{equation*}
\cN = \Big[\cN_0 \shtuka{i}{j} \cN_1\Big], \quad
\cM = \Big[\cM_0 \shtuka{i}{j} \cM_1\Big], \quad
\cQ = \Big[\cQ_0 \shtuka{i}{j} \cQ_1\Big].
\end{equation*}
Assume the following:
\begin{enumerate}
\item $\RGamma(\cM)$ and $\RGamma(\Der\cM)$ are concentrated in degree $1$.

\item $\uH^0(\cM_*) = 0$ for $* \in \{0,1\}$.

\item $\RGamma(\cQ) = 0$ and $\RGamma(\Der\cQ) = 0$.

\item $\uH^1(\cQ_*) = 0$ for $* \in \{0,1\}$.

\item $\cQ$ is linear.
%
\end{enumerate}
Then the following holds:
\begin{enumerate}
\item The natural map $\RGamma(\cN) \to \RGamma(\cM)$ is a quasi-isomorphism.

\item The natural map $\RGamma(\Der\cN) \to \RGamma(\Der\cM)$ is a
quasi-isomorphism.

\item The natural square
\begin{equation*}
\xymatrix{
\det\nolimits_{\ngF}^{\phantom{1}} \RGamma(\cN) \ar[d]_{\zeta_{\cN}}^{\rtviso}
\ar[r]^{\isosign} &
\det\nolimits_{\ngF}^{\phantom{1}} \RGamma(\cM)
\ar[d]^{\zeta_{\cM}}_{\ltviso} \\
\det\nolimits_{\ngF}^{\phantom{1}} \RGamma(\Der\cN)
\ar[r]^{\isosign} &
\det\nolimits_{\ngF}^{\phantom{1}} \RGamma(\Der\cM)
}
\end{equation*}
commutes up to the sign $(-1)^{nm}$ where
\begin{align*}
n &= \dim_{\ngF}\uH^1(\Der\cM) + \dim_F \uH^1(\cM), \\
m &= \dim_{\ngF}\uH^0(\cQ_0).
\end{align*}%
\end{enumerate}%
\end{lem}



\begin{rmk}%
Lemma \ref{zetadecomp} should certainly
hold without the conditions (1), (2) and (4).
However to establish it in this generality one should
prove a stronger version of Lemma \ref{detlemma}. Unfortunately such a proof
appears to be quite messy. We opt not to do it since the present
version of Lemma \ref{zetadecomp} is all we need to prove the main result
of this text, the class number formula.\end{rmk}

\afterall\noindent
\textit{Proof of Lemma \ref{zetadecomp}.}
In the following
we drop the subscript $\ngF$ of $\det$ and $\otimes$
to improve the legibility.
The claims
(1) and (2) are immediate consequences of the assumption (3). Let us now prove (3).
Assumptions (1) and (2) imply that we have a short exact sequence
\begin{equation*}
0 \to \uH^1(\cM) \to \uH^1(\cM_0) \xrightarrow{\,\,i-j\,\,} \uH^1(\cM_1) \to 0.
\end{equation*}
Assumption (2) also implies that $\uH^0(\cN_*) = 0$ for $* \in\{0,1\}$. 
The natural map $\RGamma(\cN) \to \RGamma(\cM)$ is a quasi-isomorphism due to
assumption (3).
We thus have a short exact sequence
\begin{equation*}
0 \to \uH^1(\cN) \to \uH^1(\cN_0) \xrightarrow{\,\,i-j\,\,} \uH^1(\cN_1) \to 0.
\end{equation*}
As $\uH^0(\cN_*) = 0$ the assumption (4) implies that the cohomology
sequence of the short exact sequence $0 \to \cN_* \to \cM_* \to \cQ_* \to 0$ has
the form
\begin{equation*}
0 \to \uH^0(\cQ_*) \xrightarrow{\,\,\delta\,\,} \uH^1(\cN_*) \to \uH^1(\cM_*)
\to 0.
\end{equation*}
Altogether we have a commutative diagram
\begin{equation*}
\xymatrix{
& & 0 \ar[d] & 0 \ar[d] & \\
& 0 \ar[r]\ar[d] & \uH^0(\cQ_0)\ar[d]_{\delta} \ar[r]^i_{\bisosign} &
\uH^0(\cQ_1) \ar[r] \ar[d]^{\delta} & 0 \\
0 \ar[r] & \uH^1(\cN) \ar[r] \ar[d]^{\rtviso} & \uH^1(\cN_0) \ar[d] \ar[r]^{i-j} &
\uH^1(\cN_1) \ar[r] \ar[d] & 0 \\
0 \ar[r] & \uH^1(\cM) \ar[d] \ar[r] & \uH^1(\cM_0) \ar[d] \ar[r]^{i-j} &
\uH^1(\cM_1) \ar[d] \ar[r] & 0 \\
& 0 & 0 & 0 &
}
\end{equation*}
The arrow $\uH^0(\cQ_0) \to \uH^0(\cQ_1)$ is labelled $i$ since the shtuka
$\cQ$ is linear by assumption (5). Let us denote $B = \uH^0(\cQ_0)$.
Taking splittings of the maps $\delta$ we
replace the diagram above with an isomorphic diagram
\begin{equation}\label{zetafuncdiag}
\vcenter{\vbox{\xymatrix{
& & B \ar@{=}[r] \ar[d] & B \ar[d] & \\
0 \ar[r] & \uH^1(\cN) \ar[r] \ar[d]^{\rtviso} & B \oplus \uH^1(\cM_0) \ar[d] \ar[r]^{i-j} &
B \oplus \uH^1(\cM_1) \ar[r] \ar[d] & 0 \\
0 \ar[r] & \uH^1(\cM) \ar[r] & \uH^1(\cM_0) \ar[r]^{i-j} &
\uH^1(\cM_1) \ar[r] & 0
}}}
\end{equation}
The same argument applied to the shtukas $\Der\cN$, $\Der\cM$, $\Der\cQ$
shows that we have a commutative diagram
\begin{equation}\label{zetafuncdiaglie}
\vcenter{\vbox{\xymatrix{
& & B \ar@{=}[r] \ar[d] & B \ar[d] & \\
0 \ar[r] & \uH^1(\Der\cN) \ar[r] \ar[d]^{\rtviso} & B \oplus \uH^1(\cM_0) \ar[d] \ar[r]^{i} &
B \oplus \uH^1(\cM_1) \ar[r] \ar[d] & 0 \\
0 \ar[r] & \uH^1(\Der\cM) \ar[r] & \uH^1(\cM_0) \ar[r]^{i} &
\uH^1(\cM_1) \ar[r] & 0
}}}
\end{equation}
Observe that the vertical arrows in the second and third column of
\eqref{zetafuncdiaglie} can be chosen to be the same as the corresponding arrows
in \eqref{zetafuncdiag}. Indeed they only depend on the chosen splittings of the
injections $\uH^0(\cQ_0) \to \uH^1(\cN_0)$ and $\uH^0(\cQ_1) \to \uH^1(\cN_1)$
and on the map $i\colon \uH^0(\cQ_0) \to \uH^0(\cQ_1)$. Here we again use the
assumption (5) that $\cQ$ is linear.

Applying Lemma \ref{detlemma} to \eqref{zetafuncdiag} we conclude that
the natural square of determinants
\begin{equation}\label{zetadetdiag}
\vcenter{\vbox{\xymatrix{
\det \uH^1(\cN) \ar[d]_{\ltviso}
\ar[r]^{\tisosign\quad\quad\quad\quad\quad\quad\quad} &
\det\big(B \oplus \uH^1(\cM_0)\big) \otimes
\det\nolimits^{-1}\big(B \oplus \uH^1(\cM_1)\big) \ar[d]^{\rtviso} \\
\det \uH^1(\cM) \ar[r]^{\tisosign\quad\quad\quad\quad} &
\det \uH^1(\cM_0) \otimes
\det\nolimits^{-1} \uH^1(\cM_1)
}}}
\end{equation}
commutes up to the sign $(-1)^{n_1 m}$ where $n_1 = \dim \uH^1(\cM)$,
$m = \dim B$. Lemma \ref{detlemma} applied to \eqref{zetafuncdiaglie} shows
that the natural square of determinants
\begin{equation}\label{zetadetdiaglie}
\vcenter{\vbox{\xymatrix{
\det \uH^1(\Der\cN) \ar[d]_{\ltviso}
\ar[r]^{\tisosign\quad\quad\quad\quad\quad\quad\quad} &
\det\big(B \oplus \uH^1(\cM_0)\big) \otimes
\det\nolimits^{-1}\big(B \oplus \uH^1(\cM_1)\big) \ar[d]^{\rtviso} \\
\det \uH^1(\Der\cM) \ar[r]^{\tisosign\quad\quad\quad\quad} &
\det \uH^1(\cM_0) \otimes
\det\nolimits^{-1} \uH^1(\cM_1)
}}}
\end{equation}
commutes up to $(-1)^{n_2 m}$ with $n_2 = \dim \uH^1(\Der\cM)$.

Now the composition of the bottom horizontal isomorphisms in
\eqref{zetadetdiag} and \eqref{zetadetdiaglie} is the $\zeta$-isomorphism of
$\cM$ by definition while the composition of the top horizontal isomorphisms is
the $\zeta$-isomorphism of $\cN$ by construction of the diagrams
\eqref{zetafuncdiag} and \eqref{zetafuncdiaglie}. As the right vertical
isomorphisms in \eqref{zetadetdiag} and \eqref{zetadetdiaglie} are the same we
get the result. \quod

\section{Elliptic shtukas}

Starting from this section we focus on the coefficient ring $\ngOF$ which is the
ring of integers of a local field $\ngF$. As before we denote $\nsOF = \Spec\ngOF$
and $\nsF = \Spec\ngF$.


\begin{dfn}\label{defglobell}\index{idx}{shtuka!elliptic}\index{idx}{ramification ideal!elliptic shtuka@of an elliptic shtuka}\index{nidx}{$\rami$, ramification ideal}%
Let $\cM$ be a shtuka on $\nsOF\times X$ and let $\rami\subset\cO_K$ be an open ideal.
We say that $\cM$ is an \emph{elliptic shtuka of ramification ideal $\rami$} if the
following holds:
\begin{enumerate}
\item $\cM$ is locally free.

\item $\cM(\ngOF/\fm_\ngF \otimes R)$ is nilpotent.

\item $\cM(\ngOF \complot\cO_K)$ is an elliptic shtuka of ramification ideal $\rami$
in the sense of Definition \lref{reg}{defell}.
\end{enumerate}\end{dfn}

\begin{rmk}We fix the ramification ideal $\rami$ throughout the rest of the chapter.  In
the following we speak simply of elliptic shtukas rather than elliptic shtukas
of ramification ideal $\rami$.\end{rmk}

\begin{lem} If a shtuka $\cM$ is elliptic
then $\Der\cM$ is elliptic.\end{lem}

\pf Follows immediately from Proposition \lref{reg}{ellder}. \quod


\begin{prp}\label{globellcoh}%
For every elliptic shtuka $\cM$ on $\nsOF\times X$ the following holds:
\begin{enumerate}
\item The natural map $\RGamma(\cM) \to \RGamma(\ngOF\complot\cO_K,\,\cM)$
is a quasi-\hspace{0pt}isomorphism. 

\item $\RGamma(\cM)$ is concentrated in degree $1$.

\item $\uH^1(\cM)$ is a free $\ngOF$-module of finite rank.%
\end{enumerate}\end{prp}

\pf Indeed $\cM(\ngOF/\fm_\ngF \otimes R)$ 
is nilpotent so Theorem \lref{complcech}{intconc} shows that
the natural map $\RGamma(\cM) \to \RGamma(\ngOF\complot\cO_K,\,\cM)$
is a quasi-isomorphism. Now $\cM(\ngOF\complot\cO_K)$ is an elliptic shtuka
of ramification ideal $\rami$ by definition. So Theorem \lref{reg}{ellcoh} shows that
the complex $\RGamma(\ngOF\complot\cO_K,\,\cM)$ is concentrated in degree $1$
and $\uH^1(\ngOF\complot\cO_K,\,\cM)$ is a free $\ngOF$-module of finite rank. \quod

\breakflow
As in the case of elliptic shtukas on $\ngOF\complot\cO_K$ we will need a
twisting construction for elliptic shtukas on $\nsOF\times X$.

\begin{dfn}\label{defglobelltwist}%
Let $\cM$ be a quasi-coherent shtuka on $\nsOF\times X$ given by a
diagram
\begin{equation*}
\cM_0 \shtuka{\quad i\quad}{j} \cM_1.
\end{equation*}
We define the shtuka $\rami\cM$ by the
diagram
\begin{equation*}
\cI\cM_0 \shtuka{\quad i\quad }{j} \cI\cM_1
\end{equation*}
where $\cI \subset \cO_{\ngS\times X}$ is the unique ideal sheaf such that
\begin{equation*}
\cI(\Spec\ngOF\otimes R) = \ngOF\otimes R, \quad
\cI(\Spec\ngOF\otimes\cO_K) = \ngOF\otimes \rami.
\end{equation*}%
The shtuka $\rami\cM$ will be called the \emph{twist} of $\cM$.

Let $n \geqslant 1$.
By construction we have a natural embedding $\rami^n\cM \hookrightarrow \cM$.
We denote $\cM/\rami^n$ the quotient $\cM/\rami^n\cM$. Observe that
$\cM/\rami^n$ coincides with the restriction of $\cM$ to the closed
affine subscheme $\Spec(\cO_F\otimes\cO_K/\rami^n)$ of $\nsOFX$.%
\end{dfn}

\begin{lem}\label{globelltwist}%
Let $\cM$ be a shtuka over $\nsOF\times X$. If $\cM$ is elliptic then
$\rami\cM$ is elliptic.\end{lem}

\pf 
The shtuka $\rami\cM$ is locally free since the ideal sheaf $\cI$ above is invertible.
Furthermore $\cM(\ngOF\otimes R) = (\rami\cM)(\ngOF\otimes R)$ so that
$(\rami\cM)(\ngOF/\fm_\ngF\otimes R)$ is nilpotent. Finally
$(\rami\cM)(\ngOF\complot\cO_K)$ is the twist of $\cM(\ngOF\complot\cO_K)$ by
$\rami$ in the sense of Definition \lref{reg}{defelltwist}. Hence
Proposition \lref{reg}{elltwist} shows that $(\rami\cM)(\ngOF\complot\cO_K)$
is an elliptic shtuka in the sense of Definition \lref{reg}{defell}. \quod

\begin{dfn}\label{defglobellreg}\index{idx}{regulator!elliptic shtuka@of an elliptic shtuka}%
An $\cO_F$-linear natural transformation
\begin{equation*}
\RGamma(\cM) \xrightarrow{\quad\quad\rho\quad\quad} \RGamma(\Der\cM)
\end{equation*}
of functors on
the category of elliptic shtukas is called a \emph{regulator} if
for every $\cM$ such that $\cM/\rami^{2n}$ is linear the diagram
\begin{equation*}
\xymatrix{
\RGamma(\cM) \ar[rr]^{\rho} \ar[d] &&\RGamma(\Der\cM)
\ar[d] \\
\RGamma(\cM/\rami^n) \ar[rr]^1 && \RGamma(\Der\cM/\rami^n)
}
\end{equation*}
is commutative.%
\end{dfn}


\begin{thm}\label{globellregexist}%
There exists a unique regulator. Moreover:
\begin{enumerate}
%
\item The regulator is a quasi-\hspace{0pt}isomorphism.

\item For every $\cM$ the square
\begin{equation}\label{globellregsquare}
\vcenter{\vbox{\xymatrix{
\RGamma(\cM) \ar[rr]^{\rho} \ar[d] && \RGamma(\Der\cM) \ar[d] \\
\RGamma(\cO_F\complot\cO_K,\,\cM) \ar[rr]^{\rho|_{\cO_F\complot\cO_K}} &&
\RGamma(\cO_F\complot\cO_K,\,\Der\cM)
}}}
\end{equation}
is commutative. Here the bottom arrow is the regulator of the elliptic shtuka
$\cM(\cO_F\complot\cO_K)$.\end{enumerate}\end{thm}

\pf
Proposition \ref{globellcoh} shows that the vertical arrows
in the diagram \eqref{globellregsquare} are quasi-\hspace{0pt}isomorphisms.
This fact has two consequences. First we can define a natural transformation
$\rho\colon\RGamma(\cM)\to\RGamma(\Der\cM)$ using \eqref{globellregsquare}.
The definition of a regulator for shtukas on $\cO_F\complot\cO_K$ then implies
that $\rho$ is indeed a regulator in the sense of Definition \ref{defglobellreg}.

Second, any regulator $\rho\colon\RGamma(\cM) \to \RGamma(\Der\cM)$ induces
an $\cO_F$-linear morphism
$\RGamma(\cO_F\complot\cO_K,\,\cM) \to \RGamma(\cO_F\complot\cO_K,\,\Der\cM)$.
Let us denote it $\rho'$.
Observe that $\rami\cM(\cO_F\complot\cO_K)$ is the twist of
$\cM(\cO_F\complot\cO_K)$ by $\rami$ in the sense of
Definition \ref{defelltwist}. So Proposition \ref{elltwistvanish} implies
that $(\rami^{2n}\cM)/\rami^{2n}$ is linear for all $n \geqslant 0$.
As a consequence the diagram
\begin{equation*}
\xymatrix{
\RGamma(\cO_F\complot\cO_K,\,\rami^{2n}\cM) \ar[rr]^{\rho'} \ar[d] &&
\RGamma(\cO_F\complot\cO_K,\,\rami^{2n}\Der\cM) \ar[d] \\
\RGamma(\cO_F\complot\cO_K,\,(\rami^{2n}\cM)/\rami^n) \ar[rr]^{1} &&
\RGamma(\cO_F\complot\cO_K,\,(\rami^{2n}\Der\cM)/\rami^n)
}
\end{equation*}
is commutative for every $n \geqslant 0$.
Theorem \ref{ellregcmp} now implies that $\rho'$ coincides with the regulator
of the elliptic shtuka $\cM(\cO_F\complot\cO_K)$ so that we get the property (2).
Now the unicity and (1) follow by Theorem \ref{ellregexist}.\quod


\begin{dfn}%
Let $\icoef \subset \ngOF$ be a nonzero ideal.
\begin{enumerate}
\item For a sheaf of modules $\cE$ on $\nsOF\times X$
we denote $\cE/\icoef$ the restriction of $\cE$
to the closed subscheme $\Spec(\ngOF/\icoef)\times X$.

\item For a shtuka $\cM$ on $\nsOFX$ we denote
$\cM/\icoef$ its restriction to the closed subscheme $\Spec(\ngOF/\icoef) \times X$.
\end{enumerate}\end{dfn}


\begin{lem}%
Let $\icoef \subset\ngOF$ be a nonzero ideal. If $\cM$ is a locally free shtuka on
$\nsOF\times X$ then the pullback map
$\RGamma(\cM)\otimes_{\ngOF}^{\mathbf L} \ngOF/\icoef \to \RGamma(\cM/\icoef)$
is a quasi-\hspace{0pt}isomorphism.\end{lem}

\pf It is an immediate consequence of Proposition
\lref{globcoeffch}{coeffchangetriiso}.\quod

\begin{lem}\label{globellregred}%
Let $\icoef = \fm_\ngF^d \subset \ngOF$ be a nonzero ideal.
Let $\cM$ be an elliptic shtuka.
If $\cM(\ngOF\otimes\cO_K/\rami^{2d})$ is linear then the following holds:
\begin{enumerate}
\item The artinian regulator $\rho_{\cM/\icoef}$ is defined for $\cM/\icoef$.

\item The natural diagram
\begin{equation*}
\xymatrix{
\RGamma(\cM)\otimes_{\ngOF}^{\mathbf L}\ngOF/\icoef \ar[rr]^{\rho_\cM} \ar[d]_{\ltviso} &&
\RGamma(\Der\cM) \otimes_{\ngOF}^{\mathbf L}\ngOF/\icoef \ar[d]^{\rtviso} \\
\RGamma(\cM/\icoef) \ar[rr]^{\rho_{\cM/\icoef}} &&
\RGamma(\Der\cM/\icoef)
}
\end{equation*}
is commutative.\end{enumerate}\end{lem}

\pf (1) Suppose that $\cM$ is given by a diagram
\begin{equation*}
\cM_0 \shtuka{i}{j}\cM_1.
\end{equation*}
As $\cM(\ngOF\otimes\cO_K/\rami^{2d})$ is linear
Theorem \lref{reg}{ellregart} shows that the artinian regulator is defined for
$\cM(\ngOF/\icoef \otimes \cO_K)$. In other words
the endomorphism $i^{-1}j$ of $\cM_0(\Spec\ngOF/\icoef\otimes K)$ preserves
the submodule $\cM_0(\Spec\ngOF/\icoef\otimes\cO_K)$.
So the artinian regulator is defined for $\cM/\icoef$.

(2) We need to prove that the square
\begin{equation*}
\xymatrix{
\RGamma(\cM) \ar[rr]^{\rho_\cM} \ar[d] &&
\RGamma(\Der\cM) \ar[d] \\
\RGamma(\cM/\icoef) \ar[rr]^{\rho_{\cM/\icoef}} &&
\RGamma(\Der\cM/\icoef)
}
\end{equation*}
is commutative. Theorem \ref{globellregexist} in combination with
Proposition \ref{globellcoh} identifies the top arrow
with the regulator of $\cM(\cO_F\complot\cO_K)$.
Similarly Lemma \ref{globartregcmp} and Proposition \ref{globartcoh} identify the bottom arrow with
the artinian regulator of $\cM(\cO_F/\icoef\otimes\cO_K)$.
So the result is a consequence of Theorem \ref{ellregart}.\quod

\breakflow
Our goal is to compare the determinant of the regulator with the
$\zeta$-isomorphism.

\begin{lem}If $\cM$ is a locally free shtuka on $\nsOF\times X$ then the
$\zeta$-isomorphism is defined for $\cM$.\end{lem}

\pf Follows from Proposition \lref{globzeta}{regzetadefined} since
$\ngOF$ is regular and $\cM$ is coherent. \quod


%

\begin{lem}\label{globellzetared}%
Let $\icoef \subset \ngOF$ be a nonzero ideal.
Let $\cM$ be a locally free shtuka on $\nsOFX$ given by a diagram
$\cM_0 \shtuka{}{} \cM_1$.
Assume that 
\begin{enumerate}
\item $\cM(\ngOF/\fm_\ngF\otimes R)$ is nilpotent,

\item $\uH^0(\cM_0/\fm_F) = 0$ and $\uH^0(\cM_1/\fm_F) = 0$. 
\end{enumerate}
Then the following holds:
\begin{enumerate}
\item The $\zeta$-isomorphism is defined for $\cM/\icoef$.

\item The natural diagram
\begin{equation*}
\xymatrix{
\big(\!\det\nolimits_{\ngOF} \RGamma(\cM)\big)/\icoef \ar[rr]^{\zeta_\cM} \ar[d] &&
\big(\!\det\nolimits_{\ngOF} \RGamma(\Der\cM)\big)/\icoef \ar[d] \\
\det\nolimits_{\ngOF/\icoef} \RGamma(\cM/\icoef) \ar[rr]^{\zeta_{\cM/\icoef}} &&
\det\nolimits_{\ngOF/\icoef} \RGamma(\Der\cM/\icoef)
}
\end{equation*}
is commutative.\end{enumerate}\end{lem}

\pf (1) Let $*\in\{0,1\}$. Lemma \ref{lfcohchar} shows that
$\uH^0(\cM_*) = 0$ and $\uH^1(\cM_*)$ is a free $\ngOF$-module of finite rank.
The base change theorem [\stacks{07VK}] then implies that
$\uH^0(\cM_*/\icoef) = 0$ and $\uH^1(\cM_*/\icoef)$ is a free $\ngOF/\icoef$-module of finite
rank. As $\cM(\ngOF/\fm_\ngF\otimes R)$ is nilpotent
Proposition \ref{globartcoh} implies that $\uH^0(\cM/\icoef) = 0$ and $\uH^1(\cM/\icoef)$
is a free $\ngOF/\icoef$-module of finite rank.
Hence the $\zeta$-isomorphism is defined for $\cM/\icoef$.
(2) Follows immediately from Proposition \lref{globzeta}{globcoeffchzeta}.\quod


%


%

\begin{lem}\label{globelltwistiso}%
Let $\cM$ be an elliptic shtuka. 
\begin{enumerate}
\item $\ngF\otimes_{\ngOF} \RGamma(\cM/\rami) = 0$.

\item The natural map
$\ngF\otimes_{\ngOF} \RGamma(\rami\cM) \to \ngF\otimes_{\ngOF} \RGamma(\cM)$
is a quasi-\hspace{0pt}isomorphism.\end{enumerate} 
\end{lem}

\pf In view of Proposition \ref{globellcoh} the result
follows from Theorem~\ref{ellcohfilt}.\quod

\begin{lem}\label{globellzetatwist}%
Let $\cM$ be an elliptic shtuka
given by a diagram $\cM_0 \shtuka{}{} \cM_1$.
If $\uH^0(\cM_0/\fm_F) = 0$ and $\uH^0(\cM_1/\fm_F) = 0$ then
then the natural square
\begin{equation*}
\xymatrix{
\ngF\otimes_{\ngOF} \det\nolimits_{\ngOF} \RGamma(\rami\cM) \ar[r]^{\isosign}
\ar[d]_{\zeta_{\rami\cM}}^{\rtviso} &
\ngF\otimes_{\ngOF} \det\nolimits_{\ngOF} \RGamma(\cM)
\ar[d]^{\zeta_{\cM}}_{\ltviso} \\
\ngF\otimes_{\ngOF} \det\nolimits_{\ngOF} \RGamma(\Der\rami\cM) \ar[r]^{\isosign} &
\ngF\otimes_{\ngOF} \det\nolimits_{\ngOF} \RGamma(\Der\cM)
}
\end{equation*}
is commutative.\end{lem}

\pf
Let $\cM^\circ$, $\rami\cM^\circ$, $\cM^\circ/\rami$ denote the pullbacks
of the respective shtukas to $\nsF\times X$. By construction we have a short exact
sequence $0 \to \rami\cM^\circ \to \cM^\circ \to \cM^\circ/\rami \to 0$.
We would like to apply Lemma \ref{zetadecomp} to this short exact sequence. To
do it we should verify that the following conditions are met:
\begin{enumerate}
\item $\RGamma(\cM^\circ)$ and $\RGamma(\Der\cM^\circ)$ are concentrated in degree $1$.

\item The sheaves of $\cM^\circ$ have cohomology concentrated in degree $1$.

\item $\RGamma(\cM/\rami) = 0$ and $\RGamma(\Der\cM/\rami) = 0$.

\item The sheaves of $\cM^\circ/\rami$ have cohomology concentrated in degree $0$.

\item $\cM^\circ/\rami$ is linear.
\end{enumerate}
(1) follows by Proposition \ref{globellcoh}, Lemma \ref{lfcohchar} implies (2)
and (3) holds by Lemma \ref{globelltwistiso}.
By construction $\cM^\circ/\rami$ is supported at a closed affine subscheme $\Spec \ngF \otimes
\cO_K/\rami$ of $\nsF \times X$. Thus the condition (4) is satisfied.
Finally the condition (5) follows since $\cM$ is elliptic.
Hence Lemma \ref{zetadecomp} demonstrates that the natural square
\begin{equation*}
\xymatrix{
\det\nolimits_{\ngF}^{\phantom{1}} \RGamma(\rami\cM^\circ) \ar[d]_{\zeta_{\rami\cM^\circ}}^{\rtviso}
\ar[r]^{\isosign} &
\det\nolimits_{\ngF}^{\phantom{1}} \RGamma(\cM^\circ)
\ar[d]^{\zeta_{\cM^\circ}}_{\ltviso} \\
\det\nolimits_{\ngF}^{\phantom{1}} \RGamma(\Der\rami\cM^\circ)
\ar[r]^{\isosign} &
\det\nolimits_{\ngF}^{\phantom{1}} \RGamma(\Der\cM^\circ)
}
\end{equation*}
commutes up to the sign $(-1)^{n m}$ where
\begin{equation*}
n = \dim_\ngF \uH^1(\cM^\circ) + \dim_\ngF \uH^1(\Der\cM^\circ).
\end{equation*}
However $\dim_\ngF \uH^1(\cM^\circ) = \dim_\ngF \uH^1(\Der\cM^\circ)$
since the shtuka $\cM$ is elliptic and so admits a regulator isomorphism
$\rho\colon \uH^1(\cM) \xrightarrow{\isosign}\uH^1(\Der\cM)$.
Thus the square above is in fact commutative.
Proposition \lref{globzeta}{globcoeffchzeta}
now implies the result. \quod


%

\section{\texorpdfstring{Euler products for $\ngOF$}{Euler products for O\_\ngF}}

In this section we define the Euler products for shtukas over $\ngOF\otimes R$.

\begin{lem}%
Let $k$ be a finite field extension of $\Fq$. Let $\cM$ be a locally
free shtuka on $\ngOF\otimes k$ given by a diagram
\begin{equation*}
M_0 \shtuka{i}{j} M_1.
\end{equation*}
If $\cM(\ngOF/\fm_\ngF \otimes k)$ is nilpotent then 
\begin{enumerate}
%
\item $i\colon M_0 \to M_1$ is an isomorphism,

\item $(1-i^{-1}j)\colon M_0 \to M_0$ is an $\ngOF$-linear isomorphism.
\end{enumerate}\end{lem}

\pf Indeed the ring $\ngOF \otimes k$ is noetherian and complete with respect to
a $\tau$-invariant ideal $\fm_\ngF \otimes k$. So Proposition
\lref{nilp}{nilpcomp} implies that
%
$\RGamma(\Der\cM) = 0$ and $\RGamma(\cM) = 0$.\quad

\begin{dfn}\label{euprodloc}%
Let $k$ be a finite field extension of $\Fq$. Let a locally free
shtuka $\cM$ on $\ngOF \otimes k$ be given by a diagram
\begin{equation*}
M_0 \shtuka{\,\,i\,\,}{j} M_1.
\end{equation*}
Assuming that $\cM(\ngOF/\fm_\ngF \otimes k)$ is nilpotent we define
\begin{equation*}
L(\cM) =
\det\nolimits_{\ngOF}(1 - i^{-1} j \mid M_0) \in \ngOF^\times.
\end{equation*}%
\end{dfn}

%
%

\begin{lem}\label{euconv}%
Let $n \geqslant 1$ be an integer.
If $\cM$ is a locally free shtuka on $\ngOF \otimes R$ such that
$\cM(\ngOF/\fm_\ngF \otimes R)$ is nilpotent then
for almost all maximal ideals $\fm \subset R$ we have
$L\big(\cM(\ngOF\otimes R/\fm)\big) \equiv 1 \pmod{\fm_\ngF^n}$
\end{lem}

\pf Let $\fm \subset R$ be a maximal ideal.
By construction
\begin{equation*}
L\big(\cM(\ngOF\otimes R/\fm)\big) \equiv
L\big(\cM(\ngOF/\fm_\ngF^n \otimes R/\fm)\big) \pmod{\fm_\ngF^n}
\end{equation*}
Applying Lemma \ref{euartok} to $\cM(\ngOF/\fm_\ngF^n \otimes R)$ we get
the result.\quod

\begin{dfn}\label{euprodglob}\index{nidx}{shtukas!$L(\cM)$}%
Let $\cM$ be a locally free shtuka on $\ngOF\otimes R$.
Assuming that $\cM(\ngOF/\fm_\ngF \otimes R)$ is nilpotent we define
\begin{equation*}
L(\cM) = \prod_\fm L(\cM(\ngOF \otimes R/\fm))^{-1} \in \ngOF^\times
\end{equation*}
where $\fm \subset R$ ranges over the maximal ideals.
This product converges by Lemma \ref{euconv}.\end{dfn}

\begin{rmk}Given a locally free shtuka $\cM$ on $\nsOFX$ such that
$\cM(\cO_F/\fm_F\otimes R)$ is nilpotent we write
$L(\cM)$ instead of $L(\cM(\ngOF\otimes R))$ to simplify the notation.
It is important to note that the closed points of $X$ in the complement of
$\Spec R$ are not taken into account.\end{rmk}

\begin{lem}%
If $\cM$ is a locally free shtuka
on $\ngOF\otimes R$ such that $\cM(\ngOF/\fm_\ngF \otimes R)$ is nilpotent
then for every $n \geqslant 1$ we have
$L(\cM) \equiv
L\big(\cM(\ngOF/\fm_\ngF^n\otimes R)\big) \pmod{\fm_\ngF^n}$.\quod
\end{lem}

\begin{prp}\label{euprodoneunit}%
If $\cM$ is a locally free shtuka
on $\ngOF\otimes R$ such that $\cM(\ngOF/\fm_\ngF \otimes R)$ is nilpotent
then $L(\cM) \equiv 1 \pmod{\fm_\ngF}$.\end{prp}

\pf Indeed $L(\cM) \equiv L(\cM(\ngOF/\fm_\ngF\otimes R)) \pmod{\fm_\ngF}$.
Since $\ngOF/\fm_\ngF$ is a field
Lemma \ref{euartone} implies that
$L(\cM(\ngOF/\fm_\ngF\otimes R)) = 1$.\quod

\section{Trace formula}

%

\begin{lem}\label{globelltracemod}%
Let $d \geqslant 1$. Let $\cM$ be an elliptic shtuka 
given
by a diagram
\begin{equation*}
\cM_0 \shtuka{}{} \cM_1.
\end{equation*}
Assume that the following holds:
\begin{enumerate}
\item $\uH^0(\cM_0/\fm_F) = 0$ and $\uH^0(\cM_1/\fm_F) = 0$,

\item $\cM/\rami^{2d}$ is linear.
\end{enumerate}
If the ramification ideal $\rami$ 
is contained in $\fm_K$ then 
$\zeta_\cM \equiv L(\cM) \cdot \det\nolimits_{\ngOF}(\rho_\cM) \pmod{\fm_\ngF^d}$.
\end{lem}

\pf Set $\ngC = \ngOF/\fm_\ngF^d$. Let $\cN = \cM/\fm_\ngF^d$ and let
\begin{equation*}
\cN_0 \shtuka{i}{j} \cN_1
\end{equation*}
be the diagram of $\cN$. We claim that the shtuka $\cN$ has the following
properties:
\begin{enumerate}
%
\item[(a)] $\uH^0(\Spec(\ngOF/\fm_F)\times X,\,\cN_0) = 0$
and $\uH^0(\Spec(\ngOF/\fm_F)\times X,\,\cN_1) = 0$,

\item[(b)] $\cN(\ngOF/\fm_\ngF \otimes R)$ is nilpotent,

\item[(c)] $\cN(\ngC \otimes \cO_K/\rami^{2d})$ is linear,

\item[(d)] $\rami^d \cdot \uH^1(\ngC \otimes \cO_K,\,\Der\cN) = 0$.
\end{enumerate}
The property (a) follows from the assumption (1), 
the property (b) holds
since $\cM$ is elliptic and (c) is a consequence of the assumption (2).
Finally, one gets (d) by applying
Lemma \lref{reg}{ellmfquot} to
the shtuka $\cM(\ngOF\complot\cO_K)/\fm_\ngF^d = \cN(\ngC \otimes \cO_K)$.

Now we apply Theorem \ref{arttrace}
to $\cN$ with $\ibase = \rami^d$ and conclude that the following is true:
\begin{enumerate}
\item[(i)] The artinian regulator $\rho_\cN$ is defined for $\cN$.

\item[(ii)] The $\zeta$-isomorphism $\zeta_\cN$ is defined for $\cN$.

\item[(iii)] $\zeta_\cN = L(\cN) \cdot \det_\ngC(\rho_\cN)$.
\end{enumerate}
The congruence $L(\cM) \equiv L(\cN) \pmod{\fm_\ngF^d}$ holds by construction.
Moreover
$\rho_\cN \equiv \rho_\cM \pmod{\fm_\ngF^d}$ by Lemma \ref{globellregred}
and $\zeta_\cN \equiv \zeta_\cM \pmod{\fm_\ngF^d}$
by Lemma \ref{globellzetared}. So the result follows.\quod 

\begin{lem}\label{globelltwistrelinv}%
Let $\cM$ be an elliptic shtuka given by a diagram $\cM_0\shtuka{}{}\cM_1$. 
Suppose that 
$\uH^0(\cM_0/\fm_F) = 0$ and $\uH^0(\cM_1/\fm_F) = 0$.
If $\alpha \in \ngOF^\times$ is the unique element such that
$\zeta_\cM = \alpha\cdot \det\nolimits_{\ngOF}(\rho_\cM)$
then
$\zeta_{\rami\cM} = \alpha \cdot \det\nolimits_{\ngOF}(\rho_{\rami\cM})$.
\end{lem}

\pf The square
\begin{equation*}
\xymatrix{
F\otimes_{\ngOF} \det\nolimits_{\ngOF}\RGamma(\rami\cM) \ar[r]^{\isosign}
\ar[d]_{\det_{\ngOF}(\rho_{\rami\cM})}^{\rtviso} &
F\otimes_{\ngOF} \det\nolimits_{\ngOF}\RGamma(\cM)
\ar[d]^{\det_{\ngOF}(\rho_{\cM})}_{\ltviso} \\
F\otimes_{\ngOF} \det\nolimits_{\ngOF}\RGamma(\Der\rami\cM) \ar[r]^{\isosign} &
F\otimes_{\ngOF} \det\nolimits_{\ngOF}\RGamma(\Der\cM)
}
\end{equation*}
is commutative by naturality of $\rho$. 
At the same time
Lemma \ref{globellzetatwist} shows that the natural square
\begin{equation*}
\xymatrix{
F\otimes_{\ngOF} \det\nolimits_{\ngOF} \RGamma(\rami\cM) \ar[r]^{\isosign}
\ar[d]_{\zeta_{\rami\cM}}^{\rtviso} &
F\otimes_{\ngOF} \det\nolimits_{\ngOF} \RGamma(\cM)
\ar[d]^{\zeta_{\cM}}_{\ltviso} \\
F\otimes_{\ngOF} \det\nolimits_{\ngOF} \RGamma(\Der\rami\cM) \ar[r]^{\isosign} &
F\otimes_{\ngOF} \det\nolimits_{\ngOF} \RGamma(\Der\cM)
}
\end{equation*}
is commutative. So we get the result.\quod

%
%

\breakflow
Finally we are ready to prove the trace
formula for regulators of elliptic shtukas.

\begin{thm}\label{elltrace}%
Let $\cM$ be an elliptic shtuka given by a diagram $\cM_0\shtuka{}{}\cM_1$. 
Suppose that 
$\uH^0(\cM_0/\fm_F) = 0$ and $\uH^0(\cM_1/\fm_F) = 0$.
If the ramification ideal $\rami$
is contained in $\fm_K$ then
\begin{equation*}
\zeta_\cM = L(\cM) \cdot \det\nolimits_{\ngOF}(\rho_\cM).
\end{equation*}\end{thm}

\breakflow
This theorem should certainly hold without the assumption on
on the cohomology of sheaves underlying $\cM$.
The bottleneck which prevents us from establishing this more natural version of
the trace formula is Lemma \ref{detlemma}. A more general variant of
this lemma is needed.
The proof of such a lemma appears to be
too messy at the moment.
Alas, the determinant theory of \cite{kmdet} is
not too user-friendly. 
Nevertheless Theorem \ref{elltrace} is still enough
to prove the class number formula for Drinfeld modules.

\afterall\noindent
\textit{Proof of Theorem \ref{elltrace}.}
%
Let $\alpha \in \ngOF^\times$ be the unique element such that
$\zeta_\cM = \alpha \cdot \det_{\ngOF}(\rho_\cM)$.
We will show that $\alpha \equiv L(\cM) \pmod{\fm_\ngF^d}$ for every $d \geqslant 1$.

Observe that $L(\cM) = L(\rami\cM)$.
Indeed the invariant $L$ depends only on
the restriction of a shtuka to $\ngOF\otimes R$ and $\cM(\ngOF\otimes R) =
\rami\cM(\ngOF\otimes R)$ by construction. At the same time Lemma \ref{globelltwistrelinv}
shows that $\zeta_{\rami\cM} = \alpha \cdot \det_{\ngOF}(\rho_{\rami\cM})$.
Furthemore the condition on the cohomology of sheaves underlying $\cM$
is preserved under twists by the ideal $\rami$.
We are thus free to replace $\cM$ by $\rami^n\cM$.

Observe that $(\rami\cM)(\ngOF\complot\cO_K)$ is the twist of
$\cM(\ngOF\complot\cO_K)$ by $\rami$ in
the sense of Definition \lref{reg}{defelltwist}.
So Proposition \lref{reg}{elltwistvanish} implies that
$(\rami^{2d}\cM)/\rami^{2d}$ is linear
for all $d \geqslant 0$.
Lemma \ref{globelltracemod} now implies the result.\quod


\chapter{The motive of a Drinfeld module}
\label{chapter:drmot}
\label{ch:drmot}
\label{ch:drinconstr}
\label{ch:coeffrings}
\label{ch:conilp}

We recall the notion of a Drinfeld module \cite{ell}, its motive as introduced by Anderson
\cite{and}, and a construction of Drinfeld \cite{nc} associating a shtuka to a Drinfeld
module.

Let $C$ be a smooth projective curve over $\Fq$ and $\infty\in C$ a closed point.
Denote $A = \Gamma(C\backslash\{\infty\}, \,\cO_C)$. Let $E$ be a Drinfeld $A$-module
over an $\Fq$-algebra $B$. The motive $M$ of $E$ is a left $A\otimes B\{\tau\}$-module
which is a locally free $A\otimes B$-module. Drinfeld's construction yields
a canonical shtuka on $C\times \Spec B$ extending $M$,
%
a ``compactification'' of $M$ in the direction of the coefficient
curve $C$.
In the subsequent chapters we will combine it with a compactification in the direction
of the base. Both compactifications are important to the proof of the class number
formula.

We should stress that all the material of this chapter is well-known,
with the possible exception of Proposition~\ref{motfingen} and the
co-\hspace{0pt}nilpotence property in Theorem~\ref{motivec}.
The latter two results are certainly known to experts.

All the nontrivial results of this chapter
are due to Drinfeld \cite{ell,nc}.
A detailed study of motives of Drinfeld modules over arbitrary, not necessarily reduced
base rings can be found in the preprint \cite{hartl} of Urs Hartl.

\section{Forms of the additive group}
\label{sec:forms}

In this section we work over a fixed $\Fq$-algebra $B$. To simplify the
exposition \emph{we assume that $B$ is reduced}. The theory which we attempt to present
here can be developed without this assumption. However it then becomes more
subtle. In applications we will only need the case of reduced $B$.

We equip $B$ with a $\tau$-ring structure given by the $q$-th power map. The
Frobenius $\tau$ and its powers are in a natural way $\Fq$-linear endomorphisms
of the additive group scheme $\bG_a$ over $B$.

\begin{lem}\label{gaend} Let $\End(\bG_a)$ be the ring of $\Fq$-linear
endomorphisms of $\bG_a$. The natural map $B\{\tau\} \to \End(\bG_a)$
is an isomorphism.\quod\end{lem}


\begin{lem}\label{gaautlin}
$B\{\tau\}^\times = B^\times$.
\end{lem}

\pf
Let $\varphi\in B\{\tau\}^\times$. If $K$ is a $B$-algebra which is a field then
the image of $\varphi$ in $K\{\tau\}$ must be an element of $K^\times$.
Therefore the constant coefficient of $\varphi$ is a unit and all other
coefficients are nilpotent.  Since $B$ is assumed to be reduced we conclude that
$\varphi\in B^\times$.\quod


%

\breakflow
Recall that an $\Fq$-vector space scheme $E$ is an abelian group scheme
equipped with a compatible $\Fq$-multiplication.

\begin{dfn} We say that an $\Fq$-vector space scheme $E$ is \emph{a form of
$\bG_a$} if it is Zariski-locally isomorphic to $\bG_a$.\end{dfn}

\breakflow
In the following we fix an $\Fq$-vector space scheme $E$ which is a form of $\bG_a$.
Our main object of study is the motive of $E$ which we now introduce.

\begin{dfn}\label{defmotive}\index{idx}{Drinfeld module!motive}\index{nidx}{Drinfeld modules!$M$, motive}%
The \emph{motive of $E$} is the abelian group $M = \Hom(E,\bG_a)$ of
$\bF_q$-linear group scheme morphisms from $E$ to $\bG_a$. The endomorphism
ring of $\bG_a$ acts on $M$ by composition making it into a left
$B\{\tau\}$-module.\end{dfn}

\begin{lem}\label{motsheaf}%
The formation of the
motive of $E$ commutes with arbitrary base change.\end{lem}

\pf For a scheme $Y$ over $X = \Spec B$ set $\cM(Y) = \Hom(E_Y, \bG_{a,Y})$ where
$E_Y$ denotes the pullback of $E$ to $Y$. Zarski descent for morphisms of
schemes implies that $\cM$ is a sheaf on the big Zariski site of $X$.
The abelian group $\cM(Y)$ carries a natural action of $\Gamma(Y,\cO_Y)$ on the
left. Together these actions make $\cM$ into an $\cO_X$-module.
The formation of $\cM$ is functorial in $E$.

If $E = \bG_a$ then
$\cM$ is the quasi-coherent sheaf defined by the left $B$-module $B\{\tau\}$.
Since $E$ is Zariski-locally isomorphic to $\bG_a$ we conclude that $\cM$ is
quasi-coherent.
Therefore the natural map
$S \otimes_B \cM(X) \to \cM(\Spec S)$
is an isomorphism for every $B$-algebra $S$. \quod

\begin{prp}\label{motzero}%
There exists a unique invertible $B$-submodule $M^0\subset M$ such that
the map
\begin{equation*}
B\{\tau\}\otimes_B M^0 \to M, \quad \varphi \otimes m \mapsto \varphi\cdot m
\end{equation*}
is an isomorphism.\end{prp}

\begin{cor}\label{motproj}%
The motive $M$ is a projective left $B\{\tau\}$-module.\quod\end{cor}

\begin{dfn}\index{nidx}{Drinfeld modules!$M^0$, degree $0$ part}%
We define the \emph{degree filtration} $M^\ast$ on $M$ in the following way.
For $n \geqslant 0$ we let
$M^n = B\{\tau\}^{n} \cdot M^0$
where $B\{\tau\}^{n} \subset B\{\tau\}$ is the submodule of
$\tau$-polynomials of degree at most $n$.
For $n < 0$ we set $M^n = 0$.
\end{dfn}

%
\afterall\noindent
\textit{Proof of Proposition \ref{motzero}.}
First let us prove unicity. If $M^0, N^0 \subset M$ are invertible $B$-submodules
such that the natural maps $B\{\tau\}\otimes_B M^0 \to M$ and 
$B\{\tau\}\otimes_B N^0 \to M$ are isomorphisms then we get an induced isomorphism
$B\{\tau\}\otimes_B M^0 \cong B\{\tau\}\otimes_B N^0$ of
left $B\{\tau\}$-modules.
Now Lemma \ref{gaautlin} implies that this isomorphism
comes from a unique $B$-linear isomorphism $M^0 \cong N^0$ which is compatible
with the inclusions $M^0 \hookrightarrow M$ and $N^0 \hookrightarrow M$.
As a consequence the submodules $M^0$ and $N^0$ of $M$ coincide.

Next let us prove the existence. If $E = \bG_a$ then we can take for $M^0$ the
submodule $B\cdot\tau^0 \subset B\{\tau\} = M$.
For an affine open subscheme $\Spec S \subset
\Spec B$ let $E_S$ be the pullback of $E$ to $\Spec S$ and let $M_S$ be the
motive of $E_S$. If $E_S$ is isomorphic to $\bG_{a,S}$ then by the remark above
we have an invertible $S$-submodule $M_S^0 \subset M_S$ satisfying the condition
of the proposition. Now the natural map $S \otimes_B M \to M_S$ is an
isomorphism by Lemma \ref{motsheaf} so the unicity part of the proposition
implies that $M_S^0$ glue to an invertible $B$-submodule $M^0 \subset M$.
The natural map $B\{\tau\}\otimes_B M^0 \to M$ is an isomorphism since it is so
after the pullback to every affine open subscheme $\Spec S\subset\Spec B$ such
that $E_S \cong \bG_{a,S}$. \quod

\begin{rmk}%
Without the assumption that $B$ is reduced the
existence part of Proposition \ref{motzero} still holds. However the submodule $M^0
\subset M$ is not unique anymore.%
\end{rmk}

\begin{lem}\label{motbasechange}%
The degree filtration on $M$ is stable under base change to an arbitrary
\emph{reduced} $B$-algebra.\end{lem}


\pf 
Let $S$ be a reduced $B$-algebra.
By Proposition~\ref{motzero} it is enough to show that the formation of
$M^0$ commutes with the base change.
Let $S$ be a $B$-algebra and let $M_S$ be the motive of $E$ over $S$.
The natural map
$S\{\tau\}\otimes_B M^0 \to S \otimes_B M$ is an isomorphism by definition of
$M^0$. Lemma \ref{motsheaf} shows that the natural map $S \otimes_B M \to M_S$
is an isomorphism. In particular $S \otimes_B M^0$ is in a natural way
an $S$-submodule of $M_S$.
Now if $S$ is reduced then Propostion~\ref{motzero}
implies that the image of $S \otimes_B M^0$ in $M_S$ is $(M_S)^0$. \quod


%
%
%
%
%
%
%

\begin{prp}\label{motzerorep}%
For every $B$-algebra $S$ the map
\begin{equation*}
E(S) \to \Hom_{B\{\tau\}}(M, S), \quad
e \mapsto \big(m \mapsto m(e)\big)
\end{equation*}
is an $\Fq$-linear isomorphism.\end{prp}

\breakflow\noindent
The left $B\{\tau\}$-module structure on $S$ is given by the $q$-power map.%

\pf Let $\cE$ be the presheaf on the big Zariski site of $\Spec B$ defined by
the functor $\Hom_{B\{\tau\}}(M, -)$. The map above defines a morphism of presheaves
$E \to \cE$.
Corollary \ref{motproj} implies that $\cE$ is a sheaf. 
By Lemma \ref{motsheaf} the formation of $M$ commutes with localization of $B$.
We can thus reduce to the case $E = \bG_a$ where the statement is clear.\quod

\begin{dfn}%
\index{nidx}{Drinfeld modules!$M^{\geqslant 1}$, positive degree part}%
\index{nidx}{Drinfeld modules!$\Omega = M/M^{\geqslant 1}$, the dual of the Lie algebra}%
We denote $M^{\geqslant 1} = B\{\tau\}\tau \otimes_B M^0$
and $\Omega = M/M^{\geqslant 1}$.
\end{dfn}

\begin{prp}\label{motpullback}%
The adjoint $\tau^\ast M \to M$ of the multiplication map $\tau\colon M \to M$
is injective with image $M^{\geqslant 1}$. \quod\end{prp}



\begin{rmk}%
We identify $\Omega$ with the $B$-module of Lie algebra homomorphisms
$\Lie_E \to \Lie_{\bG_a}$. An image of a morphism $m\in \Hom(E,\bG_a) = M$
in the quotient $\Omega = M/M^{\geqslant 1}$ is the induced morphism
of Lie algebras $dm\colon\Lie_E \to \Lie_{\bG_a}$.\end{rmk}

\begin{prp}\label{motlierep}%
For every $B$-algebra $S$ the map
\begin{equation*}
\Lie_E(S) \to
\Hom_B(\Omega, S), \quad
\varepsilon \mapsto \big(dm \mapsto dm(\varepsilon)\big)
\end{equation*}
is an isomorphism of $S$-modules.\quod\end{prp}

\section{Coefficient rings}

Let $A$ be an $\Fq$-algebra of finite type which is a Dedekind domain.
To such an algebra $A$ one can functorially associate a smooth connected
projective curve $C$ over $\Fq$ together with an open embedding $\Spec A
\subset C$.
We call $C$ the compactification of $\Spec A$.

\begin{dfn}\label{defcoeffring}\index{idx}{Drinfeld module!coefficient ring}We say
that $A$ is a \emph{coefficient ring} if the complement of $\Spec A$ in $C$
consists of a single point. This point is called the point of $A$ at
infinity and is denoted $\infty$.\end{dfn}

\begin{example*} $\Fq[t]$ is a coefficient ring.
$\Fq[t,s]/(s^2 - t^3 + 1)$ is a coefficient ring
provided $q$ is coprime to $6$.\end{example*}

\breakflow
Recall that an element $a \in A$ is called constant if it is algebraic over
$\Fq$. Since we do not assume $A$ to be geometrically irreducible there may be
constant elements not in $\Fq \subset A$.

\begin{lem}\label{coefffinflat} Let $A$ be a coefficient ring.
If $a \in A$ is not constant then
the natural map $\Fq[a] \to A$ is finite flat.
\end{lem}

\pf
Since $a$ is not constant the ring $A$ has no $\Fq[a]$-torsion
and so is flat over $\Fq[a]$.
The only nontrivial claim is that it is finite over $\Fq[a]$.

Let $C$ be the compactification of $\Spec A$ and let $\infty$ be the point in
the complement of $\Spec A$ in $C$. The inclusion $\Fq[a] \subset A$ induces a
morphism $C \to \mathbb{P}^1_\Fq$. This morphism is automatically proper. The
only point of $C$ which does not map to $\Spec\Fq[a] \subset \mathbb{P}^1_\Fq$
is $\infty$. Hence the preimage of $\Spec \Fq[a]$ in $C$ is $\Spec A$. We
conclude that the map $\Spec A \to \Spec\Fq[a]$ is proper. As a consequence it is finite
[\stacks{01WN}]. \quod

\breakflow
Attached to $A$ one has its local field at infinity $F$.
It is the completion of the function field of $C$ at the point $\infty$.
One has a canonical inclusion $A
\hookrightarrow F$.

\begin{dfn}\label{coeffdegmap}%
Let $A$ be a coefficient ring. For every nonzero $a \in A$ we define
\begin{equation*}
\deg(a) = -\nu(a)
\end{equation*}
where $\nu\colon F^\times \to \bZ$ is the normalized valuation.%
\end{dfn}

\breakflow
Observe that $\deg(a) = 0$ if and only if $a$ is a nonzero constant.

%
\begin{lem}\label{coeffdegeq} Let $A$ be a coefficient ring, $F$ the local field
at infinity, $k$ the residue field of $F$. If $a \in A$ is not constant then
\begin{equation*}
f \cdot \deg(a) = d
\end{equation*}
where $f = [k : \Fq]$, $d = [A : \Fq[a]]$.
\end{lem}

\pf Let $F_0$ be the local field
of $\Fq[a]$ at infinity. By construction $a^{-1}$ is a uniformizer of $F_0$.
Hence $\deg(a)$ equals the ramification index $e$ of $F$ over $F_0$. Moreover $f$
coincides with the inertia index of $F$ over $F_0$. Since $ef = [F : F_0]$ we
only need to prove that
$[F : F_0] = d$. Both $A$ and $\Fq[a]$ have a single point at infinity. Thus
\begin{equation*}
F_0 \otimes_{\Fq[a]} A = F.
\end{equation*}
Since the inclusion $\Fq[a] \subset A$ is finite flat it follows that $[F : F_0]
= d$. \quod

\section{Action of coefficient rings}

We continue working over the fixed $\Fq$-algebra $B$. As before we suppose that $B$
is reduced. Throughout this section we fix an $\Fq$-vector space scheme $E$ over
$B$ which is a form of $\bG_a$. We denote $M$ its motive.

We assume that $E$ is equipped with an action of a fixed coefficient ring $A$.
In other words we are given an $\Fq$-algebra homomorphism $\varphi\colon A \to \End(E)$.
The ring $A$ acts on $M = \Hom(E,\bG_a)$ on the right. As $A$ is commutative we
can view it as a left action. Thus $M$ acquires a structure of a left $A \otimes
B\{\tau\}$-module.
In this section we study how the $A \otimes B$-module structure on $M$ interacts
with the degree filtration.

\begin{lem}\label{motnotfingen} Assume that $B$ is noetherian. If $M^0$ is an $A
\otimes B$-submodule of $M$ then $M$ is \emph{not} a finitely generated $A \otimes
B$-module.\end{lem}

\pf Indeed
in this case
every $M^n \subset M$ is an
$A \otimes B$-submodule. Since the~quotients $M^n/M^{n-1}$ are not zero
we conclude
that $M$ contains an infinite increasing chain of $A \otimes B$-submodules.
Hence it is not finitely generated.\quod

\begin{lem}\label{motactfields} Let $a \in A$ and let $d \geqslant 0$. The
following are equivalent:
\begin{enumerate}
\item $M^0 a \subset M^d$ and the induced map $a\colon M^0 \to M^d/M^{d-1}$ is
an isomorphism.

\item The same holds after base change to every $B$-algebra $K$ which is a field.
\end{enumerate}\end{lem}

\pf (1) $\Rightarrow$ (2) is a consequence of Lemma \ref{motbasechange}.
(2) $\Rightarrow$ (1). Thanks to Lemma \ref{motbasechange} we may assume that $E
= \bG_a$. In this case the action of $A$ on $E$ is given by a homomorphism
$\varphi\colon A \to B\{\tau\}$. The condition (2) means that for
every $B$-algebra $K$ which is a field the polynomial $\varphi(a)$ has degree
$d$ in $K\{\tau\}$. Therefore the coefficient of $\varphi(a)$ at $\tau^d$ is a
unit while the coefficient at $\tau^n$ is nilpotent for every $n > d$. By
assumption of this section $B$ is reduced. Hence $\varphi(a)$ is of degree $d$
with top coefficient a unit. \quod

\begin{lem}\label{tmotfinconds} Assume that $A = \Fq[t]$. Let $r \geqslant 1$.
The following are equivalent:
\begin{enumerate}
\item $M^0 t \subset M^r$ and the induced map $t\colon M^0 \to M^r/M^{r-1}$ is
an isomorphism of $B$-modules.

\item The natural map $A \otimes M^{r-1} \to M$ is an isomorphism of $A \otimes
B$-modules.

\item $M$ is a locally free $A \otimes B$-module of rank $r$.
\end{enumerate}\end{lem}

\pf Thanks to Lemma \ref{motbasechange} we may assume that $E = \bG_a$.
In this case $M = B\{\tau\}$ and the degree filtration on $M$ is the filtration by
degree of $\tau$-polynomials. The action of $A$ is given by a homomorphism $\varphi\colon \Fq[t]
\to B\{\tau\}$.
We split the rest of the proof into several steps.

\breakflow
\textbf{Step 1.}
\textit{If \textup{(1)} holds then the natural map $A \otimes M^{r-1} \to M$ is
surjective.}

By assumption $\varphi(t)$ is of degree $r$ with top coefficient a unit.
Write $\varphi(t) = \psi + \alpha_r \tau^r$ with $\psi$ of degree less than
$r$. We then have
$\tau^r = \alpha_r^{-1} (\varphi(t) - \psi)$.
Multiplying both sides by $\tau^n$ on the left we obtain a relation
\begin{equation*}
\tau^{r+n} = \tau^n(\alpha_r)^{-1}(\tau^n \varphi(t) - \tau^n \psi ).
\end{equation*}
Induction over $n$ now shows that the image of the natural map
$A \otimes M^{r-1} \to M$ is the whole of $M$.

\breakflow
\textbf{Step 2.}
\textit{If $B$ is a field then \textup{(1)} implies \textup{(2)}.}

According to Step 1 the natural map $A \otimes M^{r-1} \to M$ is surjective. We
need to prove that it is injective.
Observe that for every nonzero $\alpha \in B$ and $n, d \geqslant 0$ the
$\tau$-polynomial $\alpha \tau^n \varphi(t^d)$ is of degree $r d + n$. Hence if $f
\in A \otimes B = B[t]$ is a polynomial of degree $d$ then $f \cdot \tau^n$
is of degree $r d + n$.

Now let $f_0, \dotsc,
f_{r-1} \in B[t]$.
If one of the $f_n$ is nonzero then there exists a unique $n \in
\{0,\dotsc,r-1\}$ such that $r \cdot \deg f_n + n$ is maximal. From the
observation above we deduce that $f_n \cdot \tau^n$ is of degree $r \deg f_n +
n$ while for every $m \ne n$ the element
$f_m \cdot \tau^m$ is of lesser degree. We conclude that
\begin{equation*}
f_0 \cdot 1 + f_1 \cdot \tau + \dotsc + f_{r-1} \cdot \tau^{r-1} \ne 0.
\end{equation*}

\textbf{Step 3.}
\textit{If $B$ is noetherian then \textup{(1)} implies \textup{(2)}.}

According to Step 1 the
natural map $A \otimes M^{r-1} \to M$ is surjective. We thus have a short exact
sequence
\begin{equation}\label{mottsurjses}
0 \to N \to A \otimes M^{r-1} \to M \to 0.
\end{equation}
Since $B$ is noetherian it follows that $N$ is a finitely generated $A \otimes
B$-module. 
By construction $M^{r-1}$ and $M$ are flat $B$-modules.
Therefore \eqref{mottsurjses} is a short exact sequence of flat $B$-modules.

Let $K$ be a $B$-algebra which is a field. Lemma
\ref{motbasechange} tells that the formation of $M$ commutes with base change to
$K$ and that the base change preserves the degree filtration. Therefore Step 2
shows that the second arrow of \eqref{mottsurjses} becomes an isomorphism after
base change to $K$. Since \eqref{mottsurjses} is a sequence of flat $B$-modules
we conclude that $N \otimes_B K = 0$ for every $K$.
As $N$ is a finitely generated $A \otimes B$-module Nakayama's lemma implies
that $N = 0$.

\breakflow
\textbf{Step 4.} \textit{\textup{(1)} implies \textup{(2)}.}
Write
\begin{equation*}
\varphi(t) = \alpha_0 + \alpha_1 \tau + \dotsc + \alpha_r \tau^r.
\end{equation*}
By assumption $\alpha_r$ is a unit.
Let $B_0$ be the $\Fq$-subalgebra of $B$ generated by $\alpha_i$ and
$\alpha_r^{-1}$. Let $E_0 =
\bG_{a,B_0}$ equipped with the action of $\Fq[t]$ given by $\varphi$. As
$\alpha_r$ is a unit in $B_0$ it follows that the assumption (1) holds for $E_0$.
Step 3 now implies that (2) holds for $E_0$. Lemma \ref{motbasechange}
shows that (2) holds for the base change of $E_0$ to $B$. As $E_0
\otimes_{B_0} B = E$ by construction the result follows.

\breakflow
\textbf{Step 5.} \textit{\textup{(2)} implies \textup{(3)}}.
By construction 
$M^{r-1}$ is a locally free $B$-module of rank $r$.
Therefore $A \otimes M^{r-1}$ is a locally free $A
\otimes B$-module of rank $r$.

\breakflow
\textbf{Step 6.} \textit{\textup{(3)} implies \textup{(1)}}.
Thanks to Lemma \ref{motactfields} we may suppose that
$B$ is a field. If $\varphi(t)$ is of degree $0$ then $M^0$ is an $A\otimes B$-submodule of
$M$. Lemma \ref{motnotfingen} then shows that $M$ is not a finitely generated $A
\otimes B$-module, a contradiction. Hence $\varphi(t)$ is of positive degree $d$.
Now Step 3 shows that $M$ is locally free of rank $d$ whence $d = r$. The
induced map $t\colon M^0 \to M^r /M^{r-1}$ is an isomorphism since the top
coefficient of $\varphi(t)$ is not zero. \quod

\breakflow
We now return to a general coefficient ring $A$.
Let $F$ be the local field of $A$ at infinity and let $k$ be the residue field
of $F$. We denote
\begin{equation*}
f = [ k : \Fq ]
\end{equation*}
the degree of the residue field extension at infinity.

%
%

\begin{prp}\label{motfinconds}Let $r \geqslant 1$ and let $a \in A$ be a nonconstant
element. The following are equivalent:
\begin{enumerate}
\item $M$ is a locally free $A \otimes B$-module of rank $r$.

\item $M^0 a \subset M^{f r \deg a}$ and
the induced map
\begin{equation*}
M^0 \xrightarrow{\,\,a\,\,} M^{f r \deg a}/M^{f r \deg a - 1}
\end{equation*} is an
isomorphism of $B$-modules.
\end{enumerate}\end{prp}

\pf
According to Lemma \ref{coefffinflat} the natural map $\Fq[a] \to A$ is
finite flat. Hence $M$ is locally free of rank $r$
as an $A \otimes B$-module if and only if it is locally free of rank $rd$ as an
$\Fq[a] \otimes B$-module where $d = [A:\Fq[a]]$.
Since $d = f \deg a$ by
Lemma \ref{coeffdegeq} the result follows
from Lemma \ref{tmotfinconds} applied to $t = a$. \quod

\breakflow
Assuming that the base ring $B$ is noetherian
we next show that the motive $M$ is a finitely generated $A \otimes B$-module if
and only if it is locally free. We include this result only for illustrative purposes.
It will not be used in the proof of the class number formula.

\begin{lem}Assume that $A = \Fq[t]$ and $B$ is a field. If $M$ is a finitely
generated $A \otimes B$-module then it is locally free of rank $\geqslant
1$.\end{lem}

\pf 
If $M^0 t \subset M^0$ then Lemma \ref{motnotfingen} shows that $M$ is not
finitely generated as an $\Fq[t] \otimes B$-module, t contradiction.
Therefore $M^0 t
\subset M^n$ for some $n \geqslant 1$.  Without loss of generality we may assume
that $M^0 t \not\subset M^{n-1}$. In this case the induced map $t\colon M^0 \to
M^n/M^{n-1}$ is nonzero. As $B$ is t field it is an isomorphism. Lemma
\ref{tmotfinconds} then shows that $M$ is t locally free $A \otimes B$-module of
rank $r \geqslant 1$. \quod

\begin{lem}\label{motrankineq}Assume that $A = \Fq[t]$ and $B$ is a DVR. Let
$K$ be the fraction field and $k$ the residue field of $B$.  If $M$ is a
finitely generated $A \otimes B$-module then $\rank_{A \otimes K} M \otimes_B K
\geqslant \rank_{A \otimes k} M\otimes_B k$.\end{lem}

\pf The Picard group of $B$ is trivial so by Proposition~\ref{motzero} we may
assume that $E = \bG_a$. In this case the $A$-action is given by a
homomorphism $\varphi\colon \Fq[t] \to B\{\tau\}$. Let $r_K$ be the degree of
$\varphi(t)$ in $K\{\tau\}$ and let $r_k$ be its degree in $k\{\tau\}$. Lemma
\ref{tmotfinconds} shows that $M \otimes_B K$ is a locally free $A \otimes
K$-module of rank $r_K$ while $M \otimes_B k$ is a locally free $A \otimes
k$-module of rank $r_k$. As $r_K \geqslant r_k$ by construction the result
follows.\quod

\begin{prp}\label{motfingen}%
If $B$ is noetherian and $\Spec B$ is connected then the following are
equivalent:
\begin{enumerate}
\item $M$ is a finitely generated $A \otimes B$-module.

\item $M$ is a locally free $A \otimes B$-module of constant rank.
\end{enumerate}\end{prp}

\pf (1) $\Rightarrow$ (2). If $a \in A$ is a nonconstant element then the map
$\Fq[a] \to A$ is finite flat by Lemma \ref{coefffinflat}. Hence to deduce (2)
it is enough to assume that $A = \Fq[t]$. 

Let $r\colon \Spec B \to \bZ_{\geqslant 1}$ be the function which sends a prime
$\fp \subset B$ to the rank of $M \otimes_B \Frac B/\fp$ as an $A \otimes \Frac
B/\fp$-module. We will show that $r$ is constant.

Let $\fp \subset \fq$ be primes of $B$ such that $\fp \ne \fq$. According to
[\stacks{054F}] there exists a discrete valuation ring $V$ and a morphism $B \to
V$ such that the generic point of $\Spec V$ maps to $\fp$ and the closed point
maps to $\fq$. Applying Lemma \ref{motrankineq} to the base change of $E$ to $V$
we deduce that $r(\fp) \geqslant r(\fq)$. Hence $r$ is lower semi-continuous.
However $A \otimes B$ is noetherian and $M$
is a finitely generated $A \otimes B$-module. The function $r$ is therefore also
upper semi-continuous. We conclude that it is in fact constant. Let us denote this
constant $r$.

Let $K$ be a $B$-algebra which is a field.
Lemma \ref{tmotfinconds} shows that $(M^0 \otimes_B K)t \subset (M^r \otimes_B
K)$ and the induced map $M^0 \otimes_B K \to (M^r/M^{r-1})\otimes_B K$ is an
isomorphism. Hence Lemma \ref{motactfields} shows that the same holds already on
the level of $B$. Applying Lemma \ref{tmotfinconds} again we conclude that $M$
is a locally free $\Fq[t] \otimes B$-module of rank $r \geqslant 1$. \quod

\section{Drinfeld modules}

We continue working over a fixed reduced $\Fq$-algebra $B$. Let $A$ be a coefficient
ring as in Definition \ref{defcoeffring}. 
As in the previous section $f$ denotes the degree of the residue field
extension at infinity.

\begin{dfn}\label{defdrinfeld}\index{idx}{Drinfeld module}%
A \emph{Drinfeld $A$-module} of rank $r \geqslant 1$ over $B$
is an $\Fq$-vector space scheme $E$ over $B$ equipped with an action of $A$
such that:
\begin{enumerate}
\item $E$ is a form of $\bG_a$.

\item The motive $M = \Hom(E,\bG_a)$ is a locally free $A \otimes B$-module of
rank $r$.
\end{enumerate}\end{dfn}

\breakflow\noindent
Proposition \ref{motfinconds} implies that (2) can be replaced by
\begin{enumerate}
\item[(2$'$)] There exists a nonconstant element $a \in A$ such that $M^0 a
\subset M^{f r \deg a}$ and the induced map $a\colon M^0 \to M^{f r \deg a} /
M^{f r \deg a - 1}$ is an isomorphism.\end{enumerate}
%
%
%
Using this fact it is easy to show that our definition is equivalent with Drinfeld's original definition
\cite{ell}. As in \cite{ell} the rank of our Drinfeld modules is constant on $\Spec B$.
%
%

\begin{prp}\label{motfiltprop}%
Let $E$ be a Drinfeld $A$-module of rank $r$ and let $M$ be its
motive. The degree filtration $M^\ast$ has the following properties.
\begin{enumerate}
\item $M^\ast$ is exhaustive.

\item $M^n$ is a locally free $B$-module of rank $\max(0,n + 1)$.

\item For every $n \geqslant 0$ and every nonzero $a \in A$ we have $M^n a
\subset M^{n + f r \deg a}$ and the induced map
\begin{equation*}
M^n/M^{n-1} \to M^{n + f r \deg a} / M^{n + f r \deg a - 1}
\end{equation*}
is an isomorphism.\end{enumerate}\end{prp}

\pf (1) and (2) 
are immediate from the definition of the degree filtration
and (3) follows from Proposition \ref{motfinconds}.\quod

\begin{prp} 
Let $E$ be a form of $\bG_a$ equipped with an action of $A$.
If $B$ is noetherian and $\Spec B$ is connected then
the following are equivalent:
\begin{enumerate}
\item $E$ is a Drinfeld $A$-module.

\item The motive $M = \Hom(E,\bG_a)$ is a finitely generated $A\otimes B$-module.
\end{enumerate}\end{prp}

\pf Follows instantly from Proposition \ref{motfingen}.\quod

\section{Co-nilpotence}

Before we state the main result of this chapter
let us introduce an auxillary notion. Let $R$ be a $\tau$-ring.

\begin{dfn}\label{defconilp}\index{idx}{shtuka!co-nilpotent}%
Let $M$ be an $R$-module shtuka given by a diagram
\begin{equation*}
M_0 \shtuka{\,\,i\,\,}{j} M_1.
\end{equation*}
We say that $M$ is \emph{co-nilpotent} if the adjoint
$j^a\colon \tau^\ast M_0 \to M_1$ of $j$ is an isomorphism 
and the compsition
\begin{equation*}
\tau^{\ast n}(u) \circ \dotsc \circ u
, \quad
u = (j^a)^{-1} \circ i,
\end{equation*}
is zero for $n \gg 0$.\end{dfn}

\begin{prp}\label{conilp}%
Let $M$ be an $R$-module shtuka 
given by a diagram
\begin{equation*}
M = \Big[ M_0 \shtuka{i}{j} M_1 \Big].
\end{equation*}
and let $N$ be a left $R\{\tau\}$-module.
If $M$ is co-nilpotent
then $\iHom_R(M,N)$ is nilpotent.\end{prp}

\pf Follows directly from the definition of $\iHom$
(Definition~\lref{genhomsht}{defhomsht}).
In verifying this it is convenient to
identify $M_1$ with $\tau^\ast M_0$ via
$j^a$.\quod 

\section{Drinfeld's construction}
\label{sec:drinconstr}

%
In this section we recall Drinfeld's
construction \cite{nc} of a shtuka attached to a Drinfeld module.
No originality is claimed. All nontrivial results are Drinfeld's.
%

Fix a coefficient ring $A$. Let $C$ be the projective compactification of $\Spec A$
and 
$\infty \in C$ the closed point in the complemenet of $\Spec A$.
Let $F$ be the local field of $C$ at $\infty$.
As in the previous sections $f$ denotes the degree of
the residue field of $F$ over $\Fq$.
Let $B$ be an $\Fq$-algebra. 
We equip $C \times \Spec B$ with an endomorphism $\tau$ which acts
as the identity on $C$ and as the $q$-Frobenius on $B$.
%
%
%

Fix an ample line bundle on $C$ which corresponds to the divisor $\infty$. Let
$\cO(1)$ be the pullback of this bundle to $C \times \Spec B$.

\begin{thm}\label{motivec}%
Let $E$ be a Drinfeld $A$-module of rank $r$ over $B$ and let $M$ be its motive.
The $A \otimes B$-module shtuka
\begin{equation*}
M \shtuka{\,\,1\,\,}{\tau} M
\end{equation*}
extends uniquely to a shtuka
\begin{equation*}
\cE = \Big[ \cE_{-1} \shtuka{\,\,i\,\,}{j} \cE_0 \Big]
\end{equation*}
on $C \times \Spec B$ with the following properties:
\begin{enumerate}
\item $\cE_{-1}$ and $\cE_0$ are locally free of rank $r$.

\item For every $n \in \bZ$ we have
\begin{align*}
\uH^0(C \times \Spec B, \,\cE_0(n)) &= M^{n f r}, \\
\uH^0(C \times \Spec B, \,\cE_{-1}(n)) &= M^{n f r - 1}
\end{align*}
as $B$-submodules of $M$.\end{enumerate}
Moreover the shtuka $\cE(\cO_F/\fm_F \otimes B)$ is co-nilpotent.%
\end{thm}

\begin{rmk}
The fact that $\cE(\cO_F/\fm_F\otimes B)$ is co-nilpotent is of fundamental
importance to our study. It implies that certain shtukas we construct out of
$\cE$ are nilpotent which in turn allows us to apply the theory of Chapters
\ref{chapter:reg} and \ref{chapter:trace} to Drinfeld modules.

The condition (2) can be replaced with the following pair of conditions:
\begin{itemize}
\item[(2a)] The adjoint of the $j$-arrow of $\cE(\cO_F/\fm_F \otimes B)$ is an isomorphism.

\item[(2b)] $\chi(\cE_{-1}) = 0$ for all points of $\Spec B$.
\end{itemize}
However we will use the fact that $\cE$ satisfies (2), so we opt for the
less elegant formulation.%
\end{rmk}

\afterall\noindent
\textit{Proof of Theorem \ref{motivec}.}
Uniqueness follows from (2). Let us prove the existence.
We denote
$\iota\colon \Spec(A\otimes B) \hookrightarrow C \times \Spec B$
the open immersion.
By Proposition~\ref{motfiltprop} the degree filtration $M^\ast$ has the
following properties:
\begin{itemize}
\item $M^\ast$ is exhaustive.

\item $M^n$ is a locally free $B$-module of rank $\max(0,n + 1)$.

\item For every $n \geqslant 0$ and every nonzero $a \in A$ we have $M^n a
\subset M^{n + f r \deg a}$ and the induced map
\begin{equation*}
M^n/M^{n-1} \to M^{n + f r \deg a} / M^{n + f r \deg a - 1}
\end{equation*}
is an isomorphism.\end{itemize}
We can thus invoke \cite[Corollary 1]{nc}, \cite[Proposition 3]{nc}
and conclude that there exists an increasing chain
\begin{equation*}
\dotsc \subset \cE_{-1} \subset \cE_0 \subset \cE_1 \subset \cE_2 \subset \dotsc
\end{equation*}
of locally free sheaves of rank $r$ on $C\times\Spec B$ extending $M$
such that:
\begin{enumerate}
\item $\uH^0(C \times \Spec B,\,\cE_n) = M^n$.

\item\label{defspiralspiral}$\cE_n(1) = \cE_{n + f r}$.

\item The $\tau$-multiplication map
$\tau\colon \iota_\ast M \to \tau_\ast\iota_\ast M$
sends $\cE_n$ to $\cE_{n+1}$.
%
%

\item The map $\tau^\ast(\cE_n/\cE_{n-1}) \to \cE_{n+1}/\cE_n$ induced
by the adjoint of $\tau$ is an isomorphism.
\end{enumerate}
In particular we get a shtuka
\begin{equation*}
\cE = \Big[ \cE_{-1} \shtuka{\,\,i\,\,}{j} \cE_0 \Big]
\end{equation*}
where $i\colon \cE_{-1} \to \cE_0$ is the inclusion
and $j$ is induced by $\tau$.
It satisfies all the conditions of this theorem.
What remains is to prove that $\cE(\cO_F/\fm_F\otimes B)$ is co-nilpotent.

Let $\alpha\colon\Spec(\cO_F/\fm_F\otimes B) \to C\times\Spec B$ be the closed immersion.
We denote $j^a\colon \tau^\ast\cE_n \to \cE_{n+1}$ the adjoint of $j$.
Property (4) implies that the natural map
$\tau^\ast\cE_{-1}/\cE_{-fr-1} \to \cE_0/\cE_{-fr}$
is an isomorphism.
This map coincides with $\alpha^\ast(j^a)$
since $\cE_{-fr+n} = \cE_{n}(1)$ by the property (2) above.
The same property now implies that the composition
\begin{equation*}
\tau^{\ast fr}(u) \circ \dotsc \circ u, \quad
u = \alpha^\ast(j^a)^{-1} \circ \alpha^\ast(i)
\end{equation*}
is zero.\quod

\chapter{The motive and the Hom shtuka}
\label{chapter:mothomsht}
\label{ch:hoch}
\label{ch:infbasecoh}
\label{ch:ihommotdesc}
\label{ch:ihomcoh}
\label{ch:ihomtate}

Let $A$ be a coefficient ring, $F$ the local field of $A$ at infinity
and $\omega$ the module of K\"ahler differentials of $A$ over $\Fq$.
Let $E$ be a Drinfeld $A$-module over a reduced $\Fq$-algebra $B$
and let $M = \Hom(E, \bG_a)$ be its motive. Rather than working
with $M$ directly we study the Hom shtuka
\begin{equation*}
\iHom_{A\otimes B}(M, \,\omega\otimes B).
\end{equation*}

This chapter has three groups of results. The first is purely algebraic and
deals with a Drinfeld $A$-module $E$ over a reduced $\Fq$-algebra
$B$. 
We show that
\begin{align*}
\RGamma(\iHom_{A\otimes B}(M, \,\omega \otimes B)) &\cong E(B)[-1], \\
\RGamma(\Der\iHom_{A\otimes B}(M, \,\omega \otimes B)) &\cong \Lie_E(B)[-1].
\end{align*}
This result is related to the formulas of Barsotti-Weil type
obtained by Papanikolas-\hspace{0pt}Ramachandran \cite{papa} and Taelman \cite{lt1mot}.
In some sense it is analogous to the classical Barsotti-Weil isomorphism
$\Ext^1(h(E), \unit) \cong E(k)$
for the motive $h(E)$ of
an elliptic curve $E$ over a field $k$.


The second group of results consists of two formulas.
The first applies to a Drinfeld module $E$ over a local field $K$.
We assume that the action of $A$ on $\Lie_E(K)$ extends to a continuous action
of $F$. Under this condition there is a natural map $\exp\colon \Lie_E(K) \to
E(K)$. We construct a quasi-\hspace{0pt}isomorphism
\begin{equation*}
\RGamma(\iHom_{A\otimes K}(M, \,\omega\complot K)) \cong \Big[ \Lie_E(K) \xrightarrow{\,\exp\,} E(K) \Big].
\end{equation*}
This formula is inspired by the work of Anderson \cite[Section 2]{and}.
Its proof relies on the fact that $\omega\complot K = b(F/A, K)$
is naturally the space of bounded functions $F/A \to K$ (see Definition~\ref{defbndfn}).

The second formula is a variant of the first.
It applies to a Drinfeld module $E$ over an $\Fq$-algebra $R$ which is discrete
in a finite product of local fields $K$. Assuming that the action
of $A$ on $\Lie_E(K)$ extends to a continuous action of $F$
we produce a quasi-\hspace{0pt}isomorphism
\begin{equation*}
\RGamma\!\big(\!\iHom_{A\otimes R}(M, \,\omega\complot \tfrac{K}{R})\big) \cong \Big[ \Lie_E(K) \xrightarrow{\,\exp\,} \frac{E(K)}{E(R)} \Big].
\end{equation*}
This result is crucial to the proof of the class number formula.


The third group of results relates Tate modules of a Drinfeld module over a 
field $k$ to the Hom shtuka $\iHom_{A\otimes k}(M, \,\omega \otimes k)$.
We translate certain results of David Goss
\cite[Section 5.6]{goss} concerning the motive $M$ of $E$ to the language of Hom shtukas.

%
%
%
%
%

Our theorems generalize to
Anderson modules \cite{and} in place of Drinfeld modules. Their proofs extend
without change. We limit our exposition to Drinfeld
modules since other important parts of our theory still depend on their special
properties.

Our approach 
relies on computations of
Hochschild cohomology of $A$ with coefficients in various $(A,A)$-bimodules
such as the
function spaces and the germ spaces.
First six sections of this 
chapter are devoted to such computations. Section~\ref{sec:hoch} relates 
Hochschild cohomology to shtuka cohomology.
Save for the last section the rest of the chapter is a straightforward
derivation of results on shtuka cohomology from the earlier Hochschild results.


%

In the case $A = \Fq[t]$
the shtukas isomorphic to $\iHom_{A\otimes B}(M, \,A\otimes B)$ 
have appeared in the work of Fang \cite{fang}. 
Generalizing an earlier construction of Taelman \cite{ltcarlitz}
he defined such shtukas in terms of explicit bases and matrices.
His work motivated the author of this text to study Hom shtukas.

\section{Hochschild cohomology in general}
\index{nidx}{Hochschild cohomology!$\RH(A,M)$, cohomology complex}%
\index{nidx}{Hochschild cohomology!$\uH^n(A,M)$, cohomology modules}%
Let $A$ be an associative unital $\Fq$-algebra and let $M$ be a complex of
$(A,A)$-bimodules. The Hochschild cohomology of $A$ with coefficients in $M$ is
by definition the complex
\begin{equation*}
\RHom_{A\otimes A^\circ}(A, M).
\end{equation*}
where $A^\circ$ is the opposite of $A$ and the $(A,A)$-bimodule structure on $A$
is the diagonal one. We denote this complex
$\RH(A,M)$ and $\uH^n(A, M)$ its cohomology groups. 

Even though Hochschild cohomology will figure prominently in our
computations, no nontrivial properties of it will be used.

\section{Hochschild cohomology of coefficient rings}

Let $A$ be a coefficient ring in the sense of
Definition \lref{coeffrings}{defcoeffring}.
Given an $(A,A)$-bimodule $M$ we denote
$D(A,M)$ the module of $M$-valued derivations over $\Fq$, i.e.
$\Fq$-linear maps $d\colon A \to M$ which satisfy the Leibniz identity
\begin{equation*}
d(a_0 a_1) = a_0 d(a_1) + d(a_1) a_0.
\end{equation*}
We denote $\partial\colon M \to D(A,M)$ the map which sends an element
$m \in M$ to the derivation $a \mapsto a m - m a$.

\begin{prp}\label{hochcanres}%
For every $(A,A)$-bimodule $M$ there is a natural quasi-\hspace{0pt}isomorphism
\begin{equation*}
\RH(A, M) \cong
\Big[ M \xrightarrow{\,\,\partial\,\,} D(A,M) \Big]
\end{equation*}
\end{prp}

\pf The $(A,A)$-bimodule $A$ has a canonical resolution
$0 \to I \to A \otimes A \to A \to 0$.
The $A\otimes A$-modules $I$ and $A \otimes A$ are invertible, so that
\begin{equation*}
\RHom_{A\otimes A}(A, M) =
\Big[ M \to \Hom_{A\otimes A}(I, M) \Big].
\end{equation*}
The natural map
$\Hom_{A\otimes A}(I,M) \to D(A,M)$, 
$f \mapsto (a \mapsto f(a \otimes 1 - 1 \otimes a))$,
is an isomorphism of $(A,A)$-bimodules. It identifies the complex above with the
complex in the statement of the lemma. \quod

\breakflow
We denote $\omega$ the module of K\"ahler differentials of $A$ over $\Fq$.
Let $M$ be an $A$-module.
Using the universal property of $\omega$ one easily proves the following lemma:

\begin{lem}\label{hochmiso}
The map $M \to D(A, \,\omega \otimes_A M)$
given by the formula 
$m \mapsto (a \mapsto da \otimes_A m)$ is an isomorphism.\quod\end{lem}

\begin{lem}\label{hochses}%
The sequence
\begin{equation}\label{hochseq}
0 \to \omega \otimes M \xrightarrow{\,\,\partial\,\,}
D(A, \,\omega \otimes M) \xrightarrow{\,\,-D(\mu)\,\,}
D(A, \,\omega \otimes_A M) \to 0
\end{equation}
is exact. Here 
$\mu\colon \omega \otimes M \to \omega \otimes_A M$
is the contraction map.\end{lem}

\breakflow
Observe that in \eqref{hochseq} we use the morphism $-D(\mu)$ instead
of the more natural $D(\mu)$. With this choice of sign the statements of several
results become more straightforward.
One example is Proposition \ref{hochgffsq}.

\afterall\noindent
\textit{Proof of Lemma \ref{hochses}.}
Let $0 \to I \to A \otimes A \to A \to 0$ be the canonical resolution of the
diagonal $(A,A)$-bimodule $A$. The composition with the multiplication map
$A \otimes A \to A$ induces a homomorphism
$\nu\colon D(A,A\otimes A) \to D(A,A)$.
One easily shows that the natural sequence
\begin{equation}\label{hochseqbase}
0 \to A \otimes A \to D(A,A\otimes A) \xrightarrow{-\nu} D(A,A) \to 0
\end{equation}
is exact. The tensor product of \eqref{hochseqbase}
with $\omega \otimes M$ over $A \otimes A$ coincides with the sequence
\eqref{hochseq}.
Moreover
\begin{equation*}
D(A,A) = 
\Hom_{A\otimes A}(I,A) = \Hom_A(I/I^2,A) = \Hom_A(\omega,A).
\end{equation*}
So to prove the lemma it is enough to show that 
\begin{equation*}
\Tor_1^{A\otimes A}(\Hom_A(\omega,A), \,\omega \otimes M) = 0
\end{equation*}
with the diagonal $(A,A)$-bimodule structure on $\Hom_A(\omega,A)$.
However 
\begin{equation*}
\Hom_A(\omega, A) \otimes^{\mathbf{L}}_{A\otimes A}
(\omega \otimes M) =
\Hom_A(\omega, A) \otimes^{\mathbf{L}}_A \omega
\otimes^{\mathbf{L}}_A M = M[0].\quod
\end{equation*}

\begin{prp}\label{hochtensor}%
For every $A$-module $M$ there is a natural $A$-linear
quasi-\hspace{0pt}isomorphism
$\RH(A, \,\omega \otimes M) \xrightarrow{\isosign} M[-1]$.
It is the following composition:
\begin{equation*}
\def\arraystretch{1.5}
\begin{array}{r@{}lll}
\RH(A, \,\omega \otimes M)\,\,
& \xrightarrow{\isosign} 
\big[ \omega \otimes M \xrightarrow{\partial} D(A,\,\omega \otimes M) \big]
&\textup{(Proposition \ref{hochcanres})} & \\
&\xrightarrow{\isosign} 
D(A, \,\omega \otimes_A M)[-1]
&\textup{(Lemma \ref{hochses})} & \\
& \xleftarrow{\isosign}
M[-1] & \textup{(Lemma \ref{hochmiso})} & \quod
\end{array}\end{equation*}\end{prp}

%

\section{Vanishing}

Let $A$ be a coefficient ring
and let $F$ be the local field of $A$ at $\infty$.

%

\begin{lem}$A \otimes_{A\otimes A} (A\complot F) = 0$.\end{lem}

\pf
Fix a nonconstant element $t \in A$.
The element $1 \otimes t - t \otimes 1$ acts by zero on $A$.
We will show that it becomes a unit in $\Fq[t]\complot F$
and a fortiori in $A\complot F$.

Let $\alpha$ be the image of $t$ in $F$.
Recall from Section~\ref{sec:otalgexamples} that $\Fq[t]\complot F$
is the ring of Tate series $F\langle t \rangle$ in variable
$t$ and with coefficients in $F$. The image of $1 \otimes 1 - t \otimes 1$
in $F\langle t \rangle$ is the polynomial $\alpha - t$.
It is invertible since $|\alpha| > 1$. \quod

\breakflow
We will use the following vanishing statement:

\begin{prp}\label{hochvanish}%
If $M$ is an $A\complot F$-module then $\RH(A, M) = 0$.\end{prp}

\breakflow\noindent
The $(A,A)$-bimodule structure on $M$ is given by the natural
morphism $A\otimes A \to A\complot F$.

\afterall\noindent\textit{Proof of Proposition~\ref{hochvanish}.}
Indeed $A\complot F$ is flat over $A\otimes A$ and the pullback of the
diagonal $A\otimes A$-module $A$ to $A\complot F$ is zero. So 
$\RHom_{A\otimes A}(A, M) = \RHom_{A\complot F}(0, M) = 0$
by extension of scalars.\quod

\section{Preliminaries on function spaces}%
\label{sec:hochfuncprelim}%
As before $A$ stands for a coefficient ring in the sense of
Definition \lref{coeffrings}{defcoeffring}.
From now on we denote
$F$ the local field of $A$ at $\infty$
and
$\omega$ the module of K\"ahler differentials
of $A$ over $\Fq$.
Let $\res\colon \omega\otimes_A F \to \Fq$ be the residue map at infinity
composed with the trace map to $\Fq$.

Let $U$ be one of the topological $A$-modules
$F$ or $F/A$ and let $M$ be a locally compact $A$-module.
Consider $a(U,M)$, the space of bounded locally constant $\Fq$-linear maps
$U \to M$ as in Definition \ref{defbndlc}.
The action of $A$ on the source and the target makes $a(U,M)$ into an $(A,A)$-bimodule.
In the following sections we study the Hochschild cohomology
of $a(U,M)$ for certain $A$-modules $M$.

Besides $a(U,M)$ we will also use the space of bounded $\Fq$-linear maps
$b(U,M)$ as in Definition \ref{defbndfn}
and the space of germs $g(U,M)$ as in Definition \ref{defgerm}.
The $(A,A)$-bimodule structures on such spaces are given by the action of $A$
on the source and the target.
To improve the legibility we will generally write $\RH(-)$ instead of
$\RH(A,-)$. 

\begin{dfn}\label{defphi}%
We denote
$\Res\colon \omega \otimes M \to a(F/A,\,M)$
the map given by the formula
$\eta \otimes m \mapsto (x \mapsto \res(\eta x)m)$.%
\end{dfn}

\breakflow\noindent
The map $\Res$ is a topological isomorphism of $(A,A)$-bimodules
by Corollary~\ref{funcmoddesc}.
We will view it interchangeably as a map
to $a(F/A,M)$, $a(F, M)$ and $b(F/A,M)$. 
Corollary~\ref{funcmoddesc} also shows that
$\Res$ extends uniquely
to a topological isomorphism 
of $(A,A)$-bimodules
$\omega\complot M \xrightarrow{\isosign} b(F/A,\,M)$.
We denote it $\Res$ too. The composition of $\Res$
with the natural map $b(F/A,M) \to g(F,M)$ will be
denoted $[\Res]$.

\begin{lem}\label{omegagermsesexact}%
The sequence
\begin{equation*}
0 \to \omega\otimes M \xrightarrow{\phantom{\,\,[\Res]\,\,}} \omega\complot M \xrightarrow{\,\,[\Res]\,\,} g(F,M) \to 0
\end{equation*}
is exact.\end{lem}

\begin{dfn}%
We denote
$\delta\colon \RH(g(F,M)) \to \RH(\omega\otimes M)[1]$ 
the natural morphism arising from the short exact
sequence of Lemma~\ref{omegagermsesexact}.\end{dfn}

\afterall\noindent\textit{Proof of Lemma~\ref{omegagermsesexact}. }
The quotient map $F \to F/A$ is a local isomorphism
and so induces an isomorphism $g(F/A, M) \to g(F,M)$.
Now the map $\Res$ identifies the sequence in question
with the natural sequence
$0 \to a(F/A, M) \to b(F/A, M) \to g(F/A,M) \to 0$
which is exact by Proposition~\ref{germses}.\quod
\section{Vector spaces}
\label{sec:hochfuncspf}%
We continue using the notation and the conventions
of Section~\ref{sec:hochfuncprelim}.
In this section we work with a fixed finite-\hspace{0pt}dimensional $F$-vector space $V$.
According to Proposition~\ref{hochcanres} we have
\begin{equation*}
\RH(g(F, V)) = 
\Big[ g(F, V) \xrightarrow{\partial} D(A, \,g(F, V)) \Big].
\end{equation*}
One verifies easily that the map
$v \mapsto (x \mapsto xv)$ defines a morphism
of complexes $V[0] \to \RH(g(F/A,V))$.

\begin{prp}\label{hochgffsq}%
The map
$v \mapsto (x \mapsto xv)$
induces a quasi-\hspace{0pt}isomorphism
$V[0] \xrightarrow{\isosign} \RH(g(F,V))$
which fits into a commutative square
\begin{equation*}
\xymatrix{
V[0] \ar[d]_{\ltviso} \ar@{=}[r] & V[0] \\
\RH(g(F,V)) \ar[r]^{\delta\,\,}&
\RH(\omega\otimes V)[1] \ar[u]^{\ltviso}_{\textup{Prp. }\ref{hochtensor}}
}
\end{equation*}%
\end{prp}

\breakflow
Although this statement appears innocent, we do not know an easy
proof. We split the proof into several lemmas.

\begin{lem}\label{hochaffqi}%
The map $\Res\colon \omega\otimes V \to a(F,V)$
induces a quasi-\hspace{0pt}isomorphism
$\RH(\omega\otimes V) \xrightarrow{\isosign} \RH(a(F,V))$.\end{lem}


\pf
A quick inspection reveals that the natural sequence
\begin{equation*}
0 \to a(F/A,V) \to a(F,V) \to a(A,V) \to 0
\end{equation*}
is exact. Moreover $a(A,V) = b(A,V)$ since $A$ is discrete.
Now Proposition~\ref{bndmod} shows that the $(A,A)$-bimodule structure
on $b(A,V)$ extends to an $A\complot F$-module structure.
Hence $\RH(b(A,V)) = 0$ by Proposition~\ref{hochvanish} and the result follows.\quod


%

\begin{lem}\label{hochdefmu}%
There exists a unique continuous map $\mu\colon a(F,F) \to F^*$ with the
following property. If $f\colon F \to \Fq$ is a continuous $\Fq$-linear map and
$\alpha \in F$ then the image of the map $x \mapsto f(x) \alpha$ under $\mu$ is
the map $x \mapsto f(x\alpha)$.
\end{lem}

\pf Consider the topological $\Fq$-vector space $F^* \otimes_{\textup{ic}} F$
(Definition \lref{indcot}{indcotdef}).
Let $\mu\colon F^* \otimes_{\textup{ic}} F \to F^*$ be the map which sends a
tensor $f \otimes \alpha$ to the function $x \mapsto f(\alpha x)$.
This map is easily shown to be continuous in the ind-tensor product topology.
It induces a unique continuous map $F^* \indcot F \to F^*$ by completion.
Proposition~\ref{bndlcmod} identifies $F^* \indcot F$ with
$a(F,F)$ and we get the result. \quod


\breakflow
We denote
$\rho\colon \omega \otimes_A F \xrightarrow{\isosign} F^*$ the isomorphism
induced by the residue pairing (Theorem \lref{resdual}{resdual}).

\begin{lem}\label{hochaffses}%
The diagram
\begin{equation}\label{hochaffcohdiag}
\vcenter{\vbox{\xymatrix{
0 \ar[r] & a(F,F) \ar[r]^{\partial\quad} & D(A,\,a(F,F)) \ar[r]^{-D(\mu)} & D(A,\,F^*) \ar[r] & 0 \\
0 \ar[r] & \omega \otimes F \ar[u]^{\Res} \ar[r]^{\partial\quad} & D(A,\,\omega \otimes F) \ar[u]_{D(\Res)} \ar[r]^{-D(\mu)} &
D(A,\,\omega \otimes_A F) \ar[r] \ar[u]^{\ltviso}_{\rho} & 0
}}}
\end{equation}
is commutative with exact rows.
\end{lem}

\pf
%
Let $\eta \in \omega$ and $\alpha \in F$.
By definition $\mu(\Res(\eta \otimes \alpha))$ is the function
$x \mapsto \res(x\alpha\eta)$.
It is the image of the element $\eta \otimes_A \alpha \in \omega \otimes_A F$
under $\rho$. So the right square of \eqref{hochaffcohdiag} is commutative.
The left square commutes by definition.
The bottom row is exact by Lemma~\ref{hochses}.
So Lemma~\ref{hochaffqi} implies that the top row is exact.\quod


\begin{lem}\label{mutrace}%
If $\Lambda$ is a finite-\hspace{0pt}dimensional discrete $\Fq$-vector space,
$f\colon F \to \Lambda$ and $g\colon \Lambda \to F$ continuous $\Fq$-linear maps then
$\mu(g \circ f)\colon F \to \Fq$ is the map
$x \mapsto \tr_\Lambda(f \circ x \circ g)$.\end{lem}

\pf 
If $\Lambda = \Fq$ then the claim is true by definition of $\mu$.
In general pick a splitting $\Lambda = \Lambda_1 \oplus \Lambda_2$. Let $f_1$, $f_2$ be
the compositions of $f$ with the projections to $\Lambda_1$, $\Lambda_2$ and let $g_1$,
$g_2$ be the restrictions of $g$ to $\Lambda_1$, $\Lambda_2$. We then have
\begin{equation*}
\mu(g \circ f) = \mu(g_1 \circ f_1) + \mu(g_2 \circ f_2).
\end{equation*}
At the same time
\begin{equation*}
\tr_\Lambda(f \circ x \circ g) = \tr_{\Lambda_1}(f_1 \circ x \circ g_1) +
\tr_{\Lambda_2}(f_2 \circ x \circ g_2).
\end{equation*}
So the claim follows by induction on the dimension of $\Lambda$.\quod

\afterall\noindent\textit{Proof of Proposition~\ref{hochgffsq}.}
First we extend the square in question to the right:
\begin{equation*}
\xymatrix{
V[0] \ar[d] \ar@{=}[r] & V[0] \ar@{=}[r] & V[0] \\
\RH(g(F,V)) \ar[r]^{\delta\,\,} &
\RH(\omega\otimes V)[1] \ar[u]^{\ltviso}_{\textup{Prp. }\ref{hochtensor}} \ar[r]^{\isosign} &
\RH(a(F,V))[1] \ar[u]_{\rtviso}
}
\end{equation*}
Here the bottom right arrow is the natural quasi-\hspace{0pt}isomorphism
of Lemma~\ref{hochaffqi}
and the right vertical arrow is defined by commutativity.
To prove the proposition it is enough to show that the outer rectangle
\begin{equation}\label{hochgffsquare}
\vcenter{\vbox{\xymatrix{
V[0] \ar[d] \ar@{=}[r] & V[0] \\
\RH(g(F,V)) \ar[r] &
\RH(a(F,V))[1] \ar[u]_{\rtviso}
}}}
\end{equation}
is commutative and the bottom map is a quasi-\hspace{0pt}isomorphism.

By construction the bottom arrow in \eqref{hochgffsquare}
is the map $\delta\colon \RH(g(F,V)) \to \RH(a(F,V))[1]$ arising from the
short exact sequence
\begin{equation*}
0 \to a(F,V) \to b(F,V) \to g(F,V) \to 0.
\end{equation*}
The $(A,A)$-bimodule structure on $b(F,V)$ extends to an $A\complot F$-module
structure by Proposition~\ref{bndmod}. Proposition~\ref{hochvanish}
then implies that $\RH(b(F,V)) = 0$.
So $\delta$ is a quasi-\hspace{0pt}isomorphism.
%
What remains is to show that \eqref{hochgffsquare} is commutative.

By naturality one reduces to the case $V = F$.
All the maps in \eqref{hochgffsquare} are
$F$-linear since $a(F,V)$, $b(F,V)$ and $g(F,V)$ are $(F,F)$-bimodules.
So to prove 
the commutativity
it is enough to consider
the element $1 \in F$.

Applying $\uH^0$ to \eqref{hochgffsquare} we get a square
\begin{equation*}
\xymatrix{
F \ar[d]_{\ltviso} \ar@{=}[r] & F \\
\uH^0(A,g(F,F)) \ar[r]^{\delta}_{\bisosign} &
\uH^1(A,a(F,F)) \ar[u]_{\rtviso}
}
\end{equation*}
The image of $1 \in F$ under the composition of the left, bottom and right
arrows can be described in the following way. Pick a bounded function
$f\colon F \to F$ such that $f(x) = x$ for all $x$ in an open neighbourhood of
$0$. Let $D\colon A \to a(F,F)$ be the derivation given by the formula
$D(a) = a f - f a$.
By Lemma~\ref{hochmiso} 
there exists a unique $\alpha \in F$ such that
for all $a \in A$ we have
\begin{equation*}
\mu(D(a)) = \rho(da \otimes_A \alpha).
\end{equation*}
Lemma~\ref{hochaffses} implies that
the image of $1 \in F$ is the element $-\alpha \in F$.
It is independent of the choice of $f$. 
Our goal is to prove that $\alpha = -1$.

Pick elements $a, b \in A$ such that $z = a b^{-1}$ is a uniformizer of
$F$.
Let $k$ be the residue field of $F$ so that $F = k(\!(z)\!)$.
There is a well-defined element $dz \in \omega[b^{-1}]$.

First we extend $D\colon A \to a(F,F)$ to a derivation
$D\colon F \to a(F,F)$ using the formula $D(x) = x f - f x$.
Leibniz rule for the products $z = a b^{-1}$ and $1 = b b^{-1}$ implies that
\begin{equation*}
D(z) = D(a) b^{-1} - z D(b) b^{-1}.
\end{equation*}
By assumption $\mu(D(x)) = \rho(dx \otimes_A \alpha)$ for all $x \in A$. Hence
\begin{equation*}
\mu(D(z)) = \rho(da \otimes_A \alpha b^{-1} - db \otimes_A \alpha z b^{-1}).
\end{equation*}
At the same time we have an identity
\begin{equation*}
dz \otimes_{A[b^{-1}]} \alpha = da \otimes_A \alpha b^{-1} - db \otimes_A \alpha zb^{-1}
\end{equation*}
in $\omega[b^{-1}] \otimes_{A[b^{-1}]} F = \omega \otimes_A F$.
Therefore
\begin{equation}\label{derofuni}
\mu(D(z)) = \rho(dz \otimes_{A[b^{-1}]} \alpha).
\end{equation}

Next we construct a suitable bounded function $f\colon F \to F$.
Consider the function 
defined by the formula
\begin{equation*}
f\Big( \sum_n \alpha_n z^n\Big) = 
\sum_{n \geqslant 0} \alpha_n z^n
\end{equation*}
with $\alpha_n \in k$. It is clearly bounded and satisfies
$f(x) = x$ for all $x \in \cO_F$. Moreover for every $\alpha \in k$ and
$n \in \bZ$ we have
\begin{equation*}
D(z)\colon \alpha z^n \mapsto
\left\{\begin{array}{ll}
-\alpha, &n = -1, \\
0, & n \ne -1.
\end{array}\right.
\end{equation*}
%
Applying Lemma \ref{mutrace} to the function $D(z)$, 
the vector space $\Lambda = k$ 
and the inclusion $k \hookrightarrow F$
we conclude
that $\mu(D(z)) = -\rho(dz \otimes_{A[b^{-1}]} 1)$.
In view of \eqref{derofuni} it follows that $\alpha = -1$.\quod

\section{Locally compact modules}
\label{sec:hochexp}

In this section we fix a finite-dimensional $F$-vector space
$V$, a locally compact $A$-module $M$ and a continuous
$A$-linear map $e\colon V \to M$. We assume that
$e$ is a \emph{local isomorphism} as in Definition~\ref{deflociso}.
One may view $V$ as the ``tangent space'' of $M$
at $0$ and $e\colon V \to M$ as the ``exponential map'' in the manner of
Lie theory.

\begin{dfn}%
We denote $C_e$ the $A$-module complex
$\Big[\, V \xrightarrow{\,\,e\,\,} M \,\Big]$.\end{dfn}

\afterall\noindent
Our aim is to prove that $\RH(\omega\complot M) = C_e$.

\begin{lem}\label{hochezerocoh}%
$\uH^0(C_e)$ is a finitely generated projective $A$-module.\end{lem}

\pf Since $e$ is a local isomorphism it follows that $\uH^0(C_e) = \ker e$
is a discrete $A$-submodule of $V$.
As $V$ is finite-\hspace{0pt}dimensional we conclude that $\uH^0(C_e)$
is finitely generated. It is then projective by virtue of being torsion free.\quod

\breakflow
In the following $\uD(A)$ stands for the derived category of $A$-modules.

\begin{lem}\label{hocheuniq}%
Let $h\colon M[-1] \to C_e$ be the map given by the identity in degree one
and let $C$ be a complex of $A$-modules.
If $\uH^0(C)$ is finitely generated and $\uH^n(C) = 0$ for all $n < 0$
then the map
\begin{equation*}
\Hom_{\uD(A)}(C_e, \,C) \xrightarrow{- \circ h} \Hom_{\uD(A)}(M[-1], \,C)
\end{equation*}
is injective.%
\end{lem}

\pf The map $h$ extends to a distinguished triangle
\begin{equation*}
M[-1] \xrightarrow{\,h\,} C_e \to V[0] \to M[0].
\end{equation*}
Applying $\Hom_{\uD(A)}(-,\,C)$ we get an exact sequence
\begin{equation*}
\Hom(V[0],\,C) \to \Hom(C_e,\,C) \xrightarrow{-\circ h} \Hom(M[-1],\,C).
\end{equation*}
It is thus enough to show that $\Hom(V[0],\,C) = 0$.
Since $A$ is of global dimension $1$ 
there is a non-canonical quasi-\hspace{0pt}isomorphism
\begin{equation*}
C \cong \bigoplus_n \uH^n(C)[-n].
\end{equation*}
As $\uH^n(C) = 0$ for $n < 0$
it follows that $\Hom(V[0], C) = \Hom_A(V,\,\uH^0(C))$.
However $V$ is uniquely divisible
and $\uH^0(C)$ has no divisible elements besides zero.
So the latter Hom is zero.\quod

\begin{prp}\label{hochge}%
The map $v \mapsto (x \mapsto e(xv))$ induces a quasi-\hspace{0pt}isomorphism
$V[0] \xrightarrow{\isosign} \RH(g(F,M))$ which fits into a commutative square
\begin{equation*}
\xymatrix{
V[0] \ar[r]^e \ar[d]_{\ltviso} & M[0] \\
\RH(g(F,M)) \ar[r]^\delta & \RH(\omega\otimes M)[1]\ar[u]_{\textup{Prp. }\ref{hochtensor}}^{\ltviso}
}
\end{equation*}\end{prp}

\pf The map $e\colon V \to M$ induces a morphism of short exact sequences
\begin{equation*}
\xymatrix{
0 \ar[r] & \omega \otimes V \ar[d]^{1 \otimes e} \ar[r] & \omega\complot V \ar[d]^{1\complot e} \ar[r] & g(F,V) \ar[d]^{e\circ-} \ar[r] & 0 \\
0 \ar[r] & \omega \otimes M \ar[r] & \omega\complot M \ar[r] & g(F,M) \ar[r] & 0
}
\end{equation*}
As $e$ is a local isomorphism it follows 
that the induced map
$g(F,V) \to g(F,M)$ is an isomorphism. Taking the cohomology we get a commutative
diagram
\begin{equation*}
\xymatrix{
\RH(g(F,V)) \ar[r]^{\delta} \ar[d]_{\ltviso}^{e\circ-} &
\RH(\omega\otimes V)[1] \ar[d]^{1\otimes e} \ar[rr]^{\quad\textup{Prp. }\ref{hochtensor}}_{\quad\isosign} && V[0] \ar[d]^{e} \\
\RH(g(F,M)) \ar[r]^{\delta} & \RH(\omega \otimes M)[1] \ar[rr]^{\quad\textup{Prp. }\ref{hochtensor}}_{\quad\isosign} && M[0]
}
\end{equation*}
The result now follows from Proposition~\ref{hochgffsq}.\quod

\begin{thm}\label{hochexp}%
There exists a quasi-\hspace{0pt}isomorphism
$\RH(\omega\complot M) \xrightarrow{\isosign} C_e$
with the following properties.
\nopagebreak%
\begin{enumerate}
\nopagebreak%
\item%
\nopagebreak%
It is natural in $V$, $M$, $e$.

\item
It is the unique morphism in $\uD(A)$ 
such that the square
\begin{equation*}
\xymatrix{
\RH(\omega\otimes M) \ar[rr]\ar[d]^{\rtviso}_{\textup{Proposition }\ref{hochtensor}}
&& \RH(\omega\complot M) \ar[d]^{\phantom{\textup{Proposition }\ref{hochgffsq}}} \\
M[-1] \ar[rr]^{\textup{identity in degree 1}} && C_e
}
\end{equation*}
is commutative.

\item It makes the square 
\begin{equation*}
\xymatrix{
\RH(\omega\complot M) \ar[rr] \ar[d]_{\phantom{\textup{Proposition }\ref{hochtensor}}}
&& \RH(g(F,M)) \\
C_e \ar[rr]^{\textup{identity in degree 0}} && V[0] \ar[u]^{\ltviso}_{\textup{Proposition }\ref{hochge}}
}
\end{equation*}
commutative.%
\end{enumerate}\end{thm}

\pf Consider the diagram
\begin{equation*}
\xymatrix{
\RH(\omega\otimes M)\ar[r]\ar[d]_{\textup{Proposition }\ref{hochtensor}}^{\rtviso} &
\RH(\omega\complot M)\ar[r] &
\RH(g(F,M)) \ar[r]\ar[d]^{\textup{Proposition }\ref{hochge}}_{\ltviso} & [1] \\
M[-1] \ar[r]^{\textup{id. in deg. 1}} & C_e \ar[r]^{\textup{id. in deg. 0}} & V[0] \ar[r]^e & M[0]
}
\end{equation*}
Proposition~\ref{hochge} implies that there exists a map
$\RH(\omega\complot M) \to C_e$ completing it into a morphism of distinguished triangles.
Such a map is necessarily a quasi-\hspace{0pt}isomorphism.

The square (2) commutes by construction.
$\uH^0(C_e)$ is finitely generated by Lemma~\ref{hochezerocoh}.
Applying Lemma~\ref{hocheuniq} to $C = C_e$ we conclude that
the map
\begin{equation*}
\Hom_{\uD(A)}\big(\RH(\omega\complot M), \,C_e\big) \to
\Hom_{\uD(A)}\big(\RH(\omega\otimes M), \,C_e\big)
\end{equation*}
given by composition with
$\RH(\omega\otimes M) \to \RH(\omega\complot M)$ is injective.
So the unicity part of (2) follows. It remains to prove naturality.

Suppose that we are given a commutative square
\begin{equation*}
\xymatrix{
V \ar[d]_e \ar[r]^{\xi} & W \ar[d]^f \\
M \ar[r]^u & N
}
\end{equation*}
where
\begin{itemize}
\item $W$ is a finite-\hspace{0pt}dimensional $F$-vector space,

\item $N$ is a locally compact $A$-module,

\item $f$ is an $A$-linear local isomorphism,

\item $\xi$ and $u$ are continuous and $A$-linear.
\end{itemize}
Let $C_u\colon C_e \to C_f$ be the morphism of complexes
induced by $\xi$ and $u$.
Consider the diagram
\begin{equation*}
\xymatrix{
\RH(\omega\otimes M) \ar[d]_{1\otimes u}\ar[r] & \RH(\omega\complot M) \ar[d]_{1\complot u}\ar[r]^{\quad\quad\isosign} & C_e \ar[d]^{C_u} \\
\RH(\omega\otimes N) \ar[r] & \RH(\omega\complot N) \ar[r]^{\quad\quad\isosign} & C_f
}
\end{equation*}
We want to prove that the right square is commutative.
Now $\uH^0(C_f)$ is finitely generated by Lemma~\ref{hochezerocoh}
so Lemma~\ref{hocheuniq} implies that the map
\begin{equation*}
\Hom_{\uD(A)}\big(\RH(\omega\complot M), \,C_f\big) \to
\Hom_{\uD(A)}\big(\RH(\omega\otimes M), \,C_f\big)
\end{equation*}
given by composition with $\RH(\omega\otimes M) \to \RH(\omega\complot M)$
is injective. It is thus enough to prove that the outer rectangle is commutative.
However the outer rectangle decomposes into two commutative squares
\begin{equation*}
\xymatrix{
\RH(\omega\otimes M) \ar[d]_{1\otimes u}\ar[rr]^{\textup{Prp. }\ref{hochtensor}} &&
M[-1] \ar[d]_{u}\ar[rr]^{\textup{id. in deg. }1} && C_e \ar[d]^{C_u} \\
\RH(\omega\otimes N) \ar[rr]^{\textup{Prp. }\ref{hochtensor}} &&
N[-1] \ar[rr]^{\textup{id. in deg. }1} && C_f
}
\end{equation*}
by construction.\quod

\section{The main lemma}
\label{sec:hoch}

The following simple lemma is our main tool to compute the cohomology of Hom
shtukas associated to Drinfeld modules.

Let $A$ and $B$ be associative unital $\Fq$-algebras.
For all left $A\otimes B$-modules $M$ and $N$ we have
$\Hom_{A\otimes B}(M, N) =
\Hom_{A\otimes A}(A, \,\Hom_B(M,N))$.
Taking derived functors we obtain a natural map
$\RHom_{A\otimes B}(M, N) \to \RH(A, \,\RHom_B(M,N))$.

\begin{lem}\label{hoch}%
For all left $A\otimes B$-modules $M$ and $N$
the natural map
\begin{equation*}
\RHom_{A\otimes B}(M, N) \to \RH(A, \,\RHom_B(M,N))
\end{equation*}
is a quasi-\hspace{0pt}isomorphism.%
\end{lem}

\breakflow
In applications $A$ will be
commutative in which case the quasi-isomorphism will be $A$-linear.

\afterall\noindent
\textit{Proof of Lemma \ref{hoch}}.
Let $M$ be a projective left $A \otimes B$-module and $N$ an injective left
$A \otimes B$-module.
We claim that $\Hom_B(M,N)$ is an injective $(A,A)$-bimodule.
Indeed $M$ is a flat $A$-module since it is a direct summand of a free $A\otimes
B$-module and $A \otimes B$ is $A$-flat since $B$ is $\Fq$-flat.
So
the functor
\begin{equation*}
\Hom_{A\otimes A^\circ}(-, \Hom_B(M,N)) =
\Hom_{A \otimes B}(- \otimes_A M, N).
\end{equation*}
is exact as the
composition of exact functors $-\otimes_A M$ and
$\Hom_{A\otimes B}(-, N)$.

Let $M$ and $N$ be arbitrary left $A \otimes B$-modules. 
Let $M^\bullet$ be a projective resolution of $M$ and $N^\bullet$ an injective
resolution of $N$. The argument above shows that 
$\Hom_B^\bullet (M^\bullet,N^\bullet)$ is a bounded below complex of injective
$(A,A)$-bimodules. Hence
\begin{align*}
\RHom_{A\otimes A^\circ}(A,\Hom_B^\bullet(M^\bullet,N^\bullet))
&=\Hom_{A\otimes A^\circ}(A, \Hom_B^\bullet(M^\bullet,N^\bullet)) \\
&=\Hom_{A \otimes B}^\bullet(M^\bullet,N^\bullet)\\
&=\RHom_{A \otimes B}(M, N).
\end{align*}

Next we claim that $M^\bullet$ is a projective resolution of $M$ as a
$B$-module.
Indeed if $P$ is a projective left $A \otimes B$-module then it is a direct
summand of a free $A \otimes B$-module.
However $A \otimes B$ is a projective $B$-module since
$\Hom_B(A \otimes B, -) = \Hom_{\Fq}(A,-)$. So $P$ is a projective $B$-module
and $M^\bullet$ is a projective $B$-module resolution of $M$.
It follows that
\begin{equation*}
\RHom_B(M,N) = \Hom_B^\bullet(M^\bullet,N[0]) =
\Hom_B^\bullet(M^\bullet,N^\bullet).
\end{equation*}
Thus 
$\RHom_{A\otimes A^\circ}(A,\RHom_B(M,N)) =
\RHom_{A\otimes A^\circ}(A,\Hom_B^\bullet(M^\bullet,N^\bullet))$.\quod

\breakflow
More generally this lemma applies
to associative unital algebras $A$, $B$ over an arbitrary commutative ring
$k$ provided that $A$ is a projective $k$-module and $B$ is a flat $k$-module.

\section{The Hom shtuka}
\label{sec:homshthoch}

Let $A$ be a coefficient ring, $E$ a Drinfeld $A$-module over
a reduced\footnote{The assumption that $B$ is reduced is inessential for the theory
discussed in this chapter. However in Chapter~\ref{chapter:drmot} we defined
the degree filtration on $M$ under this assumption,
so we are obliged to keep it here as well.}
$\Fq$-algebra $B$ and $M$ the motive of $E$.
As in Chapter~\ref{chapter:drmot} we denote $M^{\geqslant 1}$
the positive degree part of $M$ and $\Omega = M/M^{\geqslant 1}$.

We equip $A\otimes B$ with an endomorphism $\tau$ acting
as the identity on $A$ and as the $q$-Frobenius on $B$.
Note that $(A\otimes B)\{\tau\} = A \otimes (B\{\tau\})$
where the latter $\tau$ is the $q$-power map.
The motive $M = \Hom(E,\bG_a)$ carries a natural left action
of $B\{\tau\}$ via $\bG_a$ and a right action of $A$ via $E$.
Since $A$ is commutative we can view $M$ as a left $(A\otimes B)\{\tau\}$-module.

In this section
we study the sthuka
$\iHom_{A\otimes B}(M,N)$ for varying left $A\otimes B\{\tau\}$-modules
$N$.
To simplify the expressions we will generally write $\iHom(M,N)$ and
$\Hom(M,N)$ omitting the subscript $A\otimes B$.

\begin{prp}\label{ihommotdesc}%
The shtuka $\iHom(M,N)$ is represented by the diagram
\begin{equation*}
\Hom(M, N) \shtuka{\,\, i\,\,}{j} \Hom(M^{\geqslant 1}, N)
\end{equation*}
where $i$ is the restriction to $M^{\geqslant 1}$ and
$j$ sends an element $f$ to the map $\tau m \mapsto \tau f(m)$.
\end{prp}

\pf According to Proposition \lref{drmot}{motpullback}
the adjoint $\tau^\ast M \to M$ of the
multiplication map $\tau\colon M \to M$ is injective with image $M^{\geqslant 1}$.
So the result is a consequence of Lemma \lref{genhomsht}{rtaumodhomsht}. \quod

%

\begin{prp}\label{ihomhoch}%
There is a natural quasi-\hspace{0pt}isomorphism
\begin{equation*}
\RGamma(\iHom(M,N))
\xrightarrow{\isosign}
\RH(A,\,\Hom_{B\{\tau\}}(M,N)).
\end{equation*}
It is the following composition:
\begin{equation*}
\def\arraystretch{1.5}
\begin{array}{r@{}lll}
\RGamma(\iHom(M,N))\,\,
& \xleftarrow{\isosign} 
\RHom_{A\otimes B\{\tau\}}(M,N)
&\textup{(Theorem \ref{homshtcohrhom})} & \\
&\xrightarrow{\isosign} 
\RH(A, \,\Hom_{B\{\tau\}}(M,N))
&\textup{(Lemma \ref{hoch})} &
\end{array}\end{equation*}\end{prp}

\pf Indeed $\Hom_{B\{\tau\}}(M,N) = \RHom_{B\{\tau\}}(M,N)$ since
$M$ is a projective left $B\{\tau\}$-module by Corollary~\ref{motproj}.\quod

%
\breakflow
Consider the quotient $\Omega = M/M^{\geqslant 1}$.
Observe that $\Hom_{A\otimes B}(\Omega, N)$ is naturally
a submodule of $\Hom_{A\otimes B}(M,N)$.

\begin{lem}\label{derihomzerocoh}%
$\Hom_{A\otimes B}(\Omega,N) = \uH^0(\Der\iHom(M,N))$.\end{lem}

\pf Indeed Proposition~\ref{ihommotdesc} shows
that $\uH^0(\Der\iHom(M,N))$ consists of the
$A\otimes B$-linear maps $f\colon M \to N$ which
restrict to zero on $M^{\geqslant 1}$.\quod

\breakflow
The functor
$N \mapsto \Der\iHom(M,N)$ is exact.
So the isomorphism of Lemma~\ref{derihomzerocoh}
induces a natural map
$\RHom_{A\otimes B}(\Omega,N) \to \RGamma(\Der\iHom(M,N))$.

\begin{lem}\label{derihomcoh}%
The natural map
$\RHom_{A\otimes B}(\Omega,N) \to\RGamma(\Der\iHom(M,N))$
is a quasi-\hspace{0pt}isomorphism.\end{lem}

\pf The short exact sequence
$0 \to M^{\geqslant 1} \to M \to \Omega \to 0$
is a projective $A\otimes B$-module resolution of $\Omega$.
Therefore
\begin{equation*}
\RHom_{A\otimes B}(\Omega,N) =
\Big[ \Hom_{A\otimes B}(M,N) \xrightarrow{\,\,i\,\,}
\Hom_{A\otimes B}(M^{\geqslant 1}, N) \Big]
\end{equation*}
where $i$ is the restriction to $M^{\geqslant 1}$.
Now Theorem~\ref{shtaffcoh} in combination with
Proposition~\ref{ihommotdesc} shows that
the complex above computes $\RGamma(\Der\iHom(M,N))$.\quod

\begin{prp}\label{derihomhoch}%
There is a natural quasi-\hspace{0pt}isomorphism
\begin{equation*}
\RGamma(\Der\iHom(M,N)) \xrightarrow{\isosign}
\RH(A, \,\Hom_B(\Omega, N)).
\end{equation*}
It is the following composition:
\begin{equation*}
\def\arraystretch{1.5}
\begin{array}{r@{}lll}
\RGamma(\Der\iHom(M,N))\,\,
& \xleftarrow{\isosign} 
\RHom_{A\otimes B}(\Omega,N)
&\textup{(Lemma \ref{derihomcoh})} & \\
&\xrightarrow{\isosign} 
\RH(A, \,\Hom_{B}(\Omega,N))
&\textup{(Lemma \ref{hoch})} &
\end{array}\end{equation*}\end{prp}

\pf $\Omega$ is an invertible $B$-module
so $\Hom_B(\Omega,N) = \RHom_B(\Omega,N)$ and the result
follows.\quod

\section{Formulas of Barsotti-Weil type}

We keep the setting of the previous section.
As usual $\omega$ stands for the module of K\"ahler differentials
of $A$ over $\Fq$.
%
Given a map $m\colon E \to \bG_a$ we denote $dm$
the induced map from $\Lie_E$ to $\Lie_{\bG_a}$.
%
In the following $S$ is an arbitrary $B$-algebra.

\begin{lem}\label{barweillieind}%
The natural map
\begin{equation*}
\omega \otimes \Lie_E(S) \to
\Hom_B(\Omega, \,\omega \otimes S), \quad
\eta \otimes \alpha \mapsto
(dm \mapsto \eta \otimes dm(\alpha))
\end{equation*}
is an isomorphism of $(A,S)$-bimodules.\end{lem}

\pf The natural map
$\omega \otimes \Hom_B(\Omega, S) \to \Hom_B(\Omega,\,\omega\otimes S)$
is an isomorphism since $\Omega$ is a finitely generated $B$-module.
Proposition~\ref{motlierep} identifies
$\Hom_B(\Omega, S)$ with $\Lie_E(S)$ and the result follows.\quod

\begin{lem}\label{barweilind}%
The natural map
\begin{equation*}
\omega \otimes E(S) \to
\Hom_{B\{\tau\}}(M, \,\omega \otimes S), \quad
\eta \otimes e \mapsto (m \mapsto \eta \otimes m(e))
\end{equation*}
is an isomorphism of $(A,A)$-bimodules.\end{lem}

\pf The motive $M$ is a finitely generated $B\{\tau\}$-module
by Proposition~\ref{motzero}.
Hence the natural map
$\omega \otimes \Hom_{B\{\tau\}}(M, S) \to
\Hom_{B\{\tau\}}(M, \,\omega \otimes S)$
is an isomorphism.
Proposition~\ref{motzerorep} identifies
$\Hom_{B\{\tau\}}(M,S)$ with $E(S)$ and we get the result.\quod

\breakflow
Combining the lemmas above with the results of Section~\ref{sec:homshthoch}
and the isomorphism $\RH(A,\,\omega \otimes N) \cong N[-1]$ of
Proposition~\ref{hochtensor} we obtain the formulas of Barsotti-Weil type
annouced in the introduction.
In the case $A = \Fq[t]$ related results were obtained in
\cite{papa, lt1mot}.

\begin{thm}\label{barweilliemot}%
There is a natural $A\otimes S$-linear
quasi-\hspace{0pt}isomorphism
\begin{equation*}
\RGamma(\Der\iHom(M,\,\omega \otimes S)) \xrightarrow{\isosign} \Lie_E(S)[-1].
\end{equation*}
It is given by the composition
\begin{equation*}
\def\arraystretch{1.5}
\begin{array}{r@{}ll}
\RGamma(\Der\iHom(M,\,\omega \otimes S))\,\,
& \xrightarrow{\isosign} 
\RH(A, \,\Hom_B(\Omega,\,\omega\otimes S))
&\textup{(Proposition \ref{derihomhoch})} \\
&\xleftarrow{\isosign} 
\RH(A, \,\omega \otimes \Lie_E(S))
&\textup{(Lemma \ref{barweillieind})} \\
&\xrightarrow{\isosign}
\Lie_E(S)[-1]
&\textup{(Proposition \ref{hochtensor})} \quod
\end{array}\end{equation*}%
\end{thm}


\begin{thm}\label{barweilmot}%
There is a natural $A$-linear
quasi-\hspace{0pt}isomorphism
\begin{equation*}
\RGamma(\iHom(M,\,\omega \otimes S)) \xrightarrow{\isosign} E(S)[-1].
\end{equation*}
It is given by the composition
\begin{equation*}
\def\arraystretch{1.5}
\begin{array}{r@{}ll}
\RGamma(\iHom(M,\,\omega \otimes S))\,\,
& \xrightarrow{\isosign} 
\RH(A, \,\Hom_{B\{\tau\}}(M,\,\omega\otimes S))
&\textup{(Proposition \ref{ihomhoch})} \\
&\xleftarrow{\isosign} 
\RH(A, \,\omega \otimes E(S))
&\textup{(Lemma \ref{barweilind})} \\
&\xrightarrow{\isosign}
E(S)[-1]
&\textup{(Proposition \ref{hochtensor})} \quod
\end{array}\end{equation*}%
\end{thm}


\section{The exponential complex}

As before $A$ is a coefficient ring, $F$ the local field of $A$ at infinity
and $\omega$ the module of K\"ahler differentials of $A$ over $\Fq$.
Let $K$ be a finite product of local fields containing $\Fq$.
Fix a Drinfeld $A$-module $E$ over $K$.
We denote $M$ the motive of $E$ and
$\Omega$ the quotient $M/M^{\geqslant 1}$.
We assume that the action of $A$ on $\Lie_E(K)$ extends to a continuous action
of $F$.

\begin{dfn}\label{defexp}\index{idx}{exponential map!of a Drinfeld module}%
The \emph{exponential map} of the Drinfeld module $E$ is a map
\begin{equation*}
\exp\colon
\Lie_E(K) \to E(K)
\end{equation*}
satisfying the
following conditions:
\begin{enumerate}
\item $\exp$ is a homomorphism of $A$-modules,

\item $\exp$ is an analytic function with derivative $1$ at zero in the
following sense. Fix an $\Fq$-linear isomorphism of group schemes $E \cong
\bG_a$. It identifies $E(K)$ with $K$ while its differential identifies
$\Lie_E(K)$ with $K$. We demand that the resulting map $\exp\colon K \to K$ is
given by an everywhere convergent power series of the form
\begin{equation*}
\exp(z) = z + a_1 z^q + a_2 z^{q^2} + \dotsc
\end{equation*}
with coefficients in $K$.
\end{enumerate}
\end{dfn}

\begin{prp}%
The exponential map exists, is unique and is a local isomorphism
in the sense of Definition~\ref{deflociso}.\end{prp}

\pf For the existence and unicity see \cite[Theorem 2.1]{delhus}.
The exponential map is a local isomorphism since it has nonzero derivative at zero.\quod


\begin{dfn}The \emph{exponential complex} of $E$ is the $A$-module complex
\begin{equation*}
C_{\exp} = \Big[ \Lie_E(K) \xrightarrow{\,\,\exp\,\,} E(K) \Big].
\end{equation*}\end{dfn}

\afterall\noindent
Our goal is to show that $\RGamma(\iHom(M,\,\omega\complot K)) = C_{\exp}$.
To do it we apply the results of Section~\ref{sec:hochexp} to
the map $\exp\colon\Lie_E(K) \to E(K)$.


\begin{lem}\label{bndlieiso}%
The map $\omega\complot\Lie_E(K) \to \Hom_K(\Omega,\,\omega\complot K)$
defined by the commutative square
\begin{equation*}
\xymatrix{
\omega\complot\Lie_E(K) \ar[d]_{\Res}^{\rtviso} \ar[rrr] &&&\Hom_K(\Omega,\,\omega\complot K) \ar[d]^{\Res\circ-}_{\ltviso} \\
b(F/A,\,\Lie_E(K)) \ar[rrr]^{f \mapsto (dm \mapsto dm \circ f)} &&&\Hom_K(\Omega,\,b(F/A,K))
}
\end{equation*}
is an isomorphism of $A\complot K$-modules.
Furthermore
it is compatible with the isomorphism
$\omega\otimes\Lie_E(K) \xrightarrow{\isosign} \Hom_K(\Omega,\,\omega\otimes K)$
of Lemma~\ref{barweillieind}.%
\end{lem}

\pf The first claim follows since $\Omega$ is a free $K$-module of rank $1$.
The compatibility follows since $\Res$ identifies $\omega\otimes\Lie_E(K)$ with
the subspace $a(F/A,\,\Lie_E(K))$ of $b(F/A,\,\Lie_E(K))$.\quod


\begin{lem}\label{bndeiso}%
The map $\omega\complot E(K) \to \Hom_{K\{\tau\}}(M,\,\omega\complot K)$
defined by the commutative square
\begin{equation*}
\xymatrix{
\omega\complot E(K) \ar[d]_{\Res}^{\rtviso} \ar[rr] && \Hom_{K\{\tau\}}(M,\,\omega\complot K) \ar[d]^{\Res\circ-}_{\ltviso} \\
b(F/A,\,E(K)) \ar[rr]^{f \mapsto (m \mapsto m \circ f)\quad\quad} && \Hom_{K\{\tau\}}(M,\,b(F/A,K))
}
\end{equation*}
is an isomorphism of $(A,A)$-bimodules.
Furthermore it is compatible with the isomorphism
$\omega\otimes E(K) \xrightarrow{\isosign} \Hom_{K\{\tau\}}(M,\,\omega\otimes K)$
of Lemma~\ref{barweilind}.%
\end{lem}

\pf Follows since $M$ is a free $K\{\tau\}$-module of rank $1$.\quod

\begin{lem}\label{germlieiso}%
The map
\begin{equation*}
g(F,\,\Lie_E(K)) \to \Hom_K(\Omega,\,g(F,K)), \quad
f \mapsto (dm \mapsto dm \circ f)
\end{equation*}
is an isomorphism of $A \otimes K$-modules.\end{lem}

\pf In view of the short exact sequence
of Lemma~\ref{omegagermsesexact} the claim follows from Lemma~\ref{bndlieiso}.\quod

\begin{lem}\label{germeiso}%
The map
\begin{equation*}
g(F,\,E(K)) \to \Hom_{K\{\tau\}}(M,\,g(F,K)), \quad
f \mapsto (m \mapsto m \circ f)
\end{equation*}
is an isomorphism of $(A,A)$-bimodules.\end{lem}

\pf Follows from Lemma~\ref{bndeiso}.\quod

\begin{thm}\label{ihombvanish}%
$\RGamma(\Der\iHom(M,\,\omega\complot K)) = 0$.
\end{thm}

\pf In view of Lemma~\ref{bndlieiso}
and Proposition~\ref{derihomhoch} 
we have
\begin{equation*}
\RGamma(\Der\iHom(M,\,\omega\complot K)) =
\RH(\omega\complot\Lie_E(K)).
\end{equation*}
By assumption the action of $A$ on $\Lie_E(K)$
extends to a continuous action of $F$.
So the result follows from
Proposition~\ref{hochvanish}. \quod

%
%
%

\begin{prp}\label{ihomglie}%
The map
\begin{equation*}
\Lie_E(K) \to \Hom_{A\otimes K}(M,\,g(F,K)), \quad
\alpha \mapsto (m \mapsto (x \mapsto dm(\alpha x)))
\end{equation*}
induces a quasi-\hspace{0pt}isomorphism
$\Lie_E(K)[0] \xrightarrow{\isosign} \RGamma(\Der\iHom(M,\,g(F,K)))$.%
\end{prp}

\pf In view of Lemma~\ref{germlieiso} and Proposition~\ref{derihomhoch} 
the result is a consequence of Proposition~\ref{hochgffsq}. \quod


 
\begin{prp}\label{ihomge}%
The map
\begin{equation*}
\Lie_E(K) \to \Hom_{A\otimes K}(M,\,g(F,K)), \quad
\alpha \mapsto (m \mapsto (x \mapsto m\exp(\alpha x)))
\end{equation*}
induces a quasi-\hspace{0pt}isomorphism
$\Lie_E(K)[0] \xrightarrow{\isosign}
\RGamma(\iHom(M,\,g(F,K)))$.%
\end{prp}


\pf In view of Lemma~\ref{germeiso} and Proposition~\ref{ihomhoch} 
the result is a consequence of Proposition~\ref{hochge}. \quod


\breakflow
We state the following theorem for the sake of completeness: it will not be used in the proof
of the class number formula.

\begin{thm}\label{ihomexp}%
There exists a quasi-\hspace{0pt}isomorphism
\begin{equation*}
\RGamma(\iHom(M,\,\omega\complot K)) \xrightarrow{\isosign} C_{\exp}
\end{equation*}
with the following properties.
\begin{enumerate}
\item It is natural in $E$ and $K$.

\item It is the unique map in the derived category of
$A$-modules
such that the
square
\begin{equation*}
\xymatrix{
\RGamma(\iHom(M,\,\omega\otimes K)) \ar[r] \ar[d]^{\rtviso}_{\textup{Thm. }\ref{barweilmot}} &
\RGamma(\iHom(M,\,\omega\complot K)) \ar[d]^{\phantom{\textup{Prp. }\ref{ihomge}}} \\
E(K)[-1] \ar[r]^{\textup{identity in degree }1} & C_{\exp}
}
\end{equation*}
is commutative.

\item It makes the square
\begin{equation*}
\xymatrix{
\RGamma(\iHom(M,\,\omega\complot K)) \ar[r] \ar[d]_{\phantom{\textup{Thm. }\ref{barweilmot}}} &
\RGamma(\iHom(M,\,g(F,K))) \\
C_{\exp} \ar[r]^{\textup{identity in degree }0} & \Lie_E(K)[0] \ar[u]^{\ltviso}_{\textup{Prp. }\ref{ihomge}}
}
\end{equation*}
commutative.
\end{enumerate}
\end{thm}

\pf In view of Lemma~\ref{bndeiso} and Proposition~\ref{ihomhoch}
the result is a consequence of Theorem~\ref{hochexp}. \quod



\section{The units complex}
In this section we prove a refined variant of the formula from the previous
section. It is one of the key tools in the proof of the class number formula.

Let $K$ be a finite product of local fields and $R \subset K$ an $\Fq$-subalgebra.
We work with a Drinfeld $A$-module $E$ over $R$.
As usual the motive of $E$ is denoted $M$ and $\Omega$ stands for the quotient
$M/M^{\geqslant 1}$.
We make the following assumptions:
\begin{enumerate}
\item $R$ is discrete in $K$.

\item The action of $A$ on $\Lie_E(K)$ extends to a continuous action of $F$.
\end{enumerate}
In particular we have the exponential map
$\exp\colon\Lie_E(K) \to E(K)$ as in Definition~\ref{defexp}.

\begin{dfn}\label{defunitscomplex}%
The \emph{units complex} of $E$ is the $A$-module complex
\begin{equation*}
U_E = \Big[ \Lie_E(K) \xrightarrow{\quad\exp\quad} \frac{E(K)}{E(R)} \,\Big].
\end{equation*}\end{dfn}

\breakflow\noindent
In the following we denote $Q = K/R$. Our aim is to prove that
\begin{equation*}
\RGamma(\iHom(M,\,\omega\complot Q)) = U_E.
\end{equation*}
To imporve the legibility we will write
$\Lie_E(Q)$ for the quotient $\Lie_E(K)/\Lie_E(R)$
and $E(Q)$ for the quotient $E(K)/E(R)$.

\begin{lem}\label{codiscrexact}%
Let $V$ be a linearly topologized Hausdorff $\Fq$-vector space,
$V_0 \subset V$ a discrete subspace and $V_1 = V/V_0$.
Then the sequence
\begin{equation*}
0 \to \omega \otimes V_0 \to \omega\complot V \to \omega\complot V_1\to 0
\end{equation*}
is exact.\end{lem}

\begin{cor}\label{bnqexact}%
The sequence
$0 \to \omega\otimes R \to \omega\complot K \to \omega\complot Q \to 0$
is exact.\quod\end{cor}

\begin{cor}\label{bnqexactliee}%
The sequences
\begin{gather*}
0 \to \omega\otimes\Lie_E(R) \to \omega\complot\Lie_E(K) \to \omega\complot\Lie_E(Q) \to 0 \\
0 \to \omega\otimes E(R) \to \omega\complot E(K) \to \omega\complot E(Q) \to 0
\end{gather*}
are exact.\quod\end{cor}

\afterall\noindent\textit{Proof of Lemma~\ref{codiscrexact}.}
Let $U \subset V$ be an open $\Fq$-vector subspace such that
$V_0 \cap U = \{0\}$. It maps isomorphically onto
an open subpace $U_1 \subset V_1$. The quotients
$V/U$ and $V_1/U_1$ are discrete so that
$\omega\complot(V/U) = \omega\otimes (V/U)$ and
$\omega\complot(V_1/U_1) = \omega\otimes (V_1/U_1)$.
We then have a commutative diagram
\begin{equation*}
\xymatrix{
& & 0 \ar[d] & 0 \ar[d] \\
& 0 \ar[r]\ar[d] & \omega\complot U \ar[r]\ar[d] & \omega\complot U_1\ar[r]\ar[d] & 0 \\
0\ar[r] & \omega\otimes V_0 \ar[r]\ar[d] & \omega\complot V \ar[r]\ar[d] & \omega\complot V_1 \ar[r]\ar[d] & 0 \\
0\ar[r] & \omega\otimes V_0 \ar[r]\ar[d] & \omega\otimes V/U\ar[d] \ar[r] & \omega\otimes V_1/U_1 \ar[r]\ar[d] & 0 \\
& 0 & 0 & 0 &
}
\end{equation*}
Top and bottom rows and the left column are clearly exact.
The middle and the right columns are exact since the open
embeddings $U \to V$ and $U_1 \to V_1$ are continuously split
by Lemma~\ref{opensplit}. Hence the middle row is exact.\quod

\begin{lem}\label{bnqlieiso}%
The map $\omega\complot\Lie_E(Q) \to \Hom_R(\Omega,\,\omega\complot Q)$
defined by the commutative square
\begin{equation*}
\xymatrix{
\omega\complot\Lie_E(Q) \ar[d]_{\Res}^{\rtviso} \ar[rrr] &&&\Hom_R(\Omega,\,\omega\complot Q) \ar[d]^{\Res\circ-}_{\ltviso} \\
b(F/A,\,\Lie_E(Q)) \ar[rrr]^{f \mapsto (dm \mapsto dm \circ f)\quad} &&&\Hom_R(\Omega,\,b(F/A,Q))
}
\end{equation*}
is an isomorphism of $A\otimes R$-modules.
It is compatible with the isomorphisms 
of Lemma~\ref{bndlieiso} and of Lemma~\ref{barweillieind} for $R$ and $K$.%
\end{lem}

\pf Follows from Lemmas \ref{barweillieind} and \ref{bndlieiso}
in view of Corollary~\ref{bnqexactliee}.\quod

\begin{lem}\label{bnqeiso}%
The map $\omega\complot E(Q) \to \Hom_{R\{\tau\}}(M,\,\omega\complot Q)$
defined by the commutative square
\begin{equation*}
\xymatrix{
\omega\complot E(Q) \ar[d]_{\Res}^{\rtviso} \ar[rrr] &&&\Hom_{R\{\tau\}}(M,\,\omega\complot Q) \ar[d]^{\Res\circ-}_{\ltviso} \\
b(F/A,\,E(Q)) \ar[rrr]^{f \mapsto (m \mapsto m \circ f)\quad\quad} &&&\Hom_{R\{\tau\}}(M,\,b(F/A,Q))
}
\end{equation*}
is an isomorphism of $(A,A)$-bimodules.
It is compatible with the isomorphisms 
of Lemma~\ref{bndeiso} and of Lemma~\ref{barweilind} for $R$ and $K$.%
\end{lem}

\pf Follows from Lemmas \ref{barweilind} and \ref{bndeiso}
in view of Corollary~\ref{bnqexactliee}.\quod

\breakflow
Corollary~\ref{bnqexact} implies that the map $[\Res]\colon \omega\complot K \to g(F,K)$
factors uniquely over the natural map $\omega\complot K \to \omega\complot Q$.
We thus get a map $\omega\complot Q \to g(F,K)$.

\begin{thm}\label{derihomunit}%
There exists a unique quasi-isomorphism
\begin{equation*}
\RGamma(\Der\iHom(M, \,\omega\complot Q)) \xrightarrow{\isosign}\Lie_E(R)[0]
\end{equation*}
such that the square
\begin{equation*}
\xymatrix{
\RGamma(\Der\iHom(M,\,\omega\complot Q)) \ar[r] \ar[d]_{\phantom{\textup{Thm. }\ref{barweilmot}}} &
\RGamma(\Der\iHom(M,\,g(F,K))) \\
\Lie_E(R)[0] \ar[r] & \Lie_E(K)[0] \ar[u]^{\ltviso}_{\textup{Prp. }\ref{ihomge}}
}
\end{equation*}
is commutative.
It is natural in $R$, $K$ and $E$.%
\end{thm}

\pf 
Proposition~\ref{derihomhoch} and Lemma~\ref{bnqlieiso}
produce a natural quasi-isomorphism
\begin{equation*}
\RGamma(\Der\iHom(M,\,\omega\complot Q)) = \RH(\omega\complot\Lie_E(Q)).
\end{equation*}
Let $e\colon \Lie_E(K) \to \Lie_E(Q)$ be the natural map.
Observe that $\Lie_E(R)[0]$ is canonically the mapping fiber of $e$.
So the existence of the quasi-\hspace{0pt}isomorphism in question follows by
Theorem~\ref{hochexp} applied to the local isomorphism $e$.
The unicity part is then clear.\quod

\breakflow
Applying Theorem~\ref{barweilmot} to $R$ and $K$ we get a natural quasi-\hspace{0pt}isomorphism
\begin{equation}\label{barweilq}
\RGamma(\iHom(M,\,\omega\otimes Q)) \xrightarrow{\isosign} E(Q)[-1].
\end{equation}

\begin{thm}\label{ihomunit}%
There exists a quasi-isomorphism
\begin{equation*}
\RGamma(\iHom(M, \,\omega\complot Q)) \xrightarrow{\isosign} U_E
\end{equation*}
with the following properties:
\begin{enumerate}
\item It is natural in $R$, $K$ and $E$.

\item It is the unique map in the derived category of
$A$-modules
such that the
square
\begin{equation*}
\xymatrix{
\RGamma(\iHom(M,\,\omega\otimes Q)) \ar[r] \ar[d]^{\rtviso}_{\eqref{barweilq}} &
\RGamma(\iHom(M,\,\omega\complot Q)) \ar[d]^{\phantom{\textup{Prp. }\ref{ihomge}}} \\
E(Q)[-1] \ar[r]^{\textup{identity in degree }1} & U_E
}
\end{equation*}
is commutative.

\item It makes the square
\begin{equation*}
\xymatrix{
\RGamma(\iHom(M,\,\omega\complot Q)) \ar[r] \ar[d]_{\phantom{\eqref{barweilq}}} &
\RGamma(\iHom(M,\,g(F,K))) \\
U_E \ar[r]^{\textup{identity in degree }0} & \Lie_E(K)[0] \ar[u]^{\ltviso}_{\textup{Prp. }\ref{ihomge}}
}
\end{equation*}
commutative.%
\end{enumerate}%
\end{thm}

\pf In view of Lemma~\ref{bnqeiso} and Proposition~\ref{ihomhoch} the
result is a consequence of Theorem~\ref{hochexp} applied to
the local isomorphism $\exp\colon\Lie_E(K) \to E(Q)$.\quod

\section{Tate modules and Galois action}

As usual $A$ is a coefficient ring and $\omega$ is the module of K\"ahler
differentials of $A$ over $\Fq$. We work with a Drinfeld
$A$-module $E$ over a field $k$ containing $\Fq$. We denote $M$ the
motive of $E$. In this section we study the Hom shtuka
$\iHom(M,\omega \otimes k)$ and its relation to the Tate modules
of $E$.

The main results of this section are due to David Goss
\cite[Section 5.6]{goss}.
He stated them in terms of the motive $M$ rather than the Hom
shtuka.
We give a different argument which uses the Hom shtuka approach.

Let $k^s$ be a separable closure of $k$ and
let $G_k = \Aut(k^s/k)$ be the Galois group. 
We equip $k$ and $k^s$ with the discrete topology.
Fix a maximal ideal $\fp$ of $A$. We denote $A_\fp$ the completion of $A$ at
$\fp$ and $F_\fp$ the fraction field of $A_\fp$.
Set $\omega_\fp = \omega \otimes_A A_\fp$.
The ring $A_\fp$, the field $F_\fp$
and the module $\omega_\fp$ are assumed to carry the $\fp$-adic topologies.

We equip $\omega_\fp \complot k^s$ with an
endomorphism $\tau$ acting as the identity on $\omega_\fp$ and as the $q$-Frobenius on $k^s$.
So $\omega_\fp\complot k^s$ becomes a left
$A_\fp\complot k^s\{\tau\}$-module.

\begin{prp}\label{ihomtate}%
There exists a $G_k$-equivariant
$A_\fp$-module isomorphism
\begin{equation*}
T_\fp E \xrightarrow{\isosign}
\Hom_{A \otimes k\{\tau\}}(M,\,\omega_\fp \complot k^s)
\end{equation*}
The Galois group $G_k$ acts on the right hand side via $k^s$.%
\end{prp}


\begin{rmk}One can prove that
$\RHom_{A\otimes k\{\tau\}}(M, \,\omega_\fp\complot k^s) = T_\fp E[0]$.
We will only need the less precise result above.\end{rmk}

\afterall\noindent\textit{Proof of Proposition~\ref{ihomtate}.}
First we construct a natural isomorphism 
\begin{equation}\label{galoisequiviso}
c(F_\fp/A_\fp, \,E(k^s))
\xrightarrow{\isosign} 
\Hom_{k\{\tau\}}(M, \,\omega_\fp\complot k^s).
\end{equation}
Since $k$ is a field the motive $M$ is a free $k\{\tau\}$-module of rank $1$.
Hence the map
\begin{equation*}
c(F_\fp/A_\fp, \,\,E(k^s)) \to
\Hom_{k\{\tau\}}(M, \,c(F_\fp/A_\fp, k^s)), \quad
f \mapsto (m \mapsto m \circ f)
\end{equation*}
is an isomorphism.
Corollary~\ref{funcmoddesc} identifies $\omega_\fp\complot k^s$ with
$c(F_\fp/A_\fp, k^s)$.
Combining it with the isomorphism above
we get \eqref{galoisequiviso}.

Now \eqref{galoisequiviso} identifies 
$\Hom_{A\otimes k\{\tau\}}(M, \,\omega_\fp\complot k^s)$
with the $A_\fp$-module of 
continuous $A$-linear maps
from $F_\fp/A_\fp$ to $E(k^s)$. As $F_\fp/A_\fp$ is discrete
it follows that the latter module is
\begin{equation*}
\Hom_A(F_\fp/A_\fp, \,E(k^s)) = T_\fp E.
\end{equation*}
We get the result since
all the isomorphisms used above are $G_k$-equivariant by construction.\quod 

\begin{lem}\label{omegaplat}%
$\omega \otimes k$ is an $A\otimes k$-lattice in the $A_\fp\complot k^s$-module
$\omega_\fp \complot k^s$.\end{lem}

\pf Follows from Lemma~\ref{compllattice}.\quod

\breakflow
As $\Lie_E(k)$ is a one-dimensional $k$-vector space
the action of $A$ on $\Lie_E(k)$ determines a homomorphism
$A \to k$. The kernel of this homomorphism is called the
\emph{characteristic}\index{idx}{Drinfeld module!characteristic} of $E$.

\begin{lem}\label{ihomibiject}%
If $\fp$ is different from the characteristic of $E$
then
\begin{equation*}
\RGamma(\Der\iHom(M, \,\omega_\fp\complot k^s) = 0.
\end{equation*}%
\end{lem}

\pf Lemma~\ref{omegaplat} implies that 
\begin{equation*}
\Der\iHom(M, \,\omega \otimes k) \otimes_{A \otimes k} (A_\fp\complot k^s) =
\Der\iHom(M, \,\omega_\fp\complot k^s).
\end{equation*}
At the same time Theorem~\ref{barweilliemot} provides a
natural quasi-\hspace{0pt}isomorphism
\begin{equation*}
\RGamma(\Der\iHom(M,\,\omega\otimes k)) = \Lie_E(k)[-1].
\end{equation*}
Since $A_\fp\complot k^s$ is flat over $A\otimes k$ it is 
enough to prove that
\begin{equation*}
\Lie_E(k) \otimes_{A\otimes k}(A_\fp\complot k^s) = 0.
\end{equation*}
%
Now $\fp$ is different from the characteristic of $E$ so
there exists an element of $\fp$ which does not act on $\Lie_E(k)$ by zero.
It then acts by an automorphism and we conclude that $\Lie_E(k) \otimes_A A/\fp = 0$.
At the same time $A_\fp \complot k^s$ is the completion of the noetherian ring
$A \otimes k^s$ with respect to the ideal $\fp \otimes k^s$.
Since $\Lie_E(k)$ is a finitely generated $A\otimes k$-module Nakayama's lemma
shows that $\Lie_E(k) \otimes_{A\otimes k} (A_\fp\complot k^s)$ is zero.\quod



\breakflow
We will also need a lemma from Dieudonn\'e-Manin theory.

\begin{lem}\label{dmlattice}%
Let $N$ be a left $A_\fp\complot k^s\{\tau\}$-module which is
finitely generated free as an $A_\fp\complot k^s$-module.
If the adjoint $\tau^\ast N \to N$ of the $\tau$-multiplication map
is bijective then the natural map
$N^{\tau=1}\otimes_{A_\fp} (A_\fp\complot k^s) \to N$ is an isomorphism.\end{lem}

\pf See \cite[Proposition 4.4]{stalder}.
\quod

\begin{prp}\label{ihomtatelat}%
If $\fp$ is different from the characteristic of $E$ then
the natural map
\begin{equation*}
\Hom_{A \otimes k\{\tau\}}(M,\,\omega_\fp \complot k^s) \otimes_{A_\fp} (A_\fp \complot k^s) \to
\Hom_{A \otimes k}(M,\,\omega_\fp \complot k^s)
\end{equation*}
is an isomorphism.\end{prp}

\pf Consider the shtuka
\begin{equation*}
\iHom(M,\,\omega_\fp\complot k^s) = 
\Big[
\Hom_{A\otimes k}(M,\,\omega_\fp\complot k^s) \shtuka{\,\,i\,\,}{j}
\Hom_{A\otimes k}(\tau^\ast M,\,\omega_\fp\complot k^s)
\Big].
\end{equation*}
Lemma \ref{ihomibiject} implies that
the arrow $i$ is bijective.
We equip the $A_\fp\complot k^s$-module
\begin{equation*}
H = \Hom_{A\otimes k}(M, \,\omega_\fp\complot k^s)
\end{equation*}
with the $\tau$-linear endomorphism $i^{-1} j$.

The adjoint $\tau^\ast(\omega_\fp\complot k^s) \to \omega_\fp\complot k^s$
of the $\tau$-multiplication map
is an isomorphism. It follows from the definition of $\iHom$
that the $\tau$-adjoint of the arrow $j$ above is an isomorphism.
We can thus apply Lemma~\ref{dmlattice} and conclude that the natural
map
\begin{equation*}
H^{i^{-1} j = 1} \otimes_{A_\fp} (A_\fp\complot k^s) \to H
\end{equation*}
is an isomorphism.
Proposition~\ref{ihomzerocoh} identifies the module
$\Hom_{A\otimes k\{\tau\}}(M, \,\omega_\fp\complot k^s)$
with $\uH^0(\iHom(M,\,\omega_\fp\complot k^s))$
which is equal to $H^{i^{-1} j = 1}$ by Proposition~\ref{shtcohzerodesc}.\quod
%



\begin{prp}\label{ihomdetarith}%
Assume that $k$ is a finite extension of $\Fq$ of degree $d$.
The shtuka
\begin{equation*}
\iHom(M,\,\omega \otimes k) =
\Big[
\Hom_{A\otimes k}(M, \,\omega \otimes k) \shtuka{\,\,i\,\,}{j}
\Hom_{A\otimes k}(\tau^\ast M, \,\omega\otimes k) \Big]
\end{equation*}
has the following properties.
\begin{enumerate}
\item The arrow $j$ is a bijection.

\item 
Let
$\sigma \in G_k$ be the arithmetic Frobenius element.
For every prime $\fp$ of $A$ different from the characteristic of $E$
we have an identity
\begin{equation*}
\det\nolimits_{A_\fp}\!\big(T - \sigma \,\big|\, T_\fp E\big) =
\det\nolimits_{A \otimes k}\big(T - (j^{-1} i)^d \,\big|\, \Hom_{A\otimes k}(M, \,\omega \otimes k)\big)
\end{equation*}
in $A_\fp \otimes k[T]$.
\end{enumerate}
\end{prp}

\pf 
(1) The endomorphism $\tau$ of $\omega \otimes k$ is an automorphism
since $k$ is a finite extension of $\Fq$. From the definition of $\iHom$
it follows that the arrow $j$ is bijective.

(2) Given an $A \otimes k$-module $N$ we will
denote $\psi$ the endomorphism of the module $\Hom_{A\otimes k}(M, N)$
which sends a map $f$ to the map $m \mapsto f(\tau^d m)$.

Proposition \ref{ihomtate} identifies $T_\fp E$ 
and 
$\Hom_{A\otimes k\{\tau\}}(M, \,\omega_\fp\complot k^s)$ on
which $G_k$ acts via $k^s$.
So if $f$ is an element of the latter module then
$\sigma(f) = \tau^d \circ f$.
As $f$ is a homomorphism of left $A \otimes k\{\tau\}$-modules
it follows that $\tau^d \circ f = \psi(f)$.
We conclude that
\begin{equation*}
\det\nolimits_{A_\fp}\!\big(T - \sigma \,\big|\, T_\fp E\big) =
\det\nolimits_{A_\fp}\!\big(T - \psi \,\big|\, \Hom_{A\otimes k\{\tau\}}(M, \,\omega_\fp\complot k^s)\big).
\end{equation*}
Next, Proposition \ref{ihomtatelat} implies that
the latter polynomial is equal to
\begin{equation*}
\det\nolimits_{A_\fp\complot k^s}\big( T - \psi \,\big|\, \Hom_{A\otimes k}(M, \,\omega_\fp\complot k^s)\big).
\end{equation*}
Lemma~\ref{omegaplat} identifies this polynomial with 
\begin{equation*}
\det\nolimits_{A\otimes k}\big( T - \psi \,\big|\, \Hom_{A\otimes k}(M, \,\omega\otimes k)\big).
\end{equation*}
Now Proposition~\ref{ihommotdesc} 
implies that the endomorphism
$j^{-1} i$ of $\Hom_{A\otimes k}(M, \,\omega \otimes k)$
sends a map
$f$
to the map $m \mapsto \tau^{-1}(f(\tau m))$. As $\tau^{-d}$ is the identity automorphism
of $k$ we conclude that $(j^{-1} i)^d(f) = \psi(f)$ and the result follows.\quod

\breakflow
Let $\nfA$ be the fraction field of $A$ and let
$\nOA = \omega \otimes_A \nfA$ be the module of K\"ahler differentials
of $\nfA$ over $\Fq$.
The main result of this section is the following theorem.

\begin{thm}\label{ihomdet}%
Assume that $k$ is a finite extension of $\Fq$ of degree $d$.
The $\nfA \otimes k$-module shtuka
\begin{equation*}
\iHom(M, \,\nOA \otimes k) =
\Big[
\Hom_{A\otimes k}(M, \,\nOA \otimes k) \shtuka{\,\,i\,\,}{j}
\Hom_{A\otimes k}(\tau^\ast M, \,\nOA \otimes k) \Big]
\end{equation*}
has the following properties.
\begin{enumerate}
\item The arrow $i$ is an isomorphism.

\item
Let $\sigma^{-1} \in G_k$ be the geometric Frobenius element.
For every prime $\fp$ of $A$ different from the characteristic of $E$
we have an identity
\begin{equation*}
\det\nolimits_{A_\fp}\!\big(1 - T^d \sigma^{-1} \,\big|\, T_\fp E\big) =
\det\nolimits_{\nfA}\big(1 - T (i^{-1} j) \,\big|\, \Hom_{A\otimes k}(M, \,\nOA \otimes k)\big)
\end{equation*}
in $F_\fp[T]$. 
\end{enumerate}
\end{thm}

\begin{rmk}%
In particular
the polynomial $\det\nolimits_{A_\fp}(1 - T\sigma^{-1} \mid T_\fp E)$
is independent of the choice of $\fp$ and
has
coefficients in the subring of $F^\sharp$ consisting of elements
which are integral away from the characteristic of $E$ and
$\infty$.
Drinfeld's construction (Theorem~\ref{motivec}) implies that the coefficients
are also integral at $\infty$.\end{rmk}

\afterall\noindent
\textit{Proof of Theorem \ref{ihomdet}.}
(1) By Theorem~\ref{barweilliemot} we have
$\RGamma(\Der\iHom(M,\,\omega\otimes k) = \Lie_E(k)[-1]$.
It is thus enough to prove that $\Lie_E(k) \otimes_A \nfA = 0$.
However $\Lie_E(k)$ is a torsion $A$-module since $k$ is a finite
extension of $\Fq$.

(2) We have
\begin{equation*}
\det\nolimits_{A_\fp}\!\big(1 - T^d \sigma^{-1} \,\big|\, T_\fp E\big) =
\frac{\det\nolimits_{A_\fp}\!\big(T^d - \sigma \,\big|\, T_\fp E\big)}
{\det\nolimits_{A_\fp}\!\big( -\sigma \,\big|\, T_\fp E\big)}.
\end{equation*}
Proposition~\ref{ihomdetarith} identifies the latter fraction with
\begin{equation*}
\frac{\det\nolimits_{A \otimes k}\big(T^d - (j^{-1} i)^d \,\big|\, \Hom_{A\otimes k}(M, \,\omega \otimes k)\big)}
{\det\nolimits_{A \otimes k}\big( -(j^{-1} i)^d \,\big|\, \Hom_{A\otimes k}(M, \,\omega \otimes k)\big)}
\end{equation*}
which is in turn equal to
\begin{equation*}
\det\nolimits_{\nfA \otimes k}\big(1 - T^d (i^{-1} j)^d \,\big|\, \Hom_{A\otimes k}(M, \,\nOA \otimes k)\big).
\end{equation*}
Now \cite[Lemma 8.1.4]{bp} shows that this polynomial coincides with
\begin{equation*}
\det\nolimits_{\nfA}\big(1 - T(i^{-1} j) \,\big|\, \Hom_{A\otimes k}(M, \,\nOA \otimes k)\big).
\end{equation*}
and the result follows.\quod

\chapter{Local models}
\label{chapter:locmod}
\label{ch:refl}
\label{ch:locmod}
\label{ch:locmodcoh}
\label{ch:locmoddrinfeld}

Fix a coefficient ring $A$ as in Definition \lref{coeffrings}{defcoeffring}
and let $F$ be the local field of $A$ at infinity.
We denote $\cO_F\subset F$ the ring of integers and
$\fm_F \subset \cO_F$ the maximal ideal.
Let $K$ be a finite product of local fields containing $\Fq$.
As usual $\cO_K \subset K$ stands for the ring of integers
and $\fm_K \subset\cO_K$ denotes the Jacobson radical.
The $\tau$-ring and $\tau$-module structures used in this chapter are
as described in Section \lref{ottaustruct}{sec:ottaustruct}.

We study Drinfeld modules in a local situation.
Namely we work with a Drinfeld $A$-module $E$ over $K$ by the means of
the $F\indcot K$-module shtuka $\iHom(M,\,a(F,K))$. Here $M$ is the motive of $E$
and $a(F,K)$ is the space of locally constant bounded $\Fq$-linear maps from $F$
to $K$ as in Definition \lref{bndlc}{defbndlc}.
We introduce the notion of a local model which is an
$\cO_F\indcot\cO_K$-module subshtuka $\cM \subset \iHom(M,\,a(F,K))$ with certain
properties. Informally speaking, $\cM$ compactifies $\iHom(M,\,a(F,K))$ in the
direction of the coefficients $\cO_F \subset F$ and the base $\cO_K \subset K$.
One important result of this chapter is Theorem \ref{basecompexist} which
implies that local models exist.
Another important result is Theorem \ref{locmodell} which states that a
local model is an elliptic shtuka in the sense of Chapter \ref{chapter:reg}.

The constructions in this chapter are algebraic in nature and
the topology on the various tensor product rings plays no essential role.
Even though the (abstract) rings $\cO_F\complot\cO_K$ and $\cO_F\indcot\cO_K$
are the same we use the latter notation to stress that it is a subring
of $F\indcot K$. 

\section{Lattices}

Let $R_0 \to R$ be a homomorphism of rings and $M$ an $R$-module.
Recall that an $R_0$-submodule $M_0 \subset M$ is called a lattice
if the natural map $R \otimes_{R_0} M_0 \to M$ is an isomorphism
(Section~\ref{sec:lattices} in the chapter ``\chnotconv'').

\begin{dfn}\index{idx}{lattice}%
Let $R_0 \to R$ be a homomorphism of $\tau$-rings, $M$ an $R$-module shtuka and
$M_0 \subset M$ an $R_0$-module subshtuka. We say that $M_0$ is an
$R_0$-lattice in $M$ if the underlying modules of $M_0$ are $R_0$-lattices in
the underlying modules of $M$.\end{dfn}

\section{Reflexive sheaves}

The aim of this section is to
review some properties of reflexive sheaves on a scheme such as
$Y = \Spec k[[z,\zeta]]$ with $k$ a field and the open subscheme $U \subset Y$
which is the complement of the closed point.
The main result states that every locally free sheaf on $U$ extends uniquely to
a locally free sheaf on $Y$.
While the contents of the section is widely known, it does not seem to
appear in the literature in the form which we need.

Let $F$ be a local field containing $\Fq$ and let $K$ be a finite product of
local fields containing $\Fq$. We denote $\cO_F \subset F$ and $\cO_K \subset K$
the corresponding rings of integers.
For the rest of this section let us fix uniformizers $z \in \cO_F$ and
$\zeta\in\cO_K$. Let $Y = \Spec(\cO_F \indcot \cO_K)$ and let $U = D(z) \cup
D(\zeta) \subset Y$.

\begin{lem}\label{indcotstruct}
The ring $\cO_F \indcot \cO_K$ has the following properties:
\begin{enumerate}
\item $\cO_F \indcot \cO_K$ is a finite product of complete regular
2-dimensional local rings.

\item The maximal ideals of $\cO_F \indcot \cO_K$ are precisely the prime ideals
containing $z$ and $\zeta$.
\end{enumerate}\end{lem}

\pf By Proposition \lref{ringot}{indcotcomplot} the natural map
$\cO_F\indcot\cO_K \to \cO_F\complot\cO_K$ is an isomorphism. So the result
follows at once from Proposition \lref{otalgprops}{complotstruct}. \quod

\breakflow
Recall that a coherent sheaf $\cF$ is called \emph{reflexive} if the natural map
from $\cF$ to its double dual $\cF^{**}$ is an isomorphism. A locally free sheaf
is automatically reflexive. We use \cite{hartrefl} as the reference for the
theory of reflexive sheaves.

\begin{lem}\label{reflhartogs} If $\cF$ is a reflexive sheaf on
$Y$ then the natural map $\Gamma(Y,\cF) \to \Gamma(U,\cF)$
is an isomorphism.\end{lem}

\pf Lemma \ref{indcotstruct} implies the following
two statements.
First, $Y$ is a finite disjoint union of regular schemes. Second, all the prime
ideals in the complement of $U$ are of height $2$. Therefore the result follows
from \cite{hartrefl}, Proposition 1.6 (i) $\Rightarrow$ (iii). \quod

\breakflow
Let $\iota\colon U \to Y$ be the open embedding.

\begin{lem}\label{reflglobgen} Every coherent sheaf on $U$ is globally
generated.\end{lem}

\pf 
Let $\cF$ be a coherent sheaf on $U$.  The morphism $\iota$ is quasi-compact
quasi-separated so $\iota_\ast\cF$ is quasi-coherent. Let $f \in \cF(D(z))$. As
$\iota_\ast\cF$ is quasi-coherent and $Y$ is affine there exists an $n \gg 0$
such that $f z^n$ lifts to a global section of $\iota_\ast\cF$ or equivalently,
to a global section of $\cF$.  The same argument applies to $\cF(D(\zeta))$. We
can therefore lift all the generators of $\cF(D(z))$ and $\cF(D(\zeta))$ to
global sections of $\cF$.
\quod

\begin{lem}\label{reflext} If $\cF$ is a reflexive sheaf on $U$ then $\iota_\ast
\cF$ is reflexive.\end{lem}

\pf First we prove that $\iota_\ast\cF$ is coherent.
The sheaf $\cF^*$ is coherent. By Lemma \ref{reflglobgen} there is a surjection
$\cO_U^n \to \cF^*$ for $n \gg 0$. Dualizing it and taking the
composition with the natural isomorphism $\cF \cong \cF^{**}$ we obtain an
embedding $\cF \hookrightarrow \cO_U^n$. The sheaf $\cO_Y$ is reflexive so Lemma
\ref{reflhartogs} implies that $\iota_\ast\cO_U^n = \cO_Y^n$. Therefore
$\iota_\ast\cF$ embeds into $\cO_Y^n$.

Next we have a natural commutative diagram
\begin{equation}\label{reflsquare}
\vcenter{\vbox{\xymatrix{
\Gamma(Y,\iota_\ast\cF) \ar[r] \ar[d] & \Gamma(Y,(\iota_\ast\cF)^{**}) \ar[d] \\
\Gamma(U,\cF) \ar[r] & \Gamma(U,\cF^{**}).
}}}
\end{equation}
Here the horizontal arrows are induced by the natural maps $\cF \to \cF^{**}$
and $\iota_\ast\cF \to (\iota_\ast\cF)^{**}$, the left vertical arrow is the
restriction while the right vertical arrow is
the composition
\begin{equation*}
\Gamma(Y, (\iota_\ast\cF)^{**}) \to \Gamma(U, (\iota_\ast\cF)^{**}) \to
\Gamma(U, \cF^{**})
\end{equation*}
of the restriction map and the map which results from
applying $\Gamma(U,{-}^{**})$ to the adjunction morphism $\iota^\ast\iota_\ast\cF \to
\cF$.

Let us prove that the right vertical arrow in \eqref{reflsquare} is an isomorphism.
By \cite{hartrefl}, Corollary 1.2 the sheaf $(\iota_\ast\cF)^{**}$ is reflexive.
So the restriction map $\Gamma(Y,(\iota_\ast\cF)^{**}) \to
\Gamma(U,(\iota_\ast\cF)^{**})$ is an isomorphism by Lemma \ref{reflhartogs}.
The adjunction map $\iota^\ast\iota_\ast\cF \to \cF$ is an isomorphism since
$\iota$ is an open embedding. So our claim follows.

The right vertical arrow in \eqref{reflsquare} is an isomorphism by
construction while the bottom horizontal arrow is an isomorphism since $\cF$ is
reflexive. It follows that the top horizontal arrow is an isomorphism whence
$\iota_\ast\cF$ is reflexive. \quod

\begin{lem}\label{reflfree}
Every reflexive sheaf on $Y = \Spec\cO_F\indcot\cO_K$ is locally free.
\end{lem}

\pf According to
Lemma \ref{indcotstruct}
the ring $\cO_F \indcot \cO_K$ is a finite product of regular local
2-dimensional rings.  So the result follows from \cite{hartrefl}, Corollary 1.4.
\quod

\section{Reflexive sheaves and lattices}

Using the results of the previous section we prove several technical lemmas on
modules over the ring $\cO_F\indcot\cO_K$.


\begin{lem}\label{locmodrings} Let $z \in\cO_F$ and $\zeta\in\cO_K$ be
uniformizers.
\begin{enumerate}
\item The natural map $\Spec(F \indcot \cO_K) \to \Spec(\cO_F\indcot\cO_K)$ is
an open embedding with image $D(z)$.

\item The natural map $\Spec(\cO_F\indcot K) \to \Spec(\cO_F\indcot\cO_K)$ is an
open embedding with image $D(\zeta)$.

\item The natural map $\Spec(F\indcot K) \to \Spec(\cO_F\indcot\cO_K)$ is an
open embedding with image $D(z) \cap D(\zeta)$.
\end{enumerate}\end{lem}

\pf
Proposition \lref{otalgprops}{compindcotloc} shows that
\begin{align*}
(\cO_F \indcot \cO_K)[z^{-1}] &= F \indcot \cO_K, \\
(\cO_F \indcot \cO_K)[\zeta^{-1}] &= \cO_F \indcot K
\end{align*}
where $z \in \cO_F$ and $\zeta \in \cO_K$ are uniformizers.
Furthermore
Proposition \lref{otalgprops}{indcotloc} tells us that
\begin{equation*}
(\cO_F \indcot \cO_K)[(z\zeta)^{-1}] = F \indcot K.
\end{equation*}
So the result follows. \quod

\begin{lem}\label{ofoklathartogs}
Let $N$ be a locally free $\cO_F\indcot \cO_K$-module.
Consider the modules
\begin{equation*}
N^c = 
N \otimes_{\cO_F\indcot\cO_K}(\cO_F\indcot K), \quad
N^b = 
N \otimes_{\cO_F\indcot\cO_K}(F\indcot \cO_K)
\end{equation*}
We have $N = N^c \cap N^b$ as submodules of
$N \otimes_{\cO_F\indcot\cO_K}(F \indcot K)$.\quod\end{lem}

%

\begin{lem}\label{ofoklocfree} Let $N$ be a locally free $F \indcot K$-module of
finite rank. If $N^c \subset N$ is a locally free $\cO_F \indcot K$-lattice
and $N^b \subset N$ a locally free $F \indcot \cO_K$-lattice then $N^c \cap N^b$
is a locally free $\cO_F \indcot \cO_K$-lattice in $N^c$ and $N^b$.\end{lem}

\pf Let $z \in \cO_F$ and $\zeta\in\cO_K$ be uniformizers.
By Lemma \ref{locmodrings}
$\Spec(F \indcot \cO_K) = D(z)$ and $\Spec(\cO_F \indcot K) =
D(\zeta)$ are open subschemes in $\Spec(\cO_F \indcot \cO_K)$ 
whose intersection is $\Spec(F \indcot K)$.
By assumption $N^c$ and $N^b$ restrict to the same module $N$ on the
intersection $D(z) \cap D(\zeta) = \Spec(F \indcot K)$. Hence they
define a locally free sheaf $\cN$ on $D(z) \cup D(\zeta)$. Let $\iota\colon D(z)
\cup D(\zeta) \to \Spec(\cO_F\indcot\cO_K)$ be the embedding map. Lemmas \ref{reflext} and
\ref{reflfree} imply that $\iota_\ast\cN$ is a locally free sheaf on
$\Spec(\cO_F \indcot \cO_K)$. By construction
\begin{equation*}
\Gamma(\Spec(\cO_F \indcot \cO_K), \iota_\ast \cN) = N^c \cap N^b
\end{equation*}
so the result follows. \quod

\section{Hom shtukas}

Fix a Drinfeld $A$-module $E$ over $K$.
The motive $M = \Hom(E,\bG_a)$ of $E$
carries a natural structure of a left $A\otimes K\{\tau\}$-module with
$A$ acting on the right via $E$ and $K\{\tau\}$ acting on the left via $\bG_a$.
%
%
Let $F$ be the local field of $A$ at infinity.
In this section we list some properties of the shtuka $\iHom_{A\otimes K}(M,N)$
where $N$ is one of the function spaces
\begin{equation*}
a(F,K), \quad
b(F,K), \quad
a(F/A,K), \quad
b(F/A,K)
\end{equation*}
or the germ space $g(F,K)$.
These are immediate consequences of results in Chapters \ref{chapter:tvs} and \ref{chapter:trm}.
The $\tau$-module structures on these spaces are as described
in Section \lref{ottaustruct}{sec:ottaustruct}.


\begin{lem}%
$\iHom(M,\,a(F,K))$ is a locally free $F \indcot K$-module shtuka and
$\iHom(M,\,b(F,K))$ is a locally free $F^\#\complot K$-module shtuka.\end{lem}

\pf The $A\otimes K$-module $M$ is locally free of finite rank by definition.
Corollary~\ref{funcmodstruct} tells us that the function
space $a(F,K)$ is a free $F\indcot K$-module of rank $1$ while $b(F,K)$ is a
free $F^\#\complot K$-module of rank $1$. The result now follows
from the definition of $\iHom$.\quod

\begin{lem}\label{homablat}%
$\iHom(M,\,a(F,K))$ is an $F\indcot K$-lattice
in $\iHom(M,\,b(F,K))$.
\end{lem}

%
\pf
According to Corollary~\ref{funclat}
the function space $a(F,K)$ is an
$F\indcot K$-lattice in $b(F,K)$.
The result now follows from the definition of $\iHom$ since
$M$ is a locally free $A\otimes K$-module of finite rank. \quod

\begin{lem}\label{ihomafaklat}%
The shtuka $\iHom(M,\,a(F/A,K))$ has the following properties:
\begin{enumerate}
\item It is a locally free $A\otimes K$-module shtuka.

\item It is an $A\otimes K$-lattice in the
$F\indcot K$-module shtuka $\iHom(M,a(F,K))$.
\end{enumerate}\end{lem}

\pf By definition $M$ is a locally free $A\otimes K$-module.
By Corollaries \ref{funcmodstruct} and \ref{funclat}
the space $a(F/A,K)$ is a locally free $A\otimes K$-lattice
in the $F\indcot K$-module $a(F,K)$.
So (1) and (2) follow from the definition of $\iHom$. \quod

\begin{lem}\label{ihombfaklat}%
The shtuka $\iHom(M,\,b(F/A,K))$ has the following properties:
\begin{enumerate}
\item It is a locally free $A\complot K$-module shtuka.

\item It is an $A\complot K$-lattice in the
$F^\#\complot K$-module shtuka $\iHom(M,b(F,K))$.
\end{enumerate}\end{lem}

\pf By definition $M$ is a locally free $A\otimes K$-module.
By Corollaries \ref{funcmodstruct} and \ref{funclat}
the space $b(F/A,K)$ is a locally free $A\complot K$-lattice
in the $F^\#\complot K$-module $b(F,K)$.
So (1) and (2) follow from the definition of $\iHom$. \quod

\begin{lem}\label{ihomafakblat}%
The shtuka $\iHom(M,\,a(F/A,K))$ is a locally free
$A\otimes K$-lattice in the $A\complot K$-module shtuka
$\iHom(M,\,b(F/A,K))$.\end{lem}

\pf Follows from Corollary~\ref{funclat}.\quod

%

\section{The notion of a local model}
\label{sec:deflocmod}

%

As before we work with a fixed Drinfeld module $E$ over $K$.
We denote $M = \Hom(E,\bG_a)$ the motive of $E$ as in Definition \lref{drmot}{defmotive}.
Throughout the rest of the chapter we assume that the action of $A$ on $\Lie_E(K)$ extends
to a continuous action of $F$, the local field of $A$ at infinity.

The shtuka $\iHom(M,a(F,K))$ is a locally free $F\indcot K$-module shtuka.
In this chapter we will study models of this shtuka over various subrings
of $F\indcot K$. In particular we will introduce the notions of
\begin{enumerate}
\item the \emph{coefficient compactification} (a model over $\cO_F\indcot K$),

\item a \emph{base compactification} (a model over $F \indcot \cO_K$),

\item a \emph{local model} (over $\cO_F\indcot\cO_K$).
\end{enumerate}

Let $C$ be the projective compactification of $\Spec A$.
We denote $\cE$
the locally free shtuka on $C \times \Spec K$ constructed in Theorem \lref{drinconstr}{motivec}.

\begin{dfn}\label{locmoddrincomp}\index{idx}{compactification!coefficient compactification}%
The \emph{coefficient compactification}
$\cM^c \subset \iHom(M,\,a(F,K))$ is the
$\cO_F\indcot K$-subshtuka
\begin{equation*}
\iHom_{\cO_F\indcot K}\big(\cE(\cO_F\indcot K), \,a(F/\cO_F,K)\big).
\end{equation*}
The superscript ``c'' stands for ``coefficients''.%
\end{dfn}

\begin{lem}%
$\cM^c \subset \iHom(M,\,a(F,K))$ is a locally free
$\cO_F \indcot K$-lattice.\end{lem}

\pf
By construction $\cE(A \otimes K)$ is the shtuka given by the left $A\otimes K\{\tau\}$-module $M$.
Hence $\cE(\cO_F\indcot K)$ is a locally free
$\cO_F\indcot K$-lattice in the pullback of $M$ to $F\indcot K$.
According to Corollaries \ref{funcmodstruct} and \ref{funclat}
the function space
$a(F/\cO_F,K)$ is a locally free $\cO_F \indcot K$-lattice
in $a(F,K)$. The result now follows from the
definition of $\iHom$.\quod

\begin{lem}\label{drincompnilp}%
$\cM^c(\cO_F/\fm_F \otimes K)$ is nilpotent.\end{lem}

\pf 
By Theorem~\lref{drinconstr}{motivec} the shtuka
$\cE(\cO_F/\fm_F \otimes K)$ is co-nilpotent so the result follows by
Proposition~\ref{conilp}. \quod

\breakflow
By our assumption the action of $A$ on $\Lie_E(K)$ extends to a continuous
action of $F$ so that we have a continuous homomorphism $F \to K$. Such a
homomorphism necessarily maps $\cO_F$ to $\cO_K$.

\begin{dfn}\label{deflocmodcond}\index{idx}{ramification ideal!local model@of a local model}\index{nidx}{$\rami$, ramification ideal}%
We define the \emph{ramification ideal} $\rami \subset \cO_K$ to be the ideal
generated by $\fm_F$ in $\cO_K$.\end{dfn}

\begin{lem} The ideal $\rami$ is open and is contained in the Jacobson
radical $\fm_K$.\quod\end{lem}

\begin{dfn}\label{defbcomp}\index{idx}{compactification!base compactification}%
A \emph{base compactification} of $\iHom(M,\,a(F,K))$
is a locally free $F \indcot \cO_K$-lattice
$\cM^b$ 
such that
$\cM^b(F \otimes \cO_K/\fm_K)$ is nilpotent and
$\cM^b(F \otimes \cO_K/\rami)$ is linear.
The superscript ``b'' stands for ``base''.\end{dfn}


\begin{dfn}\label{deflocmod}\index{idx}{shtuka model!local}%
A \emph{local model} of $\iHom(M,\,a(F,K))$
is a locally free $\cO_F\indcot\cO_K$-lattice $\cM$
such that:
\begin{enumerate}
%
\item The subshtuka
$\cM(\cO_F\indcot K) \subset \iHom(M,a(F,K))$ 
coincides with $\cM^c$.

\item $\cM(F \indcot \cO_K/\fm_K)$ is nilpotent
and $\cM(F\indcot\cO_K/\rami)$ is linear.
\end{enumerate}\end{dfn}

\begin{prp}\label{locmodintersect}%
Let $\cM \subset \iHom(M,\,a(F,K))$ be an $\cO_F\indcot\cO_K$-subshtuka.
The following are equivalent:
\begin{enumerate}
\item $\cM$ is a local model.

\item The subshtuka $\cM(F \indcot \cO_K) \subset \iHom(M,\,a(F,K))$ is a base
compactification and
$\cM = \cM^c \cap \cM(F \indcot \cO_K)$
\end{enumerate}\end{prp}

\pf[\aftereq](1) $\Rightarrow$ (2).
It follows directly from the definition of a local model that
$\cM(F\indcot\cO_K)$ is a base compactification
and $\cM^c = \cM(\cO_F\indcot K)$.
Lemma \ref{ofoklathartogs} implies that
$\cM = \cM^c \cap \cM(F\indcot\cO_K)$.
(2) $\Rightarrow$ (1). Lemma \ref{ofoklocfree} shows that $\cM$ is a locally
free $\cO_F\indcot\cO_K$-lattice in $\cM^c$ so (1) follows. \quod

%

\begin{lem}\label{locmodblat}%
If $\cM \subset \iHom(M,\,a(F,K))$ is a local model then the natural map
$\cM(F^\#\complot\cO_K) \to \iHom(M,\,b(F,K))$ is an inclusion of an
$F^\#\complot\cO_K$-lattice.\end{lem}

\pf $\cM(F\indcot K) = \iHom(M,\,a(F,K))$ by definition
and $\iHom(M,\,a(F,K))$ is an $F\indcot K$-lattice in
$\iHom(M,\,b(F,K))$ by Lemma \ref{homablat}.
Therfore the natural map $\cM(F^\#\complot K) \to \iHom(M,\,b(F,K))$ is an isomorphism.
According to Proposition~\lref{otalgprops}{discrcomplotloc}
the ring $F^\#\complot K$ is a localization
of $F^\#\complot\cO_K$ at a uniformizer of $\cO_K$.
As $\cM$ is locally free it follows that
$\cM(F^\#\complot\cO_K) \to \cM(F^\#\complot K)$ is a lattice inclusion.
\quod

%
%

\breakflow
The natural map $\cO_F\indcot\cO_K \to \cO_F\complot\cO_K$
is an isomorphism by Proposition~\lref{ringot}{indcotcomplot}.
So we can view an $\cO_F\indcot\cO_K$-module
shtuka $\cM$ as an $\cO_F\complot\cO_K$-module shtuka. In particular
the constructions of Chapter~\ref{chapter:reg} apply to shtukas
on $\cO_F\indcot\cO_K$.

Let us recall the twisting construction of Section~\lref{reg}{sec:elltwist}.
Let $I \subset \cO_F\indcot\cO_K$ be a $\tau$-invariant ideal. Given
an $\cO_F\indcot\cO_K$-module shtuka 
\begin{equation*}
\cM = \Big[ M_0 \shtuka{i}{j} M_1 \Big]
\end{equation*}
we define
\begin{equation*}
I \cM = \Big[ I M_0 \shtuka{i}{j} I M_1 \Big].
\end{equation*}
The fact that $I$ is an invariant ideal guarantees that the diagram on the right
indeed defines a shtuka. We will use twists by the invariant ideal
$\cO_F \indcot \rami$. To improve legibility we will write
$\rami \cM$ in place of $(\cO_F\indcot\rami) \cM$.

\begin{prp}\label{locmodtwist}%
If $\cM$ is a local model then $\rami\cM$ is a local model.\end{prp}

\pf Indeed $(\rami\cM)(\cO_F \indcot K) = \cM(\cO_F\indcot K)$ so
$(\rami\cM)(\cO_F\indcot K)$ coincides with the coefficient compactification
of $\iHom(M,\,a(F,K))$. Lemma \lref{reg}{ellnilptwist} shows that
$(\rami\cM)(F \otimes \cO_K/\fm_K)$ is nilpotent while
Proposition \lref{reg}{elltwistvanish} implies that
$(\rami\cM)(F \otimes \cO_K/\rami)$ is linear. Since
$\rami\cM$ is a locally free $\cO_F\indcot\cO_K$-lattice in
$\iHom(M,\,a(F,K))$ by construction we conclude that
$\rami\cM$ is a local model. \quod

\begin{prp}\label{locmodmap}%
If $\cM$, $\cN$ are local models then there exists
$n \geqslant 0$ such that $\rami^n \cM \subset \cN$.\end{prp}

\pf By Proposition \ref{locmodintersect} we have
\begin{equation*}
\cM = \cM^c \cap \cM(F\indcot\cO_K), \quad
\cN = \cM^c \cap \cN(F\indcot\cO_K)
\end{equation*}
where $\cM^c \subset \iHom(M,\,a(F,K))$ is the coefficient compactification.
So to prove the proposition it is enough to find an integer $n \geqslant 0$
such that $(\rami^n\cM)(F\indcot\cO_K) \subset \cN(F\indcot\cO_K)$.
Observe that the ramification ideal ideal $\rami$ is contained in the Jacobson radical
$\fm_K$ by construction. In particular it contains a power of a uniformizer of
$\cO_K$. The result follows since
$\cM(F\indcot\cO_K)$ and $\cN(F\indcot\cO_K)$ are
$F\indcot\cO_K$-lattices in $\iHom(M,\,a(F,K))$. \quod

\section{Existence of base compactifications}

We keep the assumptions and the notation of the previous section.
In this section we prove that the $F\indcot K$-module shtuka
$\iHom(M,\,a(F,K))$ admits a base compactification over $F \indcot \cO_K$ in the sense of
the preceding section. In fact we go one step
further and show that the $A \otimes K$-module shtuka $\iHom(M,\,a(F/A,K))$ 
admits a base compactification over $A \otimes \cO_K$.
This result was inspired by the construction of extension by zero
from the theory of B\"ockle-Pink \cite[Section 4.5]{bp}.

\begin{lem}\label{reflfreea}%
If $N$ is a finitely generated reflexive module over the ring
$A\complot\cO_K$ then $N$ is locally free.\end{lem}

\pf The ring $A \otimes \cO_K$ is noetherian, regular and of Krull dimension
$2$. Hence so is its completion $A \complot \cO_K$ at the Jacobson radical
$\fm_K \subset \cO_K$. So the result follows from
\cite[Corollary 1.4]{hartrefl}. \quod

\begin{lem}\label{tatelattice}%
Every locally free $A\complot K$-module admits a locally free
$A\complot\cO_K$-lattice.\end{lem}

\pf Let $N$ be a locally free $A \complot K$-module. Clearly there exists a
finitely generated
$A\complot\cO_K$-lattice $N_0 \subset N$. A priori $N_0$ need not be
locally free. But the double dual $N_0^{**}$ is still a lattice in $N$ and is a
reflexive $A\complot\cO_K$-module by \cite[Corollary 1.2]{hartrefl}.
Now Lemma \ref{reflfreea} shows that $N_0^{**}$ is a locally free
$A\complot\cO_K$-module. \quod

\begin{prp}\label{modbcomp}%
Let $N$ be a left $A \complot K\{\tau\}$-module which is locally free of finite
rank as an $A\complot K$-module.
\begin{enumerate}
\item There exists an
$A\complot \cO_K\{\tau\}$-submodule $N_0 \subset N$ such that
$N_0$ is a locally free  $A\complot\cO_K$-lattice in the $A\complot
K$-module $N$.

\item Given an open ideal $I \subset \cO_K$ one can choose $N_0$ in such a way
as to ensure that
$\tau$ acts by zero on $N \otimes_{A\complot\cO_K} (A \otimes \cO_K/I)$.
\end{enumerate}\end{prp}

\pf
(1) By Lemma \ref{tatelattice} the $A\complot K$-module $N$ admits a locally
free $A\complot\cO_K$-lattice $N_1 \subset N$. Let $\zeta \in \cO_K$ be a
uniformizer. According to Proposition \lref{otalgprops}{discrcomplotloc}
the ring $A\complot K$ is the
localization of $A\complot\cO_K$ at $\zeta$. Hence $N_1[\zeta^{-1}] = N$.
As a consequence there exists an $n \geqslant 0$ such that
$\tau N_1 \subset \zeta^{-n} N_1$. The locally free $A\complot\cO_K$-lattice
$N_0 = \zeta^n N_1$ has the property that
$\tau N_0 \subset \zeta^{(q-1)n} N_1$.
As $q > 1$ it follows that $N_0$ is a left $A\complot\cO_K\{\tau\}$-submodule.

(2) Without loss of generality we may assume that $I = \zeta^n \cO_K$ for some
$n \geqslant 0$. Let $N_0 \subset N$ be a left
$A\complot\cO_K\{\tau\}$-submodule as in (1) and let $N_1 = \zeta^n N_0$.
We have $\tau(N_1) \subset \zeta^{(q-1)n} N_1$. The result follows since
$q > 1$. \quod

%

\begin{dfn}\label{fabasecomp}\index{idx}{compactification!base compactification}%
A \emph{base compactification}
of 
$\iHom(M,\,a(F/A,K))$ is a locally free $A\otimes\cO_K$-lattice
$\cM^b$ such that
$\cM^b(A \otimes \cO_K/\fm_K)$ is nilpotent and
$\cM^b(A \otimes \cO_K/\rami)$ is linear.
Here $\rami \subset \cO_K$ is the ramification ideal ideal of
Definition \ref{deflocmodcond}.\end{dfn}

\breakflow
Any base compactification of 
$\iHom(M,\,a(F/A,K))$
induces a base compactification of 
$\iHom(M,\,a(F,K))$ in
the sense of Definition \ref{defbcomp}.

\begin{thm}\label{basecompexist}%
The shtuka $\iHom(M,\,a(F/A,K))$ admits a base
compactification.\end{thm}

\pf The $A\otimes K$-module shtuka $\iHom(M,\,a(F/A,K))$ is locally free by
Lemma \ref{ihomafaklat} while the $A\complot K$-module shtuka
$\iHom(M,\,b(F/A,K))$ is locally free 
by Lemma \ref{ihombfaklat}.
According to Lemma \ref{ihomafakblat}
the shtuka $\iHom(M,\,a(F/A,K))$ is an $A\otimes K$-lattice in 
$\iHom(M,\,b(F/A,K))$.
Hence Beauville-Laszlo glueing theorem [\stacks{0BP2}] implies that
to prove the existence of a base compactification
it is enough to construct a locally free $A\complot\cO_K$-lattice
$\cM \subset \iHom(M,\,b(F/A,K))$ such that
$\cM(A \otimes \cO_K/\fm_K)$ is nilpotent and
$\cM(A\otimes\cO_K/\rami)$ is linear.

Now consider the shtuka
\begin{equation*}
\iHom(M,\,b(F/A,K)) =
\Big[ \Hom(M,\,b(F/A,K)) \shtuka{i}{j}
\Hom(\tau^\ast M,\,b(F/A,K)) \Big].
\end{equation*}
Theorem~\lref{ihomcoh}{ihombvanish} shows that
$\RGamma(\Der\iHom(M,\,b(F/A,K))) = 0$.
Hence
the arrow $i$ in the diagram above
is an isomorphism.
If we let $\tau$ act on 
the $A\complot K$-module $\Hom(M,\,b(F/A,K))$ via the endomorphism
$i^{-1} \circ j$ then it 
becomes a left
$A\complot K\{\tau\}$-module. By construction $\iHom(M,\,b(F/A,K))$ is isomorphic
to the shtuka defined by $\Hom(M,\,b(F/A,K))$.
Therefore applying Proposition \ref{modbcomp}
to $N = \Hom(M,\,b(F/A,K))$ with $I = \rami$ we get the result. \quod

\section{Local models as elliptic shtukas}

We keep the notation and the assumptions of the previous section.

By Lemma~\ref{homablat}
the $F\indcot K$-module shtuka $\iHom(M,\,a(F,K))$ is a lattice
in the $F^\#\complot K$-module shtuka $\iHom(M, \,b(F,K))$.
Since $g(F,K)$ is the quotient of $b(F,K)$ by $a(F,K)$ we deduce that
\begin{equation*}
\iHom(M,\,g(F,K)) = \tfrac{\iHom(M,\,a(F,K))(F^\#\complot K)}{\iHom(M,\,a(F,K))(F\indcot K)}.
\end{equation*}
As a result Proposition~\ref{locfreegermcoh} provides us with
natural quasi-\hspace{0pt}isomorphisms
\begin{align*}
\RGammag(\iHom(M,\,a(F,K))) &\xrightarrow{\isosign} \RGamma(\iHom(M,\,g(F,K)))[-1], \\
\RGammag(\Der\iHom(M,\,a(F,K))) &\xrightarrow{\isosign} \RGamma(\Der\iHom(M,\,g(F,K)))[-1].
\end{align*}

%
%

\begin{dfn}\label{ihomalie}%
We define natural quasi-\hspace{0pt}isomorphisms
\begin{align*}
\RGammag(\iHom(M,\,a(F,K))) &\xrightarrow{\isosign} \Lie_E(K)[-1], \\
\RGammag(\Der\iHom(M,\,a(F,K))) &\xrightarrow{\isosign} \Lie_E(K)[-1]
\end{align*}
as the compositions
\begin{equation*}
\xymatrix{
\RGammag(\Der\iHom(M,\,a(F,K))) \ar[d]_{\textup{Proposition }\ref{locfreegermcoh}}^{\rtviso} &
\RGammag(\iHom(M,\,a(F,K))) \ar[d]^{\textup{Proposition }\ref{locfreegermcoh}}_{\ltviso} \\
\RGamma(\Der\iHom(M,\,g(F,K)))[-1] \ar[d]_{\textup{Proposition }\ref{ihomglie}}^{\rtviso} &
\RGamma(\iHom(M,\,g(F,K)))[-1] \ar[d]^{\textup{Proposition }\ref{ihomge}}_{\ltviso} \\
\Lie_E(K)[-1] & \Lie_E(K)[-1]
}
\end{equation*}%
\end{dfn}

\begin{dfn}\label{gammamaps}%
Let $\cM \subset \iHom(M,\,a(F,K))$ be a local model.
We define a map
\begin{equation*}
\gamma\colon \RGamma(\cM) \to \Lie_E(K)[-1]
\end{equation*}
as the composition of the following maps:
\begin{itemize}
\item The natural map $\RGamma(\cM) \to \RGamma(F\indcot \cO_K,\cM)$.

\item The local germ map
$\RGamma(F\indcot\cO_K,\cM) \xrightarrow{\isosign} \RGammag(F\indcot K, \cM)$
of Definition~\lref{lcmp}{defloccmp}.

\item The isomorphism
$\RGammag(F\indcot K,\cM) = \RGammag(\iHom(M,\,a(F,K)))$ which results
from the equality $\cM(F\indcot K) = \iHom(M,\,a(F,K))$.

\item The quasi-isomorphism
$\RGammag(\iHom(M,\,a(F,K))) \xrightarrow{\isosign} \Lie_E(K)[-1]$
of Definition \ref{ihomalie}.
\end{itemize}
We define a map
\begin{equation*}
\Der\gamma\colon \RGamma(\Der\cM) \to \Lie_E(K)[-1]
\end{equation*}
in the same way.%
\end{dfn}

\breakflow
Observe that the map $\gamma$ is $\cO_F$-linear by construction
while $\Der\gamma$ is both $\cO_F$-linear and $\cO_K$-linear.

\begin{lem}\label{locqi}%
If $\cM$ is a locally free $\cO_F\indcot\cO_K$-module shtuka then
the natural map
$F \otimes_{\cO_F} \RGamma(\cM) \to \RGamma(F\indcot\cO_K,\,\cM)$
is a quasi-isomorphism.\end{lem}

\pf 
According to Proposition \lref{otalgprops}{compindcotloc} the natural map
$F \otimes_{\cO_F}(\cO_F \indcot \cO_K) \to F \indcot \cO_K$ is an isomorphism.
The differentials in the complex computing $\RGamma(\cM)$ are $\cO_F$-linear so
the result follows. \quod

\begin{cor}\label{gammaiso}%
The $F$-linear extensions
\begin{align*}
\gamma&\colon F \otimes_{\cO_F} \RGamma(\cM) \to \Lie_E(K)[-1], \\
\Der\gamma&\colon F \otimes_{\cO_F} \RGamma(\Der\cM) \to \Lie_E(K)[-1]
\end{align*}
of the maps $\gamma$ and $\Der\gamma$ are quasi-isomorphisms.\quod\end{cor}


\begin{prp}\label{locmodcohiso}%
If $\cM \hookrightarrow \cN$ is an inclusion of local models then the induced maps
\begin{align*}
F \otimes_{\cO_F} \RGamma(\cM) &\to F \otimes_{\cO_F} \RGamma(\cN), \\
F \otimes_{\cO_F} \RGamma(\Der\cM) &\to F \otimes_{\cO_F} \RGamma(\Der\cN)
\end{align*}
are quasi-isomorphisms.\end{prp}

\pf Since $\gamma$ and $\Der\gamma$ are natural the result follows from
Corollary \ref{gammaiso}. \quod

\begin{thm}\label{locmodcohdesc}%
Let $\cM$ be a local model.
\begin{enumerate}
\item $\uH^0(\cM) = 0$ and $\uH^1(\cM)$ is a finitely generated free $\cO_F$-module.

\item $\uH^0(\Der\cM) = 0$ and $\uH^1(\Der\cM)$ is a finitely generated free $\cO_F$-module.
\end{enumerate}\end{thm}

\pf (1) 
By definition $\cM(\cO_F/\fm_F \otimes K) =
\cM^c(\cO_F/\fm_F \otimes K)$ where $\cM^c \subset \iHom(M,\,a(F,K))$ is the
coefficient compactification.
Lemma \ref{drincompnilp} claims that $\cM^c(\cO_F/\fm_F \otimes K)$
is nilpotent.
Hence the result follows from Theorem
\lref{locloccoh}{llcoh}.
(2) Nilpotence is preserved under linearization. So the result follows
from Theorem \lref{locloccoh}{llcoh} as well.
\quod

\breakflow
The main result of this section is the following theorem:

\begin{thm}\label{locmodell}%
A local model $\cM$ is an elliptic shtuka of
ramification ideal $\rami$ where $\rami \subset \cO_K$ is the ideal
of Definition \ref{deflocmodcond}.\end{thm}

\pf
We verify the conditions of Definition \lref{reg}{defell} for $\cM$.
\begin{enumerate}
\setcounter{enumi}{-1}
\renewcommand{\theenumi}{E\arabic{enumi}}
\item 
Indeed $\cM$ is a locally free $\cO_F\complot\cO_K$-module shtuka by definition.

\item 
By construction $\cM(\cO_F/\fm_F \otimes K) =
\cM^c(\cO_F/\fm_F \otimes K)$ where $\cM^c \subset \iHom(M,a(F,K))$ is the
coefficient compactification.
Thus the shtuka $\cM(\cO_F/\fm_F \otimes K)$ is nilpotent by
Lemma \ref{drincompnilp}.

\item $\cM(F \otimes \cO_K/\fm_K)$ is nilpotent by definition.

\item Consider the map
$\uH^1(\Der\gamma)\colon \uH^1(\Der\cM) \to \Lie_E(K)$.
This map is $\cO_F$-linear and $\cO_K$-linear by construction.
Corollary \ref{gammaiso} in combination with Theorem \ref{locmodcohdesc}
implies that the map is injective.
Therefore
\begin{equation*}
\fm_F \cdot \uH^1(\Der\cM) = \rami \cdot \uH^1(\Der\cM)
\end{equation*}
by definition of the ramification ideal ideal $\rami$.

\item The shtuka $\cM(F \otimes \cO_K/\rami)$ is linear by definition.
As $\cM$ is locally free it follows that $\cM(\cO_F\otimes\cO_K/\rami)$ is linear.
\quod
\end{enumerate}

\breakflow
Theorem \ref{locmodell} allows us to make the following important definition:

\begin{dfn}\label{locmodreg}\index{idx}{regulator!local model@of a local model}%
The \emph{regulator}
\begin{equation*}
\rho\colon \uH^1(\cM) \xrightarrow{\,\,\isosign\,\,}
\uH^1(\Der\cM)
\end{equation*}
of a local model $\cM$ is the regulator of $\cM$ viewed as an
elliptic shtuka over $\cO_F\complot\cO_K$
(see Definition~\lref{reg}{defellreg}).%
\end{dfn}

\section{The exponential map}
\label{sec:locmodexp}

We keep the assumptions and the notation of the previous section.
In this section we describe in detail how the
maps $\gamma$ and $\Der\gamma$ act on cohomology classes
and introduce the exponential map of a local model.
Before we begin let us make a remark.
%
Let $h_1, h_2 \in b(F,K)$. Recall that the expression
\begin{equation*}
h_1 \sim h_2
\end{equation*}
signifies that the image of $h_1 - h_2$ in $g(F,K)$ is zero. In other words
there exists an open neighbourhood $U \subset F$ of $0$ such that $h_1|_U =
h_2|_U$.

\begin{prp}\label{gammadesc}%
Let $\cM = [\cM_0 \shtuka{i}{j} \cM_1]$ be a local model
and let $g \in \cM_1$.
\begin{enumerate}
\item There exists a unique $f \in \cM_0(F^\#\complot\cO_K)$ such that $(i-j)(f)
= g$.

\item Let $[g]$ be the cohomology class of $g$ in $\uH^1(\cM)$ and let $\alpha =
\gamma[g]$. The element $\alpha\in\Lie_E(K)$ is uniquely characterized by the
property that for every $m \in M^0$ one has
\begin{equation*}
f(m) \sim (x \mapsto m \exp (x\alpha)).
\end{equation*}
Here we view $f$ as an element of $\Hom(M,\,b(F,K))$ via the natural inclusion
$\cM(F^\#\complot\cO_K) \subset \iHom(M,\,b(F,K))$ of Lemma~\ref{locmodblat}.
\end{enumerate}\end{prp}

\pf (1) is a direct consequence of Proposition~\lref{lcmp}{loccmpdesc} (1).
Part (2) of this proposition implies that the image of $f$ in
$\Hom(M,\,g(F,K))$ represents the image of $[g] \in \uH^1(\cM)$ under the
local germ map
\begin{equation*}
\uH^1(\cM) \xrightarrow{\isosign} \uH^0_{\textup{g}}(F \indcot K, \cM) =
\uH^0(\iHom(M,g(F,K))).
\end{equation*}
By Proposition~\ref{ihomge} the map
\begin{equation*}
\Lie_E(K) \to \uH^0(\iHom(M,\,g(F,K))), \quad
\alpha \mapsto (m \mapsto (x \mapsto m\exp(x\alpha)))
\end{equation*}
is an isomorphism.
Hence
for every $m \in M$ we have
$f(m) \sim (x \mapsto m \exp (x\alpha))$. It remains to show that
if $\beta \in \Lie_E(K)$ is an element such that
$f(m) \sim (x \mapsto m \exp (x \beta))$ for all $m \in M^0 \subset M$ then
$\beta = \alpha$.

According to Proposition \lref{drmot}{motzero}
$M^0$ generates $M$
as a left $K\{\tau\}$-module.
As the image of $f$ in $\Hom(M,g(F,K))$ represents
a class in $\uH^0(\iHom(M,g(F,K)))$ Proposition \lref{genhomsht}{ihomzerocoh} shows that
the map $M \to g(F,K)$ induced by $f$ is in fact a homomorphism
of left $A \otimes K\{\tau\}$-modules. As a consequence
$f(m) \sim (x \mapsto m \exp (x\beta))$ for all $m \in M$. Whence the result.
\quod

\begin{prp}\label{liegammadesc}%
Let $\cM = [\cM_0 \shtuka{i}{j} \cM_1]$ be a local model
and let $g \in \cM_1$.
\begin{enumerate}
\item There exists a unique $f \in \cM_0(F^\#\complot\cO_K)$ such that $i(f) = g$.

\item
Let $[g]$ be the cohomology class of $g$ in $\uH^1(\Der\cM)$
and let $\alpha = \Der\gamma[g]$.
The element $\alpha\in\Lie_E(K)$ is uniquely characterized by
the property that for every $m \in M^0$ one has
\begin{equation*}
f(m) \sim (x \mapsto dm(x\alpha)).
\end{equation*}
Here we view $f$ as an element of $\Hom(M,\,b(F,K))$ via the natural inclusion
$\cM(F^\#\complot \cO_K) \subset \iHom(M,\,b(F,K))$ of Lemma \ref{locmodblat}.
Given $m\in M = \Hom(E,\bG_a)$ we denote $dm\colon \Lie_E \to \Lie_{\bG_a}$ the
induced map of Lie algebras.
\end{enumerate}\end{prp}

\pf Same as the proof of Proposition \ref{gammadesc}. \quod

\breakflow
In the following it would be convenient for us to
assemble the maps $\gamma$ and $\Der\gamma$ into a map acting on the cohomology
of a local model.

\begin{dfn}\label{defexpmap}\index{idx}{exponential map!of a local model}%
Let $\cM$ be a local model. We define the \emph{exponential map}
\begin{equation*}
\exp\colon F \otimes_{\cO_F} \uH^1(\Der\cM) \to F \otimes_{\cO_F} \uH^1(\cM)
\end{equation*}
as the composition $\uH^1(\gamma) \circ \uH^1(\Der\gamma)^{-1}$.\end{dfn}

\breakflow
Even though this exponential map is induced by the identity map
on the Lie algebra, one can justify its name by looking at what it does to
the auxillary elements $f \in \Hom(M,b(F,K))$ as described in Propositions
\ref{liegammadesc} and \ref{gammadesc}.

We will show in the subsequent chapters that the exponential map is nothing but
the inverse of the regulator map $\rho\colon \uH^1(\cM) \to \uH^1(\Der\cM)$
which we introduced in Definition \ref{locmodreg}. This important result appears
to be neither easy nor evident. For one thing, the regulator is a purely
shtuka-theoretic construct while the exponential map uses the arithmetic data of
the Drinfeld module in an essential way. The only proof we have at the moment is
rather technical and is based on explicit computations.


\chapter{Change of coefficients}
\label{chapter:coeffch}
\label{ch:homshtdual}
\label{ch:cclocmod}

This chapter is of a technical nature.
In it we verify that the constructions of Chapter \ref{chapter:locmod}
are compatible with restriction of the coefficient ring $A$.
This result will be used in Chapter \ref{chapter:locmodreg}
where we show that
the regulator of a local model
is the inverse of the exponential map.
We will do it by reduction to an explicit computation
in the case $A = \Fq[t]$.

\section{Duality for Hom shtukas}
\label{sec:homshtdual}

In this section we describe a general duality construction for Hom shtukas.
It will be used several times in the rest of the chapter. We begin with an
auxillary result.

\begin{lem}\label{shtmodhomsht}%
Let $R$ be a $\tau$-ring. If $M = [M_0 \shtuka{i_M}{j_M} M_1]$ is an $R$-module
shtuka and $N$ a left $R\{\tau\}$-module then
\begin{equation*}
\iHom_R(M,N) = \Big[ \Hom_R(M_1, N) \shtuka{i}{j} \Hom_R(\tau^\ast M_0, N) \Big]
\end{equation*}
where
\begin{align*}
i(f)&\colon r \otimes m \mapsto r f j_M(m), \\
j(f)&\colon r \otimes m \mapsto r \tau \cdot f i_M(m).
\end{align*}
\end{lem}

\pf Follows directly from the definition of $\iHom$ (Definition
\lref{genhomsht}{defhomsht}). \quod

\breakflow
If $\varphi\colon R \to S$ is a ring homomorphism, $M$ an $S$-module and $N$ an $R$-module then
the $R$-modules $\Hom_R(M,N)$ and $\Hom_R(S,N)$ carry natural $S$-module
structures: via $M$ in the first case and via $S$ in the second.
Furthermore the natural duality map
\begin{equation*}
\Hom_S(M,\Hom_R(S,N)) \to \Hom_R(M,N), \quad
f \mapsto [m \mapsto f(m)(1)]
\end{equation*}
is an $S$-module isomorphism. We would like to establish an analog of this
duality for a
$\tau$-ring homomorphism $\varphi\colon (R,\tau) \to (S,\sigma)$ and $\iHom$ in place
of $\Hom$. In the following it will be important to distinguish the $\tau$-endomorphisms
of $R$ and $S$. We thus denote them by different letters.

Let $\varphi\colon (R,\tau) \to (S,\sigma)$ be a homomorphism of $\tau$-rings. For
every $S$-module $M$ there is a natural base change map
\begin{equation*}
\mu_M\colon \tau^\ast M \to \sigma^\ast M, \quad
r \otimes m \mapsto \varphi(r) \otimes m.
\end{equation*}
In particular we have a base change map $\mu_S\colon \tau^\ast S \to \sigma^\ast
S = S$.

Consider the commutative diagram 
\begin{equation}\label{dualcocart}\tag{$\ast$}%
\vcenter{\vbox{\xymatrix{
S \ar[r]^{\sigma} & S \\
R \ar[u]^{\varphi} \ar[r]^{\tau} & R \ar[u]_{\varphi}
}}}
\end{equation}
For the duality statements below to work it will be necessary to assume that
\eqref{dualcocart} is cocartesian in the category of rings.
It is cocartesian if and only if the base
change map $\mu_S\colon \tau^\ast S \to S$ is an isomorphism.

\begin{prp}\label{shthomdualmod}%
Let $\varphi\colon (R,\tau) \to (S,\sigma)$ be a homomorphism of $\tau$-rings
and let $N$ be a left $R\{\tau\}$-module.  Consider the shtuka
\begin{equation*}
\iHom_R(S,N) =
\Big[ \Hom_R(S,N) \shtuka{i}{j} \Hom_R(\tau^\ast S, N) \Big].
\end{equation*}
If the
commutative diagram of rings \eqref{dualcocart}
is cocartesian then the following holds:
\begin{enumerate}
\item $i$ is an isomorphism.

\item The endomorphism $i^{-1} j$
of $\Hom_R(S,N)$ makes it into a left $S\{\sigma\}$-module.

\item If $f \in \Hom_R(S,N)$ then $g = i^{-1} j(f)$ is the unique $R$-linear map
such that
\begin{equation*}
g\big[\varphi(r) \sigma(s)\big] = r \tau \cdot f(s)
\end{equation*}
for every $r \in R$ and $s \in S$.
\end{enumerate}
\end{prp}

\pf (1) Let $\tau^a\colon \tau^\ast S \to S$ be the adjoint of the
$\tau$-multiplication map of $S$. By definition
\begin{equation*}
\tau^a(r \otimes s) = \varphi(r) \sigma(s).
\end{equation*}
At the same time
\begin{equation*}
\mu_S(r \otimes s) = \varphi(r) \otimes s = \varphi(r) \sigma(s).
\end{equation*}
Thus $\tau^a = \mu_S$ is an isomorphism and we conclude that $i$ is an
isomorphism from Lemma \ref{shtmodhomsht}.

(2) We will deduce from (3) that the endomorphism $i^{-1} j$ is $\sigma$-linear.
Let us temporarily denote $e = i^{-1} j$. Let $f \in
\Hom_R(S,N)$. If
$s_1 \in S$ then the maps
$e(f \cdot s_1)$ and $e(f) \cdot \sigma(s_1)$ satisfy
\begin{align*}
e(f\cdot s_1)\big[\varphi(r)\sigma(s)\big] &= r \tau \cdot f(s_1 s) \\
e(f)\cdot \sigma(s_1)\big[\varphi(r)\sigma(s)\big] &=
e(f)\big[\sigma(s_1)\varphi(r)\sigma(s)\big] =
r \tau \cdot f(s_1 s)
\end{align*}
for every $r \in R$, $s \in S$. Thus $e$ is
$\sigma$-linear.

(3) Let $\tau_N^a\colon \tau^\ast N \to N$ be the adjoint of the
$\tau$-multiplication map. If $f \in \Hom_R(S,N)$ then according to Lemma
\ref{shtmodhomsht}
\begin{equation*}
j(f)\colon r \otimes s \mapsto r \tau \cdot f(s)
\end{equation*}
for every $r \in R$ and $s \in S$. As we saw in (1) the map $i$ is given by the
composition with the base change map $\mu_S\colon r \otimes s \mapsto
\varphi(r)\sigma(s)$. The result is now clear. \quod

\begin{prp}\label{homshtdualiso}%
Let $\varphi\colon (R,\tau) \to (S,\sigma)$ be
a homomorphism of $\tau$-rings. Let
\begin{equation*}
M = \Big[ M_0 \shtuka{i_M}{j_M} M_1 \Big]
\end{equation*}
be an $S$-module shtuka and $N$ a
left $R\{\tau\}$-module. If the
commutative diagram of rings \eqref{dualcocart}
is cocartesian
then the maps
\begin{align*}
\Hom_S(M_1,\Hom_R(S,N)) &\xrightarrow{\,\,\quad\quad\textup{duality}\quad\quad\,\,}
\Hom_R(M_1,N), \\
\Hom_S(\sigma^\ast M_0, \Hom_R(S,N))
&\xrightarrow{\quad(\mu_{M_0})^* \circ\,\textup{duality}\quad}
\Hom_R(\tau^\ast M_0, N)
\end{align*}
define an isomorphism
\begin{equation*}
\iHom_S(M,\Hom_R(S,N)) \cong \iHom_R(M,N)
\end{equation*}
of $R$-module shtukas.
The left $S\{\sigma\}$-module structure on $\Hom_R(S,N)$
is as constructed in Proposition \ref{shthomdualmod}.
\end{prp}

\pf
Denote the duality maps
\begin{align*}
\eta_0&\colon \Hom_S(M_1,\Hom_R(S,N)) \to \Hom_R(M,N), \\
\eta_1&\colon \Hom_S(\sigma^\ast M_0, \Hom_R(S,N)) \to \Hom_R(\tau^\ast M_0, N).
\end{align*}
If \eqref{dualcocart} is cocartesian then the base change map $\mu_{M_0}$ is an
isomorphism. The inverse of $\mu_{M_0}$ is
given by the formula
\begin{equation*}
\varphi(r)\sigma(s) \otimes m \mapsto r \otimes s m.
\end{equation*}
Thus $\eta_0$ and $\eta_1$ define an isomorphism of $R$-module
shtukas provided they form a morphism of shtukas. Let us show that it is indeed
the case.

Let $i_S$, $j_S$ be the arrows of $\iHom_S(M,\Hom_R(S,N))$ and let $i_R$,
$j_R$ be the arrows of $\iHom_R(M,N)$.
If $f \in \Hom_S(M_1,\Hom_R(S,N))$ then
\begin{equation*}
\eta_0(f)\colon m \mapsto f(m)(1).
\end{equation*}
So Lemma \ref{shtmodhomsht} implies that
\begin{equation*}
j_R(\eta_0(f))\colon r \otimes m \mapsto r \tau \cdot [f i_M(m)](1).
\end{equation*}
By the same Lemma
\begin{equation*}
j_S(f)\colon s \otimes m \mapsto s \sigma \cdot f i_M(m).
\end{equation*}
If $g \in \Hom_S(\sigma^\ast M_0, \Hom_R(S,N))$ then
\begin{equation*}
\eta_1(g)\colon r \otimes m \mapsto g(\varphi(r) \otimes m)(1).
\end{equation*}
Therefore
\begin{align*}
\eta_1(j_S(f))\colon r \otimes m &
\mapsto [j_S(f)(\varphi(r) \otimes m)](1)\\
& = [\varphi(r)\sigma \cdot f i_M(m)](1)\\
& = [\sigma\cdot f i_M(m)](\varphi(r)).
\end{align*}
According to Proposition \ref{shthomdualmod}
\begin{equation*}
\sigma\cdot f i_M(m)\colon \varphi(r) \mapsto r \tau \cdot [f i_M(m)](1).
\end{equation*}
Hence
\begin{equation*}
\eta_1(j_S(f))\colon r \otimes m \mapsto r \tau \cdot [f i_M(m)](1).
\end{equation*}
We conclude that $\eta_1(j_S(f)) = j_R(\eta_0(f))$. It is easy to see that
$\eta_1(i_S(f)) = i_R(\eta_0(f))$. \quod

%
%

\section{Tensor products}


\begin{lem}\label{chdomfields}%
Let $T$ be a locally compact $\Fq$-algebra
and let $S' \subset S$ be an extension of locally compact $\Fq$-algebras.
If $S$ is finitely generated free as a topological $S'$-module then
the natural map
$S \otimes_{S'} (S' \indcot T) \to S \indcot T$
is an isomorphism.\end{lem}

\pf We rewrite the natural map in question as
\begin{equation*}
(S \otimes_{\textup{ic}} T) \otimes_{(S' \otimes_{\textup{ic}} T)} (S'\indcot T)
\to S \indcot T.
\end{equation*}
By assumption $S$ is a finitely generated free topological $S'$-module.
Therefore $S
\otimes_{\textup{ic}} T$ is a finitely generated free topological $S' \otimes_{\textup{ic}}
T$-module. The result now follows from Lemma
\lref{topmodcompl}{compllattice}.\quod

\begin{lem}\label{chdomfieldsb}%
Let $T$ be a locally compact $\Fq$-algebra and let $S' \subset S$ be an
extension of locally compact $\Fq$-algebras. If $S$ is locally free of finite
rank as an $S'$-module without topology then the natural map
$S \otimes_{S'} ((S')^\# \complot T) \to S^\#\complot T$
is an isomorphism.\end{lem}

\pf We rewrite the natural map in question as
\begin{equation*}
(S^\# \otimes_{\textup{c}} T) \otimes_{((S')^\# \otimes_{\textup{c}} T)}
((S')^\#\complot T) \to S^\# \complot T.
\end{equation*}
By assumption $S$ is locally free of finite rank as an $S'$-module without
topology.
Therefore $S^\# \otimes_{\textup{c}} T$ is a
topological direct summand of a finitely generated free
$(S')^\#\otimes_{\textup{c}} T$-module.
The result now follows from Lemma
\lref{topmodcompl}{compllattice}.\quod

\section{The setting}

We now start with the main part of this chapter. The setting is as follows.
Fix a coefficient ring $A$ as in Definition \lref{coeffrings}{defcoeffring}.
As usual we denote $F$ the local field of $A$ at infinity, $\cO_F \subset F$ the
ring of integers and $\fm_F \subset \cO_F$ the maximal ideal.
Let $K$ be a finite product of local fields containing $\Fq$.
Fix a Drinfeld $A$-module $E$ over $K$ and let $M = \Hom(E,\bG_a)$ be its motive.
Throughout the chapter we assume that
the action of $A$ on $\Lie_E(K)$ extends to a continuous action of $F$.
This is the context in which we defined and studied local models in Chapter
\ref{chapter:locmod}.

Fix an $\Fq$-subalgebra $A' \subset A$ such that $A$ is finite flat over $A'$.
Note that $A'$ is itself a coefficient ring.
We denote $F'$ the local field of $A'$ at infinity, $\cO_{F'} \subset F'$ the
ring of integers and $\fm_{F'} \subset \cO_{F'}$ the maximal ideal.

\begin{lem} The motive of $E$ viewed as a Drinfeld $A'$-module is
$M$ with its natural left $A'\otimes K\{\tau\}$-module structure.\quod\end{lem}

%

\section{Hom shtukas}

%
%
%

%

\begin{lem}\label{ccacocart}%
The commutative square of rings
\begin{equation*}
\xymatrix{
A \otimes K \ar[r] & F\indcot K \\
A'\otimes K\ar[u] \ar[r] & F'\indcot K \ar[u]
}
\end{equation*}
is cocartesian.\end{lem}

\pf Indeed $A \otimes_{A'} F' = F$ so
\begin{equation*}
A \otimes_{A'} (F'\indcot K) =
(A \otimes_{A'} F') \otimes_{F'} (F'\indcot K) =
F \otimes_{F'} (F'\indcot K).
\end{equation*}
Now Lemma \ref{chdomfields} shows that
$F \otimes_{F'} (F'\indcot K) = F\indcot K$ and the result follows.\quod

\begin{cor}\label{ccmcocart}%
The natural map
$M \otimes_{A'\otimes K}(F'\indcot K) \to
M \otimes_{A \otimes K}(F\indcot K)$
is an isomorphism of left $F'\indcot K\{\tau\}$-modules.\quod
\end{cor}

\begin{lem}\label{ccfcocart}%
The commutative square of rings
\begin{equation*}
\xymatrix{
F \indcot K \ar[r]^{\tau} & F \indcot K \\
F' \indcot K \ar[u] \ar[r]^{\tau} & F'\indcot K \ar[u]
}
\end{equation*}
is cocartesian.\end{lem}

\pf Immediate from Lemma \ref{chdomfields}.\quod

\breakflow
We are thus in position to apply the duality machinery of
Section \ref{sec:homshtdual} to the $\tau$-ring homomorphism
$F'\indcot K \to F\indcot K$. Proposition \ref{shthomdualmod}
equips the $F\indcot K$-module
\begin{equation*}
\Hom_{F'\indcot K}(F \indcot K, a(F',K))
\end{equation*}
with the structure of a left $F\indcot K\{\tau\}$-module.

\begin{lem}\label{cchoma}%
The natural map
\begin{equation*}
\Hom_{F'\indcot K}(F\indcot K,\,a(F',K)) \to \Hom_{F'}(F,\,a(F',K))
\end{equation*}
is an isomorphism of $F\indcot K$-modules.\end{lem}

\pf Follows since $F\otimes_{F'}(F'\indcot K) = F\indcot K$ by Lemma \ref{chdomfields}.\quod

\breakflow
So we get a left $F\indcot K\{\tau\}$-module structure on
$\Hom_{F'}(F,\,a(F',K))$.

\begin{lem}%
If $g \in \Hom_{F'}(F,\,a(F',K))$ then $\tau\cdot g$ maps $x \in F$ to
$\tau(g(x))$.\quod\end{lem}

\begin{lem}\label{chdomshriekiso}%
The map
\begin{equation*}
a(F,K) \to \Hom_{F'}(F,\,a(F',K)), \quad
f \mapsto [x \mapsto (y \mapsto f(yx))]
\end{equation*}
is an isomorphism of left $F\indcot K\{\tau\}$-modules.
\end{lem}

\pf 
For a finite-dimensional $F'$-vector space $V$ let $\alpha_V$ be the map
\begin{equation*}
\alpha_V\colon a(V,K) \to \Hom_{F'}(V,\,a(F',K)), \quad
f \mapsto [v \mapsto (y \mapsto f(yv))].
\end{equation*}
It is clearly an $F'\indcot K$-linear isomorphism if $V$ is of dimension $1$.
The map $\alpha_V$ is also natural in $V$. Hence $\alpha_V$ is an $F'\indcot
K$-linear isomorphism  for any finite-dimensional $F'$-vector space $V$ and in
particular for $V = F$.

The map $\alpha_F$ is also $F$-linear. By Lemma \ref{chdomfields} we have
$F\otimes_{F'}(F'\indcot K) = F\indcot K$.  Hence $\alpha_F$ is $F\indcot
K$-linear. A simple computation shows that it commutes witht the action of
$\tau$. So we get the result.\quod


\begin{prp}\label{ccihoma}%
The restriction map $a(F,K) \to a(F',K)$ induces an isomorphism
\begin{equation*}
\iHom_{A \otimes K}(M,\,a(F,K)) \cong
\iHom_{A'\otimes K}(M,\,a(F',K))
\end{equation*}
of $F'\indcot K$-module shtukas.\end{prp}

\pf We have
\begin{equation*}
\iHom_{A \otimes K}(M,\,a(F,K)) =
\iHom_{F\indcot K}(M \otimes_{A\otimes K}(F\indcot K), \,a(F,K)).
\end{equation*}
Lemma \ref{chdomshriekiso} gives us a natural isomorphism
\begin{align*}
&\iHom_{F\indcot K}(M \otimes_{A\otimes K}(F\indcot K), \,a(F,K)) \\
\cong&
\iHom_{F\indcot K}(M \otimes_{A\otimes K}(F\indcot K), \,\Hom_{F'}(F,a(F',K))).
\end{align*}
In view of Lemma \ref{ccfcocart} we can apply
Proposition \ref{homshtdualiso} to get an isomorphism
\begin{align*}
& \iHom_{F\indcot K}(M \otimes_{A\otimes K}(F\indcot K),\,\Hom_{F'}(F,a(F',K))) \\
\cong &
\iHom_{F'\indcot K}(M\otimes_{A\otimes K}(F\indcot K), \,a(F',K))
\end{align*}
of $F'\indcot K$-module shtukas.
Corollary \ref{ccmcocart} identifies
$M \otimes_{A\otimes K}(F\indcot K)$
with $M \otimes_{A'\otimes K}(F'\indcot K)$.
So we get an isomorphism
\begin{equation*}
\iHom_{A\otimes K}(M,\,a(F,K)) \xrightarrow{\isosign}
\iHom_{A'\otimes K}(M,\,a(F',K))
\end{equation*}
of $F'\indcot K$-module shtukas. A straightforward computation shows that this
isomorphism is induced by the restriction map
$a(F,K) \to a(F',K)$. \quod
 
%

\section{Coefficient compactifications}

We denote $C$ the projective compactification of $\Spec A$ and $C'$ the
projective compactification of $A'$.
Let $\rho\colon C \times \Spec K \to C' \times \Spec K$
be the map induced by the inclusion $A' \subset A$.

\begin{lem}%
The commutative square of schemes
\begin{equation*}
\xymatrix{
\Spec\cO_F\indcot K \ar[r] \ar[d] & C \times X \ar[d]^{\rho} \\
\Spec\cO_{F'} \indcot K \ar[r] & C' \times X
}
\end{equation*}
is cartesian.\end{lem}

\pf Lemma \ref{chdomfields} implies that the square
\begin{equation*}
\xymatrix{
\cO_F \otimes K \ar[r] & \cO_F\indcot K \\
\cO_{F'} \otimes K \ar[u] \ar[r] & \cO_{F'} \indcot K \ar[u]
}
\end{equation*}
is cocartesian. Whence the result. \quod

\begin{cor}\label{ccdrinbasechange}%
For every quasi-coherent sheaf $\cE$ on $C \times X$ the base change map
\begin{equation*}
(\rho_\ast \cE)(\cO_{F'} \indcot K) \to
\cE(\cO_F \indcot K)
\end{equation*}
is an isomorphism of $\cO_{F'}\indcot K$-modules.\quod\end{cor}

\begin{lem}\label{ccdrinstable}%
If $\cE$ is the shtuka on $C \times K$ which corresponds to the
left $A \otimes K\{\tau\}$-module $M$ by Theorem \lref{drinconstr}{motivec} then 
$\rho_\ast\cE$ is the shtuka on $C' \times K$ which corresponds to $M$ viewed as a
left $A' \otimes K\{\tau\}$-module.\end{lem}

\pf Suppose that $\cE$ is given by the diagram
\begin{equation*}
\Big[ \cE_{-1} \shtuka{i}{j} \cE_0 \Big] \subset
\Big[ M \shtuka{1}{\tau} M \Big].
\end{equation*}
Let $f$ be the degree of the residue field of $F$ over $\Fq$ and
let $r$ be the rank of $M$ as an $A \otimes K$-module.
By Theorem \lref{drinconstr}{motivec} the sheaves $\cE_{-1}$, $\cE_0$ are
locally free of rank $r$ and have the following property:
\begin{align*}
&\uH^0(C \times \Spec K, \,\cE_0(n)) = M^{f r n}, \\
&\uH^0(C \times \Spec K, \,\cE_{-1}(n)) = M^{f r n - 1}.
\end{align*}
Let $d = [F : F']$. Observe that
$f = d f'$ where $f'$ is the degree of the residue field of $F'$ over $\Fq$.
The morphism $\rho\colon C \times K \to C' \times K$ is finite flat of degree
$d$. Hence $\rho_\ast\cE_{-1}$, $\rho_\ast\cE_0$ are locally free of rank $dr$
and 
\begin{align*}
&\uH^0(C \times \Spec K, \,\rho_\ast\cE_0(n)) = M^{f' d r n}, \\
&\uH^0(C \times \Spec K, \,\rho_\ast\cE_{-1}(n)) = M^{f' d r n - 1}.
\end{align*}
The unicity part of Theorem \lref{drinconstr}{motivec} now implies the
result.\quod

\breakflow
We next study the function space $a(F/\cO_F,K)$.

\begin{lem}\label{ccofcocart}%
The commutative square of rings
\begin{equation*}
\xymatrix{
\cO_F \indcot K \ar[r]^{\tau} & \cO_F \indcot K \\
\cO_{F'} \indcot K \ar[u] \ar[r]^{\tau} & \cO_{F'}\indcot K \ar[u]
}
\end{equation*}
is cocartesian.\end{lem}

\pf Immediate from Lemma \ref{chdomfields}.\quod

\breakflow
So we can apply the duality constructions of
Section \ref{sec:homshtdual} to the $\tau$-ring homomorphism
$\cO_{F'}\indcot K \to \cO_F\indcot K$. Proposition \ref{shthomdualmod}
equips the $\cO_F\indcot K$-module
\begin{equation*}
\Hom_{\cO_{F'}\indcot K}(\cO_F \indcot K, \,a(F'/\cO_{F'},K))
\end{equation*}
with the structure of a left $\cO_F\indcot K\{\tau\}$-module.

\begin{lem}\label{cchomofa}%
The natural map
\begin{equation*}
\Hom_{\cO_{F'}\indcot K}(\cO_F\indcot K, \,a(F'/\cO_{F'},K)) \to
\Hom_{\cO_{F'}}(\cO_F,\,a(F'/\cO_{F'},K))
\end{equation*}
is an isomorphism of $\cO_F\indcot K$-modules.\quod\end{lem}


\breakflow\noindent
So we get a left $\cO_F\indcot K\{\tau\}$-module structure on
$\Hom_{\cO_{F'}}(\cO_F,\,a(F'/\cO_{F'},K))$.

\begin{lem}%
If $g \in \Hom_{\cO_{F'}}(\cO_F,a(F'/\cO_{F'},K))$ then $\tau\cdot g$ maps $x
\in \cO_F$ to $\tau(g(x))$.\quod\end{lem}

\begin{lem}\label{ccofshriekiso}%
The map
\begin{equation*}
a(F/\cO_F,K) \to \Hom_{\cO_{F'}}(\cO_F,\,a(F'/\cO_{F'},K)), \quad
f \mapsto [x \mapsto (y \mapsto f(yx))]
\end{equation*}
is an isomorphism of left $\cO_F\indcot K\{\tau\}$-modules.
\end{lem}

\pf We view $a(F/\cO_F,K)$ as a subspace of $a(F,K)$ consisting of functions
which vanish on $\cO_F$, and similarly for $a(F'/\cO_{F'}, K)$.
Note that
\begin{equation*}
\Hom_{F'}(F,a(F,K)) = \Hom_{\cO_{F'}}(\cO_F,\,a(F',K)).
\end{equation*}
We can thus identify
\begin{equation*}
\Hom_{\cO_{F'}}(\cO_F,\,a(F'/\cO_{F'},K))
\end{equation*}
with a submodule
of $\Hom_{F'}(F,\,a(F,K))$. In view of this remark the result
follows from Lemma \ref{chdomshriekiso}. \quod

\begin{prp}\label{ccdrincomp}%
The restriction isomorphism
of Proposition \ref{ccihoma} identifies the coefficient compactification
of $\iHom_{A\otimes K}(M,a(F,K))$ with the coefficient compactification
of $\iHom_{A'\otimes K}(M,a(F',K))$.\end{prp}

\pf Lemma \ref{ccofshriekiso} gives us a natural isomorphism
\begin{align*}
&\iHom_{\cO_F\indcot K}(\cE(\cO_F\indcot K), \,a(F/\cO_F,K)) \xrightarrow{\isosign} \\
&\iHom_{\cO_F\indcot K}(\cE(\cO_F\indcot K),\,
\Hom_{\cO_{F'}}(\cO_F,a(F'/\cO_{F'},K))).
\end{align*}
In view of Lemma \ref{ccofcocart} we can apply
Proposition~\ref{homshtdualiso} and identify the shtuka above with
\begin{equation*}
\iHom_{\cO_{F'} \indcot K}(\cE(\cO_F\indcot K), \,a(F'/\cO_{F'},K))
\end{equation*}
It is easy to see that the resulting isomorphism
\begin{align*}
&\iHom_{\cO_F\indcot K}(\cE(\cO_F\indcot K), \,a(F/\cO_F,K)) \xrightarrow{\isosign}\\
&\iHom_{\cO_{F'}\indcot K}(\cE(\cO_F \indcot K), \,a(F'/\cO_{F'},K))
\end{align*}
is induced by the restriction map
$a(F/\cO_F,K) \to a(F'/\cO_{F'},K)$.
Now Corollary \ref{ccdrinbasechange} implies that the natural map
$\rho_\ast\cE(\cO_{F'}\indcot K) \to
\cE(\cO_F\indcot K)$
is an isomorphism. Lemma \ref{ccdrinstable} implies that the shtuka
\begin{equation*}
\iHom_{\cO_{F'}\indcot K}(\rho_\ast\cE(\cO_{F'} \indcot K), \,a(F'/\cO_{F'},K))
\end{equation*}
is the coefficient compactification of $\iHom_{A'\otimes K}(M,\,a(F',K))$.
Whence the result. \quod

\section{Base compactifications}

\begin{lem}\label{ccfokcocart}%
The commutative square of rings
\begin{equation*}
\xymatrix{
F \indcot \cO_K \ar[r] & F \indcot K \\
F' \indcot \cO_K \ar[r] \ar[u] & F' \indcot K \ar[u]
}
\end{equation*}
is cocartesian.
\end{lem}

\pf Follows from Lemma \ref{chdomfields}. \quod

\begin{lem}\label{ccfokquotcocart}%
For every open ideal $I \subset \cO_K$
the commutative square of rings
\begin{equation*}
\xymatrix{
F \indcot \cO_K \ar[r] & F \otimes \cO_K/I \\
F' \indcot \cO_K \ar[r] \ar[u] & F' \otimes \cO_K/I  \ar[u]
}
\end{equation*}
is cocartesian.\end{lem}

\pf Indeed Lemma \ref{chdomfields} implies that the square
\begin{equation*}
\xymatrix{
F \otimes \cO_K \ar[r] & F \indcot \cO_K \\
F'\otimes \cO_K \ar[u] \ar[r] & F'\indcot \cO_K \ar[u]
}
\end{equation*}
is cocartesian. Whence the result. \quod

\breakflow
According to our assumptions the action of $A$ on $\Lie_E(K)$ extends to a
continuous action of $F$ so that we have a continuous homomorphism $F \to K$.
Recall that the \emph{ramification ideal} $\rami$ is the ideal generated by $\fm_F$ in
$\cO_K$ (Definition \lref{locmod}{deflocmodcond}). In the same manner we get a
ramification ideal $\rami' \subset \cO_K$ for the coefficient subring $A' \subset A$.

\begin{prp}\label{ccbasecomp}%
Let  $\cM \subset \iHom_{A\otimes K}(M,\,a(F,K))$ be a base
compactification.
If $\cM(F \otimes \cO_K/\rami')$ is linear
then the image of $\cM$ under the restriction isomorphism
of Proposition \ref{ccihoma} is a base compactification of
$\iHom_{A'\otimes K}(M,\,a(F',K))$.\end{prp}

\pf Let $\cM'$ be the image of $\cM$ in $\iHom_{A'\otimes K}(M,\,a(F',K))$.
The shtuka $\cM$ is an $F \indcot \cO_K$-lattice in $\iHom_{A\otimes K}(M,\,a(F,K))$
so Lemma \ref{ccfokcocart} implies that $\cM'$ is an $F' \indcot \cO_K$-lattice in
$\iHom_{A'\otimes K}(M,a(F',K))$.
According to Lemma \ref{ccfokquotcocart} we have natural isomorphisms
\begin{align*}
\cM'(F'\otimes \cO_K/\fm) &\cong \cM(F\otimes \cO_K/\fm), \\
\cM'(F'\otimes \cO_K/\rami') &\cong \cM(F\otimes \cO_K/\rami).
\end{align*}
By definition $\cM(F\otimes\cO_K/\fm)$ is nilpotent while
$\cM(F\otimes\cO_K/\rami')$ is linear by assumption.
It follows that $\cM'$ is a base compactification. \quod

\section{Local models}

\begin{prp}\label{cclocmod}%
Let $\cM \subset \iHom_{A\otimes K}(M,\,a(F,K))$ be a local model.
If $\cM(F \otimes \cO_K/\rami')$ is linear
then the image of $\cM$ under the restriction isomorphism
of Proposition \ref{ccihoma} is a local model of the shtuka
$\iHom_{A'\otimes K}(M,\,a(F',K))$.\end{prp}

\pf By Proposition \lref{locmod}{locmodintersect} the local model $\cM$ is the intersection of the
coefficient compactification $\cM^c$ and a base compactification $\cM^b$.
Proposition \ref{ccdrincomp} claims that $\cM^c$ is mapped isomorphically onto the
coefficient compactification of $\iHom_{A'\otimes K}(M,\,a(F',K))$. The image of
$\cM^b$ in $\iHom_{A'\otimes K}(M,\,a(F',K))$ is a base compactification
by Propostion \ref{ccbasecomp}.
According to Proposition \lref{locmod}{locmodintersect} their intersection is a
local model. \quod

%

\begin{prp}\label{ccreg}%
Let $\cM \subset \iHom_{A\otimes K}(M,\,a(F,K))$ be a local model
such that $\cM(F \otimes \cO_K/\rami')$ is linear.
Let $\cM' \subset \iHom_{A'\otimes K}(M,\,a(F',K))$ be the image of
$\cM$ under the restriction isomorphism
of Proposition \ref{ccihoma}.
The diagram
\begin{equation*}
\xymatrix{
\uH^1(\cM) \ar[d]_{\textrm{res.}}^{\rtviso} \ar[r]^{\rho} &
\uH^1(\Der\cM) \ar[d]^{\textrm{res.}}_{\ltviso} \\
\uH^1(\cM') \ar[r]^{\rho'} & \uH^1(\Der\cM')
}
\end{equation*}
is commutative. Here $\rho$ is the regulator of $\cM$
and $\rho'$ is the regulator of $\cM'$.
\end{prp}

\pf Recall that the natural maps
\begin{align*}
\cO_F\indcot\cO_K &\to \cO_F\complot\cO_K, \\
\cO_{F'}\indcot\cO_K &\to \cO_{F'}\complot\cO_K
\end{align*}
are isomorphisms by Proposition \lref{ringot}{indcotdiscrcomplot}.
We need to prove that the regulator of $\cM$ viewed as an elliptic
$\cO_F\complot\cO_K$-shtuka of ramification ideal $\rami$
coincides with the regulator of
$\cM$ viewed as an elliptic $\cO_{F'}\complot\cO_K$-shtuka of ramification ideal $\rami'$.
The field extension $F/F'$ is totally ramified of degree $d$.
As a consequence $\rami' = \rami^d$.
The result now follows from Theorem \lref{reg}{ellregcmp}. \quod

\breakflow
Next we prove that the exponential maps of local models are stable under
restriction of coefficients.

\begin{lem}\label{ccihomb}%
The restriction map $b(F,K) \to b(F',K)$ induces an isomorphism
\begin{equation*}
\iHom_{A\otimes K}(M,\,b(F,K)) \cong \iHom_{A'\otimes K}(M,\,b(F',K))
\end{equation*}
of $F^\#\complot K$-module shtukas.\end{lem}

\pf The argument is the same as in Proposition \ref{ccihoma} save for the fact
that one needs to use Lemma \ref{chdomfieldsb} in place of Lemma
\ref{chdomfields}.\quod

\begin{lem}\label{ccfsharpcocart}%
The commutative square of rings
\begin{equation*}
\xymatrix{
\cO_F \indcot \cO_K \ar[r] & F^\#\complot \cO_K \\
\cO_{F'}\indcot \cO_K \ar[u] \ar[r] & (F')^\#\complot \cO_K \ar[u]
}
\end{equation*}
is cocartesian.\end{lem}

\pf Lemma \ref{chdomfields} implies that the square
\begin{equation*}
\xymatrix{
\cO_F \otimes \cO_K \ar[r] & \cO_F \indcot \cO_K \\
\cO_{F'} \otimes \cO_K \ar[u] \ar[r] & \cO_{F'} \indcot \cO_K \ar[u]
}
\end{equation*}
is cocartesian. At the same time the square
\begin{equation*}
\xymatrix{
F \otimes \cO_K \ar[r] & F^\#\complot\cO_K \\
F' \otimes \cO_K \ar[u] \ar[r] & (F')^\#\complot\cO_K \ar[u]
}
\end{equation*}
is cocartesian by
Lemma \ref{chdomfieldsb}. Whence the result. \quod

\begin{prp}\label{ccexp}%
Let $\cM \subset \iHom_{A\otimes K}(M,\,a(F,K))$ be a local model
such that $\cM(F \otimes \cO_K/\rami')$ is linear.
Let $\cM' \subset \iHom_{A'\otimes K}(M,\,a(F',K))$ be the image of
$\cM$ under the restriction isomorphism
of Proposition \ref{ccihoma}.
The diagram
\begin{equation*}
\xymatrix{
F \otimes_{\cO_F} \uH^1(\Der\cM) \ar[d]_{\textrm{res.}}^{\rtviso} \ar[r]^{\exp} &
F \otimes_{\cO_F} \uH^1(\cM) \ar[d]^{\textrm{res.}}_{\ltviso} \\
F'\otimes_{\cO_{F'}} \uH^1(\Der\cM') \ar[r]^{\exp'} & F'\otimes_{\cO_F'} \uH^1(\cM')
}
\end{equation*}
is commutative. Here $\exp$ is the exponential map of $\cM$
and $\exp'$ is the exponential map of $\cM'$.
\end{prp}

\pf Suppose that $\cM$ is given by a diagram
\begin{equation*}
\Big[ \cM_0 \shtuka{i}{j} \cM_1 \Big].
\end{equation*}
By definition the exponential map of $\cM$ is the composition
$\uH^1(\gamma) \circ \uH^1(\Der\gamma)^{-1}$
of the maps
\begin{align*}
\uH^1(\gamma)\colon &F \otimes_{\cO_F} \uH^1(\cM) \to \Lie_E(K), \\
\uH^1(\Der\gamma)\colon & F\otimes_{\cO_F} \uH^1(\Der\cM) \to \Lie_E(K)
\end{align*}
of Definition \lref{locmod}{gammamaps}. Similarly the exponential map
of $\cM'$ is the composition of the maps
\begin{align*}
\uH^1(\gamma')\colon &F'\otimes_{\cO_{F'}} \uH^1(\cM) \to \Lie_E(K), \\
\uH^1(\Der\gamma')\colon & F'\otimes_{\cO_{F'}} \uH^1(\Der\cM) \to \Lie_E(K).
\end{align*}
To prove the proposition it is enough to show that $\uH^1(\gamma)$ and
$\uH^1(\Der\gamma)$ are compatible with the corresponding maps of $\cM'$.

Let $c \in \uH^1(\cM)$ be a
cohomology class and let $\alpha$ be the image of $c$ under $\uH^1(\gamma)$.
Let $g \in \cM_1$ be an element representing $c$.
According to Proposition \lref{locmod}{gammadesc} there exists a unique element
\begin{equation*}
f \in \cM_0(F^\#\complot\cO_K) \subset \iHom_{A\otimes K}(M,\,b(F,K))
\end{equation*}
such that $(i-j)(f) = g$. The element $\alpha \in \Lie_E(K)$ is characterized by
the fact that for all $x$ in an open neighbourhood of
$0$ in $F$ and for all $m \in M^0$ we have $f(m)(x) = m \exp (x\alpha)$.

Lemma \ref{ccihomb} identifies
$\iHom_{A \otimes K}(M,\,b(F,K))$ with
$\iHom_{A'\otimes K}(M,\,b(F',K))$ while 
Lemma \ref{ccfsharpcocart} implies that the natural map
$\cM'((F')^\#\complot\cO_K) \to \cM(F^\#\complot\cO_K)$ is an isomorphism.
Using Proposition \lref{locmod}{gammadesc} we conclude that the square
\begin{equation*}
\xymatrix{
\uH^1(\cM) \ar[d]_{\textrm{res.}}^{\rtviso}
\ar[rr]^{\uH^1(\gamma)} && \Lie_E(K) \\
\uH^1(\cM') \ar[rr]^{\uH^1(\gamma')} && \Lie_E(K) \ar@{=}[u]
}
\end{equation*}
is commutative. The same argument shows that the square 
\begin{equation*}
\xymatrix{
\uH^1(\Der\cM) \ar[d]_{\textrm{res.}}^{\rtviso}
\ar[rr]^{\uH^1(\Der\gamma)} && \Lie_E(K) \\
\uH^1(\Der\cM') \ar[rr]^{\uH^1(\Der\gamma')} && \Lie_E(K) \ar@{=}[u]
}
\end{equation*}
is commutative. So we get the result. \quod

\chapter{Regulators of local models}
\label{chapter:locmodreg}
\label{ch:locmodreg}

Let $E$ be a Drinfeld $A$-module over $K$, a finite product of local fields containing
$\Fq$. As in Chapter \ref{chapter:locmod} we assume that the action of $A$ on $\Lie_E(K)$
extends to a continuous action of $F$. In that chapter we introduced the notion of a local
model $\cM$, its regulator $\rho\colon \uH^1(\cM) \xrightarrow{\isosign} \uH^1(\Der\cM)$ and its exponential
map $\exp\colon F \otimes_{\cO_F} \uH^1(\Der\cM) \xrightarrow{\isosign} F\otimes_{\cO_F} \uH^1(\cM)$.
Our goal is to show that the square
\begin{equation*}
\xymatrix{
\uH^1(\Der\cM) \ar[r] & F \otimes_{\cO_F} \uH^1(\Der\cM) \ar[d]^{\exp}_{\ltviso} \\
\uH^1(\cM) \ar[u]^{\rho}_{\rtviso} \ar[r] & F \otimes_{\cO_F} \uH^1(\cM)
}
\end{equation*}
commutes. Unfortunately the only proof we have proceeds by reduction to the case
$A = \Fq[t]$ (using the results of Chapter \ref{chapter:coeffch})
and a brute force computation.
On the positive side, this chapter makes the abstract
machinery of Chapter \ref{chapter:locmod} more explicit.

\section{The setting}

From now on we consider the special case $A = \Fq[t]$.
As usual $F = \Fq(\!(t^{-1})\!)$ stands for the local field of $A$ at infinity
and $\cO_F = \Fq[[t^{-1}]]$ denotes the ring of integers of $F$.
Let $K$ be a local field containing $\Fq$. As before its ring of
integers is denoted $\cO_K$ and $\fm_K \subset \cO_K$ stands for the maximal
ideal. We also fix a norm on $K$ such that $|\zeta| = q^{-1}$ for a uniformizer
$\zeta \in K$.

Fix a Drinfeld $A$-module $E$ of rank $r$ over $K$.
Let $\varphi\colon A \to K\{\tau\}$
be the ring homomorphism determined by the action of $A$ on $E$.
We define the elements
$\theta,\alpha_1,\dotsc,\alpha_r \in K$ via the equation
\begin{equation*}
\varphi(t) = \theta + \alpha_1 \tau + \dotsc + \alpha_r \tau^r.
\end{equation*}
As in Chapter \ref{chapter:locmod} we assume that the action of $A$ on
$\Lie_E(K)$ extends to a continuous action of $F$.

\begin{lem}
The following are equivalent:
\begin{enumerate}
\item The action of $A$ on $\Lie_E(K)$ extends to a continuous action of $F$.

\item $|\theta| > 1$.\quod
\end{enumerate}\end{lem}

\breakflow
According to Definition \lref{locmod}{deflocmodcond} the ramification ideal $\rami
\subset \cO_K$ is the ideal generated by $\fm_F$ under the homomorphism $F \to
K$ induced by the action of $A$ on $\Lie_E(K)$.

\begin{lem} $\rami = \theta^{-1} \cO_K$. \quod\end{lem}

\breakflow
Throughout this chapter $M = \Hom(E,\bG_a)$ stands for the motive of $E$.
Without loss of generality we assume
that the underlying group scheme of $E$ is $\bG_a$ so that we can identify
$M$ with $K\{\tau\}$.

\begin{lem}\label{polymotstruct}
The motive $M = K\{\tau\}$ has the following properties.
\begin{enumerate}
\item $M$ is a free $A \otimes K$-module of rank $r$ with a basis
$1, \dotsc, \tau^{r-1}$.

\item $M^{\geqslant 1}$ is a free $A\otimes K$-module of rank $r$ with a basis
$\tau, \dotsc, \tau^r$.

\item We have a relation
\begin{equation*}
\tau^r = \alpha_r^{-1}\big((t - \theta) \cdot 1 - \alpha_1 \tau - \alpha_2 \tau^2 - \dotsc
- \alpha_{r-1} \tau^{r-1}\big).
\end{equation*}
in the $A\otimes K$-module $M$.
\end{enumerate}\end{lem}

\pf (1) follows from Lemma \lref{drmot}{tmotfinconds} (2). In view of (1) the
result (2) is a corollary of
Proposition \lref{drmot}{motpullback}. (3) is a consequence of (1). \quod


\section{Coefficient compactification}

\begin{lem}\label{ballvanish} Let $k \geqslant 1$ be an integer.
\begin{enumerate}
\item Let $f \in \Hom(M,\,a(F,K))$. The following are equivalent:
\begin{enumerate}
\item $f(\tau^n)$ vanishes on $\cO_F$ for all $n \leqslant k r$.

\item $f(1), \dotsc, f(\tau^{r-1})$ vanish on $t^{k-1} \cO_F$ and furthermore
$f(1)$ vanishes on $t^k \cO_F$.
\end{enumerate}

\medskip\item Let
$g\in \Hom(M^{\geqslant 1},\,a(F,K))$. The following are equivalent:
\begin{enumerate}
\item $g(\tau^n)$ vanishes on $\cO_F$ for all $n \leqslant kr$.

\item $g(\tau), \dotsc, g(\tau^r)$ vanish on $t^{k-1} \cO_F$.
\end{enumerate}
\end{enumerate}
\end{lem}

\pf[\afterenum](2) From the relation
\begin{equation*}
\tau^n t = \theta^{q^n} \tau^n + \alpha_1^{q^n} \tau^{n+1} + \dotsc +
\alpha_r^{q^n} \tau^{n + r}.
\end{equation*}
one concludes that if $g(\tau^n), \dotsc, g(\tau^{n+r-1})$ vanish on an open
neighbourhood $U \subset F$ of $0$ then the following are equivalent:
\begin{itemize}
\item[(i)] $g(\tau^{n+r})$ vanishes on $U$.

\item[(ii)] $g(\tau^n)$ vanishes on $tU$.
\end{itemize}
With (i) one gets (a) $\Rightarrow$ (b) and (ii) implies (b) $\Rightarrow$
(a). The argument for (1) is the same as for (2). \quod

\begin{prp}\label{indrcomp} The $\cO_F\indcot K$-modules $\cM_0^c$, $\cM_1^c$ in the
coefficient compactification
\begin{equation*}
\Big[\cM_0^c \shtuka{i}{j} \cM_1^c\Big] \subset
\iHom(M,\,a(F,K)).
\end{equation*}
admit the following description:
\begin{equation*}
\begin{array}{rll}
\cM_0^c &= \big\{ f \colon M \to a(F,K)\!\! &|\,\,\, f(\tau), \dotsc, f(\tau^{r-1})
\textrm{ \rm vanish on } t^{-1}\cO_F \\
& &\,\,\,\textrm{ \rm and } f(1) \textrm{ \rm vanishes on } \cO_F
\big\},\\ \\
\cM_1^c &= \big\{ g \colon M \to a(F,K)\!\! &|\,\,\, g(\tau), \dotsc, g(\tau^r)
\textrm{ \rm vanish on } t^{-1} \cO_F \big\}.
\end{array}
\end{equation*}\end{prp}
\pf[]Let $C$ be the compactification of $\Spec A$.
We denote
\begin{equation*}
\cE =
\Big[\cE_{-1} \shtuka{i}{j} \cE_0\Big] \subset
\Big[M \shtuka{1}{\tau} M \Big]
\end{equation*}
the shtuka on $C\times\Spec K$
produced by Theorem \lref{drinconstr}{motivec}.
By definition $\cM_0^c$ consists of those $f\colon M \to a(F,K)$
which send the submodule
\begin{equation*}
\cE_0(\cO_F \indcot K) \subset
M \otimes_{A \otimes K} (F \indcot K) 
\end{equation*}
to $a(F/\cO_F,K) \subset a(F,K)$.

Take $k \gg 0$ such that $\cE_0(k)$ is globally generated.
Theorem \lref{drinconstr}{motivec} implies that the
$\cO_F\indcot K$-submodule
\begin{equation*}
\cE_0(k)(\cO_F \indcot K) \subset 
M \otimes_{A\otimes K}(F \indcot K)
\end{equation*}
is generated by $\uH^0(C\times\Spec K,\,\cE_0(k)) = M^{kr}$.
By Lemma \ref{ballvanish}
(1) a morphism $f\colon M \to a(F,K)$ sends this submodule to
$a(F/\cO_K, K)$ if and only if $f(1)$ vanishes on $t^k \cO_F$ and $f(\tau),
\dotsc, f(\tau^{r-1})$ vanish on $t^{k-1} \cO_F$. However
\begin{equation*}
\cE_0(\cO_F\indcot K) = t^{-k}
\cE_0(k)(\cO_F\indcot K)
\end{equation*}
so we get the result for $\cM_0^c$. The case of $\cM_1^c$ is handled in a similar
manner. \quod

\section{Explicit models}
We equip the space $b(F,K)$ and its subspace $a(F,K)$ with the sup-norm induced
by the norm on $K$.

\begin{lem}\label{afkgen}
Let $\mu\colon F \to \Fq$ be a nonzero continuous $\Fq$-linear map.
\begin{enumerate}
\item $\mu$ generates $a(F,K)$ as an $F\indcot K$-module.

\item For every $\beta \in K^\times$ the subspace
$\{ f : |f| \leqslant |\beta| \} \subset a(F,K)$
is a free $F \indcot\cO_K$-submodule generated by $\beta \cdot \mu$.
\end{enumerate}
\end{lem}

\pf (1) Let $F^*$ be the continuous $\bF_q$-linear dual of $F$. According
to Theorem \lref{resdual}{resdual} the topological $F$-module $F^*$ is free of
rank one. So a nonzero function $\mu\colon F \to \Fq$ generates
$F^*$ as an $F$-module.
Now Corollary~\ref{funcomegalat} 
implies that $F^*$ is an $F$-lattice in the
$F\indcot K$-module $a(F,K)$ so we get the result.

(2)
Without loss of generality we assume that $\beta = 1$.
In this case the subspace of $a(F,K)$ in question
is $a(F,\cO_K)$.
By Corollary~\ref{funcomegalat} 
the space
$F^*$ is an $F$-lattice in the $F\indcot\cO_K$-module $a(F,\cO_K)$.
Whence the result. 
\quod

%
%

\begin{dfn} Let $b, b_1, \dotsc, b_{r-1} \in q^\bZ$ be real numbers.
We introduce conditions on $A \otimes K$-linear maps
$f\colon M \to b(F,K)$ and $g\colon M^{\geqslant 1}\to b(F,K)$:
\begin{align}
\label{EM0}\tag{EM$_0$}
&|f(1)| \leqslant |\alpha_r \theta^{-1}|b, \quad
|f(\tau)| \leqslant b_1, \quad \dotsc\quad
|f(\tau^{r-1})| \leqslant b_{r-1}, \\
\label{EM1}\tag{EM$_1$}
&|g(\tau)| \leqslant b_1, \quad \dotsc\quad
|g(\tau^{r-1})| \leqslant b_{r-1}, \quad
|g(\tau^r)| \leqslant b.
\end{align}\end{dfn}

\begin{prp}[\aftereq]\label{explbcwell} Let $b, b_1 \dotsc b_{r-1} \in q^\bZ$ be real
numbers. Consider the $F\indcot\cO_K$-submodules
\begin{align*}
\cM_0 &= \big\{ \,\, f \in \Hom(M,a(F,K)) \mid f \textit{ satisfies \uref{EM0}}\,\, \big\}, \\
\cM_1 &= \big\{ \,\, g \in \Hom(M^{\geqslant 1},a(F,K)) \mid g \textit{ satisfies \uref{EM1}} \,\, \big\}.
\end{align*}
If $b, b_1 \dotsc b_{r-1}$ satisfy 
the inequalities
\begin{align}\label{EMP1}\tag{EP$_1$}
&\Big|\frac{\alpha_n}{\alpha_r}\Big| \frac{b_n}{b} \leqslant 1,
\quad n  \in \{1, \dotsc, r -1 \} \\
\label{EMP2}\tag{EP$_2$}
\def\arraystretch{2.2}
&\Big|\frac{\alpha_r^q}{\theta^q}\Big| \frac{b^q}{b_1} \leqslant |\theta^{-1}|; \quad
\frac{b_{r-1}^q}{b} \leqslant |\theta^{-1}|; \quad
\frac{b_n^q}{b_{n+1}} \leqslant |\theta^{-1}|, \quad n \in \{1, \dotsc, r-2 \}
\end{align}
then $\cM_0$, $\cM_1$ define an $F\indcot\cO_K$-subshtuka
\begin{equation*}
\cM = [\cM_0 \shtuka{}{} \cM_1] \subset \iHom(M,\,a(F,K))
\end{equation*}
which is a base compactification in the sense of
Definition \lref{locmod}{defbcomp}.
\end{prp}
\pf 
It will be convenient for us to work with bases of
various modules. Our first goal is to construct bases for $\Hom(M,a(F,K))$ and
$\Hom(M^{\geqslant 1}, a(F,K))$.

According to Lemma \ref{polymotstruct} the elements
$1, \tau, \dotsc, \tau^{r-1}$
form
an $A \otimes K$-basis of~$M$.
Similarly the elements
$\tau, \dotsc, \tau^r$
form a basis of $M^{\geqslant 1}$.
Fix a nonzero continuous $\Fq$-linear map $\mu\colon F \to \Fq$.
Lemma \ref{afkgen} (1) shows that $\mu$ generates $a(F,K)$ as an $F\indcot K$-module.
We thus get $F\indcot K$-module bases of
$\Hom(M,a(F,K))$,
$\Hom(M^{\geqslant 1},a(F,K))$
which are dual to
the aforementioned bases of $M$.

Let $\beta, \beta_1 \dotsc \beta_{r-1} \in K^\times$ be such that
$|\beta| = b, |\beta_1| = b_1 \dotsc |\beta_{r-1}| = b_{r-1}$.
Lemma \ref{afkgen} (2) implies
that the modules
$\cM_0$, $\cM_1$ have $F\indcot\cO_K$-bases given by matrices
\begin{equation*}
\begin{pmatrix}
\alpha_r \theta^{-1}\beta & 0 & & 0 \\
0 & \beta_1 &  & 0 \\
  &   & \ddots &   \\
0 & 0 &  & \beta_{r-1} \\
\end{pmatrix}, \quad
\begin{pmatrix}
\beta_1 & & 0 & 0 \\
  & \ddots &   &   \\
0 &  & \beta_{r-1} & 0 \\
0 &  & 0 & \beta \\
\end{pmatrix}.
\end{equation*}
with respect to the fixed bases of $\Hom(M,a(F,K))$ and
$\Hom(M^{\geqslant 1}, a(F,K))$. We conclude that $\cM_0$, $\cM_1$ are
$F\indcot\cO_K$-lattices in these modules.

Let $i$ and $j$ be the arrows of the shtuka
\begin{equation*}
\iHom(M,\,a(F,K)) =
\Big[ \Hom(M,a(F,K)) \shtuka{i}{j} \Hom(M^{\geqslant 1}, a(F,K) \Big].
\end{equation*}
According to Proposition~\lref{ihommotdesc}{ihommotdesc} the map $i$ is the restriction to
$M^{\geqslant 1}$.
Hence the matrix of the map $i$ in the fixed bases of $\iHom(M,\,a(F,K))$ is
\begin{equation*}
\begin{pmatrix}
0 & 1 & & 0 \\
  &   & \ddots & \\
0 & 0 & & 1 \\
\frac{t - \theta}{\alpha_r} & -\frac{\alpha_1}{\alpha_r} &  & -\frac{\alpha_{r-1}}{\alpha_r}
\end{pmatrix}
\end{equation*}
Rewriting it in the bases of $\cM_0$, $\cM_1$ gives
\begin{equation*}
\begin{pmatrix}
0 & 1 & & 0 \\
  &   & \ddots & \\
0 & 0 & & 1 \\
t \theta^{-1} - 1 & -\frac{\alpha_1 \beta_1}{\alpha_r \beta} & &
-\frac{\alpha_{r-1}\beta_{r-1}}{\alpha_r \beta}
\end{pmatrix}
\end{equation*}
The assumptions $|\theta| > 1$ and \uref{EMP1} imply that the bottom row is in
$F\indcot\cO_K$ and therefore $i(\cM_0) \subset \cM_1$.  The determinant of the
matrix is $(-1)^r (1 - \theta^{-1} t)$. As $|\theta| > 1$ it reduces to $(-1)^r$
modulo $F \indcot\fm_K$ and in particular $i$ is an isomorphism modulo
$F\indcot\fm_K$.

According to Proposition~\lref{ihommotdesc}{ihommotdesc} the map $j$ sends $f \in \Hom(M,a(F,K))$ to
the map
\begin{equation*}
\tau m \mapsto \tau f(m).
\end{equation*}
Hence the matrix of $j$ in the bases of $\cM_0$, $\cM_1$ is
\begin{equation*}
\begin{pmatrix}
\frac{\alpha_r^q \beta^q}{\theta^q \beta_1} & 0 &  & 0 & 0 \\
0 & \frac{\beta_1^q}{\beta_2} &  & 0 & 0 \\
  &   &  \ddots &   &   \\
0 & 0 &  & \frac{\beta_{r-2}^q}{\beta_{r-1}} & 0 \\
0 & 0 &  & 0 & \frac{\beta_{r-1}^q}{\beta}
\end{pmatrix}
\end{equation*}
In our case the ramification ideal $\rami$ is the ideal
$\theta^{-1} \cO_K$.
Hence the assumptions \uref{EMP2} imply that the matrix above lies in $\rami$ whence
$j(\cM_0) \subset \cM_1$ and $j$ reduces to zero modulo $\rami$. Above we
demonstrated that $i$ reduces to an isomorphism modulo $\fm_K$.
Hence $\cM$ is a base compactification as claimed.
\quod

\begin{prp}\label{explexists} There exist $b, b_1, \dotsc, b_{r-1} \in q^\bZ$
satisfying \uref{EMP1} and \uref{EMP2}.\end{prp}

\pf A direct verification shows that for all
$\varepsilon \in q^\bZ$ small enough
the real numbers
\begin{equation*}
b = \varepsilon^2, \,\,
b_1 = \varepsilon^3, \dotsc,
b_{r-1} = \varepsilon^3
\end{equation*}
satisfy \uref{EMP1} and \uref{EMP2}. \quod

\begin{dfn}\label{defexpl} Let $b, b_1, \dotsc, b_{r-1} \in q^\bZ$ be real numbers
satisfying \uref{EMP1}, \uref{EMP2}.
The \emph{explicit model} $\cM$ of parameters $b, b_1
\dotsc b_{r-1}$ is the subshtuka in
$\iHom(M,\,a(F,K))$ defined as the intersection
\begin{equation*}
\cM = \cM^c \cap \cM^b
\end{equation*}
where $\cM^c$ is the coefficient compactification of Definition
\lref{locmoddrinfeld}{locmoddrincomp} and $\cM^b$
is the subshtuka of $\iHom(M,\,a(F,K))$ desribed in Proposition \ref{explbcwell}.
The number $b$ will be called the \emph{leading parameter}.\end{dfn}

\breakflow
It follows from Proposition \lref{locmod}{locmodintersect} that an explicit model is a local
model in the sense of Definition \lref{locmod}{deflocmod}. 
In particular we have the twists $\rami^n\cM$ of a local model
as in Proposition \lref{locmod}{locmodtwist}.

\begin{prp}\label{expltwist}
If $\cM$ is an explicit model of parameters $b, b_1, \dotsc, b_{r-1}$
then $\rami\cM$ is an explicit model of parameters
$b|\theta^{-1}|, b_1|\theta^{-1}|, \dotsc, b_{r-1}|\theta^{-1}|$.\end{prp}

\pf Follows since $\rami = \theta^{-1}\cO_K$.\quod

\breakflow
Next we prove a technical statement on explicit models as
lattices in the shtuka $\iHom(M,\,b(F,K))$.

\begin{lem}\label{afkbfkboundlat}
For every $\beta \in K^\times$
the subspace
$\{f \in a(F,K) : |f| \leqslant |\beta| \}$
is an $F\indcot\cO_K$-lattice in the
$F^\#\complot\cO_K$-module
$\{f \in b(F,K) : |f| \leqslant |\beta| \}$.
\end{lem}

\pf Without loss of generality we assume that $\beta = 1$.
In this case the spaces in question are $a(F,\cO_K)$
and $b(F,\cO_K)$
so the result follows from Corollary~\ref{funclat}.\quod

\breakflow
If $\cM \subset \iHom(M,\,a(F,K))$ is a local model then
Lemma \lref{locmod}{locmodblat} shows that the natural map
\begin{equation*}
\cM(F^\#\complot\cO_K) \to \iHom(M,\,b(F,K))
\end{equation*}
is an inclusion of an
$F^\#\complot\cO_K$-lattice. In particular we can view
$\cM(F^\#\complot\cO_K)$ as a subshtuka of $\iHom(M,\,b(F,K))$.

\begin{lem}\label{explbounddesc}
If $\cM = [ \cM_0 \shtuka{}{} \cM_1] $ is an explicit model then
\begin{align*}
\cM_0(F^\#\complot\cO_K) &=
\big\{ \,\, f \in \Hom(M,\,b(F,K)) \mid
f \textit{ satisfies \uref{EM0}} \,\, \big\}, \\
\cM_1(F^\#\complot\cO_K) &=
\big\{ \,\, g \in \Hom(M^{\geqslant 1},\,b(F,K)) \mid
g \textit{ satisfies \uref{EM1}} \,\, \big\}.
\end{align*}
\end{lem}

\pf Follows immediately from Lemma \ref{afkbfkboundlat}. \quod

\section{Exponential maps of Drinfeld modules}

As we identified the underlying group scheme of the Drinfeld module $E$ with $\bG_a$
we can view the exponential map of $E$ as a function $K \to K$.
The map $\exp\colon K \to K$ is a local
isomorphism of topological $\Fq$-vector spaces. In the following we will need a
slightly more precise version of this result.

\begin{lem}\label{expbound} There exists a constant $B_{\exp} \in \bR$ such that
\begin{enumerate}
\item $0 < B_{\exp} \leqslant \frac{1}{q}$.

\item If $|z| \leqslant B_{\exp}$ then $|\exp z - z| \leqslant |z|^q$. In
particular $|\exp z| = |z|$.

\end{enumerate}
\end{lem}

\pf Indeed $\exp z$ is given by an everywhere convergent power series
\begin{equation*}
\exp z = z + a_1 z^q + a_2 z^{q^2} + \dotsc
\end{equation*}
Therefore
$\exp z - z = z^q (a_1 + a_2 z^q + \dotsc)$.
The power series in brackets also converges everywhere and defines a continuous
function from $K$ to $K$. Hence there exists a nonzero constant $B_{\exp}$ such
that as soon as $|z| \leqslant B_{\exp}$ the values of this function have norm
less than or equal to $1$. Without loss of generality we may assume that
$B_{\exp} \leqslant \frac{1}{q}$. We thus get (1) and (2).
\quod

\breakflow
Let $B_{\exp}$ be a constant as in Lemma \ref{expbound} and
let $U \subset K$ be the ball of radius $B_{\exp}$ around $0$.
Lemma \ref{expbound} (2) implies that $\exp(U) = U$ and the induced map
$\exp\colon U \to U$ is an isomorphism of topological $\Fq$-vector spaces.

\begin{dfn}\label{deflog}
In the following
we denote $\log\colon U \to U$ the inverse of $\exp$ on $U$.
We call it the \emph{logarithmic map} of the Drinfeld module $E$.\end{dfn}


\breakflow
Denote $\varphi_t\colon K \to K$ the map given by the action of
$t$ on $K = E(K)$. In other words
$\varphi_t(z) = \theta z + \alpha_1 z^q + \dotsc + \alpha_r z^{q^r}$.

\begin{lem}\label{expllog}%
Let $z \in K$.
If $|\varphi_t(z)| \leqslant B_{\exp}$ and $|z|\leqslant B_{\exp}$ then
$\log(\varphi_t(z)) = \theta \log(z)$.\quod
\end{lem}


\section{Exponential maps of explicit models}

Let $\cM$ be an explicit model. In Section \lref{locmodcoh}{sec:locmodexp} we
introduced the exponential map of $\cM$:
\begin{equation*}
\exp\colon F \otimes_{\cO_F} \uH^1(\Der\cM) \to F \otimes_{\cO_F} \uH^1(\cM).
\end{equation*}
In this section we will describe this map on explicit representatives
of cohomology classes.


\begin{dfn}\label{spcohcl} Let $h \in b(F, K)$.
We define $A\otimes K$-linear maps
\begin{align*}
\Der g(h)&\colon M^{\geqslant 1} \to b(F,K), \\
g(h)&\colon M^{\geqslant 1} \to b(F,K)
\end{align*}
by prescribing them on the basis
$\tau, \dotsc, \tau^r$ of $M^{\geqslant 1}$ as follows:
\begin{align*}
\Der g(h)\colon &\tau^1, \dotsc, \tau^{r-1} \mapsto  0, \\
&\tau^r \mapsto \alpha_r^{-1} (ht - \theta h), \\
g(h)\colon &\tau^1, \dotsc, \tau^{r-1} \mapsto  0, \\
&\tau^r \mapsto \alpha_r^{-1} \exp \circ (ht - \theta h).
\end{align*}
where $\exp\colon K \to K$ is the exponential map of the Drinfeld module $E$.
Note that $g(h)(\tau^r)$ maps $x \in F$ to
$\alpha_r^{-1} \exp( h(tx) - \theta h(x))$ in $K$.

We define an $A\otimes K$-linear map
$\Der f(h)\colon M \to b(F,K)$ by prescribing it on the basis
$1, \dotsc, \tau^{r-1}$ of $M$ as follows:
\begin{equation*}
\Der f(h)\colon 1 \mapsto h, \quad \tau,\dotsc,\tau^{r-1} \mapsto 0.
\end{equation*}
\end{dfn}

\breakflow
We will use $g(h)$ and $\Der g(h)$ as representatives for classes in the
cohomology groups $\uH^1(\cM)$ and $\uH^1(\Der \cM)$ respectively.

\begin{lem}\label{explreplieaux} Let $i$, $j$ be the arrows of the shtuka
\begin{equation*}
\iHom(M,\,b(F,K)) =
\Big[ \Hom(M,b(F,K)) \shtuka{i}{j} \Hom(M^{\geqslant 1}, b(F,K)) \Big].
\end{equation*}
If $h \in b(F,K)$ then $i(\Der f(h)) = \Der g(h)$.
\end{lem}

\pf According to Proposition~\lref{ihommotdesc}{ihommotdesc} the map $i(\Der f(h))$ is the
restriction of $\Der f(h)$
to $M^{\geqslant 1}$. Applying $\Der f(h)$ to the relation
\begin{equation*}
\tau^r = \alpha_r^{-1} \big( (t-\theta) \cdot 1
- \alpha_1 \tau - \alpha_2 \tau^2 - \dotsc - \alpha_{r-1} \tau^{r-1} \big)
\end{equation*}
we conclude that
$\Der f(h)(\tau^r) = \alpha_r^{-1}(h t-\theta h)$. \quod

\begin{lem}\label{explreplieall}
Let $\cM$ be an
explicit model given by a diagram
\begin{equation*}
\cM_0 \shtuka{i}{j} \cM_1
\end{equation*}
and let $c \in \uH^1(\Der\cM)$ be a cohomology class. There
exists a function $h \in b(F,K)$ such that the following holds:
\begin{enumerate}
\item $\Der g(h)$ belongs to $\cM_1$ and represents $c$
in $\uH^1(\Der\cM) = \coker(i)$.

\item $\Der f(h) \in \cM_0(F^\#\complot\cO_K)$.
\end{enumerate}
\end{lem}

\pf 
Let $g \in \cM_1 \subset \Hom(M^{\geqslant 1}, a(F,K))$ be a
representative of the cohomology class $c$. According to Proposition
\lref{locmodcoh}{liegammadesc} there exists a unique $f \in
\cM_0(F^\#\complot\cO_K) \subset \Hom(M,b(F,K))$ such that $i(f) = g$.
Set $h = f(1)$.

Since $f$ is an element of $\cM_0(F^\#\complot\cO_K)$
Lemma \ref{explbounddesc} shows that
the function $h$ satisfies $|h| \leqslant |\alpha_r \theta^{-1}| b$.
The same lemma then implies that
$\Der f(h) \in \cM_0(F^\#\complot\cO_K)$ which is
the claim (2). In order to deduce (1) it is enough to show that
the difference
\begin{equation*}
\delta = f - \Der f(h)
\end{equation*}
is an element of $\cM_0$.
Indeed $i(\Der f(h)) = \Der g(h)$ by Lemma \ref{explreplieaux} so
$\Der g(h) = g - i(\delta)$ is an element of $\cM_1$ which represents the same
cohomology class as $g$.

According to Proposition~\lref{ihommotdesc}{ihommotdesc} the map $i(f) = g$ is
the restriction of $f$ to $M^{\geqslant 1}$. Therefore the difference $\delta$ acts
on the generators of $M$ as follows:
\begin{equation*}
\delta \colon 1 \mapsto 0, \,\,
\tau^n \mapsto g(\tau^n), \,\, n \in \{1, \dotsc, r-1\}.
\end{equation*}
Let $\cM^c = [ \cM^c_0 \shtuka{}{} \cM^c_1 ]$ be the coefficient compactification.
By definition of an explicit model
\begin{equation*}
\cM_0 = \cM^c_0 \cap \cM_0(F \indcot \cO_K).
\end{equation*}
We first prove that $\delta \in \cM^c_0$.
As $g \in \cM_1^c$ Proposition \ref{indrcomp} shows that the functions $g(\tau), \dotsc,
g(\tau^r)$ vanish on $t^{-1} \cO_F$. Since $\delta$ maps $1$ to $0$ the same
Proposition implies that $\delta \in \cM_0^c$.

Now we prove that $\delta \in \cM_0(F \indcot\cO_K)$.
By definition $g$ is an element of $\cM_1$ so \uref{EM1} implies that
\begin{equation*}
|\delta(\tau^n)| = |g(\tau^n)| \leqslant b_n, \quad
n \in \{1, \dotsc, r-1\}
\end{equation*}
where $b_n$ are the parameters of $\cM$.
Moreover $\delta(1) = 0$ so 
$\delta \in \cM_0(F \indcot\cO_K)$ by \uref{EM0}.\quod


\breakflow
In the following let us fix a constant $B_{\exp}$ satisfying the assumptions of
Lemma \ref{expbound}.

\begin{lem}\label{explreplog} Let $\cM$ be an
explicit model given by a diagram
\begin{equation*}
\cM_0 \shtuka{i}{j} \cM_1.
\end{equation*}
Let $h \in b(F,K)$ be a function such that
\begin{equation*}
\Der g(h), g(h) \in \cM_1, \quad
\Der f(h) \in \cM_0(F^\#\complot\cO_K).
\end{equation*}
Assume that the leading parameter $b$ of $\cM$ satisfies
$|\alpha_r| b \leqslant B_{\exp}$.
If an element
$f \in \cM_0(F^\#\complot\cO_K)$ 
satisfies $(i-j)(f) = g(h)$ then $f(1) = \exp h$.\end{lem}

\pf Set $h_1 = f(1)$. As before we denote $\varphi_t\colon K \to K$ 
the action of $t$ on $K = E(K)$. We have
\begin{equation*}
\varphi_t(z) = \theta z + \alpha_1 z^q + \dotsc + \alpha_r z^{q^r}.
\end{equation*}
We split the proof in several steps.

\medbreak
\textbf{Step 1.}
$h_1 \circ t - \varphi_t \circ h_1 = \exp \circ (h t - \theta h)$
\textit{in} $b(F,K)$.

According to Proposition~\lref{ihommotdesc}{ihommotdesc} the map $(i-j)(f)$ satisfies
\begin{equation*}
(i-j)(f)\colon \tau^{n+1} \mapsto f(\tau^{n+1}) - \tau f(\tau^n)
\end{equation*}
for all $n \geqslant 0$. Comparing with the definition of $g(h)$ we deduce that
\begin{equation}\label{explftaun}
f(\tau^n) = \tau^n f(1) = \tau^n h_1
\end{equation}
for all $n \in \{0, \dotsc, r - 1\}$ and that
\begin{equation}\label{explftaur}
f(\tau^r) = \tau^r h_1 +  \alpha_r^{-1} \exp (h t - \theta h).
\end{equation}
Appling the $A \otimes K$-linear map $f$ to the relation
\begin{equation*}
1 \cdot t - \theta \cdot 1 = \alpha_1 \tau  + \dotsc + \alpha_r \tau^r
\end{equation*}
in $M$,
we obtain what we need.

\medbreak
\textbf{Step 2.}
\textit{The function
$\varphi_t \circ h_1$ in $b(F,K)$
satisfies $|\varphi_t \circ h_1| \leqslant |\alpha_r| b$.}

We prove it estimating the expression
\begin{equation*}
\varphi_t \circ h_1 = \theta h_1 + \alpha_1 \tau h_1 + \dotsc + \alpha_r \tau^r
h_1
\end{equation*}
term by term.
As $f \in \cM_0(F^\#\complot\cO_K)$
Lemma \ref{explbounddesc} implies that
$|h_1| \leqslant |\alpha_r \theta^{-1} b|$. Hence
\begin{equation}\label{explthetah1b}
|\theta h_1|
\leqslant |\alpha_r| b.
\end{equation}
Next let $b_1, \dotsc, b_{r-1}$ be the parameters of $\cM$ and
let $n \in \{1, \dotsc, r-1\}$. According to \eqref{explftaun}
we have $f(\tau^n) = \tau^n h_1$.
The map $f$ belongs to $\cM_0(F^\#\complot\cO_K)$.
Hence Lemma \ref{explbounddesc} shows that
\begin{equation*}
|f(\tau^n)| \leqslant b_n.
\end{equation*}
However the conditions \uref{EMP1} imply that
\begin{equation*}
|\alpha_n| b_n \leqslant |\alpha_r| b.
\end{equation*}
One thus obtains the inequality
\begin{equation}\label{expltaunh1b}
|\alpha_n \tau^n h_1| \leqslant |\alpha_n| b_n \leqslant |\alpha_r| b.
\end{equation}
It remains to estimate $\alpha_r \tau^r h_1$. The equation \eqref{explftaur}
implies that
\begin{equation*}
\tau^r h_1 = f(\tau^r) - g(h)(\tau^r).
\end{equation*}
The element $f$ belongs to $\cM_0(F^\#\complot\cO_K)$ by assumption. Hence
$i(f) \in \cM_1(F^\#\complot\cO_K)$.
According to Proposition~\lref{ihommotdesc}{ihommotdesc} the map $i(f)$ is the
restriction of $f$ to $M^{\geqslant 1}$. Therefore
Lemma \ref{explbounddesc} implies that
$|f(\tau^r)| \leqslant b$.
As $g(h) \in \cM_1$ we deduce that $|g(h)(\tau^r)| \leqslant b$. Hence
\begin{equation}\label{expltaurh1b}
|\tau^r h_1| \leqslant b.
\end{equation}
Combining \eqref{explthetah1b}, \eqref{expltaunh1b}, \eqref{expltaurh1b} we get the inequality
$|\varphi_t \circ h_1|
\leqslant |\alpha_r| b$.

\medbreak
\textbf{Step 3.} \emph{$\log h_1 \in b(F,K)$ is well-defined}.

The logarithmic map is defined on the ball of radius $B_{\exp}$ around $0$
(see Definition \ref{deflog}).
The function $h_1$ satisfies $|h_1| \leqslant |\alpha_r \theta^{-1}| b$ because $f
\in \cM_0(F^\#\complot\cO_K)$.
Since $|\theta| > 1$ and $|\alpha_r| b \leqslant B_{\exp}$ by assumption we
conclude that $\log h_1$ is well-defined.

\medbreak
\textbf{Step 4.} $\Der g(\log h_1) = \Der g(h)$.

According to Step 1 we have
\begin{equation}\label{explhtlog}
h_1 \circ t - \varphi_t \circ h_1 = \exp \circ (h t - \theta h).
\end{equation}
We would like to apply $\log$ to this equation.
Step 2 shows that
\begin{equation*}
|\varphi_t \circ
h_1| \leqslant |\alpha_r| b \leqslant B_{\exp}
\end{equation*}
so 
$\log\varphi_t h_1$ is well-defined. Moreover $\log h_1$ is well-defined by
Step 3. 
Now $f(h) \in \cM_0(F^\#\complot\cO_K)$ by assumption
so Lemma \ref{explbounddesc} shows that
\begin{equation*}
|h| \leqslant |\alpha_r \theta^{-1}|b \leqslant |\theta^{-1}| B_{\exp}.
\end{equation*}
According to Lemma \ref{expbound} we have
$|\exp(z)| = |z|$ provided $|z| \leqslant B_{\exp}$.
Therefore
\begin{equation*}
|\exp \circ (ht - \theta h)| \leqslant |ht - \theta h| \leqslant B_{\exp}
\end{equation*}
and we conclude that
$\log(\exp \circ (ht - \theta h))$ is well-defined.
Applying $\log$ to \eqref{explhtlog} we get 
\begin{equation*}
\log h_1 t - \log \varphi_t h_1 = h t  - \theta h.
\end{equation*}
Lemma \ref{expllog} shows that
$\log(\varphi_t(z)) = \theta \log z$
provided $\varphi_t(z)$ and $z$ are in the domain of definition of $\log$.
As a consequence
$\log\varphi_t h_1 = \theta\log h_1$
and
\begin{equation*}
\log h_1 t - \theta\log h_1 = h t - \theta h.
\end{equation*}
It follows that
$\Der g(\log h_1) = \Der g(h)$
by definition of $\Der g$ and $g$.

\medbreak
\textbf{Step 5.} $\Der f(\log h_1) \in \cM_0(F^\#\complot\cO_K)$.

From Lemma \ref{expbound} it follows that $|\log z| = |z|$ for all $z \in K$ satisfying $|z|
\leqslant B_{\exp}$. Thus $|\log h_1| = |h_1|
\leqslant |\alpha_r \theta^{-1}| b$. Applying Lemma \ref{explbounddesc} we conclude that
$\Der f(\log h_1) \in \cM_0(F^\#\complot\cO_K)$.

\medbreak
\textbf{Step 6.} $h_1 = \exp h$. According to Lemma \ref{explreplieaux}
\begin{align*}
i(\Der f(h)) &= \Der g(h),\\
i(\Der f(\log h_1)) &= \Der g(\log h_1).
\end{align*}
Furthermore $\Der g(\log h_1) = g(h)$ by Step 4. Now $\Der f(h)$ belongs to
$\cM_0(F^\#\complot\cO_K)$ by assumption while $\Der f(\log h_1) \in
\cM_0(F^\#\complot\cO_K)$ by Step 5. Hence the unicity part of Proposition
\lref{locmodcoh}{gammadesc} implies that $\Der f(h) = \Der f(\log h_1)$. 
By definition of $\Der f$ it means that $h = \log h_1$. Therefore $h_1 = \exp h$. \quod

\begin{lem}\label{explrepdiff} Let $\cM = [\cM_0 \shtuka{}{} \cM_1]$ be an
explicit model with leading parameter $b$ satisfying $|\alpha_r| b \leqslant
B_{\exp}$. If $h \in b(F,K)$ is a function such that
$\Der g(h) \in \cM_1(F\indcot\cO_K)$ then
\begin{equation*}
\Der g(h) - g(h) \in I \cM_1(F\indcot\cO_K)
\end{equation*}
where
$I = \{ x \in K : |x| \leqslant (|\alpha_r| b)^{q-1} \} \subset \cO_K$.
\end{lem}

\pf The fact that
$\Der g(h)$ is an element of $\cM_1(F\indcot\cO_K)$ implies a bound
\begin{equation*}
|ht-\theta h| \leqslant |\alpha_r| b.
\end{equation*}
According to Lemma \ref{expbound} $|\exp z - z| \leqslant |z|^q$ as soon as
$|z| \leqslant B_{\exp}$. As $|\alpha_r| b \leqslant B_{\exp}$ we conclude that
\begin{equation*}
|\exp(h t-\theta h) - (h t-\theta h)| \leqslant |h t-\theta h|^q \leqslant
(|\alpha_r| b)^q.
\end{equation*}
Therefore
\begin{equation*}
|\Der g(h)(\tau^r) - g(h)(\tau^r)| \leqslant (|\alpha_r| b)^{q-1} \cdot b
\end{equation*}
which implies that $\Der g(h) - g(h) \in I \cM_1(F\indcot\cO_K)$. \quod

\begin{prp}\label{explmodexp} Let $\cM = [\cM_0 \shtuka{}{} \cM_1]$ be an explicit model with
leading parameter $b$. If
$|\alpha_r| b \leqslant B_{\exp}$ then for
every cohomology class $c \in \uH^1(\Der\cM)$ there exists $h \in b(F,K)$ such
that the following holds.
\begin{enumerate}
\item $\Der g(h)$ belongs to $\cM_1$ and represents $c$.

\item $g(h)$ belongs to $\cM_1$. 

\item $\exp[ \Der g(h) ] = [ g(h) ]$ in $F \otimes_{\cO_F} \uH^1(\cM)$.

\item $\Der g(h) - g(h) \in I\cM_1$ where
$I = \{ x \in K : |x| \leqslant (|\alpha_r| b)^{q-1} \} \subset \cO_K$.
\end{enumerate}
Here the brackets $[\,]$ denote cohomology
classes and $\exp$ is the exponential map of $\cM$
as in Definition \lref{locmod}{defexpmap}.
\end{prp}

\pf According to Lemma \ref{explreplieall} there exists an $h \in b(F,K)$ such
that $\Der g(h) \in \cM_1$ represents the class $c$ and $\Der f(h) \in
\cM_0(F^\#\complot\cO_K)$.

We claim that $\Der g(h) - g(h) \in I\cM_1$. To prove it we consider the
coefficient compactification $\cM^c = [\cM_0^c \shtuka{}{} \cM_1^c]$. As $\Der
g(h) \in \cM_1^c$
Proposition \ref{indrcomp} implies that $g(h) \in \cM_1^c$.
Lemma \ref{explrepdiff} shows that
\begin{equation*}
\Der g(h) - g(h) \in I \cM_1(F\indcot\cO_K)
\end{equation*}
By definition of an explicit model
\begin{equation*}
\cM_1^c \cap I \cM_1(F\indcot\cO_K) = I \cM_1.
\end{equation*}
Hence $\Der g(h) - g(h) \in I \cM_1$ and $g(h) \in \cM_1$.

It remains to prove that $\exp[\Der g(h)] = [g(h)]$.
Consider the isomorphisms
\begin{align*}
\gamma &\colon F \otimes_{\cO_F} \RGamma(\cM) \xrightarrow{\,\isosign\,} \Lie_E(K)[-1], \\
\Der\gamma &\colon F \otimes_{\cO_F} \RGamma(\Der\cM) \xrightarrow{\,\isosign\,} \Lie_E(K)[-1]
\end{align*}
of Definition \lref{locmodcoh}{gammamaps}. By Definition
\lref{locmodcoh}{defexpmap} the exponential map of $\cM$ is the composition
\begin{equation*}
\uH^1(\gamma) \circ \uH^1(\Der\gamma)^{-1}.
\end{equation*}
Let $\alpha \in \Lie_E(K)$ be such that $\uH^1(\Der\gamma)(c) = \alpha$.
Proposition \lref{locmodcoh}{liegammadesc} shows that the element $\alpha$ is
characterized by the following property:
\begin{equation*}
\Der f(h)(1) = h \sim (x \mapsto x\alpha).
\end{equation*}
Let $i$ and $j$ be the arrows of $\cM$.
According to Proposition \lref{locmod}{gammadesc} there exists a unique $f \in
\cM_0(F^\#\complot\cO_K)$ such that $(i-j)f = g(h)$.
Lemma \ref{explreplog} tells us that $f(1) = \exp h$.
Hence
\begin{equation*}
f(1) \sim (x \mapsto \exp (x\alpha)).
\end{equation*}
Thus $\uH^1(\gamma)([g(h)]) = \alpha$ by Proposition \lref{locmodcoh}{gammadesc}.
We conclude that $\exp [\Der g(h)] = [g(h)]$. \quod

\section{Regulators of explicit models}

Let $\cM$ be a local model.
By Theorem \lref{locmodcoh}{locmodcohdesc}
the cohomology modules
$\uH^1(\cM)$, $\uH^1(\Der\cM)$ are free $\cO_F$-modules of finite rank.

\begin{lem}
Let $\cM$ be an explicit model with leading parameter
$b$. If $|\alpha_r| b \leqslant B_{\exp}$ then
the exponential map
$\exp\colon F \otimes_{\cO_F} \uH^1(\Der\cM) \to
F \otimes_{\cO_F} \uH^1(\cM)$
sends the $\cO_F$-submodule $\uH^1(\Der\cM)$ to
$\uH^1(\cM)$.
\end{lem}


\pf Immediate from Proposition \ref{explmodexp}.\quod

\breakflow
As in Section \lref{reg}{sec:elltwist} we denote $\cM/\rami^n$ the quotient
$\cM/(\rami^n\cM)$.

\begin{lem}\label{explslide}%
Let $\cM$ be an explicit model with leading
parameter $b$ satisfying $|\alpha_r| b \leqslant B_{\exp}$. Let
\begin{equation*}
I = \{ x \in K : |x| \leqslant (|\alpha_r| b)^{q-1} \} \subset \cO_K.
\end{equation*}
Let $n \geqslant 0$. If $I \subset \rami^n$
and $\cM/\rami^n$ is linear
then the following square is commutative:
\begin{equation*}
\xymatrix{
\uH^1(\Der\cM) \ar[r] \ar[d]^{\exp}& \uH^1(\Der\cM/\rami^n) \ar[d]^1 \\
\uH^1(\cM) \ar[r] & \uH^1(\cM/\rami^n).
}
\end{equation*}
Here the horizontal arrows are induced by reduction modulo $\rami^n$ and the right
vertical arrow comes from the identity of the shtukas
$\Der\cM/\rami^n$ and $\cM/\rami^n$.
\end{lem}

\pf Let $\cM = [ \cM_0 \shtuka{}{} \cM_1 ]$ and let $c \in \uH^1(\Der\cM)$ be a
cohomology class. According to Proposition \ref{explmodexp} there exists a
function $h \in b(F,K)$ such that
\begin{align*}
\Der g(h), g(h) &\in \cM_1, \\
\Der g(h) - g(h) &\in I\cM_1, \\
[\Der g(h)] &= c, \\
\exp [\Der g(h)] &= [g(h)].
\end{align*}
Here the brackets $[\,]$ denote cohomology classes.
We get the result since the images of $\Der g(h)$
and $g(h)$ in $\cM_1/\rami^n$ are the same. \quod

\begin{thm}\label{explreg}%
Let $\cM$ be an explicit model with leading parameter
$b$. If $|\alpha_r| b \leqslant B_{\exp}$ then the exponential map $\exp\colon
\uH^1(\Der\cM) \to \uH^1(\cM)$ is the inverse of the regulator
$\rho\colon \uH^1(\cM) \to \uH^1(\Der\cM)$.
\end{thm}


\pf Let $n \geqslant 0$. By Proposition \ref{expltwist} the twist
$\rami^n\cM$ is an explicit model with leading parameter $b|\theta^{-n}|$.
Consider the ideal
\begin{equation*}
I_n = \{ x \in K : |x| \leqslant (|\alpha_r \theta^{-n}| b)^{q-1} \} \subset \cO_K.
\end{equation*}
Since $|\alpha_r| b \leqslant B_{\exp} < 1$ and $\rami = \theta^{-1} \cO_K$
we conclude that $I_n \subset \rami^n$.
The shtuka $(\rami^n\cM)/\rami^n$ is linear
by Proposition \lref{reg}{elltwistvanish}. So Lemma \ref{explslide} implies that
\begin{equation*}
\xymatrix{
\uH^1(\Der\rami^n\cM) \ar[r] \ar[d]^{\exp}& \uH^1(\Der(\rami^n\cM)/\rami^n) \ar[d]^1 \\
\uH^1(\rami^n\cM) \ar[r] & \uH^1((\rami^n\cM)/\rami^n).
}
\end{equation*}
Now Theorem \lref{reg}{ellregcmp} shows that $\exp$ is the inverse of the regulator map.
\quod

\section{Regulators in general}

In this section we let $A$ be an arbitrary coefficient ring
and $K$ a finite product of local fields containing $\Fq$.
As before we fix a Drinfeld $A$-module $E$ over $K$ with motive $M$.
We assume that the action of $A$ on $\Lie_E(K)$ extends to a continuous action of $F$.

We are finally ready to prove the following result:

\begin{thm}\label{regexp}
If $\cM \subset \iHom(M,\,a(F,K))$ is a local model then the diagram
\begin{equation*}
\xymatrix{
\uH^1(\Der\cM) \ar[r] & F \otimes_{\cO_F} \uH^1(\Der\cM) \ar[d]^{\exp}_{\ltviso} \\
\uH^1(\cM) \ar[u]^{\rho}_{\rtviso} \ar[r] & F \otimes_{\cO_F} \uH^1(\cM)
}
\end{equation*}
is commutative.\end{thm}

\pf We split the proof into several steps.

\medbreak
\textbf{Step 1.} \emph{Let $\cN \subset \cM$ be an inclusion of local models.
The theorem holds for $\cN$ if and only if it holds for $\cM$.}
Indeed Proposition \lref{locmodcoh}{locmodcohiso} shows that the maps
\begin{align*}
F \otimes_{\cO_F} \RGamma(\cN) &\to F\otimes_{\cO_F} \RGamma(\cM), \\
F \otimes_{\cO_F} \RGamma(\cN) &\to F\otimes_{\cO_F} \RGamma(\Der\cM)
\end{align*}
are quasi-isomorphisms. The result follows since both $\rho$ and $\exp$
are natural transformations of functors on the category of local models.

\medbreak
\textbf{Step 2.} \emph{The theorem holds in the case $A = \Fq[t]$.}

Without loss of generality we assume that $K$ is a single local field.
Let $\cM$ be a local model. According to Proposition \ref{explexists} there exists
an explicit model $\cN$ with some parameters $b$, $b_1, \dotsc, b_{r-1}$.
Proposition \lref{locmod}{locmodmap} shows that the local model $\rami^n \cN$
is a subshtuka of $\cM$ for all $n \gg 0$.
By Proposition \ref{expltwist} the twist $\rami^n \cN$ is an explicit model with leading
parameter $b|\theta^{-n}|$. So taking $n \gg 0$ we can ensure that
$|\alpha_r| b \leqslant B_{\exp}$. In this situation Theorem \ref{explreg} shows
that the result holds for $\rami^n\cN$. Step 1 implies that it holds
for $\cM$ as well.

\medbreak
Next let us fix a nonconstant element $a \in A$. We denote $A' = \Fq[a]$ and
$F'$ the local field of $A'$ at infinity and $\rami' \subset \cO_K$ the ideal
generated by $\fm_{F'}$ under the map $F' \to K$ given by the action of $A'$
on $\Lie_E(K)$.

\medbreak
\textbf{Step 3.} \emph{If $\cM$ is a local model such that
$\cM(F \otimes \cO_K/\rami')$ is linear then the theorem holds for $\cM$.}

Let $\cM'$ be $\cM$ viewed as an $\cO_{F'}\indcot\cO_K$-module shtuka.
By Proposition \lref{cclocmod}{cclocmod} the shtuka $\cM'$ is a local model of the
Drinfeld $A'$-module $E$. 
Proposition \lref{cclocmod}{ccreg} shows that the regulator of $\cM'$ is compatible
with the regulator of $\cM$ while Proposition \lref{cclocmod}{ccexp}
does the same for the exponential map. Hence the theorem holds for $\cM$ with
the coefficient ring $A$ if and only if it holds for $\cM'$ with the coefficient
ring $A'$. Applying Step 2 to $\cM'$ we get the result.

\medbreak
\textbf{Step 4.} \emph{The theorem holds for an arbitrary local model.}

Let $\cM$ be a local model.
The ramification ideal $\rami'$ is a power of $\rami$ by construction.
Hence Proposition \lref{locmod}{locmodtwist} demonstrates that
$\rami'\cM$ is a local model. 
Proposition \lref{reg}{elltwistvanish} shows that
$(\rami'\cM)(F \otimes \cO_K/\rami')$ is linear. Therefore the theorem
holds for $\rami'\cM$ by Step 3. As $\rami'\cM \subset \cM$ Step 1 implies
that the theorem holds for $\cM$. \quod

\chapter{Global models and the class number formula}
\label{chapter:cnf}

In this chapter we introduce global shtuka models
of a Drinfeld module $E$ over a Dedekind domain of finite type over $\Fq$ 
and
use them to derive the class number formula
for Drinfeld modules.

\begin{nrmk*}%
The notation of this chapter differs from the notation used
in the introduction. Namely we write $F$ in place of $F_\infty$
and $K$ in place of $K_\infty$.\end{nrmk*}

\section{The setting}

Fix a coefficient ring $A$. As before we denote $F$ the local field of $A$ at
infinity, $\cO_F\subset F$ the ring of integers and $\fm_F \subset\cO_F$ the
maximal ideal. We denote $C$ the compactification of $\Spec A$
and $\Omega_C$ the sheaf of K\"ahler differentials of $C$ over $\Fq$.
We use the following notation:
\begin{equation*}
\begin{array}{l@{\,\,}c@{\,\,}l@{\quad\quad}l@{\,\,}c@{\,\,}l}
\nSpA &=& \Spec A & \omega &=& \Gamma(\nSpA, \,\Omega_C), \\
\nsOF &=& \Spec\cO_F & \omega_{\cO_F} &=& \Gamma(\nsOF,\,\Omega_C), \\
\nsF &=& \Spec F & \omega_F &=& \Gamma(\nsF, \,\Omega_C).
\end{array}
\end{equation*}

Let $R$ be a Dedekind domain of finite type over $\Fq$.
We denote $Y$ the spectrum of $R$ and $X$ the projective curve over $\Fq$ which
compactifies $Y$. 

Let $E$ be a Drinfeld $A$-module over $R$.
As usual $M = \Hom(E,\bG_a)$ stands for the motive of $E$.
The action of $A$ on $\Lie_E$ determines a homomorphism
$A \to R$. We assume that it is \emph{finite flat}.
%
%
%

%

The $F$-algebra $K = R \otimes_A F$ is a finite product of local fields.
Note that it is not necessarily \'etale.
As before we denote $\cO_K \subset K$ the ring of
integers and $\fm_K \subset \cO_K$ the Jacobson radical.
Let $\rami$ be the ideal generated by
$\fm_F$ in $\cO_K$ under the natural map $F \to K$. We call $\rami$ the
\emph{ramification ideal}.\index{idx}{ramification ideal!global model@of a global model}\index{nidx}{$\rami$, ramification ideal}

\section{The notion of a global model}


To the Drinfeld module $E$ with motive $M$ we associate the
$A\otimes R$-module shtuka $\iHom(M,\,\omega\otimes R)$.
In order to prove the class number formula for $E$ we will extend
it to a locally free shtuka on $C \times X$. In this section we
introduce the appropriate notion of such an extension.

Let $\pi\colon C \times Y \to C$ 
be the projection map.
By construction $\pi^\ast \Omega_C$ 
is the sheaf of K\"ahler differentials of $C\times Y$
over $Y$.
We introduce the shtuka
\begin{equation*}
\Omega_{C,Y} = 
\Big[ \pi^\ast \Omega_C \shtuka{\,\,1\,\,}{j} \pi^\ast \Omega_C \Big]
\end{equation*}
where $j$ is the $\tau$-adjoint of the natural isomorphism
$\tau^\ast \pi^\ast \Omega_C \cong \pi^\ast \Omega_C$
induced by the equality $\pi \circ \tau = \pi$.

\begin{dfn}%
We set
$\cM^c = \siHom(\cE, \,\Omega_{C,Y})$
where $\siHom$ is the sheaf
Hom shtuka of Definition~\ref{defglobhomsht}
and $\cE$ is the shtuka on $C\times Y$ constructed
in Theorem~\ref{motivec}.\end{dfn}

%
%
%


\begin{lem}%
$\cM^c(A \otimes R) = \iHom(M, \,\omega \otimes R)$.\end{lem}

\pf
Indeed we have
\begin{equation*}
\cE(A\otimes R) = \Big[ M \shtuka{1}{\tau} M \Big], \quad
\Omega_{C,Y}(A\otimes R) = \Big[ \omega \otimes R \shtuka{1}{\tau} \omega \otimes R \Big]
\end{equation*}
by construction.\quod

\begin{dfn}\label{defglobmod}\index{idx}{shtuka model!global}%
A \emph{global model} of $\iHom(M,\,\omega\otimes R)$ is
a shtuka on $C\times X$ 
which has the following properties:
\begin{enumerate}
\item $\cM$ is a locally free shtuka
extending $\cM^c$ on $C\times Y$.

\item $\cM(A \otimes \cO_K/\fm_K)$ is nilpotent and
$\cM(A \otimes \cO_K/\rami)$ is linear.%
\end{enumerate}
\end{dfn}


\begin{prp}\label{globmodloc}%
If $\cM$ is a global model then $\cM(\cO_F\indcot\cO_K)$ is a local model in
the sense of Definition \lref{locmod}{deflocmod}.\end{prp}

\pf
Corollary~\ref{funclat} implies that $\omega\otimes R$ is
an $A\otimes R$-lattice in the $F\indcot K$-module $a(F,K)$.
Hence $\cM(F\indcot K) = \iHom(M,\,a(F,K))$.
The shtuka $\cM(\cO_F\indcot\cO_K)$ is a locally free
$\cO_F\indcot\cO_K$-lattice in 
$\cM(F\indcot K)$. 
To prove that
$\cM(\cO_F\indcot\cO_K)$ is a local model
we need to show the following:
\begin{enumerate}
\item $\cM(\cO_F\indcot K)$ is the coefficient compactification
in the sense of Definition~\lref{locmod}{locmoddrincomp}.

\item $\cM(F \otimes \cO_K/\fm_K)$ is nilpotent and
$\cM(F \otimes \cO_K/\rami)$ is linear.
\end{enumerate}
Corollary~\ref{funclat} implies that
$\omega_{\cO_F}\otimes R$ is an $\cO_F\otimes R$-lattice
in the $\cO_F\indcot K$-module $a(F/\cO_F, K)$.
We conclude that
$\Omega_{C,Y}(\cO_F\indcot K) = a(F/\cO_F, K)$.
Hence
\begin{equation*}
\cM(\cO_F\indcot K) = \iHom_{\cO_F\indcot K}\big(\cE(\cO_F\indcot K), \,a(F/\cO_F,K)\big)
\end{equation*}
is the coefficient compactification as in (1).
%
%
The property (2) follows directly from the definition of a global model.
\quod

\begin{prp}\label{globmodell}%
If $\cM$ is a global model then the restriction of $\cM$ to
$\nsOFX$ is an elliptic shtuka of ramification ideal $\rami$ in the sense of
Definition~\lref{globell}{defglobell}.\end{prp}

\pf We need to prove the following:
\begin{enumerate}
\item $\cM(\cO_F/\fm_F \otimes R)$ is nilpotent.

\item $\cM(\cO_F\indcot\cO_K)$ is an elliptic shtuka of ramification ideal $\rami$.
\end{enumerate}
By Theorem~\lref{drinconstr}{motivec} the shtuka $\cE(\cO_F/\fm_F \otimes R)$
is co-nilpotent.
So (1) follows from Proposition~\lref{conilp}{conilp}.
The shtuka $\cM(\cO_F\indcot\cO_K)$ is a local model
by Proposition~\ref{globmodloc}
so (2) is a corollary of
Theorem~\lref{locmod}{locmodell}. \quod

\section{Existence of global models}

\begin{lem}\label{globmodbase}%
There exists a shtuka $\cM$ on $\nSpA\times X$ with the following
properties:
\begin{enumerate}
\item $\cM$ is a locally free shtuka extending
$\iHom(M,\,\omega\otimes R)$ on $\nSpA\times Y$.


\item $\cM(A \otimes \cO_K/\fm_K)$ is nilpotent and
$\cM(A\otimes \cO_K/\rami)$ is linear.
\end{enumerate}
\end{lem}

\pf By Theorem~\lref{locmod}{basecompexist} the
shtuka $\iHom(M,\,a(F/A,K))$
admits an
$A \otimes \cO_K$-lattice $\cM^b$ 
with the following properties:
\begin{enumerate}
\item $\cM^b$ is locally free.

\item $\cM^b(A\otimes\cO_K/\fm_K)$ is nilpotent and
$\cM^b(A\otimes \cO_K/\rami)$ is linear.
\end{enumerate}
Using the Beauville-Laszlo Theorem [\stacks{0BP2}] we glue
the $A\otimes R$-module shtuka
$\iHom(M,\,\omega\otimes R)$ to
the $A\complot\cO_K$-module shtuka $\cM^b(A\complot\cO_K)$ 
over
$A\complot K$ and obtain the desired shtuka $\cM$ on
$\nSpA\times X$.\quod

\breakflow
Let $U$ denote the union of open subschemes
$C \times Y$ and $\nSpA \times X$ of $C\times X$
and let $\iota\colon U \to C\times X$ be the open immersion.
Note that the complement of $U$ in $C \times X$ consists of finitely many
points.

\begin{lem}\label{freeext}%
If $\cF$ is a locally free sheaf on $U$
then $\iota_\ast\cF$ is a locally free sheaf on $C \times X$.\end{lem}

\pf Consider the cartesian square of schemes
\begin{equation*}
\xymatrix{
V \ar[r]^{\kappa\quad\quad\quad} \ar[d]_{g} & \Spec\cO_F\indcot\cO_K \ar[d]^f \\
U \ar[r]^{\iota} & C \times X
}
\end{equation*}
where $V$ is the complement of the closed points of
$\Spec\cO_F\indcot\cO_K$
and $g$ is the natural morphism.
Observe that
$f^\ast \iota_\ast \cF = \kappa_\ast g^\ast \cF$.
Now Lemmas \lref{refl}{reflext} and \lref{refl}{reflfree} imply that
$\kappa_\ast g^\ast \cF$ is a locally free sheaf
on $\Spec\cO_F\indcot\cO_K$. Whence
$\iota_\ast\cF$ is locally free.\quod

\begin{prp}\label{globmodexist}%
The shtuka $\iHom(M,\,\omega\otimes R)$ admits a global model.\end{prp}

\pf Let $\cM^b$ be a shtuka constructed in Lemma~\ref{globmodbase}.
The shtukas $\cM^c$ and $\cM^b$ restrict to $\iHom(M,\,\omega\otimes R)$
on $\nSpA\times Y$. We thus obtain a shtuka $\cM$ on
the union $U$ of $C\times Y$ and $\nSpA\times X$.
By Lemma \ref{freeext} the shtuka
$\iota_\ast\cM$ on $C \times X$ is locally free. It is therefore a global model. \quod

\section{Cohomology of the Hom shtukas}

\begin{dfn} The \emph{complex of units} of the Drinfeld module $E$ is the $A$-module complex
\begin{equation*}
U_E =
\Big[ \Lie_E(K) \xrightarrow{\quad\exp\quad}
\frac{E(K)}{E(R)} \,\Big]
\end{equation*}
where $\exp\colon \Lie_E(K) \to E(K)$ is the exponential map.
\end{dfn}


\afterall
Our goal is to construct quasi-isomorphisms
\begin{align*}
\RGammac(\iHom(M,\,\omega\otimes R)) &\xrightarrow{\,\isosign\,} U_E[-1], \\
\RGammac(\Der\iHom(M,\,\omega\otimes R)) &\xrightarrow{\,\isosign\,} \Lie_E(R)[-1]
\end{align*}
where $\RGammac$ is the compactly supported cohomology of shtukas on
$A \otimes R$ (Definition~\lref{csc}{defcsc}).
In the next section we will use these quasi-\hspace{0pt}isomorphisms
to study the cohomology of global models.

We first study the germ cohomology
of shtukas on $A\otimes K$.
To improve legibility we will write
$\RGammag(-)$ in place of $\RGammag(A \otimes K,-)$.
The module $\omega\otimes K$ is an $A\otimes K$-lattice
in the $A\complot K$-module $\omega\complot K$.
Moreover the sequence
\begin{equation*}
0 \to \omega\otimes M \xrightarrow{\phantom{\,\,[\Res]\,\,}} \omega\complot M \xrightarrow{\,\,[\Res]\,\,} g(F,M) \to 0
\end{equation*}
is exact by Lemma~\ref{omegagermsesexact}. We thus get a natural
identification
\begin{equation*}
g(F,K) = \tfrac{(A\complot K)\otimes_{A\otimes K} (\omega\otimes K)}{\omega\otimes K}.
\end{equation*}
As a result Proposition~\ref{locfreegermcoh} provides us with
natural quasi-\hspace{0pt}isomorphisms
\begin{align*}
\RGammag(\iHom(M,\,\omega\otimes K)) &\xrightarrow{\isosign} \RGamma(\iHom(M,\,g(F,K)))[-1], \\
\RGammag(\Der\iHom(M,\,\omega\otimes K)) &\xrightarrow{\isosign} \RGamma(\Der\iHom(M,\,g(F,K)))[-1].
\end{align*}

\begin{dfn}\label{ihomfalie}%
We define natural quasi-\hspace{0pt}isomorphisms
\begin{align*}
\RGammag(\iHom(M,\,\omega\otimes K)) &\xrightarrow{\,\isosign\,} \Lie_E(K)[-1], \\
\RGammag(\Der\iHom(M,\,\omega\otimes K)) &\xrightarrow{\,\isosign\,} \Lie_E(K)[-1]
\end{align*}
as the compositions
\begin{equation*}
\xymatrix{
\RGammag(\iHom(M,\,\omega\otimes K)) \ar[d]_{\textup{Proposition }\ref{locfreegermcoh}}^{\rtviso} &
\RGammag(\Der\iHom(M,\,\omega\otimes K)) \ar[d]^{\textup{Proposition }\ref{locfreegermcoh}}_{\ltviso} \\
\RGamma(\iHom(M,\,g(F,K)))[-1] \ar[d]_{\textrm{Proposition }\lref{ihomcoh}{ihomge}\,}^{\rtviso} &
\RGamma(\Der\iHom(M,\,g(F,K)))[-1] \ar[d]^{\textrm{Proposition }\lref{ihomcoh}{ihomglie}}_{\ltviso} \\
\Lie_E(K)[-1] & \Lie_E(K)[-1]
}
\end{equation*}
\end{dfn}

\breakflow
Let $Q = K/R$. By Corollary~\ref{bnqexact} the natural sequence
\begin{equation*}
0 \to \omega\otimes R \to \omega\complot K \to \omega\complot Q \to 0
\end{equation*}
is exact.
Moreover $\omega\otimes R$ is an $A\otimes R$-lattice
in the $A\complot K$-module $\omega\complot K$.
We thus get a natural identification
\begin{equation*}
\omega\complot Q = \tfrac{(A\complot K)\otimes_{A\otimes R} (\omega\otimes R)}{\omega\otimes R}.
\end{equation*}
As a result Proposition~\ref{locfreecsc} provides us with
natural quasi-\hspace{0pt}isomorphisms
\begin{align*}
\RGammac(\iHom(M,\,\omega\otimes R)) &\xrightarrow{\isosign} \RGamma(\iHom(M,\,\omega\complot Q))[-1], \\
\RGammac(\Der\iHom(M,\,\omega\otimes R)) &\xrightarrow{\isosign} \RGamma(\Der\iHom(M,\,\omega\complot Q))[-1].
\end{align*}

\begin{dfn}\label{ihomcsc}%
We define natural quasi-\hspace{0pt}isomorphisms
\begin{align*}
\RGammac(\iHom(M,\,\omega\otimes R)) &\xrightarrow{\,\isosign\,} U_E[-1], \\
\RGammac(\Der\iHom(M,\,\omega\otimes R)) &\xrightarrow{\,\isosign\,} \Lie_E(R)[-1]
\end{align*}
as the compositions
\begin{equation*}
\xymatrix{
\RGammac(\iHom(M,\,\omega\otimes R)) \ar[d]_{\textup{Proposition }\ref{locfreecsc}}^{\rtviso} &
\RGammac(\Der\iHom(M,\,\omega\otimes R)) \ar[d]^{\textup{Proposition }\ref{locfreecsc}}_{\ltviso} \\
\RGamma(\iHom(M,\,\omega\complot Q))[-1] \ar[d]_{\textrm{Theorem }\ref{ihomunit}\,}^{\rtviso} &
\RGamma(\Der\iHom(M,\,\omega\complot Q))[-1] \ar[d]^{\textrm{Theorem }\ref{derihomunit}}_{\ltviso} \\
U_E[-1] & \Lie_E(R)[-1]
}
\end{equation*}%
Theorems \ref{ihomunit} and \ref{derihomunit} imply that these quasi-\hspace{0pt}isomorphisms
are natural in $E$.\end{dfn}

\breakflow
For every quasi-\hspace{0pt}coherent shtuka $\cM$ on $A\otimes R$
we have a natural morphism
$\RGammac(A \otimes R, \,\cM) \to \RGammag(A \otimes K, \,\cM)$
(see Definition~\ref{csctogerm}).

\begin{prp}\label{ihomcsccomm}%
The square
\begin{equation*}
\xymatrix{
\RGammac(\iHom(M,\,\omega\otimes R)) \ar[rr] \ar[d]_{\textup{Definition }\ref{ihomcsc}}^{\rtviso} &&
\RGammag(\iHom(M,\,\omega\otimes K)) \ar[d]^{\textup{Definition }\ref{ihomfalie}}_{\ltviso} \\
U_E[-1] \ar[rr]^{\textup{identity in degree }1} && \Lie_E(K)[-1]
}
\end{equation*}
is commutative.\end{prp}

\pf For every locally free $A\otimes R$-module shtuka $\cM$ the square
\begin{equation*}
\xymatrix{
\RGammac(A\otimes R,\,\cM) \ar[rr] \ar[d]_{\textup{Proposition }\ref{locfreecsc}}^{\rtviso} &&
\RGammag(A\otimes K,\,\cM) \ar[d]^{\textup{Proposition }\ref{locfreegermcoh}}_{\ltviso} \\
\RGamma\!\Big(A\otimes R,\,\tfrac{\cM(A\complot K)}{\cM(A\otimes R)}\Big)[-1] \ar[rr] &&
\RGamma\!\Big(A\otimes K,\,\tfrac{\cM(A\complot K)}{\cM(A\otimes K)}\Big)[-1]
}
\end{equation*}
is commutative by construction. So the result follows from
Theorem~\ref{ihomunit}~(3).\quod

\begin{prp}\label{derihomcsccomm}%
The square
\begin{equation*}
\xymatrix{
\RGammac(\Der\iHom(M,\,\omega\otimes R)) \ar[rr] \ar[d]_{\textup{Definition }\ref{ihomcsc}}^{\rtviso} &&
\RGammag(\Der\iHom(M,\,\omega\otimes K)) \ar[d]^{\textup{Definition }\ref{ihomfalie}}_{\ltviso} \\
\Lie_E(R)[-1] \ar[rr]^{\textup{embedding in degree }1} && \Lie_E(K)[-1]
}
\end{equation*}
is commutative.\end{prp}

\pf Same as the proof of Proposition~\ref{ihomcsccomm}.\quod

\section{Cohomology of global models}

We work with a fixed global model $\cM$.

\begin{lem}\label{globmodcsctocoh}%
The natural maps
\begin{align*}
\RGammac(\nSpA\times Y,\,\cM) &\to \RGamma(\nSpA\times X,\,\cM), \\
\RGammac(\nSpA\times Y,\,\Der\cM) &\to \RGamma(\nSpA\times X,\,\Der\cM)
\end{align*}
are quasi-\hspace{0pt}isomorphisms.\end{lem}

\pf
The shtuka
$\cM(A \otimes \cO_K/\fm_K)$ is nilpotent so 
the claim for $\RGammac(\nSpA\times Y,\,\cM)$
follows by Proposition~\lref{csc}{bcohmap}.
The same applies to $\Der\cM$.\quod

\begin{dfn}\label{defglobmodcoh}%
We define natural quasi-\hspace{0pt}isomorphisms
\begin{align*}
\RGamma(\nSpA\times X,\,\cM) &\xrightarrow{\,\isosign\,} U_E[-1], \\
\RGamma(\nSpA\times X,\,\Der\cM) &\xrightarrow{\,\isosign\,} \Lie_E(R)[-1]
\end{align*}
as the compositions
\begin{equation*}
\xymatrix{
\RGamma(\nSpA\times X,\,\cM) &
\RGamma(\nSpA\times X,\,\Der\cM) \\
\RGammac(\nSpA\times Y,\,\cM)
\ar[u]^{\textup{Proposition }\ref{bcohmap}}_{\rtviso}
\ar[d]_{\textrm{Definition }\ref{ihomcsc}\,}^{\rtviso} &
\RGammag(\nSpA\times Y,\,\Der\cM)
\ar[u]_{\textup{Proposition }\ref{bcohmap}}^{\ltviso}
\ar[d]^{\textrm{Definition }\ref{ihomcsc}}_{\ltviso} \\
U_E[-1] & \Lie_E(R)[-1]
}
\end{equation*}%
These quasi-\hspace{0pt}isomorphisms are natural in $\cM$
and $E$ by construction.%
\end{dfn}


\breakflow
According to Proposition \ref{globmodloc}
the shtuka $\cM(\cO_F\indcot\cO_K)$ is a local model.
So Lemma \lref{locmod}{locqi} shows that the maps of
Definition \lref{locmod}{gammamaps} induce $F$-linear
quasi-isomorphisms
\begin{align*}
\gamma&\colon \RGamma(F\indcot\cO_K,\,\cM) \xrightarrow{\,\isosign\,} \Lie_E(K)[-1], \\
\Der\gamma&\colon \RGamma(F\indcot\cO_K,\,\Der\cM) \xrightarrow{\,\isosign\,} \Lie_E(K)[-1]
\end{align*}

\begin{prp}\label{globmodcoh}%
The square
\begin{equation*}
\xymatrix{
\RGamma(\nSpA\times X, \,\cM) \ar[rr]^{\textup{pullback}} \ar[d]^{\rtviso}_{\textup{Definition }\ref{defglobmodcoh}} &&
\RGamma(F \indcot \cO_K, \,\cM) \ar[d]^{\gamma}_{\ltviso} \\
U_E[-1] \ar[rr]^{\textup{identity in degree }1} && \Lie_E(K)[-1]
}
\end{equation*}
is commutative.%
\end{prp}

\pf
We verify it in several steps.

\medbreak
\textbf{Step 1.}
Consider the square
\begin{equation*}
\xymatrix{
U_E[-1] \ar[rr]^{\textup{identity in degree }1} && \Lie_E(K)[-1] \\
\RGamma(\nSpA\times X,\,\cM) \ar[u]^{\textup{Definition }\ref{defglobmodcoh}}_{\rtviso} \ar[rr]^{\textup{global}} &&
\RGammag(A \otimes K,\,\cM) \ar[u]^{\ltviso}_{\textup{Definition }\ref{ihomfalie}}
}
\end{equation*}
where
the arrow labelled ``global'' is the global germ
map
\begin{equation*}
\RGamma(A \times X, \cM) \xleftarrow{\,\isosign\,}
\RGammac(A \otimes R, \cM) \to \RGammag(A \otimes K,\cM)
\end{equation*}
of Definition~\lref{csc}{defgcmp}.
Proposition~\ref{ihomcsccomm} implies that this square is commutative.

\medbreak
\textbf{Step 2.}
Consider the diagram
\begin{equation*}
\xymatrix{
\RGamma(\nSpA\times X,\,\cM) \ar[d] \ar[r]^{\textup{global}} &
\RGammag(A\otimes K,\,\cM) \ar[d] \\
\RGamma(\nsF\times X,\,\cM) \ar[r]^{\textup{global}} \ar[d] &
\RGammag(F \indcot K,\,\cM) \\
\RGamma(F \indcot \cO_K,\,\cM) \ar[r]^{\textup{local}} &
\RGammag(F \indcot K,\,\cM) \ar@{=}[u]
}
\end{equation*}
where the arrow labelled ``local'' is the local germ map
of Definition~\lref{lcmp}{defloccmp} and the unlabelled arrows are
the pullback morphisms. The top square of this diagram commutes by naturality of
the global germ map. The commutativity of the bottom square follows from
Theorem~\lref{locglobcmp}{locglobcmp}.

\medbreak
\textbf{Step 3.}
Combining Step 1 and Step 2 we get a commutative diagram
\begin{equation*}
\xymatrix{
U_E[-1] \ar[rr]^{\textup{identity in degree }1} && \Lie_E(K)[-1] \\
\RGamma(\nSpA\times X,\,\cM) \ar[u]^{\textup{Definition }\ref{defglobmodcoh}}_{\rtviso} \ar[d] \ar[rr]^{\textup{global}} &&
\RGammag(A \otimes K,\,\cM) \ar[u]^{\ltviso}_{\textup{Definition }\ref{ihomfalie}} \ar[d] \\
\RGamma(F \indcot \cO_K,\,\cM) \ar[rr]^{\textup{local}} &&
\RGammag(F \indcot K,\,\cM).
}
\end{equation*}

\medbreak
\textbf{Step 4.}
We claim that the arrow
$\RGammag(A \otimes K,\cM) \to \RGammag(F\indcot K,\cM)$ is a quasi-isomorphism
and that the
square
\begin{equation*}
\xymatrix{
\RGamma(F \indcot \cO_K,\,\cM)\ar[d]_{\phantom{Tralalalala}\textup{local}} \ar[rr]^{\gamma} &&
\Lie_E(K)[-1] \\
\RGammag(F\indcot K,\,\cM) &&
\RGammag(A\otimes K,\,\cM) \ar[u]_{\textup{Definition }\ref{ihomfalie}}^{\ltviso} \ar[ll]_{\isosign}
}
\end{equation*}
is commutative.
Together with Step 3 this claim immediately implies the
theorem.

To prove this claim we first observe that the square
\begin{equation*}
\xymatrix{
\RGammag(A \otimes K,\,\cM) \ar[rr]_{\bisosign}^{\textup{Def. }\ref{ihomfalie}}
\ar[d] && \Lie_E(K)[-1] \ar@{=}[d] \\
\RGammag(F \indcot K,\,\cM) \ar[rr]_{\bisosign}^{\textup{Def. }\lref{locmod}{ihomalie}}
&& \Lie_E(K)[-1]
}
\end{equation*}
is commutative by construction.
By definition the quasi-\hspace{0pt}isomorphism $\gamma$
is the composition
\begin{equation*}
\RGamma(F\indcot\cO_K,\,\cM) \xrightarrow{\,\textup{local}\,}
\RGammag(F\indcot K,\,\cM) \xrightarrow{\,\textup{Def. }\lref{locmod}{ihomalie}\,} \Lie_E(K)[-1].
\end{equation*}
Hence the claim and the proposition follow. \quod

\begin{prp}\label{globmodliecoh}%
The square
\begin{equation*}
\xymatrix{
\RGamma(\nSpA\times X, \,\Der\cM) \ar[rr]^{\textup{pullback}} \ar[d]^{\rtviso}_{\textup{Definition }\ref{defglobmodcoh}} &&
\RGamma(F \indcot \cO_K, \,\Der\cM) \ar[d]^{\Der\gamma}_{\ltviso} \\
\Lie_E(R)[-1] \ar[rr]^{\textup{embedding in degree }1} && \Lie_E(K)[-1]
}
\end{equation*}
is commutative.%
\end{prp}

\pf Same as the proof of Proposition~\ref{globmodcoh}. \quod


\section{Regulators}
\label{sec:globmodreg}

\begin{dfn}\index{idx}{regulator!arithmetic}The \emph{arithmetic regulator}
\begin{equation*}
\rho_E\colon F \otimes_A U_E \to \Lie_E(K)[0]
\end{equation*}
of 
$E$
is the $F$-linear extension of the morphism
$U_E \to \Lie_E(K)[0]$
given by the identity in
degree zero. 
\end{dfn}

\breakflow
Taelman \cite{ltut} demonstrated that $\rho_E$ is a quasi-isomorphism.
This will also follow from the results below.
 
We work with a fixed global model $\cM$.
By Proposition \ref{globmodell} the restriction of $\cM$
to $\nsOFX$ is an elliptic shtuka of ramification ideal $\rami$.
We can thus make the following definition.

\begin{dfn}\label{defglobmodreg}\index{idx}{regulator!global model@of a global model}%
The \emph{regulator}
\begin{equation*}
\rho\colon \RGamma(\nsFX,\,\cM) \xrightarrow{\,\,\isosign\,\,} \RGamma(\nsFX,\,\Der\cM)
\end{equation*}
of a global model $\cM$ is the $F$-linear extension of the regulator
$\rho\colon \RGamma(\nsOFX,\,\cM) \xrightarrow{\isosign} \RGamma(\nsOFX,\,\Der\cM)$
of the elliptic shtuka $\cM|_{\nsOFX}$.%
\end{dfn}


\breakflow
The $F$-linear extensions of the maps of
Definition~\ref{defglobmodcoh} give us the quasi-\hspace{0pt}isomorphisms
\begin{align}
\label{fglobmodcoh}
\RGamma(\nsFX,\,\cM) &\xrightarrow{\isosign} F\otimes_A U_E[-1], \\
\label{fderglobmodcoh}
\RGamma(\nsFX,\,\Der\cM) &\xrightarrow{\isosign} \Lie_E(K)[-1].
\end{align}

\begin{thm}\label{globmodreg}%
The square
\begin{equation*}
\xymatrix{
\RGamma(\nsFX,\,\cM) \ar[r]^{\rho} \ar[d]^{\rtviso}_{\eqref{fglobmodcoh}}  &
\RGamma(\nsFX,\,\Der\cM) \ar[d]_{\ltviso}^{\eqref{fderglobmodcoh}} \\
F \otimes_A U_E[-1] \ar[r]^{\rho_E[-1]} & \Lie_E(K)[-1]
}
\end{equation*}
is commutative.%
\end{thm}

\pf According to Proposition \ref{globmodloc}
the shtuka $\cM(\cO_F\indcot\cO_K)$ is a local model
of $\iHom(M,a(F,K))$. So it is an elliptic shtuka
of ramification ideal $\rami$ by Theorem~\lref{locmod}{locmodell}.
As such it has a regulator
\begin{equation*}
\RGamma(\cO_F\indcot\cO_K,\,\cM) \to
\RGamma(\cO_F\indcot\cO_K,\,\Der\cM).
\end{equation*}
By Lemma~\lref{locmod}{locqi} its $F$-linear extension
can be identified with a quasi-\hspace{0pt}isomorphism
\begin{equation*}
\widecheck{\rho}\colon
\RGamma(F\indcot\cO_K,\,\cM) \to \RGamma(F\indcot\cO_K,\,\Der\cM).
\end{equation*}
Theorem~\ref{globellregexist} shows that the natural square
\begin{equation*}
\xymatrix{
\RGamma(\nsFX,\,\cM) \ar[r]^{\rho} \ar[d] &
\RGamma(\nsFX,\,\Der\cM) \ar[d] \\
\RGamma(F\indcot\cO_K,\,\cM) \ar[r]^{\widecheck{\rho}} &
\RGamma(F\indcot\cO_K,\,\Der\cM)
}
\end{equation*}
is commutative.
The vertical arrows in it are
quasi-\hspace{0pt}isomorphisms by Proposition~\ref{globellcoh}.
Now consider the diagram
\begin{equation*}
\xymatrix{
& F \otimes_A \RGamma(\nSpA\times X,\,\cM) \ar[d] \ar[r]^{\isosign} &
F \otimes_A U_E[-1] \ar[d]^{\rho_E}_{\ltviso} \\
\RGamma(\nsFX,\,\cM) \ar[d]^\rho \ar[r]^{\isosign} & \RGamma(F\indcot\cO_K,\,\cM) \ar[d]_{\widecheck{\rho}} \ar[r]^{\gamma}_{\isosign} &
\Lie_E(K)[-1] \ar@{=}[d] \\
\RGamma(\nsFX,\,\Der\cM) \ar[r]^{\isosign} & \RGamma(F\indcot\cO_K,\,\Der\cM) \ar[r]^{\Der\gamma}_{\isosign} &
\Lie_E(K)[-1] \\
& F \otimes_A \RGamma(\nSpA\times X,\,\Der\cM) \ar[u] \ar[r]^{\isosign} &
F \otimes_A \Lie_E(R)[-1] \ar[u]^{\ltviso}
}
\end{equation*}
The top and bottom squares
are induced by the commutative squares of
Proposition \ref{globmodcoh} and \ref{globmodliecoh} respectively.
The middle square commutes by Theorem \lref{locmodreg}{regexp}.
The square
\begin{equation*}
\xymatrix{
\RGamma(\nsFX,\,\cM) \ar[r]^{\eqref{fglobmodcoh}}_{\isosign} \ar[d]_{\ltviso} & F\otimes_A U_E[-1] \\
\RGamma(F\indcot\cO_K,\,\cM) & F\otimes_A\RGamma(\nSpA\times X,\,\cM) \ar[u]^{\ltviso}_{\textup{Def. }\ref{defglobmodcoh}} \ar[l]_{\isosign}
}
\end{equation*}
commutes by construction.
%
An analogous observation applies to \eqref{fderglobmodcoh}.
We thus get the result. \quod

\section{Euler products}

%

As before we denote $\nfA$ the fraction field of $A$
and $\nOA$ the module of K\"ahler differentials of $\nfA$ over $\Fq$.
Given a prime $\fp$ of $A$ we denote $A_\fp$ the $\fp$-adic completion of $A$.

Let $\fm$ be a maximal ideal of $R$.
We denote $E_\fm$ the pullback of $E$ to $R/\fm$.
For each $\fm$ we fix a separable closure
of the residue field $R/\fm$.
With this choice we have for every prime $\fp$ of $A$
the $A_\fp$-adic Tate module $T_\fp E_\fm$.
Given a maximal ideal $\fm$ of $R$ and a prime $\fp$ of $A$
different from $\mu^{-1}\fm$ we define
\begin{equation*}
P_\fm(T) = \det\nolimits_{A_\fp} \!\big(1 - T \sigma^{-1}_\fm \,\big|\, T_\fp E_\fm\big) \in
A_\fp[T]
\end{equation*}
where $\sigma^{-1}_\fm$ is the geometric Frobenius element at $\fm$.

\begin{prp}%
The characteristic polynomial $P_\fm$ has coefficients in $\nfA$ and is
independent of the choice of $\fp$.\end{prp}

\pf
Theorem~\lref{ihomtate}{ihomdet} implies that
\begin{equation*}
P_\fm(T^d) = 
\det\nolimits_{\nfA} \!\big(1 - T(i^{-1} j) \,\big|\, \Hom_{A\otimes R}(M, \,\nOA \otimes R/\fm)
\big)
\end{equation*}
where $i$ and $j$ are the arrows of the shtuka
$\iHom(M, \,\nOA \otimes R/\fm)$. \quod



\begin{dfn}%
We define a formal product $L(E^*,0) \in F$ as follows:
\begin{equation*}
L(E^*, 0) = \prod_{\fm} P_\fm(1)^{-1}
\end{equation*}
where $\fm \subset R$ ranges over the maximal ideals.
\end{dfn}

\breakflow
In a moment we will see that this product converges.

\begin{lem}If $\cM$ is a global model then $\cM(\cO_F/\fm_F\otimes R)$ is nilpotent.\end{lem}

\pf Indeed the restriction of $\cM$ to $\nsOFX$ is an elliptic shtuka
by Proposition~\ref{globmodell}.\quod

\breakflow
Hence the $L$-invariants of Definitions \ref{euprodloc} and \ref{euprodglob}
make sense for $\cM$.

\begin{prp}\label{globmodeuloc}%
Let $\cM$ be a global model.
For every maximal ideal $\fm$
we have an equality
$L(\cM(\cO_F\otimes R/\fm)) = P_\fm(1)$
in $F$.
\end{prp}

\pf Suppose that $\cM$ is given by the diagram
\begin{equation*}
\cM_0 \shtuka{i}{j} \cM_1.
\end{equation*}
By definition
$L(\cM(\cO_F\otimes R/\fm))$
is the determinant
\begin{equation*}
\det\nolimits_{\cO_F}\!\big(1 - i^{-1} j \,\big|\, \cM_0(\cO_F \otimes R/\fm)\big).
\end{equation*}
As $\cM$ is a global model the shtuka $\cM(\cO_F\otimes R/\fm)$
is a locally free $\cO_F \otimes R/\fm$-lattice in the $F \otimes R/\fm$-module
shtuka $\iHom(M, \,a(F, R/\fm))$. Hence the determinant above is equal to
\begin{equation*}
\det\nolimits_F \!\big(1 - i^{-1} j \,\big|\, \Hom_{A\otimes R}(M, \,a(F, R/\fm))\big).
\end{equation*}
Next $\nOA \otimes R/\fm$ is an $\nfA \otimes R/\fm$-lattice in
the $F\otimes R/\fm$-module $a(F, R/\fm)$ so the latter determinant coincides with
\begin{equation*}
\det\nolimits_{\nfA} \!\big(1 - i^{-1} j \,\big|\, \Hom_{A\otimes R}(M,\,\nOA \otimes R/\fm))\big).
\end{equation*}
This determinant is equal to $P_\fm(1)$ by Theorem~\ref{ihomdet}.\quod

%

\begin{dfn}%
Let $\cM$ be a global model. We define $L(\cM) = L(\cM(\cO_F\otimes R))$.\end{dfn}

\begin{prp}\label{globmodeu}%
The formal product $L(E^*,0)$ has the following properties.
\begin{enumerate}
\item $L(E^*, 0)$ converges to an element of $\cO_F$.

\item $L(E^*, 0) \equiv 1 \pmod{\fm_F}$.

\item For every global model $\cM$ we have
$L(\cM) = L(E^*, 0)$.
\end{enumerate}
\end{prp}

\pf Let $\cM$ be a global model. By definition
\begin{equation*}
L(\cM) = 
L(\cM(\cO_F\otimes R)) =
\prod_\fm L(\cM(\cO_F \otimes R/\fm))^{-1}
\end{equation*}
where $\fm$ runs over all the maximal ideals.
The product defining $L(\cM)$ converges
by Lemma~\lref{globell}{euconv}.
Proposition~\ref{globmodeuloc} now implies that
$L(\cM) = L(E^*, 0)$ and we get (3).
Proposition~\lref{euprod}{euprodoneunit} shows that
$L(\cM) \equiv 1 \pmod{\fm_F}$.
Since global models exist by Proposition~\ref{globmodexist}
we get (1) and (2).\quod

\section{Trace formula}
\label{sec:globmodtrace}

%
%
%
%
%

\begin{dfn}%
Let $\cM$ be a locally free shtuka on $C \times X$ given by a diagram
\begin{equation*}
\cM_0 \shtuka{\,\,i\,\,}{\,\,j\,\,} \cM_1.
\end{equation*}
The \emph{twist} $\rami\cM$ is the shtuka
\begin{equation*}
\cI\cM_0 \shtuka{\,\,i\,\,}{\,\,j\,\,} \cI\cM_1
\end{equation*}
where $\cI \subset \cO_{C\times X}$ is the pullback to $C\times X$
of the unique ideal sheaf $\cI_0 \subset\cO_X$ satisfying
$\cI_0(\Spec R) = R$ and 
$\cI_0(\Spec\cO_K) = \rami$.
\end{dfn}

\begin{prp} If $\cM$ is a global model then $\rami\cM$ is a global model.\end{prp}

\pf The restrictions of $\cM$ and $\rami\cM$ to $C\times Y$
coincide by construction. Moreover the same argument as in
the proof of Proposition~\ref{ellbasetwist}
shows that $\rami\cM(A \otimes \cO_K/\fm_K)$
is nilpotent and $\rami\cM(A \otimes \cO_K/\rami)$ is linear.
Hence $\rami\cM$ is a global model.\quod

\begin{dfn}%
The \emph{$\zeta$-isomorphism} of a global model $\cM$ is the $\zeta$-isomorphism
\begin{equation*}
\zeta\colon \det\nolimits_A \RGamma(\nSpA\times X,\,\cM) \xrightarrow{\isosign} \det\nolimits_A \RGamma(\nSpA\times X,\,\Der\cM)
\end{equation*}
of the restriction of $\cM$ to $\nSpA\times X$.\end{dfn}
 
\begin{thm}\label{globmodtrace}%
Let $\cM$ be a global model.
For all $n\gg 0$ we have
\begin{equation*}
\zeta_{\rami^n\cM} = L(\rami^n\cM) \cdot \det\nolimits_F(\rho_{\rami^n\cM}).
\end{equation*}
\end{thm}

\begin{rmk}This formula should hold for $n = 0$ as well. See the remark
below Theorem \lref{globell}{elltrace}.\end{rmk}

\afterall\noindent
\textit{Proof of Theorem \ref{globmodtrace}.}
Let $\cM_0$ and $\cM_1$ be the underlying sheaves of $\cM$.
The ideal sheaf $\cI$ above is anti-ample relative to $C$.
Hence for all $n \gg 0$ and $*\in\{0,1\}$ we have
\begin{equation*}
\uH^0(\Spec(\cO_F/\fm_F)\times X, \, \cI^n\cM_*) = 0.
\end{equation*}
We can thus apply Theorem~\ref{elltrace}
to the restriction of $\rami^n\cM$ to $\nsOFX$.
The $\zeta$-isomorphisms are compatible with the change
of coefficients by Proposition~\ref{globcoeffchzeta}.
Hence Theorem~\ref{elltrace} implies the result.\quod

\section{The class number formula}


Recall that the complex of units of the Drinfeld module $E$ is 
\begin{equation*}
U_E =
\Big[ \Lie_E(K) \xrightarrow{\quad\exp\quad}
\frac{E(K)}{E(R)} \,\Big].
\end{equation*}
The arithmetic regulator
$\rho_E\colon F \otimes_A U_E \to \Lie_E(K)[0]$
is the $F$-linear extension of
the morphism $U_E \to \Lie_E(K)[0]$ given by the identity in
degree zero. 
Taelman \cite{ltut} observed that $U_E$ is a perfect $A$-module complex
and $\rho_E$ is a quasi-isomorphism.

By construction
$\Lie_E(K) = F \otimes_A \Lie_E(R)$
so the one-dimensional $F$-vector space $\det_F\Lie_E(K)$
contains a canonical $A$-lattice
$\det_A\Lie_E(R)$.

\begin{thm}%
The image of $\det_A U_E$ under 
$\det_F(\rho_E)$ is 
\begin{equation*}
L(E^*,0) \cdot \det\nolimits_A \Lie_E(R).
\end{equation*}
\end{thm}

\pf
Proposition \ref{globmodexist} associates a global shtuka model $\cM$
with the Drinfeld module $E$. The model comes equipped with
a $\zeta$-isomorphism
\begin{equation*}
\zeta\colon
\det\nolimits_A\RGamma(\nSpA\times X,\,\cM) \xrightarrow{\isosign}
\det\nolimits_A\RGamma(\nSpA\times X, \,\Der\cM)
\end{equation*}
and a regulator
$\rho\colon
\RGamma(\nsFX, \,\cM) \xrightarrow{\isosign} \RGamma(\nsFX,\,\Der\cM)$.

Theorem \ref{globmodtrace} shows
that after replacing $\cM$ with a twist $\rami^n\cM$ by a high enough power
of $\rami^n$ we have 
\begin{equation}\label{maintrace}
\zeta = L(\cM) \cdot \det\nolimits_F(\rho).
\end{equation}
%
Definition~\ref{defglobmodcoh}
provides us with quasi-\hspace{0pt}isomorphisms
\begin{align*}
\RGamma(\nSpA\times X,\,\cM) &\xrightarrow{\isosign} U_E[-1], \\
\RGamma(\nSpA\times X,\,\Der\cM) &\xrightarrow{\isosign} \Lie_E(R)[-1].
\end{align*}
Hence 
$\zeta$ induces an $A$-module isomorphism
\begin{equation*}
\det\nolimits_A(U_E[-1]) \xrightarrow{\isosign} \det\nolimits_A(\Lie_E(R)[-1]).
\end{equation*}
Thanks to Theorem~\ref{globmodreg}
we know that under the quasi-\hspace{0pt}isomorphisms above the
regulator $\rho$ matches with the shifted arithmetic regulator
\begin{equation*}
\rho_E[-1]\colon F \otimes_A U_E[-1] \to \Lie_E(K)[-1].
\end{equation*}
Moreover $L(\cM) = L(E^*,0)$ by Proposition~\ref{globmodeu}.
As $\det_F(\rho_E[-1])$ is the inverse of
the dual of $\det_F(\rho_E)$
we conclude that
\eqref{maintrace} implies the theorem.\quod

\breakflow\noindent
{\small
The references of the form [$wxyz$] point to tags in the Stacks
project \cite{stacks} as explained in
the chapter ``\chnotconv''.}

\Printindex{idx}{Index of terms}

\Printindex{nidx}{Index of symbols}

%
%
%

\end{document}